\documentclass[oneside]{book}
\PassOptionsToPackage{dvipsnames,table}{xcolor}
\usepackage[toc,page]{appendix}
\usepackage[utf8]{inputenc}
\usepackage[T1]{fontenc}
\usepackage{tikz-cd}
\usepackage{tikz}
\usepackage{nextpage}
\usepackage{amsfonts}
\usepackage{amsmath}
\usepackage{amssymb}
\usepackage{amsthm}
\usepackage{mathtools}
\usepackage[left=3cm,right=3cm,top=3.5cm,bottom=3.5cm]{geometry}
\usepackage[french]{babel}
\usepackage{enumitem}
\usepackage{comment}
\usepackage[ruled,vlined,french,onelanguage,algochapter]{algorithm2e}
\usepackage{makecell}
\usepackage{mathrsfs}
\usepackage{mathdots}
\usepackage[Sonny]{fncychap}
\usepackage{adjustbox}
\usepackage{thmtools}
\usepackage{textgreek}
\usepackage{textcomp}
\usepackage{float}
\usepackage[backend=bibtex, style=alphabetic,  url = false, isbn = false]{biblatex}
\DeclareFieldFormat{postnote}{#1}
\DeclareFieldFormat{multipostnote}{#1}
\usepackage{hyperref}
\usepackage{numprint}
\usepackage{arydshln}

\hypersetup{
    colorlinks,
    citecolor=blue,
    filecolor=blue,
    linkcolor=blue,
    urlcolor=blue
}
\title{}
\date{}
\DeclareMathOperator{\Tors}{Tors}
\DeclareMathOperator{\Fact}{Fact}
\DeclareMathOperator{\Fib}{Fib}
\DeclareMathOperator{\Sch}{Sch}
\DeclareMathOperator{\Set}{Set}
\DeclareMathOperator{\cd}{cd}
\DeclareMathOperator{\Sw}{Sw}
\DeclareMathOperator{\Fet}{F\acute{e}t}
\DeclareMathOperator{\coind}{coind}
\DeclareMathOperator{\Spec}{Spec}
\DeclareMathOperator{\Proj}{Proj}
\DeclareMathOperator{\Frac}{Frac}
\DeclareMathOperator{\Aut}{Aut}
\DeclareMathOperator{\Pic}{Pic}
\DeclareMathOperator{\End}{End}
\DeclareMathOperator{\rang}{rang}
\DeclareMathOperator{\res}{res}
\DeclareMathOperator{\Div}{Div}

\DeclareMathOperator{\RG}{R\Gamma}

\DeclareMathOperator{\HH}{H}
\DeclareMathOperator{\DD}{D}
\DeclareMathOperator{\colim}{colim}
\DeclareMathOperator{\id}{id}
\DeclareMathOperator{\ind}{ind}
\DeclareMathOperator{\Mod}{Mod}
\DeclareMathOperator{\im}{im}
\DeclareMathOperator{\GL}{GL}
\DeclareMathOperator{\Isom}{{Isom}}
\DeclareMathOperator{\Hom}{Hom}
\DeclareMathOperator{\Ext}{Ext}
\DeclareMathOperator{\coker}{coker}
\DeclareMathOperator{\RHom}{RHom}

\DeclareMathOperator{\Ab}{Ab}
\DeclareMathOperator{\PAb}{PAb}
\DeclareMathOperator{\rg}{rg}
\DeclareMathOperator{\Mat}{Mat}
\DeclareMathOperator{\Gal}{Gal}
\DeclareMathOperator{\tr}{tr}
\DeclareMathOperator{\Couv}{Couv}

\DeclareMathOperator{\TC}{TC}
\DeclareMathOperator{\Frob}{Frob}
\DeclareMathOperator{\Tot}{Tot}
\DeclareMathOperator{\ppcm}{ppcm}
\DeclareMathOperator{\pgcd}{pgcd}
\DeclareMathOperator{\codim}{codim}

\newcommand{\pf}{{\rm pf}}

\newcommand{\prol}{\pi_1^{{\rm pro-}\ell}}
\newcommand{\C}{\mathcal{C}}

\newcommand{\op}{{\rm op}}
\newcommand{\R}{{\rm R}}
\newcommand{\K}{\mathscr{K}}

\newcommand{\m}{\mathfrak{m}}
\newcommand{\NN}{\mathbb{N}}
\newcommand{\FF}{\mathbb{F}}
\newcommand{\CC}{\mathbb{C}}
\newcommand{\ZZ}{\mathbb{Z}}
\newcommand{\QQ}{\mathbb{Q}}
\newcommand{\PP}{\mathbb{P}}
\newcommand{\bareta}{{\bar{\eta}}}
\newcommand{\OO}{\mathcal{O}}
\newcommand{\LL}{\mathcal{L}}
\newcommand{\GG}{\mathbb{G}}
\newcommand{\ttt}{{\mathrm{t}}}

\newcommand{\F}{\mathcal{F}}
\newcommand{\G}{\mathcal{G}}
\newcommand{\I}{\mathcal{I}}
\newcommand{\T}{\mathcal{T}}

\newcommand{\OX}{\OO_X}
\newcommand{\OT}{\OO_T}
\newcommand{\OP}{\OO_{\PP^1}}
\newcommand{\A}{\mathbb{A}}
\newcommand{\tX}{{\tilde{X}}}
\newcommand{\tY}{{\tilde{Y}}}
\newcommand{\tZ}{{\tilde{Z}}}
\newcommand{\tT}{{\tilde{T}}}
\newcommand{\tL}{{\tilde{\mathcal{L}}}}
\newcommand{\tF}{{\tilde \F}}

\newcommand{\Gk}{\Gal(k |k_0)}
\newcommand\commentaire[1]{\marginpar{\textcolor{Magenta}{#1}}}
\newcommand{\sep}{{\rm sep}}
\newcommand{\parf}{{\rm parf}}

\newcommand{\cone}{\rm c\hat{o}ne}
\newcommand{\cont}{\mathrm{cont}}
\newcommand{\Xet}{X_{\rm \acute{e}t}}
\newcommand{\Yet}{Y_{\rm \acute{e}t}}
\newcommand{\Xfet}{X_{\rm f\acute{e}t}}
\newcommand{\Xlet}{X_{\rm \ell\acute{e}t}}
\newcommand{\red}{{\rm red}}
\newcommand{\sing}{{\rm sing}}
\renewcommand{\div}{\mathrm{div}}
\renewcommand{\L}{\mathcal{L}}

\SetKwInput{Comp}{Complexité}
\SetKwInput{Dep}{Algorithmes utilisés}
\SetKwInput{Source}{Référence}

\theoremstyle{definition}

\usepackage{xpatch}
\makeatletter
\xpatchcmd{\@thm}{\thm@headpunct{.}}{\thm@headpunct{}}{}{}
\makeatother

\newtheorem{df}{Définition}[section]
\newtheorem{rk}[df]{Remarque}
\newtheorem{theorem}[df]{Théorème}
\newtheorem{prop}[df]{Proposition}
\newtheorem{lem}[df]{Lemme}
\newtheorem{ex}[df]{Exemple}
\newtheorem{cor}[df]{Corollaire}

\renewcommand{\thedf}{\arabic{chapter}.\arabic{section}.\arabic{df}}
\renewcommand{\therk}{\arabic{chapter}.\arabic{section}.\arabic{rk}}
\renewcommand{\thetheorem}{\arabic{chapter}.\arabic{section}.\arabic{theorem}}
\renewcommand{\theprop}{\arabic{chapter}.\arabic{section}.\arabic{prop}}
\renewcommand{\thelem}{\arabic{chapter}.\arabic{section}.\arabic{lem}}
\renewcommand{\theex}{\arabic{chapter}.\arabic{section}.\arabic{ex}}
\renewcommand{\thecor}{\arabic{chapter}.\arabic{section}.\arabic{cor}}

\newtheorem{iprop}{Proposition}

\newtheorem{itheorem}[iprop]{Théorème}

\tikzstyle{titre} = [minimum width=5.8cm, minimum height=2cm,text centered, draw=black, fill=white!30]
\tikzstyle{titrevert} = [minimum width=2.2cm, minimum height=5.5cm,text centered, draw=black, fill=white!30]
\tikzstyle{bon} = [rectangle, rounded corners, minimum width=2.3cm, minimum height=1cm,text centered, draw=black]
\tikzstyle{final} = [rectangle, rounded corners, minimum width=4.6cm, minimum height=1cm,text centered, draw=black]
\tikzstyle{addbon} = [rectangle, rounded corners, minimum width=2.4cm, minimum height=0.5cm,text centered, draw=black]
\tikzstyle{auteur}=[rectangle, rounded corners, minimum width=2.4cm, minimum height=0.5cm,text centered, draw=black]
\tikzstyle{auteurfinal}=[rectangle, rounded corners, minimum width=4.6cm, minimum height=0.5cm,text centered, draw=black]

\title{Calcul effectif de la cohomologie des faisceaux constructibles sur le site étale d'une courbe algébrique}
\author{Christophe \textsc{Levrat}}
\date{\today}

\usepackage{tocloft}

\begin{document}



\thispagestyle{empty}

\begin{center}
\includegraphics[scale=0.2]{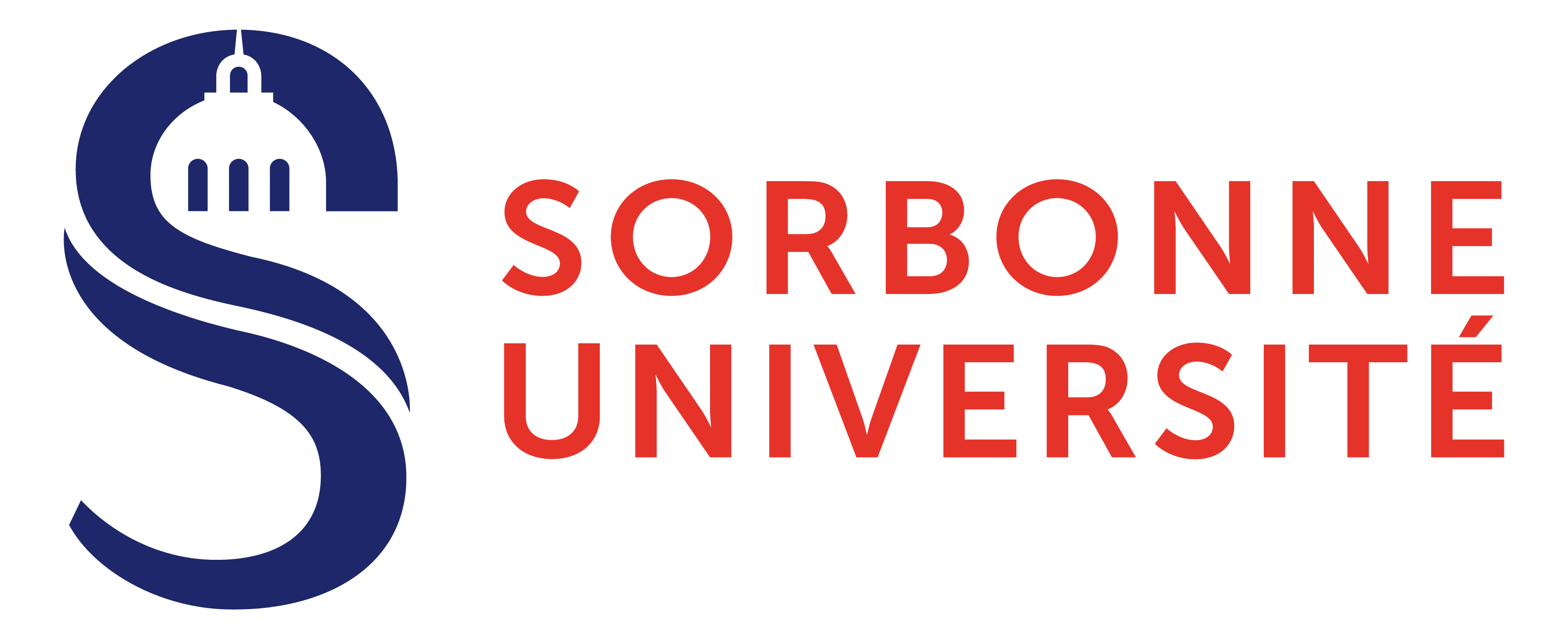}
\end{center}

\begin{center}
\vspace{\stretch{1}}
{\Large \textbf{École doctorale de Sciences Mathématiques de Paris Centre}}

\vspace{\stretch{2}}

{\Huge \textsc{Thèse de doctorat}}

\vspace{\stretch{1}}

{\LARGE Discipline : Mathématiques}

\vspace{\stretch{3}}

{\large présentée par}
\vspace{\stretch{1}}

\textbf{{\LARGE Christophe \textsc{Levrat}}}

\vspace{\stretch{2}}
\hrule
\vspace{\stretch{1}}
{\LARGE \textbf{Calcul effectif de la cohomologie des faisceaux }}
\medbreak 

{\LARGE \textbf{constructibles  sur le site étale d'une courbe }}
\medbreak

\vspace{\stretch{1}}
\hrule
\vspace{\stretch{2}}

{\Large 
dirigée par David \textsc{Madore} et Fabrice \textsc{Orgogozo}}

\vspace{\stretch{4}}

{\Large Soutenue le~30 septembre 2022 devant le jury composé de :}
\begin{large}

\vspace{\stretch{2}}
{ 
\begin{tabular}{lll}
M. Xavier \textsc{Caruso}\hspace{1cm} & Université de Bordeaux \hspace{0.2cm} & rapporteur \\
M. Alain \textsc{Couvreur} & INRIA Saclay & examinateur \\
M. Bruno \textsc{Kahn} & Sorbonne Université & président \\
M. Davide \textsc{Lombardo} & Università di Pisa & rapporteur \\
M. David \textsc{Madore} & Télécom Paris & directeur \\
M. Fabrice \textsc{Orgogozo} & Sorbonne Université & directeur \\
M. Hugues \textsc{Randriam} & ANSSI \& Télécom Paris & examinateur 
\end{tabular}
}
\end{large}

\end{center}


\newpage
\thispagestyle{empty}

\vspace*{\fill}

\noindent\begin{center}
\begin{minipage}[t]{(\textwidth - 3cm)/2}
Institut de mathématiques de Jussieu-Paris Rive gauche. UMR 7586. \\
Boîte courrier 247 \\
4 place Jussieu \\
\numprint{75252} Paris Cedex 05
\end{minipage}
\hspace{1.5cm}
\begin{minipage}[t]{(\textwidth - 3cm)/2}
 Sorbonne Université \\
 École doctorale de Sciences\\
 Mathématiques de Paris Centre.\\
 4 place Jussieu\\
 \numprint{75252} Paris Cedex 05
 \end{minipage}
\end{center}

\newpage

\frontmatter

{\topskip0pt
\vspace*{\fill}

\hfill
\begin{minipage}[t]{.67\textwidth}
„Sie meinen“, fragte Lukas, „Sie werden gar nicht kleiner, wenn Sie näher kommen? Und Sie sind auch nicht wirklich so riesengroß, wenn Sie weit entfernt sind, sondern es sieht nur so aus?“\\
„Sehr richtig“, antwortete Herr Tur Tur. „Deshalb sagte ich; ich bin ein Scheinriese.“ \\~\\
Michael Ende, \textit{Jim Knopf und Lukas der Lokomotivführer}
\end{minipage}
\vskip2cm
\vspace*{\fill}
}

\newpage

\cleartooddpage

\chapter*{Résumé/Abstract}
\vskip -1.6cm
(FR) Cette thèse porte sur la représentation algorithmique des faisceaux constructibles de groupes abéliens sur le site étale d'une variété sur un corps algébriquement clos, ainsi que sur le calcul effectif de leur cohomologie lorsque leur torsion est inversible dans le corps. Nous décrivons trois représentations de ces faisceaux sur les courbes lisses ou nodales, ainsi que des algorithmes permettant d'effectuer un certain nombre d'opérations (noyaux et conoyaux de morphismes, images directes et réciproques, $\Hom$ interne et produit tensoriel) sur ces faisceaux. Nous présentons un algorithme de calcul du complexe de cohomologie d'un faisceau localement constant constructible sur une courbe lisse ou nodale $X$, et en déduisons une description explicite du foncteur $\RG(X,-)\colon \DD^b_c(X,\ZZ/n\ZZ)\to \DD^b_c(\ZZ/n\ZZ)$, fonctorielle en le schéma $X$ et le complexe constructible considéré. En particulier, si $X$ et le faisceau $\F$ proviennent par changement de base d'un sous-corps parfait, nous décrivons l'action de Galois sur le complexe $\RG(X,\F)$ calculé. Nous donnons des bornes explicites sur le nombre d'opérations effectuées par l'algorithme calculant $\RG(X,\F)$. Nous donnons également une description explicite des cup-produits dans la cohomologie des faisceaux localement constants constructibles sur les courbes projectives lisses. Enfin, nous indiquons comment déduire de ces algorithmes une façon de calculer la cohomologie d'un faisceau constant sur une surface lisse fibrée sur la droite projective. \vskip 0.3cm
Mots-clés : cohomologie étale, cohomologie galoisienne, géométrie algébrique effective, courbe algébrique, complexité. \vskip 0.75cm

(ENG) This thesis deals with the algorithmic representation of constructible sheaves of abelian groups on the étale site of a variety over an algebraically closed field, as well as the explicit computation of their cohomology. We describe three representations of such sheaves on curves with at worst nodal singularities, as well as algorithms performing various operations (kernels and cokernels of morphisms, pullback and pushforward, internal $\Hom$ and tensor product) on these sheaves. We present an algorithm computing the cohomology complex of a locally constant constructible sheaf on a smooth or nodal curve, which in turn allows us to give an explicit description of the functor $\RG(X,-)\colon \DD^b_c(X,\ZZ/n\ZZ)\to \DD^b_c(\ZZ/n\ZZ)$. This description is functorial in the scheme $X$ and the given complex of constructible sheaves. In particular, if $X$ and the sheaf $\F$ are obtained by base change from a subfield, we describe the Galois action on the complex $\RG(X,\F)$. We give precise bounds on the number of operations performed by the algorithm computing $\RG(X,\F)$. We also give an explicit description of cup-products in the cohomology of locally constant constructible sheaves over smooth projective curves. Finally, we show how to use these algorithms in order to compute the cohomology groups of a constant sheaf on a smooth surface fibered over the projective line.\vskip 0.3cm
Keywords: étale cohomology, Galois cohomology, effective algebraic geometry, algebraic curve, complexity.

\cleartooddpage

\chapter*{Remerciements}

Ces quelques pages seront sans aucun doute les plus lues de tout ce manuscrit. Ce sont les seules à être écrites à la première personne du singulier, et c'est à partir d'elles que les chercheurs qui ne me connaissent pas se feront une image de moi. Par contre, pour vous qui me connaissez, qui m'avez accompagné pendant ces trois ans, ce n'est pas de moi qu'elles parlent, mais de vous.
Ce sont ces pages, également, qui occuperont ceux qui en auront assez de faire semblant d'écouter attentivement pendant la soutenance. Ne vous inquiétez pas, elles sont faites pour : si vous lisez ces lignes à 15h05 alors que j'ai commencé à parler à 15h, dites-vous que d'autres sont en train de faire la même chose, et que je l'aurais peut-être fait aussi si je n'étais pas occupé à soutenir ma thèse.

\paragraph{Le jury} Tout d'abord, je voudrais remercier les chercheurs qui m'ont fait l'honneur de s'intéresser à mon travail. Merci à Xavier Caruso, Alain Couvreur, Bruno Kahn, Davide Lombardo et Hugues Randriam d'avoir accepté de faire partie de mon jury. En particulier, je tiens à remercier Xavier Caruso et Davide Lombardo d'avoir lu mon manuscrit avec attention et suggéré de nombreuses corrections qui ont grandement amélioré la qualité des pages qui suivent. Merci à Alain Couvreur d'avoir toujours répondu à mes questions avec un sourire qui se lisait même dans les mails, et à Bruno Kahn pour sa présence et ses conseils bienveillants pendant ces trois années.

\paragraph{Les directeurs} David et Fabrice, je ne peux pas assez vous remercier pour votre encadrement pendant ces trois années. J'ai beaucoup appris à vos côtés, sur les maths, sur la recherche, et aussi sur mes propres qualités et défauts. Vous avez été présents tous les deux tout au long de ma thèse, et vous ne m'avez pas laissé tomber à une période où d'autres doctorants ont perdu tout contact avec leurs directeurs.
Vous avez su me parler dans les moments où j'étais en difficulté, et avez été optimistes pour moi dans les moments où je ne l'étais plus. Vous m'avez toujours communiqué votre envie que je réussisse, et j'ai fini par réussir ! En bref : merci, vous avez carrément géré.

\paragraph{Les autres chercheurs} Je tiens à remercier Matthieu Rambaud, qui m'a fait découvrir avec enthousiasme le monde de la recherche il y a maintenant plus de quatre ans, et qui m'offre désormais un avenir dans ce monde : merci d'avoir pris de mes nouvelles pendant ces trois années, et merci de croire en moi.
Je remercie également tous les autres chercheurs qui ont répondu à mes questions au cours de cette thèse : Christophe Ritzenthaler, Emmanuel Hallouin, Henri Lombardi... Sur cette liste figure en particulier le regretté Bas Edixhoven, dont une conjecture a inspiré certains des travaux de cette thèse. 

\paragraph{Les enseignants et les enseignés} La thèse, ce n'est heureusement pas que de la recherche ! Sinon, on deviendrait bien vite bien vieux.
Merci à Maxime Wolff pour ses précieux conseils sur la gestion du temps consacré à l'enseignement en début de thèse, à Frédéric Paugam qui m'a permis d'enseigner des choses que je n'aurais jamais cru avoir la chance d'expliquer à qui que ce soit, et à Adrien Deloro pour ses encouragements, mais aussi son humour. Merci enfin à tous les étudiants qui se sont intéressés de près ou de loin aux choses que j'ai eu le plaisir de leur enseigner. J'espère que vous aussi, un jour, aurez la chance d'écrire des remerciements.

\paragraph{Le SdT} Je remercie également mes camarades de l'organisation du Séminaire des Thésards. Aux permanents, Romain, Bram et Sébastien, merci pour votre soutien, votre motivation, et vos efforts pour avoir la carte du labo. Aux doctorants, merci pour ces rendez-vous (quasi-)hebdomadaires et ces discussions passionnantes. Raphaël, tu devrais aussi figurer plus bas, mais je te règle ton compte tout de suite : merci pour ta façon d'être franche et sincère, pour ton humour mais aussi ta capacité à survivre à tous ces exposés de géométrie algébrique. Je ne sais pas si tu te souviens d'une conversation qu'on a eue en sortant du 2bis l'hiver dernier : elle m'avait beaucoup rassuré. Antoine S, merci pour ta passion et ton investissement dans le fonctionnement du séminaire, et pour les conversations constructives sur les problèmes qu'on a pu rencontrer dans nos thèses respectives. Arnaud V, merci pour ta bienveillance, et pour le rôle que tu as joué dans mon intégration dans le monde des doctorants il y a maintenant trois ans.

\paragraph{Mes enseignants} Faire des mathématiques, c'est avant tout avoir de jolies images en tête. Je tiens à remercier tous ceux qui m'ont fait découvrir ces images, à commencer par Charles Vix et ses belles courbes, ainsi que Pascal Guelfi et ses tôles ondulées. Merci à Gautier Hanna qui m'a emmené à mon premier séminaire il y a déjà plus de six ans. Merci également à Nicolas Perrin qui m'a donné le goût de l'algèbre commutative et la géométrie algébrique, à Mohamed Krir qui m'a appris à faire des calculs concrets sur des courbes (je trouvais ça inutile à l'époque, qu'est-ce que je pouvais être ignorant), et à Benoît Stroh qui m'a donné par son enseignement fantastique l'envie irrésistible de suivre le M2 de Jussieu.

\paragraph{Les doctorants} 
Tout d'abord, je tiens à remercier mes camarades de bureau : Guido, qui est arrivé à Jussieu le même jour que moi et m'a enrôlé dans un groupe de travail quelques jours plus tard, et Arnaud E, qui même après son déménagement est encore venu me parler presque tous les matins pour me faire découvrir des maths fascinantes. Je vous souhaite toute la réussite au monde pour le grand avenir auquel vous êtes promis. La salle de convivialité du cinquième étage a vu passer des dizaines de doctorants pendant ces trois années, et je voudrais tous les remercier pour les conversations intéressantes que nous avons eues. Au début, il y avait la vieille garde : en particulier, Sylvain et ses innombrables histoires de voyages, grand Mathieu et ses conseils avisés sur l'enseignement, Jean-Michel avec sa prononciation soviétique de la lettre h et son goût du poivre... De ma génération, je tiens d'abord à remercier Grace pour son soutien et sa compréhension dans toutes les situations. Des cours de crypto de l'UVSQ à la fin de thèse en passant par l'enseignement à Jussieu, ta présence a toujours été synonyme de joie et de bienveillance. Merci aussi à Haowen, pas seulement pour les découvertes culinaires, mais surtout pour tous les moments de bavardage ! 
Germain, merci de m'avoir fait réviser le Erlkönig malgré moi, et de m'avoir montré qu'on pouvait trouver partout autour de soi des choses à compter et à appeler machinbidulèdre.
Jacques, merci de m'avoir fait découvrir avec ferveur les délires du monde fascinant de la naturopathie ; et même si je ne joue plus à ce jeu, j'ai toujours trouvé une familiarité réconfortante dans tes descriptions des actualités de League of Legends. Maxime, merci de faire une forme de géométrie algébrique que j'ai parfois l'impression de comprendre !
Enfin, merci à tous les autres habitants du couloir pour chaque conversation et chaque rigolade : Xavier, Ilias, Alexandre, Benoît, petit Mathieu, Mahya, Anna F, petit et grand Thomas, Christina, Perla, Thibaut, Thibault, Vadim, Raphaël, Sebastian, Antoine T, Chenyu, Fabrizio, Gabriel, Mathieu H, Adrien, Anna RS, Eva, Thiago, Nelson.

\paragraph{Les amis d'avant ou d'ailleurs} La thèse, c'est beaucoup de travail, et le travail le plus important consiste à survivre à tout ce travail. Merci à vous tous pour chaque soirée, chaque repas, chaque après-midi de jeux, chaque bavardage, chaque message de soutien. Je voudrais rappeler à chacune et chacun d'entre vous un bon souvenir que nous avons partagé pendant ces trois années.
\begin{itemize}
\item Ceux de TPS : Alexandra, David, Guillaume et Camille, Issam, Léna... merci de m'avoir intégré à votre groupe d'amis à tel point que je crois moi-même avoir vu J***** danser sur une table du Fouaille après un cocktail champagne-café. D,I,L, vous souvenez de notre pique-nique sur les quais de Seine après la première vague ? Je travaillais à l'époque sur un document détaillant ce qui est maintenant résumé dans l'annexe \ref{subsec:decprim}. A,C,G,I,L, vous vous souvenez de l'anniv de l'une d'entre nous dans un bar en sous-sol, puis un très bon resto avec en dessert une tartine très addictive ? Je travaillais sur la section \ref{subsec:X2gal} cette semaine-là.
\item Celle qui parle luxembourgeois : \'Anh-Lise, merci d'avoir bien voulu rester amie avec un piéton. Merci pour tes mots rassurants il y a un an et demi quand j'étais coincé. Quand on se revoit, rien n'a changé, et c'est génial, et je te remercie pour ça aussi ! Tu te souviens du jour d'hiver où on a pris un café place Stan avec nos mecs ? C'était la dernière année d'insouciance collective. Je travaillais à cette époque sur les algorithmes de décomposition primaire de la section \ref{subsec:decprim}.
\item Ceux du DAp : merci à Bertrand et Sylvie, à Stéphane et Magda, ainsi qu'à Tatyana pour leur accueil chaleureux dans un monde qui n'était pas le mien au départ. En plus d'être des gens super sympa, on mange bien avec vous. Vous vous souvenez de la journée chachlik ? Qu'est-ce que c'était bon ! Je travaillais à cette époque de l'année sur la proposition \ref{prop:divrat}.
\item Ceux de la L3 : Camille et Florian, et Florian aussi, on ne se voit largement pas autant qu'on s'apprécie, et je voudrais vous remercier pour tous les bons moments passés ensemble.
Vous vous souvenez sûrement de nos parties de billard et de mes difficultés à manger un cornet de glace, mais est-ce que vous vous souvenez de la fois où on s'est vus à Longwy il y a maintenant près de trois ans, et vous m'avez bien éclaté à TowerFall ? \`A cette époque, j'essayais de comprendre les algorithmes de normalisation de l'annexe \ref{subsec:decprim}.
\item Celles des vacances mais pas que : Clémence et Mélanie, merci pour toutes les belles soirées de concert, de veille ou de retour de voyage.
Est-ce que vous vous souvenez encore du premier midi à La Rochelle ? On avait bu de la Guignette sur la terrasse de notre appart, y avait du soleil et les terrasses venaient de rouvrir après 6 mois de fermeture ! La semaine d'avant, j'avais travaillé sur le contenu de la remarque \ref{rk:actgalcst}.
\item Ceux avec lesquels on n'a toujours pas joué au Cowboy Bebop : Gauthier et Mathilde, merci pour votre bonne humeur pendant les soirées passées à jouer à des jeux plus ou moins coopératifs. Vous vous souvenez de la dernière fois qu'on s'est vus ? C'était il y a quelques semaines, et en jouant à Galerapagos, Gauthier, tu as dit "Ta Ferrari n'est pas là ?" quand Mathilde a tiré la carte des clés de voiture, et ça a ouvert tout un monde de souvenirs pour moi. \`{A} ce moment-là, je travaillais sur une correction de la section \ref{subsec:pushforward}.
\item Ceux qui ont vu naître AK : Gham, Ladié, Sulu, Tris. Merci de m'avoir supporté à l'époque où j'étais insupportable, merci d'être resté en contact ces dernières années, merci pour cette alliance devenue amitié. Vous vous souvenez de notre visio sur Messenger pendant le confinement ? On comparait un peu nos situations, et on était tous un peu dans la même. Je travaillais à l'époque sur la section \ref{subsec:foncpic}.
\item Mes co-stagiaires d'un bel été : Isabella, Sarah, merci d'être toujours aussi vivantes et drôles et gentilles avec moi, et de m'avoir montré comment on finit une thèse. Isa, tu te souviens du jour très très chaud où on est allés manger des glaces chez Pozzetto avec M et on n'avait pas assez de serviettes et y en avait partout sur la table ? Je venais d'avoir le financement pour la thèse, et vous deux, vous aviez déjà presque survécu à la première année. Je lisais à l'époque un livre sur le contenu du chapitre \ref{chap:1}. Sarah, tu te souviens de notre conversation derrière le Turing au mois de mai ? Il faisait beau et tu m'as ramené au RER et on était trop contents de s'être revus. J'étais en pleine rédaction de la section \ref{sec:mo} à ce moment-là.
\item Celui qui m'a fait découvrir le goulasch au tofu : Julien, merci d'avoir pris le temps de t'intéresser à ce jeune groupe de stagiaires qui était à Saclay pendant que tu rédigeais ta thèse, et d'avoir été depuis un ami fidèle et une source de très bons conseils cinématographiques.
Tu te souviens de la fois où je suis venu te voir à Rennes, et tu m'as emmené dans une crêperie pompeuse avec des verres à vin gigantesques et des vraies serviettes dans les toilettes ? Tu savais que ça retiendrait mon attention, et je me rappelle toute cette journée comme si c'était hier. C'était le dernier jeudi avant le début de mon contrat de thèse, et je lisais à l'époque un livre sur le contenu de la section \ref{sec:gpfond}.
\item Celui qui met ses courses dans une poche : Marc, merci pour ta gentillesse, ta générosité, pour tous les déplacements que tu as faits pour venir nous voir. Tu te souviens de ta venue en décembre dernier, quand on est tombés par hasard sur une exposition de meubles Prisunic dans le Marais ? La semaine d'avant, j'avais travaillé sans succès sur ce qui est maintenant la section \ref{subsec:cohcomplisses}.
\item Celui que j'ai rencontré derrière un McDo : Marouane, tu es la seule personne qu'on puisse croiser plusieurs fois par hasard dans et autour de la fac, et je te remercie pour cette présence constante et pour toutes nos conversations stimulantes. Tu te souviens de la fois où t'es venu dans mon bureau tôt le matin pour parler de pronostics concernant l'euro de foot ? Je travaillais à l'époque sur la section \ref{repet}.
\item Celles des petits écoliers : Mathilde et Sandrine, merci pour votre gentillesse et votre bienveillance à toute épreuve, pour nos réunions pleines de joie et de bonne humeur. Est-ce que vous vous souvenez de notre pique-nique sur les quais au mois de mai, et de l'embrouille sur les tranches de jambon ? Je rédigeais cette semaine-là la section \ref{sec:mo}.
\item Celui qui a tout relu : Pascal, merci infiniment pour cet effort aussi conséquent que désintéressé. Tu te souviens de la fois en octobre 2019 où on s'est vus au Baker Street, et tu m'as donné des conseils précieux pour parler de quelque chose d'important à des gens qui m'importent ? J'essayais à l'époque de comprendre les algorithmes de normalisation évoqués dans l'annexe \ref{subsec:decprim}.
\item Le meilleur ami de mon enfance : Pierre, merci pour ton enthousiasme et ton énergie, et pour tout l'allemand que tu as retenu ! Tu te souviens du soir il y a environ un an où on a bu un verre près de la fac, et avant que Nicolas nous rejoigne, c'était la première fois depuis des lustres qu'on se voyait juste tous les deux ? On était assis en terrasse aux Rattrapages, on se racontait nos vies et on était juste super contents. La semaine en question, je travaillais sur l'annexe \ref{chap:A3}.
\item Celui qui m'a appris alors que j'étais censé lui apprendre : Quentin C, merci pour ta bonne humeur, ta philosophie, ta générosité. Merci de m'avoir laissé essayer de t'expliquer un peu de maths pendant deux ans, et d'être quand même devenu mon ami !
Tu te souviens de la soirée fin avril où tu nous as apporté une plante et tu m'as bien battu à FIFA ? Cette plante s'appelle maintenant Romy, et elle se porte très bien. Je rédigeais à ce moment-là le chapitre \ref{chap:2}.
\item Celui qui avait un trait vert sur sa télé : Quentin R, merci d'abord pour deux dessins qui figurent dans cette thèse et pour lesquels tu refuses que je te remercie plus que ça. Merci aussi pour ton soutien dans toute cette aventure, et pour tes expressions qui me font toujours bien rigoler. Et puis merci d'être persuadé que je vais réussir quand je suis persuadé que je vais échouer. Tu te souviens encore de ta soutenance il y a un an ? T'avais fait une vidéo pour expliquer aux gens comment accéder à la salle, et c'était vachement bien pensé. C'était un jeudi après-midi, et le matin, j'avais corrigé des choses dans la section \ref{sec:jin}.
\item Celui qui a connu Tigre Bois : Thomas, il y a tellement de choses pour lesquelles je dois te remercier dans ma vie que je ne vais pas tout étaler ici. Mais bien sûr, merci pour ton soutien presque quotidien, merci pour ta gentillesse, et merci de ne pas trop mal supporter les buts marqués au premier poteau. Est-ce que tu te souviens du vendredi soir juste avant mon déménagement, il y a un peu plus d'un an et demi ?  J'étais très stressé parce que ma thèse n'avançait pas, et tu avais fait ton possible pour me rassurer. Ce jour-là, je me cassais les dents sur la suite de la section \ref{subsec:gysinlisse}.
\item Celui que je n'ai pas vu, mais entendu, et que je veux remercier quand même : Tristan, merci d'être plus jeune et en même temps plus et moins sage que moi ! Un jour, on visitera le château de Fougeret ensemble.
\item Celui que je n'ai ni vu, ni entendu, mais lu, et qui aime le BK même si ce dernier devient un peu décevant : Rodolphe, merci pour ta bienveillance en toute circonstance, et merci de me rappeler fréquemment des souvenirs qui risqueraient de disparaître de nos mémoires si on ne les réactivait pas de temps en temps.
\item Ceux avec lesquels on n'a pas trop parlé de la thèse mais que je veux citer quand même, et avec chacun desquels on a fait au moins une raclette que je veux leur rappeler : merci à Kevin et Romain, à Martial et Tiago, à Thibault et Yann pour les belles soirées.
\end{itemize}

\paragraph{La famille} J'en ai une géniale. C'est trop privé, je n'en dirai pas plus sur le sujet !

\paragraph{Le meilleur pour la fin} Nicolas, merci pour tout ce qu'il y a de beau dans la vie ; ça y est, j'ai fini aussi, et je n'y serais pas arrivé sans toi.

\cleartooddpage

\tableofcontents

\chapter*{Introduction}

\addcontentsline{toc}{chapter}{Introduction}

\section*{Contexte}

Cette thèse porte sur des méthodes explicites de calcul de groupes de cohomologie étale. Les prérequis à la compréhension de cette introduction et de la suite du manuscrit sont présentées dans le chapitre \ref{chap:1}. Soit $n$ un entier. Soit $k$ un corps de caractéristique première à $n$. 
Les catégories dérivées des faisceaux de $\Lambda\coloneqq \ZZ/n\ZZ$-modules sur les $k$-schémas de type fini sont munies des six opérations $f^\star$, $\R f_\star$, $\R f_!$, $\R f^!$, $\otimes_\Lambda^L$, $\R\underline{\Hom}_\Lambda$ définies par Grothendieck. La condition raisonnable de finitude sur les faisceaux est la constructibilité : un faisceau de $\Lambda$-modules sur un schéma noethérien $X$ est dit constructible si $X$ admet une stratification telle que le faisceau soit un système local sur chaque strate. Les faisceaux constructibles sont les objets noethériens de la catégorie des faisceaux de $\Lambda$-modules sur $X$. Un résultat majeur dû à Grothendieck, Artin, Deligne et plus récemment Gabber pour la formulation la plus générale \cite[XIII, Th. 1.1.1]{travaux_gabber}, affirme que pour des schémas noethériens quasi-excellents sur lesquels $n$ est inversible, la constructibilité est stable par les six opérations. 

Supposons $k$ algébriquement clos. Soit $X$ un schéma de type fini sur $k$. Le résultat précédent implique que pour tout faisceau constructible $\F$ de $\Lambda$-modules sur $X$, les groupes de cohomologie $\HH^i(X,\F)$ sont finis. De plus, ils sont nuls dès que $i>2\dim X$. Un morphisme $f\colon Y\to X$ de $k$-schémas induit par fonctorialité un morphisme $\RG(X,\F)\to \RG(Y,f^\star\F)$. De même, si $X$ provient par changement de base d'un schéma sur un sous-corps parfait $k_0$, $\RG(X,\F)$ est muni d'une action de $\Gal(k|k_0)$. \`A défaut de pouvoir calculer les six opérations explicitement, un objectif atteignable est le calcul d'un complexe représentant $\RG(X,\F)$, ou au moins des groupes $\HH^i(X,\F)$, d'une façon qui permette de représenter ces morphismes de fonctorialité.
Si $X$ est une courbe projective lisse connexe et $\F$ est constant, les groupes $\HH^i(X,\F)$ sont bien connus ; en particulier, $\HH^1(X,\mu_n)$ est le groupe des points de $n$-torsion de la variété jacobienne de $X$. Si $X$ est de dimension $d>1$, les deux techniques de calcul usuelles procèdent par fibration et récurrence sur la dimension. D'une part, pour les variétés projectives, un pinceau de Lefschetz permet d'obtenir une fibration \[\tX\to\PP^1 \] où $\tX$ est un éclatement de $X$. D'autre part, les bons voisinages d'Artin permettent d'obtenir une fibration en courbes affines 
\[ X\to X_{d-1}\to \dots\to X_1\]
où $X_i$ est de dimension $i$. Si les suites spectrales associées à ces fibrations permettent assez rapidement de majorer le rang des $\HH^j(X,\F)$, voire le calculer pour des petites valeurs de $j$, elles ne fournissent toutefois pas immédiatement une description de ces groupes permettant de calculer les morphismes de fonctorialité évoqués ci-dessus. 

La calculabilité des groupes $\HH^i(X,\Lambda)$ a été démontrée par Poonen, Testa et van Luijk en 2015 dans le cas où $k$ est de caractéristique nulle \cite[Th. 7.9]{poonen_testa}. Une preuve en caractéristique quelconque a été donnée quelques mois plus tard par Madore et Orgogozo \cite[Th. 0.1]{mo}. Ces deux articles décrivent en détail des algorithmes permettant effectivement de calculer les $\HH^i(X,\ZZ/n\ZZ)$, mais ne donnent pas de borne sur le nombre d'opérations effectuées par ces algorithmes, qui paraissent en outre trop inefficaces pour être utilisés dans la pratique. \`A l'autre bout du spectre, Huang et Ieradi \cite[Th. 1]{huang_counting} et plus tard Couveignes \cite[Th. 1]{couveignes_linearizing} ont décrit des algorithmes calculant explicitement des classes de diviseurs sur une courbe projective lisse formant une base de la $n$-torsion de la jacobienne de cette courbe, avec des bornes de complexité très précises. Nous nous sommes intéressés au calcul du foncteur $\RG(X,-)$, lorsque $X$ est une courbe lisse ou nodale, et donnons des bornes explicites sur la complexité du calcul.

Lorsque $X$ provient par changement de base d'un corps fini $\FF_q$, le calcul des $\HH^i(X,\Lambda)$ permet également, par l'intermédiaire de la formule des traces, de compter les $\FF_q$-points de $X$. En particulier, un algorithme de calcul de $\HH^i(X,\Lambda)$ de complexité polynomiale en $n$ et $\log q$ permet de calculer le cardinal de $X(\FF_q)$ en un nombre d'opérations polynomial en $\log q$. Dans le cas où $X$ est une courbe, l'algorithme probabiliste de Huang et Ierardi atteint cette complexité. Dans le cas où $X$ est une surface, l'existence d'un tel algorithme a été conjecturée par Couveignes et Edixhoven dans \cite[Epilogue]{edixcouv} ; il y est suggéré d'employer une fibration de Lefschetz. Une grande partie du travail de cette thèse a été motivée par cette conjecture.

\section*{Contributions}

Soient $k_0$ un corps parfait et $k$ une clôture algébrique de $k_0$. Soit $n$ un entier premier à la caractéristique de $k$. Notons $\Lambda$ l'anneau $\ZZ/n\ZZ$. Soit $X_0$ une courbe lisse géométriquement intègre sur $k_0$. Notons $X=X_0\times_{k_0}k$.

\paragraph{Représentations des faisceaux constructibles sur les courbes} Nous décrivons trois représentations algorithmiques possibles des faisceaux constructibles sur $X$ : par les générateurs ou cogénérateurs classiques de la catégorie des faisceaux constructibles, et par recollement relativement à un ouvert de lissité. Nous décrivons des algorithmes (de complexité élémentaire en les entrées) permettant de passer d'une représentation à une autre. Pour la représentation par recollement, qui est la plus adaptée aux calculs de cohomologie effectués par la suite, nous présentons des algorithmes calculant les opérations suivantes sur les faisceaux : \begin{itemize}[label=$\bullet$]
\item noyau et conoyau de morphismes ;
\item somme directe, produit tensoriel, $\Hom$ interne ;
\item tiré en arrière, poussé en avant par des morphismes entre courbes lisses.
\end{itemize}
Nous adaptons également cette représentation par recollement au cas des courbes nodales.

\paragraph{Construction de revêtements galoisiens de courbes}  Nous montrons comment, à partir d'un algorithme calculant la $n$-torsion de la jacobienne de la compactification lisse de la normalisée de $X$, calculer explicitement un revêtement caractéristique $X_2$ de $X$ de groupe $\HH^1(X,\Lambda)^\vee$, ainsi qu'un revêtement caractéristique de $X_0$ de groupe $\HH^1(X_0,\Lambda)^\vee$. Nous construisons également, à partir de $X_2$, un revêtement caractéristique $X_{2,0}$ de $X_0$ tel que $X_{2,0}\times_{k_0}k\to X$ trivialise les $\Lambda$-torseurs sur $X$. Une adaptation de ces constructions au cas des courbes nodales est également présentée. La notation $X_2$ trouve son origine dans le fait que lorsque $n$ est un nombre premier $\ell$, le revêtement $X_2\to X$ correspond au quotient du complété pro-$\ell$ de $\pi_1(X)$ par le groupe numéro 2 de la série de Frattini descendante (pour plus de détails, voir \cite[§3.1]{mo}).

\paragraph{Calculabilité de la cohomologie sur les variétés} Après l'avoir passé en revue, ainsi que d'autres algorithmes de calcul de la cohomologie, nous montrons que l'algorithme de Madore et Orgogozo, qui calcule les groupes de cohomologie d'un faisceau constant sur un schéma de type fini sur $k$, est de complexité primitivement récursive dès que le schéma en entrée est lisse.

\paragraph{Cohomologie des faisceaux constructibles sur les courbes} Soit $\F_0$ un faisceau constructible de $\Lambda$-modules sur $X_0$, lisse sur un ouvert affine non vide $U_0$, de complémentaire fermé réduit $Z_0$. Notons $M$ sa fibre générique géométrique. Notons $X,U,Z,V$ les changements de base respectifs de $X_0,U_0,Z_0,V_0$ à $k$ et $\F=(\F_0)|_X$. Soit $V\to U$ un revêtement (étale connexe) galoisien trivialisant $\F|_U$. Soit $V_2\to V$ le revêtement décrit ci-avant. Nous décrivons des méthodes efficaces de calcul de ces différents objets et du groupe $G\coloneqq \Aut(V_2|U)$. La courbe affine $U$ est un $K(\pi,1)$ :  la cohomologie de $\F$ sur $U$ est la cohomologie galoisienne du $\pi_1(U)$-module $M$. Notons $C^{12}(G,M)$ le groupe des morphismes croisés $G\to M$.

\begin{iprop}(\ref{prop:RGlisse})
Le morphisme canonique \[ [M\to C^{12}(G,M)]=\tau_{\leqslant 1}\RG(G,M)\to \RG(U,\F)\] est un quasi-isomorphisme.
\end{iprop} 

Pour chaque point $z\in Z$, notons $I_z$ le groupe d'inertie d'un point de la compactification lisse de $V_2$ au-dessus de $z$, et $P_z$ le groupe d'inertie sauvage correspondant. Soit $\phi_z\colon \F_z\to M^{I_z}\subseteq M^{P_z}\xrightarrow{\sim}M_{P_z}$ le morphisme de recollement en $z$ composé avec l'isomorphisme canonique $M^{P_z}\xrightarrow{\sim}M_{P_z}$. Ici, les indices (resp. exposants) désignent les modules des coinvariants (resp. invariants) sous les groupes en question, et l'isomorphisme $M^{P_z}\xrightarrow{\sim}M_{P_z}$ associe à un élément de $M$ invariant sous $P_z$ sa classe dans le quotient $M_{P_z}$. Le morphisme $\phi_z$ fait partie des données définissant $\F$. Nous construisons dans le lemme \ref{lem:secIP} une section au morphisme d'inflation $C^{12}(I_z/P_z,M^{P_z})\to C^{12}(I_z,M)$.

\begin{itheorem}(\ref{th:compRG}) Le cône du morphisme de complexes suivant représente $\RG(X,\F)[1]$.
\[
\begin{adjustbox}{width=\textwidth}{
\begin{tikzcd}
M\oplus \bigoplus_{z\in Z}\F_z \arrow[d,"\bigoplus_{z\in Z}(\id-\phi_z)"] \arrow[r,"{(\partial_G,0)}"] & C^{12}(G,M) \arrow[d,"\bigoplus_z \res_{I_z}^G"]\arrow[r] & \bigoplus_{z\in Z}\HH^1(I_z/P_z,M_{P_z}) \arrow[d,"\id"]\arrow[r] & 0 \\
\bigoplus_{z\in Z}M_{P_z} \arrow[r,"\bigoplus_z\partial_{I_z}"] \arrow[r] & \bigoplus_{z\in Z}C^{12}(I_z/P_z,M_{P_z}) \arrow[r] & \bigoplus_{z\in Z}\HH^1(I_z/P_z,M_{P_z}) \arrow[r] & 0
\end{tikzcd}
}
\end{adjustbox}
\] 
\end{itheorem}
Lorsque $V$ provient de $k_0$, le complexe obtenu est naturellement muni d'une action de $\Gal(k|k_0)$. Nous donnons également une variante de ce résultat permettant, à partir de la donnée d'un revêtement $V_0\to U_0$ qui trivialise $\F_0$, de calculer l'action de $\Gal(k|k_0)$ sur $\RG(X,\F)$. De plus, nous adaptons ce résultat au calcul de $\RG(X,\K)$, où $\K$ est un complexe de faisceaux constructibles. 
Enfin, lorsque $k_0=\FF_q$, nous étudions la complexité de l'algorithme construisant ce complexe, basé sur des algorithmes existants de calcul de points de $n$-torsion dans la jacobienne d'une courbe. 

\begin{itheorem}(\ref{th:complexiteRG}) Notons $g$ le genre de $X$, et $r$ le nombre de points à l'infini de $U$. Soit $d$ le degré de $V\to U$. Soit $m$ un entier tel que $M$ et les fibres de $\F$ en les points de $X-U$ soient des quotients de $\Lambda^m$. 
Il existe un algorithme probabiliste (Las Vegas) qui calcule un complexe de $\Lambda[\Gk]$-modules représentant $\RG(X,\F)$ en
\[ \mathcal{P}(d,g,n,r,m,\log q)^{2^{O(d(g+r))}}\]
opérations dans $\FF_q$, où $\mathcal{P}$ est un polynôme. Si $V$ admet un modèle plan à singularités ordinaires de degré $O(g)$, cette complexité devient \[ \mathcal{P}(d,g,n,r,m,\log q)^{O\left((d(g+r))^2\right)}.\]
\end{itheorem}

\paragraph{Calcul de cup-produits} Supposons $X$ projective de genre non nul. Soit $Y\to X$ un revêtement galoisien de $X$ de degré divisible par $n$.  Nous donnons une description alternative du $\HH^2$ d'un faisceau lisse sur $X$ qui permet de calculer les cup-produits $\HH^1\times\HH^1\to \HH^2$. Soient $Y_2$ le revêtement de $Y$ de groupe $\HH^1(Y,\Lambda)^\vee$, et $Y_3$ le revêtement de $Y_2$ de groupe $\HH^1(Y_2,\Lambda)^\vee$. Notons $G_2=\Aut(Y_2|X)$ et $G_3=\Aut(Y_3|X)$. 

\begin{iprop}(\ref{prop:H2im}) Pour tout faisceau lisse $\F$ de $\Lambda$-modules sur $X$ trivialisé par $Y$ de fibre $M$, le morphisme
\[ \im(\HH^2(G_2,M)\to \HH^2(G_3,M))\to \HH^2(\pi_1(X),\F)\]
est un isomorphisme.
\end{iprop}

\begin{itheorem}(\ref{th:cupprod}) Soient $\F$ et $\G$ deux faisceaux lisses de $\Lambda$-modules sur $X$ trivialisés par $Y$, de fibres respectives $M$ et $N$. Le cup-produit 
\[ \HH^1(X,\F)\times \HH^1(X,\G) \to \HH^2(X,\F\otimes \G) \]
est réalisé par la composée
\[ \HH^1(G_2,M)\times \HH^1(G_2,N) \xrightarrow{\cup} \HH^2(G_2,M\otimes N) \rightarrow \im(\HH^2(G_2,M\otimes N)\to \HH^2(G_3,M\otimes N)).\]
\end{itheorem}

\paragraph{Cohomologie des faisceaux constants sur les surfaces} En appliquant la procédure prévue dans \cite[Epilogue]{edixcouv}, nous indiquons comment appliquer les algorithmes précédents au calcul de la cohomologie de $\Lambda$ sur une surface projective lisse sur $k$ fibrée sur $\PP^1$ par un pinceau de Lefschetz.

\section*{Organisation du manuscrit}

Le chapitre \ref{chap:1} rappelle les définitions et théorèmes célèbres de la cohomologie étale. Il permettra au lecteur peu familier avec ces notions de trouver l'ensemble des résultats utilisés dans la suite du manuscrit, et pourra aisément être laissé de côté par le lecteur expert.

Le chapitre \ref{chap:2} recense des résultats classiques sur la cohomologie et les revêtements des courbes lisses, et leurs analogues moins connus concernant les courbes nodales. En particulier, il contient la construction et l'étude d'un revêtement caractéristique utilisé dans la suite.

Le chapitre \ref{chap:3} est consacré à la description des diverses représentations explicites des faisceaux constructibles de groupes abéliens sur une courbe lisse ou nodale, ainsi qu'aux algorithmes servant à passer d'une représentation à une autre. 
Nous y décrivons un certain nombre d'opérations sur ces faisceaux dans une représentation par recollement relativement à un ouvert de la courbe sur lequel le faisceau est lisse.

Les algorithmes existants pour le calcul de la cohomologie sont décrits dans le chapitre \ref{chap:4}. Tout d'abord, nous résumons l'algorithme de Madore et Orgogozo qui calcule la cohomologie d'un faisceau sur une variété quelconque, et montrons qu'il est primitivement récursif dans le cas particulier des variétés lisses. L'algorithme de Huang et Ierardi et l'algorithme de Couveignes permettent de calculer la cohomologie d'un faisceau constant sur une courbe projective lisse ; nous expliquons comment adapter l'un et l'autre au calcul de la division par $n$ dans le groupe de Picard de la courbe. Ces deux algorithmes serviront de base à nos méthodes de calcul. Enfin, nous décrivons la méthode de Jin servant à calculer la cohomologie d'un faisceau lisse sur une courbe lisse.

Nous présentons dans le chapitre \ref{chap:5} des méthodes de calcul de la cohomologie des faisceaux constructibles sur une courbe $X$ lisse ou nodale. Le cas le plus simple, traité en premier, est celui des faisceaux constants, qui se résume à des calculs dans la jacobienne de la (compactification lisse de la normalisation de la) courbe. Vient ensuite le cas des faisceaux lisses, qui se ramène à des calculs de cohomologie galoisienne. Une section est alors consacrée aux calculs de cup-produits dans la cohomologie de ces faisceaux. Enfin, le calcul de la cohomologie des faisceaux constructibles et des complexes d'iceux s'effectue par recollement, en s'appuyant sur le cas des faisceaux lisses. Nos résultats fournissent un calcul explicite du foncteur \[\RG(X,-)\colon \DD^b_c(X,\ZZ/n\ZZ)\to \DD^b_c(\ZZ/n\ZZ)\]
où $n$ est un entier inversible sur $X$. Tous ces calculs sont fonctoriels en $X$ et en le faisceau étudié ; nous donnons des bornes de complexité précises sur les algorithmes présentés.

Le chapitre \ref{chap:6} est dédié au calcul de la cohomologie des surfaces projectives lisses. Nous montrons comment la méthode classique de fibration en courbes au moyen d'un pinceau de Lefschetz peut être utilisée en conjonction avec nos algorithmes sur les courbes pour calculer la cohomologie d'un faisceau constant sur une surface, et ainsi compter les points sur cette surface si elle provient d'un corps fini. 

Les annexes contiennent des compléments de nature algorithmique. L'annexe \ref{chap:A1} fournit des précisions sur les différentes classes de complexité ainsi que sur la notion de corps calculable, et rappelle la complexité des opérations classiques dans les anneaux de polynômes. L'annexe \ref{chap:A2} détaille la représentation des variétés algébriques utilisée par les algorithmes, et résume les opérations classiques en géométrie algébrique effective. Enfin, les algorithmes spécifiques aux courbes ainsi que leur complexité sont présentés dans l'annexe \ref{chap:A3}.

\mainmatter

\setcounter{page}{0}
\thispagestyle{empty}

\cleartooddpage

\chapter{Cohomologie étale et groupe fondamental}\label{chap:1}

\setcounter{page}{1}

Ce premier chapitre recense les définitions des objets manipulés par la suite, ainsi que les principaux résultats de la cohomologie étale utilisés dans la suite de ce texte. Afin de rendre la lecture plus abordable et citer plus rapidement des théorèmes importants, ces résultats ne sont pas présentés dans l'ordre habituel de leur démonstration. Le lecteur averti pourra passer directement au chapitre suivant.\\

Nous supposons connus les résultats classiques d'algèbre homologique, en particulier la construction des catégories et foncteurs dérivés. Une introduction brève et efficace à ces notions se trouve dans \cite{derived_working}, et une introduction plus complète dans \cite{yekutieli}.

\paragraph{Avertissement} Dans ce chapitre, le terme \textit{schéma} signifiera \textit{schéma noethérien}. Tous les schémas rencontrés dans la suite de cette thèse seront effectivement noethériens.

\section{Morphismes étales et groupe fondamental}

\subsection{Morphismes étales}

\begin{df}  Un morphisme de schémas $f\colon Y\to X$ est dit étale en un point $y$ de $Y$ s'il est plat et non ramifié en $y$. Il est dit étale s'il est étale en tout point de $Y$.  Nous noterons $\Xet$ la catégorie dont les objets sont les couples $(Y,f)$ où $f\colon Y\to X$ est un morphisme étale, et les morphismes $(Y,f)\to (Y',f')$ sont les morphismes de schémas $\phi\colon Y\to Y'$ tels que $f'\circ\phi=f$. 
Nous noterons $\Fet_X$ la sous-catégorie de $\Xet$ formée des $X$-schémas finis étales, appelés revêtements étales.
\end{df}

\begin{ex} \begin{enumerate} 
\item Toute immersion ouverte est étale.
\item L'immersion d'un fermé strict d'un schéma connexe n'est jamais étale.
\item Soient $A$ un anneau et $h\in A[t]$. Soit $g\in A[t]$ un polynôme unitaire tel que $g'$ soit inversible dans $A[t,h^{-1}]/(g)$. Alors le morphisme \[ \Spec A[t,h^{-1}]/(g)\to \Spec A\] 
est étale. Un morphisme de cette forme est appelé morphisme étale standard.
\item Si $k'/k$ est une extension galoisienne de corps et $X$ est un schéma sur $k$ alors le morphisme $X_{k'}\to X$ est étale.
\item Si $k$ est un corps et $n$ est un entier inversible dans $k$, l'endomorphisme de $\GG_m=\Spec k[t,t^{-1}]$ défini par $t\mapsto t^n$ est étale.
\item Si $E$ est une courbe elliptique sur un corps $k$ et $n$ est un entier inversible dans $k$ alors la multiplication par $n$ est un endomorphisme étale de $E$.
\item Le morphisme $\Spec k[x,y]/(y-x^2)\to \Spec k[y]$ donné par $(x,y)\mapsto y$ est plat, mais ramifié au point $(0,0)$.
\item Soit $X=\Spec k[x,y]/(y^2-x^2(x+1))$ la cubique nodale, et soit $\nu\colon\A^1\to X$ le morphisme de normalisation donné par $t\mapsto (t^2-1,t^3-t)$. Le morphisme $\nu$ est non ramifié dès que la caractéristique de $k$ est différente de 3, mais n'est pas plat.
\end{enumerate}
\end{ex}
Il existe de nombreuses caractérisations équivalentes de l'étalitude d'un morphisme (voir p. ex. \cite[02GU]{stacks}). L'une des plus explicites est la suivante.
\begin{prop} Soit $f\colon Y\to X$ un morphisme de schémas. Le morphisme $f$ est étale en $y\in Y$ si et seulement s'il existe un ouvert affine $V=\Spec B$ de $Y$ contenant $y$, un ouvert affine $U=\Spec A$ de $X$ contenant $f(y)$ et une présentation \[ A=B[x_1\dots x_n]/(f_1\dots f_n)\]
tels que $\det(\partial f_i/\partial x_j)\in B_y^\times$.
\end{prop}

L'étalitude d'un morphisme de variétés a une interprétation géométrique très simple.
\begin{prop}\cite[Prop. 2.9]{milneLEC} Soit $f\colon Y\to X$ un morphisme de variétés sur un corps algébriquement clos $k$.  Le morphisme $f$ est étale si et seulement si, pour chaque point fermé $y$ de $Y$, le morphisme induit entre les cônes tangents $\TC_{f(y)}X\to \TC_yY$ est un isomorphisme.
\end{prop}

Voici quelques propriétés classiques des morphismes étales. 
\begin{prop}\cite[02GN,02GO,03WT]{stacks}
\begin{enumerate}
\item Un morphisme étale est ouvert.
\item La composée de morphismes étales est étale.
\item Tout changement de base d'un morphisme étale est étale.
\item Soient $f\colon Y\to X$ et $g\colon X\to S$ des morphismes de schémas. Si $g$ et $g\circ f$ sont étales alors $f$ l'est également.
\end{enumerate}
\end{prop}

\subsection{Revêtements galoisiens}

\begin{df} Soit $X$ un schéma connexe. Un revêtement galoisien de $X$ est un morphisme fini étale $f\colon Y\to X$, où $Y$ est connexe, tel que le groupe d'automorphismes $\Aut(Y|X)$ agisse transitivement sur les fibres géométriques de $f$.
\end{df}

\begin{rk} Comme un morphisme fini est fermé et un morphisme étale est ouvert, tout revêtement galoisien d'un schéma connexe est surjectif.
\end{rk}

Le lemme suivant montre en quoi cette notion généralise celle d'extension galoisienne de corps ; en particulier, une extension finie $k'/k$ est galoisienne si et seulement si le morphisme $\Spec k'\to \Spec k$ est un revêtement galoisien.
\begin{lem}\cite[03SF]{stacks} Un morphisme fini étale de schémas $f\colon Y\to X$ est galoisien si et seulement si le groupe $\Aut(Y|X)$ est d'ordre $\deg(f)$.
\end{lem}

Comme pour les extensions de corps, il y a une notion de clôture galoisienne.
\begin{df} Soit $f\colon Y\to X$ un morphisme fini étale, où $X$ est un schéma connexe. Un morphisme $g\colon Z\to Y$ tel que la composée $Z\to X$ soit un revêtement galoisien est appelé clôture galoisienne de $f$ si tout autre $X$-morphisme d'un revêtement galoisien de $X$ vers $Y$ se factorise par $g$.
\end{df}
\[
\begin{tikzcd}
Z \arrow[r,"g"] & Y \arrow[r] & X \\
& M \arrow[ul,dashed,"\exists !"] \arrow[u] \arrow[ur,"{\rm gal}",swap] & 
\end{tikzcd}
\]
\begin{prop}\cite[Prop. 5.3.9]{szamuely} La clôture galoisienne existe et est unique à isomorphisme près ; si $\bar{x}$ est un point géométrique de $X$, de préimages les points géométriques $\bar{y}_1,\dots,\bar{y}_d$ de $Y$, la clôture galoisienne de $f$ est la composante connexe dans $Y\times_X \dots \times_X Y$ ($d$ facteurs) de $(\bar{y}_1,\dots,\bar{y}_d)$.
\end{prop}

\subsection{Groupe fondamental}\label{sec:gpfond}

Soit $X$ un schéma connexe. Soit $\bar x$ un point géométrique de $X$. Rappelons que $\Fet_X$ désigne la catégorie des $X$-schémas finis étales.

\begin{df} Le groupe fondamental $\pi_1(X,\bar x)$ de $X$ en $\bar x$ est le groupe des automorphismes du foncteur fibre $\Fib_{\bar x}\colon \Fet_X\to\Set, Y\mapsto Y_{\bar x}$.
\end{df}

Considérons la catégorie cofiltrante $I_{X,\bar x}$ des couples $(f_Y\colon Y\to X,\bar y)$, où $f_Y\colon Y\to X$ est un revêtement galoisien (connexe) et $y$ est un point géométrique de $Y$ vérifiant $\bar x=f_Y\circ \bar y$. Un morphisme de tels couples est un morphisme de $X$-schémas compatible avec les points géométriques. 

\begin{prop}\cite[Prop. 5.4.6]{szamuely} Le foncteur $\Fib_{\bar x}$ est pro-représenté par $\lim_{Y\in I_{X,\bar x}^{\op}} Y$, c'est-à-dire que pour tout $X$-schéma étale $Z\to X$, le morphisme
\[ \begin{array}{rcl}\colim_{(f_Y,\bar y)\in I_{X,\bar x}^{\op}}\Hom_X(Y,Z)&\longrightarrow& \Fib_{\bar x}(Z) \\ f&\longmapsto& f(\bar y)\end{array}\]
est un isomorphisme. Le groupe $\pi_1(X,\bar x)$ est isomorphe à $\lim_{(Y,\bar y)\in I_{X,\bar x}^{\op}} \Aut(Y|X)^{\op}$. C'est en particulier un groupe profini.
\end{prop}
\begin{lem} Soit $(Y,\bar y)$ un objet de $I_{X,\bar x}$. Il y a une suite exacte de groupes profinis :
\[ 1 \to \pi_1(Y,\bar y)\to \pi_1(X,\bar x)\to \Aut(Y|X)^{\op}\to 1.\]
\begin{proof} Pour tout $(Y',\bar y')\in I_{Y,\bar y}$, la clôture galoisienne de $Y'\to X$ se factorise par $Y'\to Y$. Par conséquent, si $J$ désigne la sous-catégorie de $I_{Y,\bar y}$ des $(Y',\bar y')$ tels que $Y'\to X$ soit galoisienne, $\pi_1(Y,\bar y)$ est encore isomorphe à $\lim_{(Y',\bar y')\in J^{\op}}\Aut(Y'|Y)^{\op}$. Comme $J$ est encore une sous-catégorie de $I_{X,\bar x}$, ceci définit une injection $\pi_1(Y,\bar y)\to \pi_1(X,\bar x)$, qui est le noyau de $\pi_1(X,\bar x)\to \Aut(Y|X)^{\op}$.
\end{proof}
\end{lem}
Le théorème central de la théorie de Galois-Grothendieck est le suivant.
\begin{theorem}\label{th:galsch}\cite[Th. 5.4.2]{szamuely} 	Le foncteur $\Fib_{\bar x}$ induit une équivalence entre $\Fet_X$ et la catégorie des ensembles finis munis d'une action à gauche continue de $\pi_1(X,\bar x)$. Cette équivalence fait correspondre les revêtements galoisiens de $X$ aux quotients finis de $\pi_1(X,\bar x)$.
\end{theorem}
\begin{ex}\begin{enumerate}
\item Soit $k_0$ un corps. Soit $k_0^\sep$ une clôture séparable de $k_0$. Notons $\bareta$ le point géométrique correspondant de $\Spec k_0$. Alors \[\pi_1(\Spec k_0,\bareta)=\Gal(k_0^\sep|k_0).\]
\item Soit $k$ un corps algébriquement clos de caractéristique nulle. Alors il n'y a pas de revêtement étale non trivial de $\A^n_k$, et \[\pi_1(\A^n_k,\bareta)=\pi_1(\PP^n_k,\bareta)=0.\]
\end{enumerate} 
\end{ex}
De même qu'en topologie, le choix du point-base $\bar x$ est indispensable à la fonctorialité de $\pi_1$. Soit $f\colon Y\to X$ un morphisme de schémas connexes. Soit $\bar y$ un point géométrique de $Y$ d'image $\bar x$. Notons $B_{Y\to X}=-\times_X Y\colon \Fet_X\to \Fet_Y$ le foncteur de changement de base. Alors $\Fib_{\bar x}=\Fib_{\bar y}\circ B_{Y\to X}$, ce qui induit un morphisme de groupes $f_\star\colon \pi_1(Y,\bar y)\to \pi_1(X,\bar x)$. Tout comme en topologie, le point-base n'a aucune influence sur la structure du groupe fondamental.
\begin{prop}\cite[Cor. 5.5.2]{szamuely} Soit $\bar x'$ un autre point géométrique de $X$. Il y a un isomorphisme $\pi_1(X,\bar x)\to \pi_1(X,\bar x')$, unique à un automorphisme intérieur de $\pi_1(X,\bar x)$ près.
\end{prop}
Nous pourrons donc nous permettre l'abus de notation consistant, lorsque $X$ est intègre, à noter $\pi_1(X)$ le groupe $\pi_1(X,\bareta)$, où $\bareta$ est un point générique géométrique de $X$. Dans le cas où $X$ est de surcroît normal, le groupe $\pi_1(X)$ est le groupe de Galois d'une extension de son corps des fonctions.
\begin{prop}\cite[Prop. 5.4.9]{szamuely} Soit $X$ un schéma normal intègre de corps des fonctions $K$. Soient $\bareta=\Spec \bar K\to X$ un point générique géométrique de $X$, et $K^\sep$ la clôture séparable de $K$ dans $\bar K$. Désignons par $K_X$ la composée dans $K^\sep$ des sous-extensions $L/K$ telles que la normalisation de $X$ dans $L$ soit étale sur $X$. Alors $K_X$ est une extension galoisienne de $K$, et le groupe $\Gal(K_X|K)$ est canoniquement isomorphe à $\pi_1(X,\bareta)$.
\end{prop}

\begin{theorem}\cite[IX, Th. 6.1]{sga1} Soient $k_0$ un corps, et $k$ une clôture algébrique de $k_0$. Soit $k_0^\sep$ la clôture séparable de $k_0$ dans $k$. Soit $X_0$ un $k_0$-schéma de type fini géométriquement intègre. Notons $X\coloneqq X_0\times_{k_0}k$. Soit $\bar x\colon \Spec k\to X$ un point géométrique de $X$. Le morphisme $X\to X_0$ induit une suite exacte de groupes profinis
\[ 1\to \pi_1(X,\bar x)\to \pi_1(X_0,\bar x)\to \Gal(k_0^\sep|k)\to 1.\]
\end{theorem}

\section{Le topos étale}

Introduite par Grothendieck dans les années 1960, la théorie des sites et des topos permet de généraliser la théorie usuelle des faisceaux définis sur la topologie de Zariski d'un schéma. Elle permet de considérer des topologies contenant beaucoup plus d'ouverts que la topologie de Zariski, et donne lieu à des théories cohomologiques plus fines, comme la cohomologie étale.

\subsection{Sites}

\begin{df} Soit $\C$ une catégorie. Une prétopologie sur $\C$ est la donnée d'un ensemble $\Couv(\C)$ de familles de morphismes $(u_i\colon U_i\to U)_{i\in I}$ dans $\C$, appelées familles couvrantes, vérifiant les propriétés suivantes. \begin{enumerate}
\item Si $V\to U$ est un isomorphisme dans $\C$ alors $(V\to U)\in\Couv(\C)$.
\item Si $(U_i\to U)_{i\in I}\in\Couv(\C)$ et pour tout $i\in I$, $(V_{ij}\to U_i)_{j\in J_i}\in \Couv(C)$ alors 
$(V_{ij}\to U)_{i,j}\in\Couv(C)$.
\item Si $(U_i\to U)_{i\in I}\in \Couv(\C)$ et $V\to U$ est un morphisme dans $\C$ alors les produits fibrés $U_i\times_U V$ existent et $(U_i\times_U V\to V)_{i\in I}\in\Couv(\C)$.
\end{enumerate}
Une catégorie munie d'une prétopologie sera appelée un site.
\end{df}
Cette terminologie courante  diffère de celle employée dans \cite{sga41}, où un site est défini comme une catégorie munie d'une \textit{topologie}. Ceci n'a aucune incidence sur la suite.

\begin{df} Soit $X$ un schéma. \begin{enumerate}
\item Un recouvrement de Zariski de $X$ est la donnée d'une famille d'ouverts de Zariski $(U_i)_{i\in I}$ de $X$ telle que $X=\bigcup_{i\in I}U_i$. Le site de Zariski de $X$ est la catégorie des ouverts de $X$ munie de la prétopologie dont les familles couvrantes sont les recouvrements de Zariski.
\item Un recouvrement étale de $X$ est la donnée d'une famille $(u_i\colon U_i\to X)_{i\in I}$ d'éléments de $\Xet$ telle que $\bigcup_i u_i(U_i)=X$. Le (petit) site étale sur $X$ est la catégorie $\Xet$ munie de la prétopologie dont les familles couvrantes sont les recouvrements étales. Le gros site étale sur $X$ est la catégorie des $X$-schémas, munie de la prétopologie dont les familles couvrantes sont les recouvrements étales.
\end{enumerate}
\end{df}

\subsection{Faisceaux et topos}

Soit $\C$ un site. Afin de traiter le sujet des topos en toute généralité et en évitant de rencontrer des problèmes ensemblistes, il est nécessaire de fixer en amont un univers $\mathcal{U}$ \cite[I, §0]{sga41} et ne considérer que des $\mathcal{U}$-sites \cite[II, Def. 3.0.2]{sga41}. Si la catégorie sous-jacente à $\C$ est une catégorie de schémas de type fini sur un schéma fixé, ce qui est le cas du petit site étale sur un schéma, ces complications peuvent être ignorées (voir par exemple la discussion dans \cite[II, §2, p.57]{milneEC}).

\begin{df} Un préfaisceau sur $\C$ est un foncteur $\C^\op\to \Set$. Un préfaisceau $\F$ sur $\C$ est appelé un faisceau si pour tout $U\in \C$ et toute famille couvrante $(U_i\to U)_{i\in I}$, la suite \[ \F(U)\to \prod_{i\in I}\F(U_i)\rightrightarrows \prod_{i,j\in I}\F(U_i\times_UU_j)\]
est exacte. La catégorie des faisceaux sur un site est appelée topos.
\end{df}

On définit de la même façon les (pré)faisceaux de groupes, groupes abéliens, anneaux, modules... sur $\C$. Nous noterons $\Ab(\C)$ la catégorie des faisceaux de groupes abéliens sur $\C$. Soit $X$ un schéma. Soit $\Lambda$ un anneau commutatif. Nous noterons également $\Ab(X)$ (resp. $\Mod_\Lambda(X)$) la catégorie des faisceaux de groupes abéliens (resp. de $\Lambda$-modules) sur le site $\Xet$. 

\begin{rk} \begin{enumerate}
\item \cite[II, Th. 3.4]{sga41} Comme pour la topologie de Zariski, il existe un foncteur de faisceautisation, adjoint à gauche du foncteur d'oubli de la catégorie des faisceaux sur $\C$ vers la catégorie des préfaisceaux sur $\C$.
\item Comme dans la topologie de Zariski, les opérations usuelles sur les préfaisceaux de groupes abéliens (noyau, conoyau, produit, somme directe, produit tensoriel) sont définies en effectuant les opérations concernées sur les groupes de sections. Les opérations correspondantes sur les faisceaux sont obtenues à partir de celles-ci par faisceautisation.
\end{enumerate}
\end{rk}

\begin{ex} Soit $G$ un groupe. Vu comme un groupoïde à un objet, il définit un site $\C_G$. La catégorie des $G$-ensembles est équivalente au topos des faisceaux sur $\C_G$.
\end{ex}

\begin{df} Soient $A,B$ deux topos. Un morphisme de topos $A\to B$ est la donnée d'un couple de foncteurs $(u^\star, u_\star)$ où $u^\star\colon B\to A$ commute aux limites finies et $u_\star\colon A\to B$ est adjoint à droite à $u^\star$.
\end{df}

La catégorie des faisceaux de groupes abéliens sur $\C$ est abélienne et admet suffisamment d'injectifs \cite[II, Prop. 6.7]{sga41}. Lorsque $\C$ est le site étale d'un schéma $X$, nous noterons $\DD(X)$ sa catégorie dérivée, et $\DD^b(X)$ la sous-catégorie pleine de $\DD(X)$ constituée des objets $K$ tels que $\HH^iK$ soit non nul seulement pour un nombre fini d'entiers $i$.
Le foncteur des sections globales \[ \Gamma(\C,-)\colon \Ab(\C)\to\Ab\] associe à un faisceau $\F$ le groupe des morphismes de préfaisceaux d'ensembles du préfaisceau trivial $U\mapsto \{\star\}$ vers $\F$ ; dans le cas particulier où le site $\C$ admet un objet final $X$, $\Gamma(\C,\F)=\F(X)$. Ce foncteur est exact à gauche mais pas à droite ; son foncteur dérivé est noté $\RG(\C,-)$, et les groupes de cohomologie associés sont les $\HH^i(\C,-)=\R^i\Gamma(\C,-)$. Pour tout schéma $X$, nous noterons encore $\RG(X,-)\colon \DD(X)\to \DD(X)$ le foncteur $\RG(\Xet,-)$.

\subsection{Opérations sur les faisceaux}

\begin{df} Soit $f\colon Y\to X$ un morphisme de schémas. Il donne lieu aux foncteurs suivants.\begin{enumerate}
\item Le foncteur image directe $f_\star \colon \Ab(Y)\to \Ab(X)$. Pour tout faisceau $\F\in\Ab(Y)$, $f_\star \F$ est le faisceau  $U\mapsto \F(U\times_X Y)$ sur $Y$.
\item Le foncteur image directe à support propre $f_{!}\colon \Ab(Y)\to \Ab(X)$. Pour tout faisceau $\F\in\Ab(Y)$, le faisceau $f_!\F$ est le sous-faisceau de $f_\star\F$ dont les sections sur $U\to X$ sont les $s\in f_\star\F(U)$ dont le support (c'est-à-dire le plus petit fermé sur le complémentaire duquel faisceau est nul) est propre sur $U$.
\item Le foncteur image inverse $f^\star$. Pour tout faisceau $\G\in \Ab(X)$, $f^\star \G$ est le faisceautisé du préfaisceau $U\mapsto \colim_{V}\F(V)$ sur $X$, où la colimite porte sur les $V\to X$ étales tels que la composée $U\to X$ se factorise par $V\to X$. En particulier, si $f$ est étale, il s'agit simplement de la restriction de $\F$ au site étale de $Y$.
\end{enumerate}
\end{df}

Le foncteur $f^\star$ est adjoint à gauche de $f_\star$. Le foncteur $f^\star$ est exact, et le foncteur $f_\star$ est exact à gauche ; il est également exact à droite lorsque le morphisme $f$ est fini \cite[II, Cor. 3.6]{milneEC}.\\

Soit $\F\in\Ab(X)$. \'{E}tant donné un point géométrique $\bar x\colon \Spec k\to X$, la fibre $\F_{\bar x}$ est définie comme étant $\Gamma(\Spec k,\bar{x}^\star \F)$. Un morphisme de faisceaux $\F\to \G$ est un monomorphisme (resp. épi, resp iso) si et seulement si pour tout point géométrique $\bar x$ de $X$, le morphisme induit $\F_{\bar{x}}\to\G_{\bar{x}}$ en est un.

\subsection{Exemples de faisceaux}
Soit $X$ un schéma. Donnons quelques exemples de faisceaux sur le site étale de $X$. Le premier est le faisceau structural de $X$, qui à $U\in \Xet$ associe $\OO_U(U)$ ; nous le noterons encore $\OO_X$. De même, le préfaisceau $U\mapsto \Hom_X(U,Y)$ représenté par un $X$-schéma $Y$ est un faisceau.

\begin{df} Soit $A$ un groupe abélien. Le faisceau constant associé à $A$, noté $\underline{A}_X$ ou simplement $A$, est le faisceautisé du préfaisceau $U\mapsto A$ ayant pour morphismes de restriction $\id_A$. Son groupe de sections sur $U\in\Xet$ est $\underline{A}_X(U)=A^{\pi_0(U)}$, où $\pi_0(U)$ est l'ensemble des composantes connexes de $U$. Un faisceau sur $X$ est dit constant s'il est isomorphe à un faisceau de cette forme.
\end{df}

\begin{rk}
Le foncteur faisceau constant est adjoint à gauche du foncteur des sections globales $\Gamma(X,-)\colon\F\mapsto \F(X)$.
\end{rk}

\begin{ex} Le groupe multiplicatif $\GG_{m,X}$, représenté par $\Spec \ZZ[t,t^{-1}]\times_{\ZZ}X$, associe à $U\in \Xet$ le groupe $\Gamma(U,\OO_U)^\times$. Le groupe additif $\GG_{a,X}$, représenté par $\Spec \ZZ[t]\times_{\ZZ}X$, associe à $U\in \Xet$ le groupe $\Gamma(U,\OO_U)$.
Soit $n$ un entier.
Le morphisme $[n]\colon\GG_m\to \GG_m$ défini par $t\mapsto t^n$ est un morphisme de schémas en groupes, et son noyau est $\mu_n=\Spec \ZZ[t]/(t^n-1)\times_{\ZZ}X$. 
Le faisceau représenté par $\mu_n$ est \[\mu_{n,X}\colon U\mapsto 
\{ x\in \Gamma(U,\OO_U)\mid x^n=1\}.\]
Remarquons que si $X$ est un schéma sur un corps $k$ qui contient les racines $n$-ièmes de l'unité, le choix d'une telle racine $\zeta\in\mu_n(k)$ détermine un isomorphisme de faisceaux $\ZZ/n\ZZ\to \mu_n$.
\end{ex}

\begin{prop}\cite[IX, 3.2 et 3.5]{sga43}
Soit $n$ un entier inversible sur $X$. Il y a une suite exacte dans $\Ab(X)$, appelée suite exacte de Kummer :
\[ 0\to \mu_{n,X}\to \GG_{m,X}\xrightarrow{[n]} \GG_{m,X}\to 0.\]
Supposons $X$ de caractéristique un nombre premier $p$. Il y a une suite exacte dans $\Ab(X)$, appelée suite exacte d'Artin-Schreier :
\[ 0\to \ZZ/p\ZZ\to \GG_{a,X}\xrightarrow{t\mapsto t^p-t} \GG_{a,X}\to 0.\]
\end{prop}

\begin{df} \begin{enumerate}
\item  Un faisceau sur $X$ est dit localement constant s'il existe un recouvrement étale $(U_i\to X)_{i\in I}$ tel que pour tout $i$, la restriction de $\F$ au site étale de $U_i$ soit un faisceau constant.
\item Un faisceau $\F$ sur $X$ est dit constructible si $X$ est réunion finie de parties localement fermées sur lesquelles $\F$ est localement constant à fibres finies.
\item Soit $\Lambda$ un anneau noethérien. Un faisceau de $\Lambda$-modules sur $X$ sera dit lisse s'il est localement constant et constructible, c'est-à-dire localement constant et à fibres finies.
\end{enumerate}
\end{df}
Nous noterons $\DD^b_c(X)$ la sous-catégorie triangulée pleine de $\DD^b(X)$ des objets $K$ tels que tous les faisceaux $\HH^iK$ soient constructibles. De même, nous noterons $\DD^b_c(\Lambda)$ la catégorie dérivée bornée des $\Lambda$-modules de type fini.

Voici quelques exemples de faisceaux abéliens constructibles sur $X$. 
\begin{itemize}[label=$\bullet$]
\item le faisceau représenté par un $X$-schéma étale \cite[03S8.(1)]{stacks} ;
\item si $f\colon Y\to X$ est un morphisme étale, le faisceau $f_!\Lambda$ \cite[03S8.(3)]{stacks} ;
\item si $i\colon Z\to X$ est une immersion fermée et $\F$ est un faisceau constructible sur $Z$, le faisceau gratte-ciel $i_\star \F$.
\end{itemize}
Remarquons que ce dernier point fournit un exemple de faisceau constructible qui n'est pas représentable par un $X$-schéma étale ; en effet, un tel schéma aurait une image ouverte dans $X$, et donc des sections locales non nulles en tout point d'un ouvert de $X$.

\subsection{Anneaux locaux pour la topologie étale}

L'équivalent des anneaux locaux pour la topologie étale sont les anneaux henséliens.

\begin{df} Un anneau local $(A,\m,k)$ est dit hensélien si pour tout $f\in A[t]$ et toute racine $a_0\in k$ de $\bar f\in k[t]$ telle que $\bar f'(a_0)\neq 0$, il existe $a\in A$ tel que $\bar a=a_0$ et $f(a)=0$. Il est dit strictement hensélien si le corps résiduel $k$ est séparablement clos.
\end{df}

\begin{prop}\cite[Prop. 18.8.8]{ega44} Soit $(A,\m,k)$ un anneau local. Il existe un anneau local strictement hensélien $(A^{hs},\m^{hs},k^{hs})$ muni d'un morphisme d'anneaux locaux $i^{hs}\colon A\to A^{hs}$ tel que tout morphisme d'anneaux locaux de $A$ vers un anneau strictement hensélien $(A',\m',k')$ se factorise par $A^{hs}$, et que cette factorisation est unique si l'on impose le morphisme de corps résiduels $k^{hs}\to k'$. Le couple $(A^{hs},i^{hs})$ est appelé hensélisé strict de $A$.
\end{prop}

Soit $X$ un schéma. Soit $\bar x\colon\Spec k\to X$ un point géométrique de $X$. Un voisinage étale de $\bar x$ est un morphisme étale $u\colon U\to X$ tel que $\bar x$ se factorise par $u$. La colimite des $\Gamma(U,\OO_U)$, où $U$ parcourt les voisinages étales de $X$, est alors le hensélisé strict de l'anneau local $\OO_{X,x}$. Nous noterons $X_{\bar x}$ le schéma $\Spec \OO_{X,x}^{hs}$ ; il est muni d'un morphisme canonique $X_{\bar x}\to X$.

\subsection{Recollement}\label{subsec:recollement}

Soit $X$ un schéma. Soit $j\colon U\to X$ une immersion ouverte. Soit $i\colon Z\to X$ une immersion fermée dont l'image est le complémentaire de l'image de $U$ dans $X$. Définissons une catégorie $\mathcal{C}_{U,Z}$ de la façon suivante. Ses objets sont les triplets $(\F_U,\F_Z,\phi)$ où $\F_U\in \Ab(U)$, $\F_Z\in\Ab(Z)$ et $\phi \colon \F_Z\to i^\star j_\star\F_U$. Les morphismes $(\F_U,\F_Z,\phi)\to (\F'_U,\F'_Z,\phi')$ sont les couples de morphismes $(\psi_U\colon\F_U\to\F'_U,\psi_Z\colon\F_Z\to\F'_Z)$ tels que le diagramme suivant soit commutatif.

\[\begin{tikzcd}
\F_Z \arrow[r,"\psi_Z"]\arrow[d,"\phi"] & \F'_Z \arrow[d,"\phi'"] \\
\F_U \arrow[r,"\psi_U"] & \F'_U
\end{tikzcd}\]

Remarquons que pour tout faisceau $\F$ sur $X$, l'adjonction $j^\star \dashv j_\star$ fournit un morphisme \[ \phi_\F\colon\F\to j_\star j^\star\F.\]

\begin{prop}\cite[II, Th. 3.10]{milneEC}
Le foncteur $\Ab(X)\to \mathcal{C}_{U,Z}$ qui à un faisceau $\F$ associe $(j^\star\F,i^\star\F,\phi_\F)$ et à un morphisme $f\colon \F\to \F'$ associe $(j^\star f,i^\star f)$ est une équivalence de catégories. Un quasi-inverse est donné par $(\F_U,\F_Z,\phi)\mapsto \F$, où $\F$ est défini par le diagramme cartésien suivant :
\[
\begin{tikzcd}
\F \arrow[r]\arrow[d] & i_\star \F_Z \arrow[d,"i_\star \phi"] \\
j_\star \F_U\arrow[r,"i^\star\dashv i_\star",swap] & i_\star i^\star j_\star \F_U
\end{tikzcd}
\]
\end{prop}

\begin{df} L'équivalence de catégories précédente permet de définir les foncteurs \[ j_!\colon\Ab(U)\to \Ab(X), \,\F\mapsto (\F,0,0) \] et \[i^!\colon \Ab(X)\to \Ab(Z), \,(\F_U,\F_Z,\phi)\mapsto\ker(\phi).\]
\end{df}

\begin{rk}
Le foncteur $j_!$ coïncide avec le foncteur image directe à support propre défini précédemment.
\end{rk}

\begin{prop} Les foncteurs $i_\star,i^\star,i^!,j_\star,j_!,j^\star$ vérifient les adjonctions
\[ i^\star \dashv i_\star \dashv i^! \]
et \[ j_!\dashv j^\star \dashv j_\star.\]
Les foncteurs $i^\star,i_\star,j^\star,j_!$ sont exacts. Les foncteurs $i^!$ et $j_\star$ sont exacts à gauche. Pour tout faisceau $\F$ de groupes abéliens sur $X$, il y a des suites exactes courtes
\[ 0\to j_!j^\star\F\to \F\to i_\star i^\star\F\to 0\]
et \[ 0\to i_\star i^!\F\to \F\to j_\star j^\star\F.\]
\end{prop}

\section{Torseurs, \texorpdfstring{$\HH^1$}{H\textsuperscript{1}} et fibrés en droites}

\subsection{Torseurs et premier groupe de cohomologie}
Soit $\C$ un site. Soit $\mathcal{G}$ un faisceau en groupes abéliens sur $\C$. 

\begin{df}
Un faisceau $\F$ sur $\C$ est un $\mathcal{G}$-torseur s'il est muni d'une action à gauche $\mathcal{G}\times \F\to\F$ qui est localement (pour la topologie de $\C$) isomorphe à l'action par translation $\G\times\G\to \G$.
\end{df}

Soit $\G\to \I$ un monomorphisme de $\G$ vers un faisceau injectif. Alors le morphisme \[ \HH^0(\C,\I/\G)\to \HH^1(\C,\G)\]
est surjectif. Soit $c\in \HH^1(\C,\G)$ : il est l'image d'un élément $c'\in \HH^0(\C,\I/\G)$. Considérons le faisceau $\F'\subset \I$ des antécédents de $c'$.
\begin{prop}\cite[03AJ]{stacks} L'application $\F\mapsto \F'$ définit une bijection canonique entre l'ensemble $\HH^1(\C,\G)$ et l'ensemble des classes d'isomorphisme de $\G$-torseurs sur $\C$.
\end{prop}

\begin{df} Soient $\F_1,\F_2$ deux faisceaux sur $\C$ munis d'une action à gauche de $\G$. Le produit contracté $\F_1\wedge^\G\F_2$ est le quotient du faisceau $\F_1\times \F_2$ par la relation d'équivalence définie par $(g\cdot f_1,f_2)=(f_1,g\cdot f_2)$.
\end{df}

Si $\T_1$ et $\T_2$ sont deux $\G$-torseurs, le produit contracté $\T_1\wedge^\G\T_2$ est encore un $\G$-torseur ; le produit contracté définit une loi de groupe sur l'ensemble des classes d'isomorphisme de $\G$-torseurs sur $X$, qui correspond à la loi de groupe sur $\HH^1(\C,\G)$ via la bijection ci-dessus \cite[III.4, Rem. 4.8(b)]{milneEC}. 

\subsection{Torseurs sur le site étale}

Soit $X$ un schéma. 
\begin{df} \begin{enumerate}
\item Soit $G$ un $X$-schéma en groupes. Soit $T$ un $X$-schéma muni d'une action à droite de $G$. On dit que $T$ est un $G$-torseur sur $X$ si $T\to X$ est étale et surjectif, et si le morphisme $G\times_X T\to T\times_X T, (g,t)\mapsto (t,tg)$ est un isomorphisme.
\item Soit $G$ un groupe abélien. Un faisceau $\F$ sur $X$ est un $G$-torseur si c'est un torseur sous le faisceau constant associé à $G$. Un $X$-schéma $T$ est un $G$-torseur si c'est un torseur sous le $X$-schéma en groupes $G\times X$.
\end{enumerate}
\end{df}

Le résultat de représentabilité élémentaire suivant suffira dans notre cas.

\begin{prop}\cite[Prop. 5.7.18]{fulei} Soit $X$ un schéma. Soit $G$ un schéma en groupes séparé étale de présentation finie sur $X$. Le foncteur de Yoneda $T\mapsto h_T$ induit une équivalence de catégories entre la catégorie des $X$-schémas qui sont des $G$-torseurs et la catégorie des faisceaux sur $X$ qui sont des $h_G$-torseurs.
\end{prop}

Les torseurs sous un groupe fini généralisent les revêtements galoisiens.

\begin{lem}\cite[Prop. 5.3.16]{szamuely} Soit $X$ un schéma. Soit $G$ un groupe fini. Les $G$-torseurs connexes sur $X$ sont les revêtements galoisiens de $X$ de groupe $G$.
\end{lem}

\begin{rk}[Fonctorialité] Voici comment se décrivent différents morphismes en termes de torseurs. \begin{itemize}[label=$\bullet$]
\item \'{E}tant donné un morphisme de faisceaux de groupes abéliens $u\colon\F\to\G$ sur $X$, le morphisme $u_\star \colon \HH^1(X,\F)\to \HH^1(X,\G)$ obtenu par fonctorialité de $\HH^1(X,-)$ associe à un $\F$-torseur $\T$ le $\G$-torseur $\T\wedge^\F \G$, où $\F$ agit sur $\G$ via $f\cdot g=u(f)+g$.
\item \'{E}tant donné un morphisme de schémas $\phi\colon Y\to X$, le morphisme $\phi^\star \colon \HH^1(X,\F)\to \HH^1(Y,f^\star\F)$ associe à un $\F$-torseur $\T$ le $\phi^\star\F$-torseur $\phi^\star \T$. 
\item \'{E}tant donné une suite exacte de faisceaux de groupes abéliens sur $X$ \[ 0\to \F \xrightarrow{u} \G\xrightarrow{v}\mathcal{H}\to 0 \]
le morphisme canonique $\partial \colon \HH^0(X,\mathcal{H})\to \HH^1(X,\F)$ est construit de la façon suivante \cite[Cycle, 1.1.4]{sga412}. \`{A} une section globale $s\in \HH^0(X,\mathcal{H})$, il associe le faisceau \[ v^{-1}s\colon U\mapsto \{ t\in \G(U)\mid  v_U(t)=s|_U \}. \]
La structure de $\F$-torseur sur $\partial s\coloneqq v^{-1}s$ est donnée par $f\cdot t=u(f)+t$. 
\end{itemize}
\end{rk}

Mentionnons enfin un résultat qui servira par la suite, concernant la restriction des torseurs.
\begin{lem}\label{lem:torsres} Soit $j\colon U\to X$ une immersion ouverte. Soit $\F$ un faisceau de groupes abéliens sur $U$. Alors le morphisme de restriction $\HH^1(X,j_\star\F)\to \HH^1(U,\F)$ est injectif.
\begin{proof}
Soit $\T$ un $j_\star\F$-torseur sur $X$ tel qu'il y ait un isomorphisme de $j^\star j_\star\F=\F$-torseurs $\phi\colon j^\star\T\to j^\star j_\star\F$. Comme tout morphisme de torseurs est un isomorphisme, il suffit d'exhiber un morphisme de $j_\star\F$-torseurs $\T\to j_\star\F$, ce qui est simple : la composée
\[ \T\to j_\star j^\star\T \xrightarrow{j_\star\phi} j_\star j^\star j_\star\F \simeq j_\star \F\]
convient. En effet, le diagramme
\[
\begin{tikzcd}
j_\star \F\times \T \arrow[d]\arrow[r] 
& j_\star j^\star j_\star \F\times j_\star j^\star\T \arrow[d]\arrow[r] 
&j_\star j^\star j_\star \F \times j_\star j^\star j_\star \F \arrow[d]\arrow[r] 
& j_\star \F\times j_\star \F \arrow[d] \\
 \T \arrow[r] & j_\star j^\star\T \arrow[r] & j_\star j^\star j_\star \F \arrow[r] & j_\star \F 
\end{tikzcd}
\]
est commutatif.
\end{proof}
\end{lem}

\subsection{Torseurs sous $\GG_m$ et $\mu_n$}

Soit $X$ un schéma. Fixons un entier naturel non nul $n$. Dans un premier temps, nous allons considérer les torseurs sous $\GG_m$ et $\mu_n$ sur $X$.  Ils admettent une description en termes de faisceaux inversibles sur $X$. Remarquons qu'un faisceau inversible pour la topologie étale sur $X$ définit par restriction un faisceau inversible pour la topologie de Zariski. Réciproquement, un $\OO_X$-module inversible $\LL$ pour la topologie de Zariski définit un $\OO_X$-module dont les sections sur $(U\xrightarrow{u} X)\in\Xet$ sont données par $\Gamma(U,u_{Zar}^\star\LL)$ \cite[03DV]{stacks}. Ceci définit une équivalence entre les catégories de $\OO_X$-modules inversibles pour les topologies de Zariski et étale ; il n'y a donc pas lieu de faire une distinction entre les deux.

\begin{prop}\cite[040D]{stacks} \'{E}tant donné un faisceau inversible $\LL$ sur $X$, le faisceau sur $\Xet$\[ \underline\Isom(\OO_X,\LL)\colon U\mapsto \Isom_U(\OO_U,\LL_U)\] est un $\GG_m$-torseur.
Le morphisme $\LL\mapsto \underline\Isom(\OO_X,\LL)$ définit un isomorphisme $\Pic X\to \HH^1(X,\GG_m)$ fonctoriel en $X$. L'isomorphisme inverse associe à un $\GG_m$-torseur $\T$ le faisceau inversible $\mathcal{T}\wedge^{\GG_m} \OO_X$.
\end{prop}

Considérons désormais la catégorie $\C$ dont les objets sont les couples $(\LL,\alpha)$, où $\LL$ est un faisceau inversible sur $X$ et $\alpha \colon \LL^{\otimes n}\xrightarrow{\sim} \OO_X$ est une trivialisation de $\LL^{\otimes n}$. Un morphisme entre deux couples $(\LL,\alpha),(\LL',\alpha')$ est défini comme étant un morphisme $\phi \colon\LL\to\LL'$ tel que le diagramme \[
\begin{tikzcd}
\LL^{\otimes n} \arrow[r,"\alpha"]\arrow[d,"\phi",swap] & \OO_X \arrow[d,"\id"] \\
\LL'^{\otimes n} \arrow[r,"\alpha'"] & \OO_X
\end{tikzcd}\]
soit commutatif. Soit $S$ l'ensemble des classes d'isomorphisme de tels couples. Le produit tensoriel $(\LL,\alpha)\otimes(\LL',\alpha')\coloneqq(\LL\otimes\LL',\alpha\otimes\alpha')$ définit une loi de groupe sur $S$, d'élément neutre $(\OO_X,\id)$. 

\begin{prop}\cite[040Q]{stacks}\label{prop:H1faiscinv} \'{E}tant donné un objet $(\LL,\alpha)$ de $\C$, le faisceau \[\T_\LL\colon U\mapsto \Isom_\C((\OO_U,1),(\LL|_U,\alpha|_U)) \] sur $\Xet$ est un $\mu_n$-torseur. L'application $\LL\mapsto \T_\LL$ définit un isomorphisme de groupes $S\to \HH^1(X,\mu_n)$ fonctoriel en $X$.
\end{prop}

\section{Groupe fondamental et faisceaux lisses}

\subsection{Faisceaux lisses et $\pi_1$-modules}

Soient $X$ un schéma et $\bar x$ un point géométrique de $X$. Soit $\Lambda$ un anneau fini.

\begin{prop}\cite[Prop. 5.8.1.(i)]{fulei} Soit $\F$ un faisceau à fibres finies sur $X$. Le faisceau $\F$ est localement constant si et seulement s'il est représenté par un revêtement étale de $X$. S'il l'est, il existe un morphisme fini étale surjectif $f\colon Y\to X$, avec $Y$ connexe, tel que $f^\star\F$ soit constant.
\end{prop}
En d'autres termes, le foncteur de Yoneda $Y\mapsto h_Y$ définit une équivalence entre la catégorie $\Fet_X$ des revêtements étales de $X$ et la catégorie des faisceaux localement constants constructibles sur $X$. Nous connaissons déjà une autre catégorie équivalente à $\Fet_X$ : celle des $\pi_1(X,\bar x)$-modules. Nous serons intéressés par le cas particulier des faisceaux de $\Lambda$-modules ; le résultat s'énonce alors de la façon suivante.

\begin{prop}\cite[V, Th. 4.1]{sga1} Le foncteur fibre en $\bar x$ définit une équivalence de catégories entre la catégorie des faisceaux lisses de $\Lambda$-modules sur $X$ et la catégorie des $\pi_1(X,\bar x)$-modules de type fini munis d'une structure de $\Lambda$-module qui commute à l'action de $\pi_1(X,\bar x)$.
\end{prop}

Soit $\F$ un faisceau lisse de $\Lambda$-modules sur $X$. La proposition précédente montre qu'il est uniquement déterminé par la donnée du $\pi_1(X,\bar x)$-ensemble $\F_{\bar x}$. Un faisceau localement constant $\T$ sur $X$ muni d'une action à droite de $\F$ s'identifie alors au groupe $\T_{\bar x}$, muni d'une action à droite de $\F_{\bar x}$ et d'une action continue à gauche de $\pi_1(X,\bar x)$, qui sont compatibles au sens où pour tous $s\in \pi_1(X,\bar x), t\in T_{\bar x}$ et $f\in \F_{\bar x}$ : \[ s(t\cdot f)=(st)\cdot (sf).\]
Le faisceau $\T$ est un $\F$-torseur si et seulement si $\T_{\bar x}$ est un $\F_{\bar x}$-torseur dans la catégorie des $\pi_1(X,\bar x)$-ensembles \cite[XI, §5, p231]{sga1}. Il s'en déduit un isomorphisme canonique, fonctoriel en $\F$ :
\[ \begin{array}{rcl }\HH^1(X,\F)&\xrightarrow{\sim}&\HH^1(\pi_1(X,\bar x),\F_{\bar x})\\
\T &\mapsto &\T_{\bar x}.\end{array}\]

\paragraph{Opérations sur les faisceaux lisses} Soit $f\colon Y\to X$ un revêtement galoisien de schémas connexes. Fixons des points géométriques $\bar x,\bar y$ de $X$ et $Y$ tels que $\bar x=f\circ\bar y$. 
Si $\F$ est un faisceau lisse de groupes abéliens sur $X$ de fibre $\F_{\bar x}=M$, le faisceau lisse $f^\star \F$ a pour fibre $\F_{\bar y}=M$. L'adjoint à droite de $f^\star$ est $f_\star$, qui correspond donc à l'adjoint à droite du foncteur d'oubli $\Mod_{\pi_1(X,\bar x)}\to \Mod_{\pi_1(Y,\bar y)}$, qui est la co-induction. Dans les catégories de revêtements étales de $X$ et de $Y$, il correspond à la restriction de Weil $R_{Y\to X}$ (voir annexe \ref{sec:weilres}). Le tableau suivant résume la situation :\\
\begin{center}
\begin{tabular}{|c|c|c|}\hline
Revêtements étales & Faisceaux lisses & $\pi_1$-modules  \\ \hline
$-\times_X Y$ & $f^\star$ & oubli $M\mapsto M$ \\
\hline
$R_{Y\to X}$ & $f_\star$ & $\coind_{\pi_1(X,\bar x)}^{\pi_1(Y,\bar y)}$ \\ \hline
\end{tabular} 
\end{center}

\paragraph{Revêtement trivialisant minimal} Soit $X$ un schéma. Soit $\bar x$ un point géométrique de $X$. Soit $\Lambda$ un anneau noethérien. Soit $\F$ un faisceau lisse de $\Lambda$-modules sur $X$, correspondant à un $\pi_1(X,\bar x)$-module $M$. Notons $\mathfrak{S}$ le groupe de monodromie associé : c'est l'image de $\pi_1(X,\bar x)$ dans $\Aut_\Lambda(M)$.
Soit $H$ le noyau du morphisme $\pi_1(X,\bar x)\to M$. Il correspond à un revêtement galoisien $X_{\min}\to X$ de groupe d'automorphismes $\pi_1(X,\bar x)/H\xrightarrow{\sim}\mathfrak{S}$ ; ce revêtement est minimal pour la propriété de trivialiser $\F$.

\subsection{$G$-faisceaux et descente galoisienne}

\begin{df} Soient $X$ un schéma, $\F$ un faisceau sur $X$, et $G$ un groupe d'automorphismes de $X$. Une action de $G$ sur $\F$ est une famille d'isomorphismes $(\phi_g\colon \F\to g_\star\F)_{g\in G}$ telle que $\phi_{e_G}=\id_\F$ et que pour tous $g,h\in G$, le diagramme 

\[
\begin{tikzcd}
\F \arrow[r,"\phi_g"]\arrow[d,"\phi_{gh}",swap] & g_\star\F \arrow[d,"g_\star\phi_h"] \\
(gh)_\star\F \arrow[r,"\sim"] & g_\star h_\star \F
\end{tikzcd}
\]
où la ligne du bas est l'isomorphisme canonique, soit commutatif. Un faisceau $\F$ muni d'une action de $G$ sera appelé $G$-faisceau. Un morphisme de $G$-faisceaux est un morphisme de faisceaux compatible à l'action de $G$.
\end{df}

Considérons pour la suite un revêtement galoisien $f\colon Y\to X$ de groupe $G$. Si $\F$ est un faisceau sur $X$, alors le faisceau $f^\star\F$ est muni d'une action de $G$. \'{E}tant donné un élément $g\in G$, la composition du morphisme d'adjonction $f^\star \F\to g_\star g^\star f^\star \F$ avec l'isomorphisme canonique $g_\star g^\star f^\star \F \to g_\star (f\circ g)^\star\F= g_\star f^\star \F$ fournit un morphisme $f^\star\F \to g_\star f^\star\F$.

\begin{prop}\cite[0GEZ, 0CDQ]{stacks}\label{prop:descgal}  Le foncteur $\F\mapsto f^\star\F$ définit une équivalence entre la catégorie des faisceaux sur $X$ et la catégorie des $G$-faisceaux sur $Y$. Un quasi-inverse de ce foncteur est donné par $\F\mapsto (f_\star \F)^G$.
\end{prop}

Soit $\G$ un faisceau de groupes abéliens sur $X$. Notons $\G'$ le $G$-faisceau $f^\star\G$. Soit $\F$ un $\G$-torseur sur $X$. Le $G$-faisceau $\F'\coloneqq f^\star\F$ est encore un $\G'$-torseur sur $Y$ ; le morphisme $\G'\times \F'\to \F'$ est un morphisme de $\G$-faisceaux, au sens où le diagramme suivant est commutatif pour tout $g\in G$ :
\[
\begin{tikzcd} \G'\times \F' \arrow[d]\arrow[r] & \F' \arrow[d] \\
g_\star\G'\times g_\star\F' \arrow[r] & g_\star \F'
\end{tikzcd}\]
Réciproquement, soit $\F'$ un $G$-faisceau sur $Y$ muni d'une action de $\G'$ qui en fait un $\G'$-torseur. Soit $\F$ le faisceau $(f_\star\F')^G$ sur $X$ ; il vérifie $f^\star\F=\F'$. Le morphisme $\G'\times \F'\to \F'$ descend en un morphisme $\G\times \F\to \F$ si et seulement si l'action de $\G'$ sur $\F'$ est $G$-équivariante au sens ci-dessus.

\begin{cor}\label{cor:torsrev} Soit $\mathcal{G}$ un faisceau de groupes abéliens sur $X$. Le foncteur $\F\mapsto f^\star\F$ définit une équivalence entre la catégorie des $\G$-torseurs sur $X$ et celle des $f^\star\G$-torseurs sur $Y$ munis d'une action de $G$ telle que l'action de $f^\star\G$ soit $G$-équivariante.
\end{cor}

\subsection{Schémas $K(\pi,1)$}\label{subsec:Kpi1}

Sur certains schémas, la cohomologie des faisceaux lisses est isomorphe à la cohomologie galoisienne de leur fibre : les schémas qui vérifient cette propriété sont appelés schémas $K(\pi,1)$. Une étude détaillée de leurs propriétés se trouve par exemple dans \cite{achinger_phd}.\\

Soient $X$ un schéma connexe et $\bar x$ un point géométrique de $X$. Soit $\ell$ un nombre premier inversible sur $X$. Notons $\Xfet$ le topos fini étale associé à $X$, défini par la sous-catégorie pleine des $X$-schémas finis étales munie de la topologie étale. L'équivalence de catégories donnée par le théorème \ref{th:galsch} définit un isomorphisme de topos entre $\Xfet$ et le topos $B_{\pi_1(X,\bar x)}$ des $\pi_1(X,\bar x)$-ensembles finis continus.

\begin{df}\label{df:lm} Un faisceau d'ensembles $\F$ sur $X$ est dit $\ell$-monodromique s'il est localement constant constructible et si l'image de $\pi_1(X,\bar x)$ dans $\Aut(\F_{\bar x})$ est un $\ell$-groupe. Un faisceau de groupes $\F$ sur $X$ est dit $\ell$-monodromique si ses fibres sont des $\ell$-groupes finis et s'il est $\ell$-monodromique en tant que faisceau d'ensembles.
\end{df}

Notons encore $\Xlet$ le topos défini par la sous-catégorie pleine des schémas finis étales $Y\to X$ tels que le faisceau représenté par $Y$ soit $\ell$-monodromique. Considérons les morphismes de topos $\rho\colon\Xet\to\Xfet$ et $\rho_\ell\colon\Xet\to\Xlet$ définis par la restriction des faisceaux aux sites définissant $\Xfet$ et $\Xlet$. Notons enfin $\pi_1(X,\bar x)^\ell$ le complété pro-$\ell$ de $\pi_1(X,\bar x)$.

\begin{df}\cite[Def. 9.21]{abbes_gros} Le schéma $X$ est appelé un $K(\pi,1)$ si pour tout entier $n$ inversible sur $X$ et tout faisceau $\F$ de $\ZZ/n\ZZ$-modules sur $\Xfet$, le morphisme d'adjonction \[ \F \to \R\rho_\star(\rho^\star\F)\]
est un isomorphisme. 
\end{df}

\begin{df}\cite[1.4.4]{mo}
Le schéma $X$ est un $K(\pi,1)$ pro-$\ell$ si pour tout faisceau abélien $\F$ de $\ell$-torsion sur $\Xlet$, le morphisme d'adjonction \[ \F\to \R{\rho_\ell}_\star(\rho_\ell^\star\F) \]
est un isomorphisme.
\end{df}

\begin{lem}\cite[1.4.2]{mo} Le schéma $X$ est un $K(\pi,1)$ si et seulement si, pour tout entier $n$ inversible sur $X$ et tout faisceau lisse $\F$ de $\ZZ/n\ZZ$-modules, le morphisme \[ \RG(\pi_1(X,\bar x),\F_{\bar x})\to\RG(X,\F) \]
est un isomorphisme. Le schéma $X$ est un $K(\pi,1)$ pro-$\ell$ si et seulement si, pour tout faisceau abélien $\ell$-monodromique $\F$ sur $X$, le morphisme \[ \RG(\pi_1(X,\bar x)^\ell,\F_{\bar x})\to \RG(X,\F)\]
est un isomorphisme.
\end{lem}

Si $X$ est un $K(\pi,1)$, la cohomologie d'un faisceau lisse $\F$ sur $X$ peut donc se calculer à l'aide de la cohomologie du $\pi_1(X,\bar x)$-module $\F_{\bar x}$. 
Le résultat d'effaçabilité suivant fournit une condition suffisante pour qu'un schéma soit un $K(\pi,1)$.

\begin{prop}\cite[Prop. A.3.1]{stix_thesis}\label{prop:condKpi1} Soient $X$ un schéma connexe et $\ell$ un nombre premier inversible sur $X$. Supposons que pour tout $i\geqslant 1$ et tout faisceau lisse $\F$ de torsion inversible sur $X$ (resp. de $\ZZ/\ell\ZZ$-espaces vectoriels), il existe un revêtement $\phi_i\colon Y_i\to X$ galoisien (resp. galoisien de groupe un $\ell$-groupe) tel que le morphisme $\HH^i(X,\F)\to \HH^i(Y_i,\phi_i^\star\F)$ soit nul. Alors $X$ est un $K(\pi,1)$.
\begin{proof} Un $\partial$-foncteur cohomologique effaçable étant universel \cite[Prop. 2.2.1]{tohoku}, il suffit de montrer que pour tout $i\geqslant 1$, le foncteur $\HH^i(X,\rho^\star-)\colon B_{\pi_1(X,\bar x)}\to \Ab$ est effaçable. Soit $M$ un $\pi_1(X,\bar x)$-module fini de torsion inversible sur $X$, et $\F=\rho^\star M$ le faisceau lisse associé. Soit $i\geqslant 1$. La condition de l'énoncé fournit un revêtement galoisien $\phi\colon Y\to X$ tel que $\HH^i(X,\F)\to \HH^i(Y,\phi^\star\F)$ soit nul.
Le morphisme $\F\to \phi_\star\phi^\star\F$ est injectif. Soit $\bar y$ un point géométrique de $Y$ d'image $\bar x$. Le faisceau lisse $\phi_\star\phi^\star\F$ correspond au $\pi_1(X,\bar x)$-module $N\coloneqq \coind_{\pi_1(X,\bar x)}^{\pi_1(Y,\bar y)}M$. Enfin, $\HH^1(X,\phi_\star\phi^\star\F)=\HH^1(Y,\F)$ par exactitude de $\phi_\star$. Il y a donc une injection $M\to N$ dans $B_{\pi_1(X,\bar x)}$ telle que $\HH^1(X,\rho^\star M)\to \HH^1(X,\rho^\star N)$ soit nulle.
\end{proof}
\end{prop}

\begin{rk} Sur un corps algébriquement clos, pour tout entier $m\geqslant 1$, $\pi_1(\PP^m)=0$ \cite[XI, Prop. 1.1]{sga1} mais $\HH^2(\PP^m,\Lambda)=\Lambda(-1)$ \cite[Example 16.3]{milneLEC}. Par conséquent, $\PP^m$ n'est pas un $K(\pi,1)$. Nous verrons dans la section \ref{subsec:courbKpi1} qu'à l'exception de $\PP^1$, les courbes lisses sont toutes des $K(\pi,1)$. Il existe également des résultats positifs en dimension quelconque, sur les corps finis : Achinger a montré \cite[Th. 1.1.1]{achinger_Kpi1} que tout $\FF_p$-schéma affine connexe est un $K(\pi,1)$.
\end{rk}

\begin{rk} Si $X$ est un $K(\pi,1)$ alors pour tout complexe $K$ de $\pi_1(X,\bar x)$-modules, le morphisme \[ \RG(\pi_1(X,\bar x),K)\to \RG(X,\rho^\star K)\]
est encore un isomorphisme. Notons $\mathscr{K}=\rho^\star K$. Considérons le morphisme entre les suites spectrales  $E_{2,\pi}^{pq}\coloneqq \HH^p(\pi,\HH^qK)$ et $E_{2,X}^{pq}\coloneqq \HH^p(X,\HH^q\mathscr{K})$ données par \cite[015J]{stacks}. Il est un isomorphisme dès la deuxième page puisque $X$ est un $K(\pi,1)$, et il induit donc des isomorphismes entre les aboutissements $\HH^i(\pi,K)\to \HH^i(X,\mathscr{K})$ par \cite[Th. 5.2.12]{weibel}. Par conséquent, le morphisme $\RG(\pi,K)\to\RG(X,\mathscr{K})$ est un isomorphisme dans $\DD^b_c(\Lambda)$.
\end{rk}

\subsection{Morphismes de restriction}

Cette section recense quelques lemmes utiles par la suite.

\begin{lem}
Soit $U'$ un ouvert non vide d'un schéma $U$ intègre normal. Notons ${j} \colon U'\to U$ l'inclusion. Si $\F$ est un faisceau constant de $\Lambda$-modules sur $U'$ de fibre $F$ alors $j_\star \F$ est constant sur $U$ de fibre $F$.
\begin{proof} Soit $f\colon T\to U$ un morphisme étale. Alors $T$ est encore un schéma normal \cite[025P]{stacks}. Si $T$ est connexe alors il est intègre \cite[033M]{stacks}, et son ouvert $T\times_X U'$ (qui est non vide puisque $U$' est dense) est encore irréductible, donc connexe. Par conséquent, $j_\star\F(T)=F$. Les morphismes de restriction de $j_\star\F$ sont évidemment les mêmes que ceux de $\F$.
\end{proof}
\end{lem}

\begin{lem}\label{adjiso}
Soient $U$ un schéma intègre normal, et $U'$ un ouvert non vide de $U$.  Notons ${j} \colon U'\to U$ l'inclusion. Soit $\F$ un faisceau lisse de $\Lambda$-modules sur $U$. Alors le morphisme d'adjonction $\F\to {j}_\star{j}^\star \F$ est un isomorphisme.
\begin{proof} Il existe un morphisme fini étale $p\colon V\to U$ qui trivialise $\F$, avec $V$ un schéma intègre normal. Considérons le diagramme cartésien suivant.
\[
\begin{tikzcd}
V' \arrow[r,"{j'}"]\arrow[d,"q",swap] & V \arrow[d,"p"]\\
U' \arrow[r,"{j}"] & U
\end{tikzcd}
\]
 Alors $V'$ est un ouvert non vide de $V$. Considérons le faisceau constant $F={j'}^\star p^\star \F$. Par le lemme précédent, ${j'}_\star F$ est constant, et par conséquent le morphisme $p^\star\F \to {j'}_\star{j'}^\star p^\star\F={j'}_\star q^\star {j}^\star\F=p^\star{j}_\star{j}^\star\F$ (par finitude de $p$) est un isomorphisme. Comme le morphisme $p$ est fidèlement plat et quasi-compact, le foncteur $p^\star$ reflète les isomorphismes \cite[Prop. 2.7.1]{ega42}, et par conséquent $\F\to {j}_\star{j}^\star\F$ est un isomorphisme.
\end{proof}
\end{lem}

\begin{cor}\label{cor:restiso} Soit $X$ un schéma intègre normal. Soient $U'\subset U$ deux ouverts de Zariski non vides de $X$. Soit $\F$ un faisceau lisse de $\Lambda$-modules sur $X$. Le morphisme de restriction $\F(U)\to\F(U')$ est un isomorphisme.
\end{cor}

\section{Les grands théorèmes}

\subsection{Dimension cohomologique}

\begin{df}
Soit $\ell$ un nombre premier. La $\ell$-dimension cohomologique d'un schéma $X$ est le plus petit $j\in \mathbb{N}\cup\{\infty\}$ tel que pour tout faisceau abélien $\F$ de $\ell$-torsion sur $X$ et tout $i>j$, $\HH^i(X,\F)=0$. Elle sera notée $\cd_\ell(X)$. La dimension cohomologique de $X$ est \[ \cd(X)\coloneqq \sup_\ell \cd_\ell(X)\in\NN\cup \{ \infty\}.\]
\end{df}

\begin{prop}\cite[X, Cor. 4.3]{sga43} Soient $p,\ell$ deux nombres premiers distincts. Soient $k_0$ un corps de caractéristique $p$ et $X$ un schéma affine de type fini sur $k_0$. Alors \[ \cd_\ell(X)\leqslant \cd_\ell(k_0)+2\dim(X)\]
et \[ \cd_p(X)\leqslant \dim(X)+1.\]
\end{prop}

\begin{prop}\cite[X, Th. 5.1]{sga43} Soit $X$ un schéma affine de type fini sur un corps $k_0$. Alors 
\[ \cd(X)\leqslant \dim(X)+\cd(k_0).\]
\end{prop}

\begin{rk}\begin{enumerate}
\item Lorsque $X$ est le spectre d'un corps $k_0$, sa dimension cohomologique est la dimension cohomologique de $k_0$ pour la cohomologie galoisienne. En particulier, les corps de dimension cohomologique nulle sont les corps séparablement clos. Les corps finis ainsi que les corps de fonctions de courbes sur un corps algébriquement clos sont de dimension cohomologique $1$.
\item Si $X$ est une courbe sur un corps séparablement clos et $\F$ est un faisceau abélien de torsion sur $X$ alors $\HH^i(X,\F)=0$ dès que $i\geqslant 3$. Si $X$ est de surcroît affine, $\HH^2(X,\F)$ est également nul.
\end{enumerate}
\end{rk}

\subsection{Invariance topologique}\label{subsec:invtop}

Nous nous intéresserons par la suite uniquement au calcul de la cohomologie de schémas réduits sur des corps parfaits ; les résultats suivants justifient ces restrictions.

\begin{theorem}\cite[VIII, Th. 1.1]{sga42}\label{th:homeouniv} Soit $f\colon Y\to X$ un morphisme de schémas. Si $f$ est un homéomorphisme universel alors les foncteurs $-\times_X Y\colon \Xet\to\Yet$ et $f^\star\colon\Ab(X)\to\Ab(Y)$ qu'il induit sont des équivalences de catégories. 
\end{theorem}
Supposons $X$ de dimension cohomologique finie. Avec ces notations, pour tout $K\in \DD^b(X)$, le morphisme $\RG(X,K)\to \RG(Y,f^\star K)$ est un quasi-isomorphisme. Les homéomorphismes universels étant exactement les morphismes entiers surjectifs et radiciels, ce théorème fournit les résultats ci-dessous.

\begin{cor}\begin{enumerate} \item Soit $X$ un schéma de dimension cohomologique finie, et $X_\red$ son réduit. Le morphisme $X_\red\to X$ est un homéomorphisme universel \cite[054M]{stacks} et induit donc pour tout $K\in \DD^b(X)$ un quasi-isomorphisme \[\RG(X,K)\xrightarrow{\sim} \RG(X_\red,K|_{X_\red}).\]
\item Soit $X$ un schéma sur un corps $k$. Soit $k'\to k$ une extension purement inséparable. Notons $X'=X\times_k k'$. Le morphisme $X'\to X$ est un homéomorphisme universel  et induit pour tout $K\in \DD(X)$ un quasi-isomorphisme
\[ \RG(X,K)\xrightarrow{\sim} \RG(X',K|_{X'}).\]
En particulier, cela s'applique au cas où $k'$ est une clôture parfaite de $k$.
\end{enumerate}
\end{cor}

\subsection{Cohomologie de \v{C}ech et triangle de Mayer-Vietoris}\label{subsec:mayerviet}

Soient $X$ un schéma, et $\mathcal{U}=(U_i\xrightarrow{f_i} X)_{i\in I}$ un recouvrement étale fini de $X$. Considérons le schéma $U=\coprod_{i\in I}U_i\to X$. \'{E}tant donné $i_1,\dots,i_r\in I$, notons encore $f_{i_1,\dots,i_r}\colon U_{i_1}\times_X\dots\times_X U_{i_r}\to X$. Définissons le préfaisceau \[ \ZZ_\mathcal{U}\coloneqq\coker\left[\bigoplus_{i,j\in I}(f_{ij})_!\ZZ\to \bigoplus_{i\in I}(f_i)_!\ZZ\right]\]où les extensions par zéro sont à comprendre au sens des préfaisceaux abéliens. Une résolution projective de $\ZZ_\mathcal{U}$ dans la catégorie $\PAb(X)$ des préfaisceaux de groupes abéliens sur $X$ est donnée par 
\[ K_\mathcal{U}\coloneqq \cdots\to\bigoplus_{i_1,\dots,i_r\in I}(f_{i_1,\dots,i_r})_!\ZZ\to\cdots\to\bigoplus_{i,j\in I}(f_{ij})_!\ZZ\to \bigoplus_{i\in I}(f_i)_!\ZZ.\]
Pour tout faisceau $\F\in\Ab(X)$, le complexe $\Hom(K_\mathcal{U},\F)$, où les termes non nuls de $K_\mathcal{U}$ sont placés en degrés $]-\infty,0]$, est le complexe de \v{C}ech usuel. 
Notons \[\check{\Gamma}_P(\mathcal{U},-)\coloneqq \Hom_{\PAb(X)}(\ZZ_\mathcal{U},-)\colon\PAb(X)\to \Ab.\] Soit $O\colon \Ab(X)\to\PAb(X)$ le foncteur d'oubli. Alors \[\Gamma(X,-)=\check{\Gamma}_P(\mathcal{U},-)\circ O\colon \Ab(X)\to \Ab\]
et comme tout faisceau injectif est encore un préfaisceau injectif, 
\[ \RG(X,-)=\R\check{\Gamma}_P(\mathcal{U},-)\circ \R O.\]
\'{E}tant donné un faisceau $\F\in\Ab(X)$, $\R^iO(\F)$ associe à tout $X$-schéma étale $Y$ le groupe $\HH^i(Y,\F)$. Dans le cas où $\mathcal{U}=\{ U,V\}$ est un recouvrement de $X$ par deux ouverts, la suite spectrale \[ E_2^{pq}=\check{\HH}^p(\mathcal{U},R^qO(\F))\Rightarrow \HH^{p+q}(X,\F)\]
dégénère à la seconde page, et donne une suite exacte
\[ 0\to \HH^0(X,\F)\to \HH^0(U,\F)\oplus \HH^0(V,\F)\to \HH^0(U\cap V,\F)\to \HH^1(X,\F)\to \HH^1(U,\F)\oplus \HH^1(V,\F)\to \dots.\]
Cette suite exacte, dite de Mayer-Vietoris, découle également de la proposition suivante.

\begin{prop}[Triangle de Mayer-Vietoris]\cite[0CRS]{stacks} Soit $X$ un schéma recouvert par deux ouverts $U,V$. Soit $\F\in \Ab(X)$. Alors il y a un morphisme $\RG(U\cap V,\F)\to \RG(X,\F)[1]$ qui fait de
\[ \RG(X,\F)\to\RG(U,\F)\oplus\RG(V,\F)\to \RG(U\cap V,\F)\to \RG(X,\F)[1]\]
un triangle distingué.
\end{prop}

\subsection{Théorèmes de changement de base}

Les théorèmes de changement de base permettent d'exprimer autrement des images directes dérivées ; par exemple, le théorème de changement de base propre permet d'exprimer les fibres de l'image directe dérivée d'un faisceau comme la cohomologie de ce faisceau sur les fibres du morphisme propre en question. Ils servent dans de nombreuses démonstrations, notamment celles des théorèmes de finitude présentés dans la section suivante.

Considérons le diagramme cartésien de schémas suivant :
\[\begin{tikzcd}
X' \arrow[r,"v",swap]\arrow[d,"u"] & X \arrow[d,"f"] \\
S'\arrow[r,"g"] & S
\end{tikzcd}\]
Soit $\F$ un faisceau de groupes abéliens sur $X$. Le morphisme 
\[ \R f_\star\F\to \R f_\star \R v_\star v^\star \F \]
obtenu par l'adjonction $v^\star \dashv \R v_\star$ donne via les identifications canoniques \[ \R f_\star \R v_\star =\R(f\circ v)_\star =\R(g\circ u)_\star=\R g_\star \R u_\star \]
un morphisme \[ \R f_\star \F \to \R g_\star \R u_\star v^\star \F \]
qui fournit lui-même, par l'adjonction $g^\star \dashv \R g_\star$, un morphisme dit de changement de base
\[ g^\star \R f_\star\F \to \R u_\star v^\star \F.\]

\begin{theorem}[Changement de base propre]\cite[XII, Th. 5.1]{sga43}\label{th:chbp} Si $f$ est propre et $\F$ est un faisceau de torsion alors le morphisme de changement de base est un isomorphisme.
\end{theorem}

\begin{cor} Soit $\bar{s}$ un point géométrique de $S$. Si $f\colon X\to S$ est propre et $\F$ est un faisceau abélien de torsion sur $X$ alors la fibre 
$(\R f_\star\F)_{\bar{s}}$ est canoniquement isomorphe à $\RG(X_{\bar{s}},\F)$.
\end{cor}

\begin{theorem}[Changement de base lisse]\cite[XVI, Th. 1.1]{sga43} Si $g$ est lisse, $f$ est quasi-compact et quasi-séparé et $\F$ est un faisceau dont les fibres géométriques sont de torsion d'ordre inversible sur $X$ alors le morphisme de changement de base est un isomorphisme.
\end{theorem}

\subsection{Images directes dérivées et finitude de la cohomologie}
Soit $f\colon Y\to X$ un morphisme de schémas.
Le foncteur image directe $f_\star\colon \Ab(Y)\to \Ab(X)$ est exact à gauche, et son foncteur dérivé est noté $\R f_\star$. Cette section résume comment calculer les fibres de l'image directe dérivée d'un faisceau, et sous quelles hypothèses les propriétés de constructibilité et de lissité d'un faisceau sont préservées par l'image directe dérivée.

\begin{prop}\cite[VIII, Th. 5.2]{sga42}\label{prop:fibforward} Soit $f\colon Y\to X$ un morphisme quasi-compact et quasi-séparé de schémas. Soit $K\in D^b(Y)$ un complexe de faisceaux abéliens sur $Y$. Soit $\bar x$ un point géométrique de $X$. Notons $X_{\bar x}$ le spectre de l'anneau local strictement hensélien de $X$ en $\bar x$. Il y a un isomorphisme canonique dans $D^+(\Ab)$ :
\[ (Rf_\star K)_{\bar x}\xrightarrow{\sim}\RG(Y\times_X X_{\bar x},K).\]
\end{prop}

\begin{theorem}\cite[XIII, Th. 1.1.1]{travaux_gabber} Soient $X$ un schéma quasi-excellent, $f\colon Y\to X$ un morphisme de type fini, $n\geqslant 1$ un entier inversible sur $X$ et $\F$ un faisceau constructible de $\ZZ/n\ZZ$-modules sur $Y$. Alors : \begin{enumerate}
\item Pour tout entier $q\geqslant 0$, le faisceau $\R^q f_\star\F$ est constructible.
\item Il existe un entier $N$ tel que $\R^q f_\star \F=0$ pour tout $q\geqslant N$.
\end{enumerate}
\end{theorem}

Un cas particulier de ce théorème, lorsque la cible du morphisme est un corps séparablement clos, est le suivant.

\begin{cor} Soit $X$ un schéma de type fini sur un corps séparablement clos $k$. Pour tout faisceau constructible $\F$ de groupes abéliens sur $X$ et tout entier naturel $i$, le groupe $\HH^i(X,\F)$ est fini.
\end{cor}

\begin{theorem}\cite[I, Th. 8.9]{freitag_kiehl}\label{th:ehresmann} Soit $f\colon Y\to X$ un morphisme propre et lisse de schémas. Soit $n$ un entier inversible sur $X$. Pour tout faisceau lisse de $\ZZ/n\ZZ$-modules sur $Y$ et tout entier naturel $i$, le faisceau $\R^if_\star\F$ est lisse sur $X$.
\end{theorem}

\section{Cohomologie à support dans un fermé}

\subsection{Généralités}\label{subsec:cohsuppferm}

Soit $X$ un schéma. Soient $i\colon Z\to X$ une immersion fermée, et $j\colon U\to X$ l'inclusion de l'ouvert complémentaire.

\begin{df} Le foncteur $\Gamma_Z(X,-)\colon\Ab(X)\to \Ab$ des sections à support dans $Z$ est défini par \[ \Gamma_Z(X,-)=\Gamma(Z,i^!-)\colon\F\mapsto \ker(\F(X)\to \F(U)).\]
Ce foncteur est exact à gauche, et ses foncteurs dérivés à droite sont notés $\HH^j_Z(X,-)$.
\end{df}

\begin{prop}\cite[III, Prop. 1.25]{milneEC}\label{prop:of} Soit $K\in \DD^b(X)$. La suite exacte
\[ 0\to j_!\ZZ\to \ZZ\to i_\star \ZZ\to 0\]
produit par $\RHom(-,\F)$ un triangle distingué "ouvert-fermé"
\[\RG_Z(X,\F)\to\RG(X,\F)\to\RG(U,\F)\xrightarrow{+1}.\]
\end{prop}

\subsection{Théorème de pureté et suite de Gysin}\label{subsec:gysin}

Le théorème suivant permettra de calculer la cohomologie à support dans un fermé des faisceaux lisses.

\begin{theorem}[Pureté]\cite[XVI, Th. 3.1.1]{travaux_gabber} Soit $X$ un schéma régulier. Soit $i\colon Z\to X$ une immersion fermée d'un sous-schéma régulier, de codimension $c$. Soit $\F$ un faisceau localement constant sur $X$ de torsion inversible sur $X$. Alors il y a un isomorphisme $\Lambda \to \R i^!\Lambda(c)[2c]$ dans $D^+(Z,\Lambda)$.
\end{theorem}

Nous serons intéressés par le cas où $Z$ et $X$ sont des variétés sur un corps séparablement clos $k$, et $\F$ est un faisceau lisse de $\ZZ/n\ZZ$-modules sur $X$. Le théorème fournit alors un isomorphisme 
\[ \HH^{r-2c}(Z,\F(-c))\xrightarrow{\sim}\HH^r_Z(X,\F). \]
En notant $U$ l'ouvert complémentaire de $Z$ dans $X$, la suite ouvert-fermé devient alors la \textit{suite de Gysin} :
\[ 0 \to \HH^{2c-1}(X,\F) \to \cdots \to \HH^{r-2c}(Z,\F(-c))\to \HH^r(X,\F)\to \HH^r(U,\F)\to \cdots \]

\section{Cohomologie à support compact}

\subsection{Généralités}

Soit $k$ un corps. Soit $X$ un schéma séparé de type fini sur $k$. D'après un théorème de Nagata \cite[Th. 4.3]{nagata}, il existe une immersion ouverte $j\colon X\to \bar X$, où $\bar X$ est propre sur $k$. \'{E}tant donné un complexe $K\in \DD^b(X)$, l'élément $\RG(\bar X,j_!K)\in \DD^b(\Ab)$ est indépendant du choix de $\bar X$ \cite[Prop. 18.2]{milneLEC}.

\begin{df}
Pour tout $K\in \DD^b(X,\Lambda)$, on définit 
\[ \RG_c(X,K)=\RG(\bar X,j_!K). \]
Les groupes de cohomologie à support compact de $K$, notés $\HH^i_c(X,K)$, sont les groupes de cohomologie de $\RG_c(X,K)$.
\end{df}

\begin{rk} Le foncteur $\RG_c(X,-)$ n'est pas le foncteur dérivé de $\HH^0_c(X,-)$. Prenons $X=\A^1$. Nous verrons dans la section \ref{subsubsec:courbcohsupp} que $\HH^1(\PP^1,j_!\Lambda)=\Lambda$. Par contre, \[ \HH^0(\PP^1,j_!-)=\bigoplus_{x\in |\A^1|}\HH^0_x(\A^1,-).\] Comme $\HH^1_x(\A^1,\Lambda)=\Lambda$ (voir section \ref{sec:cohfermcourb}), \[\R^1\HH^0_c(\A^1,\Lambda)\simeq \bigoplus_{x\in |\A^1|}\Lambda \neq \HH^1_c(\A^1,\Lambda).\]
\end{rk}

\subsection{Dualité de Poincaré}

Soient $k$ un corps algébriquement clos et $n$ un entier inversible dans $k$. Notons $\Lambda$ l'anneau $\ZZ/n\ZZ$. Soit $X$ un schéma connexe séparé de type fini sur $k$, lisse de dimension $d$.  Soit $\F$ un faisceau constructible de $\Lambda$-modules sur $X$. Rappelons que pour tout entier naturel $i$, $\Ext^i(\F,\Lambda(d))=\Hom_{\DD(X,\Lambda)}(\F,\Lambda(d)[i])$. Un morphisme $\F\to \Lambda(d)[i]$ définit pour tous entiers naturels $i,j$ un morphisme \[ \HH^j_c(X,\F)\to \HH^{i+j}_c(X,\Lambda(d)).\] Ceci permet de définir un accouplement \[ \HH^j_c(X,\F)\times \Ext^{i}(\F,\Lambda(d))\to \HH^{i+j}_c(X,\Lambda(d)).\]

L'énoncé général de la dualité de Poincaré se trouve dans \cite[XVIII, Th. 3.2.5]{sga43} ; le cas particulier qui nous concerne est le suivant.
\begin{theorem}\cite[VI, Th. 11.1]{milneEC}\label{th:dualpoinc} Il y a un isomorphisme canonique $\HH^{2d}_c(X,\Lambda(d))\xrightarrow{\sim} \Lambda$. Pour tout entier $j\in \{0\dots 2d\}$, l'accouplement \[ \HH^j_c(X,\F)\times \Ext^{2d-j}(\F,\Lambda(d))\to \HH^{2d}_c(X,\Lambda(d))\xrightarrow{\sim}\Lambda\]
est non dégénéré.
\end{theorem} 
Nous nous intéresserons particulièrement au cas où $X$ est une courbe et $\F$ est lisse. Le groupe $\Ext^{2-j}(\F,\Lambda(d))$ est alors canoniquement isomorphe à $\HH^{2-j}(X,\F^\vee(1))$. Il y a donc en particulier un accouplement non dégénéré
\[ \HH^1_c(X,\F)\times \HH^1(X,\F^\vee(1))\to \Lambda.\]

\section{Formule des traces et comptage de points}

Soient $k_0$ un corps fini de cardinal $q$, et $k$ une clôture algébrique de $k_0$. Pour tout entier $m\geqslant 1$, notons $k_m$ l'extension de degré $m$ de $k_0$ dans $k$. Soit $\ell$ un nombre premier inversible dans $k$.

\begin{df} Soit $X$ un schéma sur $k_0$. L'endomorphisme de Frobenius géométrique de $X$, noté $\Frob_X$, est l'endomorphisme défini par l'identité sur l'espace topologique sous-jacent, et par la mise à la puissance $q$ sur le faisceau structural $\OO_X$.
\end{df}

En particulier, si $X_0$ est une variété sur $k_0$, le morphisme de Frobenius géométrique sur $X(k)$ correspond à la mise à la puissance $q$ des coordonnées des points. Notons $X=X_0\times_{k_0}k$. Les $k_m$-points de $X$ sont alors les points fixes de $\Frob_X^m$ sur $X_0(k)$. La formule des traces permet de calculer le nombre d'intersection dans $X\times X$ du graphe de $\Frob_X$ avec la diagonale ; cette intersection étant transverse, ce nombre d'intersection est exactement le nombre de $k_0$-points de $X_0$. Nous n'aurons pas besoin de l'énoncé général portant sur la cohomologie $\ell$-adique, mais simplement de la formulation suivante.

\begin{theorem}\cite[Rapport, Th. 3.2]{sga412}\label{th:trace}
Soit $X$ un schéma de type fini sur $k_0$ de dimension $d$. Notons $X=X_0\times_{k_0}k$. Soit $\ell$ un nombre premier inversible dans $k_0$. Alors pour tous entiers $m,n\geqslant 1$,
\begin{eqnarray*} \# X(k_m)&\equiv &\sum_{i=1}^{2d}\tr((\Frob_X^m)^\star\mid \HH^i_c(X,\ZZ/\ell^n\ZZ)) \mod \ell^n.
\end{eqnarray*}
\end{theorem}

\cleartooddpage

\chapter{Revêtements et cohomologie des courbes}\label{chap:2}

Dans tout ce chapitre, le mot \textit{courbe} désigne un schéma équidimensionnel de dimension 1 sur un corps. Nous préciserons à chaque fois que cela sera nécessaire s'il s'agit d'une courbe connexe, intègre, lisse, affine, projective...
Ce chapitre décrit explicitement les groupes de cohomologie des faisceaux constants sur les courbes intègres lisses ou nodales, ainsi que leurs revêtements galoisiens. En particulier, il contient en \ref{sec:revcar} la construction d'un revêtement caractéristique de ces courbes qui servira par la suite. \\

\paragraph{Notations} Rappelons que la cohomologie d'un schéma sur un corps non parfait est canoniquement isomorphe à celle du changement de base de ce schéma à la clôture parfaite du corps (voir section \ref{subsec:invtop}). Dans tout ce chapitre, $k_0$ désigne un corps parfait de caractéristique $p\geqslant 0$, et $k$ une clôture algébrique de $k_0$. Le groupe $\Gal(k|k_0)$ est noté $\mathfrak{G}_0$. Nous fixons un nombre premier $\ell$ distinct de $p$ ainsi qu'un entier naturel $n$ premier à $p$ et notons $\Lambda$ l'anneau $\ZZ/n\ZZ$. 

\section{Groupe fondamental des courbes}

\subsection{Revêtements de courbes}

Nous serons principalement intéressés par les morphismes entre courbes normales (i.e. régulières par \cite[IV, Th. 11]{serre}), qui sont déterminés par des extensions de corps.

\begin{theorem}\cite[Th. 4.3.10, Prop. 4.4.5]{szamuely}\label{th:revcourb} Soit $X$ une courbe intègre régulière sur un corps. Soient $\mathcal{C}$ la catégorie des courbes intègres régulières $Y$ munies d'un morphisme fini $f\colon Y\to X$, et $\mathcal{D}$ la catégorie des extensions finies du corps des fonctions $k(X)$ de $X$. Le foncteur $\mathcal{C}^{\op}\to \mathcal{D}$ qui à $(Y,f)$ associe l'extension $k(X)\xrightarrow{f^\star}k(Y)$ est une équivalence de catégories.
\end{theorem}
Ceci n'est \textit{jamais} le cas pour les courbes singulières, puisque le morphisme de normalisation induit un isomorphisme entre les corps de fonctions. Pour les morphismes vers une courbe régulière, la platitude est garantie par le résultat ci-après.
\begin{lem}\cite[0CCK]{stacks} Soit $f\colon Y\to X$ un morphisme non constant de courbes intègres sur un corps. Si $X$ est normale alors $f$ est plat.
\end{lem}
Par conséquent, il est aisé de vérifier si un morphisme $f\colon Y\to X$ entre courbes régulières sur un corps est étale en un point $y$ de $Y$. L'anneau local $\OO_{X,f(y)}$ est un anneau de valuation discrète ; choisissons une uniformisante $\pi$ de cet anneau. Alors $f$ est étale si et seulement si $f^\star\pi$ est une uniformisante de $\OO_{Y,y}$. La proposition suivante affirme qu'une extension séparable de corps de fonctions correspond à un morphisme \textit{génériquement étale}. 
\begin{prop}\cite[Prop. 4.5.9]{szamuely} Soit $f\colon Y\to X$ un morphisme de courbes intègres sur un corps. Si l'extension de corps de fonctions correspondante est séparable alors il existe un ouvert $U$ de $Y$ tel que $f|_U$ soit étale.
\end{prop}

\begin{ex}\begin{enumerate}
\item Soit $p$ un nombre premier. Soit $X=\PP^1_{\FF_p}$. Le morphisme $Y\to X$ défini par l'extension purement inséparable $\FF_p(t)\to \FF_p(\sqrt[p]{t})$ n'est étale en aucun point de $Y$ \cite[0CCY]{stacks}.
\item Considérons le morphisme de courbes lisses \[ f\colon Y=\Spec \CC[x,y]/(y^3-y+x)\to \A^1_{\CC}=\Spec \CC[x] \] défini par $(x,y)\mapsto x$. L'extension de corps de fonctions $\CC(x)\to \CC(x)[y]/(y^3-y+x)$ 
est séparable, car son discriminant $-4+27x^2$ est non nul. Par conséquent, le morphisme $f$ est génériquement étale. Elle n'est pas normale, car comme $-4+27x^2$ n'est pas un carré dans $\CC(x)$, le polynôme minimal $T^3-T+x$ de $y$ n'a pas toutes ses racines dans $\CC(Y)$. 
Le morphisme est ramifié au-dessus des racines de $-4+27x^2$. Considérons l'ouvert
\[ U=\Spec k[x,(27x^2-4)^{-1}] \]
de $\A^1$, et sa préimage \[ V=\Spec k[x,(27x^2-4)^{-1},y]/(y^3-y+x)\] dans $X$. Le morphisme $V\to U$ est alors fini étale, mais pas galoisien.
\end{enumerate}
\end{ex}

\subsection{Groupe fondamental des courbes}

Il est généralement très difficile de calculer le groupe fondamental d'un schéma ; toutefois, le cas des courbes sur un corps algébriquement clos est bien étudié. Soit $X$ une courbe intègre lisse sur un corps algébriquement clos $k$ de caractéristique $p$. Lorsque $k=\CC$, la théorie des extensions de corps de $k(X)$ revient à l'étude des extensions du corps des fonctions méromorphes sur la surface de Riemann associée à $X$, et le groupe fondamental de $X$ est le complété profini du groupe fondamental topologique de cette surface de Riemann. Ce résultat se transpose ensuite à tout corps de caractéristique nulle. La preuve du résultat correspondant en caractéristique positive, dû à Grothendieck, utilise des techniques bien plus profondes. Elle consiste à considérer $X$ comme la fibre spéciale d'un schéma sur un anneau de valuation discrète de corps des fractions de caractéristique nulle et de corps résiduel $k$, puis à conclure par un théorème de spécialisation du groupe fondamental \cite[X, Th. 3.8]{sga1}. 
\begin{theorem}\cite[Th. 5.7.13]{szamuely} Soit $r\in\NN$. Soit $X$ une courbe intègre propre lisse sur $k$ de genre $g$. Soient $x_1,\dots,x_r$ des points fermés de $X$. Le plus grand quotient pro-$p'$ (c'est-à-dire limite de groupes finis d'ordre premier à $p$) du groupe fondamental de $U\coloneqq X-\{x_1,\dots,x_r\}$ est le groupe pro-$p'$ donné par la présentation suivante. \[ \pi_1^{(p')}(U)=\langle a_1,b_1,\dots,a_{g},b_{g},c_1,\dots,c_r\mid [a_1,b_1]\cdots [a_g,b_g ]c_1\cdots c_r=1\rangle_{p'}\]
\end{theorem}
Ce théorème donne en particulier une condition nécessaire pour qu'un groupe fini $G$ apparaisse comme le groupe d'automorphismes d'un revêtement galoisien d'une courbe de genre $g$ : il faut pour cela que le plus grand quotient de $G$ d'ordre premier à $p$ puisse être engendré par une famille d'au plus $2g+r$ éléments. La conjecture d'Abhyankar, démontrée par Harbater en 1994 \cite[Th. 6.2]{harbater}, affirme la réciproque.

\begin{rk}
En caractéristique positive, le groupe fondamental lui-même devient gigantesque. Par exemple, celui de $\A^1_{\overline{\FF_p}}$ n'est pas topologiquement de type fini. En effet, il y a pour toute puissance $p^j$ de $p$ un revêtement de groupe $(\ZZ/p\ZZ)^j$ : le revêtement d'Artin-Schreier $x\mapsto x^{p^j}-x$. Ceci permet de construire, pour chaque entier $j$ positif, $p^j$ morphismes continus $\pi_1(\A^1_{\overline{\FF_p}})\to \ZZ/p\ZZ$. Il y a donc une infinité de tels morphismes continus, ce qui serait impossible si le groupe fondamental était topologiquement de type fini.
\end{rk}

Le cas des courbes singulières a été étudié récemment ; par normalisation, le groupe fondamental d'une courbe singulière s'exprime comme produit libre du groupe fondamental d'une courbe lisse par un groupe profini libre.

\begin{theorem}\cite[Th. 1.1]{das_galois}
Soit $X$ une courbe projective connexe sur $k$ à $s$ composantes irréductibles. Soit $\nu\colon\tX=C_1\sqcup\dots\sqcup C_s\to X$ sa normalisation. Notons \[ \delta\coloneqq 1-s+\sum_{x\in X}\left(|\nu^{-1}(x)|-1\right).\]
Notons $F_\delta$ le groupe libre à $\delta$ générateurs, et $\widehat{F_\delta}$ son complété profini. Il y a un isomorphisme de groupes profinis
\[ \pi_1(X)\xrightarrow{\sim} \pi_1(C_1)\star\dots \star\pi_1(C_s)\star \widehat {F_\delta}.\]
\end{theorem}

\subsection{Revêtements génériquement étales}\label{subsec:inert}

\begin{df} Soit $f\colon Y\to X$ un morphisme fini de courbes intègres sur $k_0$. Nous dirons que $f$ est génériquement galoisien si l'extension de corps de fonctions correspondante est galoisienne.
\end{df}
\begin{prop}\label{prop:galcorps} Soit $f\colon Y\to X$ un morphisme fini étale de courbes intègres normales sur un corps. Notons $L/K$ l'extension de corps de fonctions correspondante. Alors $\Aut(Y|X)^{\op}=\Aut(L|K)$, et le morphisme $f$ est un revêtement galoisien si et seulement s'il est génériquement galoisien.
\begin{proof} Le théorème \ref{th:revcourb} assure que $\Aut(Y|X)^{\op}=\Aut(L|K)$. Si $f$ est étale alors l'extension $L/K$ est séparable. Par conséquent, le groupe $\Aut(L|K)$ est d'ordre $\deg(f)=[L:K]$ si et seulement si l'extension $L/K$ est galoisienne. 
\end{proof}
\end{prop}

Un morphisme génériquement galoisien de courbes est étale, et donc galoisien, sur un ouvert. En particulier, si $f\colon Y\to X$ est un revêtement galoisien de courbes affines connexes, le morphisme $\bar f\colon \bar Y\to\bar X$ entre les compactifications lisses de ces courbes est un revêtement génériquement galoisien. Lorsque les courbes en question sont régulières, l'étude locale du morphisme $Y\to X$ au-dessus d'un point $x$ de $X$ revient à l'étude d'une extension de l'anneau de valuation discrète $\OO_{X,x}$. Nous nous placerons donc pour le reste de cette section dans la situation 
\[
\begin{tikzcd}
 B \arrow[r,hook] &L \\
 A \arrow[u,hook]\arrow[r,hook] &K\arrow[u,hook]
\end{tikzcd}
\]
où $A$ est un anneau de valuation discrète, $L$ est une extension galoisienne de son corps des fractions $K$, et $B$ est la normalisation de $A$ dans $L$. Le groupe $G\coloneqq \Aut(L|K)$ agit transitivement sur les idéaux maximaux de $B$ au-dessus de $\m_A$. Soit $\m_B$ un idéal maximal de $B$ au-dessus de $\m_A$. Notons $k_A\coloneqq A/\m_A$.
\begin{df} Le groupe de décomposition de $\m_B$ est le groupe \[ D_{\m_B}=\{ \sigma\in G\mid \sigma(\m_B)=\m_B\}.\]
Le groupe d'inertie associé est le groupe \[ I_{\m_B}=\ker(D_{\m_B}\to \Aut((B/\m_B)|k_A)).\]
\end{df}

Supposons désormais $L/K$ finie. Notons $\m_1,\dots,\m_d$ les idéaux maximaux de $B$ au-dessus de $\m_A$, et $k_1,\dots,k_d$ leurs corps résiduels. Il existe \cite[09EB]{stacks} des entiers $e$ et $f$, appelés respectivement indice de ramification et degré résiduel de $B/A$, tels que pour tout $i\in \{1\dots d\}$, $\m_AB_{\m_i}=\m_i^e B_{\m_i}$ et $[k_i:k_A]=f$. L'extension est dite sauvagement ramifiée au-dessus de $\m_A$ si la caractéristique de $k_A$ divise $e$, et modérément ramifiée sinon. De plus, $[L:K]=def$. 

Fixons un indice $i\in \{1\dots d\}$, posons $\m=\m_i$ et $k_\m=k_i$.
Notons $\pi_A$ une uniformisante de $A$, et $\pi_B$ une uniformisante du localisé $B_{\m}$. Alors pour tout $\sigma\in I_\m$, il existe un unique élément $u_\sigma\in B_{\m}^\times$ tel que $\pi_B=\sigma(\pi_B)u_\sigma$. 
\begin{prop}\cite[09EE]{stacks}\label{prop:gpinert} La composée
\[ I_{\m}\xrightarrow{\sigma\mapsto u_\sigma}B_{\m}^\times\to k_{\m}^\times\]
est à image dans le groupe $\mu_e(k_\m)$ des racines $e$-ièmes de l'unité dans $k_\m$ et définit un morphisme de groupes surjectif \[I_\m\to \mu_e(k_\m)\]
dont le noyau $P_{\m}$ est nul si $k_A$ est de caractéristique nulle, et un $p$-groupe si la caractéristique de $k_A$ est $p>0$. Le groupe quotient $I_{\m}/P_{\m}$ est cyclique d'ordre le plus grand diviseur de $e$ premier à $p$.
\begin{proof}
Vérifions simplement la première assertion. Il existe un unique $u\in A^\times$ tel que $\pi_B^e=u\pi_A$. Comme $u_\sigma=\sigma(\pi_B)\pi_B^{-1}$, $u_\sigma^e=\sigma(u)u^{-1}$. Comme $\sigma\in I_{\m}$, cet élément est congru à 1 modulo $\m_B$.
\end{proof} 
\end{prop}

Nous serons amenés, étant donné un tel revêtement, à calculer la cohomologie de son groupe d'inertie. 

\begin{lem} Soit $m$ un entier. Soit $C$ un groupe cyclique d'ordre divisible par $m$. Soit $M$ un $(\ZZ/m\ZZ)[C]$-module sur lequel le sous-groupe de $C$ d'ordre $m$ agit trivialement. Alors le choix d'un générateur de $C$ détermine un isomorphisme $\HH^1(C,M)\to M_C$. La notation $M_C$ désigne le module des coinvariants, c'est-à-dire le quotient de $M$ par le sous-module engendré par les éléments de la forme $\sigma\cdot m-m$ où $\sigma\in C$ et $m\in M$.
\begin{proof} Soit $t$ un générateur de $C$. Notons $c$ l'ordre de $C$, et $N=\sum_{j=0}^{c-1}t^j\in \End(M)$. 
Les résultats classiques sur la cohomologie des groupes cycliques \cite[III.1, Ex. 2]{brown} montrent alors qu'il y a un isomorphisme $\HH^1(C,M)\to \ker(N)/(t-1)M$, qui associe à un morphisme croisé $f\colon C\to M$ la classe de l'élément $f(t)\in\ker(N)$. 
Comme $\langle t^{c/m}\rangle$ agit trivialement sur $A$, \[ N=\sum_{i=0}^{m-1}\sum_{j=0}^{c/m-1}t^{ie/m+j}=m\sum_{j=0}^{c/m-1}t^j\]
est l'endomorphisme nul puisque $M$ est de $m$-torsion. Par conséquent, $\HH^1(C,M)$ est isomorphe à\[ M/(t-1)M=M_C.\]
\end{proof}
\end{lem}

La preuve du résultat suivant est une adaptation de celle de \cite[Prop. 8.1.4]{fulei}. Rappelons que $\Lambda=\ZZ/n\ZZ$, où $n$ est un entier inversible dans le corps algébriquement clos $k$.
\begin{cor}\label{cor:cohinert} Soit $f\colon Y\to X$ un revêtement galoisien de courbes lisses sur $k$ de groupe $G$. Soient $y$ un point fermé de $Y$ et $x=f(y)$. Supposons que l'indice de ramification de $f$ en $y$ soit divisible par $n$. Il y a pour tout $\Lambda[G]$-module $M$ un isomorphisme canonique \[ \HH^1(I_y,M)\xrightarrow{\sim} M_I(-1)\]
où $M_I(-1)$ désigne $M_I\otimes_\Lambda\mu_n(k)^\vee=M_I\otimes_\Lambda\Hom(\mu_n,\Lambda)$.
\begin{proof} D'après la proposition \ref{prop:gpinert}, le groupe $I_y/P_y$ est canoniquement isomorphe à $\mu_e(k)$, où $e$ est l'indice de ramification en $y$. Fixons un générateur $\zeta$ de $\mu_e(k)$. Comme $P_y$ est un $p$-groupe et $M$ est de $n$-torsion, le morphisme $M^{P_y}\to M_{P_y}$ est un isomorphisme et les groupes de cohomologie $\HH^i(P_y,M)$ sont nuls dès que $i\geqslant 1$ \cite[Prop. 6.1.10]{weibel}.
Le lemme précédent appliqué au $I_y/P_y$-module $M^{P_y}$ assure qu'il y a un isomorphisme \[ \HH^1(I_y/P_y,M^{P_y})\xrightarrow{\sim} (M^{P_y})_{I_y/P_y}\simeq (M_{P_y})_{I_y/P_y}\simeq M_{I_y} \] qui à un morphisme croisé $f\colon I_y/P_y\simeq \mu_e(k)\to M$ associe l'élément $f(\zeta)$. Le morphisme \[ \HH^1(I_y/P_y,M)\otimes \mu_e(k)\to M_{I_y} \] qui à un morphisme croisé $f\colon I_y/P_y\simeq \mu_e(k)\to M_{I_y}$ et un élément $t\in \mu_n(k)$ associe $f(t)$ est encore un isomorphisme, canonique cette fois-ci. Par conséquent, il y a des isomorphismes canoniques \[\HH^1(I_y/P_y,M)=M_{I_y}\otimes\mu_e(k)^\vee=M_{I_y}\otimes\mu_n(k)^\vee.\]
 La suite spectrale de Hochschild-Serre fournit alors la suite exacte
\[ 0\to \HH^1(I_y/P_y,M_{P_y})\to \HH^1(I_y,M)\to \HH^0(I_y/P_y,\HH^1(P_y,M))=0\]
qui permet de conclure.
\end{proof}
\end{cor}

\begin{df}[Troncature] Soit \[ K= \cdots \xrightarrow{d^{i-1}} K^i \xrightarrow{d^i} K^{i+1}\xrightarrow{d^{i+1}}\cdots \]
 un complexe dans une catégorie abélienne $\mathcal{A}$. Soit $r$ un entier. Nous noterons $\tau_{\leqslant r}K$ et appellerons tronqué de $K$ en degré $\leqslant r$ le complexe :
\[ \cdots \xrightarrow{d^{r-2}}K^{r-1}\xrightarrow{d^{r-1}} \ker(d^r)\to 0\to 0\to\cdots \]
Le morphisme évident $\tau_{\leqslant r}K\to K$
induit un isomorphisme $\HH^i(\tau_{\leqslant r} K)\to \HH^iK$ pour tout entier $i\leqslant r$. Ce foncteur de troncature définit un endofoncteur de la catégorie dérivée $D(\mathcal{A})$, que nous noterons encore $\tau_{\leqslant r}$.
\end{df}

Pour tout groupe $H$ agissant sur $M$, notons $C^{12}(H,M)$ le groupe des morphismes croisés $H\to M$. Le complexe $M\xrightarrow{\partial} C^{12}(H,M)$ représente alors $\tau_{\leqslant 1}\RG(H,M)$.

\begin{lem}\label{lem:secIP} Reprenons les notations et hypothèses du corollaire précédent. Le morphisme canonique $C^{12}(I_y/P_y,M^{P_y})\to C^{12}(I_y,M)$ admet une section. 
\begin{proof} Soit $u\colon I_y\to M$ un morphisme croisé.
Considérons le diagramme commutatif
\[
\begin{tikzcd}
I_y \arrow[d]\arrow[r,"u"] & M \arrow[r,"q"] & M_{P_y} \\
I_y/P_y \arrow[r,dashed] & M^{P_y}\arrow[ur,"\alpha"]\arrow[u] & 
\end{tikzcd}
\]
et notons $f=q\circ u$. Pour tout $x\in P_y$ et tout $g\in I_y$, la définition de $M_{P_y}$ assure que $f(x g)=f(x)+q(x\cdot u(g))=f(x)+f(g)$. Par conséquent, pour tout $x\in P_y$, $f(x^{|P_y|})=|P_y|f(x)$ doit être nul ; comme la multiplication par $|P_y|$ est un automorphisme de $M$, cela signifie que $f(x)=0$. Par conséquent, $f$ passe au quotient en $\bar f\colon I_y/P_y\to M_{P_y}$. Notons $\bar u=\alpha^{-1}\circ \bar{f} \colon I_y/P_y\to M^{P_y}$. L'application $u\mapsto\bar u$ est évidemment linéaire. De plus, $\bar u$ est encore un morphisme croisé car $\bar{u}(\bar g\bar h)=\alpha^{-1}f(gh)=\alpha^{-1}f(g)+q(g\cdot \alpha^{-1}u(h))$. La $I_y$-linéarité de $\alpha^{-1}q$ conclut.
\end{proof}
\end{lem}

\begin{rk}
Ce lemme assure que le quasi-isomorphisme des complexes de cochaînes usuels \[ \tau_{\leqslant 1}\RG(I_y/P_y,M^{P_y})\to \tau_{\leqslant 1} \RG(I_y,M) \]
admet un inverse dans $\DD^b_c(\Lambda)$ qui est donné par un vrai morphisme de complexes.
\end{rk}

Les revêtements étales d'une courbe affine modérément ramifiés à l'infini sont eux aussi classifiés par un groupe profini : le groupe fondamental modéré.
\begin{df} Soient $X$ une courbe projective intègre lisse sur un corps séparablement clos $k$, et $K$ son corps des fonctions. Soit $U$ un ouvert de $X$. Soient $K^\sep$ une clôture séparable de $K$, et $\bareta$ le point générique géométrique correspondant de $U$. Soit $K^\ttt$ la composée dans $K^\sep$ des extensions finies $L/K$ telles que la normalisation de $X$ dans $L$ soit étale sur $U$ et modérément ramifiée au-dessus de $X-U$. Le groupe fondamental modéré de $U$, noté $\pi_1^\ttt(U,\bareta)$ ou simplement $\pi_1^\ttt(U)$, est le groupe $\Gal(K^\ttt|K)$.
\end{df}
De la même façon que $\pi_1^{(p')}(U)$, le groupe $\pi_1^\ttt(U)$ est topologiquement de type fini \cite[XIII, Cor. 2.12]{sga1}. Un résultat plus fort a été récemment démontré par Esnault, Shusterman et Srinivas.
\begin{theorem}\cite[Th. 1.2]{esnault} Soit $X$ une courbe intègre lisse sur un corps algébriquement clos de caractéristique positive. Le groupe $\pi_1^\ttt(X)$ est un groupe profini de présentation finie. Si $X$ est affine alors $\pi_1^\ttt(X)$ est projectif, c'est-à-dire isomorphe à un sous-groupe d'un groupe profini libre.
\end{theorem}

\subsection{Ramification : théorie locale}

Soit $X$ une courbe intègre lisse sur $k$. Notons $K$ le corps des fonctions de $X$, et $K^\sep$ une clôture séparable de $K$. Notons $G=\Gal(K^\sep|K)$. Soit $\bar x$ un point fermé de $X$. Notons $K_{\bar x}$ le corps des fractions de l'anneau strictement hensélien $\OO_{X,\bar x}\subset K$. Soit $K_{\bar x}^\sep$ une clôture séparable de $K_{\bar x}$. Le choix d'un plongement $K^\sep\to K_{\bar x}^\sep$ détermine une place $\m$ de $K^\sep$. Le groupe de décomposition $D_{\m}\subset G$ de $\m$ s'identifie à $\Gal(K_{\bar x}^\sep|K_{\bar x})$ ; c'est aussi, puisque $k$ est algébriquement clos, le groupe d'inertie $I_{\m}$, que nous noterons encore $I$. Rappelons que $G=\lim_\lambda G_\lambda$, où les $G_\lambda$ sont les quotients finis de $G$, c'est-à-dire les groupes d'automorphismes des revêtements finis génériquement galoisiens de $X$. Le résultat suivant découle alors de la proposition \ref{prop:gpinert}.

\begin{prop}\cite[0BUA]{stacks} Il y a un morphisme canonique surjectif $I\to\lim_{p\nmid n}\mu_n(k)$. Son noyau $P$ est trivial si $p=0$, et un pro-$p$-groupe sinon.
\end{prop}
Les groupes $I$, $P$, $I_\ttt\coloneqq I/P$ sont encore appelés groupe d'inertie (resp. d'inertie sauvage, resp. d'inertie modérée) en $x$.  La proposition suivante est analogue au corollaire \ref{cor:cohinert}.

\begin{prop}\cite[Prop. 8.1.4]{fulei}\label{prop:cohinertinf} Avec les notations ci-dessus, soit $M$ un $\Lambda[I]$-module. Alors \[ \HH^i(I,M)=\left\lbrace\begin{array}{lcl}
M^I & \qquad&\text{ si }i=0\\
M_I(-1) & & \text{ si }i=1\\
0 & & \text{ si }i\geqslant 2.
\end{array}\right.\]
\end{prop}

\begin{cor} Soit $f\colon Y\to X$ un revêtement galoisien de courbes lisses sur $k$. Soient $y$ un point fermé de $Y$ et $x=f(y)$. Supposons que l'indice de ramification de $f$ en $y$ soit divisible par $n$. Notons $I\subset \Gal(K^\sep|K)$ le groupe d'inertie en $x$, et $I_y\subset \Aut(Y|X)$ le quotient correspondant. Soit $M$ un $\Lambda[I]$-module. Alors le morphisme composé
\[ \tau_{\leqslant 1}\RG(I_y,M)\to \tau_{\leqslant 1}\RG(I,M)\to \RG(I,M)\]
est un isomorphisme dans $\DD^b_c(X,\Lambda)$.
\begin{proof}
Les morphismes en degrés 0 et 1 sont des isomorphismes d'après les propositions précédentes ; en degré supérieur, les groupes de cohomologie de $\RG(I,M)$ sont nuls.
\end{proof}
\end{cor}

\section{Groupe de Picard et variété jacobienne}

Nous verrons dans la suite que le premier groupe de cohomologie d'un faisceau constant sur une courbe lisse sur $k$ est isomorphe à un groupe de points de torsion de la jacobienne de cette courbe.

\subsection{Foncteur de Picard et variété jacobienne}\label{subsec:foncpic}

\begin{df} Le groupe de Picard d'un schéma $X$ est le groupe des classes d'isomorphisme de $\OO_X$-modules (pour la topologie de Zariski) inversibles, avec pour loi de groupe le produit tensoriel. Soit $f\colon X\to S$ un morphisme propre de schémas. Le foncteur de Picard $\Pic_{X/S}$ (resp. $\Pic^0_{X/S}$) est le faisceau associé au préfaisceau $T\mapsto \Pic(X\times_S T)$ (resp. $T\mapsto \Pic^0(X\times_S T)$) sur le gros site étale de $S$.
\end{df}
Si $f$ admet une section alors pour tout $S$-schéma $T$, les sections de ce faisceau sur $T$ sont données par $\Pic_{X/S}(T)=\Pic(X\times_S T)/\Pic(T)$ \cite[8.1, Th. 4]{blr}.

\begin{theorem}\cite[8.2, Th. 1]{blr} Soit $f\colon X\to S$ un morphisme projectif de présentation finie. Supposons $f$ plat à fibres géométriques intègres. Alors le foncteur $\Pic_{X/S}$ est représentable par un $S$-schéma séparé localement de présentation finie sur $S$.
\end{theorem}
Dans le cas des courbes projectives lisses sur un corps, le schéma qui représente le foncteur de Picard est encore une variété projective.
\begin{theorem}\cite[III, Th. 1.6]{milneAV} Soit $X_0$ une courbe projective lisse de genre $g$ sur $k_0$ munie d’un point rationnel. Alors le foncteur de Picard $\Pic^0_{X/k_0}$ est représenté par une variété abélienne $J_{X_0}$ de dimension $g$ appelée variété jacobienne de $X_0$. 
\end{theorem}

Soit $X_0$ une courbe projective lisse connexe de genre $g$ sur $k_0$. La première construction de la variété jacobienne $J_{X_0}$ est due à Weil \cite{weil_jacobienne}. Elle consiste à construire une loi de groupe birationnelle sur la puissance symétrique $X_0^{(g)}$ exprimant l'addition des diviseurs, puis à montrer qu'il existe une variété abélienne $G$ et un morphisme birationnel $X_0^{(g)}\to G$ compatible aux lois de groupe. 
Une deuxième construction, plus directe en ce qu'elle évite l'intermédiaire de la loi de groupe birationnelle, a été donnée par Chow en 1954 \cite{chow_jacobian}. Enfin, Anderson propose dans \cite{anderson} une construction différente de la variété jacobienne, sans donner de bornes sur le nombre d'étapes de la construction.\\

Il existe donc plusieurs algorithmes calculant des équations de la jacobienne de n'importe quelle courbe lisse ; cependant, pour aucun d'entre eux, nous ne disposons de bornes sur sa complexité. Dans le cas des courbes hyperelliptiques, une description explicite de la jacobienne est connue \cite[§2]{mumford_tata}. La complexité du calcul représente un obstacle conséquent à cette construction : par exemple, la jacobienne d'une courbe de genre 2 est décrite par Cassels et Flynn comme intersection de 72 quadriques dans $\PP^{15}$ \cite[Ch. 2, §3]{cassels_flynn}. Cependant, il est beaucoup plus simple de calculer avec des classes de diviseurs sans se soucier de la structure de variété de la jacobienne. Les algorithmes effectuant ces calculs sont présentés dans l'annexe \ref{sec:divrr}.
Nous n'aurons besoin que de résultats précis concernant par exemple la $n$-torsion de la jacobienne. La proposition suivante est un cas particulier d'un résultat général sur les variétés abéliennes \cite[I, Th. 7.2]{milneAV}.

\begin{prop} Soit $X$ une courbe intègre projective lisse sur $k$ de genre $g$. Soit $n$ un entier inversible dans $k$. La multiplication par $n$ sur $J_X$ est une isogénie étale de degré $n^{2g}$, et les $k$-points de son noyau $J_X[n]$ forment un $\ZZ/n\ZZ$-module libre de rang $2g$.
\end{prop}

\subsection{Jacobienne généralisée}

Soient $X_0$ une courbe lisse géométriquement connexe sur $k_0$, et $X=X_0\times_{k_0}k$. Notons $K$ le corps des fonctions de $X$. Soit $\mathfrak{m}=\sum_P m_PP\in \Div(X)$ un diviseur invariant sous l'action du groupe $\mathfrak{G}_0=\Gal(k|k_0)$. Notons $|\m|$ son support. \'{E}tant donné une fonction $f\in K^\times$, on dit que $f\equiv 1\mod\m$ si pour tout $P\in |\m|$, $v_P(1-f)\geqslant m_P$. Deux diviseurs $D,D'\in \Div(X)$ sont dits $\m$-équivalents s'il existe $f\in K^\times$ telle que $D'=D+\div(f)$ et $f\equiv 1\mod\m$. Notons $\Div_\m(X)\coloneqq\Div(X-|\m|)$. 

\begin{df}\label{def:Picm} Le groupe de Picard $\Pic_\m(X)$ de $X$ relativement à $\m$ est le quotient de $\Div_\m(X)$ par le sous-groupe des diviseurs de fonctions congrues à 1 modulo $\m$. Notons encore $\Pic^0_\m(X)$ le sous-groupe de $\Pic_\m(X)$ des classes de diviseurs de degré 0.
\end{df}

Comme le diviseur $\m$ est défini sur $k_0$, les groupes $\Pic_\m(X)$ et $\Pic^0_\m(X)$ sont naturellement munis d'une action du groupe $\mathfrak{G}_0$. Remarquons que la classe d'équivalence modulo $\m$ d'un diviseur $D$ est incluse dans sa classe d'équivalence linéaire usuelle, ce qui permet de définir un morphisme surjectif $\Pic_\m(X)\to\Pic(X)$. Nous supposons dans toute la suite que $\m$ est réduit, c'est-à-dire que tous les coefficients non nuls de $\m$ soient égaux à 1. Au vu de sa définition, le groupe $\Pic_\m(X)$ décrit alors les faisceaux inversibles sur $X$ triviaux sur $|\m|$. Nous supposons dans toute la suite de cette section que $X$ est projective.

\begin{lem}\label{lem:Picm}\cite[V.13, Prop. 7]{serrega} Rappelons que $\m$ est réduit. Notons $P_1,\dots,P_s$ les points de son support. Il y a des suites exactes de groupes abéliens
\[ 0 \to \frac{(k^\times)^{s}}{k^\times}\to \Pic_\m(X)\to \Pic(X)\to 0\]
et 
\[ 0 \to \frac{(k^\times)^{s}}{k^\times}\to \Pic^0_\m(X)\to \Pic^0(X)\to 0\]
où $k^\times$ agit par la diagonale sur $(k^\times)^{s}$, et la flèche de gauche associe à $(\lambda_1,\dots,\lambda_s)$ le diviseur d'une fonction $f$ vérifiant $f(P_i)=\lambda_i$ pour tout $i$.
\begin{proof} Un diviseur sur $X-\m$ est un diviseur principal sur $X$ si et seulement s'il est de la forme $\div(f),f\in K^\times$ où $f$ n'a ni zéro ni pôle sur le support de $\m$ ; il est alors de degré 0. Deux diviseurs $\div(f),\div(g)$ de cette forme sont $\m$-équivalents si et seulement si $f(P_i)=g(P_i)$ pour tout $i\in \{1\dots s\}$, d'où l'exactitude au milieu. De même, le diviseur d'une fonction $f$ est $\m$-équivalent à 0 si et seulement s'il prend la même valeur en tous les points de $\m$, d'où l'injectivité à gauche. La surjectivité du morphisme de droite découle du fait que tout diviseur sur $X$ est équivalent à un diviseur de support disjoint de $\m$ (voir annexe \ref{subsec:moving} ou \cite[§1.3]{shafarevich1}). Par construction, les morphismes de cette suite sont compatibles à l'action de $\Gk$.
\end{proof}
\end{lem}
La proposition suivante permettra de calculer le sous-groupe de $n$-torsion de $\Pic^0_\m(X)$.

\begin{prop}\label{prop:H1m} Il y a une suite exacte courte de $\Lambda$-modules \[ 0\to \frac{\mu_n(k)^{|\m|}}{\mu_n(k)}\to \Pic^0_\m(X)[n] \to \Pic^0(X)[n]\to 0.\]
\begin{proof} Le diagramme commutatif à lignes exactes 
\[
\begin{tikzcd}
0 \arrow[r]& \frac{(K^\times)^{|\m|}}{K^\times}\arrow[r]\arrow[d,"n"] & \Pic^0_\m(X) \arrow[r]\arrow[d,"n"] & \Pic^0(X) \arrow[r]\arrow[d,"n"] & 0 \\
0 \arrow[r]& \frac{(K^\times)^{|\m|}}{K^\times}\arrow[r] & \Pic^0_\m(X) \arrow[r] & \Pic^0(X) \arrow[r] & 0
\end{tikzcd}
\]
fournit, par surjectivité de la mise à la puissance $n$ sur $(K^\times)^{|\m|}/K^\times$, la suite exacte des noyaux :
\[ 0\to \frac{\mu_n(k)^{|\m|}}{\mu_n(k)}\to \Pic^0_\m(X)[n] \to \Pic^0(X)[n]\to 0.\]
\end{proof}
\end{prop}

\begin{rk} Il existe un groupe algébrique $J_\m$ sur $k_0$ muni d'une application rationnelle $X_0\to J_\m$ régulière sur $X-|\m|$ qui induit un isomorphisme $\Pic^0_\m(X)\to J_\m(k)$ : c'est la jacobienne généralisée de $X$ relativement à $\m$ \cite[V.9, Th. 1]{serrega}. Il paraît hors de portée pour l'instant de déterminer explicitement en temps raisonnable des équations définissant cette jacobienne généralisée ; la proposition précédente montre  toutefois comment calculer dans son groupe de $n$-torsion.
\end{rk}

\subsection{Les $\mu_n$-torseurs sur une courbe lisse}

\begin{prop}\label{prop:muntorseurslisses} Soit $X$ une courbe intègre lisse sur $k$, de corps des fonctions $K$.
Le groupe $\HH^1(X,\mu_n)$ est canoniquement isomorphe au quotient du groupe \[ \{ (D,f)\in \Div(X)\times K^\times\mid nD=\div(f)\} \] par le sous-groupe des $(D,f)$ où $f\in (K^\times)^n$.
\begin{proof}
Soit $G$ le groupe quotient en question. Remarquons que si un couple $(D,f)$ est nul dans $G$ alors $D$ est un diviseur principal. Rappelons que d'après la proposition \ref{prop:H1faiscinv}, le groupe $\HH^1(X,\mu_n)$ est canoniquement isomorphe au groupe des classes d'isomorphisme de couples $(\mathcal{L},\alpha)$ où $\mathcal{L}$ est un faisceau inversible sur $X$ et $\alpha\colon \mathcal{L}^{\otimes n}\xrightarrow{\sim}\OO_X$.
Considérons le morphisme $F\colon G\to \HH^1(X,\mu_n)$ défini par $F(D,f)=(\OO_X(D),\OO_X(nD)\xrightarrow{m_f} \OO_X)$ où $m_f$ désigne la multiplication par $f$.  \begin{itemize}[label=$\bullet$]
\item Injectivité : si $F(D,f)$ est égal dans $\HH^1(X,\mu_n)$ à $(\OO_X,\id)$ alors il existe un isomorphisme $\phi \colon \OO_X(D)\to\OO_X$ tel que $\phi^n$ soit la multiplication par $f$. Un tel isomorphisme étant nécessairement la multiplication par un élément de $K^\times$, il en résulte que $f$ est une puissance $n$-ième dans $K^\times$.
\item Surjectivité : soit $(\LL,\alpha)\in S$. Il existe un diviseur $D\in\Div(X)$ et un isomorphisme de faisceaux $\phi \colon \LL\to\OO_X(D)$. En particulier, ceci livre un isomorphisme $\alpha\circ\phi^{-n}\colon \OO_X(nD)\to \OO_X$, qui est la multiplication par une section globale $f\in K^\times$ de $\OO_X(nD)$. Alors $(\LL,\alpha)$ est égal dans $S$ à $F(D,f)$.
\end{itemize}
\end{proof}
\end{prop}

\begin{cor}\label{cor:H1gen} Sous les mêmes hypothèses, $\HH^1(X,\mu_n)$ est isomorphe au quotient du groupe \[ \{ (A,f)\in \Pic(X)\times K^\times \mid \exists D\in \Div(X)\colon [D]=A\text{ et }  nD=\div(f)\}\] par le sous-groupe des $(A,f)$ avec $f\in (K^\times)^n$. Si de plus $X$ est projective alors $\HH^1(X,\mu_n)$ est isomorphe au sous-groupe de $n$-torsion $\Pic(X)[n]$ de $\Pic(X)$.
\begin{proof}
Le morphisme $(D,f)\mapsto ([D],f)$ défini par passage au quotient est évidemment surjectif et a pour noyau les couples $(D,f)$ tels que $f\in (K^\times)^n$. Si $X$ est projective alors l'application surjective $\HH^1(X,\mu_n)\to\Pic(X)[n]$ qui à un couple $(A,f)$ associe $A$ est également injective, puisque deux fonctions rationnelles sur $X$ de même diviseur sont égales à un élément de $k^\times$ près.
\end{proof}
\end{cor}

\begin{lem}\label{lem:H1aff} Supposons $X$ projective. Soient $U$ un ouvert de $X$, et $Z$ le fermé réduit complémentaire. Notons $\Div^0_Z(X)$ le sous-groupe de $\Div^0(X)$ formé des diviseurs à support dans $Z$.
Considérons les triplets $(A,D',f)$ où $A\in\Pic^0(X)$, $D'\in\Div^0_Z(X)$, $f\in K^\times$ et il existe un diviseur $\bar D$ de classe $A$ tels que $n\bar{D}=\div(f)+D'$ dans $\Div^0(X)$. Ces triplets forment un sous-groupe de $\Pic^0(X)\times\Div^0_Z(X)\times K^\times$. Soit $H$ le quotient de ce sous-groupe par celui des $([D'],nD',f)$. Alors $\HH^1(U,\mu_n)$ est canoniquement isomorphe à $H$.
\begin{proof}
Considérons le morphisme \[ \begin{array}{rcl} H&\longrightarrow& \HH^1(U,\mu_n) \\ (A,D',f)&\longmapsto&(A|_{U},f) \end{array}\] où $\HH^1(U,\mu_n)$ est identifié au groupe décrit dans le corollaire précédent. Notons que comme $n\bar{D}=\div(f)+D'$, la restriction à $U$ donne bien $n\bar{D}_{|U}=\div(f)|_U$. De plus, si $f$ est une puissance $n$-ième alors il existe $h\in K^\times$ tel que $n\bar{D}=\div(h^n)+D'$, autrement dit $D'=n(\bar{D}-\div(h))$, et $([D'],D',f)=([\bar{D}-\div(h)],n(\bar{D}-\div(h)),f)$. Reste à montrer la surjectivité de ce morphisme. Soit donc $([D],f)\in \HH^1(U,\mu_n)$, avec $D\in \Div(U)$. Alors $nD=\div(f)|_{U}$. Soit $\bar{D}\in \Div^0(X)$ un diviseur de restriction $D$. Alors $E\coloneqq n\bar{D}-\div(f)$ est un diviseur à support dans $Z$, de degré nul. Par conséquent, $([D],f)$ est l'image de $([\bar{D}],E,f)\in H$.
\end{proof}
\end{lem}

\section{Groupe de Picard des courbes nodales}

\subsection{Groupe de Picard des courbes singulières}\label{subsec:picsing}

Soit $X_0$ une courbe intègre sur $k_0$, et $X=X_0\times_{k_0} k$. Soit $\pi \colon\tX\to X$ sa normalisation, qui est finie \cite[035S]{stacks}. 
Supposons que $X$ a un unique point singulier $P$, qui est alors nécessairement défini sur $k_0$, et que les points de $\tX$ au-dessus de $P$ sont également définis sur $k_0$. 
Considérons la suite exacte de faisceaux sur $X$ :
\[ 0\to \GG_{m,X} \to \pi_\star \GG_{m,\tX} \to (\pi_\star \GG_{m,\tX})/\GG_{m,X}\to 0.\]
Comme $\pi$ est un isomorphisme en-dehors du lieu singulier de $X$, le support du faisceau quotient $\mathcal{Q}\coloneqq(\pi_\star \GG_{m,\tX})/\GG_{m,X}$ est $\{ P\}$. Par conséquent, la suite exacte longue en cohomologie de cette suite exacte courte est :
\[ 0\to \HH^0(X,\GG_m) \to \HH^0(\tX,\GG_m) \to \HH^0(X,\mathcal{Q}) \to \Pic(X)\to \Pic(\tX) \to 0.\]
Lorsque $X$ est projective, $\HH^0(X,\GG_m)=\HH^0(\tX,\GG_m)=k^\times$, ce qui fournit une suite exacte : \[ 0\to \HH^0(X,\mathcal{Q})\to \Pic(X)\to \Pic(\tX)\to 0.\]

\begin{rk} La normalisée de la cubique nodale $\Spec k[x,y]/(y^2-x^2(x+1))$ ou de la cubique cuspidale $\Spec k[x,y]/(y^2-x^3)$ est la droite affine $\Spec k[t]$. Dans ces deux cas, $\HH^0(X,\GG_m)=\HH^0(\tX,\GG_m)=k^\times$, et la suite \[ 0\to \HH^0(X,\mathcal{Q})\to \Pic(X)\to \Pic(\tX)\to 0\]
est exacte.
\end{rk}

\begin{ex} Voici un exemple de courbe affine $X$ telle que $\HH^0(\tX,\GG_m)$ contienne strictement $\HH^0(X,\GG_m)$. Considérons l'ouvert $\{x\neq 0\}$ de la cubique nodale précédente, c'est-à-dire \[ X=\Spec k[x,y,z]/(y^2-x^2(x+1),zx-1).\]
Sa normalisée est l'ouvert de $k[t]$ situé au-dessus de $\{ x\neq 0\}$, c'est-à-dire $\tX=\Spec k[t,s]/(st^2-1)$. Le morphisme de normalisation est donné par $x\mapsto t^2, y\mapsto t^3, z\mapsto s$.
La fonction $t=\frac{y}{x}$ est inversible sur $\tX$ d'inverse $st$, alors que $\frac{y}{x}\not\in \HH^0(X,\OO_X)^\times$.
\end{ex}

\subsection{Courbes nodales}

Soit $X_0$ une courbe sur $k_0$. Notons $X=X_0\times_{k_0}k$, où $k$ désigne toujours la clôture algébrique de $k$.

\begin{df}\cite[0C47]{stacks} Un point fermé $P_0$ de $X_0$ est dit nodal s'il existe un point $P$ de $X$ d'image $P_0$ tel que le complété de l'anneau local $\OO_{X,P}$ soit isomorphe à $k[[x,y]]/(xy)$. Une courbe sera dite nodale si elle est singulière et tous ses points singuliers sont nodaux.
\end{df}

\begin{rk} Le critère jacobien montre que si $X$ possède un unique point singulier $P$ alors $X_0$ possède un unique point singulier $P_0$, qui est un $k_0$-point. Si $P$ est nodal alors $P_0$ l'est aussi.
\end{rk}

Voyons comment reconstruire une courbe nodale à partir de sa normalisation.
\begin{df}\'{E}tant donné une courbe lisse $Y$ sur $k$ et des points deux à deux distincts \[ Q_1,R_1,\dots,Q_s,R_s\in Y(k)\] nous noterons $Y_{Q_1=R_1,\dots,Q_s=R_s}$, ou $Y_{\underline Q=\underline R}$, la courbe obtenue en identifiant pour tout $i\in \{1\dots s\}$ les deux points $Q_i$ et $R_i$. Sa construction est décrite dans \cite[IV,§4]{serrega}. 
\end{df}

\begin{prop}\cite[IV.3, Prop. 2) et IV.4, Exemple b]{serrega} Avec les notations de la définition, si $Y$ est connexe, la courbe $Y_{Q_1=R_1,\dots,Q_s=R_s}$ est irréductible, de normalisation $Y$. Elle possède $s$ points singuliers nodaux.
\end{prop}

\begin{rk}\label{rk:conetg} Pour $i\in \{1\dots s\}$, notons $P_i$ l'image dans $Y_{\underline Q=\underline R}$ de $Q_i$ et $R_i$. Les morphismes $C_{Q_i}Y\to C_{P_i}(Y_{\underline Q=\underline R})$ et $C_{R_i}Y\to C_{P_i}(Y_{\underline Q=\underline R})$ entre cônes tangents induisent un isomorphisme \[ (C_{Q_i}Y\sqcup
C_{R_i}Y)_{Q_i=R_i}\xrightarrow{\sim} C_{P_i}(Y_{\underline Q=\underline R})\]
où $\sqcup$ désigne le coproduit de schémas.
\end{rk}

\subsection{Groupe de Picard d'une courbe nodale}

Soit $X_0$ une courbe sur $k_0$, et $X=X_0\times_{k_0}k$.
Supposons que $X$ soit intègre, nodale et possède un unique point singulier $P$. Notons $K$ le corps des fonctions de $X$. Soient $\pi \colon\tX\to X$ sa normalisation, et $j\colon\tX\to \bar X$ la complétion projective lisse de $\tX$. Il y a exactement deux points $Q,R$ de $\tX$ au-dessus de $P$ \cite[0CBW]{stacks}.
Remarquons que pour un ouvert $U$ de $X$ contenant $P$, \[ \OO_X(U)=\{ f\in \OO_{\tX}(U\times_X \tX)\mid f(Q)=f(R) \}\] (voir \cite[IV.4,Exemple b)]{serrega}). Ainsi, un fibré en droites sur $X$ est défini par la donnée d'un fibré en droites $\tL$ sur $\tX$ muni d'un isomorphisme entre les fibres $\tL(Q)$ et $\tL(R)$.\\

\begin{lem}\label{lem:picnod} Le groupe $\Pic^0 X$ est isomorphe au groupe $\Pic^0_{\m}(\tX)$ de la définition \ref{def:Picm}, où $\m$ est le diviseur effectif $Q+R$.
\begin{proof} Soit $D\in \Div(X-P)$. Considérons le faisceau inversible $\tL_D\coloneqq \OO_\tX(D)$, où l'on considère $D$ comme un diviseur sur la normalisée $\tX$ de $X$ via l'isomorphisme $X-P\xrightarrow{\sim} \tX-\{Q,R\}$. Comme la valuation de $D$ en $P$ est nulle, $\tL_D(Q)$ est canoniquement isomorphe à la fibre en $Q$ de $\OO_\tX$, elle-même canoniquement isomorphe à $k$. Il en est de même pour $\tL_D(R)$. Les fibres $\tL_D(Q)$ et $\tL_D(R)$ sont identifiées par ces isomorphismes canoniques, ce qui permet de définir canoniquement un fibré en droites $\L_D$ sur $X$.
Le morphisme de groupes $\Div(X-P)\to \Pic(X)$ qui à un diviseur $D$ associe $\L_D$ a pour noyau l'ensemble des diviseurs $D$ tels que $\tL_D$ soit isomorphe à $\OO_\tX$, c'est-à-dire les diviseurs de la forme $\div(f)\in \Div(X-P)$, où $f\in K^\times$ est définie et n'a ni zéro ni pôle en $P$. Ceci définit par passage au quotient un morphisme injectif $\Pic^0_\m(\tX)\to \Pic^0(X)$. Montrons qu'il est surjectif.
Un faisceau inversible $\L$ sur $X$ donne un faisceau inversible $\tL=\pi^\star \L$ sur $\tX$ dont les fibres en $Q$ et $R$ sont canoniquement isomorphes. Le faisceau $\tL$ est isomorphe à un faisceau $\tL_D$, où $D\in \Div(\tX)$ ; quitte à lui ajouter le diviseur d'une fonction définie en $P$, on peut supposer que $D\in \Div(\tX-\{ Q,R\})=\Div(X-P)$. 
\end{proof}
\end{lem}

\begin{lem}\label{lem:sesPic} Nous retrouvons ainsi explicitement la suite exacte courte de groupes abéliens \[ 0\to \HH^0(X,\GG_m)\to \HH^0(\tX,\GG_m)\to (k^\times\times k^\times)/k^\times \to \Pic(X)\to \Pic(\tX) \to 0 \]
déjà évoquée dans la section \ref{subsec:picsing}. 
\begin{proof}
On utilise la description de $\Pic(X)$ donnée dans le lemme \ref{lem:picnod}. La flèche \[ \HH^0(\tX,\GG_m)\to (k^\times\times k^\times)/k^\times \] associe à une fonction $f$ inversible sur $\tX$ la classe du couple $(f(Q),f(R))\in k^\times\times k^\times$ ; ce couple appartient à la diagonale si et seulement si $f$ est définie en $P$, c'est-à-dire $f\in \HH^0(X,\GG_m)$. La flèche $k^\times\times k^\times \to \Pic(X)$ associe à $(a,b)$ le diviseur d'une fonction $f$ telle que $f(Q)=a$ et $f(R)=b$. Ce diviseur est celui d'une fonction définie et non nulle en $P$ si et seulement si $a=b$. La flèche de droite est le tiré en arrière $\pi^\star$, qui à la classe de $D\in \Div(X-P)$ dans $\Pic(X)$ associe la classe de $D\in \Div(\tX-\{Q,R\})$ dans $\Pic(\tX)$. 
\end{proof}
\end{lem}

\begin{rk}
Décrivons l'action de $\mathfrak{G}_0$ sur $\HH^0(X,\pi_\star\GG_m/\GG_m)=(k^\times\times k^\times)/k^\times$. Soit $\sigma\in \mathfrak{G}_0$. Si $\sigma$ échange $Q$ et $R$ alors pour tout $(a,b)\in k^\times\times k^\times$, $\sigma\cdot [(a,b)]=[(\sigma(b),\sigma(a))]$, où les crochets désignent la classe dans le quotient par $k^\times$. Sinon, $\sigma\cdot (a,b)=(\sigma(a),\sigma(b))$. Dans toute la suite de ce chapitre, l'action de $\mathfrak{G}_0$ sur $(k^\times\times k^\times)/k^\times$ sera celle-ci.
\end{rk}

\begin{rk} Tous les résultats précédents se généralisent immédiatement au cas où une courbe projective $X$ a plusieurs singularités nodales toutes définies sur $k_0$. Notons $P_1,\dots,P_r$ ces points. Au-dessus de chaque $P_i$ se trouvent deux points $Q_i,R_i$ de $\tX$. Alors $\pi_\star\GG_m/\GG_m$ est supporté en $P_1,\dots,P_r$. Il y a de même une suite exacte \[ 0\to \left(\frac{k^\times\times k^\times}{k^\times}\right)^r\to \Pic(X)\to\Pic(\tX)\to 0.\]
Le groupe $\Pic(X)$ se décrit comme le quotient du groupe des diviseurs sur $X-\{P_1,\dots, P_r\}$ par les diviseurs de fonctions $f\in\bigcap_i \OO_{X,P_i}^\times$. Le premier morphisme de la suite exacte associe à $(a_i,b_i)_{1\leqslant i\leqslant r}$ le diviseur d'une fonction telle que $f(Q_i)=a_i$ et $f(R_i)=b_i$. On obtient enfin la suite exacte
\[ 0\to \left(\frac{\mu_n(k)\times\mu_n(k)}{\mu_n(k)}\right)^r\to \HH^1(X,\mu_n)\to \HH^1(\tX,\mu_n)\to 0.\]
\end{rk}

\subsection{Description des $\mu_n$-torseurs sur une courbe nodale}

Soient $X_0$ une courbe sur $k_0$, et $X=X_0\times_{k_0}k$.
Supposons que $X$ soit intègre, nodale et possède un unique point singulier $P$.
Soit $\bar X$ la compactification lisse de $\tX$. Rappelons (voir lemme \ref{lem:H1aff}) que $\HH^1(X,\mu_n)$ s'identifie au groupe des classes d'isomorphisme de couples $(\L,\alpha)$ où $\L$ est un fibré en droites sur $X$ et $\alpha \colon \L^{\otimes n}\to \OO_X$ est un isomorphisme. 
Soit $G$ le quotient du groupe \[\{(D,f)\in \Div(\tX)\times (\OO_{\tX,Q}^\times\cap \OO_{\tX,R}^\times)\mid nD=\div(f)|_\tX, f(Q)=f(R)\} \]
par le sous-groupe des couples $(D,f)$ tels que $f$ soit la puissance $n$-ième dans $K^\times$ d'une fonction $g$ vérifiant $g(Q)=g(R)$. Ici, les groupes d'inversibles des anneaux locaux $\OO_{\tX,Q}$ et $\OO_{\tX,R}$ sont vus comme sous-groupes de $K^\times$. Construisons une application $\Phi\colon G\to \HH^1(X,\mu_n)$. Soit $(D,f)\in G$. Alors comme $nD=\div(f)$ et $f$ est inversible en $Q$ et $R$, les points $Q$ et $R$ n'appartiennent pas au support de $D$. Par conséquent, les fibres en $Q$ et $R$ du faisceau inversible $\tL_D\coloneqq\OO_\tX(D)$ sur $\tX$ sont canoniquement isomorphes à celles de $\OO_\tX$, elles-mêmes canoniquement isomorphes à $k$. L'isomorphisme $\tL_D(Q)\to \tL_D(R)$ défini par le diagramme

\[\begin{tikzcd}
\tL_D(Q) \arrow[r,"\sim"]\arrow[d] & \OO_{\tX}(Q)\arrow[r,"\sim"] & k \arrow[d,"1_k"] \\
\tL_D(R) \arrow[r,"\sim"] & \OO_{\tX}(R) \arrow[r,"\sim"] & k
\end{tikzcd}\]
permet de définir un faisceau inversible $\L_D$ sur $X$. La multiplication par $f$ définit un isomorphisme $\tL_D^{\otimes n}\to\OO_\tX$. La condition $f(Q)=f(R)$ assure que cet isomorphisme de $\OO_\tX$-modules définit un isomorphisme de $\OO_X$-modules $m_f\colon\L_D^{\otimes n}\to\OO_X$. L'application $\Phi\colon G\to \HH^1(X,\mu_n)$ associe au couple $(D,f)$ le couple $(\L_D,m_f)$.

\begin{prop}\label{prop:H1nodbis} L'application $\Phi\colon G\to \HH^1(X,\mu_n)$ construite ci-dessus est un isomorphisme de $\ZZ/n\ZZ$-modules.
\begin{proof} Elle est bien définie car elle associe à un couple $(D,f)$ où $D=\div(f)$ et $f$ est la puissance $n$-ième d'une fonction $g\in K^\times$ vérifiant $g(Q)=g(R)$, le couple $(\L_D,m_{f})$, qui est isomorphe au couple $(\OO_X,1)$ via la multiplication par $g$ (qui définit bien un isomorphisme $\L_D\to\OO_X$ car $g(Q)=g(R)$). La linéarité de $\Phi$ découle directement de sa construction ; passons à l'injectivité. Soit $(D,f)\in G$ tel qu'il existe un isomorphisme de couples $(\L_D,m_f)\to (\OO_X,1)$. Il y a encore un isomorphisme $(\tL_D,f)\to (\OO_\tX,1)$, qui est nécessairement la multiplication par une fonction $g$ telle que $D=\div(g)$. Comme l'isomorphisme$(\tL_D,f)\to (\OO_\tX,1)$ défini par $g$ provient d'un isomorphisme $(\L_D,m_f)\to (\OO_X,1)$, la fonction $g$ vérifie encore $g(Q)=g(R)$ ; par conséquent, la classe du couple $(D,f)$ dans $G$ est $0$. Enfin, soit $(\L,\alpha)\in \HH^1(X,\mu_n)$. Le couple $(\tL,\tilde\alpha)\in \HH^1(\tX,\mu_n)$ qui s'en déduit est isomorphe à $(\OO_\tX(D),g)$ pour un diviseur $D$ sur $\tX$ et une fonction $g\in K^\times$ telle que $nD=\div(g)$. Supposons, quitte à ajouter à $D$ un diviseur principal, que le support de $D$ ne contient pas $Q$ et $R$. Par conséquent, $(\L,\alpha)$ est égal dans $\HH^1(X,\mu_n)$ à $(\L_D,m_g)$.
\end{proof}
\end{prop}

Supposons maintenant $X$ affine. Soit $Z$ le fermé réduit complémentaire de $\tX$ dans $\bar X$. Notons $H$ le quotient du groupe 
\[ \{ (\bar D,D',f)\in \Div^0(\bar X)\times \Div^0_Z(X)\times (\OO_{\bar X,Q}^\times\cap\OO_{\bar X,R}^\times)\mid n\bar D+D'=\div(f),f(Q)=f(R)\} \]
par le sous-groupe des $(\bar D,D',f)$ où $f$ est la puissance $n$-ième dans $K^\times$ d'une fonction $g$ vérifiant $g(Q)=g(R)$. Considérons l'application $\Psi\colon H\to \HH^1(X,\mu_n)$ qui à $(\bar D,D',f)$ associe $(\bar D|_\tX,f)$.

\begin{prop} L'application $\Psi\colon H\to \HH^1(X,\mu_n)$ définie ci-dessus est un isomorphisme de $\Lambda[\mathfrak{G}_0]$-modules.
\begin{proof} La linéarité et l'injectivité de $\Psi$ découlent immédiatement de sa construction. Soit maintenant $(D,f)\in \HH^1(X,\mu_n)$. Ce couple vérifie $nD=\div(f)|_\tX$. Soit $\bar D$ un diviseur de degré 0 sur $\bar X$ de restriction $D$. Posons $D'=n\bar D-\div (f)$. C'est un diviseur de degré 0 et de support inclus dans $Z$. Par conséquent, $(\bar D,D',f)\in H$ a pour image $(D,f)$ par $\Psi$.
\end{proof}
\end{prop}

\section{Cohomologie des faisceaux constants sur les courbes}

\subsection{Courbes lisses sur un corps algébriquement clos}

\subsubsection{Courbes projectives}

Soit $X$ une courbe intègre lisse sur un corps algébriquement clos $k$. Rappelons que $n$ désigne un entier premier à la caractéristique de $k$.
\begin{prop}\cite[Arcata, Prop. 3.1]{sga412} \label{prop:H1Gm} Il y a des isomorphismes canoniques
\[
\HH^i(X,\GG_m)=\left\lbrace \begin{array}{lcl}
\Gamma(X,\OO_X)^\times &\text{ si }& i=0\\
\Pic(X) &\text{ si }&i=1\\
0& \text{ si }&i\geqslant 2.
\end{array}\right.
\]
\end{prop}

Supposons désormais $X$ projective ; ceci entraîne que le morphisme $\GG_m\xrightarrow{f\mapsto f^n}\GG_m$ est, sur les sections globales, la mise à la puissance $n$ sur $k^\times$. La suite exacte de Kummer 
\[ 0\to \mu_n\to \GG_m \xrightarrow{n} \GG_m \to 0\]
donne alors la suite exacte longue
\[ 0\to \HH^1(X,\mu_n)\to\Pic(X)\xrightarrow{n}\Pic(X) \to \HH^2(X,\mu_n)\to 0.\]
La théorie des variétés abéliennes nous enseigne que la multiplication par $n$ sur la jacobienne $J_X$ est surjective. La suite exacte
\[ 0\to J_X(k)\to \Pic(X)\xrightarrow{\deg}\ZZ \to 0\]
permet d'en déduire que le conoyau de la multiplication par $n$ sur $\Pic(X)$ est $\ZZ/n\ZZ$. 
\begin{theorem}\cite[03RQ]{stacks}\label{th:cohproj} Soit $X$ une courbe intègre projective lisse sur $k$.  Les groupes de cohomologie de $X$ à valeurs dans $\mu_n$ sont canoniquement isomorphes aux groupes suivants.
\[
\HH^i(X,\mu_n)=\left\lbrace \begin{array}{lcl}
\mu_n(k) &\text{ si }& i=0\\
J_X[n] &\text{ si }&i=1\\
\ZZ/n\ZZ &\text{ si }&i=2\\
0& \text{ si }&i\geqslant 3.
\end{array}\right.
\]
Si $f\colon Y\to X$ est un morphisme de courbes intègres projectives lisses sur $k$ alors le morphisme induit $f^\star\colon \HH^2(X,\mu_n)\to \HH^2(Y,\mu_n)$ est la multiplication par $\deg(f)$.
\end{theorem}

\begin{lem}\label{lem:corpsdef} Supposons que $X$ provient par changement de base d'une courbe définie sur le sous-corps parfait $k_0$ de $k$. L'action de $\Gal(k|k_0)$ sur $\HH^1(X,\mu_n)$ se factorise par un quotient d'ordre $n^{O(g^2)}$, où $g$ désigne le genre de $X$. Par conséquent, si $k_0$ est fini, il existe une extension $k_1/k_0$ de degré $n^{O(g^2)}$ et des diviseurs $k_1$-rationnels $D_1,\dots,D_{2g}$ formant une base de $\HH^1(X,\Lambda)$.
\begin{proof} La première assertion vient par passage au quotient du morphisme \[\Gal(k|k_0)\to \Aut_\Lambda(\HH^1(X,\mu_n))\] sachant que $\HH^1(X,\mu_n)$ est un $\Lambda$-module libre de rang $2g$. 
La deuxième phrase découle du fait qu'une classe de diviseurs $k_1$-rationnelle contient toujours un diviseur $k_1$-rationnel (voir lemme \ref{lem:PicDiv} en annexe).
\end{proof} 
\end{lem}

\subsubsection{Action du Frobenius et comptage de points}

Supposons que $k_0$ soit un corps fini $\FF_q$, et que $X$ provienne par changement de base d'une courbe $X_0$ sur $k_0$. Alors $X$ est munie de l'action du morphisme de Frobenius géométrique $\Frob_q$. L'action de $\Frob_q$ sur $\mu_n(k)$ est la mise à la puissance $q$, et celle sur $\HH^2(X,\mu_n)$ est l'identité puisqu'un automorphisme ne change pas le degré des diviseurs. Ainsi, en notant $t$ la trace de l'endomorphisme de Frobenius sur $\HH^1(X,\mu_n)$, la formule des traces (théorème \ref{th:trace}) assure que le nombre de $k_0$-points de $X$ est donné par :
\[ \# X(k_0)=1+q-t\mod n.\]
La fonction zêta de $X_0$ est alors \[ Z_{X_0}(t)=\frac{L(t)}{(1-t)(1-qt)}\]
où $L(t)=\det(\id-t\Frob_q^\star\mid \HH^1(X,\QQ_\ell))\in \ZZ[t]$ est un polynôme de degré $2g$ dont les racines complexes ont pour module $\sqrt{q}^{-1}$. Ceci permet également de compter le nombre de $\FF_q$-points de la jacobienne $J_X$ de $X$. En effet, comme montré dans \cite[VIII, Cor. 6.3]{lorenzini} :\[ \#J_X(\FF_q)=\# \ker(\id-t\Frob_q)=\det(1-\Frob_q^\star\mid \HH^1(C,\QQ_\ell))=L(1). \]
En particulier, \[ \#J_X(\FF_{q^r})\sim_{r\to\infty} q^{rg}.\]

\subsubsection{Courbes affines lisses : suite de Gysin}\label{subsubsec:Gyscst}

Soit $X$ une courbe projective connexe lisse sur $k$. Soient $U$ un ouvert de $X$, et $Z$ le fermé réduit complémentaire. Seuls les groupes $\HH^0(U,\mu_n)=\mu_n(k)$ et $\HH^1(U,\mu_n)$ sont non nuls. Rappelons (voir lemme \ref{lem:H1aff}) qu'un élément de $\HH^1(U,\mu_n)$ est représenté par un triplet $(D,D',f)\in \Div^0(X)\times \Div^0_Z(X)\times K^\times$ vérifiant $nD=\div(f)+D'$. La suite exacte de Gysin est alors la suite exacte de $\Lambda$-modules libres \cite[03RR]{stacks}\[
0\to \HH^1(X,\mu_n)\xrightarrow{\phi} \HH^1(U,\mu_n)\xrightarrow{\psi} \HH^0(Z,\Lambda)\xrightarrow{\Sigma} \Lambda\to 0
\]
où les flèches sont décrites, avec les notations précédentes, par : \begin{itemize}[label=$\bullet$]
\item $\phi([D],f)=([D],0,f)$
\item $\psi([D],D',f)=D'\mod n$
\item $\Sigma((\alpha_P)_{P\in Z})=\sum_{P\in Z}\alpha_P$.
\end{itemize}

On peut également décrire explicitement la fonctorialité en la paire $(X,U)$ de cette suite. Soient $\phi \colon X'\to X$ un morphisme de courbes projectives lisses, $U'=X'\times_X U$ et $Z'=X'\times_X Z$. Alors le morphisme $\phi^\star \colon \HH^1(U,\mu_n)\to \HH^1(U',\mu_n)$ défini par $\phi^\star([D],D',f)=([\phi^\star D],\phi^\star D',\phi^\star f)$, où $\phi^\star$ désigne encore les tirés en arrière usuels, s'insère dans le diagramme suivant.
\[\begin{tikzcd}
0 \arrow[r]& \HH^1(X,\mu_n)\arrow[d,"\phi^\star"]\arrow[r]& \HH^1(U,\mu_n)\arrow[r]\arrow[d,"\phi^\star"]& \HH^0(Z,\Lambda)\arrow[r]\arrow[d,"\phi^\star"]& \HH^2(X,\mu_n)\arrow[d,"\phi^\star"] \arrow[r]& 0 \\
0 \arrow[r] & \HH^1(X',\mu_n)\arrow[r]& \HH^1(U',\mu_n)\arrow[r]& \HH^0(Z',\Lambda)\arrow[r]& \HH^2(X',\mu_n) \arrow[r]&  0
\end{tikzcd}\]
Ici, le morphisme $\HH^2(X,\mu_n)\to \HH^2(X',\mu_n)$ est la multiplication par le degré de $\phi$ \cite[0AMB]{stacks}.

\subsection{Dualité de Poincaré pour les courbes lisses}\label{subsubsec:courbcohsupp}

Soit $X_0$ une courbe projective lisse sur $k_0$. Soit $j_0\colon U_0\to X_0$ l'inclusion d'un ouvert strict, et $i_0\colon Z_0\to X_0$ l'inclusion du fermé réduit complémentaire. 
Notons $U,X,Z,i,j$ leurs changements de base à $k$.
La suite exacte de faisceaux sur $X$ \[ 0\to j_!\mu_n\to \mu_n\to i_\star\mu_n\to 0 \]
montre que $\HH^0_c(U,\mu_n)=\ker(\mu_n(k)\to \mu_n(k)^{|Z|})=0$ et $\HH^2_c(U,\mu_n)=\HH^2(X,\mu_n)=\Lambda$. Dans le cas des faisceaux constants, le seul calcul à effectuer est donc celui de $\HH^1_c(U,\mu_n)$. \\

Soient $P_1,\dots,P_r$ les points fermés de $Z$. Le diviseur $\mathfrak{m}=\sum_iP_i$ sur $X$ est $\mathfrak{G}_0$-invariant puisque $Z$ provient de $Z_0$. Considérons le faisceau $\GG_{m,Z}\coloneqq\ker(\GG_m\to i_\star\GG_m)$ des fonctions congrues à 1 modulo $\m$.
Le groupe $\HH^1(X,\GG_{m,Z})$ classifie les faisceaux inversibles sur $X$ trivialisés sur $Z$ \cite[Arcata, §2.3]{sga412} et est donc isomorphe à $\Pic_\m(X)$. La suite longue en cohomologie qui se déduit de cette suite exacte courte est la suite \[ 0\to \frac{\HH^0(Z,\GG_m)}{\HH^0(X,\GG_m)}\to \Pic_\m(X)\to \Pic(X)\to 0\]
du lemme \ref{lem:Picm}.

\begin{lem}\cite[Arcata, 2.3.(a)]{sga412} \label{lem:H1c} Il y a un isomorphisme de $\Lambda[\mathfrak{G}_0]$-modules \[ \HH^1_c(U,\mu_n)\xrightarrow{\sim} \Pic^0_\m(X)[n].\]
\begin{proof}
La suite exacte de faisceaux sur $X$  \[ 0\to j_!\mu_n\to \GG_{m,Z}\xrightarrow{n} \GG_{m,Z}\to 0.\]
donne la suite exacte longue de $\Lambda[\mathfrak{G}_0]$-modules \[ \HH^1_c(U,\mu_n) \to \HH^1(X,\GG_{m,Z}) \xrightarrow{n} \HH^1(X,\GG_{m,Z})\to 0.\]
Comme le groupe $H^0(X,\GG_{m,Z})$ est trivial, le groupe $\HH^1_c(U,\mu_n)$ est donc isomorphe à $\HH^1(X,\GG_{m,Z})[n]=\Pic_\m(X)[n]=\Pic^0_\m(X)[n]$.
\end{proof}
\end{lem}

Dans le cas où $U=X$, l'accouplement
\[ \HH^1(U,\mu_n)\times \HH^1_c(U,\mu_n) \to \HH^2_c(U,\mu_n^{\otimes 2})\xrightarrow{\sim} \mu_n(k)\]
provient de l'autodualité de la jacobienne \cite[Arcata, §2.3]{sga412} et est appelé accouplement de Weil. Sa construction explicite se trouve dans \cite[§20, p184]{mumfordAV}. La construction suivante généralise cette dernière au cas où $U$ est affine. Notons $K$ son corps des fonctions.\\

\begin{df} Soit $f\in K$. Soit $D\in\Div(X)$ tel que $f$ n'ait ni zéro, ni pôle sur le support de $D$. L'évaluation de $f$ en $D$ est définie par
\[ f(D)=\prod_{P\in |D|}f(P)^{v_P(D)}.\]
\end{df}
La loi de réciprocité suivante, due à Weil, servira dans la construction de l'accouplement.
\begin{prop}\cite[III.4, Prop. 7]{serrega} Soient $f,g\in K^\times$ deux fonctions de diviseurs disjoints. Alors \[ f(\div g)=g(\div f).\]
\end{prop}

Souvenons-nous qu'un élément du groupe $\HH^1(U,\mu_n)$ est la classe d'un triplet $([D],D',f)$ avec $[D]\in\Pic^0(X)[n]$, $D'\in\Div^0_Z(X)$ et $f\in K^\times$ vérifie $nD=D'+\div(f)$ (voir lemme \ref{lem:H1aff}). D'autre part, $\HH^1_c(U,\mu_n)$ est la $n$-torsion du groupe $\Div(U)/\{\div(f), f\equiv 1\mod Z\}$ (voir lemme \ref{lem:H1c}). Fixons une uniformisante $t$ de $X$ en $P_0$.
Soient $u_1=([D_1],D'_1,f_1)\in \HH^1(U,\mu_n)$ et $u_2=[D_2]\in \HH^1_c(U,\mu_n)$. Soit $f_2\in K$ telle que $nD_2=\div(f_2)$ et $f_2\equiv 1\mod Z$.  On peut supposer, quitte à ajouter à $D_1$ le diviseur d'une fonction $g$ et multiplier $f_1$ par $g^n$, que les supports de $\div(f_1)$ et $D_2$ sont disjoints. On suppose également, quitte à multiplier $f_1,f_2$ par des éléments de $K^\times$, que $(t^{-v_{P_0}(f_i)}f_i)(P_0)=1$ pour $i=1,2$.

\begin{lem} Avec ces notations, \[ \frac{f_1(D_2)}{f_2(D_1)}\in \mu_n(k).\]
\begin{proof}
\begin{align*}
\left(\frac{f_1(D_2)}{f_2(D_1)}\right)^n &= \frac{f_1(nD_2)}{f_2(nD_1)} & \\
&= \frac{f_1(\div f_2)}{f_2(\div f_1+D_1')}& \\
&= \frac{f_1(\div f_2)}{f_2(\div(f_1))f_2(D_1')}&\\
&= \frac{1}{f_2(D_1')}& \text{par  réciprocité de Weil}\\
&= 1 & \text{car }f_2\equiv 1\mod Z.
\end{align*}
\end{proof}
\end{lem}

\begin{df}
Avec ces notations, nous appellerons accouplement de Weil généralisé l'application $e_n\colon \HH^1(X,\mu_n)\times \HH^1_c(U,\mu_n)\to\mu_n$ définie par \[ e_n(u_1,u_2)\coloneqq \frac{f_1(D_2)}{f_2(D_1)}.\]
\end{df}

Décrivons explicitement comment ces flèches permettent de réaliser la suite de Gysin comme la duale de la suite exacte de cohomologie à support propre.
Considérons les isomorphismes \[\Lambda^\vee\xrightarrow{\tr^\vee} \HH^2(X,\mu_n)^\vee\xrightarrow{u} \HH^0(X,\Lambda)=\Lambda\] et $\HH^0(Z,\Lambda)^\vee\xrightarrow{v} \HH^0(Z,\Lambda)$ donnés respectivement par $u\colon \alpha\mapsto \alpha(1)$ et $v\colon \alpha\mapsto (\alpha(i_P))_{P\in Z}$, où $(i_P)_{P\in Z}$ désigne la base canonique de $\HH^0(Z,\Lambda)$. Ici, $\tr$ est l'isomorphisme $\HH^2(X,\mu_n)\to\Lambda$ du théorème \ref{th:dualpoinc}.
Le diagramme 
\[ \begin{tikzcd} 
0 \arrow[r]& \HH^2(X,\mu_n)^\vee(1) \arrow[r]\arrow[d,"u(1)"] &\HH^0(Z,\Lambda)^\vee(1) \arrow[d,"v(1)"] \arrow[r]& \HH^1(U,\mu_n)^\vee(1)\arrow[r]\arrow[d]& \HH^1(X,\mu_n)^\vee(1) \arrow[r]\arrow[d]& 0 \\
0 \arrow[r]& \HH^0(X,\mu_n) \arrow[r]& \HH^0(Z,\mu_n) \arrow[r]& \HH^1_c(U,\mu_n) \arrow[r] & \HH^1(X,\mu_n)\arrow[r] & 0
\end{tikzcd}
\]
où la ligne supérieure est la duale de la suite de Gysin décrite dans la section \ref{subsec:gysin}, la suite du bas est celle décrite dans la proposition \ref{prop:H1m}, et les flèches verticales sont celles décrites précédemment, est commutatif.
En particulier, l'orthogonal de $\HH^0(Z,\mu_n)/\HH^0(X,\mu_n)$ dans $\HH^1(U,\mu_n)$ est l'image de $\HH^1(X,\mu_n)$.

\subsection{Courbes nodales sur un corps algébriquement clos}

\begin{prop} Soit $X$ une courbe sur un corps algébriquement clos $k$, de normalisée $\tX$. Le morphisme canonique $\HH^2(X,\mu_n)\to \HH^2(\tX,\mu_n)$ déduit de $\tX\to X$ par fonctorialité est un isomorphisme. 
\begin{proof} Notons $S$ le lieu singulier de $X$ et $\nu\colon \tX\to X$ la normalisation de $X$. Comme $\nu$ est fini, le morphisme $\mu_n\to\nu_\star\nu^\star\mu_{n}$ est injectif ; comme $\nu$ est un isomorphisme en-dehors de $S$, le quotient $\mathcal{Q}$ est supporté sur $S$. La suite exacte \[ 0\to \mu_n \to \nu_\star\mu_n \to \mathcal{Q}\to 0 \]
donne la suite exacte longue \[ \HH^1(X,\mathcal{Q})\to \HH^2(X,\mu_n)\to \HH^2(X,\nu_\star \mu_n)\to \HH^2(X,\mathcal{Q}) \]
dont les extrémités sont nulles car $\mathcal{Q}$ est un faisceau gratte-ciel. Par conséquent, $\HH^2(X,\mu_n)$ est isomorphe à $\HH^2(X,\nu_\star\mu_n)$, lui-même isomorphe à $\HH^2(\tX,\mu_n)$ par exactitude de $\nu_\star$. 
\end{proof}
\end{prop}

\subsubsection{Courbes nodales projectives}

Soit $X_0$ une courbe projective nodale sur $k_0$, ayant un unique point nodal $P$. Rappelons que $n$ est premier à la caractéristique de $k_0$. Notons $X=X_0\times_{k_0}k$, et $\tX$ la normalisée de $X$ ; c'est le changement de base à $k$ de la normalisée de $X_0$. Soient $Q,R$ les points de $\tX$ au-dessus de $P$. Soit $\m$ le diviseur $Q+R$ sur $\tX$.

\begin{prop}\label{prop:H1mnod} Il y a une suite exacte courte de $\Lambda[\Gal(k|k_0)]$-modules \[ 0\to \frac{\mu_n(k)\times\mu_n(k)}{\mu_n(k)}\to  \HH^1(X,\mu_n)\to \HH^1(\tX,\mu_n)\to 0.\]
\begin{proof} Comme $P\in X(k_0)$, le diviseur $\m$ sur $X$ est $\mathfrak{G}_0$-invariant. La suite de Kummer fournit la suite exacte 
\[ 0\to \HH^1(X,\mu_n)\to \HH^1(X,\GG_m)\xrightarrow{n} \HH^1(X,\GG_m).\]
Le lemme \ref{lem:picnod} fournit alors un isomorphisme canonique \[ \HH^1(X,\mu_n)\xrightarrow{\sim}\Pic^0_\m(\tX)[n].\]
Le résultat est maintenant une conséquence directe de la proposition \ref{prop:H1m}.
\end{proof}
\end{prop}

\begin{cor} Les groupes de cohomologie de $X$ à valeurs dans $\mu_n$ sont les suivants.
\[
\HH^i(X,\mu_n)=\left\lbrace \begin{array}{lcl}
\mu_n(k) &\text{ si }& i=0\\
\Pic_{\m}^0(X)[n] &\text{ si }&i=1\\
\ZZ/n\ZZ &\text{ si }&i=2\\
0& \text{ si }&i\geqslant 3.
\end{array}\right.
\]
\end{cor}

\subsubsection{Courbes nodales affines}

Reprenons les notations précédentes, en supposant $X$ affine ; soit $\bar X$ la compactification lisse de la normalisation $\tX$ de $X$, et $Z$ le fermé réduit complémentaire de $\tX$ dans $\bar X$. Soit $Y$ la courbe construite à partir de $\bar X$ en identifiant $Q$ et $R$.
C'est une courbe projective qui contient $X$, est lisse en-dehors de $P$ et a pour normalisation $\bar X$. Elle se construit par exemple de la façon suivante. Soit $X\to\A^1\to \PP^1$ le morphisme donné par une coordonnée. Alors la normalisation de $\PP^1$ dans $X$ est une courbe projective contenant un ouvert isomorphe à $X$ \cite[03GT]{stacks}, et lisse en-dehors de cet ouvert. Rappelons que $\HH^1(X,\mu_n)$ est le groupe des classes de triplets $([D],D',f)\in \Pic^0(\bar X)\times \Div^0_Z(\bar X)\times k(X)^\times$ tels que $\div(f)=nD+D'$ et $f(Q)=f(R)\neq 0$. Le groupe $\HH^1(Y,\mu_n)$ est le sous-groupe de $\HH^1(X,\mu_n)$ constitué des triplets $(\bar D,D',f)$ tels que $D'=0$.

\begin{cor} Il y a une suite exacte de $\Lambda[\mathfrak{G}_0]$-modules \[ 0 \to \HH^1(Y,\mu_n) \xrightarrow{\phi} \HH^1(X,\mu_n) \xrightarrow{\psi} \Div^0_Z(\bar X)\otimes_\ZZ\Lambda  \to 0\]
dont les flèches sont décrites, via les isomorphismes précédents, par : \begin{itemize}[label=$\bullet$]
\item $\phi(D,f)=(D,0,f)$
\item $\psi(\bar D,D',f)=D'\mod n$
\end{itemize}
\begin{proof}
L'exactitude de la suite en $\HH^1(Y,\mu_n)$ et $\HH^1(X,\mu_n)$ est immédiate. Vérifions la surjectivité à droite. Soit $D'$ un zéro-cycle sur $Z$ vu comme diviseur sur $X$. Soit $\bar D\in \Div(\bar X)$ un diviseur tel que dans $\Pic(\bar X)$, $n[\bar D]=[D']$. Quitte à ajouter un diviseur principal à $\bar D$, on peut supposer que son support est disjoint de $\{ Q,R\}$. Il existe donc une fonction $f\in \OO_{\bar X,Q}^\times\cap\OO_{\bar X,R}^\times$ telle que $n\bar D=D'+\div(f)$. Alors $D'$ est l'image de $([\bar D],D',f)$. 
\end{proof}
\end{cor}

\subsection{Courbes lisses sur un corps fini}

Supposons ici que $k_0$ soit un corps fini à $q$ éléments. Soit $X_0$ une courbe intègre projective lisse sur $k_0$. Notons $f\colon X_0\to \Spec k_0$ le morphisme structural, et $K_0$ le corps des fonctions de $X_0$. Notons encore $X=X_0\times_{k_0}k$. Soit $\mathfrak{G}_0$ le groupe $\Gk$. 

\subsubsection{Cohomologie de $\GG_m$ et $\mu_n$}

La suite spectrale de Leray associée à l'isomorphisme canonique \[ \RG(X_0,\GG_m)=\RG(\Spec k_0,\R f_\star\GG_m)=\RG(\mathfrak{G}_0,\RG(X,\GG_m)) \] fournit, grâce à la proposition \ref{prop:H1Gm} et quelques calculs de cohomologie galoisienne, le résultat suivant. Une démonstration détaillée se trouve dans \cite{thuses}.

\begin{theorem} Les groupes de cohomologie de $X_0$ à valeurs dans $\GG_m$ sont les suivants.
\[
\HH^i(X_0,\GG_m)=\left\lbrace \begin{array}{lcl}
\Gamma(X_0,\OO_{X_0})^\times &\text{ si }& i=0\\
\Pic(X_0) &\text{ si }&i=1\\
\QQ/\ZZ &\text{ si }&i=3\\
0& \text{~~~~\, sinon}.
\end{array}\right.
\]
\end{theorem}

Comme dans le cas des corps algébriquement clos, la suite exacte de Kummer permet alors de calculer la cohomologie de $\mu_n$.

\begin{theorem} Les groupes de cohomologie de $X_0$ à valeurs dans $\mu_n$ sont les suivants.
\[
\HH^i(X_0,\mu_n)=\left\lbrace \begin{array}{lcl}
\mu_n(k_0) &\text{ si }& i=0\\
\Tors_{X_0}(\mu_n) &\text{ si }&i=1\\
\Pic(X_0)/n\Pic(X_0) &\text{ si }&i=2\\
\ZZ/n\ZZ& \text{ si }&i= 3\\
0& \text{ si }&i\geqslant 4.
\end{array}\right.
\]
\end{theorem}

Nous avons vu que \[ \HH^1(X_0,\mu_n)=\frac{ \{ ([D],f)\in \Pic(X_0)\times K_0^\times\mid nD=\div f\}}{\{ \{(D,f)\mid f\in (K_0^\times)^n\}}.\]
Rappelons que $\mathfrak{G}_0\simeq\hat \ZZ$ est de dimension cohomologique 1, et que $\HH^1(\mathfrak{G}_0,\mu_n(k))=k_0^\times/(k_0^\times)^n$ par le théorème de Hilbert 90. La suite spectrale de Hochschild-Serre associée à l'isomorphisme canonique de foncteurs dérivés \[\RG(X_0,\mu_n)=\RG(\mathfrak{G}_0,\RG(X,\mu_n))\] fournit alors une suite exacte courte de groupes abéliens
\[ 0\to k_0^\times/(k_0^\times)^n\to \HH^1(X_0,\mu_n)\to \HH^1(X,\mu_n)^{\mathfrak{G}_0}\to 0\]
dont la première flèche associe à $\alpha$ le couple $(0,\alpha)$. Le morphisme de droite est surjectif car $k_0$ est de dimension cohomologique 1.

\begin{rk}
Le $\Lambda$-module $\HH^1(X_0,\mu_n)$ est libre si et seulement si $k_0^\times/(k_0^\times)^n$ l'est, ce qui est le cas si $\pgcd(n,q-1)\in \{1,n\}$.
\end{rk}

\subsubsection{Cup-produit sur les corps finis}\label{sec:bleher}

Soit $\ell$ un nombre premier divisant $q-1$. L'étude de l'accouplement \[ \HH^1(X_0,\mu_\ell)\times \HH^1(X_0,\mu_\ell)\to \HH^2(X_0,\mu_\ell^{\otimes 2})=\Pic(X_0)\otimes \mu_\ell(k) \]
est réalisée dans \cite{chinburg_bleher}. Les résultats qui y sont obtenus ne permettent le calcul de tous les cup-produits que dans le cas des courbes de genre 1. De même que sur les corps algébriquement clos, il est nécessaire de fixer un point $P\in X_0(k)$ ; une fonction $f\in K_0^\times$ est dite normalisée si le coefficient dominant de sa série de Laurent en $P$ est une puissance $\ell$-ième dans le corps résiduel $k_0(P)$ \cite[§1]{chinburg_bleher}. Ceci est indépendant du choix d'une uniformisante en $P$. Toute fonction s'écrit alors comme produit d'une fonction normalisée avec un élément de $k_0^\times/(k_0^\times)^n$. \'{E}tant donné une fonction $f\in K_0^\times$ dont la classe du diviseur appartient à $\ell\Pic(X)$, notons $[f]$ son image dans $\HH^1(X_0,\mu_\ell)$. Le calcul du cup-produit $[a]\cup [b]$, où $a,b\in K_0^\times$, se ramène alors à deux situations : $a$ et $b$ sont toutes les deux normalisées, ou $a$ est normalisée et $b\in k_0^\times$.

\begin{theorem}\cite[Th. 1.1, Th 1.2]{chinburg_bleher}  Notons $\langle-,-\rangle$ l'accouplement de Weil sur $\HH^1(X,\mu_\ell)$. Si $X$ est de genre 1 alors pour toutes fonctions $a,b\in K_0$ normalisées et de diviseur dans $\ell\Pic(C)$, \[ [a]\cup [b]=\frac{1}{[k_0(P):k]}\left( [P]\otimes \left\langle \left[\frac{\div(a)}{\ell}\right],\left[\frac{\div(b)}{\ell}\right]\right\rangle\right)\in \Pic(X_0)\otimes \mu_\ell(k).\]
Si $X$ est de genre $\geqslant 2$, cette formule est encore valable si et seulement si l'image dans $\HH^2(X_0,\mu_\ell^{\otimes 2})$ par le cup-produit des classes de fonctions normalisées est un sous-espace vectoriel de dimension au plus 1. 
\end{theorem}
Dans la section \ref{sec:cupprod}, nous décrirons une façon de calculer explicitement les cup-produits dans la cohomologie des faisceaux lisses sur $X$.

\section{Cohomologie à support dans un fermé}\label{sec:cohfermcourb}

Soit $X$ une courbe intègre sur un corps algébriquement clos $k$. Soient $Z$ un sous-schéma fermé réduit strict de $X$, et $U$ l'ouvert complémentaire. Comme $Z$ est zéro-dimensionnel, $\HH^0_Z(X,-)=\bigoplus_{z\in Z} \HH^0_z(X,-)$, ce qui ramène le problème du calcul des $\HH^i_Z(X,-)$ au cas où $Z$ est un point fermé de $X$. Notons $i\colon z\to X$ l'inclusion de ce point, et $j\colon U\to X$ l'inclusion de l'ouvert $X-\{ z\}$.
Dans le cas où $X$ est lisse et $\F$ est un faisceau lisse sur $X$, le théorème de pureté affirme que $\HH^j_z(X,\F)=\HH^{j-2}(z,\F(-1))$. Ces groupes sont donc tous nuls, sauf lorsque $j=2$ ; le groupe $\HH^2_z(X,\F)$ est isomorphe par pureté à $\HH^0(z,\F(-1))$. Le résultat suivant concerne le prolongement par zéro à $X$ des faisceaux lisses sur $U$, sans hypothèse de régularité sur $X$. \\

\begin{lem}\label{lem:cohsuppz} Soit $\F$ un faisceau lisse sur $U$. \begin{itemize}[label=$\bullet$]
\item $\HH^0_z(X,j_!\F)=\HH^0_z(X,j_\star\F)=\HH^1_z(X,j_\star\F)=0$
\item $\HH^1_z(X,j_!\F)=\HH^0(z,\F)$
\item $\HH^2_z(X,j_!\F)=\HH^2_z(X,j_\star\F)$
\item Pour tout $i\geqslant 3$, $\HH^i_z(X,j_!\F)=\HH^i_z(X,j_\star\F)=0$.
\end{itemize}
\begin{proof} Comme $\HH^0(X,j_\star\F)\to \HH^0(U,\F)$ est un isomorphisme et $\HH^1(X,j_\star\F)\to \HH^1(U,\F)$ est un monomorphisme d'après le lemme 
\ref{lem:torsres}, la suite exacte de cohomologie à support pour $j_\star\F$ assure que $\HH^0_z(X,j_\star\F)=\HH^1_z(X,j_\star\F)=0$. 
De plus, pour tout $j\geqslant 3$, les groupes $\HH^{j-1}(U,j_\star\F)$ et $\HH^j(X,\F)$ sont nuls, donc $\HH^j_z(X,j_\star\F)=0$. 
Rappelons que $\HH^i_z(X,i_\star-)=\HH^i(Z,-)$. La suite exacte longue des $\HH^i_z(X,-)$ associée à la suite exacte courte \[ 0\to j_!\F\to j_\star\F\to i_\star i^\star j_\star\F\to 0\]
assure que $\HH^0_z(X,j_!\F)=0$, $\HH^1_z(X,j_!\F)=\HH^0(Z,\F)$ et $\HH^2_z(X,j_!\F)=\HH^2_z(X,j_\star\F)$. 
De même, le groupe $\HH^i_z(X,j_!\F)$ est nul dès que $i\geqslant 3$.\\
\end{proof}
\end{lem}

Le seul groupe restant à calculer est donc $\HH^2_z(X,j_!\F)=\HH^2_z(X,j_\star\F)$.

\subsection{Calcul de $\HH^2_z(X,j_!\F)$ lorsque $z$ est régulier}

Voici la situation : $X$ est une courbe intègre sur un corps algébriquement clos, $z$ est un point fermé régulier de $X$ de complémentaire l'ouvert $U$, et $\F$ est un faisceau lisse sur $U$. Notons $X_{\bar z}$ le spectre de l'hensélisé strict de l'anneau local de $X$ en $z$. Notons $j'\colon\eta\to X_{\bar z}$ l'inclusion du point générique.

\begin{lem}\label{lem:cohsuppferme} Il y a un isomorphisme canonique $\HH^2_z(X,j_!\F)\to \HH^1(\eta,\F_\eta)$.
\begin{proof}
Par excision, pour tout voisinage étale affine de $Y\to X$ tel que la préimage de $z$ soit réduite à un point $y$, le groupe $\HH^2_z(X,j_!\F)$ est canoniquement isomorphe à $\HH^2_y(Y,(j_!\F)|_Y)$ \cite[III, Cor. 1.28]{milneEC}. Comme la limite de ces voisinages est $X_{\bar z}$, on a $\HH^2_z(X,j_!\F)=\HH^2_z(X_{\bar z},j'_!\F_\eta)$. La suite exacte de cohomologie sur $X_{\bar z}$ à support sur $z$ montre alors que $\HH^2_z(X_{\bar z},j'_!\F_\eta)=\HH^1(\eta,\F_\eta)$ \cite[II, Prop. 1.1]{milneADT}. 
\end{proof}
\end{lem}

Notons $K=\Frac(\OO_{X,\bar z})$. Rappelons que le groupe d'inertie $I_z$ est aussi le groupe $\Gal(K^{\rm sep}|K)$. La cohomologie de $\eta$ est donc la cohomologie galoisienne de $I_z$, et $\HH^1(\eta,\F_\eta)=\HH^1(I_z,M)$, où $M$ désigne le $\pi_1(U,u)$-module $\F_\eta$. Rappelons que $\HH^1(I_z,M)$ est isomorphe à $M_{I_z}(-1)$ par la proposition \ref{prop:cohinertinf}. 

Il est maintenant possible de comprendre plus explicitement la suite exacte ouvert-fermé de la proposition \ref{prop:of} pour le faisceau $j_!\F$. C'est la suite \[ 0\to \HH^0(U,\F)\to \HH^1_Z(X,j_!\F)\to \HH^1(X,j_!\F)\to \HH^1(U,\F)\to \HH^2_Z(X,j_!\F)\to \HH^2(X,j_!\F)\to 0.\]
Soit $G$ le groupe d'automorphismes d'un revêtement galoisien $V\to U$ qui trivialise $\F$.
Le groupe $\HH^2(X,j_!\F)$ est isomorphe par dualité de Poincaré à $\HH^0(U,\F^\vee(1))^\vee=((M^\vee(1))^G)^\vee=M_G(-1)$. La suite se réécrit donc
\[ 0\to M^G \to \bigoplus_z M^{I_z}\to \HH^1(X,j_!\F)\to \HH^1(U,\F) \to \bigoplus_z M_{I_z}(-1)\to M_G(-1)\to 0.\]
Remarquons que l'inclusion $j_!\F\to j_\star\F$ fournit le diagramme commutatif
\[
\begin{tikzcd}
\HH^1(X,j_!\F) \arrow[r]\arrow[d] &\HH^1(U,j^\star j_!\F)\arrow[d] \arrow[r,"\sim"]& \HH^1(U,\F)\arrow[d,"\id"] \\
\HH^1(X,j_\star\F) \arrow[r] &  \HH^1(U,j^\star j_\star \F)\arrow[r,"\sim"] & \HH^1(U,\F)
\end{tikzcd}
\]
La flèche $\HH^1(X,j_!\F)\to \HH^1(U,\F)$ est donc la composée du morphisme $\HH^1(X,j_!\F)\to \HH^1(X,j_\star\F)$ déduit de $j_!\F\to j_\star\F$  avec le morphisme injectif $\HH^1(X,j_\star \F)\to \HH^1(U,j^\star j_\star \F)=\HH^1(U,\F)$ obtenu par restriction. En termes de torseurs, cette flèche associe à un $j_!\F$-torseur $\T$ sur $X$ la restriction à $U$ du $j_\star\F$-torseur $\T\wedge^{j_!\F}j_\star\F$. La première moitié de la suite ci-dessus s'identifie à la suite exacte issue de la suite exacte longue associée à \[ 0\to j_!\F\to j_\star \F\to i_\star i^\star j_\star\F\to 0.\]
La deuxième moitié de la suite s'écrit
\[ 0\to \HH^1(X,j_\star\F)\to \HH^1(U,\F)\to \bigoplus_z M_{I_z}(-1)\to M_G(-1)\to 0.\]
Nous verrons dans la section \ref{subsubsec:X2} comment construire un revêtement trivialisant $V\to U$ tel que $\HH^1(U,\F)=\HH^1(G,M)$, ce qui permet de comprendre le morphisme $\HH^1(U,\F)\to \HH^2_z(X,j_!\F)$ comme le morphisme de restriction $\HH^1(G,M)\to \HH^1(I_z/(I_z\cap \pi_1 V),M)$. 

\subsection{Calcul de $\HH^2_z(X,j_!\F)$ lorsque $z$ est nodal}

Considérons désormais la situation suivante : $X$ est une courbe intègre sur un corps algébriquement clos, $z$ est un point fermé nodal de $X$ de complémentaire l'ouvert $U$. Notons $\nu\colon\tX\to X$ la normalisation de $X$, et $x,y$ les points de $\tX$ au-dessus de $z$. Notons $X_{\bar z}$ le spectre de l'hensélisé strict de l'anneau local de $X$ en $z$. Il est constitué de trois points : un point fermé $z'$, et deux idéaux premiers minimaux $x',y'$ correspondant aux branches de $X$ en $z$. Notons $U_{\bar z}\coloneqq X_{\bar z}\times_X U=\{ x'\}\sqcup \{y'\}$ ; ce schéma est canoniquement isomorphe au coproduit des points génériques des hensélisés de $\tX$ aux points $x$ et $y$ (voir section \ref{sec:constrnod} pour les détails). La situation est résumée par le diagramme cartésien suivant.
\[
\begin{tikzcd}
U_{\bar z}\arrow[d,"g'"]\arrow[r,"j'"] & X_{\bar z}\arrow[d,"g"] &\arrow[l,"i'",swap] z' \arrow[d] \\
U \arrow[r,"j"] & X & \arrow[l,"i",swap] z
\end{tikzcd}
\]

\begin{lem} Soit $\F$ un faisceau sur $U_{\bar z}$. Pour tout entier naturel $q$, le groupe $\HH^q(X_{\bar z},j'_!\F)$ est nul.
\begin{proof} Cette preuve adapte celle de \cite[II, Prop. 1.1]{milneADT} au cas des courbes nodales. L'affirmation est vraie pour $q=0$ car $\HH^0(X_{\bar z},j'_!\F)$ est le noyau de $\HH^0(X_{\bar z},j'_\star\F)\to \HH^0(X_{\bar z},i'_\star i'^\star j'_\star\F)$, qui n'est rien d'autre que le morphisme identité de $\HH^0(U_{\bar z},\F)$.
Montrons d'abord que pour tout faisceau injectif $J$ sur $U_{\bar z}$, le faisceau $j'_!J$ sur $X$ est acyclique. Pour ce faire, commençons par prouver que la suite exacte courte \[ 0\to j'_!J\to j'_\star J\to i'_\star i'^\star j'_\star J\to 0\] est une résolution injective de $j'_!J$. Fixons des clôtures séparables de $k(x')$ et $k(y')$, et notons respectivement $I_{x'}$ et $I_{y'}$ les groupes de Galois associés. Le foncteur $i^\star j_\star$ s'identifie au foncteur \[ \begin{array}{rcl} \Mod_{I_{x'}}\times\Mod_{I_{y'}}&\to &\Ab \\ (M,N)&\mapsto& M^{I_{x'}}\times N^{I_{y'}}\end{array}\]
qui admet pour adjoint à gauche le foncteur associant à un groupe abélien $M$ le couple $(M,M)$ muni des actions triviales de $I_{x'}$ et $I_{y'}$. Cet adjoint à gauche étant exact, le foncteur $i'^\star j'_\star$ préserve les injectifs. Comme $i'_\star$ et $j'_\star$ préservent également les injectifs, la suite exacte ci-dessus est bien une résolution injective de $j'_!J$.
La suite exacte longue en cohomologie associée à cette suite exacte courte montre alors que $\HH^q(X_{\bar z},j'_\star J)=0$ pour tout entier $q\geqslant 1$. Soit $\F$ un faisceau sur $U_{\bar z}$, et $J^\bullet$ une résolution injective de $\F$. Alors $j'_!J^\bullet$ est une résolution acyclique de $j'_!\F$, et $\HH^q(X_{\bar z},j'_!\F)=\HH^q(\Gamma(X_{\bar z},j'_!J^\bullet))$. Ce dernier groupe est l'image de $\F$ par le $q$-ième foncteur dérivé de $\Gamma(X_{\bar z},j'_!-)$, qui est nul. 
\end{proof}
\end{lem}

\begin{lem}\label{lem:j!nod} Soit $\LL$ un faisceau lisse sur $U_{\bar z}$. Il y a un isomorphisme canonique 
\[ \HH^1(x',\LL_{x'})\times \HH^1(y',\LL_{y'})\xrightarrow{\sim}\HH^2_z(X,j_!\LL).\]
\begin{proof} Il y a toujours par excision un isomorphisme canonique \[\HH^2_z(X,j_!\LL)\xrightarrow{\sim}\HH^2_{z'}(X_{\bar z},g^\star j_!\LL)\xrightarrow{\sim}\HH^2_{z'}(X_{\bar z},j'_!g'^\star\LL).\]
La suite exacte de cohomologie sur $X_{\bar z}$ à support sur $z'$ pour le faisceau $j'_!g'^\star\LL$ s'écrit :
\[ \HH^1(X_{\bar z},j'_!g'^\star \LL)\to \HH^1(U_{\bar z},g'^\star\LL)\to \HH^2_{z'}(X_{\bar z},j'_!g'^\star\LL)\to \HH^2(X_{\bar z},j'_!g'^\star\LL).\]
Le lemme précédent assure alors que
\[ \HH^1(U_{\bar z},g'^\star\LL)\to \HH^2_{z'}(X_{\bar z},j'_!g'^\star\LL) \]
est un isomorphisme. Comme $U_{\bar z}$ est le coproduit des schémas $\{ x'\}$ et $\{ y'\}$, le groupe $\HH^1(U_{\bar z},g'^\star\LL)$ est simplement le produit de $\HH^1(x',\LL_{x'})$ et $\HH^1(y',\LL_{y'})$.
\end{proof}
\end{lem}

Soit $Z$ un fermé reduit zéro-dimensionnel de $X$, et $U$ l'ouvert complémentaire. Notons $\nu\colon\tX\to X$ la normalisation de $X$, et $\tilde{Z}$ la préimage de $Z$ dans $\tX$. Il y a un morphisme de foncteurs \[ \Gamma_{Z}(X,-)\to \Gamma_{\tilde Z}(\tX,\nu^\star -).\]

\begin{cor}\label{cor:cohsuppnodiso} Soit $\LL$ un faisceau lisse sur $U$. Le morphisme $\RG_Z(X,j_\star\LL)\to\RG_{\tilde Z}(\tX,\nu^\star j_\star\LL)$ est un quasi-isomorphisme.
\begin{proof} D'après le lemme \ref{lem:cohsuppz}, il suffit de montrer que $\HH^2_Z(X,j_\star\LL)\to\RG_{\tilde Z}(\tX,\nu^\star j_\star\LL)$ est un isomorphisme. Notons $Z_{\rm reg}$ (resp. $Z_{\rm sing}$) l'ensemble des points de $Z$ qui sont des points réguliers (resp. singuliers) de $X$. Soient $z_1,\dots,z_r$ les points de $Z_{\rm sing}$ ; pour chaque $i\in \{1\dots r\}$, notons $x_i,y_i$ les antécédents de $z_i$ dans $\tilde{Z}$. En un point $z$ de $Z_{\rm reg}$, le morphisme $\HH^2_z(X,j_\star\LL)\to \HH^2_z(\tX,\nu^\star j_\star\LL)$ est clairement un isomorphisme. Pour tout entier $i$, le morphisme \[\HH^2_{z_i}(X,j_\star\LL)\to \HH^2_{x_i}(\tX,j_\star\LL)\oplus \HH^2_{y_i}(\tX,j_\star\LL)\] n'est autre que le morphisme identité de \[ \HH^1(x_i',\LL_{x_i'})\oplus\HH^1(y_i',\LL_{y_i'})\to \HH^1(x_i',\LL_{x_i'})\oplus\HH^1(y_i',\LL_{y_i'})\]
où $x_i',y_i'$ désignent les points génériques respectifs des hensélisés stricts de $\tX$ en $x_i,y_i$. Le morphisme $\HH^2_Z(X,j_\star\LL)\to\HH^2_{\tilde Z}(\tX,\nu^\star j_\star\LL)$ est simplement la somme directe de tous ces isomorphismes.
\end{proof}
\end{cor}

\section{Revêtements cycliques de courbes}

\subsection{Revêtements cycliques de courbes lisses}

Fixons dans cette section un corps parfait $k_0$ et un entier $n$ inversible dans $k_0$. Notons $\Lambda=\ZZ/n\ZZ$. Supposons également que $k_0$ contient une racine primitive $n$-ième de l'unité. Soit $X_0$ une courbe intègre lisse sur $k_0$. Notons $K$ son corps des fonctions.

\begin{lem}\cite[Lem. 5.8.2]{szamuely} Soient $D\in\Div(X_0)$ et $f\in K$ tels que $\div(f)=nD\in\Div(X)$. Alors la normalisation $Y_0$ de $X_0$ dans $K(\sqrt[n]{f})$ est finie étale sur $X_0$. Si de plus $D$ est d'ordre $n$ dans $\Pic(X_0)$ alors $Y_0\to X_0$ est galoisien de groupe $\Lambda$.

\begin{proof} Soient $P$ un point fermé de $X_0$, et $t\in K$ une uniformisante en $P$. Si $v_P(f)=0$ alors la fibre $Y_{0,P}$ est isomorphe à $k[x]/(x^n-f(P))$, qui est une algèbre étale sur $k$ puisque $n$ est inversible dans $k$. Si $v_P(f)\neq 0$, comme $\div(f)=nD$, il existe $i\in\NN$ et $u\in \OO_{X_0,P}^\times$ tels que $f=ut^{ni}$. Alors $K(\sqrt[n]{f})=K(\sqrt[n]{u})$, ce qui nous ramène au cas précédent puisque $v_P(u)=0$. La finitude de $Y_0\to X_0$ est une propriété de la normalisation \cite[Prop. 12.44]{goertz_wedhorn}.\\
Supposons désormais que la classe de $\div(f)$ est d'ordre $n$ dans $\Pic(X_0)$. Notons $g=\sqrt[n]{f}\in L\coloneqq K(\sqrt[n]{f})$. Montrons que pour tout diviseur strict $d$ de $n$, la fonction $g^d\in L$ n'est pas un élément de $K$, ce qui suffit d'après \cite[Th. 6.2.(ii)]{lang} à prouver que $L/K$ est galoisienne de groupe $\Lambda$. Supposons qu'il existe un tel $d<n$ tel que $g^d\in K$. Alors le diviseur $D'\coloneqq \div(g^d)-dD\in \Div X_0$ vérifie $\frac{n}{d}D'=0\in\Div(X_0)$, ce qui implique que $D'=0$. Par conséquent, $dD=\div(g^d)$ est principal, ce qui est absurde puisque $D$ est d'ordre $n$ dans $\Pic(X_0)$.
\end{proof}
\end{lem}

\begin{cor} Rappelons que $k_0$ contient une racine primitive $n$-ième de l'unité. Un morphisme $Y_0\to X_0$ est un revêtement étale galoisien de groupe $\Lambda$ si et seulement s'il existe $f\in K^\times$ tel que $Y_0$ soit la normalisation de $X_0$ dans $K(\sqrt[n]{f})$.
\begin{proof} Soit $L$ le corps des fonctions de $Y_0$. Si $Y_0\to X_0$ est galoisien de groupe $\Lambda$, c'est encore le cas de $L/K$. La théorie de Kummer assure alors qu'un tel $f$ existe \cite[09DX]{stacks}. La réciproque est le lemme précédent.
\end{proof}
\end{cor}

\subsection{Revêtements de courbes nodales}

\subsubsection{Revêtements irréductibles}\label{subsubsec:revirred}

Soit $X$ une courbe nodale sur $k$, de points nodaux $P^{(1)},\dots,P^{(r)}$. Soit $\tX$ sa normalisée. Pour $i\in \{1\dots r\}$, notons $Q^{(i)},R^{(i)}$ les antécédents de $P^{(i)}$ dans $\tX$.\\

Soit $\tY$ un revêtement galoisien de $\tX$ de groupe $G$ d'ordre $s$. Au-dessus de $Q^{(i)}\in \tX$ se trouvent $s$ points $Q^{(i)}_1,\dots,Q^{(i)}_s\in \tY$. Pour $i\in \{1\dots r\}$, soit $R_1^{(i)}$ un point de $\tY$ au-dessus de $R^{(i)}$. Notons, pour $j\in \{1\dots s\}$, $\sigma_j^{(i)}$ le $\tX$-automorphisme de $\tY$ qui envoie $Q_1^{(i)}$ sur $Q_j^{(i)}$, et notons $R^{(i)}_j$ le point $\sigma_j^{(i)}(R_1^{(i)})$. Considérons la courbe $Y\coloneqq (\tY)_{Q_j^{(i)}=R_j^{(i)}}$ où $i$ parcourt $\{ 1\dots r\}$ et $j$ parcourt $\{1\dots s\}$. 
Elle possède $rs$ points singuliers $P_j^{(i)}$.

\begin{lem}\label{lem:revencoregal} La courbe $Y$ est un revêtement galoisien de $X$ de groupe $G$.
\begin{proof} L'étalitude du morphisme $Y\to X$ au point $P^{(i)}$ est montrée par le diagramme commutatif suivant, déduit de la remarque \ref{rk:conetg}.
\[
\begin{tikzcd}
(C_{Q_i}\tY\sqcup C_{R_i}\tY)_{Q_i=R_i} \arrow[r,"\sim"] \arrow[d,"\sim"] & C_{P_i}Y \arrow[d] \\
(C_Q\tX\sqcup C_R\tX)_{Q=R} \arrow[r,"\sim"]& C_{P}X
\end{tikzcd}
\]
Comme les points $Q^{(i)}_j, R^{(i)}_j$ ont été numérotés de façon compatible à l'action de $G$ sur $\tY$, tout $\tX$-automorphisme $\sigma_j$ de $\tY$ induit un $X$-automorphisme de $Y$, vérifiant $\sigma_j(P_1^{(i)})=P_j^{(i)}$. Par conséquent, $Y\to X$ est un revêtement galoisien de groupe $G$.
\end{proof}
\end{lem}

\begin{figure}[H]\centering \includegraphics[scale=0.5]{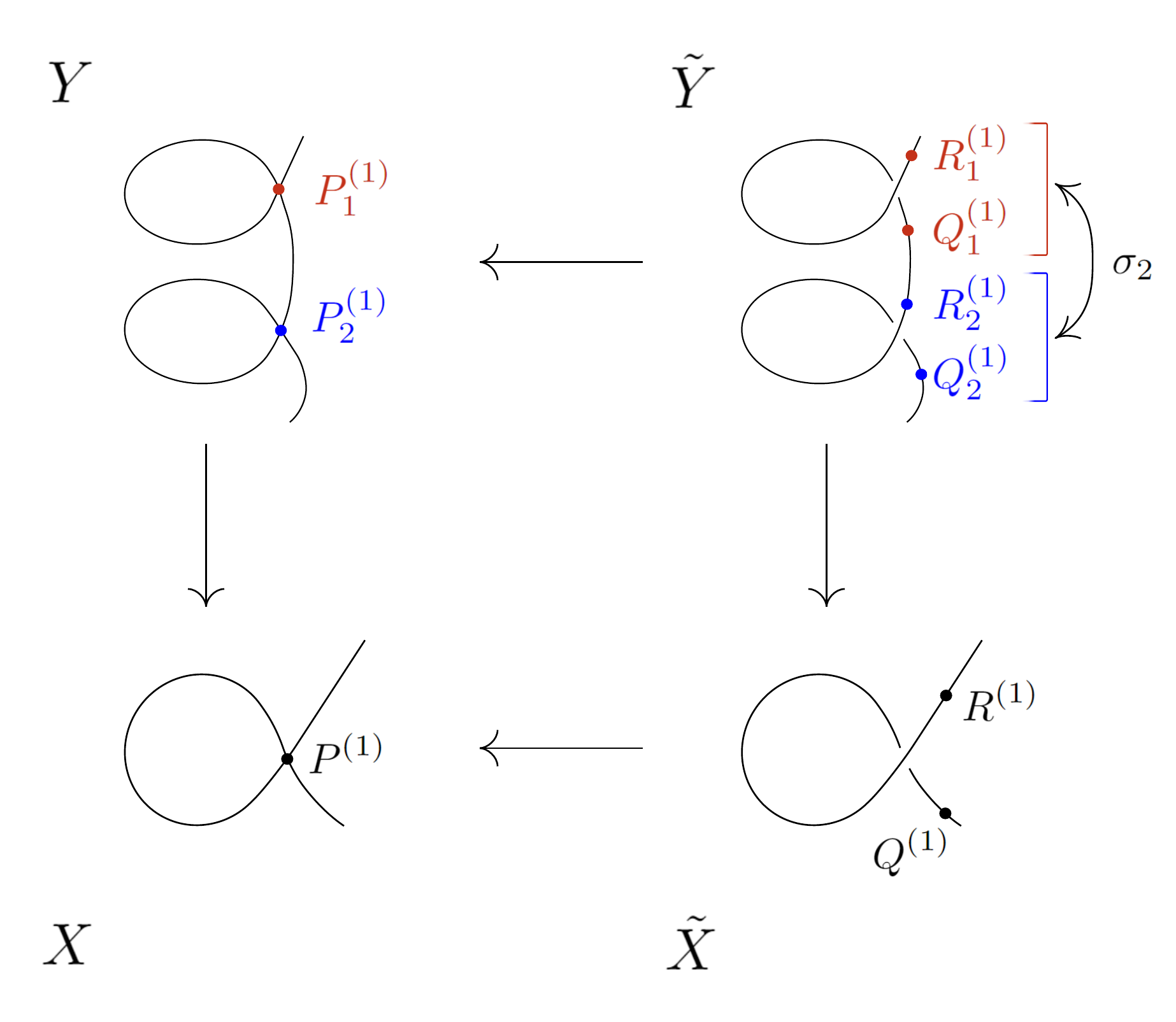}
\caption{Le revêtement $Y\to X$ lorsque $r=1$}
\end{figure}

\begin{cor} Soit $\C$ la catégorie des revêtements galoisiens de $X$ de groupe $G$. Soit $\tilde{C}$ la catégorie des revêtements galoisiens (connexes !) de $\tX$ de groupe $G$. Le foncteur $F\colon C\to\tilde{C}$ qui à un revêtement $W\to X$ associe $\tilde{W}\to\tX$, où $\tilde W$ est la normalisée de $W$, est une équivalence de catégories.
\begin{proof} Montrons que le foncteur $F'\colon \tilde{C}\to C$ construit au début de cette section \ref{subsubsec:revirred} est un quasi-inverse de $F$. Soit $W\to X$ un revêtement galoisien. Le morphisme $\tilde W\to W$ passe au quotient en $F'(\tilde W)\to W$, qui est un isomorphisme. D'autre part, soit $Y\to \tX$ un revêtement galoisien. Le morphisme $Y\to F'(Y)$ est fini, et  la propriété universelle de la normalisation assure qu'il se factorise par la normalisation de $F'(Y)$. Le morphisme $Y\to F(F'(Y))$ ainsi obtenu est un isomorphisme.
\end{proof}
\end{cor}

\subsubsection{Revêtements cycliques non irréductibles}\label{subsubsec:revnonirred}
Les courbes connexes lisses sur $k_0$ sont toutes irréductibles, ce qui n'est pas le cas des courbes connexes singulières. Voici comment construire des revêtements connexes cycliques non irréductibles de courbes nodales irréductibles. \\

Soit $X$ une courbe intègre sur $k_0$. Supposons $X$ nodale, de points singuliers $P^{(1)},\dots,P^{(r)}$. Soit $\tX$ sa normalisée. Notons, pour $i\in \{1\dots r\}$, $Q^{(i)}$ et $R^{(i)}$ les points de $\tX$ au-dessus de $P^{(i)}$. Nous supposerons dans la suite de la construction que les points $Q^{(i)},R^{(i)}$ sont tous définis sur $k_0$.
Dans $\tX\times\Lambda$, notons $Q^{(i)}_j,R^{(i)}_j$ les points de la $j$-ième composante $\tX$ au-dessus de $P^{(i)}$. Considérons la courbe $W_i$ obtenue à partir de $\tX\times \Lambda$ en identifiant, pour tout $j\in \Lambda$ et tout $m\neq i$,  $Q^{(i)}_j$ à $R^{(i)}_{j+1}$ et $Q^{(m)}_j$ à $R^{(m)}_j$. La courbe $W_i$ est connexe, de normalisée $\tX\times\Lambda$, et possède $n$ composantes irréductibles, images de celles de $\tX\times\Lambda$. Notons encore $P^{(i)}_1,\dots,P^{(i)}_n$ les antécédents de $P^{(i)}$ dans $W_i$. Dans l'illustration ci-dessous, les couleurs permettent de distinguer les deux composantes irréductibles de $W_1$.
\begin{figure}[H]\centering \includegraphics[scale=0.7]{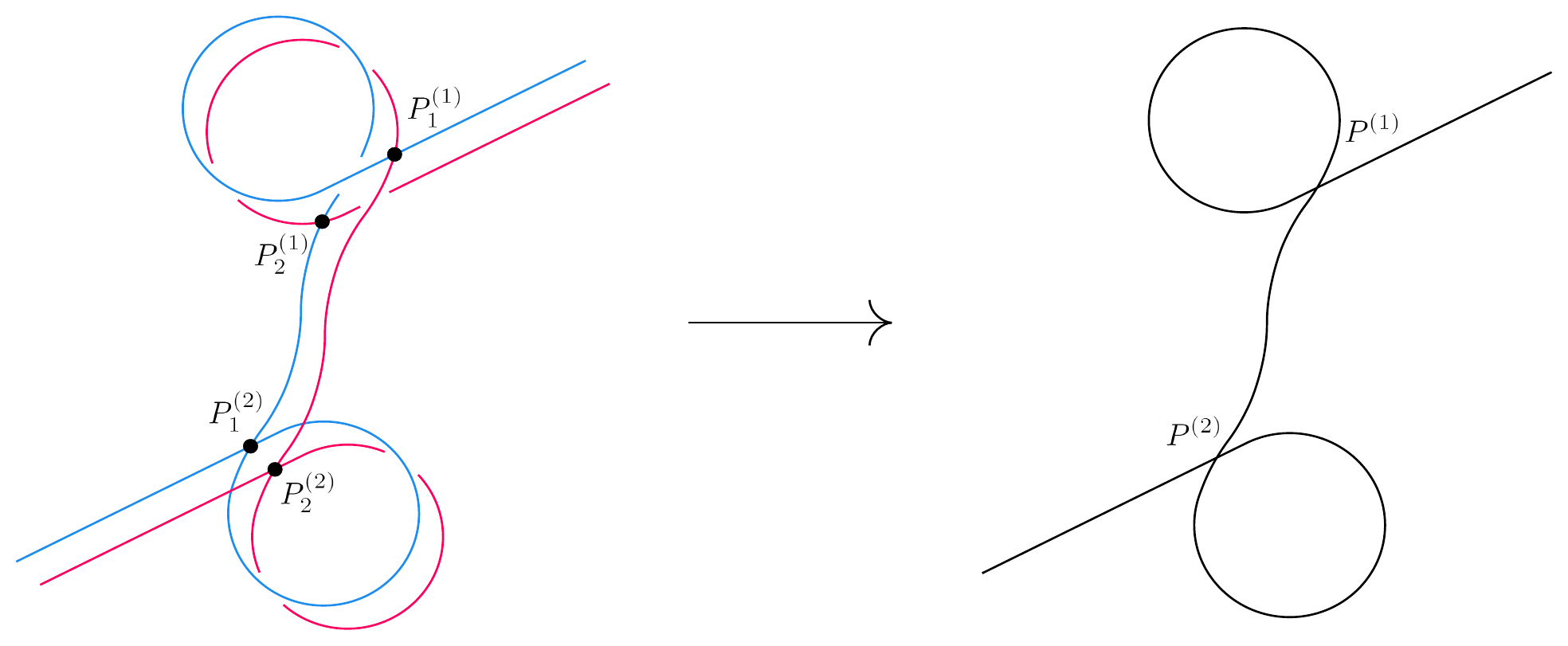}
\caption{Le revêtement $W_1\to X$ pour $n=2$}
\end{figure}
\begin{prop} Le morphisme $f_i\colon W_i\to X$ est un revêtement galoisien de groupe $\Lambda$.
\begin{proof} Le morphisme $f_i$ est clairement étale en dehors de $f_i^{-1}(P^{(1)},\dots,P^{(r)})$. Soit $m\neq i$. Au-dessus du point $P^{(m)}$ pour $m\neq i$, le revêtement $W_i\to X$ est localement coproduit de $n$ revêtements obtenus à partir de $\tX\to X$ en identifiant les deux points $Q^{(m)}$ et $R^{(m)}$ : il est donc étale en ce point. \'Etudions les antécédents de $P^{(i)}$. 
Pour tout $j\in\Lambda$, la remarque \ref{rk:conetg} fournit un diagramme commutatif 
\[
\begin{tikzcd}
(C_{Q^{(i)}_j}\tX\sqcup C_{R^{(i)}_{j+1}}\tX)_{Q^{(i)}_j=R^{(i)}_{j+1}} \arrow[r,"\sim"] \arrow[d,"\sim" {anchor=north, rotate=90}] & C_{P^{(i)}_j}W \arrow[d,"f_i"] \\
(C_{Q^{(i)}}\tX\sqcup C_{R^{(i)}}\tX)_{Q^{(i)}=R^{(i)}} \arrow[r,"\sim"]& C_{P^{(i)}}X
\end{tikzcd}
\]
qui montre que $f_i$ est étale en $P^{(i)}_j$. De plus, il est propre car $\tX\times\Lambda\to X$ l'est \cite[0AH6]{stacks}. Comme ses fibres sont finies, $f_i$ est donc fini \cite[02OG]{stacks}. Le degré de $f_i$ est clairement $n$. Exhibons un monomorphisme $\Lambda\to \Aut(W|X)$. Pour $i\in \Lambda$, le $X$-automorphisme $\sigma_i\colon \tX\times \Lambda\to \tX\times \Lambda$ qui envoie la $j$-ième copie de $\tX$ sur la $i+j$-ième induit un $X$-automorphisme de $W_i$. Par conséquent, $f_i$ est un revêtement galoisien de groupe $\Lambda$.
\end{proof}
\end{prop}

\begin{rk} Le revêtement $W_i\to X$ est en particulier un $\Lambda$-torseur (non trivial car connexe) sur $X$. Remarquons que $W_i\times_X \tX$ est la normalisation de $W$ \cite[07TD, 0CDV]{stacks}, c'est-à-dire $\tX\times\Lambda$. Le $\Lambda$-torseur non trivial $W_i\to X$ devient donc trivial après changement de base à $\tX$.
\end{rk}

Considérons désormais le revêtement $W\coloneqq W_1\times_X\dots\times_X W_r$ de $X$. 

\begin{lem} Le morphisme $f\colon W\to X$ est un revêtement galoisien de groupe $\Lambda^r$.
\begin{proof} Le morphisme $f$ est clairement fini étale de degré $n^r$. Montrons que $W$ est connexe. Comme $W\to W_1$ est lisse, la normalisée de $W$ est \cite[07TD]{stacks} : \[
\begin{array}{rcl}
 \tilde W &=& \tilde W_1 \times_{W_1} (W_1\times_X\dots\times_X W_r)
 \\ &=&(\tX\times \Lambda)\times_{W_1} W_1\times_X W_2\times\dots\times_X W_r\\
&=&(\tX\times\Lambda)\times_X W_2\times_X\dots\times_X W_r \\
&=& (\tilde W_2\times \Lambda)\times_X\dots \times_X W_r \\
&=& (\tX\times\Lambda^2)\times_X W_3\times_X\dots\times_X W_r\\
&\vdots & \\
&=& \tX\times\Lambda^r.
\end{array}\]
Le point $(Q^{(i)},j_1,\dots,j_r)\in \tX\times \Lambda^r$ a pour image \[ T^{(i)}_{j_1,\dots,j_r}\coloneqq (f_1(Q^{(i)}_{j_1}),\dots,f_r(Q^{(i)}_{j_r}))\in W_1\times_X\times_X\dots\times_X W_r.\] Les images dans $W$ des copies numéro $(j_1,\dots,j_i,\dots,j_r)$ et $(j_1,\dots, j_i+1,\dots,j_r)$ de $\tX$ dans $\tilde W$ ont pour point commun $T^{(i)}_{j_1,\dots,j_r}$, qui est l'image de $(Q^{(i)},j_1,\dots,j_r)$ et $(R^{(i)},j_1,\dots,j_i+1,\dots,j_r)$. Par conséquent, $W$ est connexe. Le morphisme $\prod_{i=1}^{r-1}\Aut(W_i|X)\to \Aut(W|X)$ est injectif car chaque $W_i\to X$ est surjectif, et par conséquent $W\to X$ est galoisien de groupe $\Lambda^r$.
\end{proof}
\end{lem}

\subsection{Les courbes sont des $K(\pi,1)$} \label{subsec:courbKpi1}

Rappelons que le corps $k$ est algébriquement clos.

\begin{prop}\label{prop:Kpi1} Soit $X$ une courbe intègre sur $k$. Si $X$ remplit l'une des conditions suivantes alors $X$ est un $K(\pi,1)$ et un $K(\pi,1)$ pro-$\ell$.\begin{enumerate}
\item $X$ est affine
\item $X$ est projective lisse de genre non nul
\item $X$ est projective nodale et de genre géométrique non nul
\end{enumerate}

\begin{proof} Soit $\F$ un faisceau lisse de $\Lambda=\ZZ/n\ZZ$-modules sur $X$. Notons $\rho$ le morphisme de topos $\Xet\to\Xfet$.
Comme $X$ est une courbe, $\HH^i(X,\F)$ est nul dès que $i> 2$. De plus, il est toujours vrai que $\id\to R\rho_\star\rho^\star$ est un isomorphisme en degrés $\leqslant 1$. La proposition est donc démontrée dans le cas où $X$ est affine. Supposons désormais qu'elle est projective, de genre géométrique non nul.
D'après la proposition \ref{prop:condKpi1}, il suffit donc de montrer qu'il existe un revêtement $\pi\colon Z\to X$ tel que le morphisme $\HH^2(X,\F|)\to \HH^2(Z,\F|_Z)$ soit nul. Notons $\nu\colon\tX\to X$ la normalisation de $X$. Soit $Y\to \tX$ un revêtement connexe trivialisant le faisceau lisse $\nu^\star\F$. Considérons le revêtement $Y_2$ de $Y$ de groupe $\HH^1(Y,\Lambda)^\vee$. Soit $Z$ le revêtement de $X$ obtenu comme décrit dans la section \ref{subsubsec:revirred} ; sa normalisation est $Y_2$. Le morphisme $\HH^2(Y,\F)\to \HH^2(Y_2,\F)$ est nul, car c'est la multiplication par $\deg(Z\to Y)$, qui est un multiple de $n$ puisque le genre de $Y$ est non nul. Le diagramme commutatif
\[
\begin{tikzcd} 
\HH^2(Z,\F) \arrow[r,"\sim"] & \HH^2(Y_2,\F) \\
\HH^2(X,\F) \arrow[r,"\sim"] \arrow[u] & \HH^2(\tX,\F) \arrow[u,"0",swap]
\end{tikzcd} 
 \]
conclut. 
La courbe $X$ est donc un $K(\pi,1)$. La même preuve convient pour les $K(\pi,1)$ pro-$\ell$ : un faisceau $\ell$-monodromique sur $\tX$ est trivialisé par un revêtement $Y\to \tX$ de groupe un $\ell$-groupe, et le revêtement $Z\to X$ obtenu à partir de $Y_2\to Y$ est encore un $\ell$-revêtement.
\end{proof}
\end{prop}

\begin{cor} Soit $X_0$ une courbe sur $k_0$. Notons $X=X_0\times_{k_0}k$. Si $X$ vérifie l'une des propriétés de la proposition précédente alors $X_0$ est un $K(\pi,1)$ et un $K(\pi,1)$ pro-$\ell$.
\begin{proof}
Soit $\F$ un faisceau lisse sur $X_0$. L'isomorphisme canonique $\mathfrak{G}_0$-équivariant \[ \RG(\pi_1(X),\F_\bareta)\xrightarrow{\sim}\RG(X,\F|_X) \] induit un isomorphisme canonique
\[ \RG(\pi_1(X_0),\F_{\bareta})=\RG(\mathfrak{G}_0,\RG(\pi_1(X),\F_\bareta))\xrightarrow{\sim}\RG(\mathfrak{G}_0,\RG(X,\F|_X))=\RG(X_0,\F).\]
\end{proof}
\end{cor}

\section{Un revêtement caractéristique}\label{sec:revcar}

\subsection{Sur un corps algébriquement clos}

\subsubsection{La construction} 
\label{subsubsec:X2}
\label{subsec:X2}

\begin{prop}\label{prop:sgcar} Soit $G$ un groupe topologiquement de type fini. Considérons le groupe abélien $\Lambda=\ZZ/n\ZZ$, muni de l'action triviale de $G$. Il existe un unique sous-groupe distingué $H$ de $G$ tel que $G/H$ soit isomorphe au $\Lambda$-dual $\HH^1(G,\Lambda)^\vee$ de $\HH^1(G,\Lambda)$. Ce sous-groupe est l'adhérence de $G^n[G,G]$ ; il est caractéristique dans $G$.
\begin{proof} Notons $S$ l'adhérence de $G^n[G,G]$. C'est un sous-groupe caractéristique de $G$ car tout automorphisme de $G$ préserve $G^n$ et $[G,G]$.
Comme $\Lambda$ est un groupe abélien de $n$-torsion, il y a un isomorphisme canonique \[ \Hom_{\cont}(G/S,\Lambda)\xrightarrow{\sim}\Hom_\cont(G,\Lambda)=\HH^1(G,\Lambda).\]
Notons que comme $G$ est topologiquement de type fini et $G/S$ est un $\Lambda$-module, $G/S$ a un ouvert dense fini. Comme un ouvert d'un groupe topologique compact est d'indice fini, $G/S$ est lui-même un $\Lambda$-module de type fini, et $\Hom_{\cont}(G/S,\Lambda)=\Hom(G/S,\Lambda)$.
Un $\Lambda$-module de type fini n'étant rien d'autre qu'une somme directe de $\ZZ/n_i\ZZ$ avec $n_i|n$, tout $\Lambda$-module de type fini est canoniquement isomorphe à son bidual. Par conséquent, l'isomorphisme ci-dessus devient par passage au dual \[\HH^1(G,\Lambda)^\vee \xrightarrow{\sim} G/S.\]
De plus, comme $\HH^1(G,\Lambda)^\vee$ est un groupe abélien de $n$-torsion, tout sous-groupe $H$ de $G$ tel qu'il y ait un isomorphisme $G/H\to \HH^1(G,\Lambda)^\vee$ contient $S$ ; il y est même égal puisque les quotients ont le même cardinal. 
\end{proof}
\end{prop}

Soit $X$ une courbe intègre lisse ou nodale sur $k$. Notons $K$ son corps des fonctions. Comme le groupe $\pi_1^{(p')}(X)$ est topologiquement de type fini, il existe un revêtement galoisien $X_2\to X$, modérément ramifié à l'infini, de groupe $\HH^1(\pi_1^{(p')}(X),\Lambda)^\vee=\HH^1(\pi_1(X),\Lambda)^\vee$. Décrivons comment construire ce revêtement $X_2\to X$. Nous avons vu dans la proposition \ref{prop:muntorseurslisses} que le groupe \[ \HH^1(\pi_1(X),\mu_n(k))\simeq \HH^1(X,\mu_n)\simeq\{ (D,g)\in \Div(X)\times K^\times \mid nD=\div g\} / \{(D,g)\mid g\in (K^\times)^n\} \] est un $\ZZ/n\ZZ$-module libre. Notons $r$ son rang. 
Soit $(D_i,g_i)_{1\leqslant i\leqslant r}$ une base de $\HH^1(X,\mu_n)$. Considérons, pour $i=0,\dots,r$, les extensions $L_i=K(\sqrt[n]{g_1},\dots,\sqrt[n]{g_{i}})$ de $K$ ; soit $\phi_i\colon Y_i\to X$ la normalisation de $X$ dans $L_i$. Notons $L=L_{r}$, et $\phi\colon Y\to X$ le revêtement correspondant de $X$.

\begin{lem} Le morphisme $\phi \colon Y\to X$ ainsi construit est isomorphe à $X_2\to X$.
\begin{proof} Vérifions par récurrence sur $j\in \{1,\dots,r\}$ que $Y_j\to Y_{j-1}$ est galoisien de groupe $\Lambda$. C'est vrai pour $j=1$.
Pour tout $i\in \{ 1,\dots,j-1\}$ l'extension $k(Y_i)=k(Y_{i-1})(\sqrt[n]{g_i})/k(Y_{i-1})$ est galoisienne de groupe $\Lambda$ par hypothèse de récurrence. La suite spectrale de Hochschild-Serre donne une suite exacte \[ 0\to \HH^1(\Lambda,\Lambda)\to \HH^1(Y_{i-1},\Lambda)\to \HH^1(Y_i,\Lambda). \] Par conséquent, le noyau de $\phi_i^\star\colon \HH^1(X,\Lambda)\to \HH^1(Y_i,\Lambda)$ est $\Lambda [D_1]\oplus\dots\oplus\Lambda [D_i]\simeq \Lambda^i$, et $[D_j]$ n'y appartient pas ; l'élément $\phi_i^\star [D_j]$ est encore d'ordre $n$ dans $\HH^1(Y_i,\mu_n)$. Ainsi, $Y\to X$ est fini étale d'ordre $n^{r}$. Le corps $k$ étant algébriquement clos, l'extension $L/K$ est le corps de décomposition des polynômes $T^n-g_1,\dots,T^n-g_r$, elle est donc galoisienne. Par la proposition \ref{prop:galcorps}, le morphisme $Y\to X$ est donc un revêtement galoisien. Un élément du groupe $\Aut(Y|X)$ est un automorphisme défini par $(\sqrt[n]{g_1}\mapsto \zeta_1\sqrt[n]{g_1},\dots,\sqrt[n]{g_{r}}\mapsto \zeta_{r}\sqrt[n]{g_{r}})$, où les $\zeta_i$ sont des racines $n$-ièmes de l'unité dans $k$ ; le groupe $\Aut(Y|X)$ est donc canoniquement isomorphe à $\Hom_\Lambda(\HH^1(X,\mu_n),\mu_n)=\HH^1(X,\Lambda)^\vee$. La proposition précédente assure que $Y$ est isomorphe à $X_2$.
\end{proof}
\end{lem}

\begin{rk}\label{rk:choixfct} Si $X$ est affine de compactification lisse $\bar X$, alors il est possible pour tout point $P\in Z\coloneqq\bar X-X$ de choisir les fonctions $g_1,\dots,g_r$ de façon à ce qu'elles soient toutes de valuation positive en $P$. En effet, les éléments d'une base $(D_i,g_i)$ de $\HH^1(\bar X,\mu_n)$ peuvent être choisis de façon à ce que le support de $D_i$ évite $Z$ (voir annexe \ref{subsec:moving}). Cette base peut être complétée en une base de $\HH^1(X,\mu_n)$ en ajoutant des antécédents de la base $(Q-Q_0)_{Q\in Z}$ de $\Div_Z^0(X)\otimes\Lambda$ pour un $Q_0\in Z-\{ P\}$ fixé.
\end{rk}

\begin{lem}\label{lem:revtriv} Pour tout $\Lambda$-module de type fini $F$, le revêtement $\phi \colon X_2\to X$ trivialise tous les $F$-torseurs sur $X$.
\begin{proof} Par construction, le revêtement $X_2\to X$ trivialise les $\Lambda$-torseurs sur $X$. Le $\Lambda$-module de type fini $F$ n'est qu'un produit de $\ZZ/n_i\ZZ$ où $n_i$ divise $n$, et $\HH^1(X,\ZZ/n_i\ZZ)\subset \HH^1(X,\Lambda)$. Par conséquent, le morphisme $\HH^1(X,\ZZ/n_i\ZZ)\to \HH^1(X_2,\ZZ/n_i\ZZ)$ est nul, et $\HH^1(X,F)\to \HH^1(X_2,F)$ est nul.
\end{proof}
\end{lem}

\begin{rk} Dans le cas où $X$ est une courbe projective lisse sur $k$ de jacobienne $J$, le choix d'un point $P\in X(k)$ détermine un plongement $i_P\colon X\to J$. Considérons le revêtement $Y\to X$ défini par le diagramme cartésien :
\[
\begin{tikzcd}
Y \arrow[d] \arrow[r] & J \arrow[d,"{[n]}"] \\
X \arrow[r,"i_P"] & J
\end{tikzcd}
\]
Le groupe d'automorphismes de l'isogénie $[n]$ est isomorphe à $J[n]=\HH^1(X,\mu_n)$. Le revêtement $Y\to X$, qui est de même degré que $[n]$, est encore galoisien de groupe d'automorphismes isomorphe à $\HH^1(X,\mu_n)$.
\end{rk}

\begin{rk}\label{rk:X2fin}
La même construction s'applique au cas de la cohomologie des courbes sur les corps finis. Supposons $k_0$ fini. Soit $X_0$ une courbe projective lisse géométriquement connexe sur $k_0$ de corps des fonctions $K_0$. Le morphisme $\HH^1(X_0,\mu_n)\to \HH^1(X,\mu_n)^{\mathfrak{G}_0}$ est surjectif ; en effet, comme $\mathfrak{G}_0$ est de dimension cohomologique 1, le terme suivant dans la suite spectrale de Hochschild-Serre est nul. Le lemme \ref{lem:PicDiv} assure alors que tout classe $\mathfrak{G}_0$-invariante de $\Pic^0(X)[n]$ contient un diviseur $\mathfrak{G}_0$-invariant $D$, et il existe une fonction $f\in K_0$ telle que $\div(f)=nD$. Un algorithme permettant de déterminer ces éléments $\mathfrak{G}_0$-invariants, et donc de construire une base de $\HH^1(X,\mu_n)^{\mathfrak{G}_0}$, est donné dans la proposition \ref{prop:divrat}. Il suffit ensuite de compléter cette base par la classe d'un générateur de $k_0^\times$, qui engendre le $\Lambda$-module $k_0^\times/(k_0^\times)^n$, pour obtenir une famille génératrice $(([D_0],f_0),\dots,([D_r],f_r))$ de $\HH^1(X_0,\mu_n)$. Par la même preuve que ci-dessus, la normalisation de $X_0$ dans $K_0(\sqrt[n]{f_0},\dots,\sqrt[n]{f_r})$ est alors un revêtement galoisien de $X_0$ de groupe $\HH^1(X_0,\Lambda)^\vee$.
\end{rk}

\subsubsection{Composition avec un autre revêtement}

Considérons maintenant le cas d'un revêtement galoisien de courbes intègres lisses $f\colon Y\to X$.  Soit $\bar f\colon\bar Y\to \bar X$ le morphisme fini entre les compactifications lisses induit par $f$. Alors $\bar f^{-1}(\bar X-X)=\bar Y-Y$, et en particulier, tout élément $\tau\in\Aut(Y|X)$ induit une permutation de $\bar Y-Y$. L'image de $(D_i,g_i)\in \HH^1(X,\mu_n(k))$ par $\tau^\star$ est donc simplement $(\tau^\star D_i,\tau^\star g_i)$. Considérons le revêtement $Y_2\to Y$ construit précédemment. 

\begin{lem}\label{lem:revcompgal} Le morphisme composé $Y_2\to Y\to X$ est encore un revêtement galoisien.
\begin{proof}  
Ceci découle directement du fait que $\pi_1(Y_2)$ est caractéristique dans $\pi_1(Y)$.
\end{proof}
\end{lem}

\begin{lem} Soit $\F$ un faisceau lisse de $\Lambda$-modules sur $X$ trivialisé par $Y$. Alors le revêtement $Y_2\to X$ trivialise tous les $\F$-torseurs.
\begin{proof} 
Notons $F=\HH^0(Y,\F)$. Le morphisme $\HH^1(X,\F)\to \HH^1(Y_2,F)$ se factorise par la flèche $\HH^1(Y,F)\to \HH^1(Y_2,F)$, qui est nulle par le lemme \ref{lem:revtriv}.
\end{proof}
\end{lem}

\subsubsection{Ramification}\label{subsubsec:ram}

Soit $X$ une courbe affine lisse sur $k$, de compactification lisse $\bar X$ et de genre géométrique $g_X$. 
Notons $P_0,\dots,P_r$ les points de $\bar X-X$, et $\overline{X_2}$ la compactification lisse du revêtement $X_2\to X$ de groupe $\HH^1(X,\Lambda)^\vee$. Le morphisme $\overline{X_2}\to \bar X$ est fini \cite[II, Prop. 6.8]{hartshorne}, et \[\overline{X_2}\times_{\bar X}(\bar X-X)=\overline{X_2}-X_2.\] \'{E}tudions les antécédents des $P_i$ dans $\overline{X_2}-X_2$ et leur ramification. Rappelons que $\HH^1(\bar X,\Lambda)^\vee$ est un quotient de $\HH^1(X,\Lambda)^\vee$ ; la compactification lisse $(\bar X)_2$ du revêtement de $X$ correspondant est étale au-dessus de $\bar X$. Il suffit donc d'étudier le revêtement $\overline{X_2}\to (\bar X)_2$. Une $\Lambda$-base de $\Div^0_{\bar X-X}(\bar X)\otimes\Lambda$ est donnée par $P_1-P_0,\dots,P_r-P_0$. Considérons des fonctions $g_1,\dots,g_r$ telles que \[\div(g_i)=nD_i+(P_i-P_0)\] où $D_i\in\Div^0(X-\{P_0,\dots,P_r\})$. Le revêtement $\overline{X_2}\to (\bar X)_2$ a pour corps de fonctions \[k((\bar X)_2)(\sqrt[n]{g_1},\dots,\sqrt[n]{g_r}).\] Soit $i\in \{1\dots r\}$. L'extension $k((\bar X)_2)(\sqrt[n]{g_j},j\neq i)$ de $k((\bar X)_2)$ fournit un revêtement $Y_i\to (\bar X)_2$ non ramifié au-dessus de $P_i$ puisque $v_{P_i}(g_j)=0$. Le revêtement $\overline{X_2}\to Y_i$ y est, quant à lui, ramifié ; soit $Q_i$ l'un des points de $Y_i$ au-dessus de $P_i$. Comme $v_{Q_i}(g_i)=1$, la fibre $(\overline{X_2})_{Q_i}$ est isomorphe à $k[x]/(x^n)$, et l'indice de ramification est $n$. En résumé, il y a au-dessus de $P_i$ exactement $\frac{1}{n}|\HH^1(X,\Lambda)|$ points de $\overline{X_2}$, tous d'indice de ramification $n$ au-dessus de $P_i$. Soit $R_i$ un antécédent de $Q_i$ dans $\overline{X_2}$.

Le sous-groupe d'inertie $I_{R_{i}|P_i}\subset \Aut(X_2|X)$ du point $R_i$ au-dessus de $\bar X$ s'insère dans la suite exacte
\[ 0\to I_{R_{i}|Q_i} \to I_{R_{i}|P_i} \to I_{Q_i|P_i} \to 0\]
(voir \cite[0BU7]{stacks} pour la surjectivité). Comme $I_{Q_i|P_i}=0$, il y a des isomorphismes \[I_{R_i|P_i}=I_{R_i|Q_i}=\Aut(\overline{X_2}|Y_i)\simeq\Lambda.\]
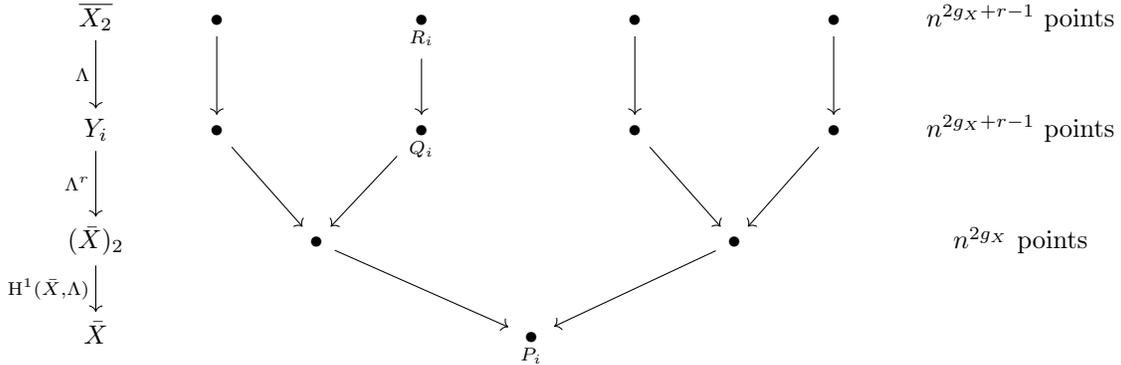
\begin{figure}[H]
\[
\begin{tikzcd}
\overline{X_2} \arrow[d,"\Lambda",swap] & \bullet \arrow[d]& &\underset{R_i}\bullet\arrow[d] & & \bullet\arrow[d] & & \bullet\arrow[d] & n^{2g_{X}+r-1}\text{ points} \\
Y_i\arrow[d,"\Lambda^r",swap] & \bullet\arrow[dr]& &\underset{Q_i}\bullet\arrow[dl] & & \bullet\arrow[dr] & & \bullet\arrow[dl] & n^{2g_{X}+r-1}\text{ points}\\
(\bar X)_2\arrow[d,"{\HH^1(\bar X,\Lambda)}",swap] & &  \bullet\arrow[drr] & & & & \bullet\arrow[dll] & &n^{2g_{X}}\text{ points}  \\
\bar X & & & & \underset{P_i}{\bullet}  & & & & 
\end{tikzcd}
\]
\caption{Ramification à l'infini de $X_2\to X$}
\end{figure}
En tant que sous-groupe de $\Aut(X_2|X)$, le groupe $I_{R_i|P_i}$ est celui engendré par $\sqrt[n]{g_i}\mapsto \zeta\sqrt[n]{g_i}$, où $\zeta$ est une racine $n$-ième primitive de l'unité. Tout ceci s'applique encore au point $P_0$, qui avait été choisi arbitrairement au début. Avec les notations ci-dessus, le sous-groupe $I_{R_0|P_0}$ est engendré par $(\sqrt[n]{g_1},\dots,\sqrt[n]{g_r})\mapsto (\zeta\sqrt[n]{g_1},\dots,\zeta\sqrt[n]{g_r})$.

\begin{lem} Le genre géométrique de $X_2$ est donné par \[ g_{X_2}=1+\frac{1}{2}n^{2g_{ X}-1+r}\left(n(2g_X-2)-1 \right).\]
\begin{proof}
L'application de la formule de Riemann-Hurwitz au revêtement ramifié de courbes projectives lisses $\overline{X_2}\to\bar X$ fournit l'égalité
\[ 2g_{{X_2}}-2=|\HH^1(X,\Lambda)|(2g_{\bar X}-2)+\sum_{P\in \overline{X_2}-X_2}(e_P-1)\]
où $e_P$ désigne l'indice de ramification de $X_2\to X$ en $P$. Rappelons que $|\HH^1(X,\Lambda)|=n^{2g_{X}+r}$. De plus, il y a au-dessus de chaque point $Q\in \bar X-X$ exactement $n^{2g_{X}-1+r}$ points, tous d'indice de ramification $n$. Par conséquent,
\[ 2g_{X_2}=2+n^{2g_{X}+r}(2g_{X}-2)+(n-1)n^{2g_{X}-1+r} \]
et le résultat s'en déduit immédiatement.
\end{proof}
\end{lem}

\subsection{Un revêtement semblable défini sur $k_0$}\label{subsec:X2gal}

Soit $Y_0$ une courbe connexe (mais pas nécessairement géométriquement connexe) lisse sur $k_0$. Notons $Y=Y_0\times_{k_0}k$. Le but de cette section est de déterminer un revêtement galoisien caractéristique $Y_{2,0}\to Y_0$ tel que $Y_{2,0}\times_{k_0}k\to Y$ trivialise tous les $\mu_n$-torseurs sur $Y$. Soient $X_0$ une courbe lisse géométriquement connexe sur $k_0$, et $f\colon Y_0\to X_0$ un revêtement galoisien. Notons encore $X=X_0\times_{k_0}k$.

\begin{rk} Considérons le cas particulier où $Y_0$ possède un $k_0$-point $y_0$ ; dans ce cas, $Y$ est connexe \cite[04KV]{stacks}. Soient $\bar y$ un point géométrique de $Y$ d'image $y_0$, et $\bar x$ son image par $f$. Soit $\bar y_2$ un point géométrique de $Y_2$ d'image $\bar y$. Le groupe $\pi_1(Y,\bar y_2)$ est caractéristique dans $\pi_1(Y,\bar y)$, qui est lui-même distingué dans $\pi_1(X,\bar x)$. Ainsi, l'action de $\mathfrak{G}_0$ par automorphismes sur $\pi_1(X,\bar x)$ induit encore une action sur $\pi_1(Y_2,\bar y_2)$, et par passage au quotient une action sur $\Aut(Y_2|X)$. 
\end{rk}

Cependant, la courbe $Y$ n'est pas nécessairement connexe : soient $Y^{(1)},\dots,Y^{(t)}$ ses composantes connexes.

\paragraph{Résumé de l'idée} Commençons par construire le revêtement $Y_2^{(1)}\to Y^{(1)}$ de groupe $\HH^1(Y^{(1)},\Lambda)$, puis le schéma $W=\coprod_{\sigma}(Y_2^{(1)})^\sigma$, où $\sigma$ parcourt les automorphismes d'une extension galoisienne suffisamment grande de $k_0$. Le revêtement $W\to Y$ provient d'un $k_0$-revêtement $Y_{2,0}\to Y_0$, qui vérifie la propriété recherchée.

\paragraph{Construction} Soit $k_1$ la clôture algébrique de $k_0$ dans $k_0(Y_0)$. C'est une extension séparable car $Y_0\to X_0$ est étale. Soit $k_2$ l'extension minimale de $k_0$ par laquelle se factorise l'action de $\mathfrak{G}_0$ sur $\HH^1(Y_1,\mu_n)$.  Soit $L$ la clôture galoisienne de la sous-extension de $k$ engendrée par $\mu_n(k),k_1,k_2$.  Soit $\alpha$ un élément primitif de l'extension séparable $L/k_1$, et $m\in k_1[t]$ son polynôme minimal. \'{E}crivons $k_0(Y_0)=k_0(x)[y]/(f)$.  Soient $(D_1,g_1),\dots,(D_r,g_r)$ des couples diviseur-fonction qui forment une base de $\HH^1(Y^{(1)},\mu_n)$. \'{E}crivons $g_i=g'_i(\alpha,x,y)$ avec $g'_i\in k_0(x)[t,y]/(m(t),f(x,y))$. Le morphisme $Y_{2,0}\to Y_0$ est alors défini par l'extension \[ k_0(Y_{2,0})=k_0(\alpha,x,y)(\sqrt[n]{g'_1},\dots,\sqrt[n]{g'_r}) \]
de $k_0(Y_0)$. La courbe $Y_{2,0}\times_{k_0}k$ est alors isomorphe à l'orbite sous $\Gal(L|k_0)$ de $Y_2^{(1)}$, et possède $[L:k_0]$ composantes connexes.

\begin{lem}\label{lem:degext} Le degré de $Y_{2,0}\to X_0$ est $n^r[L:k_1]\deg(Y_0\to X_0)$.
\begin{proof} Le degré de $Y_{2,0}\to (Y_0\times_{k_1}L)$ est $n^r$ et le degré de $Y_0\times_{k_1}L\to Y_0$ est $[L:k_1]$. 
\end{proof}
\end{lem}

\paragraph{Calcul de $\Aut(Y_{2,0}|X_0)$} Le morphisme $Y_{2,0}\to Y_0$ est de degré $n^r[L:k_1]$. Posons $z_i=\sqrt[n]{g'_i}$ avec les notations ci-dessus. Le groupe $\Aut(k_0(Y_{2,0})|k_0(Y_0))$ est constitué des automorphismes définis par $t\mapsto \sigma(t),z_i\mapsto \zeta_i z_i$ où $\sigma\in \Gal(L|k_1)$ et $\zeta_i\in\mu_n(L)$. Il y en a $\deg(Y_{2,0}\to Y_0)=n^{r}[L:k_1]$ car $\mu_n(k)\subset L$ : le revêtement $Y_{2,0}\to Y_0$ est donc galoisien.  
Calculons les automorphismes de $Y_{2,0}\to X_0$. Ce sont les \[ (t,x,y,z_1,\dots,z_r)\mapsto( \sigma(t),x',y',z'_1,\dots,z'_r) \]
où $\sigma\in \Gal(L|k_0)$, où $(x',y')$ est l'image de $(x,y)$ par un $X_0$-automorphisme de $Y_0$ induisant le même élément de $\Gal(k_1|k_0)$ que $\sigma$, et $z'_i\in k_0(Y_{2,0})$ vérifie ${z'_i}^n=\phi(g'_i)$. 
Il y a comme attendu $\deg(Y_{2,0}|X_0)=n^r[L:k_1]\deg(Y_0\to X_0)$ automorphismes de $Y_{2,0}\to X_0$, qui est donc un revêtement galoisien. 

\subsection{Adaptation aux courbes nodales}\label{subsec:X2nod}

Soit $X_0$ une courbe nodale sur $k_0$. Notons $X=X_0\times_{k_0}k$, et $\nu\colon\tX\to X$ sa normalisation. Soient $P_1,\dots,P_r$ les points nodaux de $X$. Notons $s=|H^1(\tX,\Lambda)|$.
Soit $\tX_2\to \tX$ le revêtement galoisien de $\tX$ de groupe $\HH^1(\tX,\Lambda)^\vee$. Soit $X'\to X$ le revêtement de $X$ de groupe $\HH^1(\tX,\Lambda)^\vee$ obtenu par la construction de la section \ref{subsubsec:revirred}. C'est une courbe qui a $rs$ points nodaux. Considérons également le revêtement $W\to X$ non irréductible de groupe $\Lambda^r$ construit dans la section \ref{subsubsec:revnonirred}. Considérons la courbe $Z\coloneqq W\times_X X'$. Le diagramme à carrés cartésiens ci-dessous, dont les flèches sont étiquetées par le degré des morphismes, résume la situation.
\[
\begin{tikzcd}
\tX_2 \arrow[r]\arrow[d,"s"] & X' \arrow[d,"s"] & \arrow[l,"n^r"] Z\arrow[d,"s"] \\
\tX \arrow[r] & X & \arrow[l,"n^r"] W
\end{tikzcd}
\]
\begin{lem}\label{lem:Zestconn} Le schéma $Z$ est connexe.
\begin{proof} Le morphisme $X'\to X$ étant lisse, le morphisme $Z\to W$ l'est encore ; d'après \cite[07TD]{stacks}, la normalisation de $Z$ est donc \[ \tilde Z=\tilde W\times_W Z=\tilde  W\times_W W\times_XX'=\tilde W\times_X X'=(\tX\times\Lambda^r)\times_X X'=\tX_2\times \Lambda^r.\]
Les points de $\tX_2$ au-dessus de $P^{(i)}$ sont $Q_1^{(i)},\dots,Q_s^{(i)},R_1^{(i)},\dots,R_s^{(i)}$. Soient $P^{(i)}_1,\dots,P^{(i)}_s$ leurs images respectives dans $X'$.
Notons $T^{(i)}_{j_1,\dots,j_r}$, où $(j_1,\dots,j_r)\in\Lambda^r$, les points de $W$ au-dessus de $P^{(i)}$. Les notations des antécédents de $P^{(i)}$ dans $\tX$, $\tX_2$, $X'$ et $W$ sont résumées dans le tableau ci-après, où $a\in \{1\dots s\}$ et $(j_1,\dots,j_r)\in \Lambda^r$.

\begin{table}[H]
\centering
{\renewcommand{\arraystretch}{2}
\begin{tabular}{|c|c|c|c|}
\hline
$\tX_2$               & $X'$        &      $Z=X'\times_X W$                                               & $\tilde Z=\tX_2\times\Lambda^r$                       \\ \hdashline
$Q^{(i)}_a,R^{(i)}_a$ & $P^{(i)}_a$ &  $(P^{(i)}_a,T^{(i)}_{j_1,\dots,j_r})$ & $(Q^{(i)}_a,j_1,\dots,j_r),(R^{(i)}_a,j_1,\dots,j_r)$ \\ \hline
$Q^{(i)},R^{(i)}$     & $P^{(i)}$   & $T^{(i)}_{j_1,\dots,j_r}$ & $(Q^{(i)},j_1,\dots,j_r),(R^{(i)},j_1,\dots,j_r)$     \\  \hdashline
$\tX$                 & $X$         & $W$                                                 & $\tilde W=\tX\times \Lambda^r$                        \\ \hline
\end{tabular}}
\end{table}

Souvenons-nous que le morphisme $\tilde W= \tX\times \Lambda^r\to W$ associe au couple $(Q^{(i)},j_1,\dots,j_r)$ le point $T^{(i)}_{j_1,\dots,j_r}$ et à $(R^{(i)},j_1,\dots,j_r)$ le point $T^{(i)}_{j_1,\dots,j_i-1,\dots,j_r}$.
Le morphisme \[ \tilde Z=\tX_2\times_X \Lambda^r\to Z=X'\times_X W \] associe aux points $(Q^{(i)}_a,j_1,\dots,j_r)$ et $(R^{(i)}_a,j_1,\dots,j_i+1,\dots,j_r)$ le couple $(P^{(i)}_a,T^{(i)}_{j_1,\dots,j_r})$. Les composantes irréductibles de $Z$ sont les images des $n$ composantes connexes de $\tilde Z$, toutes isomorphes à $\tX_2$ ; le point $(P^{(1)}_1,T^{(1)}_{j_1,\dots,j_r})$ appartient à l'image dans $Z$ de la $(j_1,\dots,j_r)$-ième et de la $(j_1,\dots,j_i+1,\dots,j_r)$-ième composante de $\tilde Z$. Ainsi, deux composantes irréductibles $C,C'$ de $\tilde Z$ sont toujours jointes par une suite \[ (C=C_0, C_1,\dots,C_m=C')\]
telle que pour tout $i$, l'intersection $C_i\cap C_{i+1}$ soit non vide.
\end{proof}
\end{lem}

\begin{prop} Le morphisme $Z\to X$ est un revêtement galoisien de groupe isomorphe à $\HH^1(X,\Lambda)^\vee$.
\begin{proof}
Il est fini étale de degré $n^rs$ car composée de morphismes finis étales de degrés respectifs $n^r$ et $s$, et connexe d'après le lemme précédent. Comme $Z=W\times_X X'$, il y a un morphisme $\Aut(W|X)\times\Aut(X'|X)=\Lambda^r\times \HH^1(\tX,\Lambda)^\vee\to \Aut(Z|X)$, qui est injectif car $W\to X$ et $X'\to X$ sont surjectifs. Le groupe de gauche étant d'ordre $\deg(Z\to X)$, ceci prouve que $Z\to X$ est galoisien. Rappelons que le groupe $\HH^1(X,\Lambda)$ est isomorphe à $\Lambda^r\times \HH^1(\tX,\Lambda)^\vee$, ce qui conclut.
\end{proof}
\end{prop}

\begin{cor} Le revêtement $Z\to X$ est caractéristique et trivialise tous les $\Lambda$-torseurs sur $X$.
\begin{proof} Le revêtement est caractéristique car son groupe est isomorphe à $\HH^1(X,\Lambda)^\vee$ (voir proposition \ref{prop:sgcar}). Notons $G=\Aut(Z|X)$. La suite spectrale de Hochschild-Serre pour $Z\to X$ donne une suite exacte \[ 0\to \HH^1(G,\Lambda)\to \HH^1(X,\Lambda) \to \HH^1(Z,\Lambda).\]
Sachant que le groupe $G\simeq \HH^1(X,\Lambda)$ est un $\Lambda$-module libre et que l'action de $G$ sur $\Lambda$ est triviale, $\HH^1(G,\Lambda)=\Hom_\Lambda(G,\Lambda)=\Lambda^{\rg_\Lambda G}=\HH^1(X,\Lambda)$. Par conséquent, le morphisme $\HH^1(X,\Lambda)\to \HH^1(Z,\Lambda)$ est nul. 
\end{proof}
\end{cor}

\paragraph{Construction du revêtement défini sur $k_0$} Soit $k'/k_0$ l'extension minimale de $k_0$ sur laquelle sont définis les antécédents dans $\tX$ des points singuliers de $X$. Soit $X_1$ la normalisation de $X_0$ dans $k_1(X_0)$. Considérons le revêtement galoisien $\tilde{Z}_1\to \tX_1$ de la section \ref{subsec:X2gal}. Construisons comme dans la section \ref{subsubsec:revirred} le revêtement galoisien $Z_1$ de $X_1$ correspondant. De même, construisons le revêtement non irréductible $W$ de $X_1$ de groupe $\Lambda^r$ défini dans la section \ref{subsubsec:revnonirred}. Posons enfin $X_{2,0}=Z_1\times_{X_1}W$. La connexité de $X_{2,0}$ se montre comme dans le lemme \ref{lem:Zestconn}. Si $X_0\to Y_0$ est un revêtement galoisien, le même argument que précédemment montre que $\Aut(X_{2,0}|Y_0)$ est isomorphe à $\Lambda^r\times \Aut(Z_1|Y_0)$ ; or $Z_1\to Y_0$ est galoisien car $\tilde{Z}_1\to\tilde{Y}_0$ l'est, donc $X_{2,0}\to Y_0$ l'est encore.

\subsection{Un exemple détaillé}\label{subsec:exdetrev}

Prenons $n=2$. Supposons que $-1$ n'est pas un carré dans $k_0$. Notons $V=\PP^1-\{ 0,\pm 1,\infty \}$ et $U=\PP^1-\{ 0,1,\infty\}$. Considérons le revêtement étale de degré 2
\[f\colon \begin{array}[t]{rcl} V&\longrightarrow& U \\ y&\longmapsto& y^2\end{array}\]
de groupe d'automorphismes engendré par $\tau\colon y\mapsto -y$. Le faisceau $\F\coloneqq f_\star\Lambda$ est un faisceau lisse sur $U$, trivialisé par le revêtement $f\colon V\to U$ puisque $f^\star f_\star \Lambda\simeq \Lambda^2$. Il correspond au $\Aut(V|U)$-module $\Lambda^2$, où l'élément non trivial de $\Aut(V|U)$ intervertit les deux copies de $\Lambda$.  \\

\paragraph{Calcul de $V_2$} Le groupe $\HH^1(V,\mu_2)\simeq\Lambda^3$ est engendré par les couples diviseur-fonction $(0-\infty,x),(1-\infty,x-1),(-1-\infty,x+1)$. Le revêtement $V_2\to V$ de groupe $\HH^1(V,\Lambda)^\vee$ correspond à l'extension de corps $k(\sqrt{x},\sqrt{x-1},\sqrt{x+1})/k(x)$. Le revêtement de $\bar V=\PP^1$ correspondant est le morphisme
\[ \begin{array}{ccc} \Proj k[y,z,t,h]/(z^2-(y^2-h^2),t^2-(y^2+h^2)) &\longrightarrow& \Proj k[y,h] \\ (y:z:t:h)&\longmapsto& (y^2:h^2)\end{array}\]
ramifié au-dessus de $0,\pm 1,\infty$. 
\paragraph{Calcul de $\Aut(V_2|U)$} Le groupe d'automorphismes $G\coloneqq \Aut(V_2|U)$ est d'ordre 16 ; il suffit pour le déterminer entièrement d'y trouver un antécédent du générateur $\tau$ de $\Aut(V|U)$. Un tel antécédent est $\gamma \colon (y:z:t:h)\mapsto (\sqrt{-1}y:\sqrt{-1}t:\sqrt{-1}z:h)$. Notons $\sigma_1\colon y\mapsto -y$, $\sigma_2\colon z\mapsto -z$, $\sigma_3\colon t\mapsto -t$ les générateurs évidents de $\Aut(V_2|V)\triangleleft G$. Alors $\gamma\sigma_2=\sigma_3\gamma$ et $\gamma\sigma_3=\sigma_3\gamma$, ce qui implique que $\langle \sigma_2,\sigma_3\rangle$ est distingué dans $G$. On vérifie aisément que la composée \[\langle \gamma \rangle \to G\to G/\langle \sigma_2,\sigma_3\rangle\] est un isomorphisme ; par conséquent, \[ G=\langle \sigma_2,\sigma_3 \rangle \rtimes \langle \gamma\rangle.\]

\paragraph{Ramification} Notons $Z=\bar U-U$, $W=\bar V-V$ et $W'=\overline V_2-V_2$. Le tableau ci-dessous résume la situation.

\begin{table}[H]
\resizebox{\textwidth}{!}
{\begin{tabular}{|l|clll|clllclll|clll|}
\hline
Points de $Z$                                                                                  & \multicolumn{4}{c|}{$0$}                                                                                  & \multicolumn{8}{c|}{$1$}                                                                                                                                                                                                     & \multicolumn{4}{c|}{$\infty$}                                                                       \\ \hline
\begin{tabular}[c]{@{}l@{}}Antécédents dans $W$ \\ Ramification \end{tabular}  & \multicolumn{4}{c|}{\begin{tabular}[c]{@{}c@{}}$0$\\ indice $2$\end{tabular}}                                    & \multicolumn{4}{c|}{\begin{tabular}[c]{@{}c@{}}$-1$\\ indice $1$\end{tabular}}                                           & \multicolumn{4}{c|}{\begin{tabular}[c]{@{}c@{}}$1$\\ indice $1$\end{tabular}}                                   & \multicolumn{4}{c|}{\begin{tabular}[c]{@{}c@{}}$\infty$\\ indice $2$\end{tabular}}                         \\ \hline
\begin{tabular}[c]{@{}l@{}}Antécédents dans $W'$ \\ Ramification \end{tabular} & \multicolumn{4}{c|}{\begin{tabular}[c]{@{}c@{}}4 points \\indice 4\end{tabular}} & \multicolumn{4}{c|}{\begin{tabular}[c]{@{}c@{}}4 points \\ indice 2 \end{tabular}} & \multicolumn{4}{c|}{\begin{tabular}[c]{@{}c@{}}4 points \\indice 2 \end{tabular}} & \multicolumn{4}{c|}{\begin{tabular}[c]{@{}c@{}}4 points \\ indice 4 \end{tabular}} \\ \hline
\begin{tabular}[c]{@{}l@{}}Un antécédent dans $W'$ \\ Son groupe d'inertie\end{tabular} & \multicolumn{4}{c|}{\begin{tabular}[c]{@{}c@{}}$P_0=(0,\sqrt{-1},1)$ \\ $\langle \gamma\sigma_2\rangle\simeq \mu_4(k)$\end{tabular}} & \multicolumn{4}{c|}{\begin{tabular}[c]{@{}c@{}}$P_{-1}=(\sqrt{-1},\sqrt{-2},0)$
 \\ $\langle \sigma_3\rangle\simeq \mu_2(k)$ \end{tabular}} & \multicolumn{4}{c|}{\begin{tabular}[c]{@{}c@{}}$P_1=(1,0,\sqrt{2})$
 \\ $\langle \sigma_2\rangle\simeq \mu_2(k)$ \end{tabular}} & \multicolumn{4}{c|}{\begin{tabular}[c]{@{}c@{}}$P_\infty=(1:1:1:0)$
 \\ $\langle \gamma\rangle\simeq \mu_4(k)$ \end{tabular}} \\ \hline
\end{tabular}}
\end{table}

L'isomorphisme canonique $I_{P_0}\to \mu_4(k)$ est obtenu explicitement de la façon suivante. Une uniformisante de $\bar V_2$ en $P_0=(0,\sqrt{-1},1)$ est $y$. L'orbite de $y$ sous l'action de $I_{P_0}=\langle \gamma\sigma_2\rangle$ est $\{ \pm y,\pm \sqrt{-1}y\}$. L'ensemble des $\frac{\sigma(y)}{y}(P_0)$ où $\sigma$ parcourt $I_{P_0}$ est donc exactement $\mu_4(k)$. \`{A} un élément $\sigma\in I_{P_0}$, l'isomorphisme $I_{P_0}\to \mu_4(k)$ associe $\frac{\sigma(y)}{y}(P_0)$.

Le générateur $\sqrt{-1}$ de $\mu_4(k)$ échange les deux copies de $\Lambda$ dans $M=\Lambda^2$. Le $\Lambda$-module des morphismes croisés $\mu_4(k)\to M$ est isomorphe à $\Lambda^2$, et $\tau_{\leqslant 1}\RG(I_{P_0},M)$ est représenté par le complexe suivant. 

\[ \begin{array}{rcl} \Lambda^2 &\to& \Lambda^2\\
(a,b)&\mapsto& \left[\sqrt{-1}\mapsto (a+b,a+b)\right] \end{array} \] 
Le groupe $\F_0=\HH^0(I_{P_0},M)$ est engendré par $(1,1)$, et $\HH^2_0(X,j_\star\F)=\HH^1(I_{P_0},M)$ est engendré par la classe de $(0,1)$. Le calcul de $\tau_{\leqslant 1}\RG(I_{P_\infty},M)$ est très semblable. Le groupe $I_{P_1}$ est, quant à lui, canoniquement isomorphe à $\mu_2(k)$, et agit trivialement sur $M$. Par conséquent, $\tau_{\leqslant 1}\RG(I_{P_1},M)$ est représenté par le complexe suivant.
\[ \begin{array}{rcl}\Lambda^2 &\to& \Lambda^2 \\ (a,b)&\mapsto &0 \end{array} \]
Nous calculerons $\RG(U,\F)$ dans la section \ref{subsec:exdetlis}.

\cleartooddpage

\chapter{Algorithmique des faisceaux constructibles}\label{chap:3}

Fixons un corps $k_0$, et une clôture algébrique $k$ de $k_0$. Soit $n$ un entier naturel non nul. Notons $\Lambda$ l'anneau $\ZZ/n\ZZ$. Dans toute cette section, les faisceaux constructibles considérés seront des faisceaux de $\Lambda$-modules.\\

L'objectif de ce chapitre est de donner diverses représentations et opérations sur les faisceaux lisses sur les $k_0$-schémas de type fini, puis des faisceaux constructibles sur les courbes lisses sur $k$. Nous donnons dans le cas des courbes lisses des algorithmes permettant de passer d'une représentation à une autre, ainsi que des algorithmes permettant d'effectuer des opérations (images directes et réciproques, noyaux et conoyaux de morphismes, $\Hom$ interne et produit tensoriel...) sur ces faisceaux. Nous montrons également comment effectuer ces opérations dans le cas général des faisceaux constructibles après avoir décrit comment calculer les poussés en avant de faisceaux lisses par des morphismes entre variétés régulières de même dimension.
Nous nous assurons que tous les algorithmes présentés sont de complexité élémentaire (voir annexe \ref{sec:compl}) en les entrées. Les schémas sont décrits comme recollement de schémas affines (voir annexe \ref{sec:repsch}). Pour cette représentation, il existe des algorithmes de complexité élémentaire calculant la normalisation ou la décomposition primaire d'une variété. Les courbes projectives lisses sont décrites par des produits de corps de fonctions (voir annexe \ref{subsec:modplan}), et les courbes affines lisses comme un ouvert d'un modèle plan de leur compactification lisse.

\section{Faisceaux lisses}

\subsection{Représentations des faisceaux lisses}\label{subsec:replisses}

Soit $X$ un schéma intègre de type fini sur $k_0$. Soit $\F$ un faisceau lisse sur $X$, correspondant à un $\pi_1(X)$-module $M$. Nous nous intéresserons à deux façons de définir explicitement $\F$ : \begin{enumerate}
\item par un $X$-schéma en groupes fini étale $F$ qui le représente ;
\item par un revêtement galoisien $Y\to X$ qui le trivialise ainsi que le $\Aut(Y|X)$-module $M$.
\end{enumerate}

\paragraph{Passage de la première à la deuxième représentation} Supposons $\F$ défini par un morphisme fini étale $T\to X$ de degré $d$, ainsi qu'une application $T\times_X T\to T$ définissant sa loi de groupe. Le calcul d'un revêtement trivialisant, décrit par exemple dans \cite[Prop. 5.8.1.(i)]{fulei}, se fait de la façon suivante : trouver une composante connexe $T'$ de $T$ telle que $T'\to X$ soit de degré >1, et changer de base à $T'$. Recommencer cette opération avec $T'\times_X T\to T'$ (toujours de degré $d$), jusqu'à obtenir un schéma $Y$ avec $d$ composantes connexes. \\

Remarquons qu'à chaque étape, le nombre de composantes connexes de $T$, et donc le nombre d'éléments de $\F(T)$, augmente. Il y a dans cet algorithme au plus $d-1$ appels récursifs ; comme $\F$ est un faisceau de groupes abéliens, c'est même $\log_2(d)$ puisqu'à chaque étape, $\F(T)$ est un sous-groupe strict de $\F(T')$. Chacune de ces étapes consiste en une décomposition primaire, puis le calcul d'un produit fibré. \\

Une fois cette opération effectuée, il reste à calculer la clôture galoisienne $Z\to X$ de $Y\to X$, qui est une composante connexe de $Y\times_X \dots\times_X Y$. Celle-ci se calcule d'une façon semblable à la clôture galoisienne d'une extension de corps \cite[§2]{huang_closure}. L'action de $\sigma\in\Aut(Z|X)$ sur $\F(Z)$ est donnée par la permutation des composantes connexes de $T\times_X Z\simeq \sqcup^d Z$ induite par \[T\times_X Z\xrightarrow{\id \times \sigma} T\times_X Z. \]
Rappelons que le degré de $Z\to X$ est majoré par $\deg(Y\to X)!$. Le cardinal de $\F(Z)$ est, quant à lui, égal au degré de $T\to X$.

\paragraph{Passage de la deuxième à la première représentation} Supposons $\F$ défini par un revêtement galoisien $f\colon Y\to X$ de groupe $G$, et le $G$-module $M=\HH^0(Y,f^\star\F)$. Rappelons que $\F=(f_\star f^\star\F)^G$. Le faisceau $f_\star f^\star\F$ est représenté par la restriction de Weil  $R\coloneqq R_{Y\to X}(M\times Y)$, et est encore muni d'une action de $G$ qui permute ses composantes connexes. Le faisceau $(f_\star f^\star \F)^G$ est représenté par l'intersection schématique $\bigcap_{g\in G}\ker(g-\id_R)$.

\paragraph{Calcul du revêtement trivialisant minimal} Une fois calculé un revêtement galoisien $Z\to X$ qui trivialise $\F$, un revêtement minimal est donné par $Z/H$, où $H$ est le noyau de $\Aut(Z|X)\to \Aut(\HH^0(Z,\F))$. Le degré de $Z/H\to X$ est le cardinal du groupe de monodromie, image de $\Aut(Z|X)$ dans $\Aut(\HH^0(Z,\F))$.

\paragraph{Simplifications dans le cas des courbes intègres lisses} Soit $X$ une courbe intègre lisse sur $k_0$. Un revêtement $Y$ de $X$ est simplement défini par l'extension $L/K$ de corps de fonctions correspondante. Le groupe $\Aut(Y|X)$ est $G=\Aut(L|K)$, et si le revêtement est galoisien, le faisceau lisse $\F$ n'est rien d'autre qu'un $\Lambda[G]$-module $F$. Le revêtement minimal est alors donné par $L^H$, où $H=\ker(G\to\Aut_\Lambda(M))$ ; c'est un simple calcul d'algèbre linéaire sur le $K$-espace vectoriel $L$.

\paragraph{Complexité du calcul d'un revêtement trivialisant} Tous nos algorithmes utiliseront la représentation par fibre générique et revêtement trivialisant. Soit $X$ une courbe intègre lisse sur $k_0$. Soit $\F$ un faisceau lisse sur $X$, représenté par un schéma en groupes $F=\bigsqcup_{i=1}^r F_i\to X$ où les $F_i$ sont connexes et étales sur $X$. Supposons $X$ et les $F_i$ décrites par un modèle plan (voir annexe \ref{subsec:modplan}). Pour chaque $i\in \{1\dots r\}$, notons $d_i$ le degré de $F_i$ et $f_i$ le degré de $F_i\to X$. Notons $f=f_1+\dots+f_r$ le degré de $F\to X$, c'est-à-dire le cardinal de la fibre générique de $\F$. Soit $d=\max(d_1,\dots,d_r,f)$. D'après l'annexe \ref{sec:prodfib}, le calcul de $F_i\times_X F_j$ nécessite $d_i^{4d_i}f_j^{8f_j}$ opérations, et la courbe obtenue est de degré $O(d_if_j^2)$. L'algorithme peut commencer par la composante de $F$ correspondant à la section nulle : quitte à la remplacer par $X$, son degré est celui de $X$. Comme il y a au plus $\log_2 f$ étapes de récursion, le degré de la courbe trouvée à la fin est $O(d^{1+f})$. Le nombre de calculs à effectuer est $O(d^{13fd^{1+f}})$.

\subsection{Morphismes, noyaux et conoyaux}

Soient $\F,\F'$ deux faisceaux lisses de $\Lambda$-modules sur un schéma intègre $X$. Soit $Y\to X$ un revêtement galoisien de groupe $G$ et de degré $d$ qui trivialise $\F$. Définissons de même $Y',G',d'$ pour $\F'$. Un revêtement galoisien $W$ de $X$ ayant pour corps de fonctions la composée de ceux de $Y$ et $Y'$ trivialise $\F$. Son degré est borné par $dd'$. 
Soit $H=\Aut(W|X)$. Notons $M$ et $M'$ les $\Lambda[H]$-modules $\HH^0(W,f^\star\F)$ et $\HH^0(W,f^\star\F')$. Le revêtement $W$ est l'une des composantes connexes de $Y\times_X Y'$ (qui est encore galoisienne sur $X$). Dans le cas général, le calcul de $W$ est donc de complexité élémentaire en les entrées (voir annexe \ref{subsec:decprim}). Dans le cas des courbes, le produit fibré $Y\times_X Y'$ se détermine comme décrit dans l'annexe \ref{sec:prodfib}. 
 
Un morphisme $\alpha\colon\F\to \F'$ peut être représenté par un morphisme de $X$-schémas en groupes $T\to T'$, ou par un morphisme de $\Lambda[H]$-modules $M\to M'$. Le faisceau $\coker \alpha$ est encore lisse puisque  $\phi^\star\coker\alpha=\coker(\phi^\star \alpha)$ est le conoyau d'un morphisme de faisceaux constants et est donc constant \cite[093J]{stacks}. Le noyau et le conoyau de $\alpha$ se calculent alors dans la catégorie $\Mod_{\Lambda[H]}$. Soient $m,m'$ les nombres de générateurs donnés de $M$ et $M'$. Les calculs de $\ker\alpha$ et $\coker\alpha$ se font donc en $O(\max(m,m')^3)$ opérations une fois que $W$ est construit.

\subsection{Faisceaux lisses sur les courbes nodales}\label{subsec:replisnod}

Voici comment seront représentés les faisceaux lisses sur les courbes nodales. Soient $X$ une courbe intègre nodale, $\tX$ sa normalisée, $\F$ un faisceau lisse sur $X$, et $Y\to X$ un revêtement trivialisant de $\F$. La courbe $X$ est représentée par la donnée d'un modèle plan $\tX_P$ de $\tX$ dont les singularités sont images de points de $\tX$ d'image régulière dans $X$, et des couples $(x,y)\in \tX_P(k)^2$ de points marqués qui sont les antécédents des points nodaux de $X$. L'utilité de distinguer les points de $\tX$ au-dessus des points nodaux de $X$ deviendra claire dans la section \ref{sec:constrnod} ; un tel modèle s'obtient par transformations quadratiques à partir d'un modèle plan de $X$ en éclatant chacun des points nodaux.
La courbe $Y$ est décrite par sa normalisée $\tY=Y\times_X\tX$, qui n'est peut-être pas connexe : c'est une réunion disjointe $\tY_P$ de courbes planes (dont les points singuliers ne sont pas au-dessus de points singuliers de $Y$) avec des couples de points marqués (n'appartenant pas nécessairement à une même composante connexe). Le faisceau $\F$ est donné par $\tY_P\to \tX_P$ et le $\Aut(\tY|\tX)$-module $F=\HH^0(Y,\F)$.

\section{Images directes de faisceaux lisses}\label{sec:opfaisclisses}

Soit $f\colon Y\to X$ un morphisme quasi-fini de variétés régulières intègres de même dimension sur $k$. Nous allons décrire comment, étant donné un faisceau lisse sur $Y$, calculer les faisceaux constructibles $f_\star\F$ et $f_!\F$ sur $X$. Rappelons que par le théorème principal de Zariski \cite[02LR]{stacks}, $f$ peut être décomposé en
\[
\begin{tikzcd}
Y \arrow[r,"j"]\arrow[rd,"f",swap] & X' \arrow[d,"\nu"] \\
& X
\end{tikzcd}
\]
où $j$ est une immersion ouverte quasi-compacte et la normalisation $X'\to X$ de $X$ dans $Y$ est un morphisme fini. Comme $X$ et $Y$ sont régulières, le théorème de "platitude miraculeuse" \cite[00R4]{stacks} assure que $X'\to X$ est plat, c'est-à-dire localement libre \cite[02KB]{stacks}. Par conséquent, il suffit de considérer deux cas particuliers : celui où $f$ est une immersion ouverte, et celui où $f$ est fini localement libre.\\

Voici un argument ad hoc pour calculer une telle décomposition $Y\to X'\to X$.
Comme $f$ est quasi-fini, il est génériquement fini \cite[03I1]{stacks}. Le morphisme $Y\to X$ correspond (localement sur $X$) à une extension $A\to B$ de $k$-algèbres, et une extension $K\to L$ de corps de fonctions. Soient $x_1,\dots,x_r$ des générateurs de $B$ comme $A$-algèbre, et $f_1,\dots,f_r\in K[t]$ leurs polynômes minimaux respectifs. Pour tout $i\in \{1\dots r\}$, soit $g_i$ le produit des dénominateurs de $f_i$. Alors en posant $U=\Spec A[(g_1\cdots g_r)^{-1}]$, le morphisme \[X'=U\times_X Y\to U\] est fini.\\

\begin{rk} Dans le cas où $X$ et $Y$ sont des courbes intègres régulières, la normalisation de $X$ dans $Y$ est simplement l'image inverse de $X$ dans la compactification régulière de $Y$.
\end{rk}

\subsection{Image directe}

\paragraph{Par une immersion ouverte}

Le résultat permettant de calculer le poussé en avant d'un faisceau lisse par une immersion ouverte est le suivant.

\begin{prop}\cite[Corollary 5.8]{jinbi_jin}\label{prop:jstarnorm} Soit $k$ un corps. Soient $X$ un schéma de type fini sur $k$, et $j\colon U\to X$ une immersion ouverte telle que la normalisation de $X$ dans $U$ soit $X$ (ce qui est toujours le cas si $X$ est normal). Soit $\F$ un faisceau fini localement constant d'ensembles sur $U$, représenté par un $U$-schéma fini étale $F$. Alors $j_\star\F$ est représenté par le lieu étale sur $X$ de la normalisation de $X$ dans $F$.
\end{prop}

Si $j\colon U\to X$ est une immersion ouverte de variétés normales et $\F$ est un faisceau sur $U$ représenté par un $U$-schéma $F$, le faisceau $j_\star \F$ se détermine donc en calculant la normalisation de $X$ dans $F$, puis le lieu étale (c'est-à-dire non ramifié) de cette dernière au-dessus de $X$.

\begin{rk} Si $\F$ est un faisceau de groupes abéliens, voici comment se calcule la loi de groupe sur le schéma représentant $j_\star\F$. Notons que comme \[j_\star(\F\times \F)=j_\star\F\times j_\star \F \] il y a un isomorphisme canonique de schémas \[j_\star (F\times_U F)\xrightarrow{\sim}j_\star F\times_X j_\star F.\] La loi de groupe $F\times_U F\to F$ est donnée. Notons $X'$ la normalisation de $X$ dans $F$, et $X''$ la normalisation de $X$ dans $F\times_U F$. Alors le morphisme $F\times_U F\to F\to X'$ donne lieu par propriété universelle de la normalisation relative \cite[035I]{stacks} à un morphisme de $X$-schémas $X''\to X'$. La restriction de celui-ci au lieu étale de $X'\to X$ donne un morphisme $j_\star (F\times_U F)\to j_\star F$.
\end{rk}

\paragraph{Par un morphisme fini localement libre}

Supposons que $f\colon Y\to X$ soit fini localement libre. Soit $\F$ un faisceau sur $Y$ représenté par un $Y$-schéma $F$. Alors par \cite[7.6, Th. 4]{blr}, la restriction de Weil $R_{Y\to X}(F)$ existe et représente $f_\star\F$. L'algorithme qui calcule la restriction de Weil se trouve dans la section \ref{sec:weilres}.

\subsection{Image directe à support propre}

\paragraph{Par une immersion ouverte}

\begin{rk} Soient $U$ un schéma, et $f\colon G\to U$ un $U$-schéma en groupes fini étale. Considérons la section nulle $e\colon U\to G$. Alors $fe=\id_U$ est fini étale, et comme $f$ est fini étale, $e$ l'est aussi. En particulier, $e(U)$ est ouvert-fermé dans $G$, c'est donc une composante connexe de $G$. De plus, $e$ a un inverse à gauche, c'est donc un isomorphisme sur son image. Le schéma $G$ s'écrit donc $G=U\sqcup G'$, et $G(U)=\{ \id \} \sqcup G'(U)$.
\end{rk}

\begin{lem} Soit $j\colon U\to X$ une immersion ouverte de schémas connexes. Soit $\F$ un faisceau lisse de $\Lambda$-modules sur $U$, représenté par un $U$-schéma en groupes $G\to U$. Notons $G=U\sqcup G'$, où $U$ est la section nulle. Alors le schéma $Y\coloneqq X\sqcup G'$, muni de la loi de groupe induite par celle de $U\sqcup G'$, représente le faisceau $j_!\F$.
\begin{proof} Soit $\alpha \colon T\to X$ étale. Supposons $T$ connexe. Alors $Y(T)=X(T)\sqcup {G'}(T)=\{\alpha \}\sqcup {G'}(T)$. Remarquons que si l'image de $T\to X$ est incluse dans $U$, alors $\Hom_X(T,G')=\Hom_U(T,G')$, et sinon $\Hom_X(T,G')=\varnothing$. Ainsi, $Y(T)=\Hom_U (T,G)=\F(T)$ dans le premier cas, et $Y(T)=\{\alpha\}$ dans le second, ce qui implique que $Y$ représente $j_!\F$.
\end{proof}
\end{lem}
\begin{figure}[H]
\begin{center}
\begin{tikzpicture}
\draw (-5,0) -- (-1,0);
\draw (1,0) -- (5,0);
\draw (-4,1.5) -- (-2,1.5);
\draw (-4,2) -- (-2,2);
\draw (-4,3) -- (-2,3);
\draw (-4,2.5) -- (-2,2.5);

\draw (1,1.5) -- (5,1.5);
\draw (2,3) -- (4,3);
\draw (2,2.5) -- (4,2.5);
\draw (2,2) -- (4,2);

\draw[->] (-3,1) -- (-3,0.5);
\draw[->] (3,1) -- (3,0.5);

\draw (-4,-0.1) -- (-4, 0.1);
\draw (-2,-0.1) -- (-2, 0.1);

\draw (2,-0.1) -- (2, 0.1);
\draw (4,-0.1) -- (4, 0.1);

\draw (-3,0) node[below]{$U$};
\draw (3,0) node[below]{$U$};

\draw (-5,2.25) node[left]{$G$};
\draw (5,2.25) node[right]{$j_!G$};

\draw (-1,0) node[right] {$X$};
\draw (5,0) node[right] {$X$};

\end{tikzpicture}
\caption{Prolongement par zéro d'un faisceau lisse}
\end{center}
\end{figure}
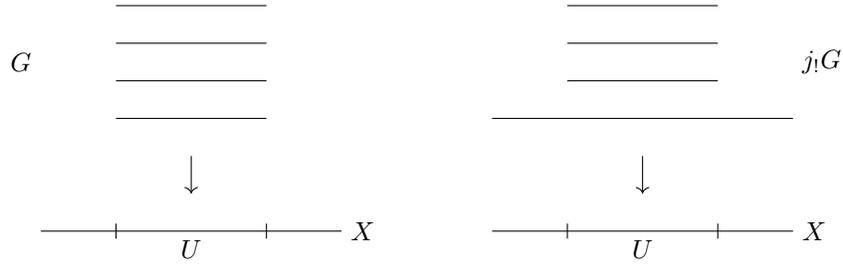

\paragraph{Par un morphisme fini localement libre} Pour un morphisme fini $f$, les foncteurs $f_!$ et $f_\star$ sont égaux ; la représentabilité de $f_\star\F$ a été traitée ci-dessus.\\

En résumé, comme la détermination des composantes connexes d'un schéma et le calcul d'une restriction de Weil sont de complexité élémentaire (voir annexes \ref{subsec:decprim} et \ref{sec:weilres}), nous avons montré le résultat suivant.

\begin{prop}
Il existe un algorithme de complexité élémentaire qui, étant donné un morphisme de $k$-variétés régulières $f\colon Y\to X$ et un faisceau lisse $\F$ sur $Y$ représenté par un $Y$-schéma en groupes étale, calcule des $X$-schémas en groupes représentant les faisceaux $f_\star\F$ et $f_!\F$.
\end{prop}

\section{Faisceaux constructibles sur les courbes lisses}

Nous allons maintenant décrire trois représentations différentes des faisceaux constructibles de $\Lambda$-modules sur une courbe intègre lisse sur $k$. Tout d'abord, un faisceau constructible peut toujours être exprimé comme conoyau d'un morphisme de la forme $f_!\Lambda\to g_!\Lambda$, où $f$ et $g$ sont étales. De même, il peut être exprimé comme noyau d'un morphisme $p_\star M\to q_\star N$, où $p$ et $q$ sont finis et $M,N$ sont des $\Lambda$-modules de type fini. Enfin, il peut être défini par recollement relativement à un ouvert sur lequel il est lisse.

\subsection{Représentation comme conoyau "$(!)$"}

\begin{prop}\cite[IX, Prop. 2.7]{sga43}\label{prop:coprod} Soit $A$ un anneau noethérien. Soient $X$ un schéma quasi-compact et quasi-séparé et $\F$ un faisceau de $A$-modules sur $X$. Pour que $\F$ soit constructible, il faut et il suffit qu’il soit isomorphe au conoyau d’un morphisme  $f_!A\to g_!A$, où $f\colon X_1\to X$ et $g\colon X_2\to X$ sont deux morphismes étales de présentation finie.
\end{prop}

La preuve classique de cette proposition consiste à considérer pour chaque point $x\in X$ un point géométrique $\bar x$ d'image $x$, et chaque $s\in \F_{\bar x}$ un schéma affine $\phi_{\bar x,s}\colon T_{\bar x,s}\to X$ tel que $s$ appartienne à l'image de $\F(T_{\bar x,s})\to \F_{\bar x}$ ; cela signifie qu'il y a un morphisme $(\phi_{\bar x,s})_!\Lambda\to \F$ dont l'image contient $s$. La constructibilité du faisceau permet alors d'appliquer un argument de noethérianité à $\bigoplus_{\bar x,s}\phi_{\bar x,s}$. En particulier, si $X$ est une courbe intègre sur $k$, cet argument de finitude peut être rendu explicite : si $U$ est un ouvert de lissité de $\F$, il suffit de considérer la fibre générique de $U$, ainsi que les fibres en les points du complémentaire zéro-dimensionnel de $U$. Pour la fibre générique, si $\phi\colon V\to U$ est un revêtement trivialisant, le morphisme $\F(V)\to \F_{\bareta}$ est un isomorphisme.

\paragraph{La représentation $(!)$} Soit $X$ une courbe intègre sur $k$. Un faisceau constructible $\F$ sur $X$ sera représenté par deux morphismes étales $f\colon X_1\to X$ et $g\colon X_2\to X$ ainsi que d'un morphisme $u\colon f_!\Lambda\to g_!\Lambda$ tel que $\F=\coker(u)$ ; ce morphisme sera décrit explicitement de l'une des deux façons suivantes.

\paragraph{Représentation d'un morphisme $f_!\Lambda\to g_!\Lambda$} Soient $f\colon X_1\to X$ et $g\colon X_2\to X$ deux morphismes étales. Les faisceaux $f_!\Lambda$ et $g_!\Lambda$ sont représentables, et nous savons d'après \ref{sec:opfaisclisses} déterminer explicitement des schémas en groupes $Y_1,Y_2$ qui les représentent. Un morphisme $f_!\Lambda\to g_!\Lambda$ est alors simplement un morphisme de schémas en groupes $Y_1\to Y_2$. Un tel morphisme peut également être représenté de façon plus succincte : par adjonction, il est défini par la donnée d'un morphisme $\Lambda \to f^\star g_!\Lambda$, c'est-à-dire d'une section globale de $f^\star g_!\Lambda$ pour chaque composante connexe de $X_1$. Comme $f$ est étale, cela revient à se donner un élément de $g_!\Lambda(X_1)$. Or \cite[p209]{fulei} $g_!\Lambda(X_1)$ est l'ensemble des $s\in \Lambda(X_1\times_X X_2)=\Lambda^{\pi_0(X_1\times_X X_2)}$ de support propre sur $X_1$. Le support d'une telle section est une réunion de composantes connexes de $X_1\times_X X_2$ : il suffit de déterminer quelles composantes connexes sont propres (c'est-à-dire finies puisque $f$ est étale, voir \cite[02OG]{stacks}) sur $X_1$.

\paragraph{Calcul du morphisme $f_!\Lambda\to \F$ associé à une section de $f^\star \F$}
\label{morphsec}
Soit $f\colon Y\to X$ un morphisme étale. Il se décompose en $f=gj$, où $j\colon Y\to X'$ est une immersion ouverte et $g\colon X'\to X$ est fini localement libre. Soit $\F$ un faisceau lisse sur $X$, représenté par un schéma en groupes fini étale $F\to X$. Soit $s\in \HH^0(Y,f^\star \F)$, qui correspond à un morphisme de $Y$-schémas en groupes $\coprod_\Lambda Y\to F\times_X Y$.
Par adjonction, $s$ définit un morphisme $j_!\Lambda\to g^\star \F$. Le faisceau $j_!\Lambda$ est représenté par $X'\sqcup\coprod_{\Lambda - \{0\}}Y$.
Les $n-1$ morphismes $Y\to F\times_X X'$ sont simplement les composées $Y\to F\times_X Y\to F\times_X X'$. Le morphisme $X'\to F\times_X X'$ est la section nulle, déduite de la section nulle $X\to F$ par changement de base. Ce morphisme $j_!\Lambda \to g^\star \F$ donne à nouveau par adjonction un morphisme $f_!\Lambda=g_\star j_!\Lambda \to \F$, qui se calcule explicitement puisque $g_\star j_!\Lambda$ est la restriction de Weil $R_{X'\to X}(j_!\Lambda)$.

\subsection{Représentation comme noyau "$(\star)$"}

\begin{prop}\cite[IX, Prop. 2.14.(ii)]{sga43}\label{prop:rep*} Soit $X$ un schéma de type fini sur un corps ou sur $\ZZ$. Soit $A$ un anneau noethérien. Soit $\F$ un faisceau constructible de $A$-modules sur $X$. Il existe des schémas noethériens intègres $X_i, i=1,\dots,m$, des morphismes finis $p_i\colon X_i\to X$ et des $A$-modules de type fini $M_i$ tels qu'il y ait un monomorphisme \[ \F \hookrightarrow \bigoplus_{i=1}^m {p_i}_\star M_i.\]
\end{prop}

La preuve classique de cette proposition consiste (dans le cas irréductible) à considérer une décomposition $X=\bigcup_i U_i$ en parties localement fermées telles que chaque $\F|_{U_i}$ soit lisse. Il existe des morphismes finis étales $U_i'\to U_i$ tels que chaque $\F|_{U_i'}$ soit constant. La normalisation $X_i$ de $U_i'$ fournit le morphisme $p_i\colon X_i\to X$ recherché. Notons que si $X$ est normal, les $X_i$ le sont aussi.\\

\paragraph{Description d'un morphisme $p_\star M\to q_\star N$} Soit $X$ une courbe intègre sur $k$. Soient $p\colon Y\to X$ et $q\colon Z\to X$ des morphismes finis, avec $Y,Z$ intègres. Soient $M,N$ deux groupes abéliens finis. Montrons, selon les dimensions de $Y$ et $Z$, comment décrire un morphisme $p_\star M\to q_\star N$, c'est-à-dire par adjonction un morphisme $q^\star p_\star M\to N$. \begin{enumerate}
\item Si $Y$ et $Z$ sont deux courbes, $p_\star M\to q_\star N$ est un morphisme de $X$-schémas entre les restrictions de Weil $R_p(Y\times M)\to R_q(Z\times N)$.
\item Si $Z$ est un point fermé de $X$, $q^\star p_\star M$ est le faisceau de fibre $(p_\star M)_q=\bigoplus_{Y_q}M$ sur le point. Un morphisme $q^\star p_\star M\to N$ est donc simplement un morphisme de $\Lambda$-modules $\bigoplus_{Y_q}M\to N$. 
\item Si $Y$ est un point fermé et $Z$ est une courbe, montrons que le seul morphisme $p_\star M\to q_\star N$ est le morphisme nul. Notons $j$ l'inclusion d'un ouvert de $X$ ne contenant pas $Y$. La transformation naturelle $\id\to j_\star j^\star$ fournit un diagramme commutatif :
\[
\begin{tikzcd}p_\star M \arrow[r]\arrow[d] & j_\star j^\star p_\star M \arrow[d] \\
q_\star N \arrow[r] & j_\star j^\star q_\star N
\end{tikzcd}
\]
D'une part, $j_\star j^\star p_\star M=0$. D'autre part, notons $q'\colon Z\times_X U\to U$ et $j'\colon Z\times_X U\to Z$. Le morphisme $q_\star N\to j_\star j^\star q_\star N$ est un isomorphisme, car il s'identifie via les isomorphismes canoniques $j^\star q_\star N\xrightarrow{\sim} q'_\star j'^\star N$ et $N\xrightarrow{\sim} j'_\star j'^\star N$ au morphisme identité de $q_\star N$. Par conséquent, le morphisme $p_\star M\to q_\star N$ est nul.
\end{enumerate}

\paragraph{La représentation $(\star)$} Reprenons les notations de la proposition \ref{prop:rep*}, en supposant que $X$ est une courbe intègre sur $k$. En appliquant la proposition au conoyau $\G$ de $\F\to \bigoplus_{i}{p_i}_\star M_i$, le faisceau $\F$ s'exprime comme le noyau d'un morphisme $\bigoplus_{i}{p_i}_\star M_i\to \bigoplus_j {q_j}_\star N_j$.
La représentation $(\star)$ de $\F$ est la donnée de ces morphismes ${p_i}_\star M_i\to {q_j}_\star N_j$, décrits explicitement comme dans le paragraphe précédent.

\subsection{Représentation par recollement "$(\sqcup)$"}

Soient $X$ une courbe intègre sur $k$ et $\F$ un faisceau constructible sur $X$. Soit $U$ un ouvert de $X$ sur lequel $\F$ est localement constant. Notons $Z$ le fermé réduit complémentaire. Notons $j\colon U\to X, i\colon Z\to X$ les inclusions, et $z_1,\dots,z_r$ les points fermés de $Z$. Par recollement (voir section \ref{subsec:recollement}), $\F$ est uniquement déterminé par les données suivantes :
\begin{itemize}[label=$\bullet$]
\item le faisceau lisse $\mathscr{L}=j^\star\F$ ;
\item le faisceau $\F_Z=i^\star\F$, défini par les groupes abéliens finis $\F_{z_1},\dots,\F_{z_r}$ ;
\item le morphisme de recollement $\phi \colon \F_Z\to i^\star j_\star\mathscr{L}$.
\end{itemize}
Nous appellerons représentation par recollement ou représentation $(\sqcup)$ de $\F$ par rapport à $(U,Z)$ la donnée du triplet $(\mathscr{L},\F_Z,\phi)$. Toute donnée de cette forme définit un faisceau constructible sur $X$.
Le carré suivant, où les flèches non étiquetées désignent les unités d'adjonction, est alors cartésien.
\[
\begin{tikzcd}
\F \arrow[r]\arrow[d,] & i_\star \F_Z \arrow[d,"i_\star \phi"] \\
j_\star \mathscr{L}\arrow[r,"c"] & i_\star i^\star j_\star \mathscr{L}
\end{tikzcd}
\]
Soit $Y\to X$ étale. Le diagramme cartésien ci-dessus montre que \[ \F(Y)= \{ (s,t)\in i_\star\F_Z(Y)\times j_\star \mathscr{L}(Y)\mid i_\star\phi(s)=c(t) \}.\]
D'une part, \[i_\star\F_Z(Y)=\F_Z(Y|_Z)=\bigoplus_{i=1}^r\bigoplus_{y\in Y_{z_i}(k)}\F_{z_i}.\] D'autre part, $j_\star\mathscr{L}(Y)=\mathscr{L}(Y|_U)$ et le morphisme \[ j_\star\mathscr{L}(Y)\to i_\star i^\star j_\star\mathscr{L}(Y)=\bigoplus_{i=1}^r\bigoplus_{y\in Y_{z_i}(k)}(j_\star\mathscr{L})_{z_i}\] envoie une section $s\in \mathscr{L}(Y)$ sur $(s_{z_i})_{i\in \{ 1\dots r\},y\in Y_{z_i}(k)}$. En résumé :

\[ \F(Y)=\{ (s,t)\in (\bigoplus_{i=1}^r\bigoplus_{y\in Y_{z_i}(k)}\F_{z_i}) \times j_\star\mathscr{L}(Y)\mid \forall i\in \{1\dots r\},\forall y\in Y_{z_i}(k) ,\phi(s_{z_i,y})=t_{z_i}\}.\]

\begin{rk} Nous disposons de deux moyens pour déterminer la fibre en chaque $z_i$ de $j_\star \mathscr{L}$. D'une part, si $\mathscr{L}$ est donné comme un $U$-schéma en groupes fini étale, le schéma $F$ représentant $j_\star\mathscr{L}$ se calcule comme décrit en \ref{prop:jstarnorm}, et la fibre de $j_\star\mathscr{L}$ en $z_i$ est simplement $F\times_X z_i$. Afin de déterminer cette fibre explicitement comme groupe fini, il est nécessaire de déterminer ses points. D'autre part, si $X$ est lisse et $\mathscr{L}$ est donné par un revêtement trivialisant $X'\to X$ et une action de $\Aut(X'|X)$ sur $F\coloneqq \F(X')$,  alors $(j_\star\mathscr{L})_{z_i}=F^{I_{z_i}}$ , où $I_{z_i}\subset \Aut(X'|X)$ désigne le stabilisateur dans $\Aut(X'|X)$ d'un antécédent de $z_i$. Il suffit donc de déterminer ces groupes d'inertie.
\end{rk}

\begin{rk} Le faisceau constructible défini par $(\mathscr{L},\F_Z,\phi)$ est lisse si et seulement si la fibre $M$ de $\mathscr{L}$ est invariante sous les groupes d'inertie $I_z,z\in Z$ et le morphisme $\phi$ est un isomorphisme.
\end{rk}

\subsection{Représentation $(\sqcup)$ d'un faisceau constructible représentable}\label{repet}

Soit $X$ une courbe intègre sur $k$.
Soit $\F$ un faisceau constructible sur $X$ représenté par un schéma étale $F\to X$. Montrons comment représenter $\F$ par recollement. Comme $f$ est étale, il est quasi-fini. Soit $j\colon U\to X$ un ouvert tel que $F'\coloneqq F\times_X U$, qui représente alors $j^\star\F$, soit fini sur $U$. Notons $X'$ le normalisé de $X$ dans $F$, qui est encore le normalisé de $X$ dans l'ouvert $F'$ de $F$ : la situation est résumée par le diagramme commutatif suivant.
\[
\begin{tikzcd}
F' \arrow[r] \arrow[d]& F \arrow[r]\arrow[d] & X'\arrow[dl] \\
U \arrow[r] & X &
\end{tikzcd}\]
Comme $F\to X$ est étale, $F$ est inclus dans le lieu étale de $X'\to X$, qui représente $j_\star j^\star\F$. Ceci fournit directement la flèche injective $\F\to j_\star j^\star \F$, qui sur les fibres en les points de $X-U$ donne le morphisme de recollement. 
\\
Un morphisme $f\colon F\to G$ de schémas finis étales sur $X$, peut également être représenté par recollement, en choisissant pour $U$ un ouvert sur lequel les deux schémas sont finis. Le morphisme de faisceaux lisses $f|_U$ sur $U$ est alors simplement un morphisme de schémas en groupes, et le morphisme sur la fibre en un point $z\colon\Spec k\to X-U$ est simplement $F_z\to G_z$.

\subsection{\'{E}quivalence entre $(\star)$ et $(\sqcup)$}

Soit $X$ une courbe intègre lisse sur $k$. Soit $\F$ un faisceau constructible sur $X$. 

\paragraph{$(\star)\rightarrow (\sqcup)$} Voici comment donner une représentation par recollement d'un faisceau défini par une représentation $(\star)$. Dans cette représentation, le faisceau est somme directe de noyaux de morphismes de la forme $p_\star M\to q_\star N$ où $p,q$ sont des morphismes finis de cible $X$.
Il suffit donc de savoir représenter par recollement un faisceau de la forme $p_\star M$ avec $p\colon Y\to X$ fini, ainsi que les morphismes entre de tels faisceaux. Il y a deux cas à considérer :
\begin{enumerate}
\item Si $Y$ est une courbe lisse, $p_\star M$ est la restriction de Weil $R_p(Y\times M)$, qui est étale sur $X$ \cite[Prop. 4.9]{scheiderer}. Nous avons décrit dans la section \ref{repet} comment représenter par recollement les (morphismes de) $X$-schémas étales.
\item Si $Y$ est l'inclusion d'un point fermé de $X$, le faisceau $p_\star M$ est constant sur $X-Y$ de valeur 0, a pour fibre $M$ en $Y$, et son morphisme de recollement est nul.
\end{enumerate}

\paragraph{$(\sqcup)\rightarrow (\star)$} Supposons $\F$ représenté par recollement relativement à un ouvert de lissité $j\colon U\to X$ et son fermé réduit complémentaire $i\colon Z\to X$. Soit $V$ un revêtement galoisien de $U$ trivialisant le faisceau lisse $j^\star\F$ de fibre $M$. Nous avons décrit dans la section \ref{subsec:replisses} comment calculer un $U$-schéma $F$ représentant $j^\star\F$. La normalisation $p\colon V'\to X$ de $X$ dans $V$ fournit un premier morphisme fini ; les inclusions $z\colon\Spec k\to X$ des points de $Z$ fournissent les autres. Décrivons par recollement l'injection \[ \phi\colon\F \longrightarrow p_\star M \oplus \bigoplus_{z\in Z} \F_z.\]
Sur $U$, c'est le morphisme \[ j^\star\F\longrightarrow j^\star p_\star M =j^\star {p}_\star p^\star \F=p{|_U}_\star {p|_U}^\star j^\star \F \]
obtenu en appliquant l'unité d'adjonction $\id\to p{|_U}_\star p|_U^\star$ au faisceau lisse $j_\star\F$. En termes de schémas, c'est le morphisme $F \to R_{Y\to X}(F\times_U V)$ qui est calculable explicitement (voir annexe \ref{sec:weilres}). Sur $z\in Z$, le morphisme $\phi$ est l'injection $0\oplus \F_z\to M^{I_z}\oplus \F_z$. Cette description de $\phi$ par recollement permet de calculer $\coker\phi$, et d'obtenir par la même procédure une injection \[ \coker\phi \to q_\star N \oplus \bigoplus_{w\in W} w_\star N_w \] où $q$ est la normalisation de $X$ dans un schéma étale sur $X$, et $W$ est un fermé zéro-dimensionnel de $X$. Le morphisme \[ p_\star M\oplus\bigoplus_z z_\star \F_z \to q_\star N\oplus\bigoplus_w w_\star N_w \]
qui s'en déduit a pour noyau $\F$. Ce morphisme est représenté de la façon suivante. Le morphisme ${p}_\star M\to q_\star N$ est la composée \[ {p}_\star M\longrightarrow j^\star\coker\phi \longrightarrow {q}_\star N. \] Il est décrit par le morphisme correspondant entre restrictions de Weil. Les morphismes \[ {p}_\star M\oplus\bigoplus_z z_\star\F_z \to \bigoplus_{w}w_\star N_w\] sont décrits fibre à fibre.

\subsection{\'{E}quivalence entre $(!)$ et $(\sqcup)$}

Soit $X$ une courbe intègre lisse sur $k$. Soit $\F$ un faisceau constructible sur $X$. 

\paragraph{$(!)\rightarrow (\sqcup)$} \'{E}tant donné une représentation de $\F$ comme $\coker (u\colon f_!\Lambda\to g_!\Lambda)$, où $f$ est $g$ sont étales de type fini, il est possible de calculer explicitement le morphisme de $X$-schémas étales représenté par $u$. Nous avons décrit dans la section \ref{repet} comment donner une représentation par recollement d'un tel morphisme. 

\paragraph{$(\sqcup)\rightarrow (!)$} Supposons $\F$ défini par recollement relativement à un couple ouvert-fermé 
\[
\begin{tikzcd}
U \arrow[rr,hook,"j"]& & X \arrow[rr,hookleftarrow,"i"]& & Z 
\end{tikzcd}
\]
par la donnée d'un faisceau lisse $\mathscr{L}$ sur $U$ trivialisé par un revêtement étale $V$, des fibres $\F_z$ en les points de $Z$ et un morphisme $\phi\colon i^\star\F\to i^\star j_\star \mathscr{L}$.
Remarquons que la preuve de la proposition \ref{prop:coprod} est constructive, mis à part pour la détermination, pour $z\in Z$ et $s\in \F_z$, d'un morphisme étale $f\colon T\to X$ tel que l'image de $\F(T)\to\F_z$ contienne $s$. \\

Soient donc $z\in Z$ et $s\in \F_z$. Rappelons que $\F=j_\star \mathscr{L}\times_{i_\star i^\star j_\star \mathscr{L}} i_\star i^\star\F$. Par conséquent, l'image de $\F(T)\to\F_z$ contient $s$ si et seulement si l'image de $s$ par le morphisme de recollement $\phi_T\colon i^\star\F(T)=\bigoplus_z\F_z\to i^\star j_\star \mathscr{L}(T)=\bigoplus_z (j_\star \mathscr{L})_z$ a un antécédent par $j_\star \mathscr{L}(T)\to \bigoplus_z (j_\star \mathscr{L})_z$. Notons $G$ le schéma représentant $\mathscr{L}$. Le faisceau $\F'\coloneqq j_\star \mathscr{L} $ est représentable par un schéma étale $F$ sur $X$ et la fibre $\F'_z$ est simplement $F\times_X z$. Le diagramme commutatif suivant résume les notations.\\

\[
\begin{tikzcd}
G\arrow[r]\arrow[d]& F \arrow[d]& \arrow[l]F_z\arrow[d] \\
U\arrow[r] & X &\arrow[l] z
\end{tikzcd}
\]

Dans le cas où $\F$ est lisse sur $X$, l'algorithme classique consiste à construire par changements de base successifs un revêtement étale $Y\to X$ qui trivialise $\F$ ; le morphisme $\F(Y)\to \F_z$ est alors un isomorphisme.
L'algorithme suivant, qui s'inspire de celui-ci, permettra d'obtenir le schéma souhaité pour les faisceaux constructibles. \\

Posons $U_0=U$, $X_0=X$, $F_0=F$, $G_0=G$ et $z_0=z$. 
L'algorithme construit une suite $(X_i,U_i,F_i,G_i,z_i)_{i\geqslant 0}$ de la façon suivante. Soit $i\in\NN$.
Soit $C$ une composante connexe de $G_i$. Notons $Y$ le lieu étale de la normalisation de $X_i$ dans $C$.\begin{itemize}[label=$\bullet$]
\item Si $\deg(C\to U_i)>1$ et la fibre $Y_{z_i}$ est non vide, soit $z_{i+1}\in Y_{z_i}$. Posons alors $U_{i+1}=C$, $X_{i+1}=Y$, $G_{i+1}=G_i\times_{U_i} C$ et $F_{i+1}=F_i\times_{X_i}Y$. 
\item Sinon, on revient au choix de $C$ ; si pour toute composante connexe $C$ de $G_i$, $\deg(C\to U_i)>1$ ou $Y_{z_i}$ est vide, l'algorithme s'arrête et renvoie $Y_f\coloneqq X_i$. 
\end{itemize}

Cet algorithme s'arrête au bout d'un nombre fini d'opérations : en effet, pour chaque entier $i\in\{0\dots i_f-1\}$, $G_{i+1}$ a strictement plus de composantes connexes que $G_{i}$ alors que $\deg(G_{i+1}\to U_{i+1})=\deg(G\to U)$. 
Notons $i_f$ la valeur de l'indice $i$ lorsque l'algorithme termine. Nous allons prouver que $\F(Y_f)\to\F_z$ est surjectif.  

\begin{lem}  Pour tout entier $i\in \{0\dots i_f\}$, $j_{i,\star} (G_i)=F_i$.
\begin{proof} Montrons-le par récurrence sur l'entier $i$. L'assertion concernant $F_0$ vient directement de sa définition. Soit désormais $i$ tel que $j_{i,\star} (G_i)=F_i$. Soient $C$ la composante connexe de $G_i$ choisie et $Y$ la normalisation de $X_i$ dans $C$. Il y a un diagramme commutatif : 
\[
\begin{tikzcd}
& G_{i+1} \arrow[dl] \arrow[dd] \arrow[rr] & & F_{i+1}\arrow[dl]\arrow[dd]  \\
G_i \arrow[dd]\arrow[rr,crossing over] & & F_i& \\
& C \arrow[dl]\arrow[rr] & & Y \arrow[dl]  \\
U_{i} \arrow[rr,"j_i"] & & X_i\arrow[from=uu, crossing over] & 
\end{tikzcd}
\]
Nous savons que $j_{i+1,\star}(G_{i+1})$ est le lieu étale de la normalisation $\tY$ de $Y$ dans $G_i$. De plus, par hypothèse de récurrence, $j_{i,\star} G_i=F_i$ est le lieu étale de la normalisation $\tX_i$ de $X_i$ dans $G_i$. Comme le morphisme $Y\to X_i$ est lisse, le changement de base $-\times_{X_i} Y$ commute à la normalisation \cite[03GV]{stacks} et $\tY=\tX_i\times_{X_i} Y$. Enfin, comme $Y$ est plat sur $X_i$ et $\tX_i\to X_i$ est de présentation finie, le lieu étale de $\tX_i\times_{X_i} Y\to Y$ est le changement de base à $Y$ du lieu étale de $\tX_i\to X_i$ \cite[0476]{stacks}. Cela signifie que le lieu étale de $\tY\to Y$, qui représente $j_{i+1,\star} (G_{i+1})$, est $F_{i+1}$.
\end{proof}
\end{lem}

\begin{lem}[3.3.6.0.2] Si, pour toute composante connexe $C$ de $G_i$ telle que $\deg(C_i\to U_i)>1$, le lieu étale de la normalisation de $X_i$ dans $C$ ne contient aucun point au-dessus de $z_i$, alors le morphisme $F_i(X_i)\to F_{i,z_i}$ est un isomorphisme.

\begin{proof} \'{E}crivons $G_i=\bigsqcup_\alpha C_\alpha$. Alors $j_{i,\star} G_i$ est représenté par $F_i=\bigsqcup_\alpha X_{i,\alpha}$, où $X_{i,\alpha}$ est le lieu étale de la normalisation de $X_i$ dans $C_\alpha$. La normalisation de $X_i$ dans $U_i$ étant $X_i$,  \[ j_{i,\star} G_i =\coprod_{\deg(C_\alpha\to U_i)=1} X \sqcup \coprod_{\deg(C_\alpha\to U_{i})>1}X_{i,\alpha}.\] Par conséquent, l'hypothèse assure que \[ |F_{i,z_i}|=|\{ \alpha \mid \deg(C_\alpha\to U_i)=1\}|=:d.\]
De plus, pour tout $T\to X_i$ voisinage étale de $z_i$, $F_i(T)=G_i(T\times_X U)$ contient au moins $d$ éléments (les flèches $T\times_{X_i} U_i\to C_\alpha$ avec $\deg(C_\alpha\to U)=1$), qui sont préservés par les flèches de restriction $F_i(T)\to F_i(T')$ et sont donc encore distincts dans $F_{i,z_i}=\colim_{(T,t)\to (Y,z)}F_i(T)$. Enfin, $F_i(X_i)=G_i(U_i)$ contient exactement $d$ éléments, puisque comme $G_i\to U_i$ est fini étale, les sections $U_i\to G_i$ sont en bijection avec les composantes connexes de $G_i$ de degré 1 sur $U_i$.
\end{proof}
\end{lem}

\begin{prop} Le morphisme $\F(Y_f)\to \F_z$ est un isomorphisme.
\begin{proof}
\`{A} chaque étape de la boucle, $F_{i,z_i}=F_{z}$ et comme $F_i\to F$ est étale, $F_i(X_i)=F(X_i)$. Par conséquent, lorsque l'algorithme s'est arrêté, $F(X_i)\to F_z$ est un isomorphisme. Choisissons un ouvert $Y'$ de $Y_f$ ne contenant qu'un seul point au-dessus de $Z$, d'image $z$. Alors $F(Y')=F(Y_f)$ et le morphisme $F(Y')\to F_z$ est encore un isomorphisme. 
Rappelons que \[\F_z=\colim_{T\to (X,z)} [\F|_Z(T\times_X Z)\times_{F(T\times_X Z)}F(T)]\] et que la catégorie des voisinages étales connexes $T$ de $(X,z)$ n'ayant au-dessus de $Z$ qu'un seul point d'image $z$ est cofinale dans celle des voisinages étale de $(X,z)$. Par conséquent, $\F_z=\colim_T \F_z\times_{F_z}F(T)$. Le morphisme $F(Y')\to F_z$ étant bijectif, le morphisme $\F(Y')\to \F_z$ l'est encore. Par conséquent, le morphisme $\F(Y_f)\to \F_z$ l'est aussi.
\end{proof}
\end{prop}

\section{Opérations sur les faisceaux dans la représentation $(\sqcup)$}

Dans toute cette section, $X$ désigne une courbe intègre lisse sur $k$. Nous décrivons comment effectuer des opérations sur les faisceaux constructibles sur $X$. La représentation par recollement est celle utilisée par les algorithmes de calcul de cohomologie du chapitre \ref{chap:5}, c'est pourquoi nous nous efforçons de décrire toutes les opérations dans cette représentation, même si elles admettent une expression plus simple dans l'une des autres représentations (ce que nous mentionnons le cas échéant).

\subsection{Restriction à un ouvert plus petit}

Il sera utile par la suite, étant donné une représentation par recollement d'un faisceau constructible $\F$ par rapport à un ouvert de lissité $U$, de déterminer la représentation par recollement de ce même faisceau $\F$ par rapport à un ouvert plus petit que $U$.\\

Soient $U$ un ouvert de $X$, et $Z$ le fermé complémentaire.
\'{E}tant donné un faisceau constructible $\F$ lisse sur $U$, défini comme ci-dessus par le triplet $(\mathscr{L},\F_Z,\phi)$, et un autre ouvert $U'\subset U$ de complémentaire réduit $Z'$ dans $X$, comment calculer la représentation de $\F$ relativement au couple $(U',Z')$ ?

\[ U'\xrightarrow{\alpha} U\xrightarrow{j} X \xleftarrow{i'} Z'\xleftarrow{\beta} Z \]

Le morphisme donné est $\phi \colon \beta^\star i'^\star \F \to \beta^\star i'^\star j_\star j^\star \F$. Calculons $\psi \colon i'^\star \F \to i'^\star j_\star \alpha_\star \alpha^\star j^\star \F$, défini par ses fibres en les points $z\in Z'$. 

Pour les points $z\in Z$, ce morphisme est le morphisme déja donné. En effet, d'après le lemme \ref{adjiso}, le morphisme $j^\star\F\to \alpha_\star\alpha^\star j^\star\F$ est un isomorphisme. Par conséquent, le morphisme \[ i'^\star j_\star j^\star\F \to i'^\star j_\star\alpha_\star\alpha^\star j^\star\F \] est encore un isomorphisme, et le morphisme cherché \[i'^\star\F\to i'^\star j_\star\alpha_\star\alpha^\star j^\star\F \] s'identifie au morphisme \[i'^\star\F\to i'^\star j_\star j^\star \F.\] Ses fibres en chaque point de $Z$ sont celles du morphisme donné \[\phi \colon \beta^\star i'^\star\F\to \beta^\star i'^\star j_\star j^\star \F.\] 
Pour un point $w\in Z'-Z=U-U'$, la flèche cherchée \[\F_w\to (j_\star\alpha_\star\alpha^\star j^\star\F)_w\] s'identifie à la flèche $\F_w\to (j_\star j^\star\F)_w$, qui est elle-même un isomorphisme puisque $w$ est un point de $U$.
\\ 

\subsection{Morphismes et somme directe}

Soient $\F_1$ et $\F_2$ deux faisceaux constructibles sur $X$. Soient $U_1,U_2$ les ouverts de lissité de $\F_1$ et $\F_2$ donnés dans leur définition. Posons $U=U_1\cap U_2$ et $Z=X-U$. Notons $j\colon U\to X$ et $i\colon Z\to X$ les inclusions. Alors $U$ est un ouvert de lissité de $\F_1$ et $\F_2$, et nous savons déterminer les représentations $(\mathscr{L}_1,\F_{1,Z},\phi_1)$ et $(\mathscr{L}_2,\F_{2,Z},\phi_2)$ de $\F_1$ et $\F_2$ relativement à $(X,U,Z)$. Un morphisme $f \colon\F_1\to \F_2$ est alors déterminé par la donnée suivante : \begin{itemize}[label=$\bullet$]
\item un morphisme de faisceaux lisses $\alpha \colon \mathscr{L}_1\to \mathscr{L}_2$ ;
\item un morphisme $\beta \colon\F_{1,Z}\to \F_{2,Z}$ défini sur les fibres en les $z\in Z$ ;
\item un morphisme $\gamma\colon i^\star j_\star \mathscr{L}_1 \to i^\star j_\star \mathscr{L}_2$ vérifiant $\gamma\circ \phi_1=\phi_2\circ\beta$.
\end{itemize}
Pour qu'une telle donnée définisse correctement un morphisme, il faut et il suffit que $\gamma=i^\star j_\star\alpha$.

\paragraph{Noyau} Par exactitude à gauche de $i^\star$, $j^\star$ et $j_\star$, le noyau d'un tel morphisme est défini par le triplet $(\ker \alpha,\ker \beta, \phi|_{\ker\beta})$. \\

\paragraph{Conoyau} De même, l'exactitude à droite de $i^\star$ et $j^\star$ montre que \[\coker \alpha=j^\star\coker f \text{~~et~~}\coker\beta=i^\star\coker f.\]  Enfin, la flèche $j_\star j^\star \F_2\to j_\star j^\star \coker(f)$ se factorise par $\coker(j_\star j^\star f)$ ; la flèche composée \[\F_2\to j_\star j^\star \F_2\to \coker j_\star j^\star f \xrightarrow{u} j_\star j^\star\coker(f)\] passe au quotient en \[\coker f \to \coker (j_\star j^\star \F_1)\xrightarrow{u} j_\star j^\star\coker(f).\] 
On en déduit que la flèche de recollement $i^\star\coker(f)\to i^\star j_\star j^\star\coker (f)$ est la composée du morphisme $\coker(\beta)\to \coker(\gamma)$ (induit par la flèche de recollement $\phi_2$) suivi de $i^\star u\colon \coker \gamma\to i^\star j_\star j^\star\coker (f)$.
\\
Remarquons que la dernière flèche se décrit explicitement. Soit $z\in Z(k)$. Soient $F_1, F_2$ les fibres des faisceaux lisses $\mathscr{L}_1,\mathscr{L}_2$, et notons $f_\bareta\colon F_1\to F_2$ le morphisme induit par $f$. Notons $I_z$ le groupe d'inertie en $z$. Alors la flèche $u_z$ est simplement $F_2^{I_z}/f_\bareta(F_1^{I_z})\to (F_2/f_\bareta(F_1))^{I_z}$.\\

\paragraph{Somme directe} De même, la somme directe de $\F_1$ et $\F_2$ est simplement définie sur $(X,U,Z)$ par $(\mathscr{L}_1\oplus \mathscr{L}_2,\F_{1,Z}\oplus\F_{2,Z},\phi_1\oplus\phi_2)$.\\

\subsection{Produit tensoriel et Hom interne}

\begin{lem}\label{lem:jstartens} Soit $Y$ un schéma. Soit $j\colon U\to Y$ un morphisme étale. Alors pour tous faisceaux $\F,\F'$ de groupes abéliens sur $Y$, il y a des isomorphismes canoniques $j^\star\F\otimes j^\star\F'\xrightarrow{\sim} j^\star(\F\otimes \F')$ et $j^\star\underline{\Hom}(\F,\F')\xrightarrow{\sim}\underline{\Hom}(j^\star\F,j^\star\F')$.
\begin{proof} Le foncteur $j^\star$ étant simplement la restriction du site étale de $X$ à celui de $U$, l'assertion pour le $\Hom$ interne est claire, puisque les deux faisceaux en question ont pour sections sur $V\to U$ le groupe $\Hom(\F|_V,\F'|_V)$. Quant au produit tensoriel, les isomorphismes canoniques suivants, valables pour tout faisceau $\G$ de groupes abéliens sur $X$, permettent de conclure.
\begin{align*}
\Hom(j^\star\F\otimes j^\star\F',\G) &= \Hom(j^\star\F,\underline{\Hom}(j^\star \F',\G)) &\text{\cite[3.19]{milneEC}} \\ &= \Hom(\F,j_\star \underline{\Hom}(j^\star\F',\G)) \\ &= \Hom(\F,\underline{\Hom}(\F',j_\star \G))& \text{\cite[3.22.a)]{milneEC}} \\ &= \Hom(\F\otimes\F',j_\star\G) \\ &= \Hom(j^\star(\F\otimes\F'),\G).
\end{align*}
\end{proof}
\end{lem}

Soit $U$ un ouvert de $X$, de complémentaire réduit $Z$. Soient $\F$ et $\F'$ deux faisceaux constructibles de $\Lambda$-modules sur $X$ lisses sur $U$, définis par recollement par les données $(\mathscr{L},\F_Z,\phi)$ et $(\mathscr{L}',\F_Z',\phi')$. Notons $F,F'$ les fibres respectives des faisceaux lisses $\mathscr{L},\mathscr{L}'$.

\paragraph{Produit tensoriel} Comme le produit tensoriel commute au tiré en arrière, $j^\star(\F\otimes \F')=\mathscr{L}\otimes \mathscr{L}'$ ; le faisceau $\mathscr{L}\otimes\mathscr{L}'$ est encore un faisceau lisse sur $U$ \cite[093V]{stacks}. Notons $M=\HH^0(V,\mathscr{L})$ et $M'=\HH^0(V,\mathscr{L}')$. Soit $V\to U$ un revêtement qui trivialise $\mathscr{L}$ et $\mathscr{L}'$. Pour $z\in Z$, notons $I_z$ le sous-groupe des éléments de $\Aut(V|U)$ fixant un même antécédent de $z$ dans la compactification lisse $\bar V$ de $V$. La fibre de $\F\otimes\F'$ en un point $z$ de $Z$ est encore $\F_z\otimes \F_z'$. Il reste à décrire les fibres en les points de $Z$ du morphisme $\F\otimes\F'\to j_\star j^\star(\F\otimes\F')$, qui provient par adjonction des foncteurs $j^\star$ et $j_\star$ du morphisme identité de $\F\otimes\F'$. Notons $g\colon \F\otimes \F'\to j_\star j^\star\F\otimes j_\star j^\star \F'$ le morphisme déduit de $\id\to j_\star j^\star$. 
Le diagramme commutatif
\[
\begin{tikzcd}
j^\star(\F\otimes \F') \arrow[r,"\id"] \arrow[d,"j^\star g",swap,"\sim"'{anchor=north, rotate=90}] & j^\star(\F\otimes \F') \\
j^\star (j_\star j^\star\F\otimes j_\star j^\star\F') \arrow[ur,"(j^\star g)^{-1}",swap] & 
\end{tikzcd}
\]
produit par adjonction le diagramme commutatif
\[
\begin{tikzcd}
\F\otimes \F' \arrow[r] \arrow[d,"g",swap] & j_\star j^\star(\F\otimes \F') \\
j_\star j^\star\F\otimes j_\star j^\star\F' \arrow[ur] & 
\end{tikzcd}\]
qui montre que le morphisme de recollement $i^\star(\F\otimes\F')\to i^\star j_\star j^\star (\F\otimes\F')$ est la composée 
\[ \F_Z\otimes \F_Z' \to i^\star j_\star\mathscr{L}\otimes i^\star j_\star\mathscr{L}' = i^\star (j_\star\mathscr{L}\otimes j_\star\mathscr{L}')\to i^\star j_\star (\mathscr{L}\otimes \mathscr{L}').\]
Sa fibre en $z$ est la composée
\[ \F_z\otimes \F'_z \xrightarrow{\phi_z\otimes\phi_z'} M^{I_z}\otimes M'^{I_z} \to (M\otimes M')^{I_z}.\]

\paragraph{Hom interne} Comme vu dans le lemme \ref{lem:jstartens}, $j^\star\underline{\Hom}(\F,\F')$ est le faisceau lisse $\underline{\Hom}(\mathscr{L},\mathscr{L}')$ dont la fibre est $\Hom_\Lambda(F,F')$. 

\begin{prop} Soit $z$ un point de $Z$. Soit $V\to U$ un revêtement galoisien trivialisant $\F$ et $\F'$. Notons $I\triangleleft \Aut(V|U)$ le sous-groupe distingué engendré par les groupes d'inertie en tous les points de $\bar V$ au-dessus de $Z$. Alors \[\underline{\Hom}(\F,\F')_{z}=\{(\alpha,\beta)\in \Hom_{\Lambda}(F,F')^{I}\times \Hom_\Lambda(\F_{z},\F'_{z})\mid \alpha\phi_{\F}=\phi_{\F'}\beta\}.\]

\begin{proof} Soit $V\to U$ un revêtement galoisien trivialisant les faisceaux $\mathscr{L}$ et $\mathscr{L}'$. Soit $V_1$ le sous-revêtement de $V$ de groupe $\Aut(V_1|U)=G/I$. Notons $X'=X-(Z-\{ z\})$. Soient $V'$ et $V_1'$ les normalisées respectives de $X'$ dans $V$ et $V_1$. Le morphisme $V_1'\to X'$ est étale, puisque $G/I$ agit librement sur $(V_1')_z$ ; c'est le sous-revêtement non ramifié maximal de $V'\to X'$.
Par définition, 
\[ \underline{\Hom}(\F,\F')_{z}=\colim_{(Y,y)}\Hom(\F|_Y,\F'|_Y)\]
où les couples $(Y,y)$ sont des voisinages étales de $(X,z)$. La sous-catégorie des voisinages étales $(Y,y)$ vérifiant : \begin{itemize}[label=$\bullet$]
\item $Y$ est connexe
\item $Y_z$ est réduit à un point
\item $Y\to X$ se factorise par $V_1'$
\item $k(Y)/k(X)$ est galoisienne
\end{itemize}
est cofinale dans celle des voisinages étales de $(X,z)$. 
Considérons un tel $f\colon (Y,y)\to (X,z)$ et le diagramme cartésien suivant.
\[ \begin{tikzcd}
Y_U\arrow[dd]\arrow[r] & Y \arrow[d] & \arrow[l]\arrow[dd] y\\
 & V_1' \arrow[d] &\\
U \arrow[r,"j"] & X' & \arrow[l] z
\end{tikzcd}
\] 
Montrons que $\Hom(\F|_Y,\F'|_Y)$ est le groupe décrit dans l'énoncé. Par recollement, 
\[ \Hom(\F|_Y,\F'|_Y)=\{ (\alpha,\beta)\in \Hom(\F|_{Y_U},\F'|_{Y_U})\times \Hom_{\Lambda}(\F_z,\F'_z)\mid \alpha\phi_\F=\phi_{\F'}\beta\}.\]
Il suffit désormais de montrer que $\Hom(\F|_{Y_U},\F'|_{Y_U})$ est égal à $\Hom_{\Lambda}(F,F')^I$. C'est le groupe des sections sur $Y_U$ du faisceau lisse $\underline{\Hom}(\mathscr{L},\mathscr{L}')$ sur $U$ de fibre $\Hom(F,F')$. Par conséquent, il est canoniquement isomorphe à $\Hom(F,F')^{\pi_1(Y_U)}$. Or l'action de $\pi_1(Y_U)$ sur $\Hom(F,F')$ se factorise de la façon suivante :
\[
\begin{tikzcd}
\pi_1(Y_U)\arrow[dr] \arrow[r] & \pi_1(V_1)\arrow[d] \arrow[r] & \pi_1(U) \arrow[d]\arrow[dr] &\\ & I \arrow[r]& G\arrow[r]&\Aut_\Lambda(\Hom(F,F'))
\end{tikzcd}
\]
Notons $J$ l'image de $\pi_1(Y_U)$ dans $I$. Le morphisme $Y_U\to V_1$ se factorise par le sous-revêtement $W$ de $V\to V_1$ défini par $\Aut(V|W)=J$ ; par conséquent, $Y\to V_1'$ se factorise encore par la normalisation $W'$ de $V_1'$ dans $W$. Comme $k(Y)/k(X)$ est galoisienne, $Y\to X$ se factorise alors par la clôture galoisienne de $W'\to X$, qui ne saurait être étale en aucun point au-dessus de $z$, sauf si $W'=V_1'$, c'est-à-dire $J=I$. L'étalitude de $Y\to X$ entraîne donc la surjectivité de $\pi_1(Y_U)\to I$. Par conséquent, \[ \Hom_{\Lambda}(F,F')^{\pi_1(Y_U)}=\Hom_\Lambda(F,F')^I.\] 
\end{proof}
\end{prop}
Enfin, le morphisme de recollement $\underline{\Hom}(\F,\F')_z\to (j_\star \underline{\Hom}(\mathscr{L},\F'_U))_z$ se calcule en considérant les sections sur chaque $Y\to X$ ; c'est simplement la projection
\[ \underline{\Hom}(\F,\F')_z\subseteq \Hom_\Lambda(F,F')^I\times \Hom_\Lambda(\F_z,\F'_z) \to \Hom_\Lambda(F,F')^{I_z}\]
où $I_z$ est le groupe d'inertie d'un point de $V'$ au-dessus de $Z$.

\subsection{Tiré en arrière}

Soit $f\colon Y\to X$ un morphisme de courbes intègres lisses sur $k$. Comme vu en \ref{sec:opfaisclisses}, $f$ se factorise en $Y\xrightarrow{s} Y'\xrightarrow{\nu} X$, où $s$ est une immersion ouverte et $\nu$ est un morphisme fini localement libre. Il suffit donc de traiter séparément le cas où $f$ est une immersion ouverte, et le cas où $f$ est fini.\\

Soit $\F$ un faisceau constructible sur $X$, défini par un ouvert de lissité $j\colon U\to X$, le fermé réduit complémentaire $i\colon Z\to X$, le faisceau lisse $\mathscr{L}= j^\star\F$ sur $U$, le faisceau $\F_Z=i^\star\F$ sur $Z$ et le morphisme de recollement $\phi \colon i_\star\F_Z\to i_\star i^\star j_\star\mathscr{L}$.

\paragraph{Le cas d'une immersion ouverte} Supposons que $f\colon Y\to X$ soit une immersion ouverte. Notons $V=U\cap Y, W=Z\cap Y$. Alors $f^\star \F$ est clairement lisse sur $V$, les fibres de $f^\star\F$ en les points de $W$ sont égales à celles de $\F$ en leurs images dans $Z$, et le morphisme de recollement de $f^\star\F$ relativement à $V,W$ est simplement la restriction à $W$ de celui de $\F$.

\paragraph{Le cas d'un morphisme fini} Supposons que $f\colon Y\to X$ soit fini.  Considérons le diagramme suivant, dont les carrés sont cartésiens.

\[
\begin{tikzcd}
V \arrow[r,"j'"]\arrow[d,"p",swap] & Y \arrow[d,"f"] & \arrow[l,"i'",swap]\arrow[d,"q"] W \\
U \arrow[r,"j"] & X & \arrow[l,"i",swap] Z
\end{tikzcd}
\]

Le faisceau $j'^\star f^\star \F=p^\star j^\star \F$ est lisse sur $V$. De plus, la commutativité du carré de droite permet de calculer simplement les fibres de $f^\star\F$ en les points de $W$. Déterminons le morphisme d'adjonction $i'^\star f^\star\F\to i'^\star j'_\star j'^\star f^\star \F=i'^\star j'_\star p^\star \mathscr{L}$. La flèche $i^\star \F\to i^\star j_\star j^\star \F$ est connue ; en lui appliquant $q^\star$, on obtient $i'^\star f^\star \F\to i'^\star f^\star j_\star j^\star \F$. Reste à calculer la flèche de comparaison $f^\star j_\star j^\star \F \to j'_\star j'^\star f^\star \F=j'_\star p^\star j^\star \F$.
Notons que $f^\star j_\star j^\star \F$ est représenté par un schéma étale sur $Y$ qui se calcule explicitement par une normalisation puis un produit fibré. Il en est de même du morphisme de schémas qui représente la flèche d'adjonction $j^\star\F \to p_\star p^\star j^\star \F$, ainsi que celui représentant $f^\star j_\star j^\star\F \to f^\star j_\star p_\star p^\star j^\star \F$. Remarquons que le faisceau de droite est canoniquement isomorphe à $f^\star f_\star j'_\star p^\star j^\star \F$ ; il suffit maintenant de lui appliquer l'unité d'adjonction $f^\star f_\star\to \id$ (pour sa description, voir l'annexe \ref{sec:weilres}) pour obtenir la composée $f^\star j_\star j^\star \F \to j'_\star p^\star j^\star \F$. En résumé, il s'agit de calculer les fibres en les points de $W$ de la composée :
 
 \[ f^\star \F \to f^\star j_\star j^\star \F \to f^\star j_\star p_\star p^\star j^\star \F \xrightarrow{\sim} f^\star f_\star j'_\star p^\star j^\star \F \to j'_\star p^\star j^\star \F \xrightarrow{\sim} j'_\star j'^\star f^\star \F.
 \]
\begin{rk}
Les représentations $(!)$ et $(\star)$ sont conceptuellement bien mieux adaptées à cette tâche, car le foncteur $f^\star$, qui se calcule sur les schémas par un simple produit fibré, commute aux noyaux et conoyaux. L'avantage de la méthode décrite ci-dessus est de ne pas avoir recours au calcul possiblement coûteux de l'une de ces deux représentations. 
\end{rk}
 
\subsection{Poussé en avant}\label{subsec:pushforward}

Soit $f\colon Y\to X$ un morphisme de courbes intègres lisses sur $k$. Nous allons montrer comment calculer les foncteurs $\R^i f_\star$, $i\geqslant 0$. La notation $f_\star$ seule désignera toujours le foncteur non dérivé. Lorsque $f$ est fini, l'exactitude de $f_\star$ assure que $\R^if_\star=0$ dès que $i\geqslant 0$. 

\paragraph{Le cas d'une immersion ouverte} Soit $j\colon U\to X$ une immersion ouverte de courbes intègres lisses sur $k$. Soit $Z$ le fermé réduit complémentaire de $U$ dans $X$. Soit $\F$ un faisceau constructible de $\Lambda$-modules sur $U$, lisse sur un ouvert $V$ de $U$, de fibre générique géométrique $M$. Comme $j^\star j_\star \F\to \F$ est un isomorphisme, le faisceau $j_\star\F$ est encore lisse sur $V$. De plus, pour tout point géométrique $z$ de $Z$, la proposition \ref{prop:fibforward} assure que
\[ \begin{array}{rcl} (\R j_\star \F)_z&=& \RG(U\times_X X_{(\bar z)},\F) \\ 
&=& \RG(V\times_X X_{(\bar z)},\F) \\
&=& \RG(\eta_z,\F) \\
&=& \RG(I_z,M)  \end{array}\]
où $\eta_z$ est le point générique de l'anneau local strictement hensélien $X_{\bar z}=\Spec \OO_{X,z}^{sh}$, le point $\bareta_z$ est un point générique géométrique de $X_{(\bar z)}$ et $I=\Gal(\bareta_z|\eta_z)$. Le morphisme de recollement $j_\star\F\to j_\star j^\star j_\star\F$ est un isomorphisme. D'autre part, \[ j^\star \R^1j_\star=\R^1(j^\star j_\star)=0\]
puisque $j^\star j_\star$ est un isomorphisme. Le faisceau $\R^1 j_\star\F$ est donc supporté sur $Z$, et $j_\star j^\star \R^1j_\star\F=0$. Comme la cohomologie de $\RG(I,M)$ est concentrée en degrés 0 et 1, le foncteur $\R^ij_\star$ est nul pour tout $i\geqslant 2$. 

\paragraph{Le cas d'un morphisme fini} Soit $f\colon Y\to X$ un morphisme fini de courbes intègres lisses. Soit $\F$ un faisceau constructible de $\Lambda$-modules sur $Y$ défini par recollement relativement à $V\xrightarrow{j} Y\xleftarrow{i} W$, où $V$ est un ouvert de $Y$ tel que $\F|_V$ soit lisse, trivialisé par un revêtement étale $T\to V$. Notons $\phi \colon i^\star\F\to i^\star j_\star j^\star \F$ le morphisme de recollement. Décrivons par recollement le faisceau $f_\star \F$, qui est encore constructible d'après \cite[IX, Prop. 2.14.(i)]{sga43}. Quitte à remplacer $f$ par sa complétion projective lisse et les faisceaux en question par leur prolongement par zéro, nous pouvons supposer $X$ et $Y$ projectives.
Si la caractéristique de $k$ est nulle, le morphisme $f$ est génériquement étale. 
Si la caractéristique de $k$ est $p>0$, le morphisme $f$ se factorise en composée d'une succession de Frobenius relatifs suivi d'un morphisme génériquement étale \cite[0CD2]{stacks}. Il suffit donc de traiter séparément les cas du Frobenius et d'un morphisme génériquement étale.\\

Le morphisme de Frobenius relatif $\phi\colon Y\to Y^{(p)}$ étant un homéomorphisme universel, le théorème \ref{th:homeouniv} assure que le foncteur $\phi^\star$, dont $\phi_\star$ est l'adjoint à droite, est une équivalence de catégories. Voici comment la donnée de recollement de $\phi_\star\F$ se déduit de celle de $\F$.
Le faisceau $\phi_\star\F$ est lisse sur l'ouvert $V^{(p)}$ de $Y^{(p)}$. Son fermé complémentaire réduit $Z$ a autant de points fermés que $W$, et étant donné un point fermé $w$ de $W$, la fibre du faisceau $\phi_\star\F$ en $\phi(w)$ est $\F_w$ car $\phi$ est radiciel. La description de $\phi_\star\F$ sur $X$ est donc la suivante : il est lisse sur l'ouvert $V^{(p)}$, de même fibre générique géométrique que $\F$ et trivialisé par $T^{(p)}$. Ses fibres sur $Z$ sont les mêmes que sur $W$. En notant $j'$ l'inclusion de $U$ dans $X$, la finitude de $\phi$ assure que $j'_\star j'^\star\phi_\star\F=\phi_\star j'_\star j'^\star \F$ a les mêmes fibres que $j_\star j^\star \F$, et le morphisme de recollement de $\phi_\star\F$ est le même que celui de $\F$.\\

Supposons désormais $f$ génériquement étale.
Soit $Z$ un fermé réduit de $X$ contenant $f(W)$ ainsi que tous les points de $X$ au-dessus desquels $f$ n'est pas étale. Soit $U$ l'ouvert complémentaire. Considérons le diagramme suivant, dont les deux carrés sont cartésiens :

\[
\begin{tikzcd}
V' \arrow[r,"j'"]\arrow[d,"p",swap] & Y \arrow[d,"f"] & \arrow[l,"i'",swap]\arrow[d,"q"] W' \\
U \arrow[r,"j"] & X & \arrow[l,"i",swap] Z
\end{tikzcd}
\]

Comme $V'\subset V$, le faisceau $j'^\star\F$ est lisse ; le morphisme $p$ étant fini étale, $p_\star j'^\star\F$ est lisse sur $U$. Il est représenté par le $U$-schéma fini étale $R_{V'\to U}(j'^\star \F)$. 
Par finitude de $f$, les morphismes de changement de base des deux carrés de ce diagramme sont des isomorphismes \cite[Corollary 5.3.9]{fulei}.
Par conséquent, le faisceau $i^\star f_\star \F$ est isomorphe à $q_\star i'^\star \F$, qui se calcule explicitement comme un ensemble de $\Lambda$-modules indexé par les points de $Z$. Le morphisme $\phi \colon i'^\star \F\to i'^\star j'_\star j'^\star \F$ est connu. 
Alors comme $j_\star j^\star f_\star \F=j_\star p_\star j'^\star \F=f_\star j'_\star j'^\star \F$, la flèche $q_\star \phi$ est un morphisme $i^\star f_\star \F \to i^\star j_\star j^\star f_\star \F$. Il s'agit du morphisme d'adjonction $f_\star \F\to j_\star f^\star f_\star \F$, défini pour tout $X$-schéma étale $T$ par le morphisme $\F(T\times_X Y)\to \F(T\times_X V')$ obtenu par fonctorialité à partir de l'inclusion $T\times_X V'\to T\times_X Y$.

\begin{rk}
Au moins pour le cas d'un morphisme fini, les représentations $(\star)$ et $(!)$ se prêtent conceptuellement mieux au calcul de $f_\star$, qui s'effectue alors par restriction de Weil.
\end{rk}

\section{Un exemple détaillé}

Considérons sous une perspective différente l'exemple de la mise au carré sur la droite affine. Dans cet exemple, $\Lambda=\ZZ/2\ZZ$ et $k=\bar \QQ$. Considérons $X=\A^1_k=\Spec k[x]$ et $Y=\Spec k[x,y]/(x-y^2)$. Notons $A=k[x], B=k[x,y]/(x-y^2)$. Soit $f\colon Y\to X$ le morphisme défini par $x\mapsto x$. Remarquons que $B=A\cdot 1\oplus A\cdot y$ : le morphisme $f$ est fini et libre. Il est étale au-dessus de $U=\GG_m=\Spec k[x^{\pm 1}]$, et ramifié au-dessus de 0 où la fibre est $k[y]/(y^2)$. Notons $U'=Y|_{\GG_m}$.\\
Considérons le faisceau constant $\F=\Lambda\simeq \mu_2$ sur $Y$, représenté par le $Y$-schéma $F\coloneqq \Spec B[t]/(t^2-1)$, et calculons le faisceau constructible $f_\star\F$.

\begin{figure}[H]\centering \includegraphics[scale=0.3]{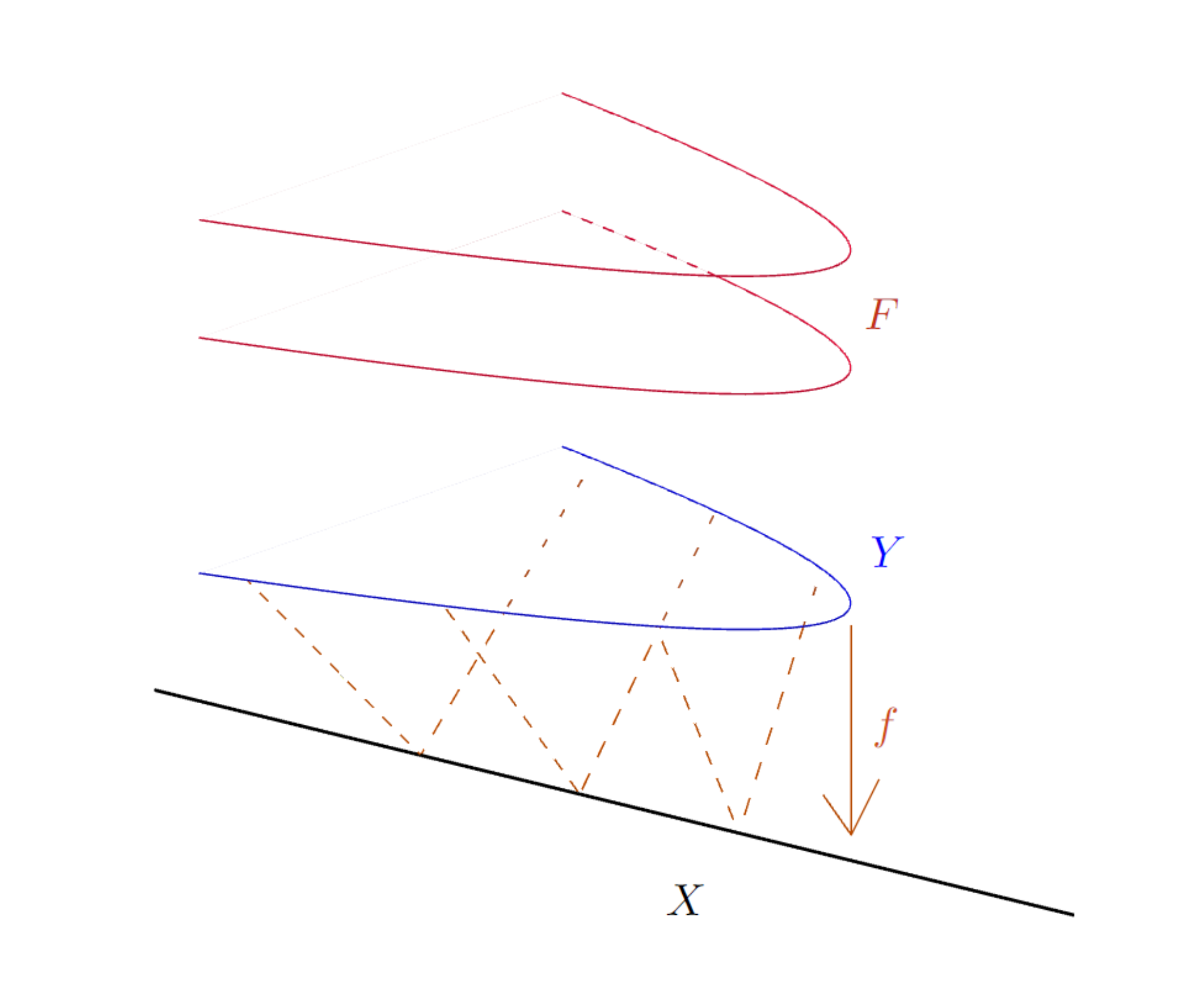}\includegraphics[scale=0.3]{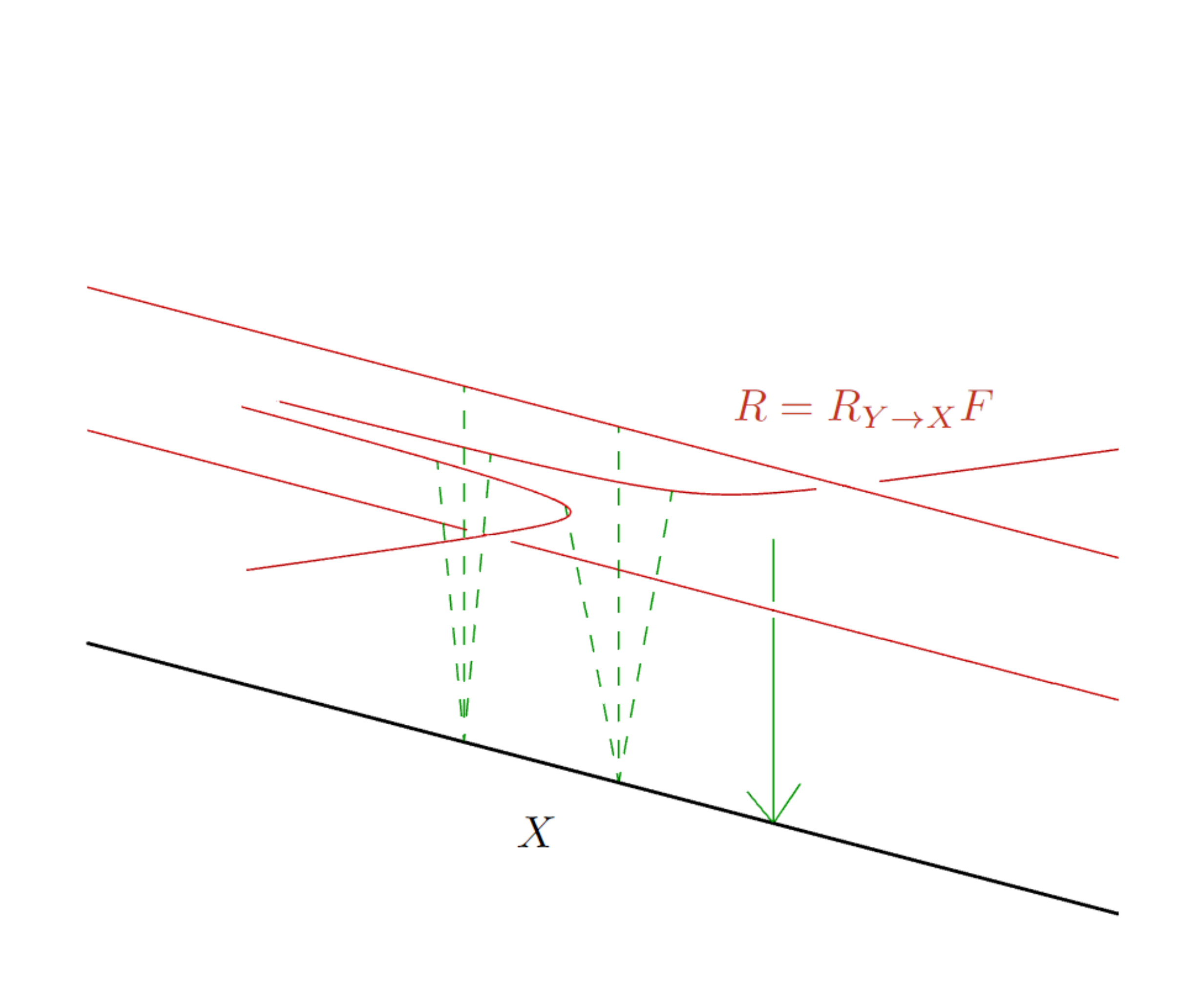}
\caption{Le $Y$-schéma $F$ et sa restriction de Weil}
\end{figure}

\paragraph{Calcul de $f_\star \F$} Le faisceau $f_\star \F$ est représenté par la restriction de Weil $R\coloneqq R_{Y\to X}F$. En notant $\phi=t^2-1\in A[t]$,  \[\phi(\alpha \cdot 1+\beta\cdot y)=(\alpha ^2+\beta^2x-1)\cdot 1+ 2\alpha\beta \cdot y.\] Par conséquent, \[ R=\Spec A[\alpha,\beta]/(\alpha ^2+\beta^2x-1,\alpha\beta).\]
Il est prévisible que $R\to X$ soit étale, mais pas fini car $f$ est ramifié. L'étalitude se vérifie rapidement : le déterminant jacobien de la présentation ci-dessus est égal à $\alpha^2-x\beta ^2$, qui vaut encore $2\alpha^2-1=(1-2x\beta^2)^{-1}$ dans $\HH^0(R,\OO_R)$. Il n'est effectivement pas fini puisque $\beta$ n'est pas racine d'un polynôme unitaire à coefficients dans $k[x]$. Cependant, $R$ doit être fini au-dessus de $U$ (où $f$ est étale), ce qui est le cas puisque $\alpha^3=\alpha$ et $\beta^3=x^{-1}\beta$. Calculons la loi de groupe sur $R$. La loi de groupe sur $F$ est donnée par \[ \begin{array}{rcl} B[t]/(t^2-1)&\longrightarrow& B[a,b]/(a^2-1,b^2-1)\\ t&\longmapsto& ab. \end{array}\] En écrivant $t=\alpha+y\beta$, $a=a_1+ya_2$, $b=b_1+yb_2$, la loi de groupe sur $R$ est alors donnée par \[  \begin{array}{rcl}
A[\alpha,\beta]/(\alpha^2+x\beta^2-1,\alpha\beta)&\to& A[a_1,b_1,a_2,b_2]/(a_1^2+xa_2^2-1,a_1a_2,b_1^2+xb_2^2-1,b_1b_2) \\ \alpha&\mapsto& a_1a_2+xb_1b_2 \\ \beta&\mapsto& a_1b_2+b_1a_2.\end{array}\]

\paragraph{Calcul de $(f_\star \F)_0$} La fibre de $R$ en $x=0$ est $k[\alpha,\beta]/(\alpha^2-1,\beta)$ et contient deux points. La fibre générique géométrique de $R\to X$ est $\Spec \overline{k(x)}[\alpha,\beta]/(\alpha^2+x\beta^2-1,\alpha\beta)=\{ (0,\pm\sqrt{x^{-1}}),(\pm 1,0)\}$. La loi de groupe sur cette fibre se déduit des formules déterminées ci-dessus ; par exemple, $(0,\sqrt{x^{-1}})\cdot (-1,0)=(0\cdot (-1)+x\sqrt{x^{-1}}\cdot 0,0\cdot 0 +\sqrt{x^{-1}}\cdot(-1))=(0,-\sqrt{x^{-1}})$. C'est un groupe isomorphe à $\Lambda^2$, d'élément neutre $(1,0)$.\\

\paragraph{Calcul de $j_\star j^\star f_\star\F$} Notons $j\colon U\to X$. Le faisceau $j_\star j^\star R$ est représenté par le lieu étale sur $\A^1$ du normalisé de $\A^1$ dans $R'\coloneqq R|_U=\Spec k[x^{\pm 1},\alpha,\beta]/(\alpha^2+\beta^2x-1,\alpha\beta)$. Il est temps de s'intéresser à la structure de $R'$. Le morphisme $R'\to U$ admet au moins deux sections (les sections globales de $f_\star\F$, correspondant à $\beta=0$ et $\alpha=\pm 1$), dont les images sont des composantes connexes de $R'$ isomorphes à $U$. Un calcul rapide donne $R'=U\sqcup U \sqcup \Spec k[x^{\pm 1},\beta]/(x\beta^2-1)$. D'une part, la normalisation de $\A^1$ dans $\GG_m$ est $\A^1$. Calculons d'autre part la clôture intégrale de $k[x]$ dans $S=k[x^{\pm 1},\beta]/(x\beta^2-1)$. Le couple $(1,x\beta)$ est une $k(x)$-base de $\Frac S$. Le sous-anneau $k[x,x\beta]\simeq k[x,s]/(s^2-x)$ de $S$ est entier sur $k[x]$, et son corps des fractions est $\Frac S$ ; comme il est normal, c'est la clôture intégrale de $k[x]$ dans $S$. Par conséquent, $j_\star R'$ est le lieu étale de $\A^1\sqcup \A^1\sqcup \Spec k[x,s]/(s^2-x)\to \A^1$, c'est-à-dire $\A^1\sqcup \A^1\sqcup \Spec k[x^{\pm 1},s]/(s^2-x)$. La fibre en $x=0$ de $j_\star R'$ est donc $\Spec(k[x]\times k[x])$ : elle contient deux $k$-points. La flèche de recollement $R_0\to(j_\star R')_0$ est, comme précisé plus haut, l'identité. 

\paragraph{Remarques sur $R'$ et $(j_\star R')_0$} Nous avons vu que le faisceau lisse $j^\star f_\star \Lambda$ sur $U$ était représenté par $R'=U\sqcup U \sqcup V$ où $V= \Spec k[x^{\pm 1},\beta]/(x\beta^2-1)$. Le morphisme $V\to U$ est fini étale, et galoisien car de degré 2. Par conséquent, $V\times_ U V=V\sqcup V$, et $R'\times_U V=\sqcup^4 V$. Le revêtement $V\to U$ trivialise donc $j^\star f_\star \Lambda$, et le groupe $G=\Aut(V|U)$ est isomorphe à $\ZZ/2\ZZ$. Le morphisme non trivial dans $G$ est $\sigma \colon\beta\mapsto -\beta$. Considérons les complétions projectives lisses $\bar U=\PP^1=\Proj k[X,Y]$ et $\bar V=\Proj k[X,B,Y]/(XB^2-Y)$ de $U$ et $V$. Au-dessus du point $0=(0:1)$ de $\PP^1$ se trouve uniquement $(0:1:0)\in \bar V$, qui est donc invariant sous $\sigma$. Par conséquent, l'action sur $(j^\star f_\star \Lambda)(V)=\Hom_U(V,U\sqcup U\sqcup V)=\{ V\to U^{(1)},V\to U^{(2)},\id_V,\sigma \}$ du groupe d'inertie $I_0$ est l'action (par précomposition) de $G=\{ \id_V,\sigma \}$. L'élément $\sigma\in G$ agit trivialement sur les deux morphismes $V\to U$, mais échange $\id_V$ et $\sigma$. Par conséquent, $R'(V)^{I_0}=\{ V\to U^{(1)},V\to U^{(2)} \}\simeq \ZZ/2\ZZ$ est bien isomorphe au groupe $(j_\star R')_0$ trouvé précédemment.

\section{Faisceaux constructibles sur les courbes nodales}\label{sec:constrnod}

Soit $X$ une courbe intègre nodale sur $k$ de corps des fonctions $K$. Soient $j\colon U\to X$ et $i\colon Z\to X$ les immersions d'un ouvert et d'un fermé réduit complémentaires. Afin de décrire par recollement un faisceau constructible sur $X$ lisse sur $U$, il est nécessaire de savoir décrire le faisceau $j_\star j^\star\F$. Notons $\mathscr{L}$ le faisceau $j^\star\F$, et $F$ sa fibre générique. Soit $\bar z$ un point géométrique de $Z$. Notons $K_{\bar z}$ le corps des fractions de $\OO_{X,\bar z}$, et $I_{\bar z}$ le groupe d'inertie associé à un choix de plongement $K^\sep\to  K_{\bar z}^\sep$. Notons $V$ l'intersection de $U$ avec le lieu non singulier de $X$. Calculons $(j_\star\mathscr{L})_{\bar z}$.

\paragraph{Si $\bar z$ est non singulier} 
Rappelons que $(j_\star\mathscr{L})_{\bar z}=\HH^0(U\times_{X}X_{\bar z},\mathscr{L})$. Or $U\times_X X_{\bar z}=V\times_X X_{\bar z}$, qui est réduit au point générique de $X_{\bar z}$, dont le groupe fondamental est $\Gal(K_{\bar z}^\sep|K^\sep)=I_{\bar z}$. Par conséquent, \[(j_\star\mathscr{L})_{\bar z}=\HH^0(I_{\bar z},F).\]

\paragraph{Si $\bar z$ est singulier}  Notons $\tX$ et $\widetilde{X_{\bar z}}$ les normalisations respectives de $X$ et $X_{\bar z}=\Spec \OO_{x,\bar z}$. Le schéma $X_{\bar z}$ est constitué de trois points : un point fermé, et deux idéaux premiers minimaux $p,q$ correspondant aux branches de $X$ en $\bar z$ \cite[06DT]{stacks}.
D'après \cite[0CBM]{stacks}, $\widetilde{X_{\bar z}}=\tX\times_X X_{\bar z}$. Par conséquent, en notant $P,Q$ les antécédents de $\bar z$ dans $\tX$,  \[ \widetilde{X_{\bar z}}=\tX_{P}\sqcup \tX_{Q}.\]
Les schémas $\tX_{P}$ et $\tX_{Q}$ sont des traits dont les points fermés ont pour image celui de $X_{\bar z}$, et les points génériques ont pour images les points $p,q$ de $X_{\bar z}$. Notons $\eta_P,\eta_Q$ les points génériques respectifs de $\tX_{P}, \tX_{Q}$.
En particulier, il y a des isomorphismes de schémas \[ U\times_X X_{\bar z}=\{p,q\} \]
et \[ U\times_X X_{\bar z}\times_{X_{\bar z}}\widetilde{X_{\bar z}}=\eta_P\sqcup \eta_Q \xrightarrow{\sim} \{ p,q\}.\]
Par conséquent, \begin{eqnarray*} (j_\star\mathscr{L})_{\bar z}&=& \HH^0(\{ p,q\},\mathscr{L}) \\ &=&\HH^0(\eta_P,\mathscr{L})\times \HH^0(\eta_Q,\mathscr{L}) \\ &=& \HH^0(I_{P},F)\times \HH^0(I_{Q},F). \end{eqnarray*}

\paragraph{Représentation $(\sqcup)$ pour les courbes nodales} Voici comment représenter par recollement un faisceau constructible $\F$ sur une courbe nodale $X$ sur $k$. La courbe $X$ est représentée par sa normalisée $\tX$ avec des couples de points marqués. L'ouvert de lissité $U$ -- dont nous supposons, quitte à ce qu'il ne soit pas maximal, qu'il ne contient pas les points nodaux -- est défini par sa normalisée $j\colon \tilde U\to \tX$. De même, le fermé $Z$ est donné par $\tilde{Z}\coloneqq Z\times_X\tX\to \tX$.
Le faisceau $\F$ est défini par : \begin{enumerate}
\item un faisceau lisse $\mathscr{L}$ sur $\tilde U$ de fibre $F$, représenté comme dans la section \ref{subsec:replisnod} à l'aide d'un revêtement trivialisant $\tY\to \tX$;
\item des fibres $M_z$ en les points de $\tilde Z$ -- une seule fibre $M_{(x,y)}$ par couple $(x,y)$ de points marqués ;
\item des morphismes $M_z\to F^{I_z}$ pour $z$ non marqué, et $M_{(x,y)}\to F^{I_x}\times F^{I_y}$ pour tout couple de points marqués $(x,y)$, où $I_z,I_x,I_y$ sont des sous-groupes du groupe de $X$-automorphismes d'une composante connexe de $\tY$, déterminés par le choix d'antécédents quelconques de $z,x,y$ dans une même composante connexe de $\tY$.
\end{enumerate}

\section{Constructions sur les surfaces}

Voyons brièvement comment représenter les faisceaux constructibles en dimension supérieure.  Afin de simplifier l'exposition, considérons le cas des surfaces, où apparaissent déjà quelques complications. Soit $X$ une surface lisse sur $k$. Soit $\F$ un faisceau constructible sur $X$. Il existe un ouvert $j\colon U\to X$ sur lequel $\F$ est lisse ; soit $i\colon Z\to X$ le fermé réduit complémentaire.

\paragraph{La représentation $(\sqcup)$} Les faisceaux $i^\star\F$ et $i^\star j_\star j^\star \F$ sont constructibles sur $Z$, qui est une courbe : nous savons représenter ces faisceaux ainsi que le morphisme de recollement entre eux. Cet argument s'adapte récursivement au cas où $X$ est de dimension supérieure.

\paragraph{La représentation $(!)$} La représentation comme conoyau s'adapte immédiatement en dimension quelconque : nous avons décrit la représentabilité de $f_!\Lambda$ pour tout morphisme étale $f$ entre $k$-variétés.

\paragraph{La représentation $(\star)$} Il s'agit ici de décrire des morphismes de faisceaux $\alpha\colon p_\star M\to p'_\star N$, où $p\colon Y\to X$ et $p'\colon Y'\to X$ sont des morphismes finis et $M,N$ sont des $\Lambda$-modules de type fini. La difficulté est le cas (inexistant si $X$ est une courbe) où $\dim Y=\dim Y'=1$ : le morphisme $\alpha$ peut être non nul, par exemple si $Y=Y'$, et les faisceaux $p_\star\Lambda$ et $q_\star\Lambda$ ne sont pas représentables par des schémas sur $X$.

\cleartooddpage

\chapter{Calculabilité de la cohomologie et algorithmes existants}\label{chap:4}

Soit $X$ un schéma de type fini sur un corps algébriquement clos $k$. Soit $n$ un entier inversible dans $k$. Les groupes $\HH^i(X,\ZZ/n\ZZ)$ sont des $\ZZ/n\ZZ$-modules de type fini, dont la question du calcul effectif a été extensivement étudiée. D'abord prouvée en caractéristique nulle par Poonen, Testa et van Luijk \cite[Th. 7.9]{poonen_testa}, la calculabilité de ces groupes a été démontrée en 2015 par Madore et Orgogozo \cite[Th. 0.1]{mo}. L'algorithme décrit dans leur article est résumé dans la section \ref{sec:mo} ; nous y montrerons que sa complexité est primitivement récursive. Les seuls algorithmes pour lesquels des bornes de complexité sont connues concernent le cas particulier des courbes projectives lisses. Dans ce cas, seul le groupe $\HH^1(X,\ZZ/n\ZZ)$, qui s'identifie d'après le théorème \ref{th:cohproj} au groupe de $n$-torsion de la jacobienne de $X$, présente un intérêt. La meilleure complexité pour le calcul de ce groupe est atteinte par l'algorithme probabiliste de Couveignes \cite{couveignes_linearizing}, décrit dans la section \ref{sec:couv}. Cependant, cet algorithme se cantonne au cas des courbes définies sur un corps fini, et nécessite d'avoir calculé au préalable la fonction zêta de la courbe en question. L'algorithme de Huang et Ierardi \cite{huang_counting}, pensé pour le comptage de points des courbes sur les corps finis, s'adapte quant à lui à tous les corps calculables munis d'un algorithme de factorisation. Une légère modification, détaillée dans la section \ref{sec:huang}, permet même de l'adapter au calcul de la division par $n$ dans $\Pic^0(X)$. Enfin, plus récemment, Jin a proposé une méthode de calcul de la cohomologie des faisceaux lisses sur les courbes lisses, ainsi que des bornes de complexité explicites \cite{jinbi_jin}. Son algorithme est expliqué dans la section \ref{sec:jin}.

\section{Calculabilité : l'algorithme de Madore et Orgogozo}\label{sec:mo}

\subsection{Représentation des objets et résumé de l'algorithme}

Soient $k$ un corps algébriquement clos et $\ell$ un nombre premier inversible dans $k$. Notons $\Lambda$ l'anneau $\ZZ/\ell\ZZ$. Les schémas de type fini sur $k$ sont représentés par recollement de schémas affines comme décrit dans l'annexe \ref{sec:repsch}. 

\begin{theorem}\cite[Th. 0.1]{mo} Il existe un algorithme calculant les groupes de cohomologie $\HH^i(X,\Lambda)$ d'un schéma $X$ de type fini sur $k$, ainsi que l'application $\HH^i(X,\Lambda)\to \HH^i(Y,\Lambda)$ déduite par fonctorialité d'un morphisme $Y\to X$. 
\end{theorem}

Soit $X$ un schéma de type fini sur $k$. Notons $d$ la dimension de $X$. L'algorithme procède de la façon suivante. Il commence par calculer les $r$ premiers étages d'un hyperrecouvrement $X_\bullet\to X$ dont les composantes sont des $K(\pi,1)$ pro-$\ell$ et tel que la cohomologie de $X$ soit celle du topos total $\Tot X_\bullet$ associé à l'hyperrecouvrement  (voir début de la section \ref{subsec:toposapp}). Dès que $r>d$, le morphisme canonique \[ \RG(\Tot X_{\bullet\leqslant r},\Lambda)\to \RG(\Tot X_{\bullet},\Lambda)\]
est un quasi-isomorphisme.
La cohomologie de $\Tot X_{\bullet \leqslant r}$ est calculée à l'aide d'approximations $X_{\bullet\leqslant r}^{(\lambda)}$ de $X_{\bullet\leqslant r}$, constituées de $\ell$-revêtements galoisiens des composantes connexes de $X_{\bullet\leqslant r}$, vérifiant
\[ \HH^i(\Tot X_{\bullet \leqslant r},\Lambda)=\colim_\lambda \HH^i(\Tot X_{\bullet\leqslant r}^{(\lambda)},\Lambda).\]
La cohomologie de $\Tot X_{\bullet\leqslant r}^{(\lambda)}$ est déterminée explicitement à partir d'une résolution de Godement tronquée du faisceau $\Lambda$. Il reste à déterminer des entiers $\alpha,\beta$ tels que \[ \HH^i(X,\Lambda)=\im (\HH^i(\Tot X_{\bullet\leqslant r}^{(\alpha)},\Lambda)\to \HH^i(\Tot X_{\bullet\leqslant r}^{(\beta)},\Lambda)).\]
La cohomologie de $\Tot X_{\bullet\leqslant r}^{(\lambda)}$ est également l'aboutissement d'une suite spectrale dont les coefficients de la première page sont les $\HH^j(X_i^{(\lambda)},\Lambda)$.
Chaque composante connexe de $X_i$ est un $K(\pi,1)$ pro-$\ell$, et les $X_i^{(\lambda)}$ correspondent à des quotients $\pi^{(\lambda)}$ du pro-$\ell$ groupe fondamental $\pi$ de cette composante ; il y a pour tout entier naturel $j$ un isomorphisme canonique \[ \HH^j(\pi,\Lambda)=\colim_\lambda \HH^j(\pi^{(\lambda)},\Lambda).\]
Un résultat sur les systèmes inductifs de suites spectrales permet de déterminer des entiers $a,b$ tels que \[ \HH^j(\pi,\Lambda)=\im (\HH^j(\pi^{(a)},\Lambda)\to \HH^j(\pi^{(b)},\Lambda))\]
et par conséquent les entiers $\alpha,\beta$ cherchés.

\subsection{Fibrations en courbes élémentaires}\label{subsec:fibel}

\begin{df} Soit $f\colon X\to S$ un morphisme de schémas. Il est appelé courbe élémentaire sur $S$ s'il peut être plongé dans un diagramme commutatif \[
\begin{tikzcd}
X\arrow["f",rd,swap]\arrow[r,"j"] & Y\arrow[d,"\bar{f}"] & D \arrow[l,"i",swap]\arrow[ld,"g"] \\
& S & 
\end{tikzcd}
\]
où $X=Y- D$, le morphisme $\bar{f}$ est une courbe relative projective lisse à fibres géométriquement connexes, et $g$ est un revêtement étale à fibres non vides. Une courbe $\ell$-élémentaire est une courbe élémentaire $f\colon X\to S$ telle que le faisceau $\R^1f_\star\ZZ/\ell\ZZ$ soit $\ell$-monodromique (voir définition \ref{df:lm}).
Un morphisme de schémas est appelé polycourbe ($\ell$-)élémentaire s'il admet une factorisation en courbes ($\ell$-)élémentaires.
\end{df}

Le résultat suivant, dû à M. Artin, affirme que tout schéma lisse sur un corps algébriquement clos est, localement pour la topologie de Zariski, une polycourbe élémentaire. 

\begin{prop}\cite[XI, Prop. 3.3]{sga43} Soit $X$ un schéma lisse sur un corps algébriquement clos $k$. Soit $x\in X$ un $k$-point. Il existe un ouvert $U$ de $X$ contenant $x$ tel que $U\to \Spec k$ soit une polycourbe élémentaire.
\end{prop}

Afin d'appliquer au schéma $X$ la théorie des pro-$\ell$-groupes, il suffit de montrer que l'on peut supposer, après changement de base étale, que l'ouvert $U$ de la proposition est une polycourbe $\ell$-élémentaire. Ceci se montre de la façon suivante : si $U\to k$ se factorise par une courbe élémentaire $g\colon U\to V$, et si $V'$ désigne un revêtement étale de $V$ trivialisant le faisceau $\R^1g_\star \Lambda$ (qui est lisse par \cite[XIII, Cor. 2.9]{sga1}), le changement de base $U\times_VV'\to V'$ est une courbe $\ell$-élémentaire.  Une récurrence sur la dimension permet de conclure.
Voici comment construire $V'$ : soit $\eta$ le point générique de $V$. Le $\Gal(\bareta|\eta)$-module $\HH^1(U_\bareta,\Lambda)$ se calcule comme un groupe de cohomologie d'une courbe (par exemple avec l'un des algorithmes des sections \ref{sec:huang} et \ref{sec:jin}) ; si $\eta'\to\eta$ est une extension finie séparable par laquelle se factorise l'action de $\Gal(\bareta|\eta)$, un ouvert étale sur $V$ de la normalisation de $V$ dans $\eta'$ convient pour $V'$.\\

Enfin, les polycourbes $\ell$-élémentaires sur $k$ sont des $K(\pi,1)$ pro-$\ell$. En effet, si $f\colon Y\to X$ est une courbe $\ell$-élémentaire et $\bareta$ est un point générique géométrique de $X$, il y a une suite exacte
\[ 1\to \prol(Y_\bareta)\to \prol(Y)\to \prol (X)\to 1.\]
Par récurrence sur la dimension relative de $f$, partant du cas des courbes affines lisses sur $k$ déjà vu dans la proposition \ref{prop:Kpi1}, on peut supposer que $X\to\Spec k$ et $Y_\bareta\to\bareta$ sont des $K(\pi,1)$ pro-$\ell$. Ainsi, pour tout faisceau $\ell$-monodromique $\F$ sur $Y$,
\[ \RG(Y,\F)=\RG(X,\R f_\star \F)=\RG(\prol(X),\RG(\prol(X_\bareta,\F_\bareta))=\RG(\prol(Y),\F_\bareta).\]

\paragraph{Complexité de la construction lorsque $X$ est lisse} Soit $x\in X$. Montrons que la construction d'un voisinage étale de $(X,x)$ qui est une polycourbe $\ell$-élémentaire sur $\Spec k$ est primitivement récursive. D'après les explications ci-avant, il suffit de construire un voisinage de Zariski de $(X,x)$ qui est une polycourbe élémentaire. Cette construction est décrite explicitement dans \cite[XI, Prop. 3.3]{sga43} : supposons $X$ plongé dans un espace affine $\A^r$. La construction commence par considérer un plongement projectif de la normalisation $\bar X$ de l'adhérence de $X$ dans $\PP^r$. Le point clé est la construction d'hyperplans en position générale $H_1,\dots,H_{d-1}$, où $d=\dim X$, qui coupent $\bar X$ et $Y$ transversalement et contiennent $x$. Ceci est décrit dans l'annexe \ref{subsec:hyplis}. Le reste de la construction se résume à un éclatement et à la détermination d'un ouvert de lissité par critère jacobien. Comme toutes ces opérations sont primitivement récursives, étant donné $x\in X(k)$, la construction d'un voisinage de $X$ qui est une polycourbe $\ell$-élémentaire est primitivement récursive. L'annexe \ref{subsec:recprimrec} assure alors qu'il existe un algorithme primitivement récursif qui recouvre $X$ par des polycourbes $\ell$-élémentaires.

\subsection{Calcul explicite de la filtration de Frattini itérée}

Soit $X$ une courbe intègre lisse sur $k$ de genre non nul. Notons $\pi$ le complété pro-$\ell$ du groupe fondamental de $X$. D'après \cite[1.20]{dixon}, le sous-groupe de Frattini $\Phi(\pi)$ de $\pi$ est alors le groupe $\pi^\ell [\pi,\pi]$, et le quotient $\pi/\Phi(\pi)$ est canoniquement isomorphe à $\HH^1(\pi,\FF_\ell)^\vee$ comme vu à la proposition \ref{prop:sgcar}. 

Définissons comme dans \cite[3.3]{mo}, la filtration de Frattini descendante par $\pi_1=\pi$, puis $\pi^{[\lambda+1]}=(\pi^{[\lambda]})^\ell[\pi^{[\lambda]},\pi^{[\lambda]}]$, et considérons les quotients $\pi^{(\lambda)}=\pi/\pi^{[\lambda]}$. En particulier, $\pi^{(\lambda)}$ est un quotient de $\pi^{(\lambda+1)}$. 
D'après \cite[1.116.(iii)]{dixon}, les $\pi^{[\lambda]}$ forment une base de voisinages de 1 dans $\pi$, ce qui entraîne que le morphisme $\pi\to\lim_{\lambda\geqslant 2}\pi^{(\lambda)}$ est un isomorphisme. Soit $X^{[\lambda]}$ un revêtement de $X$ correspondant à $\pi^{[\lambda]}$. \\

Nous avons vu dans la section \ref{subsubsec:X2} comment calculer $X^{[2]}$, alors noté $X_2$.
\'{E}valuons la complexité de cette construction. Afin de simplifier la présentation, supposons d'abord $X$ projective ; notons $g$ son genre. Nous avons vu dans la section \ref{subsubsec:ram} que \[ g(X^{[2]})=1+\ell^{2g}(g-1) \sim_{\ell,g} g\ell^{2g}.\]
Le genre de $X^{[\lambda]}$ est alors 
\[ g(X^{[\lambda]})=O_{g,\ell}\left(g\ell^{2g+2g\ell^{2g}+\dots +2g\ell^{2g+2g\ell^{2g}+\dots +2g\ell^{\dots\iddots^{2g\ell^{2g}}}}}\right)\]
où le nombre de termes horizontaux et de termes verticaux est $\lambda$. D'autre part, si $X=\bar X-\{ P_1,\dots,P_r\}$, la courbe $\bar X_2$ a $\ell^{r-2}$ points au-dessus de chaque $P_i$, et \[ |\HH^1(X^{(\lambda)},\Lambda)|\geqslant \ell^{\ell^{\iddots^{\ell^{r-2}}}}.\]
Ainsi, dans tous les cas sauf $\ell=2$ et $X=\GG_m$, la complexité du calcul de $X^{(\lambda)}$ est une exponentielle à $\lambda$ étages : ce n'est pas une fonction élémentaire en $(g_X,r,\lambda)$.

\begin{rk} La filtration considérée dans \cite{mo} est plus fine que la filtration de Frattini ; elle consiste à itérer, pour un sous-groupe $H$ de $\pi$, $H\mapsto \HH^\ell[H,\pi]$, et non $\HH^\ell[H,H]$ comme dans la filtration de Frattini. Cependant, une partie des bornes explicites obtenues sur les constructions s'exprime uniquement en termes de la filtration de Frattini.
\end{rk}

\subsection{Topos $\ell$-étale $\lambda$-approché d'un schéma simplicial}\label{subsec:toposapp}

Soit $X$ un schéma. Soit $\lambda$ un entier naturel non nul. 
Le topos $\ell$-étale $\lambda$-approché $X^{(\lambda)}$ de $X$ est le topos des faisceaux sur $X$ trivialisés par le revêtement $X^{[\lambda]}\to X$. C'est le topos associé au site $X_\lambda$ dont les objets sont les quotients (finis étales) de $X^{[\lambda]}\to X$. Un faisceau constructible sur $X^{(\lambda)}$ n'est rien d'autre qu'un $\pi^{(\lambda)}$-module. 

De même, il est possible d'associer à un schéma simplicial $X_\bullet$ un topos $X_\bullet^{(\lambda)}$. Le topos total associé est noté $\Tot X_\bullet^{(\lambda)}$ ; un objet de ce topos est la donnée, pour tout $i$ et tout ouvert $U$ du site $X_{i,\lambda}$, d'un ensemble $\F_i(U)$, fonctoriellement en $i$ et en $U$. On définit de même les topos $\Tot X_\bullet$  (resp. $\Tot X_{\bullet \ell {\rm \acute{e}t}}$), en remplaçant le site $X_{i,\lambda}$ par le site étale (resp. $\ell$-étale, voir section \ref{subsec:Kpi1}) du schéma $X_i$.

Soit $\F_\bullet$ un faisceau abélien du topos $\Tot X_\bullet^{(\lambda)}$. La construction pour chaque $i$ d'une résolution flasque $\F_i^\bullet$ de $\F_i$ permet de représenter $\RG(\Tot X_\bullet^{(\lambda)},\F_\bullet)$ par le complexe total associé au complexe double $(\Gamma(X_i^{(\lambda)},\F_i^j)_{i,j}$. Une telle résolution flasque se détermine en choisissant un point géométrique de chaque composante connexe de $X_0$, puis en calculant toutes ses images dans $X_0,\dots,X_r$ par les applications de bord et de dégénérescence. Ceci fournit un ensemble fini de points, qui forment un topos simplicial discret $P_{\bullet \leqslant r}$. Le morphisme $u\colon P_{\bullet\leqslant r}\to X^{(\lambda)}_{\bullet\leqslant r}$ est à image inverse conservative, et pour tout faisceau $\F_\bullet$ sur $\Tot X_{\bullet}^{(\lambda)}$, le morphisme $\F_{\bullet \leqslant r}\to u_\star u^\star\F_{\bullet \leqslant r}$ est le début d'une résolution flasque de $\F_{\bullet \leqslant r}$. \'{E}tant donné un schéma simplicial en groupes abéliens représentant $\F_{\bullet \leqslant r}$, il est possible de calculer explicitement un schéma simplicial représentant $u_\star u^\star \F_{\bullet\leqslant r}$ : c'est simplement un coproduit de composantes connexes de $X^{[\lambda]}$ \cite[4.2.1]{mo}. \\

Par descente cohomologique, il y a un isomorphisme dans $\DD^b_c(X,\Lambda)$ : \[ \RG(X,\Lambda)\xrightarrow{\sim} \RG(\Tot X_\bullet,\Lambda).\]
Si de plus chaque $X_i$ est un $K(\pi,1)$ pro-$\ell$, il y a un isomorphisme :
\[ \RG(X,\Lambda)\xrightarrow{\sim} \RG(\Tot X_{\bullet \ell{\rm \acute{e}t}},\Lambda).\]
Pour tout entier $j\geqslant 0$, il y a un isomorphisme :
\[ \colim_\lambda \HH^j(\Tot X_{\bullet \ell{\rm \acute{e}t}}^{(\lambda)},\Lambda) \xrightarrow{\sim} \HH^j(\Tot X_{\bullet \ell{\rm \acute{e}t}},\Lambda).\]
La cohomologie de $\Tot X_{\bullet \ell{\rm \acute{e}t}}^{(\lambda)}$ est calculée par la suite spectrale
\[ E_1^{i,j}=\HH^j(X_i^{(\lambda)},\Lambda)=\HH^j(\prol(X_i)^{(\lambda)},\Lambda)\Rightarrow \HH^{i+j} (\Tot X_{\bullet \ell{\rm \acute{e}t}}^{(\lambda)},\Lambda).\]

Les résultats de la section suivante permettent de déterminer $\HH^j(X,\Lambda)$ comme l'image d'un morphisme \[ \HH^j(\Tot X_{\bullet \ell{\rm \acute{e}t}}^{(\alpha)},\Lambda)\to \HH^j(\Tot X_{\bullet \ell{\rm \acute{e}t}}^{(\beta)},\Lambda).\]

\subsection{Systèmes essentiellement constants et cohomologie des polycourbes $\ell$-élémentaires}

Soit $\C$ une catégorie abélienne dans laquelle les colimites filtrantes existent. Soit $(A_i)_{i\in\NN}$ un système inductif d'objets noethériens dans $\C$. Si la colimite $A_\infty$ est noethérienne, il existe des entiers $j,k$ tels que la restriction du morphisme $A_k\to A_\infty$ à l'image de $A_j\to A_k$ soit un isomorphisme. 

\begin{df} Soient $N$ un entier naturel et $\phi \colon \NN\to \NN$ une application. Le système inductif $A$ est dit $(N,\phi)$-essentiellement constant (ou $(N,\phi)$-ec) s'il vérifie les deux conditions suivantes : \begin{enumerate}
\item pour tout entier $j$ et tout $k\geqslant \phi(j)$, $\ker(A_j\to A_{\phi(j)})\to \ker(A_j\to A_k)$ est un isomorphisme ;
\item pour tout $j\geqslant N$, $\im(A_N\to A_{\phi(j)})\to \im(A_j\to A_{\phi(j)})$ est un isomorphisme.\\
\end{enumerate}
\end{df}
Un tel système inductif $A$ est dit explicitement essentiellement constant (ou eec) s'il existe $(N,\phi)$ calculables tels que $A$ soit $(N,\phi)$-ec.
L'essentielle constance des systèmes inductifs se comporte bien vis-à-vis des suites exactes.
\begin{prop}\cite[Prop. 5.6]{mo} Soit \[ 0\to A'\to A\to A''\to 0\]
une suite exacte de systèmes inductifs dans $\C$ indexés par $\NN$. \begin{enumerate}
\item Si $A$ est $(N,\phi)$-ec et $A''$ est $(N'',\phi'')$-ec alors $A'$ est $(\phi''(N),\phi)$-ec.
\item Si $A$ est $(N,\phi)$-ec et $A'$ est $(N',\phi')$-ec alors $A''$ est $(N,\max(\phi,N'))$-ec.
\item Si $A'$ est $(N',\phi')$-ec et $A''$ est $(N'',\phi'')$-ec alors $A$ est $(\max(N',N''),\phi'\circ\phi'')$-ec.
\end{enumerate}
\end{prop}
Ceci implique qu'étant donné un système inductif de suites spectrales $(E_{2,\lambda}^{p,q}\Rightarrow \HH^{p+q}_\lambda)_\lambda$, si chaque $(E_{2,\lambda}^{p,q})_\lambda$ est eec et s'il existe un entier $r$ tel que $E_{2,\lambda}^{p,q}=0$ dès que $q>r$ alors $(\HH^{p+q}_\lambda)_\lambda$ est eec. Le résultat suivant précise cet énoncé dans un cas particulier.\\

\begin{lem}\label{lem:E2ec} Soit $(E_{2,\lambda}^{pq}\Longrightarrow H_\lambda^{p+q})_\lambda$ un système inductif de suites spectrales. Supposons que pour tous entiers $p,q$, le système inductif $(E_{2,\lambda}^{pq})_\lambda$ soit $(N,\phi)$-essentiellement constant. Alors le système $(E_{r,\lambda}^{pq})_\lambda$ est $(\phi^{r-2}(N),\max(\phi,\phi^{r-3}(N))$-essentiellement constant.
\begin{proof} Ceci se démontre par récurrence sur $r$. Le groupe $E_{r+1,\lambda}^{pq}$ est un quotient de la forme $\ker(f_\lambda)/\im(g_\lambda)$, où $f$ et $g$ sont des flèches de la page $E_r$. Par la proposition précédente et l'hypothèse de récurrence, le système $(\ker f_\lambda)_\lambda$ est $(\phi^{r-1}(N),\max(\phi,\phi^{r-3}(N))$-essentiellement constant, et le système $(\im g_\lambda)_\lambda$ est $(\phi^{r-2}(N),\max(\phi,\phi^{r-3}(N))$-essentiellement constant. Une nouvelle application de la proposition précédente montre que le système $(E_{r+1,\lambda}^{pq})_{\lambda}$ est  $(\phi^{r-1}(N),\max(\phi,\phi^{r-2}(N))$-essentiellement constant.
\end{proof}
\end{lem}

Soit \[ X=X_d\to X_{d-1}\to \dots \to X_1\to \Spec k\]
une polycourbe $\ell$-élémentaire (affine) sur $k$. Notons $\pi$ le groupe fondamental pro-$\ell$ de $X$ et $\bareta$ un point générique géométrique de $X_1$. Il y a, comme vu précédemment, une suite exacte
\[ 1\to \prol(X_\bareta) \to \pi\to \prol(X_1)\to 1\]
où le groupe de droite est un pro-$\ell$-groupe libre, et le groupe de gauche est (par récurrence sur la dimension) une extension itérée de tels groupes. Pour tout niveau d'approximation $\lambda$, cette suite fournit une nouvelle suite 
\[ 1\to \pi''_\lambda\to \pi^{(\lambda)}\to \prol(X_1)^{(\lambda)}\to 1\]
où $\pi''_\lambda$ est simplement défini comme étant le noyau de $\pi^{(\lambda)}\to \prol(X_1)^{(\lambda)}$. Les groupes de cette suite se calculent explicitement. Il est montré dans \cite[Cor. 6.4]{mo} que si $(\prol(X_\bareta)^{(\lambda)})_\lambda$ est eec, il en est de même pour $\pi''_\lambda$. Notons qu'un couple $(N,\phi)$ explicitant ce fait est exhibé dans \cite[Prop. 6.3] {mo} : sa détermination est primitivement récursive.
Par récurrence sur la dimension, il suffit pour montrer que le système $(\HH^j(\pi^{(\lambda)},\Lambda))_\lambda$ est eec, de montrer que c'est le cas lorsque $X$ est une courbe. 
C'est évident lorsque $j=0$. Pour $j=1$, nous verrons dans la section \ref{subsec:cohlisG} que le système est $(2,\lambda\mapsto \lambda+1)$-essentiellement constant. Pour $j\geqslant 2$, il est possible de trouver $\lambda_0$ tel que pour tout $\lambda > \lambda_0$, $\HH^j(\pi^{(\lambda_0)},\Lambda)\to \HH^j(\pi^{(\lambda)},\Lambda)$ soit nulle. Il suffit pour cela que $\HH^1(\pi^{[\lambda_0]},\Lambda)\to \HH^j(\pi^{[\lambda]},\Lambda)$ soit nulle \cite[Rem. 6.6]{mo}, ce qui est le cas dès que $\lambda_0\geqslant 2$ (voir lemme \ref{lem:revtriv}). Le système $(\HH^j(\pi^{(\lambda)},\Lambda))_\lambda$ est donc également $(2,\lambda\mapsto \lambda+1)$-ec. \\

\subsection{Adaptation de l'algorithme à des cas plus généraux}

\paragraph{Variétés singulières} Soit $X$ un schéma de type fini sur $k$. Le théorème de résolution des singularités de de Jong assure qu'il est localement lisse pour la topologie des altérations. Ceci permet de se ramener par changement de base au cas des schémas lisses, mais pose problème en termes de complexité : à moins d'étudier en détail chaque construction de \cite{dejong_sing}, une recherche non bornée est nécessaire pour trouver une altération convenable.

\paragraph{Images directes de faisceaux constructibles} 
Pour calculer les groupes de cohomologie d'un faisceau constructible $\F$, il suffit de plonger $\F$ dans un faisceau $\tilde \F$ tel que pour tout $i>0$, le morphisme \[\R^if_\star\F\to \R^i f_\star\tilde \F\] soit nul. 
En effet, $\R^j f_\star\F$ se calcule alors récursivement comme \[\R^j f_\star\F=\coker(\R^{j-1}f_\star \tilde \F\to \R^{j-1}f_\star (\tilde{\F}/\F)).\] 
La calcul d'un tel $\tilde \F$ est décrit en détail dans \cite[11.4.4]{mo}. Nous donnerons une construction explicite dans le cas des courbes dans la section \ref{sec:coheff}.

\paragraph{Cohomologie d'un complexe de faisceaux constants sur une courbe} Soit $X$ une courbe affine lisse sur $k$. Soit $K=[K^0\to\dots\to K^N]$ un complexe de $\Lambda$-modules. Considérons toujours les quotients de Frattini $\pi^{(\lambda)}$ du groupe fondamental pro-$\ell$ de $X$. Rappelons qu'il y a une suite spectrale \cite[0AVG]{stacks} de deuxième page 
\[ E_2^{p,q}=\HH^p(\pi,\HH^qK)\Longrightarrow \HH^{p+q}(\pi,K).\]
Rappelons que pour tout $p\geqslant 2$ et tout $q\in\{ 0\dots N\}$, le système $(\HH^p(\pi^{(\lambda)},\HH^qK))_\lambda$ est $(2,\lambda\mapsto \lambda+1)$ essentiellement constant. 
\begin{lem} Pour tout entier $i\geqslant 0$, le système inductif $(\HH^i(\pi^{(\lambda)},K))_{\lambda\geqslant 2}$ est $(N+2,\lambda\mapsto \lambda+N+2)$-essentiellement constant.
\begin{proof}
Considérons la suite spectrale $E_{2,\lambda}^{pq}=\HH^p(\pi^{(\lambda)},\HH^q(K))\Longrightarrow \HH^{p+q}(\pi^{(\lambda)},K)$. Le terme $E_{2,\lambda}^{pq}$ est nul dès que $q>N$. Par conséquent, toutes les flèches de la page $E_{N+2}$ sont nulles, et la suite spectrale converge : $E_{\infty,\lambda}=E_{N+2,\lambda}$. Pour tous $p,q$, le système $(E_{2,\lambda}^{pq})_{\lambda\geqslant 2}$ est $(2,\lambda\mapsto \lambda+1)$-essentiellement constant. Le lemme \ref{lem:E2ec} assure alors que le système $(E_{N+2,\lambda}^{pq})_\lambda$ est $(N+2,\lambda\mapsto \max(\lambda,N)+1))$-essentiellement constant.\\
Le groupe $\HH^i(\pi^{(\lambda)},K)$ est extension itérée des $E_{N+2,\lambda}^{p,q}$ avec $p+q=i$. Plus précisément, il y a une suite exacte courte de groupes abéliens \[ 0\to A_{i,\lambda}\to \HH^i(\pi^{(\lambda)},K) \to E_{N+2,\lambda}^{0,i} \to 0\]
puis, pour $1\leqslant r\leqslant i-1$, une suite exacte courte
\[ 0\to A_{i-r,\lambda}\to A_{i-r+1,\lambda}\to E_{N+2,\lambda}^{r,i-r}\to 0\]
et enfin une suite exacte courte
\[ 0\to E_{N+2,\lambda}^{i,0}\to A_{1,\lambda}\to E_{N+2,\lambda}^{i-1,1}\to 0.\]
D'après les résultats précédents, le système inductif $(A_{1,\lambda})_\lambda$ est $(N+2,\max(\lambda\mapsto \lambda+2,N+1))$-essentiellement constant. Une récurrence immédiate donne alors que $\HH^i(\pi^{(\lambda)},K)$ est $(N+2,\lambda\mapsto \max(N+1,\lambda+N+2))$-essentiellement constant, c'est-à-dire $(N+2,\lambda\mapsto \lambda+N+2)$-essentiellement constant. 
\end{proof}
\end{lem}

\subsection{Remarques sur la complexité de l'algorithme}

\begin{prop} Soit $X$ un schéma lisse de type fini sur un corps algébriquement clos $k$. Soit $\ell$ un nombre premier inversible dans $k$. L'algorithme décrit dans \cite{mo} qui calcule les groupes de cohomologie $\HH^j(X,\ZZ/\ell\ZZ)$ pour $j\in \{0\dots 2\dim X\}$ est primitivement récursif.
\begin{proof}
Nous avons montré ci-dessus que le calcul d'un hyperrecouvrement de $X$ par des $K(\pi,1)$ pro-$\ell$ était primitivement récursif, de même que le calcul des $r$ premiers étages d'une résolution flasque du faisceau $\ZZ/\ell\ZZ$ dans chaque topos $\Tot X_{\bullet\leqslant r}^{[\lambda]}$. Nous avons vu qu'il en est de même, pour tout entier $j\leqslant 2\dim X$, des fonctions explicitant le fait que le système inductif $(\HH^j(\Tot X_{\bullet\leqslant r}^{(\lambda)},\Lambda))_\lambda$ est essentiellement constant. 
\end{proof}
\end{prop}

\paragraph{De la difficulté de montrer qu'il est élémentaire} Considérons la forme la plus simple de l'algorithme, qui calcule les $\HH^i(X,\Lambda)$ pour un schéma $X$ intègre lisse de type fini sur $k$. Notons $d$ sa dimension. Soit $X_\bullet$ un hyperrecouvrement étale de $X$. Notons $s_i$ (resp. $D_i$) le nombre (resp. le degré maximal sur $X$) des composantes connexes de $X_i$. Alors \[ s_{i+1}\leqslant s_i^2D_i \text{~~et~~} D_{i+1}\leqslant D_i^2.\]
Remarquons que dès que l'une des composantes connexes de $X_0$ est galoisienne de degré $b$ sur son image dans $X$, le nombre de composantes connexes de $X_r$ est supérieur à $b^r$. 
Déjà lorsque $X$ est une courbe, les systèmes $(\HH^q(X_p^{(\lambda)},\Lambda))_\lambda$  pour $2\leqslant q\leqslant r$ sont $(2,\lambda\mapsto \lambda+1)$-essentiellement constants pour la filtration de Frattini, et nous n'avons pour l'instant pas de borne pour une autre filtration. Dans le cas d'un schéma de dimension supérieure, nous n'arriverons certainement pas à obtenir mieux. Il faut donc s'attendre à ce que la première page \[E_1^{pq}=\HH^q(X_p^{(\lambda)},\Lambda)\]
de la suite spectrale calculant $\HH^j(X,\Lambda)$ soit constituée d'au moins $r\times r$ systèmes inductifs, dont au moins $(r-2)\times r$ ne sont pas constants.
Le raisonnement du lemme \ref{lem:E2ec} montre alors que par exemple $(\HH^{r}(\Tot X_{\bullet\leqslant r}^{(\lambda)},\Lambda))_\lambda$ est $(r-4,\lambda\mapsto \lambda+r-5)$-essentiellement constant. Comme l'algorithme prend pour $r$ un nombre strictement supérieur à la dimension $d$ de $X$, il faut s'attendre à calculer au moins le $(d-4)$-ième revêtement de Frattini d'une courbe, dont la complexité est une exponentielle à $d-4$ étages, qui n'est pas une fonction élémentaire en $d$. Il faudrait donc, afin d'obtenir une complexité élémentaire, trouver des bornes sur l'essentielle constance des systèmes ci-dessus pour une filtration plus fine, comme celle proposée dans \cite[§3]{mo}.

\section{L'algorithme de Couveignes}\label{sec:couv}

\subsection{Données et principe de l'algorithme}

Soient $\ell$ un nombre premier, et $\FF_q$ un corps fini de caractéristique différente de $\ell$. Soit $j$ un entier naturel. Soit $X$ une courbe projective intègre lisse sur $\FF_q$ de genre $g$, décrite par un modèle plan et ses branches singulières comme dans l'annexe \ref{subsec:modplan}. L'algorithme décrit dans \cite{couveignes_linearizing} a pour but de calculer des diviseurs représentant les éléments de $J(\FF_q)[\ell^j]$. Nous nous contentons ici de résumer cet algorithme et de décrire en \ref{subsec:racJX} le travail restant à faire pour généraliser l'algorithme à la division par $n$ dans $J(\overline{\FF_q})$.\\

Soit $i$ un entier naturel. Soit $Q$ une puissance de $q$ assez grande pour que $|J(\FF_Q)[\ell^{2g}]|=\ell^{2g}$. La procédure décrite dans la section \ref{subsec:tiral} permet de tirer aléatoirement des diviseurs non principaux dans $J(\FF_Q)$ qui engendrent un sous-groupe de $J(\FF_Q)$ de petit indice. L'application de Kummer décrite dans la section \ref{subsec:appkummer} permet de leur associer des éléments d'un sous-groupe $H\subset J(\FF_Q)[\ell^i]$, où $i$ est un entier naturel. En choisissant $i$ assez grand, il est possible de faire en sorte que le sous-groupe $H$ contienne $J(\FF_Q)[\ell^j]$. Le couplage de Weil permet alors de retrouver $J(\FF_q)[\ell^j]$. La complexité de l'algorithme de Couveignes est donnée par le résultat suivant.

\begin{theorem}\cite[Th. 1]{couveignes_linearizing} Il existe un algorithme probabiliste (Monte-Carlo) qui prend en entrée une courbe projective lisse géométriquement intègre $X$ de genre $g$ sur $\FF_q$ définie par un modèle plan de degré $d$, un diviseur $\FF_q$-rationnel $O=D^+-D^-$ de degré 1 sur $\FF_q$ avec $\deg D^+=O(g)$, un nombre premier $\ell$ ne divisant pas $q$, un entier naturel $j$ et la fonction zêta de $X$, et renvoie une $\ZZ/\ell^j$-base de $\Pic(X)(\FF_q)[\ell^j]$ dont les éléments sont des diviseurs de la forme $G-gO$ où $G$ est effectif. Le nombre d'opérations arithmétiques effectué par l'algorithme est polynomial en $d,g,\log q$ et $\ell^j$ ; l'algorithme renvoie un sous-groupe de $\Pic(X)(\FF_q)[\ell^j]$, qui est le groupe entier avec probabilité $\geqslant \frac{1}{2}$.
\end{theorem}

\subsection{Tirage aléatoire de diviseurs}\label{subsec:tiral}

L'algorithme de Couveignes construit un sous-groupe assez grand de $\Pic^0(X)(\FF_q)$ en tirant aléatoirement des diviseurs $\FF_q$-rationnels sur $X$. L'idée est la suivante : le modèle plan donné de $X$ fournit un morphisme $x\colon X\to\PP^1$. L'image par $x$ de l'ensemble $\mathcal{P}(r,q)$ de $\FF_q$-places de $X$ de degré $r$ est incluse dans l'ensemble $\mathcal{U}(r,q)$ des polynômes unitaires irréductibles de degré $r$ sur $\FF_q$. Le tirage des diviseurs est réalisé en tirant aléatoirement selon une distribution uniforme un polynôme unitaire de degré $r$ à coefficients dans $\FF_q$, puis en testant son irréductibilité, avant de calculer le cas échéant ses antécédents dans $\mathcal{P}(r,q)$. Le choix d'un diviseur effectif $\Omega$ de degré $r$ sur $X$ détermine une application $\mathcal{P}(r,q)\to \Pic^0(X)(\FF_q), \alpha\mapsto [\alpha]-\Omega$.
 Il existe une unique mesure $\mu$ sur $\mathcal{P}(r,q)$ telle que les fibres non vides de $\mathcal{P}(r,q)\to \mathcal{U}(r,q)$ soient de même mesure, et tous les points d'une même fibre aient la même mesure. La mesure de l'ensemble des places dans $\mathcal{P}(r,q)$ d'image non nulle dans $\Pic^0(X)(\FF_q)$ est alors supérieure à $\frac{1}{2d}$, et au bout de $2d$ tirages, la probabilité d'avoir obtenu un diviseur non trivial est supérieure à $\frac{1}{2}$.\\
Cette méthode permet de construire par tirages successifs des diviseurs engendrant avec probabilité $\geqslant \frac{1}{2}$ un sous-groupe $G$ de $\Pic^0(X)(\FF_q)$ d'indice inférieur à une borne $\iota$ fixée. Pour que la complexité de l'algorithme reste polynomiale en les entrées, il n'est pas possible d'obtenir ainsi $\Pic^0(X)(\FF_q)$ tout entier ; il est toutefois loisible de choisir $\iota$ linéaire en $d$ et $g$ (dans l'article, $\iota=\max\{48g,24d,720\}$).

\subsection{Groupes $\ell$-divisibles}

Soit $\chi_\ell$ la réduction mod $\ell$ du polynôme caractéristique de l'endomorphisme de Frobenius sur $J_X$. Il se factorise en
\[ \chi_\ell(t)=(t-1)^b{\bar f}^\perp(t)\]
où $\bar{f}^\perp(t)$ est premier à $t-1$. Le lemme de Hensel permet de relever cette factorisation en une unique factorisation en produit de polynômes unitaires
\[ \chi(t)=f(t)f^\perp(t)\in \ZZ_\ell[t].\]
Soient $e,e^\perp$ les idempotents de $\ZZ_\ell[t]/\chi$ correspondant respectivement aux éléments $(1,0)$ et $(0,1)$ de $\ZZ_\ell[t]/f\times\ZZ_\ell[t]/f^\perp$. Les éléments $e_1(\phi_q)$ et $e_1^\perp(\phi_q)$ définissent des endomorphismes du groupe $J_X[\ell^\infty]$, et leurs images respectives sont les sous-groupes $\ell$-divisibles $\GG_1$ et $\GG_1^\perp$ de $J_X[\ell^\infty]$. Ceci définit une décomposition \[ J_X[\ell^\infty]=\GG_1\times\GG_1^\perp \] en produit de deux groupes $\ell$-divisibles.

\subsection{Une surjection $J(\FF_q)\to J[\ell^j](\FF_q)$}

\subsubsection{Détermination d'une extension de corps adéquate}\label{subsubsec:extcorpscouv}

Soit $A$ une variété abélienne sur $\FF_q$. Soit $\GG$ un sous-groupe de $A[\ell^\infty]$. Pour une extension $L$ de $\FF_q$, notons $\GG(L)=\GG\cap A(L)$.
Soit $\chi\in\ZZ_\ell[t]$ le polynôme caractéristique de l'automorphisme de Frobenius sur $A$. Voici comment déterminer une extension $L/\FF_q$ telle que $\GG[\ell^j](L)=\GG[\ell^j](\overline{\FF_q})$. Considérons l'application \[ (\ZZ/\ell^j\ZZ)[t]/(\chi~{\rm mod}~ \ell^j)\to \End(\GG[\ell^j])\]
qui envoie $t$ sur l'automorphisme de Frobenius $\phi_q$. L'ordre de $\phi_q$ dans $\Aut(\GG[\ell^j])$ est le degré de la plus petite extension de $\FF_q$ sur laquelle sont définis tous les points de $\GG[\ell^j]$. Il divise l'ordre de $t$ dans $\left( (\ZZ/\ell^j\ZZ)[t]/\chi\right)^\times$. \\

Dans le cas où $\GG$ est le groupe $\GG_1$ de la section précédente, voici comment déterminer cet ordre. Soit $b$ la valuation $(t-1)$-adique de $\chi\in\FF_\ell[t]$. L'entier $\gamma=\lceil \log_{\ell} b\rceil $ vérifie $(t-1)^{\ell^\gamma}=0\in \FF_\ell[t]/\chi$, et l'ordre de $t$ dans $(\FF_\ell[t]/\chi)^\times$ divise $\ell^\gamma$. Par conséquent, l'ordre de $t$ dans $(\ZZ/\ell^j\ZZ[t]/\chi)^\times$ divise $\ell^{\gamma+j-1}$. Dans le cas où $\GG=J_X$, l'ordre de $t$ divise\[ A_k\coloneqq\ell^{\lceil \log_\ell(2g)\rceil +j-1}\prod_i (\ell^{f_i}-1)\] où les $f_i$ sont les degrés des facteurs irréductibles de $\chi$ modulo $\ell$.

\subsubsection{L'application de Kummer}\label{subsec:appkummer}

Soit $A$ une variété abélienne sur $\FF_q$. Soit $\GG$ un sous-groupe $\ell$-divisible de $A[\ell^\infty]$. Supposons que $\GG[\ell^j](\FF_q)=\GG[\ell^j](\overline{\FF_q})$. Pour chaque $P\in \GG(\FF_q)$, choisissons un point $R_P\in \GG(\overline{\FF_q})$ tel que $\ell^j R_P=P$. Associons à $P$ l'élément \[ K_{\ell^j,q}(P)\coloneqq R_P^{\phi_q}-R_P\in \GG[\ell^j](\overline{\FF_q})=\GG[\ell^j](\FF_q).\] Ceci définit un isomorphisme de groupes abéliens \[ K_{\ell^j,q}\colon\GG(\FF_q)/\ell^j\GG(\FF_q)\to \GG[\ell^j](\FF_q).\]
Il n'est évidemment pas question de calculer explicitement l'élément $R_P$ : cet isomorphisme se calcule de façon plus efficace. Comme vu dans la partie précédente, il est possible de calculer un entier $a$ tel que $t^a=1$ dans $(\ZZ/\ell^j\ZZ)[t]/\chi$. En relevant cette identité à $\ZZ_\ell$, il existe un unique $M_{j,a}\in\ZZ_\ell[t]/\chi$ tel que $t^a-1=M_{j,a}(t)\ell^j$. L'isomorphisme de Kummer ci-dessus s'identifie alors à $P\mapsto M_{j,a}(\phi_q)(P)$.

\subsection{Calcul de la $\ell$-torsion du groupe de Picard}

Soit $Q=q^{(\ell-1)\ell^{\gamma+j-1}}$. Alors $\GG_1[\ell](\FF_Q)=\GG_1[\ell](\overline{\FF_q})$. On commence par tirer des éléments de $J(\FF_{Q})$ qui engendrent un sous-groupe d'indice inférieur à $\iota=\max(48g,24d,720)$. L'image de ce sous-groupe par la composition de morphismes surjectifs 
\[ J(\FF_Q)\overset{e(\phi_q)}{\longrightarrow}\GG_1(\FF_Q)\overset{K_{n,q}}{\longrightarrow}\GG_1[\ell^j](\FF_Q) \]
est un sous-groupe $H$ de $\GG_1[\ell^j](\FF_Q)$ d'indice au plus $\iota$. Dès que $j$ est supérieur à $\delta\coloneqq\lceil \log_\ell \iota\rceil$, le sous-groupe $H$ contient $\GG_1[\ell^{j-\delta}](\FF_Q)=\GG_1[\ell^{j-\delta}](\overline{\FF_q})$. Calculer l'ordre des éléments de $H$ permet alors de déterminer $\GG_1[\ell^{j-\delta}](\FF_Q)$. Le groupe $J[\ell^{j-\delta}](\FF_q)$ s'en déduit comme le noyau de l'endomorphisme $\phi_q-\id$ de $\GG_1[\ell^{j-\delta}](\FF_Q)$.

\subsection{Racines $n$-ièmes d'éléments non nuls de $J_X(\FF_q)$}\label{subsec:racJX}

Soit $n$ un entier naturel premier à $q$. L'algorithme de Couveignes produit les éléments de $n$-torsion de $J_X(\FF_q)$. Une question plus générale, et qui sert dans le calcul de la cohomologie des courbes affines, est la suivante : étant donné un élément $D\in J_X(\FF_q)$, déterminer $D'\in J_X(\overline{\FF_q})$ tel que $nD'=D$. L'endomorphisme $\phi_q-\id$ de $J_X$ étant surjectif, il existe $D_1\in\Pic^0(X)$ tel que $D=D_1^{\phi_q}-D_1$. Notons $\tau_{D_1}$ la translation par $D_1$ sur $J_X$. Alors l'application
\[ K_{\ell^j,q}\circ \tau_{D_1}\colon\GG_1(\FF_Q)\to \{ E\in \GG_1(\FF_Q)\mid nE=D\}\]
est surjective. La même procédure que précédemment permet alors de tirer un élément non trivial de $J(\FF_Q)$, de calculer son image dans $\GG_1(\FF_Q)$ puis par translation un antécédent de $D$ par $[n]_{\GG_1(\FF_Q)}$. Ici, connaissant déjà $J[n]$, il suffit de trouver un seul antécédent.
\'{E}tant donné $D$, le problème du calcul dans $J(\overline{\FF_q})$ d'un antécédent par la multiplication par $n$ se réduit donc au calcul d'un antécédent par $\phi_q-\id$.

En fonction des objectifs de complexité, cette réduction peut ou non apporter un bénéfice : l'isogénie de multiplication par $n$ est de degré $n^{2g}$, alors que $\phi_q-\id$ est de degré $O(q^g)$. Cependant, il n'y a à la connaissance de l'auteur de ces lignes pas de méthode efficace pour calculer un antécédent par $\phi_q-\id$. Une possibilité serait d'employer un analogue de l'algorithme de Huang-Ierardi présenté dans la section suivante, qui consisterait étant donné un diviseur $D$ à chercher explicitement des diviseurs $F$ tels que $F^{\phi_q}-F$ soit équivalent à un diviseur très simple (voir définition \ref{df:tressimple}) équivalent à $D$. Cependant, la complexité de cet algorithme serait au moins linéaire en $q$.

\section{L'algorithme de Huang et Ierardi}\label{sec:huang}

Dans toute cette section, $C_0=\Proj k_0[x,y,z]/(f)$ désigne une courbe projective plane intègre sur $k_0=\FF_q$, et $X_0$ désigne sa normalisée. Notons $d$ le degré de $f$ et $g$ le genre de $X_0$. Notons $X=X_0\times_{k_0}k$. Soit $J_X$ la jacobienne de $X$. 

\subsection{Structure de l'algorithme}

\subsubsection{Données et hypothèses}

La courbe $C_0$ est décrite par le polynôme homogène $f\in k_0[x,y,z]$ de degré $d$. 
Pour le calcul de $C_0$ étant donné $X_0$, et vice-versa, voir les annexes \ref{subsec:modplan} et \ref{subsubsec:desing}. 
Notons $Q_1,\dots,Q_r$ les points singuliers de $C$, et $m_1,\dots,m_r$ leurs multiplicités respectives. Rappelons que $r\leqslant { d-1\choose 2}$. Les singularités de $C$ sont supposées ordinaires ; ceci s'obtient par des transformations quadratiques \cite[§1.7]{kollar_resolution}, et demande éventuellement de remplacer $C$ par une courbe $C'$ de degré $d'=O(2^{d^2})$. Supposons $C$ remplacée par $C'$, et $d$ par $d'$. 
\\
Pour chaque $i\in \{1\dots r\}$ est donnée (voir \cite[§3.3]{huang_counting}) une courbe $C'_{Q_i}$ birationnelle à $C$ sur laquelle $Q_i$ a $m_i$ antécédents $Q_{ij}$ ; elle est obtenue en éclatant $C$ en $Q_i$.
Il est possible que la courbe $C$ contienne un point $Q_T$ dit "terrible" (voir \cite[Appendix A, p113]{fulton}, qui est alors unique. Il peut être régulier ou singulier. Notons $m_T$ sa multiplicité. Il est toujours possible de construire, en au plus deux transformations birationnelles \cite[§3.3]{huang_counting}, une courbe $C'_T$ birationnelle à $C$ sur laquelle $Q_T$ a $m_T$ antécédents réguliers et non terribles. Si $T$ est régulier, il sera parfois traité de la même façon que les points singuliers, ce que nous préciserons le cas échéant.
\\

Le corps fini $k_0$ est supposé assez grand ($d^6<|k_0|$) pour pouvoir énumérer un nombre suffisant d'éléments distincts dans $k_0$ dans les diverses constructions, ainsi que pour contenir les coordonnées des $Q_i$ et des $Q_{ij}$ (ce qui nécessite de remplacer le corps de départ par une extension de degré au plus $d'^3!$), ainsi que d'un point régulier de $C$, qui sera noté $\infty$.
Dans la description de l'algorithme, nous supposerons que $k_0$ a été remplacé par une telle extension.

\subsubsection{Le résultat}

L'algorithme développé par Huang et Ierardi avait pour but principal le comptage de points sur les courbes, d'où le premier énoncé. Cependant, les calculs de complexité présentés dans \cite[§5.4]{huang_counting} démontrent également le second résultat.

\begin{theorem}\cite[Theorem 1.1]{huang_counting}\label{th:huang}
Soit $d\in\NN$. Il existe un réel $\alpha>0$ et un algorithme probabiliste Las Vegas qui prend en entrée une puissance $q$ d'un nombre premier, un entier $n$ et une courbe projective plane $C_0\subset\PP^2_{\FF_q}$ de degré $d$ n'ayant que des singularités ordinaires, de normalisée $X_0$, et qui renvoie le $\Lambda[\mathfrak{G}_0]$-module $\HH^1(X_0\times_{\FF_q}\overline{\FF_q},\mu_n)$ en $O((d\log(q))^{d^{\alpha}})$ opérations dans $\FF_q$ dès lors que $n=O(d^2\log q)$. La probabilité de succès de l'algorithme est supérieure à $\frac{1}{2}$.
\end{theorem}

\begin{prop} Soient $n,d$ deux entiers naturels. Soit $q$ une puissance d'un nombre premier $p$. Il existe des entiers $\alpha,\beta >0$ et un algorithme déterministe qui prend en entrée une courbe projective plane $C_0\subset\PP^2_{\FF_q}$ de degré $d$ n'ayant que des singularités ordinaires, de normalisée $X_0$, et qui renvoie le $\Lambda[\mathfrak{G}_0]$-module $\HH^1(X_0\times_{\FF_q}\overline{\FF_q},\mu_n)$ en $O_{n,d,q}(p^{\alpha d^3!} (dn)^{d^{\beta}})$ opérations dans $\FF_q$.
\end{prop}

L'énoncé correspondant pour une courbe plane quelconque s'obtient en remplaçant $d$ par $2^{d}$ dans la complexité. Nous nous servirons par la suite du résultat suivant, qui se déduit des précédents en adaptant l'algorithme de Huang et Ierardi en suivant les remarques \ref{rk:huangrac}, \ref{rk:huangrac2} et \ref{rk:huangrac3}, à l'aide des calculs de complexité effectués dans les sections \ref{subsubsec:detrep} et \ref{subsubsec:calculJx}.\\

\begin{prop}\label{th:huangrac}
Soit $d\in\NN$. Il existe un algorithme probabiliste (Las Vegas) qui prend en entrée une puissance $q$ d'un nombre premier, un entier $n$, une courbe projective plane $C_0\subset\PP^2_{\FF_q}$ de degré $d$ et de genre $g$ n'ayant que des singularités ordinaires, de normalisée $X_0$, et un diviseur $F=F^+-F^-\in\Div^0(X_0)(\FF_q)$ avec $\deg(F^+)\leqslant g$, et qui renvoie le $\Lambda[\mathfrak{G}_0]$-module \[ \{ [D]\in\Pic^0(X_0\times_{\FF_q}\overline{\FF_q})\mid n[D]=F\} \] en $\mathcal{P}(d,g,n,\log q)^{O(g^4)}$ opérations dans $\FF_q$, où $\mathcal{P}$ est un polynôme. La probabilité de succès de l'algorithme est supérieure à $\frac{1}{2}$.
\end{prop}

\subsubsection{Résumé de l'algorithme}

L'algorithme de Huang et Ierardi est essentiellement l'algorithme de Brill-Noether (voir annexe \ref{subsec:rr}) appliqué à un diviseur générique. Il consiste à considérer un diviseur $D$ dont les coordonnées des points du support sont des indéterminées, calculer un diviseur équivalent à $nD$, et construire comme dans l'algorithme de Brill-Noether une fonction dont $nD$ est le diviseur. L'existence d'une telle fonction se traduit par des équations sur les indéterminées ; il suffit alors de trouver des points dans le schéma défini par ces équations. \\

\begin{df}\label{df:tressimple}
Un diviseur $D=\sum a_iP_i\in\Div(X)$ est dit très simple si les $P_i$ sont des points réguliers non-terribles deux à deux distincts de $C$ et si $a_i=\pm 1$ pour tout $i$.
\end{df}

Soit $D\in \Div^0(X)$. Rappelons qu'un point $\infty\in X$ a été fixé. Le théorème de Riemann-Roch assure que l'espace $\LL(D+g\infty)$ est non vide ; il existe donc un diviseur effectif $E$ de degré $g$ sur $X$ tel que $E-g\infty$ soit équivalent à $D$. Par conséquent, les diviseurs de degré zéro seront représentés par un diviseur équivalent de la forme $E-g\infty$ avec $E\geqslant 0$.\\

Tout diviseur effectif $D$ de degré $g$ sur $C$ s'écrit de façon unique  $D=D_1+D_2$, où le support de $D_1$ est constitué de points de $X$ au-dessus de points singuliers ou terribles de $C$ et $\deg(D_1)\leqslant g$.
Pour chaque $D_1$ (décrit explicitement comme $\sum_i a_iQ_i$), l'algorithme cherche les $D_2$ (décrits par des indéterminées) de degré $g-\deg(D_1)$ tels que $D_1+D_2-g\infty$ soit un diviseur de $n$-torsion. Soit donc $D=D_1+D_2$ un diviseur effectif de degré $g$ sur $X$. 
L'algorithme construit un diviseur très simple $D'$ (dont les coordonnées s'expriment en fonction de celles de $D_2$) équivalent à $n(D-g\infty)$, puis vérifie si $D'$ est principal. Il construit un polynôme homogène $G_\infty$ de façon à ce que $D'$ soit principal si et seulement s'il existe un numérateur $G_0$ tel que $D'=\div(G_0/G_\infty)$. Une condition nécessaire et suffisante d'existence d'un tel numérateur est donnée par des équations portant sur les coordonnées de $D_2$.\\

Ceci fournit un schéma affine $S_{D_1}$ dont chaque $k$-point paramètre un diviseur $D_2$ tel que $n(D_1+D_2)$ soit principal. Soit $S$ le coproduit des $S_{D_1}$. Alors chaque classe de diviseurs de $n$-torsion contient au moins un diviseur décrit par un point de $S$. Le morphisme $S\to J[n]$ qui à un diviseur associe sa classe est constant sur chaque composante irréductible ; il suffit donc de trouver un point dans chaque composante pour obtenir un diviseur de chaque classe.

\begin{rk}\label{rk:huangrac}
Cet algorithme peut également s'adapter, comme nous le verrons, à la recherche des racines $n$-ièmes dans $\Pic^0(X)$ d'un diviseur de degré zéro $F$ fixé : il convient de construire un diviseur très simple équivalent à $n (D-g\infty)-F$, où $D=D_1+D_2$ avec $D_2$ décrit par des indéterminées.
\end{rk}

\subsection{Construction d'un diviseur très simple}

\subsubsection{Pour un diviseur connu}

Voici comment déterminer un diviseur très simple équivalent à un diviseur donné $n(D-g\infty)$. Soit $P\in |D|$ : c'est un point de $X$, dont l'image dans $C$ est peut-être singulière. Soit $f\colon C'_P\to C$ un morphisme birationnel tel que $P$ corresponde au point régulier $(0,0)$ de $C'_P$ ; il se construit par transformée quadratique. Choisissons des éléments $r_1,\dots,r_n$ deux à deux distincts de $k$. Supposons donné un ensemble fini $S_P$ de points "à éviter" de $C'_P$. Soit $a_P\in k$ tel que les droites $L_i$ d'équation $y=(r_i+a_P)x$ vérifient : \begin{itemize}[label=$\bullet$]
\item Chaque $L_i$ intersecte $C'_P$ seulement en des points réguliers et non terribles dont l'image par $f$ a un unique antécédent
\item Chaque droite $L_i$ intersecte $C'_P$ en $\deg(C'_P)$ points affines
\item $L_i\cap  S_P=\varnothing$
\end{itemize}

Les deux premiers points assurent que le diviseur $\div(L_i)-(0,0)$ est très simple sur $C'_P$, et que son image dans $C$ est encore très simple. Le troisième point assure que le support de ce diviseur sur $C'_P$ ne contient aucun point de $S_P$. Un tel $a_P$ vérifie que le diviseur du produit $L_P\coloneqq L_1\cdots L_n$ sur la partie affine de $C'_P$ est égal à $nP+A_P$, où $A_P$ est un diviseur effectif très simple dont le support est disjoint de $S_P\cup \{ P\}$.

Pour construire un diviseur très simple équivalent à $n(D-g\infty)$ où $D=D_1+D_2$ comme précédemment, Huang et Ierardi commencent par appliquer cette procédure à $D_1$. Posons $T=|D_2|$ et pour un premier point $P\in |D_1|$, $S_P=(C'_P\to C)^{-1}T$. Ceci donne un produit $L_P$ tel que $\div(L_P)=nP+A_P\in \Div(C'_{P}\cap\A^2)$ avec $A_P$ effectif très simple. Pour chaque $P$, ils appliquent ensuite la procédure au point $\infty\in C'_P$, en remplaçant $S$ par $(S\cup |A_P|)-\{\infty\}$. Ceci fournit un polynôme $L'_P$ de diviseur $n\infty+B_P$. Le diviseur de la fonction rationnelle $\frac{L_P}{L'_P}$ est alors $nP-n\infty+A_P-B_P$. Le diviseur très simple $A_P-B_P$ est équivalent à $nP-n\infty$, et son support disjoint de celui de $D$. Pour passer au point $P$ suivant, il suffit de remplacer $T$ par $T\cup |A_P-B_P|$. Cette construction donne un diviseur $D_1'=D_1'^+ - D_1'^-$ équivalent à $n(D_1-\deg(D_1)\infty)$. \\

Pour $D_2$, la procédure est plus simple : nul besoin de courbes $C'_P$, puisque les points du support de $D_2$ sont déjà réguliers. L'ensemble $S$ considéré pour chaque $P\in |D_2|$ est $|D_1'|\cup|D_2|-|P|$. Pour le point $\infty$, la même procédure est appliquée avec $S=|D_1'|\cup|D_2|-\{\infty\}$. \\

\begin{rk}\label{rk:huangrac2}
Soit $F$ un diviseur de degré 0 sur $X$. Supposons qu'il est de la forme $D_F-g\infty$.
Afin de trouver un diviseur très simple équivalent à $n(D-g\infty)+F$, la procédure ci-dessus peut être appliquée à $D_F$ (en cherchant une seule droite pour chaque $P\in |D_F|$ et non pas $n$) pour le rendre très simple, puis à $D$ en ajoutant à chaque ensemble $S_P$ les points du support de $F$.
\end{rk}

\subsubsection{Pour un diviseur indéterminé}\label{subsubsec:divind}

Fixons un diviseur $D_1$ de degré $d_1\leqslant g$ dont le support ne contient que des points singuliers ou terribles de $C$. Notons $d_2=g-d_1$. L'algorithme cherchant les $D_2$ tels que $D_1+D_2-g\infty$ soit de $n$-torsion procède de la façon suivante. Soit $|D_2|$ un ensemble de $d_2$ éléments, qui représente moralement le support du diviseur $D_2$ cherché. Fixons deux familles $(r_{P,i})_{P\in |D_2|,1\leqslant i\leqslant n}\in k^{nd_2}$ et $(r_{P,\infty,i})_{P\in |D_2|,1\leqslant i\leqslant n}\in k^{nd_2}$ d'éléments de $k$ deux à deux distincts. Fixons également un ensemble $\mathcal{A}\subset k^{d_2+1}$ de familles $a=(a_P)_{P\in |D_2|\cup \{\infty\}}$. Ces deux ensembles serviront à paramétrer des droites passant par les points de $|D_2|$. Soit $a\in\mathcal{A}$. Soient $d_2$ triplets d'indéterminées $(x_P,y_P,z_P)_{P\in |D_2|}$, qui représentent les points de $D_2$. Nous noterons $k(D_2)$ le corps $k((x_P,y_P,z_P)_{P\in |D_2|})$.\\

Notons, pour tout $i\in\{1\dots n\}$ et tout $P\in |D_2|$, $L_{P,i}=(y-y_Pz)-(r_{P,i}+a_P)(x-x_Pz)\in k(D_2)[x,y,z]$, et $L_{P,\infty,i}=(y-y_\infty z)-(r_{P,\infty,i}+a_\infty)(x-x_\infty z)$. Supposons, quitte à changer $\mathcal{A}$, que pour tout $Q\in |D_1'|\cup \mathrm{Sing}(C)$, $L_{P,\infty,i}(x_Q,y_Q,z_Q)\neq 0$.
Considérons les conditions suivantes sur ces indéterminées :
\begin{itemize}[label=$\bullet$]
\item pour tout point singulier ou terrible $Q=(x_Q,y_Q,z_Q)$ de $C$, pour tout $i\in \{1\dots n\}$, pour tout $P\in |D_2|$, $L_{P,i}(x_Q,y_Q,z_Q)\neq 0$ ;
\item pour tout point $Q\in |D_1'|$ (déjà calculé explicitement), pour tout $i\in \{1\dots n\}$, pour tout $P\in |D_2|$, $L_{P,i}(x_Q,y_Q,z_Q)\neq 0$ ;
\item pour tout $P\in |D_2|$, pour tout $i\in\{1\dots n\}$, $L_{P,\infty,i}(x_P,y_P,z_P)\neq 0$ ;
\item pour tous $i,j\in \{1\dots n\}$, pour tous $P\neq P'\in |D_2|$,  le point d'intersection des droites d'équations $L_{P,i}$ et $L_{P',j}$ n'est pas un zéro du polynôme $f$ qui définit $C$ ;
\item  pour tout $i\in \{1,\dots,n\}$ et tout $P\in |D_2|$, la droite d'équation $L_{P,i}$ n'est pas tangente à $C$. 
\end{itemize}
La dernière condition se vérifie de la façon suivante : les abscisses des points d'intersection de la droite $L_{P,i}$ avec $C$ dans le plan affine $\A^2=\Spec k(D_2)[x,y]$ sont les racines du polynôme $f(x,y_P+(r_{P,i}+a_P)(x-x_P))\in k(D_2)[x]$. Les points $P$ pour lesquels $L_{P,i}$ est tangente à $C$ en un point sont ceux tels que le discriminant de ce polynôme s'annule. La condition est donc $\mathrm{discr}(f(x,y_P+(r_{P,i}+a_P)(x-x_P))\neq 0$. Une condition semblable est imposée pour l'ouvert affine $\Spec k(D_2)[x,z]$.

Ces inéquations définissent, pour chaque $a\in\mathcal{A}$, un ouvert $S_a$ de $\A_k^{3d_2}$ qui paramètre les diviseurs $D_2$ tels que pour ce $a$ fixé, les droites $L_{P,i}$ d'équation $y-y_Pz=(r_{P,i}+a_P)(x-x_Pz)$ vérifient les trois conditions mentionnées précédemment, et permettent donc de remplacer $n(D-g\infty)$ par un diviseur très simple qui lui est équivalent. En effet, soient $L_0=\prod_{P,i}L_{P,i}$ et $L_\infty=\prod_{P,i}L_{P,\infty,i}$ : alors le diviseur  $\div(\frac{L_0}{L_\infty})-n(D_2-d_2\infty)$ est un diviseur très simple $B=B^+-B^-$. Par conséquent, $A=B+D_1'$ est linéairement équivalent à $n(D-g\infty)$. Notons $A=A^+-A^-$. Alors $A^+=B^++D_1'^+$, et $A^-=B^-+D_1'^-$. 

\begin{rk}
Si l'ensemble $\mathcal{A}$ est assez grand (de taille $O(n^2g^2d)$, voir \cite[4.2, p12]{huang_counting}), alors chaque diviseur de $\Div^0X$ est paramétré par un point de l'un des $S_a,a\in\mathcal{A}$.
\end{rk}

\begin{rk}\label{rk:huangrac3}
La méthode décrite dans cette section fonctionne tout aussi bien, étant donné un diviseur $F$ de degré nul, pour remplacer $n(D_1+D_2-g\infty)+F$ par un diviseur très simple. Construisons d'abord $F'$ très simple équivalent à $F$, comme décrit précédemment. Ajoutons $|F'|$ à l'ensemble $S$ utilisé pour construire les équations. Finalement, ajoutons $F'$ au diviseur $A$. Ceci donne un diviseur très simple équivalent à $n(D-g\infty)+F$. Cette construction ajoute $O(n)$ équations et inéquations et n'altère donc pas la complexité globale de l'algorithme.
\end{rk}

\subsection{Construction d'une fonction de diviseur très simple}

La partie précédente a permis de déterminer une fonction $\frac{L_0}{L_\infty}$ telle que $A=\div(\frac{L_0}{L_\infty})-n(D-g\infty)$ soit très simple. Il reste à tester s'il existe une fonction $\frac{G_0}{G_\infty}$ dont $A$ est le diviseur. 
Pour ce faire, la procédure classique de l'algorithme de Brill-Noether commence par construire un dénominateur $G_\infty$ (explicitement, en fonction des indéterminées $x_P,y_P,z_P$ pour $P\in |D_2|$) qui vérifie que $A$ est principal si et seulement s'il existe un numérateur $G_0$ tel que $\div(\frac{G_0}{G_\infty})=A$. Des conditions nécessaires et suffisantes pour l'existence de ce numérateur sont ensuite déterminées. Rappelons que $Q_1,\dots,Q_r$ sont les points singuliers de $C$. Notons $m_1,\dots,m_r$ leurs multiplicités respectives. Pour chaque $Q_i$, notons $Q_{i,1},\dots,Q_{i,m_i}$ les points de $X$ au-dessus de $Q_i$. Définissons pour la suite le diviseur 
\[ E=\sum_{i=1}^r\sum_{j=1}^{m_i}(m_i-1)Q_{ij}.\]

\subsubsection{Construction du dénominateur}

Soit $H\in k[x,y,z]$ un polynôme homogène tel que $\div(H)-E$ soit effectif et très simple. Un tel polynôme se construit explicitement en choisissant, pour chaque point singulier $Q_i$, $m_i-1$ droites passant par $Q_i$ avec multiplicité $m_i$ qui intersectent $C$ en $d-m_i$ autres points, et en considérant le produit de toutes les formes linéaires définissant ces droites. Pour les détails, voir \cite[§4.3,p13]{huang_counting}. Notons $A_e=\div(H)-E$. Supposons avoir calculé $H$ avant la section précédente, et avoir ajouté $|\div(H)|$ à chaque ensemble $S$ de points à éviter dans la section précédente. Par construction, le degré de $H$ est inférieur à $\deg(E)d\leqslant d^4$ : ceci n'altère pas la complexité globale de l'algorithme, et permet d'assurer que $|\div(H)|$ soit disjoint de $|A|$.
Supposons fixés les $a_P,r_{P,i}\in k$, les droites $L_{P,i}\in k(D_2)[x,y,z]$ et enfin les polynômes $L_0,L_\infty \in k(D_2)[x,y,z]$ comme ci-dessus. Rappelons que $\div(L_\infty)=B^-+nd_2\infty$, où $B^-$ est la partie négative du diviseur $\div(\frac{L_0}{L_\infty})-n(D_2-d_2\infty)$. \\

Soit $U\in k(D_2)[u_x,u_y,u_z]$ le $u$-résultant de $f$ et $L_\infty$ (voir annexe \ref{subsec:repdiv}). Il vérifie  \[ U=\prod_{P\in |B^-|\cup \{\infty\}}(x_Pu_x+y_Pu_y+z_Pu_z)=(x_\infty u_x+y_\infty u_y+z_\infty u_z)^{nd_2}R_B(u_x,u_y,u_z).\]
Comme les coordonnées du point $\infty$ sont connues, il est possible de calculer explicitement $R_B$, et en particulier $R_B(-z,0,x)=\prod_{P\in |B^-|}(x-x_Pz)$. La même procédure appliquée à $D_1'^-$ fournit $R_D(-z,0,x)=\prod_{P\in |D_1'^-|}(x-x_Pz)$. Posons $R(x,y,z)=R_B(-z,0,x)R_D(-z,0,x)$.
Supposons que $\div(R)$ soit très simple et que son support soit disjoint de $A^++A_e$, et définissons $G_\infty\coloneqq HR$. Alors le diviseur de $G_\infty$ est supérieur à $A^-+E+A_e$, et s'écrit \[ \div(G_\infty)=A^-+E+A_e+M\] avec $M+A_e$ très simple.\\

\begin{prop}\cite[Th. 4.1]{huang_counting}
Soient $D,D'\in\Div(C)$ deux diviseurs. Soit $G_\infty\in k[x,y,z]$ un polynôme homogène tel que $\div(G_\infty)\geqslant E$. Notons $F=\div(G_\infty)-D-E$. Alors $D$ est équivalent à $D'$ si et seulement s'il existe un polynôme homogène $G_0\in k[x,y,z]$ tel que $\div(G_0)\geqslant E$ et $\div(G_0)=D'+E+F$.
\end{prop}

Rappelons que dans notre situation, $A=A^+-A^-$ est un diviseur équivalent à $n(D-g\infty)$, et le polynôme $G_\infty $ vérifie $\div(G_\infty)=A^-+E+A_e+M$. Le diviseur $n(D-g\infty)$ est principal si et seulement si $A^+$ est équivalent à $A^-$. 
La proposition affirme, en prenant $F=A_e+M$, $D=A^+$  et $D'=A^-$, que $n(D-g\infty)$ est principal si et seulement s'il existe un polynôme homogène $G_0\in k[x,y,z]$ tel que $\div(G_0)\geqslant E$ et $\div(G_0)=A^++E+A_e+M$.

\begin{rk}\label{rk:phi}
Ci-dessus, le diviseur de $R$ était supposé très simple. Il peut toujours être rendu très simple par un automorphisme de $\PP^2$ dont la restriction à $\A^2=\Spec k[x,y]$ est de la forme $(x,y)\mapsto (ax+by,y)$. Un tel automorphisme se trouve toujours en considérant un ensemble $\Phi$ d'automorphismes assez grand (de taille $O(dgn)$, cf \cite[4.3, p14]{huang_counting}).
\end{rk}

\subsubsection{Existence d'un numérateur}

Il reste à donner une condition nécessaire et suffisante sur $(x_P,y_P,z_P)_{P\in |D_2|}$ pour qu'il existe un polynôme homogène $G_0\in k[x,y,z]$ de même degré que $G_\infty$ tel que $\div(G_0)\geqslant E$ et \[ \div(G_0)=A^++E+A_e+M. \] Le polynôme $G_0$ est défini par ses coordonnées dans la base canonique du $k$-espace vectoriel des polynômes homogènes de degré $\deg G_\infty$. Il suffit, au vu de la contrainte sur son degré, de montrer que $G_0$ est supérieur à la fois à $E$, à $A^+$ et à $A_e+M$.

\paragraph{Condition 1 : $\div(G_0)\geqslant E$} Cette condition est équivalente à ce que la multiplicité de $G_0$ en chaque point singulier $Q_i$ soit supérieure à $m_i-1$ \cite[Lem. 4.2]{huang_counting}, ce qui se décrit par l'annulation de certaines formes linéaires en les coefficients de $G_0$ \cite[Lem. 4.3]{huang_counting}. 

\paragraph{Condition 2 : $\div(G_0)\geqslant A^+$ } Souvenons-nous que $A^+=B^++D_1'^+$, où $D_1'$ est un diviseur explicitement connu, et $B^+=\div(L_0)-nD_2$ est un diviseur indéterminé, où $L_0$ est un produit de polynômes homogènes de degré 1 dans $k(D_2)[x,y,z]$.
L'annulation de $G_0$ en les points du support de $D_1'^+$ se traduit par des équations linéaires en ses coefficients. Soit maintenant $L_{P,i}$ une des formes linéaires qui divisent $L_0$ ; pour simplifier l'exposition, supposons qu'il s'agit de la droite $y=mx$. Le diviseur de $L_{P,i}$ a $\deg(f)$ points sur la partie affine de $C$. Alors le diviseur de $G_0$ est supérieur à $\div(L_{P,i})-P$ si et seulement si les polynômes $G_0(t,mt,1)$ et $f(t,mt,1)$ ont $\deg(f)-1$ zéros en commun, c'est-à-dire s'ils ont un facteur commun de degré $d-1$. Ceci revient à demander qu'il existe un polynôme $h(t)\in k(D_2)[t]$ de degré $\deg(G_0)-\deg(f)+1$ tel que $tG_0(t,mt,1)-h(t)f(t,mt,1)=0$, ce qui se traduit par un système linéaire en les coefficients de $G_0$. 

\paragraph{Condition 3 : $\div(G_0)\geqslant A_e+M$} Un paramétrage rationnel du support de $M$ (resp. $A_e$) permet de s'assurer que $G_0$ s'annule en chaque point de $M$ (resp. $A_e$). Il s'obtient en construisant un polynôme $Q\in k(D_2)[t]$ et des fonctions rationnelles $r,s\in k(D_2)(t)$ tels que les $\overline{k(D_2)}$-points de $M$ (resp. $A_e$) soient les $(r(\theta),s(\theta))$ où $\theta\in\overline{k(D_2)}$ parcourt les racines de $Q$. Pour les détails, voir \cite[Lem. 4.4]{huang_counting} et \cite[Lem. 2.2]{canny}. Alors $G_0$ s'annule sur le support de $M$  si et seulement si $G_0 (r(t),s(t))$ s'annule sur les zéros de $Q$. Réduisons les fractions $r$ et $s$ au même dénominateur, et considérons leurs numérateurs respectifs $r',s'$. Notons $G_0'(t)=G_0\circ(r',s')\in k(D_2)[t]$. La fonction $G_0$ s'annule sur $|M|$ si et seulement s'il existe $I'\in k(D_2)[t]$ de degré $\deg(G_0)-\deg(Q)$ tel que $G_0'=I'Q$. Ceci se traduit encore par un système linéaire en les coefficients de $I'$ et $G_0$. De même, la condition $\div(G_0)\geqslant A_e$ se traduit par l'existence d'un $I''$ solution d'un certain système linéaire.

\paragraph{Résumé} La mise bout à bout de ces systèmes linéaires fournit une matrice $T$ à coefficients dans $k(D_2)$ telle que $G_0$ soit un numérateur convenable si et seulement si $(G_0,I',I'')$ est un élément non trivial du noyau de $T$. Par conséquent, le diviseur $D_2$ vérifie que $n(D_1+D_2-g\infty)$ est principal si et seulement si le rang de $T$ n'est pas maximal, ce qui se traduit par des équations en les $(x_P,y_P,z_P)_{P\in |D_2|}$.\\

En résumé, l'algorithme détermine pour chaque diviseur $D_1$ de degré inférieur à $g$ de support inclus dans ${\rm Sing}(C)$, chaque transformation projective $\phi\in\Phi$ (voir remarque \ref{rk:phi}) et chaque $a\in \mathcal{A}$ (voir début de la section \ref{subsubsec:divind}), une partie constructible $S_{D_1,\phi,a}$ de l'espace affine $\A^{3(g-\deg D_1)}=\Spec k[x_P,y_P,z_P ; P\in |D_2|]$ et un morphisme \begin{align*} S_{D_1,\phi,a}&\to J_X[n] \\ D_2&\mapsto D_1+D_2-g\infty.\end{align*} 
Chaque diviseur de $n$-torsion est équivalent à un diviseur de la forme $D_1+D_2-g\infty$, où $D_2$ appartient à l'un des $S_{D_1}$. Par conséquent, en notant $S=\bigsqcup_{D_1,\phi,a}S_{D_1,\phi,a}$, il y a un morphisme surjectif $f\colon S\to J_X[n]$ qui s'insère dans le diagramme commutatif
\[
\begin{tikzcd}
S \arrow[r]\arrow[d,"f",swap] & C^g \arrow[r] & C^{(g)} \arrow[d,"D\mapsto D-g\infty"]\\
J_X[n] \arrow[rr] & & J_X
\end{tikzcd}
\]

\subsection{Détermination de $J_X[n]$}

\subsubsection{Détermination de représentants des classes de diviseurs}\label{subsubsec:detrep}

Pour chaque diviseur $D_1$ de degré $g$, chaque transformation projective $\phi\in\Phi$ (voir remarque \ref{rk:phi}) et chaque $a\in \mathcal{A}$ (voir début de la section \ref{subsubsec:divind}), il y a une application $S_{D_1,\phi,a}\to J_X[n]$, $D_2\mapsto D_1+D_2-g\infty$. Chaque diviseur de $n$-torsion est équivalent à un diviseur de la forme $D_1+D_2$, avec $D_2\in S_{D_1}$. Ceci définit une application surjective $f\colon \bigsqcup_{D_1}S_{D_1}\to J_X[n]$.

\begin{lem}
Soit $f\colon X\to Y$ une application continue entre espaces topologiques. Si $Y$ est discret alors $f$ est constante sur chaque composante connexe.
\begin{proof}
L'image d'une partie connexe par une application continue est connexe. Comme $Y$ est discret, ses composantes connexes sont des points.
\end{proof}
\end{lem}

L'application $f\colon S=\bigcup_{D_1,\phi,a}S_{D_1,\phi,a}\to J_X[n]$ est donc surjective et constante sur chaque composante irréductible. Par conséquent, pour trouver tous les éléments de $J_X[n]$, il suffit de trouver au moins un point de chaque composante irréductible de $S$.

Dans la suite, la notation $\mathcal{P}(d,g,n)$ sera employée génériquement pour remplacer des polynômes en $d,g,n$ à coefficients réels : le polynôme exact caché derrière la notation $\mathcal{P}$ pourra varier. Chaque $S_{D_1,\phi,a}\subset \A^{3g}$ est défini par un nombre d'équations et d'inéquations polynomial en $d,g,n$ ; le degré de chacune de ces (in)équations est également polynomial en $d,g,n$. Par conséquent, la complexité totale du calcul d'un élément de chaque composante irréductible à l'aide de l'algorithme décrit dans l'annexe \ref{subsubsec:ptscompconstr} est $\mathcal{P}(d,g,n)^{O(g^3)}+3g \Fact(\mathcal{P}(d,g,n))$, où $\Fact$ désigne la complexité de la factorisation absolue d'un polynôme de $k_0[t]$ de degré donné. De plus, le nombre de points trouvés est $\mathcal{P}(d,g,n)^{O(g^2)}$, et chaque point est défini sur une extension de $k_0$ de degré $n^{O(g^3)}$.
Enfin, le nombre de valeurs de $\Phi$ et de $a$ à considérer est polynomial en $d,g,n$, et le nombre de diviseurs $D_1$ à considérer est $O(n^{2g})$. Au total, la complexité de cet algorithme est donc $\mathcal{P}(d,g,n)^{O(g^3)}+\mathcal{P}(d,g,n) \Fact(k_0,\mathcal{P}(d,g,n))$.

\subsubsection{Calcul de $J_X[n]$ et de sa loi de groupe}\label{subsubsec:calculJx}

L'algorithme précédent a permis de trouver un ensemble $T$ de $\mathcal{P}(d,g,n)^{O(g^2)}$ diviseurs, tel que chaque élément de $\Pic^0(X)[n]$ soit la classe d'un diviseur de $T$. L'algorithme choisit $D^1\in T$ et pose $J=\{ D^1\}$. Puis, pour chaque diviseur $D\in T$, il vérifie s'il est équivalent à un diviseur $D'\in J$ ; si non, il remplace $J$ par $J\cup \{ D'\}$. Au total, il vérifie pour $|T|$ diviseurs s'ils sont équivalents à $O(n^{2g})$ diviseurs. Il y a donc $\mathcal{P}(d,g,n)^{O(g^2)}$ vérifications à faire. Chacune de ces vérifications nécessite seulement de construire une extension de corps sur laquelle sont définis les deux diviseurs. \\

Afin de vérifier si deux diviseurs de degré zéro $D,D'$ sont équivalents, il suffit de calculer l'espace de Riemann-Roch $\mathcal{L}(D-D')=\HH^0(X,\OO_X(D-D'))$ et de tester s'il est nul. Comme les diviseurs en question sont de la forme $D-g\infty$, ce calcul se fait en temps $O(g^7d^{14})$ (voir annexe \ref{subsec:rr}). Ce processus fournit exactement un représentant de chaque classe de $J_X[n]$ en temps $\mathcal{P}(d,g,n)^{O(g^2)}$.

Notons $D^1,\dots,D^{n^{2g}}$ les diviseurs obtenus. L'addition dans $J$ de deux classes de diviseurs $D^i$ et $D^j$ se fait en calculant le diviseur $D^i+D^{j}$ puis en vérifiant pour chaque indice $\ell$ si $D^i+D^j$ est équivalent à $D^\ell$, c'est-à-dire si $\mathcal{L}(D^i+D^j-D^\ell)$ contient un élément non nul. La complexité totale de cette opération est $O(n^{4g}g^7d^{14})$.

\subsection{Complexité sur $\FF_q(t)$}

La construction du schéma $S$ paramétrant les diviseurs de $n$-torsion est parfaitement indépendante du corps de base $k_0$, et peut être réalisée sur n'importe quel corps calculable. Le calcul des points dans les composantes connexes demande simplement de disposer d'un algorithme de factorisation absolue des polynômes à coefficients dans $k_0$.\\

Intéressons-nous à la complexité de cet algorithme dans le cas où $k_0=\FF_q(t)$ et les points $Q_{ij}$ au-dessus des points singuliers ou terribles du modèle plan de $X$ sont tous définis sur $\FF_q(t)$. C'est une restriction très forte, mais elle permet d'obtenir rapidement une évaluation de la complexité de l'algorithme. \'{E}tant donné $f=\frac{a}{b}\in k_0$ avec $a,b\in \FF_q[t]$ premiers entre eux, définissons sa hauteur par $h(x)=\max (\deg_t a,\deg_t b)$. Notons $D$ la hauteur maximale des coefficients de l'équation $F$ de la courbe $C$, et des points $Q_{ij}$. Une transformation quadratique ne change pas la hauteur des coefficients. En supposant toujours $q$ assez grand ($d^6<q$ et $n=O(g\log q)$), il est possible de choisir pour l'ensemble $\mathcal{A}$ une partie de $\FF_q$. 
Les coefficients des équations de $S$ sont obtenus par produit d'un nombre polynomial en $d,g,n$ d'éléments de $\mathcal{A}$, de coordonnées des points $Q_{ij}$ et de coefficients de $F$ ; la hauteur de ces coefficients est alors polynomiale en $D,d,g,n$. Par conséquent, le calcul d'éléments de chaque composante irréductible de $S$ à l'aide de l'algorithme probabiliste de factorisation dans $\FF_q(t)[s]$ présenté dans l'annexe \ref{susubsec:factmult} nécessite $\mathcal{P}(D,d,g,n)^{O(g^3)}+\mathcal{P}(D,d,g,n)\log(q)$ opérations.

\section{L'algorithme de Jin}\label{sec:jin}

\subsection{Données et structure de l'algorithme}

Cet algorithme calcule les groupes de cohomologie d'un faisceau lisse sur une courbe lisse sur un corps algébriquement clos. C'est le seul algorithme existant qui effectue cette tâche pour les faisceaux lisses en un nombre d'opérations borné explicitement. Cependant, il n'est pas utilisable dans la pratique en raison du grand nombre de variables en jeu.

\subsubsection{Données}

Soient $k_0$ un corps parfait et $k$ une clôture algébrique de $k_0$.
Soit $X$ une courbe projective lisse sur $k_0$. Elle est représentée comme décrit dans l'annexe \ref{subsub:OP1alg}, c'est-à-dire par une $\OP$-algèbre $\mathcal{E}$ telle qu'il existe un morphisme $\phi \colon X\to\PP^1$ avec $\phi_\star\OO_X=\mathcal{E}$. Cette algèbre est définie par des entiers $a_1,\dots,a_r$ tels que $\mathcal{E}\simeq\OP(a_1)\oplus\dots\oplus\OP(a_r)$, la matrice $M\in \Mat_{r\times r^2}(k_0[x,y])$ qui décrit la multiplication sur $\mathcal{E}$, et la matrice $I\in\Mat_{1\times r}(k_0[x,y])$ qui décrit le morphisme structural $\OP\to\mathcal{E}$.

Le faisceau lisse $\F$ est décrit par la donnée d'un revêtement galoisien $f\colon Y\to X$ (lui-même représenté comme $\OP$-algèbre) qui le trivialise, du groupe $G$ du revêtement, et du $G$-module $F\coloneqq \HH^0(Y,f^\star\F)$. D'après le corollaire \ref{cor:torsrev}, la catégorie des $\F$-torseurs sur $X$ est équivalente à celle des $F$-torseurs sur $Y$ munis d'une action de $G$ telle que l'action de $F$ soit $G$-équivariante.

Les $F$-torseurs $T\to Y$ cherchés sont décrits de la même manière, par des entiers $b_1,\dots,b_s$, une matrice de multiplication $N\in\Mat_{s\times s^2}(k_0[x,y])$ et une matrice $J\in\Mat_{1\times s}(k_0[x,y])$ avec quelques données supplémentaires. Le morphisme $T\to Y$ est donné par une matrice $S\in \Mat_{s\times r}(k_0[x,y])$. L'action du groupe $G$ est décrite par des matrices $\Phi_g\in\Mat_{s\times s}(k_0[x,y])$, où $g$ parcourt $G$. De même, l'action du groupe $F$ est décrite par des matrices $\Psi_f\in\Mat_{s\times s}(k_0[x,y])$, où $f$ parcourt $F$. Ces matrices vérifient des égalités traduisant le fait qu'elles définissent une action, ainsi que la $G$-équivariance de l'action de $F$ (voir section \ref{subsec:doneq}).

\subsubsection{Résumé de l'algorithme}
\label{subsubsec:resalg}

L'algorithme construit un schéma paramétrant tous les $\F$-torseurs sur $X$, de façon à ce que les composantes connexes de ce schéma soient en bijection avec les classes d'isomorphisme de $\F$-torseurs sur $X$.

Il commence par sélectionner un nombre fini de $b=(b_i)$ possibles pour les torseurs $T$ envisagés, déterminés par les conditions $b_i\leqslant 0$ \cite[Lem. 6.17]{jinbi_jin} et $\sum_i b_i=|G|\cdot \sum_j a_j$ (voir le corollaire \ref{cor:det}).

Pour chacun de ces $b$, l'algorithme écrit les équations que la $\OP$-algèbre $\OO_T$, munie des actions de $F$ et $G$, vérifie si et seulement si $T\to Y$ est un $F$-torseur $G$-équivariant. Ceci donne un fermé $\mathcal{U}_b$ d'un espace affine sur $k_0$. Il se trouve que les composantes connexes de $\mathcal{U}_b$ sont irréductibles, et en bijection $\Gal(k|k_0)$-équivariante avec les classes d'isomorphisme de $F$-torseurs sur $Y$ (voir la proposition \ref{thprinc}). Il suffit alors de trouver un point dans chaque composante irréductible de chaque $\mathcal{U}_b$ afin d'obtenir $\HH^1(X,\F)$ comme $\Gal(k|k_0)$-ensemble.

\subsection{Conditions pour être un torseur}\label{subsec:condtors}

Rappelons qu'un $F$-torseur sur un schéma $Y$ est un $Y$-schéma étale $T$ tel que le morphisme $F\times T\to T\times T$ soit un isomorphisme. 
Les morphismes $T\to Y$ considérés sont toujours plats, c'est-à-dire localement libres. L'étalitude de $T\to Y$ se vérifie en deux étapes : tout d'abord, une condition nécessaire pour qu'un morphisme localement libre soit étale est que son rang soit constant. Ensuite, une condition nécessaire et suffisante portant sur les morphismes de rang constant garantit l'étalitude.

\subsubsection{Rang constant}

\begin{lem}\cite[Lem. 3.49]{jinbi_jin_thesis} Soit $S$ un schéma. Soit $X$ un $\PP^1_S$-schéma fini localement libre, et lisse sur $S$. Soit $T$ un $X$-schéma fini localement libre sur $\PP^1_S$. Alors $T$ est un $X$-schéma fini localement libre.
\end{lem}

\begin{lem}\label{surjcon}
Soient $S$ un schéma, et $f\colon X\to \PP^1_S$ un morphisme fini localement libre de $S$-schémas. Alors la restriction de $f$ à toute composante connexe de $X$ est surjective.
\begin{proof} Ceci se vérifie fibre à fibre sur $S$ : il suffit de traiter le cas où $S$ est le spectre d'un corps $K$. Soit $Y$ une composante connexe de $X$. Le morphisme $f\vert_Y\colon Y\to\PP^1_K$ est fini localement libre, donc $\dim_K Y=\dim_K X=1$. De même, par finitude de $f$, $\dim f(Y)=\dim Y=1$. Enfin, comme $f$ est fini, il est fermé, donc $f(Y)$ est un fermé de $\PP^1_K$ de dimension 1 : $f(Y)=\PP^1_K$. Par conséquent, $f\vert_Y$ est surjectif.
\end{proof}
\end{lem}

\begin{lem} Soient $S$ un schéma et $f\colon T\to X$, $g\colon X\to \PP^1_S$ deux morphismes de $S$-schémas finis localement libres. Supposons que le morphisme $f_0\colon T_0\to X_0$ entre les fibres au-dessus de $0\in\PP^1_S(S)$ est fini localement libre de rang constant $r$. Alors $f$ est fini localement libre de rang constant $r$.
\begin{proof} Soit $Y$ une composante connexe de $X$. Alors $f\vert_Y\colon T|_Y Y \to Y$ est localement libre de rang constant un entier $d$. Montrons que $d=r$. 
Notons $Y_0=Y\times_{\PP^1_S}0$ : c'est un schéma non vide par le lemme précédent. Alors $T_0\times_XY_0$ peut s'écrire comme le tiré en arrière de $T_0\to X_0$ par $Y\to X$, mais aussi comme le tiré en arrière de $f\vert_Y\colon T\vert_Y\to X$ par $0\to \PP^1_S$, ce qui implique que $r=d$.
\end{proof}
\end{lem}

Les morphismes $T\to Y$ qui sont de rang constant $m=|G|$ sont donc définis par la donnée supplémentaire d'un isomorphisme de $\OO_{Y_0}$-modules $\OO_{T_0}\xrightarrow{\sim}\OO_{Y_0}^m$. L'évaluation de la matrice de multiplication de $\OO_T$, à coefficients dans $k_0[x,y]$, en $x=0$ et $y=1$ permet de déduire $\OO_{T_0}$ de $\OO_T$. En particulier, ceci fournit la matrice $N_0$ de la multiplication sur le $k_0$-espace vectoriel $\OO_{T_0}$. La matrice $M_0$ de multiplication sur $\OO_{Y_0}$ se détermine de la même façon. La structure de $\OO_{Y_0}$-module sur $\OO_{T_0}$ est obtenue en évaluant de même la matrice du morphisme $\OO_Y\to\OO_T$. Un isomorphisme de $\OO_{X_0}$-modules $\OO_{T_0}\xrightarrow{\sim}\OO_{Y_0}^m$ est un isomorphisme de $k_0$-espaces vectoriels compatible à l'action de $\OO_{Y_0}$. Il est donc défini par une matrice $B\in\Mat_{sm\times sm}(k_0)$ à coefficients dans $k_0$, où $s=\deg(Y\to\PP^1)$. Notons $F_0$ la matrice du morphisme $\OO_{Y_0}\to\OO_{T_0}$, obtenue en évaluant la matrice de $\OO_X\to\OO_T$ en $x=0$ et $y=1$. La condition de $\OO_{Y_0}$-linéarité est la commutativité du diagramme suivant.
\[
\begin{tikzcd}
\OO_{Y_0}\otimes_{k_0}\OO_{Y_0}^m \arrow[d,"I_n\otimes M_0",swap]\arrow[r,"I_s\otimes B"] & \OO_{Y_0}\otimes_{k_0} \OO_{T_0} \arrow[d,"N_0(F_0\otimes B)"] \\
\OO_{Y_0}^m\arrow[r,"B",swap] & \OO_{T_0}
\end{tikzcd}
\]

\subsubsection{Discriminant}

Rappelons que le déterminant d'un module localement libre $\mathcal{E}$ de rang $r$ sur un schéma $X$ est le faisceau inversible $\wedge^r \mathcal{E}$. Un morphisme de modules localement libres de même rang $E\to E'$ induit un morphisme $\det \mathcal{E}\to\det \mathcal{E}'$, qui est localement la multiplication par une section de $\OO_X$.

\begin{df}
Soit $f\colon Y\to X$ un morphisme fini localement libre. Le morphisme composé $f_\star\OO_Y\otimes_{\OO_X}f_\star\OO_Y\to f_\star\OO_Y\xrightarrow{Tr}\OO_X$ définit un morphisme $\phi\colon f_\star\OO_Y\to \Hom_{\OO_X}(f_\star\OO_Y,\OO_X)$. Le déterminant de $\phi$ est un faisceau d'idéaux principaux sur $X$ appelé discriminant de $f$ et noté $\Delta_f$.  
\end{df}

\begin{prop}\cite[0BVH, 49.3.1]{stacks} Un morphisme fini localement libre $f\colon Y\to X$ de schémas est étale si et seulement si $\Delta_f \simeq \OO_X$.
\end{prop}

\begin{prop}\cite[Cor. 6.8]{jinbi_jin} Soient $g\colon T\to Y$, $f\colon Y\to W$ des morphismes finis localement libres de rang constant. Notons $r$ le rang de $g$. Alors $g$ est étale si et seulement si $\det_{\OO_W}\OO_T\simeq (\det_{\OO_W}\OO_Y)^{\otimes r}$ et $\Delta_{f\circ g}$ et $\Delta_f^{\otimes r}$ diffèrent d'un élément de $\OO_W^\times$.
\end{prop}

Dans la situation qui nous concerne, $W=\PP^1_S$, le schéma $Y$ est une courbe projective lisse, et $T$ est un $Y$-schéma dont l'étalitude est à vérifier.
Voici comment calculer le discriminant $\Delta_f$ d'un morphisme $f\colon Y\to\PP^1_{k_0}$, représenté par la $\OP$-algèbre $\mathcal{E}=f_\star\OO_Y$. C'est le déterminant de la trace de la multiplication $\mathcal{E}\otimes_{\OP}\mathcal{E}\to\mathcal{E}$. Pour le calculer, plaçons-nous sur un ouvert de $\PP^1$, disons $U_0=\Spec k_0[x]$. Alors $\mathcal{E}(U_0)$ est un $k_0[x]$-module libre, et la matrice de la multiplication sur $U_0$ dans une base $(e_1,\dots,e_s)$ est fournie. 
Pour chaque couple $(i,j)\in \{1\dots s\}^2$, soit $t_{ij}$ la trace de la matrice à coefficients dans $k_0[x]$ dans la base $(e_1,\dots,e_s)$ de la multiplication par $e_ie_j$ dans $\mathcal{E}(U_0)$. Le déterminant $\det(\mathcal{E})(U_0)$ est le déterminant de la matrice $(t_{ij})_{ij}$ est un élément de $k_0[x]$. La même procédure s'applique à l'ouvert $U_1$.
Algorithmiquement, il suffit de calculer une seule matrice à coefficients dans $k_0[x,y]$, puis d'évaluer son déterminant successivement en $x=1$ et $y=1$. Un calcul direct donne le résultat suivant.

\begin{lem}\label{lem:discr} Soit $N\in \Mat_{s\times s}(k_0[x,y])$ la matrice de la multiplication sur $\mathcal{E}$ avec les notations ci-dessus. Alors 
\[ {\rm Tr}(e_ie_j\cdot)=\sum_{\alpha,\beta=1}^sN_{\alpha, (i-1)s+j}N_{\beta,(\alpha-1)s+j}.\]
\end{lem}

\begin{ex} Considérons la courbe elliptique $E=\Proj k[X,Y,Z]/(X^2Z-Y^3+YZ^2)$ munie du morphisme $f\colon (X:Y:Z)\mapsto (X:Z)$ vers $\PP^1=\Proj k[X,Z]$. Soit $U_0$ l'ouvert $\Spec k[x]$ de $\PP^1$. Alors $f_\star \OO_E(U_0)=k[x][y]/(x^2-y^3+y)$, de $k[x]$-base $1,y,y^2$. 
Calculons par exemple la matrice de la multiplication par $y\cdot y=y^2$. On a $y^2\cdot 1=y^2, y^2\cdot y=x^2+y, y^2\cdot y^2=x^2y+y^2$. La matrice est donc \[ \left(\begin{matrix}
0 & x^2 & 0 \\
0 & 1 & x^2 \\
1 & 0 & 1
\end{matrix}\right)\]
et sa trace vaut 2. De même, la matrice la multiplication par $y^3$ est \[\left(\begin{matrix}
x^2 & 0 & x^2 \\
1 & x^2 & 1 \\
0 & 1 & x^2
\end{matrix}\right) \]
car $y^5=y^2(x^2+y)=x^2+y+x^2y^2$, et sa trace vaut $3x^2$.
Après avoir fait tous les calculs, on obtient
la matrice $({\rm tr} (y^iy^j\cdot ))_{0\leqslant i,j\leqslant 2}$ : \[ \left(\begin{matrix}
3 & 0 & 2 \\
0 & 2 & 3x^2 \\
2 & 3x^2 & 2
\end{matrix}\right).\]
Le déterminant de cette matrice est $4-9x^4$. Par conséquent, $\Delta_{E\to\PP^1}(U_0)=4-9x^4$, et le morphisme $E\to\PP^1$ n'est pas étale.
\end{ex}

Le corollaire suivant sert à borner le nombre de $\OP$-algèbres à considérer au début de l'algorithme décrit en \ref{subsubsec:resalg}.

\begin{cor}\label{cor:det}
Soit $S$ un schéma. Soit $T\to Y$ est un morphisme fini localement libre de rang $m$ de $\OO_{\PP^1_S}$-modules localement libres. Notons $(T\to\PP^1)_\star\OO_T\simeq \bigoplus_{j=1}^s\OP(b_j)$ et $(Y\to\PP^1_S)_\star\OO_Y\simeq\bigoplus_{i=1}^r\OP(a_i)$. Si $f$ est étale alors $\sum_j b_j=m\sum_i a_i$.
\begin{proof}
Cela découle immédiatement du fait que $\det_{\OP}\OO_T\simeq(\det_{\OP}\OO_X)^{\otimes m}$. 
\end{proof}
\end{cor}

\subsubsection{Vérification que $T$ est un $F$-torseur}

\begin{lem}\cite[Lem. 6.14]{jinbi_jin} Soit $T$ un schéma muni d'une action d'un groupe $G$ d'ordre $r$. Soit $f\colon T\to X$ un morphisme étale $G$-équivariant de rang constant $r$. Alors le lieu de $X$ au-dessus duquel $T$ est un torseur est ouvert et fermé dans $X$.
\end{lem}

\begin{cor} Soit $X$ un $\PP^1_S$-schéma fini localement libre. Soient $T$ un schéma muni d'une action d'un groupe $G$ d'ordre $r$, et $f\colon T\to X$ un morphisme étale $G$-équivariant de rang constant $r$.  Notons $X_0,T_0$ les fibres de $X$ et $T$ au-dessus de $0\in\PP^1_S$. Si $T_0$ est un $G$-torseur sur $X_0$ alors $T$ est un $G$-torseur sur $X$.
\begin{proof}
Le lemme \ref{surjcon} montre que $X_0$ rencontre toutes les composantes connexes de $X$. Par le lemme précédent, le lieu de $X$ où $T$ est un $G$-torseur est ouvert et fermé, et il contient $X_0$, c'est donc tout $X$.
\end{proof}
\end{cor}

Un morphisme étale $T\to X$ est donc un $F$-torseur si et seulement si $F\times_{Y_0} T_0\to T_0\times_{Y_0}T_0, (f,x)\mapsto (x,f\cdot x)$ est un isomorphisme de schémas. Cela revient à demander que $\OO_{T_0}\otimes_{\OO_{Y_0}}\OO_{T_0}\to \OO_{T_0}^{|F|}, x\otimes y\mapsto (f\cdot xy)_{f\in F}$ soit un isomorphisme de $\OO_{Y_0}$-modules. Concrètement, la matrice du morphisme $\OO_{T_0}\otimes_{\OO_{Y_0}}\OO_{T_0}\to \OO_{T_0}^{|F|}, x\otimes y\mapsto (f\cdot xy)_{f\in F}$ se calcule comme le produit $\Psi_{F,0}N_0$ où $\Psi_{F,0}$ est la concaténation (verticale) des matrices $\Psi_{f,0},f\in F$ qui représentent l'action des éléments $f\in F$ sur $\OO_{T_0}$. La matrice obtenue est à coefficients dans $k$, et demander à ce qu'elle soit inversible revient à ajouter une variable $W$ et exiger que $\det(\Psi_{F,0}N_0)W-1=0$.

\subsection{Résumé des données et équations}\label{subsec:doneq}

Soit $Y$ une courbe projective lisse munie d'un morphisme $\phi\colon Y\to\PP^1$ tel que $\phi_\star\OO_Y\simeq \OP(a_1)\oplus\dots\oplus\OP(a_r)$. Elle est représentée par des matrices $I\in\Mat_{r\times 1}(k_0)$ et $M\in\Mat_{r\times r^2}(k_0)$ qui définissent la structure de $\OP$-algèbre sur $\phi_\star\OO_Y$.
Soient $G$ un groupe fini d'ordre $n$ donné par une famille génératrice $(g_1,\dots,g_\alpha)$, et $F=\{f_1,\dots,f_\beta\}$ un $\Lambda$-module.
Un $Y$-schéma $T$ de type $b=(b_1,\dots,b_s)$, avec $s=rn$, est défini par les indéterminées suivantes : \begin{itemize}[label=$\bullet$]
\item Une matrice $J\in\Mat_{1\times s}(k_0[x,y])$ et une matrice $N\in\Mat_{s\times s^2}(k_0[x,y])$ qui définissent la structure de $\OP$-algèbre de $\OT$
\item Des matrices $\Phi_{g_1},\dots,\Phi_{g_\alpha}\in\Mat_{s\times s}(k_0[x,y])$ qui définissent l'action de $G$ sur $\OT$
\item Des matrices $\Psi_{f_1},\dots,\Psi_{f_\beta}\in\Mat_{s\times s}(k_0[x,y])$ qui définissent l'action de $F$ sur $\OT$
\item Une matrice $S\in\Mat_{r\times s}(k_0[x,y])$ qui définit le morphisme de $\OP$-algèbres $T\to Y$
\item Une matrice $B\in\Mat_{s\times s}(k_0)$ qui définit un isomorphisme de $\OO_{Y_0}$-modules $\OO_{Y_0}^n\to\OO_{T_0}$ 
\item Des indéterminées supplémentaires $V_1,\dots, V_6$
\end{itemize}

Pour une matrice $A$ à coefficients dans $k_0[x,y]$, notons $A_0$ la matrice $A$ évaluée en $x=0$ et $y=1$. Notons également \[ \Psi_F=\left(\begin{matrix} \Psi_{f_1} \\ \vdots \\\Psi_{f_\alpha}\end{matrix}\right).\] 
Soient $\Delta_{T,1},\Delta_{T,2},\Delta_{X,1},\Delta_{X,2}$ les discriminants respectifs des morphismes $T\to\PP^1$ et $Y\to\PP^1$ sur les deux ouverts standard de $\PP^1$, calculés à partir des matrices $N$ et $M$ grâce aux formules du lemme \ref{lem:discr}. Pour que le schéma $T\to Y$ défini par les données ci-dessus soit un torseur, il faut et il suffit que ces données vérifient les équations suivantes :\begin{itemize}[label=$\bullet$]
\item $N(N\otimes I_s)=N(I_s\otimes N)$ (associativité de la multiplication)
\item $N=N{0 ~ I_s \choose I_s ~ 0}$ (commutativité de la multiplication)
\item $N(J\otimes I_s)=I_s$ (unité)
\item $\Phi_{1_G}=I_s$ et pour tous $i,j\in \{1\dots \alpha\}$ : $\Phi_{g_i}\Phi_{g_j}=\Phi_{g_ig_j}$ (action de $G$)
\item $\Psi_{0_F}=I_s$ et pour tous $i,j\in \{1\dots \beta\}$ : $\Psi_{f_i}\Psi_{f_j}=\Psi_{f_if_j}$ (action de $F$)
\item Pour tous $i\in \{1\dots \alpha\}, j\in\{1\dots\beta\}$ : $\Phi_{g_i}\Phi_{f_j}\Phi_{g_i}^{-1}=\Phi_{g_i\cdot f_j}$ ($G$-équivariance de l'action de $F$)
\item $N(S\otimes S)=SM$ (compatibilité à la multiplication du morphisme $\OO_X\to\OO_T$)
\item $B(I_n\otimes M_0)=N_0(S_0\otimes B)$ ($\OO_{X_0}$-linéarité de $\OO_{X_0}^n\xrightarrow{\sim}\OO_{T_0}$)
\item $V_2\Delta_{T,1}-V_1\Delta_{X,1}=0$, $V_1V_2-1=0$ (non-ramification de $T\to X$)
\item $V_4\Delta_{T,2}-V_3\Delta_{X,2}=0$, $V_3V_4-1=0$ (non-ramification de $T\to X$)
\item $\det (B)V_5-1=0$ (bijectivité du morphisme $\OO_{Y_0}^n\to\OO_{T_0}$)
\item $\det (\Psi_{F,0}M_0)V_6-1=0$ (bijectivité du morphisme $G\times_{X_0}T_0\to T_0\times_{X_0}T_0$) 
\end{itemize}

Le schéma $\mathcal{U}_b$ défini par ces équations est un sous-schéma fermé d'un espace affine $\A^N_{k_0}$ ; chacun de ses $k$-points définit un $F$-torseur $T$ sur $Y$. Un schéma $\mathcal{R}_b$ qui paramètre les morphismes entre deux torseurs se construit de la même façon. En particulier, il y a deux morphismes "source" et "but" $s_b,t_b\colon\mathcal{R}_b\to \mathcal{U}_b$. La classe d'isomorphisme d'un torseur défini par un point $x\in\mathcal{U}_b(k)$ est $s_b(t_b^{-1}x)$. Une étude détaillée de la construction de $\mathcal{U}_b$ et $\mathcal{R}_b$ donne le résultat suivant.

\begin{prop}\cite[Prop. 7.2]{jinbi_jin}\label{fibirred}
Les morphismes $s_b,t_b\colon\mathcal{R}_b\rightrightarrows \mathcal{U}_b$ sont lisses à fibres géométriques irréductibles. De plus, les schémas $t_b(s_b^{-1}x)$, pour $x\in \mathcal{U}_b(k)$, ont tous la même dimension.
\end{prop}
Définissons enfin $\mathcal{U}\coloneqq \bigsqcup_b\mathcal{U}_b$, $\mathcal{R}\coloneqq \bigsqcup_b\mathcal{R}_b$ ; les morphismes $s_b$, $t_b$ définissent encore des morphismes $s,t\colon\mathcal{U}\rightrightarrows\mathcal{R}$.

\subsection{Le schéma en groupoïdes qui paramètre les torseurs}

\subsubsection{Catégories fibrées et schémas en groupoïdes}

\begin{df}  Soit $p\colon \mathcal{D}\to\mathcal{C}$ un foncteur. \begin{enumerate}
\item Soit $f\colon y\to x$ un morphisme dans $\mathcal{D}$. Le morphisme $f$ est dit fortement cartésien si pour tout objet $z$ de $\mathcal{D}$, l'application $\Hom_\mathcal{D}(z,y)\to \Hom_\mathcal{D}(z,x)\times_{\Hom_\mathcal{C}(p(z),p(x))}\Hom_\mathcal{C}(p(z),p(y)), \phi\mapsto(f\circ\phi,p(\phi))$ est bijective. 
\item La catégorie $\mathcal{D}$ est dite fibrée au-dessus de $\mathcal{C}$ si pour tout objet $x$ de $\mathcal{D}$, et tout morphisme $g\colon c\to p(x)$ dans $\mathcal{C}$, il existe un morphisme fortement cartésien $f\colon y\to x$ dans $\mathcal{D}$ tel que $p(f)=g$.
\item Un morphisme entre deux catégories fibrées $p\colon \mathcal{D}\to\mathcal{C}$, $p'\colon\mathcal{D}'\to\mathcal{C}$ est un foncteur $F\colon\mathcal{D}\to\mathcal{D}'$ tel que $p'\circ F=p$, et qui préserve la forte cartésianité des morphismes.
\end{enumerate}
\end{df}

\begin{df} Soit $p\colon\mathcal{D}\to\mathcal{C}$ une catégorie fibrée.
\'{E}tant donné un objet $c$ de $\mathcal{C}$, notons $\mathcal{D}(c)$ la catégorie des objets $d$ de $\mathcal{D}$ tels que $p(d)=c$, avec pour morphismes les $f\colon d'\to d$ tels que $p(f)=\id_c$. La catégorie $\mathcal{D}$ est dite fibrée en groupoïdes au-dessus de $\mathcal{C}$ si pour tout objet $c$ de $\mathcal{C}$, la catégorie $\mathcal{D}(c)$ est un groupoïde.
\end{df}

\begin{rk} Soit $p\colon\mathcal{D}\to\mathcal{C}$ une catégorie fibrée.\begin{enumerate}
\item \'{E}tant donné un objet $x$ de $\mathcal{D}$ et un morphisme $g\colon c\to p(x)$ dans $\mathcal{C}$, un relèvement fortement cartésien $f\colon y\to x$ de $g$ est unique à isomorphisme près. La notation $y=g^\star x$ est donc sans ambiguïté. Fixons désormais un tel élément $g^\star x'$ pour chaque morphisme $x'\to x$ dans $\mathcal{D}(p(x))$. Alors, pour un tel morphisme $\alpha \colon x'\to x$, il existe un unique morphisme $g^\star \alpha \colon g^\star x'\to g^\star x$ dans $\mathcal{D}(c)$ tel que le diagramme suivant soit commutatif.

\[
\begin{tikzcd}
g^\star x' \arrow[d]\arrow[r,"g^\star \alpha"] & g^\star x\arrow[d] \\
x'\arrow[r,"\alpha"] & x
\end{tikzcd}
\]
Ceci définit un foncteur $g^\star \colon \mathcal{D}(p(x))\to\mathcal{D}(c)$ \cite[02XJ, Def. 4.33.6]{stacks}. De plus, si $f$ et $g$ sont des morphismes composables dans $\mathcal{C}$ alors il y a un unique isomorphisme $(f\circ g)^\star\xrightarrow{\sim} g^\star\circ f^\star$ \cite[02XJ, Lem. 4.33.7]{stacks}. De plus, un morphisme de catégories fibrées préserve (à isomorphisme près) les tirés en arrière. 
\item En particulier, étant donné un schéma $S$, il y a dans toute catégorie fibrée au-dessus de $\Sch/S$ des foncteurs de changement de base : si $f\colon Y\to X$ est un morphisme de $S$-schémas, il y a un foncteur $f^\star=\times_X Y \colon\mathcal{D}(X)\to\mathcal{D}(Y)$.
\item 
Lorsque $S=\Spec K$ est le spectre d'un corps et $L$ est une extension galoisienne de $K$, l'ensemble $\pi_0(\mathcal{D}(L))$ des classes d'isomorphisme dans $\mathcal{D}(L)$ est muni d'une action à droite de $\Gal(L|K)$, c'est-à-dire un morphisme de groupes $\Gal(L|K)^{\rm op}\to\mathfrak{S}(\pi_0(\mathcal{D}(L))$ défini par $\sigma\mapsto\sigma^\star$.
\end{enumerate}
\end{rk}

\begin{df}\label{cattors}
Soient $S$ un schéma, $X$ un $S$-schéma et $G$ un $X$-schéma en groupes. La catégorie fibrée sur $\Sch/S$ des $G$-torseurs sur $X$ est la catégorie $\mathcal{T}$ dont les objets sont les $(S',T)$ où $S'\in\Sch/S$ et $T$ est un $G_{S'}$-torseur sur $X_{S'}$. Les morphismes entre $(S'',T'')$ et $(S',T')$ sont les morphismes de torseurs si $S''=S'$, et il n'y en a pas si $S''\neq S'$. \'{E}tant donné $S''\to S'$, le foncteur de changement de base $\mathcal{T}(S')\to\mathcal{T}(S'')$ est le foncteur $\times_{S'}S''$. Le foncteur vers $\Sch/S$ est le foncteur d'oubli, qui à un torseur associe le schéma sous-jacent, et à un morphisme $(S',T)\to (S',T')$ associe $\id_{S'}$.
\end{df}

\begin{df} Un schéma en groupoïdes sur un schéma $S$ est la donnée d'un couple $(R,U)$ de $S$-schémas muni de morphismes "source" et "but" $s,t\colon R\rightrightarrows U$ et "composition" $\circ\colon R\times_{U,s,t}R\to R$ tels que pour tout $S$-schéma $T$, la catégorie d'objets $U(T)$, de morphismes $R(T)$ dont la source et le but sont donnés par $s_T$ et $t_T$, avec la composition définie par $\circ_T$, soit un groupoïde.
\end{df}

\begin{df}\label{catschgrp} Soient $s,t\colon R\rightrightarrows U$ un schéma en groupoïdes sur un $k_0$-schéma $S$. La catégorie fibrée au-dessus de $\Sch/S$ associée est définie comme suit.
Ses objets sont les couples $(S',x)$ avec $S'\to S$ et $x\in U(S')$. Ses morphismes sont définis par $\Hom((S',x'),(S',x))=\{ f\in R(S')\mid s(f)=x', t(f)=x\}$, et le foncteur vers $\Sch/S$ est le foncteur d'oubli $(S',x)\mapsto S'$ qui envoie tout morphisme sur $\id_{S'}$. Ceci définit une catégorie fibrée en groupoïdes sur $\Sch/S$, dont la fibre au-dessus de $T\in\Sch/k_0$ est le groupoïde $U(T)$. \'{E}tant donné un morphisme $\alpha \colon S''\to S'$ dans $\Sch/S$, le foncteur de changement de base $\alpha^\star \colon U(S')\to U(S'')$ est la composition à gauche par $\alpha$.
\end{df}

\subsubsection{Le résultat principal}

La proposition suivante est une variante de \cite[Prop. 7.5]{jinbi_jin} avec des hypothèses moins contraignantes ; la preuve est essentiellement inchangée.

\begin{prop}\label{thprinc}
Soit $U$ un $k_0$-schéma de type fini.
Soit $\mathcal{T}\to \Sch/k_0$ une catégorie fibrée en groupoïdes.
Soit, pour chaque $k_0$-schéma $S$, une application $F_S\colon U(S)\to \mathcal{T}(S)$.  Supposons que ces données vérifient les propriétés suivantes : \begin{enumerate}
\item pour tout morphisme $\phi\colon S''\to S'$ de $k_0$-schémas et tout $x\in U(S')$, $F_{S''}(x\circ \phi)=\phi^\star(F_{S'}(x))$ ;
\item l'application $F_k\colon U(k)\to\mathcal{T}(k)$ est essentiellement surjective ;
\item pour tout $k_0$-schéma $S$ et tous objets $x,y\in\mathcal{T}(S)$, le foncteur $\Sch/S\to\Set$, $(S'\xrightarrow{\alpha}S)\mapsto \Isom_{\mathcal{T}(S')}(\alpha^\star x,\alpha^\star y)$ est représenté par un schéma $Y\to S$ ouvert et fermé.
\end{enumerate}
Alors l'application surjective et $\Gal(k|k_0)$-équivariante $U(k)\to\pi_0(\mathcal{T}(k))$ induite par $F_k$ se factorise en une application encore surjective et $\Gal(k|k_0)$-équivariante $\pi_0^{sch}(U_k)\to \pi_0(\mathcal{T}(k))$.
\end{prop}

\begin{proof}
La Galois-équivariance de l'application vient de la compatibilité de $F$ au changement de base. Soient $x\in\mathcal{U}(k)$, et $C$ la composante connexe de $U_k$ contenant l'image de $x$. Notons $f\colon C\to\Spec k$ le morphisme structural, et $\bar x\colon \Spec k\to C$ la factorisation de $x$ par $C$. La situation est résumée par le diagramme suivant (qui devient commutatif en retirant la flèche $f$) :
\[
\begin{tikzcd}
\Spec k \arrow[r,"x"]\arrow[d,"\bar x",shift left=1] &U\\
\arrow[u,"f",shift left=1] C \arrow[r,"j",swap]& U_{k}\arrow[u,"i",swap]
\end{tikzcd}
\]

Soient $\phi_x=x\circ f=f^\star x^\star\id_U$ et $\phi_C=i\circ j=j^\star i^\star\id_U\in U(C)$. Notons $\psi_x,\psi_C$ les images par $F_C$ de $\phi_x,\phi_C$ dans $\mathcal{T}(C)$. 
Le foncteur $\Sch/C\to \Set$, $(S\xrightarrow{\alpha}C)\mapsto{\rm Isom}_{\mathcal{T}(C)}(\alpha^\star\psi_C,\alpha^\star\psi_x)$ est par hypothèse représentable par un schéma $Y\to C$ ouvert et fermé. Par construction, $Y(\bar x)$ est non vide. En effet, \[\bar{x}^\star\psi_x=\bar{x}^\star F_C(\phi_x)=F_k(\phi_x\circ\bar x)=F_k(x\circ f\circ\bar x)=F_k(x) \]
par compatibilité au changement de base, et car $\bar x$ est une section de $f$. De même, $\bar{x}^\star\psi_C=F_k(x)$ car $x=i\circ j\circ\bar{x}$. Ainsi, $\id_{F_k(x)}\in Y(\bar x)$.
Le morphisme $Y\to C$ est donc un morphisme ouvert et fermé d'un schéma non vide vers un schéma connexe : il est surjectif. Par conséquent, pour tout autre point $\bar x'\in C(k)$, $Y(\bar x')\neq 0$. 
Il y a donc pour tous $\bar x,\bar x'\in C(k)$ un isomorphisme $\bar{x}'^\star\psi_x\simeq \bar{x}'^\star\psi_C$ dans $\mathcal{T}(k)$. Or $\bar x'^\star\psi_x=\bar x'^\star F_C(x\circ f)=F_k(x\circ f\circ\bar{x}')=F_k(x)$ ; de même, $\bar{x}'^\star\psi_{x'}=F_k(x')$. 

Par conséquent, $F_k(x)$ et $F_k(x')$ sont isomorphes dans $\mathcal{T}(k)$. La surjection $\mathcal{U}(k)\to \mathcal{T}(k)$ se factorise donc en une application $\pi_0^{sch}(\mathcal{U}_k)\to \pi_0(\mathcal{T}(k))$, qui est encore surjective. 
\end{proof}

Rappelons \cite[0478]{stacks} que tout schéma de type fini $X$ sur $k$ est un schéma de Jacobson  : les points fermés sont denses dans tout fermé de $X$. En particulier, les composantes irréductibles (resp. connexes) de $X$ sont en bijection canonique avec celles de $X(k)$ muni de sa topologie de Zariski naïve.

\begin{prop}\label{prop:compconnirred} Mêmes notations et hypothèses que la proposition précédente. Si, de plus, les préimages par $F_k$ des classes d'isomorphisme de $\mathcal{T}(k)$ sont connexes dans $U(k)$  alors l'application $\pi_0^{sch}(U_k)\to \pi_0(\mathcal{T}(k))$ induite par $F_k$ est bijective. Enfin, si ces préimages sont irréductibles dans $U(k)$ alors les composantes connexes de $U_k$ sont irréductibles.

\end{prop}
\begin{proof} Considérons $x,x'\in U(k)$ tels que $F_k(x)$ soit isomorphe à $F_k(x')$ dans $\mathcal{T}(k)$. Par hypothèse, ils appartiennent à une même partie connexe de $|U_k|$. L'image dans $U_k$ de cette partie connexe est encore connexe puisque l'injection $|U_k|\to U_k$ est continue, et les images de $x,x'$ dans $U_k$ appartiennent donc à une même composante connexe de $U_k$. Pour le second point, servons-nous de la bijection croissante entre les composantes connexes (resp. irréductibles) de $U_k$ et celles de son ensemble de points fermés $|U_k|=U(k)$. Considérons deux points $x,x'\in U(k)$ dans une même composante connexe de $U(k)$. Alors d'après la proposition \ref{thprinc}, leurs images par $F_k$ sont isomorphes. Par conséquent, ils appartiennent à une même composante irréductible de $U(k)$, ce qui conclut.
\end{proof}

\subsubsection{Application}

Le théorème est appliqué à la catégorie $\mathcal{T}$ des $G$-torseurs sur une courbe $Y$ (voir la définition \ref{cattors}), et au schéma $\mathcal{U}$ construit à la fin de la section \ref{subsec:doneq}.\\

\'{E}tant donné un $k$-schéma $S$, un point de $\mathcal{U}(S)$ définit encore un torseur sur $X_S$. En effet, un morphisme de fibrés vectoriels $\bigoplus_{i=1}^s\OO_{\PP^1_S}(a_i)\to\bigoplus_{j=1}^t\OO_{\PP^1_S}(b_j)$ est défini par une matrice $t\times s$ à coefficients dans $\HH^0(S,\OO_S)[x,y]$, dont l'élément en position $(i,j)$ est un polynôme homogène de degré $b_j-a_i$. La donnée de matrices (multiplication, action du groupe $F$) à coefficients dans $\HH^0(S,OO_S)$ définit donc encore un schéma fini localement libre sur $\PP^1_S$. Les conditions suffisantes pour être un torseur ont été démontrées dans un cadre relatif (voir la section \ref{subsec:condtors}), et s'appliquent donc encore ici. Remarquons également que toute cette construction commute au changement de base : étant donnés un morphisme $f\colon S'\to S$ et un point $x\in U(S)$ définissant un torseur $T\to X$, le torseur défini par $f^\star x=x\circ f\in U(S')$ est le torseur $T\times_S S'$. Enfin, chaque $F$-torseur sur $Y$ est isomorphe à un torseur défini par l'un des $k$-points de $\mathcal{U}$. La proposition \ref{thprinc} s'applique donc à $\mathcal{U}$ : pour trouver un représentant de chaque classe d'isomorphisme de $F$-torseurs sur $Y$, il suffit de trouver un $k$-point dans chaque composante connexe de $\mathcal{U}_k$.

\begin{cor} Le schéma $\mathcal{U}_b$ est équidimensionnel.
\begin{proof}
Les composantes irréductibles de $\mathcal{U}_b$ sont d'après le théorème principal les classes d'isomorphisme $t(s^{-1}x)$ pour $x\in \mathcal{U}_b(k)$. D'après le lemme \ref{fibirred}, elles sont toutes de même dimension.
\end{proof}
\end{cor}

Nous avons expliqué la construction du schéma $\mathcal{U}$ paramétrant les $\F$-torseurs sur $X$. D'après la discussion précédente, il suffit désormais pour calculer $\HH^1(X,\F)$ de déterminer un point dans chaque composante connexe de $X$, c'est-à-dire dans chaque composante irréductible de $X$ d'après la proposition \ref{prop:compconnirred}.

\subsection{Calcul de représentants des classes de torseurs}

Voici comment déterminer au moins un point de chaque composante irréductible du schéma affine $\mathcal{U}$ calculé précédemment. Comme ce dernier est équidimensionnel, la procédure décrite dans la section \ref{subsubsec:ptsequidim}, qui consiste à calculer la fibre en 0 d'une normalisation de Noether $\nu \colon\mathcal{U}\to \A^D$, convient. Elle fournit une liste de points, dont plusieurs appartiennent peut-être à une même composante. Deux points $x,y$ de cette liste appartiennent à la même composante si et seulement si $s^{-1}x\times_{\mathcal{R}} t^{-1}y$ est non vide : il suffit d'en déterminer des équations, puis de tester si elles engendrent l'idéal unité. 
Une fois déterminé l'ensemble $\HH^1(X,\F)$, il reste à calculer sa loi de groupe ; celle-ci est donnée par le produit contracté, et s'exprime en termes de corps de fonctions par des méthodes d'algèbre linéaire.

Ceci conclut la description de l'algorithme de Jin pour les faisceaux lisses sur les courbes projectives. Nous expliquons dans la section suivante comment une petite modification de cet algorithme permet de calculer également la cohomologie d'un faisceau lisse sur une courbe affine.

\subsection{Courbes affines}

L'algorithme de Jin est exposé dans \cite{jinbi_jin} dans un cadre plus général, qui englobe le calcul de $\HH^1(U,\F)$ et $\HH^1_c(U,\F)$ où $U$ est une courbe affine lisse sur $k$. Considérons la situation suivante : $U$ est une courbe lisse, $X$ sa complétion projective lisse munie d'un morphisme $X\to \PP^1_k$ tel que $U=X\times_{\PP^1}\A^1$. La construction d'un tel morphisme est décrite dans \cite[Prop. 9.8]{jinbi_jin}. Le faisceau lisse $\F$ sur $U$ de fibre $F$ est trivialisé par un revêtement galoisien $V\to U$ de groupe $G$ ; le schéma $Y$ est la normalisation de $X$ dans $V$, et $W$ est la fibre de $Y$ au-dessus de $0\in \PP^1$.

\[
\begin{tikzcd}
V\arrow[d,"g"]\arrow[r,"j'"] & Y \arrow[d,"f"] & \arrow[l,"i'",swap] W\arrow[d,"h"]\\
U\arrow[r,"j"]\arrow[d] &X\arrow[d]& \arrow[l,"i",swap] Z \arrow[d]\\
\A^1_{k} \arrow[r] & \PP^1_{k} &\arrow[l] 0
\end{tikzcd}
\]

\begin{lem} Soient $C$ la catégorie des $\F$-torseurs sur $U$, et $D$ la catégorie des schémas $T\to Y$ lisses sur $k$, munis d'une action $G$-équivariante de $F$ tels que $T|_V\to V$ soit un $\F$-torseur. La composition de $g^\star$ avec la complétion projective lisse définit une équivalence de catégories $C\to D$, de quasi-inverse $g_\star^{G}\circ(-\times_{\PP^1}\A^1$).
\begin{proof} D'après le corollaire \ref{cor:torsrev}, $g^\star$ induit une équivalence de la catégorie $C$ avec la catégorie des $\F|_V$-torseurs $G$-équivariants sur $V$. Ensuite, un schéma $T$ est un $\F|_V$-torseur $G$-équivariant sur $V$ si et seulement si sa complétion projective lisse est un schéma lisse sur $k$ muni d'un morphisme $G$-équivariant vers $Y$ qui en fait un $\F|_V$-torseur sur $V$.
\end{proof}
\end{lem}

Afin de calculer $\HH^1(U,\F)$, il suffit de modifier deux choses dans la construction du schéma paramétrant les $F$-torseurs $G$-équivariants sur $Y$. D'une part, l'étalitude de $T\to Y$ doit être remplacée par celle de $T\times_{\PP^1}\A^1\to V$. Cela revient \cite[Prop. 6.12]{jinbi_jin} à demander que le discriminant de $\Delta_{T\to\PP^1}$ diffère de $\Delta_{Y\to\PP^1}$ non pas d'un élément de $k^\times$, mais d'un élément de $k^\times$ fois une puissance de $y$ qui se calcule à partir des données de $\OO_T$ et $\OO_X$. D'autre part, il faut s'assurer que le schéma $T$ est lisse au-dessus de $0\in\PP^1$ : une condition nécessaire et suffisante se traduisant par des équations est donnée dans \cite[Prop. 6.16]{jinbi_jin}.

\subsection{Complexité}

Soit $X$ une courbe lisse sur $k_0$ de genre géométrique $g$. Supposons donné un morphisme $\phi\colon X\to\PP^1$, ainsi que la $\OP$-algèbre $\phi_\star\OX\simeq\OP(a_1)\oplus\dots\oplus \OP(a_r)$. Soit $\F$ un faisceau lisse de groupes abéliens sur $X$, de fibre $F$, décrit explicitement par un revêtement galoisien trivialisant $Y\to X$ de groupe $G$ et le groupe abélien $F=\HH^0(Y,\F)$ d'ordre $m$.
Les torseurs $T\to X$ considérés par l'algorithme vérifient $(T\to\PP^1)_\star\OT\simeq\OP(b_1)\oplus\dots\oplus\OP(b_s)$ avec $s=mr$.

\begin{prop}\label{lem:compjin}\cite[Prop. 6.21,6.24]{jinbi_jin} Le nombre de types $b=(b_1,\dots,b_s)$ vérifiant $b_j\leqslant 0$, $s=mr$ et $\sum_j b_j=m\sum_j a_j$ est $O((rgm)^{rm})$, où $g$ est le genre de $X$. Pour chaque type $b$, le schéma $\mathcal{U}_b$ calculé par l'algorithme est donné par $O(s^5m^4)$ polynômes de degré au plus $sn$ en $O(s^5n^4)$ variables. 
\end{prop}

\begin{prop}\cite[Cor. 8.12]{jinbi_jin} 
L'algorithme de Jin calcule un $k$-schéma en groupes zéro-dimensionnel représentant $\HH^1(X,\F)$ en $\exp(O(r^{15}m^{12}g^3\log(rm)^3))$ opérations dans $k_0$.
\end{prop}

La complexité exponentielle en $m\log(m)$ de l'algorithme global vient donc directement de la construction du schéma $\mathcal{U}$, et non d'un problème de nature algorithmique lié au calcul des composantes connexes de ce schéma : cette complexité ne peut pas être améliorée par des astuces algorithmiques.

\begin{rk}
Afin d'obtenir cette complexité, Jin fait usage d'un algorithme de Khuri-Makdisi \cite[§7]{khuri_makdisi}, qui suppose qu'il existe un algorithme de factorisation des polynômes de $k_0[t]$ en temps polynomial. Ceci est vrai lorsque $k_0$ est fini, et encore vrai si $k_0$ est un corps de nombres ; cependant, dans ce cas, il faut prendre en compte dans le calcul de la complexité la hauteur des coefficients des polynômes concernés (et non pas seulement leur degré), ce qui n'est pas étudié dans \cite{jinbi_jin}. Ce travail paraît fastidieux, mais pas difficile, à réaliser.
\end{rk}

\section{Aspects pratiques}

Voici quelques remarques sur les aspects pratiques des algorithmes mentionnés ci-dessus, dont aucun n'a été implémenté jusqu'ici. Le seul algorithme qui paraît implémentable avec une quantité d'efforts raisonnable est celui de Couveignes. De plus, au vu de sa complexité, il serait certainement utilisable dans la pratique pour de petits paramètres. Les algorithmes de Huang et Ierardi et de Jin sont très semblables : ils construisent tous les deux un grand schéma $S$ paramétrant les objets dont on cherche des classes d'isomorphisme, puis cherchent des points dans les composantes irréductibles de ce schéma. Les deux paraissent en l'état cauchemardesques à implémenter. Cependant, une fois implémenté, l'algorithme de Huang et Ierardi aurait un avantage de rapidité d'exécution, car la dimension de l'espace dans lequel est plongé le schéma $S$, linéaire en le genre de la courbe, augmente bien moins vite en fonction des paramètres que dans le cas de l'algorithme de Jin, où elle dépend non seulement de la représentation de la courbe mais aussi du cardinal de la fibre du faisceau (voir lemme \ref{lem:compjin}). La complexité de la recherche de points étant exponentielle en cette dimension (voir annexe \ref{subsubsec:ptscompconstr}), cet avantage serait conséquent. Enfin, l'algorithme de Madore et Orgogozo nécessite en premier lieu l'utilisation d'un algorithme de calcul du $\HH^1$ des faisceaux constants sur les courbes, dont aucun n'a encore été implémenté ; de plus, il ne serait pas utilisable dans la pratique en raison de sa complexité conséquente.

Dans le chapitre suivant, nous décrivons des algorithmes de calcul de la cohomologie des faisceaux constructibles sur les courbes lisses ou nodales. Nous nous baserons sur les algorithmes de Couveignes et de Huang et Ierardi pour nos résultats ; en l'absence d'une implémentation de ces algorithmes, nous nous contenterons dans nos exemples de courbes où le calcul de points de torsion de la jacobienne peut se faire par des moyens plus élémentaires.

\cleartooddpage

\chapter{Calcul effectif de la cohomologie}\label{chap:5}
Dans l'ensemble de ce chapitre, $k_0$ désigne un corps parfait, $k$ une clôture algébrique de $k_0$, et $n$ un entier inversible dans $k$. Nous noterons $\Lambda$ l'anneau $\ZZ/n\ZZ$, et $\mathfrak{G}_0$ le groupe $\Gal(k|k_0)$.\\

L'objectif de ce chapitre est de donner une description explicite, étant donné un complexe $\K$ de faisceaux constructibles de $\Lambda$-modules sur une $k$-courbe lisse $X$, d'un complexe de $\Lambda$-modules représentant $\RG(X,\K)\in \DD^b_c(\Lambda)$. Lorsque $X$ provient par changement de base de $k_0$, nous décrivons l'action de $\mathfrak{G}_0$ sur un tel complexe. Nous estimerons également la complexité du calcul de $\RG(X,\K)$. Ce calcul repose sur celui de la cohomologie des faisceaux lisses sur une courbe affine lisse, qui se déduit lui-même du cas des faisceaux constants. Nous prouverons en particulier le résultat suivant.

\begin{theorem}\label{th:compRG} Soit $X$ une courbe intègre lisse sur $k$. Soit $\F$ un faisceau constructible de $\Lambda$-modules sur $X$. Soient $U$ un ouvert de $X$ sur lequel $\F$ est lisse de fibre générique géométrique $M$, et $Z$ le fermé réduit complémentaire. Soient $V\to U$ un revêtement galoisien qui trivialise $\F|_U$, et $V_2\to V$ le revêtement de $V$ de groupe $H^1(V,\Lambda)^\vee$. Notons $G=\Aut(V_2|U)$. Pour chaque point $z\in Z$, notons $I_z$ le groupe d'inertie d'un point de la compactification lisse de $V_2$ au-dessus de $z$, et $P_z$ le groupe d'inertie sauvage correspondant. Notons $\phi_z\colon \F_z\to M^{I_z}\subseteq M^{P_z}\xrightarrow{\sim}M_{P_z}$ le morphisme de recollement en $z$ composé avec l'isomorphisme canonique $M^{P_z}\xrightarrow{\sim}M_{P_z}$, qui à un élément $P_z$-invariant de $M$ associe sa classe dans le module des coinvariants $M_{P_z}$. Alors $\RG(X,\F)[1]$ est le cône du morphisme de complexes suivant, où $C^{12}(G,M)$ désigne le groupe des morphismes croisés $G\to M$.
\[
\begin{adjustbox}{width=\textwidth}{
\begin{tikzcd}
M\oplus \bigoplus_{z\in Z}\F_z \arrow[d,"\bigoplus_{z\in Z}(\id-\phi_z)"] \arrow[r,"{(\partial_G,0)}"] & C^{12}(G,M) \arrow[d,"\bigoplus_z \res_{I_z}^G"]\arrow[r] & \bigoplus_{z\in Z}\HH^1(I_z/P_z,M_{P_z}) \arrow[d,"\id"]\arrow[r] & 0 \\
\bigoplus_{z\in Z}M_{P_z} \arrow[r,"\bigoplus_z\partial_{I_z}"] \arrow[r] & \bigoplus_{z\in Z}C^{12}(I_z/P_z,M_{P_z}) \arrow[r] & \bigoplus_{z\in Z}\HH^1(I_z/P_z,M_{P_z}) \arrow[r] & 0
\end{tikzcd}
}
\end{adjustbox}
\]
\end{theorem}
Nous étudierons en détail la complexité des algorithmes présentés, et montrerons qu'elle est doublement exponentielle en le genre de la courbe dans le cas général. Les différentes méthodes de calcul de la cohomologie et leurs interdépendances sont résumées dans le diagramme ci-après. Ici, $X$ est une courbe intègre lisse sur $k$, $\F$ est un faisceau constructible sur $X$ et $\mathscr{K}$ est un complexe de faisceaux constructibles sur $X$. La colonne de gauche indique la nature de $\F$ (resp. des termes de $\mathscr{K}$).

\begin{tikzpicture}[node distance=2cm, every text node part/.style={align=center}]
\node (proj) [titre,xshift=-3cm] {$X$ projective lisse};
\node (aff) [titre, right of=proj,xshift=4.5cm] {$X$ affine lisse};
\node (const) [titrevert, left of = proj,xshift=-2.5cm, yshift=-4.5cm] {faisceau $\F$ \\ constant};
\node (lis) [titrevert, below of = const,yshift=-4cm] {faisceau $\F$ \\ lisse};
\node (constr) [titrevert, below of = lis, yshift=-4cm] {faisceau $\F$ \\ constructible};

\node (H1cstcouv) [bon, below of=proj, yshift=-1cm] {$\HH^1(X,\F)$ \\ \ref{sec:couv}};
\node (couv) [auteur,above of=H1cstcouv,yshift=-1.2cm]{Couveignes};
\node (H1caff) [bon, below of=aff,yshift=-1cm] {$\HH^1_c(X,\F)$ \\ \ref{sec:cohaff}};
\draw [-{Stealth[scale=1.8]}] (H1cstcouv)--(H1caff);
\node (H1cst) [bon, below of=proj, yshift=-3.5cm] {$\HH^1(X,\F)$ \\ \ref{sec:huang}};
\node (huang) [auteur,above of=H1cst,yshift=-1.2cm]{Huang-Ierardi};
\node (H1cstaff) [bon, below of=aff,yshift=-3.5cm] {$\HH^1(X,\F)$ \\ \ref{sec:cohaff}};
\draw [-{Stealth[scale=1.8]}] (H1cst)--(H1cstaff);
\node (RGlisaff) [bon, below of=H1cstaff,yshift=-1.5cm] {$\RG(X,\F)$ \\ \ref{subsec:cohlisG}};
\draw [-{Stealth[scale=1.8]}] (H1cstaff)--(RGlisaff);
\node (RGlis) [bon, below of=H1cst, yshift=-1.5cm] {$\RG(X,\F)$ \\ rem. \ref{rk:mayerRG}};
\draw [-{Stealth[scale=1.8]}] (RGlisaff)--(RGlis);
\node (H1lisjin) [final, below of=RGlis,xshift = 3.25cm,yshift=-0.5cm] {$\HH^1(X,\F)$ \\ \ref{sec:jin}};
\node (H1constrjin) [final, below of=H1lisjin,yshift=-1cm] {$\HH^1(X,\F)$ \\ \ref{sec:jinconstr}};
\draw [-{Stealth[scale=1.8]}] (H1lisjin)--(H1constrjin);
\node (jin) [auteurfinal,above of=H1lisjin,yshift=-1.2cm]{Jin};
\node (RGconstr) [final, below of=H1constrjin] {$\RG(X,\F)$ \\ \ref{subsec:algconstr}};
\node (RGconstrcomp) [final, below of=RGconstr] {$\RG(X,\K)$ \\ \ref{subsec:algconstrcomp}};
\draw [-{Stealth[scale=1.8]}] (RGlisaff)--(3.5,-16.5)--(RGconstr);
\draw [-{Stealth[scale=1.8]}] (RGconstr)--(RGconstrcomp);
\end{tikzpicture}

\section{Scindage explicite des suites exactes courtes}\label{sec:scindage}

Le principe ci-dessous sera souvent utilisé dans la suite. Soit $G$ un groupe fini. Considérons une suite exacte courte de $\Lambda[G]$-modules de type fini
\[ 0\to A\xrightarrow{\phi} B\xrightarrow{\psi} C\to 0 \]
où $C$ est libre. Le but est de calculer $B$ en connaissant $A$ et $C$. Supposons que l'on dispose d'une description explicite de $A$ et $C$ comme $\Lambda[G]$-modules, c'est-à-dire de présentations de $A$ et $C$ comme $\Lambda$-modules et pour tout $\sigma\in G$, des endomorphismes de $A$ et $C$ définis par $\sigma$.
Soient $(a_1,\dots,a_r)$ et $(c_1,\dots,c_m)$ des familles génératrices respectives de $A$ et $C$.
Supposons également qu'il existe des éléments $\lambda_{ij,g}\in \Lambda, a_{i,g}\in A$ (que l'on sait calculer) tels qu'il existe une section $\Lambda$-linéaire $s$ de $\psi$ vérifiant pour tous $g\in G, i\in \{1,\dots, m\}$ : \[ g\cdot s(c_i)=\sum_{j=1}^m \lambda_{ij,g}s(c_j)+a_{i,g}.\tag{$\star$}\]
Décrivons comment calculer le $\Lambda[G]$-module $B$ explicitement. Il est représenté en tant que $\Lambda$-module par $A\oplus C$ ; l'action de $g\in G$ sur $a+c\in B$ se calcule à l'aide de $(\star)$.

Voyons comment cette description se comporte vis-à-vis des morphismes de suites exactes. \'{E}tant donné un morphisme de suites exactes de $\Lambda[G]$-modules :
\[
\begin{tikzcd}
0 \arrow[r] & A \arrow[r,"\phi"] \arrow[d,"\alpha"] & B \arrow[r,"\psi"] \arrow[d,"f"] & C\arrow[d,"\gamma"]\arrow[r] & 0 \\
0 \arrow[r] & A' \arrow[r,"\phi'"] & B' \arrow[r,"\psi'"] & C' \arrow[r] & 0
\end{tikzcd}
\]
où l'on dispose explicitement des données ci-dessus pour les deux suites exactes, ainsi que des morphismes $\alpha$ et $\gamma$, il est possible de calculer le morphisme $f$ dès que les sections $\Lambda$-linéaires $s$ de $\psi$ et $s'$ de $\psi'$ sont compatibles au sens où $s'\gamma=fs$. Soit $b\in B$. \'{E}crivons $b=\phi(a)+s(c)$, où $s$ est la section donnée de $\psi$. Alors la commutativité du carré de gauche donne $f(b)=\phi' \alpha(a)+ f(b-\phi(a))=\phi'\alpha(a)+fs(c)=\phi'\alpha(a)+s'\gamma(c)$.

\section{Faisceaux constants sur les courbes affines}\label{sec:cohaff}

Soit $X_0$ une courbe projective lisse sur $k_0$.  Soient $Z_0$ un fermé réduit de $X_0$, et $U_0$ l'ouvert complémentaire. Notons $X,Z,U$ leurs changements de base à $k$. Cette section porte sur le calcul de $\HH^1(U,\Lambda)$, où $\Lambda$ est le faisceau constant de valeur $\ZZ/n\ZZ$. Le choix d'une racine primitive $n$-ième de l'unité dans $k$ détermine un isomorphisme $\Lambda\to\mu_n$. Il y a un isomorphisme canonique
\[ \HH^1(U,\Lambda)\otimes\mu_n\xrightarrow{\sim}\HH^1(U,\mu_n)\]
qui permet de déduire $\HH^1(U,\Lambda)$ avec son action de $\mathfrak{G}_0=\Gk$ de $\HH^1(U,\mu_n)$ ; cet isomorphisme se construit par exemple aisément sur les groupes de cohomologie de \v{C}ech.
 Nous expliquons d'abord comment représenter et manipuler les éléments de $\HH^1(U,\mu_n)$, puis comment calculer les groupes $\HH^1(U,\mu_n)$ et $\HH^1_c(U,\mu_n)$, puis comment, étant donné un morphisme de courbes affines $V\to U$, décrire le morphisme canonique $\HH^1(U,\mu_n)\to \HH^1(V,\mu_n)$.

\subsection{Calcul dans $\HH^1(U,\mu_n)$}
Rappelons (cf lemme \ref{lem:H1aff}) que $\HH^1(U,\mu_n)$ est en bijection avec le groupe des classes d'équivalence de triplets $([D],D',f)$ où $D\in\Div^0X$, $D'\in \Div^0X$ est à support dans $Z$ et $nD=\div(f)+D'$.

\begin{rk}\label{rk:egalH1} Décrivons précisément comment calculer dans ce groupe quotient isomorphe à $\HH^1(U,\mu_n)$. Deux triplets $([D_1],D'_1,f_1)$ et $([D_2],D'_2,f_2)$ sont égaux dans le quotient si et seulement si $D_1'-D_2'\in n\Div^0_Z(X)$ et $(D_1-D_2)-\frac{1}{n}(D_1'-D_2')$ est principal. Ceci se teste algorithmiquement, en vérifiant si les coefficients de $D'_1-D'_2$ sont multiples de $n$, puis en calculant le cas échéant l'espace de Riemann-Roch associé à $(D_1-D_2)-\frac{1}{n}(D_1'-D_2')$ à l'aide de l'un des algorithmes présentés dans l'annexe \ref{subsec:rr}. 
\end{rk}

La proposition suivante montre comment la décomposition d'un élément de $\HH^1(U,\mu_n)$ dans une base de ce $\Lambda$-module libre se traduit en termes de corps de fonctions. Elle servira dans la section \ref{subsec:cohlisalgcomp}.
\begin{lem}\label{lem:H1affdec} Notons $A_i=([D_i],D'_i,f_i),i=1,\dots,r$ les éléments d'une base de $\HH^1(U,\mu_n)$ avec les notations ci-dessus. Soit $A=([D],D',f)\in \HH^1(U,\mu_n)$.  Soient $\alpha_1,\dots,\alpha_r\in \Lambda$ tels que
\[ A=\sum_{j=1}^r\alpha_{i}A_i\in \HH^1(U,\mu_n).\]
Alors il existe une fonction $h\in k(X)^\times$ et des entiers $a_i$ relevant les $\alpha_i$ tels que \[ f=h^n\prod_{i=1}^rf_i^{a_i}.\]
\begin{proof} D'après la remarque \ref{rk:egalH1}, il existe $h\in k(X)^\times$, des entiers $a_i$ relevant les $\alpha_i$ et $D''\in \Div^0_Z(X)$ tels que \[D=\sum_{i=1}^r a_iD_i+D''+\div(h), \qquad D'=\sum_{i=1}^ra_iD_i'+nD''\quad \text{et}\quad \div(f)=nD-D''.\] Sur la courbe projective lisse $X$, il y a donc une égalité de diviseurs :
\begin{align*}
\div(f)&=\sum_{i=1}^r a_{i}(nD_i-D_i')+\div(h^n) \\
&=\sum_{i=1}^r a_{i}\div(f_i)+\div(h^n) \\
&=\div\left(h^n\prod_{i=1}^r f_j^{a_{i}}\right).
\end{align*}
Par conséquent, il existe $C\in k^\times$ tel que $Ch^n\prod_i f_i^{a_i}=f$. Comme $k$ est algébriquement clos, il existe $c\in k^\times$ tel que $c^n=C$ ; il suffit de remplacer $h$ par $ch$ pour conclure.
\end{proof}
\end{lem}

\begin{rk}\label{rk:compdec} \'{E}tant donné les éléments $A,A_1,\dots,A_r,\alpha_1,\dots,\alpha_r$ du lemme, le diviseur $D''$ de la démonstration se calcule comme $\frac{1}{n}(D'-\sum a_iD_i')$. Le calcul de la fonction $h$ nécessite le calcul de l'espace de Riemann-Roch du diviseur $nD-D''$ (différence de deux diviseurs effectifs de même degré inférieur à $g+n$), ainsi que le calcul dans $k$ d'une racine $n$-ième de $C$ (voir annexe \ref{subsec:racn} dans le cas des corps finis).
\end{rk}

\subsection{Calcul de la cohomologie}\label{subsec:calculcohaff}

Voyons comment calculer explicitement les groupes $\HH^1(U,\mu_n)$ et $\HH^1_c(U,\mu_n)$.
Le lemme \ref{lem:H1c} assure que $\HH^1_c(U,\mu_n)$ est isomorphe au groupe des classes d'équivalence (modulo $X-U$) de diviseurs sur $U$. Rappelons également qu'il y a des suites exactes courtes fonctorielles en $(X,U)$ :
\[
0\to \HH^1(X,\mu_n)\to \HH^1(U,\mu_n)\to \Div_Z^0(X)\otimes\Lambda\to 0 
\]
dont les flèches sont données respectivement par $([D],f)\mapsto ([D],0,f)$ et $([D],D',f)\mapsto D'$, et 
\[
0\to \frac{\bigoplus_{z\in Z}\mu_n(k)}{\mu_n(k)}\to \HH^1_c(U,\mu_n)\to \HH^1(X,\mu_n)\to 0 
\]
dont la flèche de gauche associe à $(\zeta_P)_{P\in Z}$ la classe du diviseur d'une fonction $f$ vérifiant $f(P)=\zeta_P$ pour tout $P\in Z$, et celle de droite associe à un diviseur $D$ le couple $([D],f)$ où $f$ est n'importe quelle fonction de diviseur $nD$. Comme les termes de gauche et de droite de ces suites exactes sont des $\Lambda$-modules libres, $\HH^1(U,\mu_n)$ et $\HH^1_c(U,\mu_n)$ sont des $\Lambda$-modules libres.
Ces deux suites sont duales l'une de l'autre par l'accouplement de Weil, comme vu à la fin de la section \ref{subsubsec:courbcohsupp}. Soit $\Gamma$ un quotient de $\mathfrak{G}_0=\Gk$ tel que l'action de $\mathfrak{G}_0$ sur $\HH^1(X,\mu_n)$ et $\Div_Z^0(X)\otimes\Lambda$ se factorise par $\Gamma$. Il suffit, afin de disposer d'une description complète de $\HH^1(U,\mu_n)$ comme $\Lambda[\Gamma]$-module, de calculer explicitement une section $\Lambda$-linéaire de $\HH^1(U,\mu_n)\to\Div_Z^0(X)\otimes\Lambda$, ainsi que l'action de $\Gamma$ sous la forme décrite dans la section \ref{sec:scindage}. Une telle section existe toujours car les $\Lambda$-modules en question sont libres. Cela revient, étant donné $P,Q\in Z$, à calculer un diviseur $D\in \Div^0(X)$ tel que $nD$ soit linéairement équivalent à $E\coloneqq (P)-(Q)$, ce qui s'effectue à l'aide de l'adaptation de l'algorithme de Huang-Ierardi décrite dans la remarque \ref{rk:huangrac}. 

\begin{rk}\label{rk:actgalcst} Voici comment calculer l'action de $\Gamma$ sur $\HH^1(U,\mu_n)$ en suivant la méthode décrite dans la section \ref{sec:scindage}. Soit $(([D_1],f_1),\dots,([D_{2g}],f_{2g}))$ une $\Lambda$-base de $\HH^1(X,\mu_n)$. Soit également $(D'_{2g+1},\dots,D'_{2g+r})$ une $\Lambda$-base de $\Div^0_Z(X)\otimes\Lambda$. Pour chaque $i\in\{2g+1\dots 2g+r\}$, soit $([D_i],D'_i,f_i)$ un représentant d'un antécédent de $D'_i$ dans $\HH^1(U,\mu_n)$. La famille des classes de triplets
\[ ([D_1],0,f_1),\dots,([D_{2g}],0,f_{2g}),([D_{2g+1}],D'_{2g+1},f_{2g+1}),\dots,([D_{2g+r}],D'_{2g+r},f_{2g+r})\]
forme alors une base de $\HH^1(U,\mu_n)$. Soit $([D],D',f)\in \HH^1(U,\mu_n)$. L'action d'un $\sigma\in\Gamma$ sur le triplet $([D],D',f)$ se calcule de la façon suivante : l'élément $\sigma\cdot D'\in \Div^0_Z(X)\otimes\Lambda$ se décompose comme combinaison linéaire \[\sigma\cdot D'=\sum_{i=2g+1}^{2g+r}\alpha_iD'_i.\]
L'élément $\sigma\cdot([D],D',f)-\sum_i\alpha_i([D_i],D'_i,f_i)$ appartient à l'image de $\HH^1(X,\mu_n)$ : il suffit maintenant de décomposer $(\sigma\cdot [D]-\sum_i\alpha_i [D_i],(\sigma\cdot f)/\prod_if_i^{\alpha_i})$ dans la base $([D_i],f_i)_{1\leqslant i\leqslant 2g}$ de $\HH^1(X,\mu_n)$ pour obtenir la décomposition complète de $\sigma\cdot ([D],D',f)$ dans la base $([D_i],D'_i,f_i)_{1\leqslant i\leqslant 2g+r}$.
\end{rk}

\begin{ex} Soit $E$ une courbe elliptique sur $k$. Soit $C=E-Z$ où $Z=\{ P_1,\dots,P_r\}$. Le $\Lambda$-module $\Div^0_Z(E)\otimes\Lambda$ est engendré par les diviseurs de la forme $(P_i)-(P_r)$ pour $i=1,\dots,r-1$. Pour chaque $i\in \{1,\dots,r\}$, soit $Q_i$ un point tel que $nQ_i=P_i$. Ces points s'obtiennent concrètement de la façon suivante. La multiplication par $n$ sur $E$ est donnée sur les abscisses des points par le polynôme de $n$-division $\phi_n\in k[x]$. Il suffit alors de résoudre $\phi_n(x)=x_{P_i}$ ; choisissons une solution $x$. L'un des deux points de $E$ ayant $x$ pour abscisse convient. Ayant ainsi obtenu des points $Q_1,\dots,Q_r$, posons $D_i=(Q_r)-(Q_i)$. Le diviseur $(P_i)-(P_r)+nD_i$ est alors principal, car de degré et de somme nuls. C'est un antécédent dans $\HH^1(C,\mu_n)$ de $(P_i)-(P_r)$.
\end{ex}

Voyons comment déterminer explicitement le $\Lambda[\mathfrak{G}_0]$-module $\HH^1_c(U,\mu_n)$. Supposons $X_0$ décrite par un modèle plan birationnel $C_0$, tel que la restriction de $X\to C_0\times_{k_0}k$ à $Z$ soit un isomorphisme. Décrivons d'abord l'inclusion de $(\bigoplus_{z\in Z}\mu_n(k))/\mu_n(k)$ dans $\HH^1_c(U,\mu_n)$. Notons $P_i=(x_i,y_i)$, $i=1,\dots,r$ les points de $Z\subset C$, que nous pouvons tous supposer dans la carte affine $z\neq 0$. Soit $(\zeta_{P_1},\dots,\zeta_{P_r})\in \mu_n(k)^Z$. La fonction \[ f(x,y)=\sum_{i=1}^r \zeta_{P_i}\prod_{j\mid x_j\neq x_i} \frac{x-x_j}{x_i-x_j}\prod_{j\mid y_j\neq y_i}\frac{y-y_j}{y_i-y_j} \tag{$\diamond$}\]
vérifie $f(P_i)=\zeta_{P_i}$ pour tout $i\in \{1\dots r\}$. L'image de $(\zeta_{P_1},\dots,\zeta_{P_r})$ dans $\HH^1_c(U,\mu_n)$ est alors $\div(f)$. D'autre part, un antécédent dans $\HH^1_c(U,\mu_n)$ d'un élément de $\HH^1(X,\mu_n)$ représenté par un diviseur $D$ se détermine de la façon suivante. Commençons par calculer un diviseur $D'$ équivalent à $D$ de support disjoint de $Z$, comme décrit dans l'annexe \ref{subsec:moving}. Le diviseur $nD'$ est le diviseur d'une fonction $f$. Calculons, pour $i=1,\dots,r$, une racine $n$-ième $\lambda_i$ de $f(P_i)$ dans $k$, puis une fonction $g$ telle que $g(P_i)=\lambda_i$ pour tout $i$. Le diviseur $D'-\div(g)$ est un antécédent de $D$ dans $\HH^1_c(U,\mu_n)$.  En utilisant l'algorithme de Couveignes décrit dans la section \ref{sec:couv}, nous obtenons le résultat suivant.

\begin{prop}
Il existe un algorithme probabiliste (Monte-Carlo) qui, étant donné une courbe projective lisse $X_0$ de genre $g$ sur $\FF_q$ représentée par un modèle plan de degré $d$, un entier $n$ premier à $q$, un diviseur $\FF_q$-rationnel de degré 1 sur $X$, le polynôme caractéristique de la fonction zêta de $X_0$ et un ensemble non vide $Z=\{P_1,\dots, P_r\}\subset X_0(\FF_q)$ calcule un ensemble de diviseurs $D_1,\dots,D_{2g+r-1}$ de degré 0 sur $(X_0-Z)_{\overline{\FF_q}}$ dont les classes forment une base de $\HH^1_c((X_0-Z)_{\overline{\FF_q}},\mu_n)$. Le nombre d'opérations effectué est polynomial en $d$, $\log q$, $n^{2g}$ et $r$. 
\begin{proof}
Quitte à passer à une extension de  degré $d$ de $\FF_q$ pour trouver un point $P_0\in X(\FF_q)$, l'algorithme de Couveignes renvoie des classes $D_1,\dots,D_{2g}$ de diviseurs de degré 0 sur $X$ de la forme $G_i-gP_0$, où les $G_i$ sont définis sur une extension de $\FF_q$ de degré $O(gn^{2g})$, en temps polynomial en $d$, $n^{2g}$ et $\log q$. Le calcul d'un diviseur équivalent à $G_i-gP_0$ de support disjoint de $Z$ se fait donc en temps polynomial en $d$, $n^{2g}$, $\log q$ et $r$ d'après l'annexe \ref{subsec:moving}. Les diviseurs $D_{2g+1},\dots,D_{2g+r-1}$ sont obtenus en calculant le diviseur de fonctions données par la formule $(\diamond)$, qui sont de degré $O(r)$.
\end{proof}
\end{prop}

\begin{rk} Le groupe $\HH^1_c(U,\mu_n)$ se détermine très facilement une fois que l'on connaît une description de $\HH^1(X,\mu_n)$ en termes de diviseurs. Ce n'est pas le cas de $\HH^1(U,\mu_n)$, dont le calcul nécessite de déterminer des racines $n$-ièmes dans $\Pic(X)$. Si la dualité de Poincaré permet de décrire l'un des groupes à partir de l'autre, elle ne permet pas de déterminer une fonction $f$ telle que $\div(f)$ soit un antécédent dans $\HH^1(U,\mu_n)$ d'un diviseur de degré 0 supporté sur $Z$.
\end{rk}

\subsection{Fonctorialité en $U$}

Soit $U'$ une courbe affine lisse sur $k$ munie d'un $k$-morphisme $f\colon U'\to U$. Soient $X'$ la compactification lisse de $U'$, et $Z'$ le fermé réduit complémentaire de $U'$ dans $X'$. Considérons également $V'=U\times_{X}X'$, et son fermé réduit complémentaire $W'$ dans $X'$. 
Le diagramme commutatif suivant, dont les deux rectangles sont cartésiens, résume la situation.

\[
\begin{tikzcd}
U'\arrow[r]\arrow[dr,"f",swap]& V'\arrow[r]\arrow[d]& X'\arrow[d]& \arrow[l] Z'& \arrow[l]W'\arrow[d] \\
& U \arrow[r] & X & &\arrow[ll] Z
\end{tikzcd}
\]

Décrivons comment calculer le morphisme $\HH^1(U,\Lambda)\to \HH^1(U',\Lambda)$ induit par $f$. 
La fonctorialité de la suite exacte précédente permet de calculer le morphisme $\HH^1(U,\Lambda)\to \HH^1(V',\Lambda)$, c'est-à-dire la composée des morphismes de $\Lambda$-modules \[ \HH^1(U,\Lambda) \xrightarrow{\sim}\HH^1(X,\Lambda)\oplus \left(\Div^0_Z(X)\otimes\Lambda\right)\to \HH^1(X',\Lambda)\oplus (\Div^0_{W'}(X')\otimes\Lambda)\xrightarrow{\sim} \HH^1(V',\Lambda)\]  \\
où le calcul de l'isomorphisme de droite repose sur celui d'une section $\Lambda$-linéaire $s'$ de \[\HH^1(V',\Lambda)\to \Div^0_{W'}(X')\otimes\Lambda\] compatible à $f$. 
Pour faire cela, notons $s$ une section déjà calculée de $\HH^1(X,\mu_n)\to \Div^0_Z(X)\otimes\Lambda$ et $\gamma \colon  \Div^0_Z(X)\to \Div^0_{W'}(X')$. 
Pour chaque diviseur $w=\gamma(v)\in\im(\gamma)$, définissons $s'(w)\coloneqq f^\star s(v)$. Pour les $w$ qui ne sont pas dans l'image de $\gamma$, définissons $s'(w)$ comme étant une racine $n$-ième de $w$ dans $\Pic^0(X')$.\\

Remarquons que la flèche $\Div^0_Z(X)\otimes\Lambda\to \Div^0_{W'}(X')\otimes\Lambda$ dépend des indices de ramification du morphisme : notons $|Z|=\{z_1,\dots,z_n\}$, $|W|=\{w_1,\dots,w_d\}$ avec $f(w_i)=\colon z_{\delta_i}$. L'image de $(f_0,\dots,f_s)\in \HH^0(Z,\Lambda)$ est alors $(e_{w_i}(f)z_{\delta_i})_{1\leqslant i\leqslant d}$.  \\

Le morphisme $\HH^1(U,\Lambda)\to \HH^1(U',\Lambda)$ se factorise par $\HH^1(U,\Lambda)\to \HH^1(V',\Lambda)$. Il reste à déterminer la flèche $\HH^1(V',\Lambda)\to \HH^1(U',\Lambda)$.  Remarquons qu'avec la description explicite de $\HH^1(V',\Lambda)$ et $\HH^1(U',\Lambda)$ dont nous disposons, le morphisme $\HH^1(V',\Lambda)\to \HH^1(U',\Lambda)$ se calcule explicitement et est simplement une inclusion. \`{A} un triplet $([D],D',f)$ où $[D]\in \Pic (X')[n], D'\in \Div^0_{Z'}(X')\otimes\Lambda$ et $nD=D'+\div(f)$, il associe encore $([D],D',f)$. Il est donc possible de calculer explicitement un complémentaire de $\HH^1(V',\Lambda)$ dans $\HH^1(U',\Lambda)$, et ainsi représenter le morphisme composé \[ \HH^1(U,\Lambda)\to \HH^1(V',\Lambda)\hookrightarrow \HH^1(U',\Lambda). \]

\section{Cohomologie des faisceaux lisses}

Cette section est dédiée au calcul du complexe de cohomologie d'un faisceau lisse sur une courbe lisse ou nodale sur un corps algébriquement clos. Dès que son genre géométrique est non nul, une telle courbe est un $K(\pi,1)$, et notre méthode consiste à calculer la cohomologie des faisceaux lisses comme cohomologie du groupe d'automorphismes d'un revêtement galoisien de la courbe. Nous commençons par décrire cette méthode, puis nous détaillons les algorithmes employés pour la mettre en œuvre et calculons leur complexité. Nous montrons également comment calculer explicitement les morphismes de la suite de Gysin. Nous donnons ensuite deux exemples de calcul de groupes de cohomologie. Nous terminons par la description de la généralisation de cette méthode à la détermination de la cohomologie d'un complexe de faisceaux lisses sur une telle courbe.

\subsection{Calcul de la cohomologie}\label{subsec:cohlisG}
Rappelons que $k$ est un corps algébriquement clos, et $\Lambda=\ZZ/n\ZZ$ avec $n$ premier à la caractéristique de $k$.
Soit $U$ une courbe intègre lisse ou nodale de genre géométrique non nul sur $k$. Soit $\F$ un faisceau lisse de $\Lambda$-modules sur $U$, de fibre générique géométrique $M$. D'après la proposition \ref{prop:Kpi1}, $U$ est un $K(\pi,1)$ : le morphisme canonique
\[ \RG(\pi_1(U),M)\to \RG(U,\F) \]
est un isomorphisme dans $\DD^b_c(\Lambda)$. Soit $f\colon V\to U$ un revêtement qui trivialise $\F$. Considérons un revêtement caractéristique $W\to V$ trivialisant tous les $f^\star\F$-torseurs sur $V$, par exemple celui construit dans la section \ref{subsubsec:X2}. Le revêtement $W\to U$ est encore galoisien. Notons $G=\Aut(W|U)$. 

\begin{prop}\label{prop:RGlisse} Le morphisme canonique \[ \tau_{\leqslant 1}\RG(G,M)\to \tau_{\leqslant 1}\RG(\pi_1(U),M)\xrightarrow{\sim}\tau_{\leqslant 1}\RG(U,\F) \]
dans $\DD^b_c(\Lambda)$ induit par le quotient $\pi_1(U)\to G$ est un isomorphisme dans $\DD^b_c(\Lambda)$.
\begin{proof}
La suite spectrale de Hochschild-Serre fournit la suite exacte courte \[ 0\to \HH^1(G,M)\to \HH^1(U,\F)\to \HH^0(G,\HH^1(W,\F|_{W})).\]
Comme $W\to V$ trivialise tous les $M$-torseurs, le morphisme $\HH^1(U,\F)\to \HH^1(W,\F|_{W})$ est nul. Par conséquent, le morphisme 
\[ \HH^1(G,M)\to \HH^1(U,\F) \]
induit par le quotient $\pi_1(U)\to G$ est un isomorphisme. De même, le morphisme $\HH^0(G,M)\to \HH^0(U,\F)$ est un isomorphisme. Ceci implique que le morphisme composé
\[ \RG(G,M)\to \RG(\pi_1(U),M)\xrightarrow{\sim}\RG(U,\F) \]
dans $\DD^b_c(\Lambda)$ est un quasi-isomorphisme en degrés 0 et 1. 
\end{proof}
\end{prop}

\begin{rk}\label{rk:mayerRG}
Dans le cas où $U$ est affine, $\HH^i(U,\F)=0$ dès que $i\geqslant 2$. Par conséquent, le morphisme \[ \tau_{\leqslant 1}\RG(G,M)\to \tau_{\leqslant 1}\RG(\pi_1(U),M)\to  \tau_{\leqslant 1}\RG(U,\F) \to \RG(U,\F)\]
est un isomorphisme dans $\DD^b_c(\Lambda)$.

Dans le cas où $U$ est projective, choisissons deux ouverts $U_1,U_2$ de $U$. Le triangle de Mayer-Vietoris (voir section \ref{subsec:mayerviet}) assure alors que
\[\RG(U,\F)=\cone\left((\RG(U_1,\F)\oplus \RG(U_2,\F)\to \RG(U_1\cap U_2,\F)\right)[-1].\] 
\end{rk}

\paragraph{Fonctorialité sur $\Spec k$}  Soit $\phi\colon U'\to U$ un morphisme de courbes affines intègres lisses sur $k$. Soit $\F$ un faisceau lisse sur $U$. Expliquons comment calculer le morphisme $\RG(U,\F)\to \RG(U',\phi^\star\F)$ déduit de $\phi$ par fonctorialité. Soit $f\colon V\to U$ un revêtement galoisien de $U$ qui trivialise $\F$. Ici, nous supposons que $W$ est le revêtement $V_2$ de $V$ de groupe $H^1(V,\Lambda)^\vee$ construit dans la section \ref{subsubsec:X2}. Calculons une décomposition primaire de $V\times_U U'$ (voir annexe \ref{sec:prodfib}). 
Considérons une composante connexe $V'$ de $V\times_U U'$. Notons $K$, $K'$ les corps de fonctions respectifs de $V$, $V'$. Quitte à le remplacer par sa clôture galoisienne sur $U$, supposons $V'\to U$ galoisien.
Le schéma $W$ est la normalisation de $V$ dans $K(\sqrt[n]{g_1},\dots,\sqrt[n]{g_r})$, où les $g_i$ ont un diviseur multiple de $n$ dans $\Pic(X)$ ; considérons de même le revêtement $W'$ de $V'$ de groupe $\HH^1(V',\Lambda)^\vee$. Notons $L$, $L'$ les corps de fonctions respectifs de $W$, $W'$. Soit $T'\to V'$ la normalisation de $V'$ dans $K'(\sqrt[n]{g_1},\dots,\sqrt[n]{g_r})\subseteq L'$. 
Le diagramme commutatif suivant, dont les deux flèches verticales composées sont des revêtements galoisiens, résume la situation.
\[ \begin{tikzcd}
W' \arrow[d] & \\
T' \arrow[d]\arrow[r] & W \arrow[d] \\
V' \arrow[d] \arrow[r] & V \arrow[d,"f"]\\
U' \arrow[r,"\phi"] & U
\end{tikzcd} \]
Le morphisme $W'\to W$ ainsi obtenu induit un morphisme $u\colon \Aut(W'|U')\to\Aut(W|U)$. En effet, comme $L/k(U)$ est normale, tout élément de $\Aut(L'|k(U'))\subseteq\Aut(L'|k(U))$ donne par restriction un élément de $\Aut(L|k(U))=\Aut(W|U)$. Ce morphisme $\Aut(W'|U')\to \Aut(W|U)$ donne par fonctorialité le morphisme cherché entre les complexes de $\Lambda$-modules représentant \[\RG(\Aut(W|U),M)\to \RG(\Aut(W'|U'),M).\]

\begin{rk}
Le morphisme $\HH^1(U,\F)\to \HH^1(V,\F|_V)$ est en particulier très simple à calculer : c'est simplement le morphisme $\HH^1(G,M)\to \HH^1(\HH^1(V,\Lambda)^\vee,M)$ déduit de l'inclusion $\HH^1(V,\Lambda)^\vee\subset G$.
\end{rk}

\paragraph{Action de $\mathfrak{G}_0=\Gk$} Supposons désormais que $U,\F$ proviennent par changement de base d'une courbe géométriquement connexe $U_0$ sur $k_0$ et d'un faisceau lisse $\F_0$ sur $U_0$ trivialisé par $f_0\colon V_0\to U_0$. Notons $V=V_0\times_{k_0}k$. \begin{enumerate}
\item Traitons d'abord le cas où $V$ est connexe et possède un $k_0$-point $y_0$. Soit $\bar y_0$ un point géométrique de $V$ au-dessus de $y_0$, d'image un point géométrique $\bar x_0$ de $U$. Soit toujours $W=V_2$ le revêtement de $V$ de groupe $\HH^1(V,\Lambda)^\vee$. Notons $\bar z_0$ un point géométrique de $W$ d'image un point géométrique $\bar y_0$ au-dessus de $y_0$.  Pour tout $\sigma\in \mathfrak{G}_0$, l'automorphisme \[ \sigma_\star\colon \pi_1(U,\bar x_0)\to \pi_1(U,\bar x_0) \] se restreint en un automorphisme de $\pi_1(V,\bar y_0)$, et donc encore en un automorphisme de $\pi_1(W,\bar z_0)$ puisque ce dernier est caractéristique dans $\pi_1(V)$. Par conséquent, il induit par passage au quotient un automorphisme de $\Aut(W|U)$. L'action de $\mathfrak{G}_0$ sur $\tau_{\leqslant 1}\RG(\Aut(W|U),M)$ découle alors de son action sur $\Aut(W|U)$ par fonctorialité. Cette méthode sera illustrée dans la section \ref{subsec:exdetlis}.
\item Supposons désormais que $V$ n'est pas connexe, c'est-à-dire que $k_0(V_0)$ contient une extension finie non triviale de $k_0$. Construisons le revêtement $V_{2,0}\to V_0$ défini sur $k_0$ décrit dans la section \ref{subsec:X2gal}. Notons $W=V_{2,0}\times_{k_0}k$. Alors $W\to U$ est encore un $G=\Aut(V_{2,0}|U_0)$-torseur, et le même raisonnement permet de calculer $\tau_{\leqslant 1}\RG(\Aut(V_{2,0}|U_0),\HH^0(W,\F))$. Remarquons que les composantes connexes de $W$ sont permutées par un groupe $\Gal(L'|k_0)$, où $L'$ est une extension galoisienne de $k_0$. Ceci permet de calculer l'action de $\mathfrak{G}_0$ sur $\HH^0(W,\F)$. Le groupe $\mathfrak{G}_0$ agit ici trivialement sur $\Aut(V_{2,0}|U_0)$, puisque $V_{2,0}$ est défini sur $k_0$. Soient $V'$ une composante connexe de $V$ et $W'$ une composante connexe de $W$ d'image $V'$. Il y a une suite exacte \[ 1\to \Aut(W'|V')\to G\to \Gal(L'|k_0)\to 1\]
et \[ \HH^0(W,\F)=\ind_{\Aut(W'|V')}^G \HH^0(W',\F).\]
Le complexe $\RG(G,\HH^0(W,\F))$ ainsi calculé, muni d'une action naturelle de $\mathfrak{G}_0$, représente bien $\RG(U,\F)$, car le lemme de Shapiro \cite[Th. 4.19]{neukirch_cft} assure qu'il y a un isomorphisme \[\RG(G,\HH^0(W,\F))=\RG(\Aut(W'|V'),\HH^0(W',\F)).\]
\end{enumerate}

\paragraph{Mise en garde} Il serait tentant, lorsque $U$ est affine, de faire le raisonnement suivant. Le morphisme $V\to U$ étant galoisien, $\RG(U,\F)=\RG(\Aut(V|U),\RG(V,\F|_V))$. Comme $V$ est un $K(\pi,1)$, un complexe représentant $\RG(V,\F|_V)$ est le complexe de cochaînes usuel calculant $\RG(\pi_1(V,y),M)$. Il peut même être tronqué en degré $\leqslant 1$ puisque $V$ est affine. Le morphisme \[\tau_{\leqslant 1}\RG(\Aut(V_2|V),M)\to \tau_{\leqslant 1}\RG(\pi_1(V,y),M)\] est encore un quasi-isomorphisme par les arguments ci-dessus.
Comme $\F$ est constant sur $V$, l'action de $\pi_1(V,y)$ sur $M$ est triviale ; par conséquent, la première flèche du complexe calculant $\RG(\Aut(V_2|V),M)$ est nulle. Il en découle que $\RG(V,\F|_V)$ est représenté par le complexe \[\HH^0(V,\F|_V)\xrightarrow{0} \HH^1(V,\F|_V).\] L'action de $\Aut(V|U)$ sur $\Aut(V_2|V)=\HH^1(V,\Lambda)^\vee$, et donc sur ce complexe, se calcule aisément. Ceci permettrait de déterminer $\RG(\Aut(V|U),\RG(V,\F|_V))$ sans avoir à calculer de revêtement de $V$. L'erreur est la suivante : le complexe $\HH^0(V,\F|_V)\xrightarrow{0} \HH^1(V,\F|_V)$ représente $\RG(V,\F|_V)$ dans $\DD^b_c(\Lambda)$, et ses groupes de cohomologie sont munis d'une action naturelle de $\Aut(V|U)$. Mais ce complexe \textit{ne représente en général pas} $\RG(V,\F|_V)$ dans $\DD^b_c(\Lambda[\Aut(V|U)])$. Lorsque $[V:U]$ est divisible par $n$, ce n'est pas le cas. Prenons l'exemple de $V=U=\GG_m$ et $f\colon V\to U, x\mapsto x^n$. Considérons le faisceau constant $\Lambda$ sur $U$. Alors \[\RG(V,\Lambda)=[\Lambda\xrightarrow{0} \Lambda] \] et l'action de $\Aut(V|U)\simeq \mu_n(k)$ sur $\HH^0(V,\Lambda)$ et $\HH^1(V,\Lambda)$ est triviale puisque $\mu_n(k)$ fixe les points $0$ et $\infty$. Par conséquent, $\RG(V,\Lambda)$ est le complexe de $\mu_n(k)$-modules triviaux $\Lambda[0]\oplus\Lambda[-1]$, et pour tout entier $i\geqslant 0$, $\HH^i(U,\Lambda)=\HH^i(\mu_n(k),\Lambda)\oplus \HH^{i-1}(\mu_n(k),\Lambda)$. Ce résultat est absurde puisque ces groupes doivent être nuls en degré $i\geqslant 2$, et que $\HH^1(U,\Lambda)$ est de rang 1 et non 2.

\subsection{Algorithmes et complexité}\label{subsec:cohlisalgcomp}

Nous allons décrire ici les algorithmes permettant de calculer explicitement la cohomologie d'un faisceau lisse par la méthode décrite dans la section précédente.
Rappelons les notations : $U_0$ est une courbe lisse géométriquement connexe sur le corps parfait $k_0$, et $V_0\to U_0$ est un revêtement galoisien. Supposons d'abord par simplicité que $V_0$ est géométriquement connexe. Les courbes $U,V$ sont les changements de base de $U_0,V_0$ à la clôture algébrique $k$ de $k_0$. Le revêtement caractéristique $V_2\to V$ défini dans la section \ref{subsubsec:X2} a pour groupe $H^1(V,\Lambda)^\vee$.

Pour les résultats de complexité, nous noterons $C(g,q,n,r)$ la complexité du calcul de $\HH^1(U,\mu_n)$, où $U$ est une courbe intègre lisse sur $\overline{\FF_q}$ provenant de $\FF_q$, de compactification lisse $X$ de genre $g$, avec $r=|X-U|$. Notons également $D(g,n,r)$ le degré de la plus petite extension de $\FF_q$ sur laquelle sont définis les éléments de $\HH^1(U,\mu_n)$ obtenus. Une majoration de ces deux entiers est donnée à la fin de la section \ref{subsubsec:detrep}.

\begin{rk}\label{rk:Dgnr}
Il est possible de majorer $D(g,n,r)$ par $|\Aut_{\Lambda}\HH^1(U,\mu_n)|=\OO( n^{(2g+r)^2})$. En effet, l'action de $\mathfrak{G}_0=\Gk$ sur $\HH^1(U,\mu_n)$ se factorise par un quotient $\mathfrak{G}_1$ d'ordre $|\Aut_{\Lambda}\HH^1(U,\mu_n)|$, c'est-à-dire le groupe de Galois d'une extension de $\FF_q$ d'ordre $\OO(n^{(2g+r})^2)$ que l'on peut fixer. Les classes dans $\Pic(X)$ des diviseurs obtenus représentant des éléments de $\HH^1(U,\mu_n)$ sont toutes $\mathfrak{G}_1$-invariantes, et l'algorithme de la proposition \ref{prop:divrat} permet alors de trouver des diviseurs $\mathfrak{G}_1$-invariants qui leur sont équivalents ; la complexité de cette opération est négligeable devant celle des algorithmes de calcul du $\HH^1$.
\end{rk}

\paragraph{Calcul de $\Aut(V_2|U)$} Soit $\tau\in\Aut(k(V)/k(U))$. Notons $t=2g+r$. Il y a $n^t$ automorphismes de $k(V_2)$ de restriction $\tau$ ; voici comment en construire un. Soit \[ \mathcal{B}=(([D_1],D'_1,g_1),\dots ,([D_t],D'_t,g_t))\] une $\Lambda$-base de $\HH^1(V,\mu_n)$. Fixons $f_i=\sqrt[n]{g_i}$. Un automorphisme $\sigma$ de $k(V_2)$ de restriction $\tau$ vérifie nécessairement $\sigma(f_i)^n=\tau(g_i)$. Calculons une racine $n$-ième de $\tau(g_i)$. Décomposons l'élément $\tau^\star([D_i], D'_i, g_i)$ de $\HH^1(V,\mu_n)$ dans la base $\mathcal{B}$ : cela revient à calculer des éléments $\alpha_{ij}\in \Lambda$ tels que \[ (\tau^\star [D_i],\tau^\star D'_i,\tau(g_i))=\sum_{j=1}^t \alpha_{ij}([D_j],D'_j,g_j)\] dans $\HH^1(V,\mu_n)$.
Il existe d'après le lemme \ref{lem:H1affdec} une fonction $h_i\in k(V)^\times$ et des entiers $a_{ij}$ relevant les $\alpha_{ij}$ tels que \[ \left(h_i\prod_{j=1}^t f_j^{a_{ij}}\right)^n=\tau(g_i)\] dans $k(V)$. Une racine $n$-ième de $\tau(g_i)$ dans $k(V_2)$ s'obtient alors sous la forme \[ f_{\tau,i}\coloneqq h_i\prod_{j=1}^t f_j^{a_{ij}}.\]
Soit $(b_1,\dots,b_s)$ une $k(U)$-base de $k(V)$. Alors \[ (b_if_1^{a_1}\dots f_t^{a_t})_{1\leqslant i \leqslant s, 0\leqslant a_j\leqslant n-1}\] est une $k(U)$-base de $k(V_2)$. L'endomorphisme $\sigma$ de $k(V_2)$ défini par \[ \sigma(b_if_1^{a_1}\dots f_t^{a_t})=\tau(b_i)f_{\tau,1}^{a_1}\dots f_{\tau,t}^{a_t}\] est un élément de $\Aut(k(V_2)|k(U))$ dont la restriction à $k(V)$ vaut $\tau$. Les autres antécédents de $\tau$ dans $\Aut(k(V_2)|k(U))$ sont obtenus en composant celui-ci avec un élément de $\Aut(k(V_2)|k(V))\simeq \Hom_\Lambda(\HH^1(V,\mu_n),\mu_n)$.

L'algorithme pour calculer $\sigma$ est donc le suivant : pour chaque $i\in \{ 1\dots t\}$, calculer $\tau^\star ([D_i],D'_i,g_i)$, puis déterminer $h_i$ grâce à la preuve du lemme \ref{lem:H1affdec}, et en déduire $f_{\tau,i}$.
Rappelons (voir remarque \ref{rk:compdec}) que le calcul de $h_i$ nécessite des calculs d'espaces de Riemann-Roch de diviseurs sur $V$ qui sont chacun différence de deux diviseurs effectifs de même degré $g+n$, ainsi que le calcul d'une racine $n$-ième dans $k$. La complexité du calcul de $\Aut(V_2|U)$, qui est d'ordre $\deg(V\to U)n^t$, est donc dominée par celle du calcul des $f_{\tau,i}$, c'est-à-dire de $\HH^1(V,\mu_n)$.\\

Considérons désormais le cas général du calcul de $\RG(U,\F)$ : nous ne supposons plus $V_0$ géométriquement connexe. Le revêtement galoisien calculé, que nous noterons $V_{2,0}\to U_0$, est celui décrit dans la section \ref{subsec:X2gal}. Le diagramme suivant, dont le carré est cartésien, résume les notations.
\[
\begin{tikzcd}
 & V_{2,0} \arrow[d] \\
V \arrow[d]\arrow[r] & V_0 \arrow[d] \\
U \arrow[r]& U_0
\end{tikzcd}
\]

\paragraph{Ordre de $\Aut(V_{2,0}|U_0)$} Nous avons vu dans la section \ref{subsec:X2gal} comment le déduire du calcul de $\Aut(V_2|U)$. Calculons son ordre lorsque $k_0=\FF_q$. Notons $r=|X-U|$, où $X$ est la compactification lisse de $U$. Rappelons que le corps de fonctions de $V_{2,0}$ contient une extension $L$ de $k_0$, qui est la composée de la clôture algébrique $k_1$ de $k_0$ dans $k_0(V_0)$, de $k_0(\mu_n)$ et de l'extension de $k_0$ par laquelle se factorise l'action de $\Gk$ sur $\HH^1(V_0\times_{k_0}k,\Lambda)$. Soit $V^{(1)}$ une composante connexe de $V_{2,0}\times_{k_0}k$ ; d'après la formule de Riemann-Hurwitz, son genre est $O((g+r)[V^{(1)}:U])=O((g+r)\frac{[V_0:U_0]}{[k_1:k_0]})$, où $g$ est le genre de $X$. Notons $d=\frac{[V_0:U_0]}{[k_1:k_0]}$.
D'après les lemmes \ref{lem:degext} et \ref{lem:corpsdef}, l'ordre de $G=\Aut(W_0|U_0)$ est donné par 
\[ \begin{array}{rcl} 
|G|&=&n^{2g(V^{(1)})}\times [L:k_1]\times [V_0:U_0] \\ 
&=&O\left(g(V^{(1)})nD(g(V^{(1)}),q,n,dr)\right)[V_0:U_0] \\
&=& O\left((g+r)dn[V_0:U_0] D(d(2g+r),n,dr)\right).
\end{array}
\] 
Rappelons que $[V_0:U_0]$ est l'ordre du groupe de monodromie $\mathfrak{S}\subseteq \Aut_\Lambda(M)$. De plus, la remarque \ref{rk:Dgnr} assure que $D(g,n,r)$ est majoré par $n^{2g+r}$. En particulier, si $M=\Lambda^j$, $|\Aut_\Lambda(M)|=O(n^{j^2})$ et \[ |G|\leqslant (2g+r) n^{2j^2+3(2g+r)n^{j^2}}.\]
Cette majoration donne une bonne idée de la complexité de l'algorithme dans le pire cas, mais cache le fait qu'à taille constante de groupe de monodromie, cette complexité est polynomiale en $j^3n^{2^g}$. Si la courbe $V$ admet un modèle plan à singularités ordinaires de degré $O(g)$ (ce qui est le cas pour une courbe générale), cette complexité est $O(j^3n^{g^4})$.

\paragraph{Calcul de $\RG(U,\F)$} En supposant que l'on ait déjà construit le revêtement $V_{2,0}\to U_0$ et calculé son groupe de Galois, il ne reste plus qu'à calculer le complexe $\RG(G,M)$ où $G=\Aut(V_{2,0}|U_0)$ grâce à la résolution libre usuelle (résolution bar) de $\Lambda$ comme $\Lambda[G]$-module. Le module $M=\HH^0(V,\F)$ est représenté comme quotient d'un module libre : $M=\Lambda^m/N$. Rappelons que nous souhaitons calculer la cohomologie de $G$ à valeurs dans un module $M'$ induit à partir de $M$, isomorphe comme $\Lambda$-module à $M^{[k_1:k_0]}$. Notons $s=[k_1:k_0]\leqslant [V_0:U_0]$.
Rappelons que l'ordre de $G$ est $d=[L:k_1]n^r[V_0:U_0]$, où $r=\rg_\Lambda(\HH^1(V^{(1)},\Lambda))$. Le complexe représentant $\tau_{\leqslant 1}\RG(U,\F)$ est donc une troncature de la dernière ligne du diagramme commutatif à colonnes exactes suivant.

\[
\begin{tikzcd} 0\arrow[d] & 0\arrow[d] & 0\arrow[d] \\
N \arrow[r]\arrow[d] & N^{d} \arrow[r]\arrow[d] & N^{d^2}\arrow[d] \\
\Lambda^{sm} \arrow[r] \arrow[d]& \Lambda^{dsm} \arrow[r]\arrow[d] & \Lambda^{d^2sm}\arrow[d] \\
M^s \arrow[r]\arrow[d] & M^{sd} \arrow[r]\arrow[d] & M^{sd^2}\arrow[d] \\
0 & 0 & 0
\end{tikzcd}
\]
Les calculs effectués pour déterminer $\RG(U,\F)$ étant simplement des calculs de noyaux et conoyaux de morphismes de $\Lambda$-modules engendrés par au plus $sd^2m$ éléments, la complexité du calcul de $\RG(U,\F)$ est $O((sd^2m)^3)$ opérations dans $\Lambda$.\\

\begin{theorem}\label{th:complis} Soit $U_0$ une courbe lisse géométriquement connexe sur $\FF_q$. Notons $U=(U_0)_{\overline{\FF_q}}$. Notons $X$ sa compactification lisse, $g$ son genre, et $r=|X-U|$. Soit $\F$ un faisceau lisse de $\Lambda$-modules sur $U_0$ de fibre générique géométrique un $\Lambda$-module $M$ à $m$ générateurs ; notons $\mathfrak{S}$ l'image de $\pi_1(U)$ dans $M$, et $d$ son ordre. Supposons donné un revêtement galoisien $V_0\to U_0$ de groupe $\mathfrak{S}$ qui trivialise $\F$. Il existe un algorithme calculant un complexe de $\Lambda[\mathfrak{G}_0]$-modules représentant $\RG(U,\F)$ en \[ O\left(C(d(g+r),q,n,dr)+(md^{3}n^{4d(g+2r)}D(d(2g+r),n,dr))^3\right)\]
opérations dans $\FF_q$.
\begin{proof} L'algorithme consiste à calculer le revêtement caractéristique $V_{2,0}$ de $V_0$ défini à partir de $\HH^1(V,\Lambda)$, son groupe d'automorphismes $G$ ainsi que $\RG(G,M')$ où $M'=\HH^0(V_{2,0}\times_{k_0}k,\F)$. Le rang de $\HH^1(V^{(1)},\Lambda)$ est inférieur à $2d(g+r)+dr-1=2d(g+2r)-1$. L'ordre du groupe $G=\Aut(V_{2,0}|U_0)$ est donc inférieur à $dn^{2d(g+2r)}D(g,n,r)$. L'estimation précédente du coût du calcul de $\tau_{\leqslant 1}\RG(G,M')$ à partir de $G$ permet de conclure.
\end{proof}
\end{theorem}

\begin{cor} Avec les notations et hypothèses du théorème, il existe un algorithme probabiliste (Las Vegas) qui calcule $\RG(U,\F)$ en
\[ \mathcal{P}(d,g,n,r,m,\log q)^{2^{O\left((d(g+r))^2\right)}}\]
opérations dans $\FF_q$, où $\mathcal{P}$ est un polynôme. Si $V$ admet un modèle plan ordinaire de degré $O(g)$, cette complexité devient \[ \mathcal{P}(d,g,n,r,m,\log q)^{O\left((d(g+r))^4\right)}.\]
\begin{proof}
Nous avons utilisé les majorations de $C(g,q,n,r)$ et $D(g,n,r)$ obtenues dans la section \ref{subsubsec:detrep} : ils sont tous les deux bornés par $\mathcal{P}(d,g,n)^{O(g^4)}$ pour une courbe à singularités ordinaires, et donc $\mathcal{P}(d,g,n)^{2^{O(4g^2)}}$ pour une courbe quelconque.
\end{proof}
\end{cor}

\subsection{Suite de Gysin}\label{subsec:gysinlisse}

Soit $f\colon Y\to X$ un revêtement galoisien de courbes intègres projectives lisses sur $k$. Soit $\F$ un faisceau lisse de $\Lambda$-modules sur $X$, trivialisé par $f$. Notons $F=\HH^0(Y,f^\star\F)$. Soit $Z$ un fermé réduit non vide de $X$. Notons $U$ son ouvert complémentaire, $W=Z\times_X Y$ et $V=U\times_X Y$. Les notations sont résumées par le diagramme commutatif suivant, dont les deux carrés sont cartésiens.
\[
\begin{tikzcd}
V \arrow[r]\arrow[d] & Y\arrow[d] & \arrow[l]\arrow[d] W \\
U \arrow[r] & X & \arrow[l] Z
\end{tikzcd}
\]
Les suites exactes de Gysin pour $\F$ sur $X$ et pour le faisceau constant $F$ sur $Y$ fournissent le diagramme commutatif
\[
\begin{tikzcd} 0 \arrow[r]& \HH^1(X,\F)\arrow[r]\arrow[d]&  \HH^1(U,\F) \arrow[r]\arrow[d]&  \HH^0(Z,\F(-1)) \arrow[r]\arrow[d]&  \HH^0(X,\F(-1)^\vee)^\vee \arrow[r]\arrow[d]&  0 \\
0 \arrow[r] & \HH^1(Y,F) \arrow[r] & \HH^1(V,F) \arrow[r] & \HH^0(W,F(-1))\arrow[r] & F(-1) \arrow[r] & 0
\end{tikzcd}\]
dont la ligne inférieure est celle de la section \ref{subsubsec:Gyscst}.\\

Le morphisme $\HH^0(Z,\F(-1))=\HH^0(Z,F(-1))\to \HH^0(W,F(-1))$ s'écrit \[ F(-1)^Z\to F(-1)^W, (a_z)_{z\in Z}\mapsto ((a_z)_{w\in W_z})_{z\in Z} \]
et admet une rétraction $r$ donnée par des projections. Ceci permet de calculer le morphisme \[ \HH^1(U,\F)\to \HH^0(Z,\F(-1))\] comme la composée
\[ \HH^1(U,\F)\to \HH^1(V,F)\to \HH^0(W,F(-1))\xrightarrow{r} \HH^0(Z,\F(-1)).\]
La flèche $\HH^0(Z,\F(-1))\to \HH^0(X,\F(-1)^\vee)^\vee$ associe à $(a_z)_{z\in Z}$ le morphisme $\HH^0(X,\F(-1)^\vee)\to \Lambda$ qui à $u \colon\F(-1)\to \Lambda$ associe $\sum_{z\in Z} u_z(a_z)$. La flèche verticale de droite associe à $u \colon \HH^0(X,\F(-1)^\vee)\to \Lambda$ l'élément $\deg(f)u|_Y \colon F(-1)^\vee\to \Lambda$ de $F(-1)^{\vee\vee}=F(-1)$.

\subsection{Trois exemples détaillés}\label{subsec:exdetlis}

\subsubsection{Un ouvert de $\PP^1$}

Posons $n=2$. Supposons que $-1$ n'est pas un carré dans le corps parfait $k_0$. Notons $V=\PP^1-\{ 0,\pm 1,\infty \}$ et $U=\PP^1-\{ 0,1,\infty\}$. Considérons le revêtement étale de degré 2
\[f\colon \begin{array}[t]{rcl} V&\longrightarrow& U \\ y&\longmapsto& y^2\end{array}\]
de groupe d'automorphismes engendré par $\tau\colon y\mapsto -y$. Le faisceau $\F\coloneqq f_\star\Lambda$ est un faisceau lisse sur $U$, trivialisé par le revêtement $f\colon V\to U$ puisque $f^\star f_\star \Lambda\simeq \Lambda^2$. Il correspond au $\Aut(V|U)$-module $\Lambda^2$, où l'élément non trivial de $\Aut(V|U)$ intervertit les deux copies de $\Lambda$. Comme $f$ est fini, $\R f_\star \Lambda=(f_\star\Lambda)[0]$ et il y a un isomorphisme canonique \[\RG(U,f_\star \Lambda)=\RG(V,\Lambda).\] Nous devrions donc trouver\[\HH^1(U,\F)=\HH^1(V,\Lambda)\simeq\Lambda^3.\] \\

Nous avons déjà calculé dans l'exemple de la section \ref{subsec:exdetrev} le revêtement $V_2\to V$ de groupe $\HH^1(V,\Lambda)^\vee$. C'est la normalisation de $V$ dans 
\[ \begin{array}{rcl} \Proj k[y,z,t,h]/(z^2-(y^2-h^2),t^2-(y^2+h^2)) &\longrightarrow& \Proj k[y,h]=\bar V \\ (y:z:t:h)&\longmapsto& (y^2:h^2).\end{array}\]
Le groupe $G=\Aut(V_2|U)$ est d'ordre 16, engendré par $\gamma(y,z,t,h)=(\sqrt{-1}y,\sqrt{-1}t,\sqrt{-1}z,h)$ et trois éléments $\sigma_1,\sigma_2,\sigma_3$ d'ordre 2. 

\paragraph{Calcul de $\RG(U,\F)$} Un complexe à deux termes représentant $\RG(U,\F)$ est 
\[ \Lambda^2 \to \{ f\colon G\to \Lambda^2\mid \forall g,h\in G,f(gh)=f(g)+gf(h)\}. \]
Un morphisme croisé $f\colon G\to \Lambda^2$ est déterminé par les images de $\sigma_1,\sigma_2,\sigma_3,\gamma$. En exploitant les relations $\gamma\sigma_1=\sigma_1\gamma, \gamma\sigma_2=\sigma_3\gamma$ et $\gamma^2=\sigma_1\sigma_2\sigma_3$, il vient qu'un tel morphisme croisé est uniquement déterminé par un quadruplet $(a,a_1,a_2,a_3)\in \Lambda^4$ ; le morphisme croisé correspondant à ce quadruplet vérifie $f(\sigma_1)=(a_1,a_1)$, $f(\sigma_2)=(a_2,a_3)$, $f(\sigma_3)=(a_3,a_2)$ et $f(\gamma)=(a,a+a_1+a_2+a_3)$. Parmi ceux-ci, les morphismes croisés principaux sont les images de $(0,0,0,0)$ et $(1,0,0,0)$. Le complexe calculé est donc isomorphe à
\[ \begin{array}{rcl}\Lambda^2 &\longrightarrow&\Lambda^4 \\ (a,b) & \longmapsto &(a+b,0,0,0) \end{array} \] et ses groupes de cohomologie sont $\Lambda$ et $\Lambda^3$, ce qui est le résultat attendu. 

\paragraph{Action de Galois} L'action de $\mathfrak{G}_0=\Gal(k|k_0)$ sur $\Aut(Y|X)$ se factorise manifestement par le quotient $\Gal(k_0(\sqrt{-1})|k_0)$. Ce groupe est engendré par la conjugaison $\sqrt{-1}\mapsto -\sqrt{-1}$. L'automorphisme $\phi$ agit trivialement sur $\sigma_1,\sigma_2,\sigma_3$, et $\phi\cdot \gamma=\sigma_1\sigma_2\sigma_3\gamma\colon (y,z,t)\mapsto (-\sqrt{-1}y,-\sqrt{-1}z,-\sqrt{-1}t)$. L'action de $\phi$ sur $\Lambda^2$ est triviale, et son action sur le groupe $\Lambda^4$ des morphismes croisés est $(\phi\cdot f)(x)=\phi f(\phi^{-1}x)=f(\phi^{-1}x)$. Explicitement, comme $\phi\cdot\gamma=\sigma_1\sigma_2\sigma_3\gamma$, cette action est \[\phi\cdot (a,a_1,a_2,a_3)=(a+a_1+a_2+a_3,a_1,a_2,a_3).\]

\subsubsection{Une courbe elliptique}

Cet exemple illustre la possibilité que les groupes de cohomologie d'un faisceau lisse dont la fibre générique géométrique est un $\Lambda$-module libre ne soient pas des $\Lambda$-modules libres. Soit $E$ une courbe elliptique sur un corps algébriquement clos de caractéristique impaire. Posons $n=4$. L'isogénie de multiplication par 2 sur $E$ est un revêtement galoisien de groupe $(\ZZ/2\ZZ)^2$ ; notons $E'\to E$ ce revêtement. Considérons le faisceau lisse $\F$ de $\Lambda=\ZZ/4\ZZ$-modules de fibre $M=\Lambda^2$ sur $E$ trivialisé par $E'\to E$, correspondant à la représentation :
\[ \begin{array}{rcl}
(\ZZ/2\ZZ)^2&\longrightarrow& \GL_2(\Lambda) \\
(1,0),(0,1) &\longmapsto& \begin{pmatrix}
1 & 2 \\
0 & 1
\end{pmatrix}
\end{array} \]
Le groupe $\HH^0(E,\F)=\HH^0((\ZZ/2\ZZ)^2,M)$ est $2\Lambda\times \Lambda\simeq (\ZZ/2\ZZ)\times(\ZZ/4\ZZ)$.
Le revêtement $E'_2\to E'$ de groupe $\HH^1(E',\Lambda)^\vee$ est la multiplication par 4 sur $E'$ ; le revêtement composé $E'_2\to E$ est la multiplication par 8 sur $E$, de groupe d'automorphismes $(\ZZ/8\ZZ)^2$. Le groupe $\HH^1(E,\F)$ est donc isomorphe à $\HH^1((\ZZ/8\ZZ)^2,M)$, où l'action de $(\ZZ/8\ZZ)^2$ sur $M$ est définie par le morphisme composé $(\ZZ/8\ZZ)^2\to (\ZZ/2\ZZ)^2\to \GL_2(\Lambda)$. Ce groupe est isomorphe au quotient du $\Lambda$-module libre $\Lambda^4$ par le sous-module engendré par $(2,0,0,0)$. Plus précisément, il est engendré par les classes des cocycles $c_1,c_2,c_3,c_4$, où $2c_1$ est de classe nulle, définis comme suit.
\begin{center}
\begin{tabular}{|c|c|c|c|c|c|c|}\hline
Cocycle $c$ & $c(1,0)$ & $c(0,1)$ \\ \hline
$c_1$ & $(1~0)$  & $(1~0)$ \\ \hline
$c_2$ & $(0~0)$ & $(1~0)$ \\ \hline
$c_3$ & $(0~0)$ & $(0~1)$ \\ \hline
$c_4$ & $(0~1)$ & $(0~0)$ \\ \hline
\end{tabular}
\end{center}
Le cocycle $2c_1$ associe à $(0,1)$ et $(1,0)$ le vecteur $(2,0)$ ; c'est $\partial c'$, où $c'$ est le 0-cocycle de valeur $(0,1)$.

\subsubsection{Un ouvert d'une courbe elliptique}

Ce troisième exemple a nécessité l'utilisation d'un logiciel de calcul formel, à la fois pour les opérations dans la jacobienne d'une courbe de genre 2 et pour la cohomologie des groupes. Les algorithmes du chapitre \ref{chap:4} n'ayant pas été implémentés, nous nous servons des spécificités des courbes hyperelliptiques pour la détermination de la 2-torsion de la jacobienne, et suivons pas à pas notre algorithme précédemment décrit pour tout le reste. Les fonctions pour calculer dans la jacobienne d'une courbe hyperelliptique sont disponibles dans \textsc{SageMath} et \textsc{Magma}.\\

\paragraph{Le revêtement considéré} Dans ce deuxième exemple, nous prenons toujours $n=2$. Le corps $k_0$ considéré est $\FF_{11}$, et tous les calculs sont effectués sur $\FF_{121}=\FF_{11}(a)$, où $a$ est un générateur de $\FF_{121}^\times$ vérifiant $a^2+7a+2=0$. Soit $\bar E$ la courbe elliptique sur $k=\overline{\FF_{121}}$ définie par l'équation de Weierstrass affine $y^2=(x-1)(x-2)(x-2)$. 
Soit $\bar C$ la courbe projective lisse de genre 2 sur $k$ d'équation affine $y^2=(x^2-1)(x^2-2)(x^2-3)$. La courbe $\bar C$ possède deux points $\infty_+,\infty_-$ qui ne sont pas situés sur l'ouvert affine donné par cette équation. Considérons le revêtement $\bar f\colon \bar C\to \bar E$ de degré 2 donné par $(x,y)\mapsto (x^2,y)$. Il est ramifié en les points affines $P=(0,4)$ et $Q=(0,7)$ de $\bar C$ d'images respectives $(0,4)$ et $(0,7)$. Notons $C=\bar C-\{ P,Q\}$ et $E=\bar E-\{ \bar f(P),\bar f(Q)\}$. Notons $f\colon C\to E$ le revêtement étale galoisien obtenu. Remarquons que les courbes $C$ et $E$ proviennent de courbes $C_0$ et $E_0$ sur $k_0=\FF_{11}$.

\paragraph{Calcul de $\HH^1(C,\mu_2)$} Notons
\[P_1^\pm=(\pm 1,0), \qquad P_2^\pm=(\pm a^6,0),\qquad P_3^\pm =(\pm 5,0)\]
les points de $C$ d'ordonnée nulle. Une base de la $2$-torsion de la jacobienne $J_{\bar C}$ est donnée par les classes des diviseurs \[ D_1\coloneqq P_1^+-P_1^-,\qquad D_2\coloneqq P_2^+-P_2^-,\qquad D_3\coloneqq P_2^+-P_3^+,\qquad D_4\coloneqq P_1^+-P_3^-.\]
Les fonctions \[ f_1\coloneqq \frac{x-1}{x+1},\qquad f_2\coloneqq \frac{x-a^6}{x+a^6},\qquad f_3\coloneqq\frac{x-a^6}{x-5},\qquad f_4\coloneqq\frac{x-1}{x+5}\]
vérifient $2D_i=\div(f_i)$. Notons $D_5$ le diviseur $P-Q$. Le double du diviseur \[ \bar D_5\coloneqq (a^{41},a^{29})+(-a^{41},a^{29})-(\infty_++\infty_-)\]
est équivalent à $D_5$. En l'absence d'algorithme implémenté à cet effet, nous l'avons obtenu en parcourant les points de $J_C(\FF_{121})$. La fonction
 \[ f_5=\frac{1}{xy}+\frac{a^8x^2+7}{x}\]
vérifie $\div(f_5)=2\bar D_5-D_5$. Une $\FF_2$-base de $\HH^1(C,\mu_2)$ est donc donnée par les triplets \[([D_1],0,f_1),\dots,([D_4],0,f_4),([\bar D_5],D_5,f_5).\]

\paragraph{Le revêtement $C_2\to E$} Le revêtement $C_2\to C$ de groupe $\HH^1(C,\Lambda)^\vee$ a pour corps de fonctions $k(C)(z_1,\dots,z_5)$ où $z_i^2=f_i$. 
Le groupe $G=\Aut(C_2|E)$ est d'ordre 64 ; il contient le sous-groupe distingué $H=\Aut(C_2|C)\simeq (\ZZ/2\ZZ)^5$ engendré par les $\gamma_i\colon z_i\mapsto -z_i$. Voici comment déterminer un élément de $G$ d'image le générateur $\sigma\colon (x,y)\mapsto (-x,y)$ de  $\Aut(C|E)$. Le calcul des diviseurs $\sigma^\star D_i$ fournit le résultat suivant.
\[
\begin{array}{rllll}
\sigma^\star D_1 & = & -D_1 & &\\
\sigma^\star D_2 & = & -D_2 & &\\
\sigma^\star D_3 & = & D_1+D_3+\div(h_3) & \text{où} & h_3=\frac{y}{x^3+a^{58}x^2+a^2x+a^{54}}\\
\sigma^\star D_4 & = & D_2+D_4+\div(h_4) & \text{où} & h_4=\frac{y}{x^3+a^{80}x^2+a^{103}x+a^{114}}\\
\sigma^\star D_5 & = & D_5
\end{array}
\]
Remarquons que $\sigma^\star f_3=h_3^2f_1f_3$, $\sigma^\star f_4=h_4^2f_2f_4$ et $\sigma^\star f_5=-f_5$. Comme $a^{30}$ est une racine carrée de $-1$ dans $\FF_{121}$, un antécédent de $\sigma$ est l'automorphisme $\delta$ défini par \[ 
x\mapsto -x, \quad y\mapsto y,\quad z_1\mapsto \frac{1}{z_1},\quad z_2\mapsto \frac{1}{z_2},\quad z_3\mapsto h_3(x,y)z_1z_3,\quad z_4\mapsto h_4(x,y)z_2z_4,\quad z_5\mapsto a^{30}z_5.\]
En particulier, comme $h_3(x,y)h_3(-x,y)=h_4(x,y)h_4(-x,y)=-1$, l'automorphisme $\delta^2$ est donné par \[ 
x\mapsto x, \quad y\mapsto y,\quad z_1\mapsto z_1,\quad z_2\mapsto z_2,\quad z_3\mapsto -z_3,\quad z_4\mapsto -z_4,\quad z_5\mapsto -z_5.\] Par conséquent, $\delta^2=\gamma_3\gamma_4\gamma_5$ et $\delta$ est d'ordre 4 dans $G$. De plus, $\delta\gamma_1=\gamma_1\gamma_3\delta$ et $\delta\gamma_2=\gamma_2\gamma_4\delta$. 
Les éléments $\gamma_3,\gamma_4,\gamma_5$ engendrent le centre de $G$. Son groupe dérivé est $\langle \gamma_3=[\gamma_1,\delta],\gamma_4=[\gamma_2,\delta]\rangle$. 
Une fois ce groupe calculé, il est possible par de simples méthodes d'algèbre linéaire implémentées dans tous les logiciels de calcul formel courants de calculer la cohomologie de n'importe quel faisceau lisse de $\Lambda$-modules sur $E$ trivialisé par $C$. Il est également possible de calculer l'action de $\Gal(\FF_{121}|\FF_{11})=\colon\langle\phi\rangle$ sur $G$. L'action de $\phi$ sur les éléments de $G$ se calcule simplement sur les équations qui les définissent : il fixe les $\gamma_i$ et envoie $\delta$ sur $\gamma_5\delta$.
 \\
\paragraph{Un premier calcul de cohomologie} Soit $\F_0$ le faisceau lisse de $\Lambda=\ZZ/2\ZZ$-modules sur $E_0$ trivialisé par $C_0$ de fibre générique géométrique un $\Lambda$-module libre $M$ de dimension 2, défini par la représentation : \[ \begin{array}{rcl}
\Aut(C_0|E_0)&\longrightarrow& \GL_2(\Lambda) \\
\sigma&\longmapsto& \begin{pmatrix}
0 & 1 \\
1 & 0
\end{pmatrix}
\end{array} \]
Déterminons la cohomologie du faisceau $\F=(\F_0)|_E=f_\star\Lambda$. Le complexe \[ M\to C^{12}(G,M)\] représentant $\tau_{\leqslant 1}\RG(G,M)$ est de la forme 
\[ \Lambda^2\to \Lambda^6.\]
Le groupe $C^{12}(G,M)$ est engendré par des 1-cocycles $c_1,\dots,c_5,c'\colon G\to M$. Ils sont donnés dans le tableau suivant.

\begin{center}
\begin{tabular}{|c|c|c|c|c|c|c|}\hline
Cocycle $c$ & $c(\gamma_1)$ & $c(\gamma_2)$ & $c(\gamma_3)$ & $c(\gamma_4)$ & $c(\gamma_5)$ & $c(\delta)$ \\ \hline
$c_1$ & $(1~0)$  & $(0~0)$ & $(1~1)$  & $(0~0)$ & $(0~0)$ & $(0~1)$\\ \hline
$c_2$ & $(0~1)$ & $(0~0)$ & $(1~1)$ & $(0~0)$ & $(0~0)$ & $(0~1)$\\ \hline
$c_3$ & $(0~0)$ & $(1~0)$ & $(0~0)$  & $(1~1)$ & $(0~0)$ & $(0~1)$\\ \hline
$c_4$ & $(0~0)$ & $(0~1)$ & $(0~0)$ & $(1~1)$ & $(0~0)$ & $(0~1)$\\ \hline
$c_5$ & $(0~0)$ & $(0~0)$ & $(0~1)$ & $(0~0)$ & $(1~1)$ & $(0~1)$\\ \hline
$c'$ &  $(0~0)$ & $(0~0)$ & $(0~0)$ & $(0~0)$ & $(0~0)$ & $(1~1)$ \\\hline
\end{tabular}
\end{center}

Ici, $c'$ est l'image des 0-cocycles $\alpha=(1~0)$ et $\beta=(0~1)$ par l'application $\partial\colon M\to C^1(G,M)$.
En particulier, nous retrouvons le résultat attendu $\HH^1(E,f_\star \Lambda)=\HH^1(C,\Lambda)\simeq \Lambda^5$. L'action du groupe $\Gal(\FF_{121}|\FF_{11})=\langle \phi\rangle$ sur ce complexe fixe $c_1,\dots,c_4,c'$ et envoie $c_5$ sur $\phi^\star c_5=c_5+c'$. Par conséquent, comme attendu, son action sur $\HH^1(E,f_\star\Lambda)$ est triviale.

\paragraph{Un second calcul de cohomologie} Considérons désormais le faisceau lisse $\F_0$ sur $E_0$ trivialisé par $C_0$ de fibre $M=\Lambda^3$ défini par la représentation :
\[ \begin{array}{rcl}
\Aut(C_0|E_0)&\longrightarrow& \GL_3(\Lambda) \\
\sigma&\longmapsto& \begin{pmatrix}
0 & 0  & 1\\
0 & 1 & 0 \\
1 & 0 & 0
\end{pmatrix}
\end{array} \]
Nous trouvons $\HH^0(E,\F)\simeq\Lambda^2$ et $\HH^1(E,\F)\simeq\Lambda^8$. Un calcul fournit des 1-cocycles $c_1,\dots,c_8$ formant une base de $\HH^1(G,M)$ et le 1-cocycle principal $c'$ défini par $(1~0~0)\in M$.
\begin{center}
\begin{tabular}{|c|c|c|c|c|c|c|}\hline
Cocycle $c$ & $c(\gamma_1)$ & $c(\gamma_2)$ & $c(\gamma_3)$ & $c(\gamma_4)$ & $c(\gamma_5)$ & $c(\delta)$ \\ \hline
$c_1$ & $(1~0~0)$ & $(0~0~0)$ & $(1~0~1)$  & $(0~0~0)$ & $(0~0~0)$ & $(0~0~1)$\\ \hline
$c_2$ & $(0~1~0)$ & $(0~0~0)$ & $(0~0~0)$ & $(0~0~0)$ & $(0~0~0)$ & $(0~0~0)$\\ \hline
$c_3$ & $(0~0~1)$ & $(0~0~0)$ & $(1~0~1)$  & $(0~0~0)$ & $(0~0~0)$ & $(0~0~1)$\\ \hline
$c_4$ & $(0~0~0)$ & $(1~0~0)$ & $(0~0~0)$ & $(1~0~1)$ & $(0~0~0)$ & $(0~0~1)$\\ \hline
$c_5$ & $(0~0~0)$ & $(0~1~0)$ & $(0~0~0)$ & $(0~0~0)$ & $(0~0~0)$ & $(0~0~0)$\\ \hline
$c_6$ & $(0~0~0)$ & $(0~0~1)$ & $(0~0~0)$ & $(1~0~1)$ & $(0~0~0)$ & $(0~0~1)$\\ \hline
$c_7$ & $(0~0~0)$ & $(0~0~0)$ & $(0~0~0)$ & $(0~0~0)$ & $(1~0~1)$ & $(0~0~1)$\\ \hline
$c_8$ & $(0~0~0)$ & $(0~0~0)$ & $(0~0~0)$ & $(0~0~0)$ & $(0~0~0)$ & $(0~1~0)$\\ \hline
$c'$ &  $(0~0~0)$ & $(0~0~0)$ & $(0~0~0)$ & $(0~0~0)$ & $(0~0~0)$ & $(1~0~1)$ \\\hline
\end{tabular}
\end{center}
L'action de $\phi\in\Gal(\FF_{121}|\FF_{11})$ sur $C^{12}(G,M)$ fixe tous ces cocycles sauf $c_7$, qui vérifie $\phi^\star c_7=c_7+c'$. L'action sur $\HH^1(G,M)$ est donc triviale.

\subsection{Cohomologie des complexes de faisceaux lisses}\label{subsec:cohcomplisses}

Soit $X$ une courbe intègre lisse sur $k$. Supposons d'abord $X$ affine. Soit \[\K=[\K^0\to \cdots \to \K^m] \] un complexe borné de faisceaux lisses sur $X$ de fibre générique géométrique un complexe $K$ de $\pi_1(X)$-modules. Nous allons montrer comment utiliser la méthode décrite en \ref{subsec:cohlisG} pour déterminer un complexe de $\Lambda$-modules représentant $\RG(X,\K)$.
Soit $f\colon Y\to X$ un revêtement galoisien de $X$ qui trivialise tous les $\mathscr{K}^i$. Notons $\pi=\pi_1(X)$.
Pour tout groupe profini $H$, notons $P_H(\Lambda)$ la résolution projective usuelle (résolution bar) de $\Lambda$ comme $\Lambda[[H]]$-module. Nous souhaitons calculer $\RG(\pi,K)=\RHom_{\Lambda[[\pi]]}(\Lambda,K)$, qui est représenté par le complexe \[ \Hom^\bullet_{\Lambda[[\pi]]}(P_\pi(\Lambda),K)=\Tot(A^{\bullet,\bullet}) \] où $A^{p,q}=\Hom_{\Lambda[[\pi]]}(P_\pi^{-q}(\Lambda),K^p)$.
Soit $Y_2\to Y$ le revêtement de $Y$ de groupe $\HH^1(Y,\Lambda)^\vee$. Notons de plus $G=\Aut(Y_2|X)$.
Considérons le complexe double $B^{p,q}=\Hom_{\Lambda[G]}(\tau_{\geqslant -1}P_{G}^{-q}(\Lambda),K^p)$. 

\begin{prop} Le complexe $\Tot B^{\bullet,\bullet}$ représente $\RG(X,\K)$.
\begin{proof}
Le morphisme $B^{\bullet,\bullet}\to A^{\bullet,\bullet}$ induit par le quotient $\pi\to G$ définit encore un morphisme entre les suites spectrales associées à ces quotients.
Rappelons que pour tout faisceau lisse $\F$ sur $X$ trivialisé par $Y$ de fibre $F$, le morphisme $\tau_{\leqslant 1}\Hom_{\Lambda[G]}(P_G(\Lambda),F)\to 
\Hom_{\Lambda[[\pi]]}(P_{\pi}(\Lambda),F)$ est un quasi-isomorphisme. 
Remarquons que $\Hom_{\Lambda[G]}(\tau_{\geqslant -1}P_G(\Lambda),F)=\tau_{\leqslant 1}\Hom_{\Lambda[H]}(P_G(\Lambda),F)$ par exactitude à gauche de $\Hom(-,M)$.
Le morphisme $B^{\bullet,\bullet}\to A^{\bullet,\bullet}$ est donc, sur chaque colonne, un quasi-isomorphisme de complexes. Par conséquent, il induit un isomorphisme entre les premières pages des suites spectrales (pour l'orientation verticale) associées à $B^{\bullet,\bullet}$ et $A^{\bullet,\bullet}$ : il s'agit en position $(p,q)$ de l'isomorphisme $\HH^q(H,K^p)\to \HH^q(\pi,K^p)$ si $q\leqslant 1$, et $0\to \HH^q(\pi,K^p)=0$ sinon. Par conséquent, le morphisme entre les aboutissements de ces deux suites spectrales est un isomorphisme, c'est-à-dire que le morphisme $\Tot(B^{\bullet,\bullet})\to \Tot(A^{\bullet,\bullet})$ est un quasi-isomorphisme.
\end{proof}
\end{prop}

\section{Cup-produits dans la cohomologie des faisceaux lisses}\label{sec:cupprod}

Nous décrivons dans cette section comment calculer des cup-produits dans la cohomologie de faisceaux lisses sur des courbes lisses sur un corps fini ou algébriquement clos. Dans un premier temps, nous décrivons un calcul du $\HH^2$ d'un faisceau lisse sur une courbe projective lisse comme groupe de cohomologie d'un quotient du groupe fondamental de la courbe.
Nous indiquons ensuite comment calculer le cup-produit de la dualité de Poincaré dans un cas simple. Nous nous servons enfin du calcul du $\HH^2$ exposé au début de la section pour décrire les cup-produits $\HH^1\times \HH^1\to \HH^2$ entre groupes de cohomologie de faisceaux lisses.

\subsection{Un autre calcul du $\HH^2$}

Soit $X$ une courbe intègre propre lisse de genre non nul sur $k$. Soit $\F$ un faisceau lisse de $\Lambda$-modules sur $k$ de fibre $M$, trivialisé par un revêtement galoisien $Y\to X$. Nous allons montrer comment calculer $\HH^2(X,\F)$ comme le $\HH^2$ d'un quotient de $\pi_1(X)$ à valeurs dans $M$.
Supposons (quitte à passer à un revêtement caractéristique de $Y$) que le degré de $Y\to X$ est un multiple de $n$.
Soient $Y_2\to Y$ le revêtement caractéristique de groupe $\HH^1(Y,\Lambda)^\vee$ et $Y_3\to Y_2$ le revêtement caractéristique de groupe $\HH^1(Y_2,\Lambda)^\vee$. Notons $\pi,\pi_Y,\pi_Y^{[2]},\pi_Y^{[3]}$ les groupes fondamentaux respectifs de $X,Y,Y_2,Y_3$. Notons $G=\pi/\pi_Y=\Aut(Y|X)$, $G_2=\pi/\pi_Y^{[2]}=\Aut(Y_2|X)$ et $G_3=\pi/\pi_Y^{[3]}=\Aut(Y_3|X)$. 

\begin{prop}\label{prop:H2im} La restriction de $\HH^2(G_3,M)\to \HH^2(\pi,M)$ à l'image de $\HH^2(G_2,M)$ est un isomorphisme.
\begin{proof}
Rappelons que comme les morphismes \[ \HH^1(X,\F)\to \HH^1(Y_2,\F) \to \HH^1(Y_3,\F)\] sont nuls, les morphismes \[ \HH^1(G_2,M)\to \HH^1(G_3,M)\to \HH^1(\pi,M)=\HH^1(X,\F) \] sont des isomorphismes.
Considérons les morphismes entre les suites spectrales de Hochschild-Serre associés aux morphismes d'extensions de groupes :
\[ 
\begin{tikzcd}
1 \arrow[r] &\pi_Y^{[3]} \arrow[r] \arrow[d] &\pi \arrow[r]\arrow[d] & G_3 \arrow[r]\arrow[d] & 1 \\
1 \arrow[r] &\pi_Y^{[2]} \arrow[r]\arrow[d] & \pi \arrow[r]\arrow[d] & G_2 \arrow[r]\arrow[d] & 1 \\
1 \arrow[r] &\pi_Y \arrow[r] & \pi \arrow[r] & G \arrow[r]& 1 \\
\end{tikzcd}
\]
Ils donnent les morphismes de suites exactes de bas degré :
\[
\begin{adjustbox}{width=\textwidth}{
\begin{tikzcd}0 \arrow[r] & \HH^1(G,M)\arrow[r] \arrow[d] & \HH^1(\pi,M) \arrow[r] \arrow[d] & \HH^0(G,\HH^1(\pi_Y,M))\arrow[r] \arrow[d] & \HH^2(G,M)\arrow[r] \arrow[d] & \ker(\HH^2(\pi,M)\to \HH^2(\pi_Y,M))\arrow[r] \arrow[d] & \HH^1(G,\HH^1(\pi_Y,M))\arrow[d] \\
0 \arrow[r] & \HH^1(G_2,M)\arrow[r,"\sim"] \arrow[d] & \HH^1(\pi,M) \arrow[r] \arrow[d] & \HH^0(G_2,\HH^1(\pi_Y^{[2]},M))\arrow[r] \arrow[d] & \HH^2(G_2,M)\arrow[r] \arrow[d] & \ker(\HH^2(\pi,M)\to \HH^2(\pi_Y^{[2]},M))\arrow[r] \arrow[d] & \HH^1(G_2,\HH^1(\pi_Y^{[2]},M))\arrow[d] \\
0 \arrow[r] & \HH^1(G_2,M)\arrow[r,"\sim"] & \HH^1(\pi,M) \arrow[r] & \HH^0(G_3,\HH^1(\pi_Y^{[3]},M))\arrow[r] &  \HH^2(G_3,M)\arrow[r]  & \ker(\HH^2(\pi,M)\to \HH^2(\pi_Y^{[3]},M))\arrow[r]& \HH^1(G_3,\HH^1(\pi_Y^{[3]},M))
\end{tikzcd}
}
\end{adjustbox}
\]
Comme les degrés des revêtements $Y\to X$ et $Y_2\to X$ sont des multiples de $n$, les morphismes \[ \HH^2(\pi,M)\to \HH^2(\pi_Y,M) \text{~~et~~} \HH^2(\pi,M)\to \HH^2(\pi_Y^{[2]},M) \] sont nuls. De plus, comme $\F$ est trivial sur $Y$, le morphisme $\HH^1(\pi_Y,M)\to \HH^1(\pi_Y^{[2]},M)$ est nul. Le carré supérieur droit du diagramme ci-dessus est donc
\[ \begin{tikzcd}
\ker(\HH^2(\pi,M)\to \HH^2(\pi_Y,M)) \arrow[r]\arrow[d,"\sim"{anchor=north, rotate=90}] & \HH^1(G,\HH^1(\pi_Y,M)) \arrow[d,"0"] \\
\HH^2(\pi,M) \arrow[r] & \HH^1(G_2,\HH^1(\pi_Y^{[2]},M))
\end{tikzcd}
\]
ce qui démontre la nullité du morphisme $\HH^2(\pi,M)\to \HH^1(G_2,\HH^1(\pi_Y^{[2]},M))$. Le même raisonnement peut être appliqué au morphisme $Y_3\to Y_2$ pour obtenir le diagramme commutatif suivant.
\[ 
\begin{tikzcd}
0 \arrow[r] & \HH^0(G_2,\HH^1(\pi_Y^{[2]},M)) \arrow[r] \arrow[d,"0"] & \HH^2(G_2,M) \arrow[r] \arrow[d] & \HH^2(\pi,M) \arrow[r]\arrow[d,"\sim"{anchor=north, rotate=90}] & 0\\
0 \arrow[r] & \HH^0(G_3,\HH^1(\pi_Y^{[3]},M)) \arrow[r] & \HH^2(G_3,M) \arrow[r] & \HH^2(\pi,M) \arrow[r] & 0
\end{tikzcd} 
\]
Le morphisme $\HH^2(G_2,M)\to \HH^2(G_3,M)$ passe au quotient en un morphisme \[ \HH^2(G_2,M)/\HH^0(G_2,\HH^1(\pi_Y^{[2]},M)) \to \HH^2(G_3,M)\]
et le cardinal de son image est au plus $|\HH^2(\pi,M)|$. Par conséquent, la restriction du morphisme $\HH^2(G_3,M)\to \HH^2(\pi,M)$ à l'image de $\HH^2(G_2,M)\to \HH^2(G_3,M)$, que l'on savait déjà surjective, est un isomorphisme.
\end{proof}
\end{prop}

\begin{rk}
Le même raisonnement s'applique mot pour mot au cas des courbes projectives lisses géométriquement connexes sur les corps finis. La description du revêtement $Y_2$ se trouve dans la remarque \ref{rk:X2fin}. 
\end{rk}

\begin{rk}
L'hypothèse de divisibilité par $n$ du degré de $Y\to X$ est nécessaire pour notre méthode de démonstration, mais ne semble pas nécessaire en général pour obtenir ce résultat. Par exemple, soit $E$ une courbe elliptique sur $k$. Alors $G_2=\Aut(E_2|E)\simeq \Lambda^2$ et $G_3=(\ZZ/n^2\ZZ)^2$. Le groupe $\HH^2(G_2,\Lambda)$ est un $\Lambda$-module libre de rang 3 par la formule de Künneth. Or il contient $\Ext^1(G_2,\Lambda)$ qui est de rang 2. Il y a donc une extension centrale de $G_2$ par $\Lambda$ qui n'est pas abélienne, et qui fournit par produit fibré une extension centrale non abélienne de $G_3$ par $\Lambda$, c'est-à-dire un élément non nul de $\HH^2(G_3,\Lambda)$. De plus, le morphisme $\Ext^1(G_2,\Lambda)\to \Ext^1(G_3,\Lambda)$ est nul : l'image du morphisme $\HH^2(G_2,\Lambda)\to \HH^2(G_3,\Lambda)$ est de rang 1, c'est donc $\HH^2(E,\Lambda)$.
\end{rk}

\subsection{La dualité de Poincaré dans un cas particulier}

Soit $X$ une courbe projective lisse sur $k$. Soit $\F$ un faisceau lisse de $\Lambda$-modules sur $X$, trivialisé par un revêtement galoisien $f\colon Y\to X$ de groupe $G$. Nous souhaitons calculer le cup-produit
\[ \HH^1(X,\F)\times \HH^1(X,\F^\vee)\to \HH^2(X,\Lambda).\]
Dans le cas particulier où $n$ est premier à l'ordre de $G$, le morphisme $\HH^2(X,\Lambda)\to \HH^2(Y,\Lambda)$ est un isomorphisme. Le diagramme commutatif suivant
\[
\begin{tikzcd}
\HH^1(Y,f^\star\F)\times \HH^1(Y,f^\star\F^\vee) \arrow[r] & \HH^2(Y,\Lambda) \\
\HH^1(X,\F)\times \HH^1(X,\F^\vee) \arrow[r]\arrow[u] & \HH^2(X,\Lambda)\arrow[u]
\end{tikzcd}
\]
permet alors de se ramener au cas particulier du cup-produit dans la cohomologie des faisceaux constants, qui est l'accouplement de Weil décrit dans la section \ref{subsubsec:courbcohsupp}.

\subsection{Calcul des cup-produits en général}

Soit $X$ une courbe projective lisse sur le corps algébriquement clos $k$.
Soient $\F,\G$ deux faisceaux lisses de $\Lambda$-modules sur $X$ de fibres respectives $M$ et $N$. Nous souhaitons calculer le cup-produit 
\[ \HH^1(X,\F)\times \HH^1(X,\G)\to \HH^2(X,\F\otimes \G).\]
Soit $f\colon Y\to X$ un revêtement galoisien qui trivialise à la fois $\F$ et $\G$. Supposons que le degré de $Y\to X$ soit divisible par $n$. Soient $Y_2\to Y$ le revêtement caractéristique de groupe $\HH^1(Y,\Lambda)^\vee$ et $Y_3\to Y_2$ le revêtement caractéristique de groupe $\HH^1(Y_2,\Lambda)^\vee$. Notons $G_2=\Aut(Y_2|X)$ et $G_3=\Aut(Y_3|X)$.\\

\begin{theorem}\label{th:cupprod} Le cup-produit 
\[ \HH^1(X,\F)\times \HH^1(X,\G) \to \HH^2(X,\F\otimes \G) \]
est réalisée par la composée
\[ \HH^1(G_2,M)\times \HH^1(G_2,N) \xrightarrow{\cup} \HH^2(G_2,M\otimes N) \rightarrow \im(\HH^2(G_2,M\otimes N)\to \HH^2(G_3,M\otimes N)).\]
\begin{proof} Rappelons que les morphismes
\[ \im (\HH^2(G_2,M\otimes N) \to \HH^2(G_3,M\otimes N))\to \HH^2(\pi,M\otimes N)\to \HH^2(X,\F\otimes \G) \]\\
sont des isomorphismes d'après la proposition \ref{prop:H2im}. Le diagramme commutatif
\[ 
\begin{tikzcd}
\HH^1(G_2,M)\times \HH^1(G_2,N) \arrow[d,"\sim"{anchor=south, rotate=90}] \arrow[r] & \HH^2(G_2,M\otimes N) \arrow[d] \\
\HH^1(G_3,M)\times \HH^1(G_3,N) \arrow[d,"\sim"{anchor=south, rotate=90}] \arrow[r] & \HH^2(G_3,M\otimes N) \arrow[d] \\
\HH^1(\pi,M) \times \HH^1(\pi,N) \arrow[d,"\sim"{anchor=south, rotate=90}] \arrow[r] & \HH^2(\pi,M\otimes N)\arrow[d,"\sim"{anchor=south, rotate=90}] \\
\HH^1(X,\F)\times \HH^1(X,\G) \arrow[r] & \HH^2(X,\F\otimes \G)
\end{tikzcd}
\]
où les trois premières lignes sont les cup-produits de cohomologie des groupes,
montre comment calculer le cup-produit de la dernière ligne : c'est la composée
\[ \HH^1(G_2,M)\times \HH^1(G_2,N) \xrightarrow{\cup} \HH^2(G_2,M\otimes N) \rightarrow \im(\HH^2(G_2,M\otimes N)\to \HH^2(G_3,M\otimes N)).\]
\end{proof}
\end{theorem}

\begin{rk}
Ceci s'applique de la même façon au cas des courbes projectives lisses géométriquement connexes sur un corps fini.
\end{rk}

\section{Cohomologie des faisceaux constructibles : calcul du \texorpdfstring{$\HH^1$}{H\textsuperscript{1}}}

Cette section est dédiée aux méthodes de calcul direct du premier groupe de cohomologie d'un faisceau constructible sur une courbe lisse, en employant des méthodes déjà proposées dans la littérature. 

\subsection{Avec l'algorithme de Jin}\label{sec:jinconstr}

Soit $X$ une courbe intègre lisse sur $k$.
Soit $\F$ un faisceau constructible sur $X$. Supposons $\F$ défini par recollement par rapport à un couple ouvert-fermé $(j\colon U\to X,i\colon Z\to X)$ par le triplet $(\mathscr{L},\F_Z,\phi)$. La stratégie suivante, proposée par Jin dans \cite[§9.4]{jinbi_jin}, permet de calculer $\HH^1(X,\F)$.
La suite exacte de faisceaux sur $X$ 
\[ 0\to j_!\mathscr{L} \to \F\to i_\star i^\star\F \to 0\]
donne d'une part un isomorphisme
\[ \HH^2_c(U,j^\star \F)\to \HH^2(X,\F) \]
et d'autre part la suite exacte
\[ 0\to \HH^0(X,j_!\mathscr{L})\to \HH^0(X,\F)\to \HH^0(Z,\F_Z)\xrightarrow{\delta} \HH^1(X,j_!\mathscr{L})\to \HH^1(X,\F)\to 0 .\]
Il suffit donc de savoir calculer explicitement la flèche $\delta$ afin de déterminer $\HH^1(X,\F)$.
D'après \cite[Lem. 5.3]{jinbi_jin}, $\HH^1(X,j_!\mathscr{L})$ est le groupe des classes d'isomorphismes de couples $(T,s)$ où $T$ est un $\mathscr{L}$-torseur sur $U$ et $s\in i^\star j_\star \mathscr{L}(Z)$. La flèche $\delta$ associe à $s\in \HH^0(Z,\F_Z)$ le $\mathscr{L}$-torseur trivial sur $U$, et la section $\phi(s)\in  \HH^0(Z,i^\star j_\star \mathscr{L})$. \\

Soit maintenant $g\colon V\to U$ un revêtement galoisien de groupe $G$ qui trivialise $\mathscr{L}$. Notons $F\coloneqq g^\star\mathscr{L}$. Soit $Y$ la normalisation de $X$ dans $V$, et $i'\colon W\to Y$ le changement de base de $i\colon Z\to X$ par rapport à $Y\to X$. Notons également $j'_\star\colon U'\to Y$ le changement de base de $j\colon U\to X$. L'algorithme de Jin calcule des représentants des éléments de $\HH^1(X,j_!\mathscr{L})$ sous forme de couples $(T',s')$ où $T'$ est un $F$-torseur $G$-équivariant sur $Y$, et $s'\in \HH^0(G,i'^*T(W))$ \cite[Lem. 5.10]{jinbi_jin}. Le couple $(T',s')$ correspondant à $(T,s)$ est $(j'_\star g^\star T,s)$. Comme $g^\star T$ est un $G$-torseur sur $V$, le morphisme $g^\star T\to V$ est fini étale. La proposition \ref{prop:jstarnorm} indique alors comment calculer le torseur $j'_\star g^\star T$ : c'est le lieu étale de la normalisation de $Y$ dans $g^\star T$.

\subsection{Par effacement}\label{sec:coheff}
 
Soit $X$ une courbe lisse sur $k$. Soit $\F$ un faisceau constructible de $\Lambda$-modules sur $X$. La méthode suivante, proposée en dimension quelconque dans \cite[§11.4]{mo}, permet de calculer $\HH^1(X,\F)$. Supposons $\F$ lisse sur un ouvert $U$, de complémentaire réduit $Z$. Soit $V\to U$ un revêtement galoisien trivialisant $\F|_U$. Construisons un monomorphisme $\F\to \G$ de faisceaux constructibles tel que $\HH^1(X,\F)\to \HH^1(X,\G)$ soit nul. Dans ce cas, le groupe $\HH^1(X,\F)$ est isomorphe au conoyau du morphisme $\HH^0(X,\F)\to \HH^0(X,\coker (\F\to \G))$.\\

Plongeons $\F$ dans $p_\star M$, où $M$ est un $\Lambda$-module de type fini et $p$ est un morphisme fini. Explicitement, prenons $p\colon Y'=X'\sqcup Z\to X$ où $X'$ est la normalisation de $X$ dans $V$. Considérons $(X')_2$, le revêtement de $X'$ de groupe $\HH^1(X',\Lambda)^\vee$ qui trivialise les $\Lambda$-torseurs sur $X'$, et $Y''=(X')_2\sqcup Z$. D'une part, la composée $\HH^1(X,\F)\to \HH^1(U,\F|_U)\to \HH^1(X'',\F|_{X'})$ est nulle ; d'autre part, $\HH^1(X,\F)\to \HH^1(Z,\F|_Z)$ est nulle car $\HH^1(Z,\F|_Z)=0$. Par conséquent, le morphisme $\HH^1(X,\F)\to \HH^1(Y'',\F)$ est nul. Notons $q\colon Y''\to X$. Alors comme $Y''\to Y'$ est surjectif, $\F \hookrightarrow p_\star M\hookrightarrow q_\star M$.\\

Les schémas $X'$ et $X''$ se calculent explicitement, de même que le poussé en avant d'un faisceau constant par un morphisme fini. Ainsi, il est possible d'obtenir une description par recollement du faisceau $q_\star M$, du morphisme $p_\star M\to q_\star M$ et finalement du conoyau de la composée $\F\to q_\star M$. Les sections globales d'un faisceau constructible décrit par recollement étant simplement un produit fibré de groupes abéliens, le morphisme $\HH^0(X,\F)\to \HH^0(X,\coker (\F\to q_\star M))$ se déduit immédiatement de $\F\to q_\star M$.

\section{Calcul de \texorpdfstring{$\RG(X,-)\colon \DD^b_c(X,\Lambda)\to \DD^b_c(\Lambda)$}{RGamma:Dbc(X,Lambda)->Dbc(Lambda)}}\label{sec:algconstr}

Le but de cette section est de donner un algorithme qui calcule explicitement le foncteur dérivé des sections globales sur une courbe lisse ou nodale sur un corps algébriquement clos. Dans un premier temps, nous décrivons une méthode pour calculer le $\RG$ d'un faisceau constructible sur une courbe lisse. Nous donnons ensuite des bornes précises sur la complexité de cet algorithme. Nous exposons ensuite la généralisation de cette méthode au calcul du $\RG$ d'un complexe de faisceaux constructibles, et terminons par un exemple détaillé.

\subsection{Calcul de $\RG(X,\F)$ pour $\F$ constructible sur $X$ lisse}\label{subsec:algconstr}

Soit $X$ une courbe intègre lisse sur le corps algébriquement clos $k$. Nous notons toujours $\Lambda=\ZZ/n\ZZ$, avec $n$ inversible dans $k$. Soit $\F$ un faisceau constructible de $\Lambda$-modules sur $X$, lisse sur un ouvert affine $j\colon U\hookrightarrow X$. Nous allons montrer comment calculer $\RG(X,\F)$. Notons $\mathscr{L}$ le faisceau lisse $j^\star\F$. Soit $i\colon Z\to X$ l'inclusion du fermé réduit complémentaire. Notons $\phi \colon i^\star\F \to i^\star j_\star \mathscr{L}$ le morphisme de recollement.

\subsubsection{Réduction à un cas particulier}

 Considérons le morphisme 
\[ f\colon \F\to j_\star \mathscr{L}\oplus i_\star i^\star \F \]
défini par adjonction. Sa fibre en $z\in Z$ est le morphisme injectif \[ f_z=(\phi_z,\id)\colon \F_z\to (j_\star \mathscr{L})_z\oplus\F_z. \]
De plus, $f|_U$ est simplement l'isomorphisme canonique $\mathscr{L}\to j^\star j_\star \mathscr{L}$.
Notons $Q_Z$ le conoyau du morphisme $f_Z=\bigoplus_{z\in Z}f_z$. Les morphismes \[(\id,-\phi_z)\colon (j_\star \mathscr{L})_z\oplus\F_z\to (j_\star \mathscr{L})_z \]
induisent un isomorphisme \[Q_Z\xrightarrow{\sim} \bigoplus_{z\in Z}(j_\star \mathscr{L})_z.\]
Comme $f$ est un isomorphisme sur $U$, le conoyau de $f$ est $i_\star Q_Z$. La suite exacte
\[ 0\to \F\to j_\star \mathscr{L}\oplus i_\star i^\star\F\to i_\star Q_Z \to 0\]
donne le triangle distingué
\[ \RG(X,\F)\to \RG(X,j_\star \mathscr{L}) \oplus \RG(X,i_\star i^\star \F)\to \RG(X,i_\star Q_Z)\xrightarrow{+1} \]
qui fait apparaître $\RG(X,\F)[1]$ comme le cône du morphisme de droite.
Pour calculer $\RG(X,\F)$, il suffit donc désormais,  de calculer $\RG(X,j_\star \mathscr{L})$ et le morphisme $\RG(X,j_\star \mathscr{L})\to \RG(Z,Q_Z)$.

\subsubsection{Calcul de $\RG(X,j_\star\mathscr{L})$ pour $\mathscr{L}$ lisse sur $U$}

Soit $\mathscr{L}$ un faisceau lisse sur la courbe affine lisse $U$, trivialisé par le revêtement galoisien $V\to U$. Notons $V_2\to V$ le revêtement caractéristique de $V$ de groupe $\HH^1(V,\Lambda)^\vee$, et $\overline{V_2}$ sa complétion projective lisse. Notons $G=\Aut(V_2|U)$, et $M=\HH^0(V,\mathscr{L}|_V)$. Rappelons que $\RG_Z$ et $\HH^i_z$ désignent les foncteurs de cohomologie à support sur $Z$ définis dans la section \ref{subsec:cohsuppferm}. Le triangle distingué 
\[ \RG_Z(X,j_\star\mathscr{L})\to \RG(X,j_\star \mathscr{L})\to \RG(U,\mathscr{L})\xrightarrow{+1}\]
fait apparaître $\RG(X,j_\star\mathscr{L})[1]$ comme le cône de $\RG(U,\mathscr{L})\to \RG_Z(X,j_\star\mathscr{L})[1]$. Nous avons vu dans la section \ref{subsec:cohlisG} comment calculer un complexe représentant $\RG(U,\mathscr{L})=\tau_{\leqslant 1}\RG(G,M)$. Soit $z\in Z$. Comme $\HH^i_z(X,j_\star \mathscr{L})=0$ pour tout $i\neq 2$ d'après le lemme \ref{lem:cohsuppz}, il y a dans $\DD^b_c(X,\Lambda)$ un isomorphisme \[\RG_Z(X,j_\star\mathscr{L})\xrightarrow{\sim}\HH^2_z(X,j_\star \mathscr{L})[-2].\]
Rappelons (cf lemme \ref{lem:cohsuppferme}) qu'il y a un isomorphisme $\HH^2_z(X,j_\star\mathscr{L})=\HH^1(I_z/P_z,M_{P_z})$ où $I_z\subseteq G$ est le groupe d'inertie en un point de $\overline{V_2}$ d'image $z$, et $P_z$ son sous-groupe d'inertie sauvage. Rappelons également que $I_z/P_z$ est canoniquement isomorphe à $\mu_e(k)$, où $e$ est l'indice de ramification de $\overline{V_2}\to X$ en $z$. Le morphisme $\RG(U,\mathscr{L})\to \RG_Z(X,\mathscr{L})[1]$ n'est donc autre que \[ \tau_{\leqslant 1}\RG(G,M)\to \bigoplus_{z\in Z}\HH^1(I_z/P_z,M_{P_z})[-1].\] Ce morphisme est nul en tout degré sauf 1, et en degré 1 c'est le morphisme qui à un 1-cocycle $c\colon G\to M$ associe l'image dans $\HH^1(I_z,M)$, puis dans $\HH^1(I_z/P_z,M_{P_z})$, de la classe de $c$ dans $\HH^1(G,M)$. La description du morphisme $\HH^1(I_z,M)\to \HH^1(I_z/P_z,M_{P_z})$ se trouve dans le lemme \ref{lem:secIP}.

\subsubsection{Conclusion de l'argument}
Reprenons le fil : $X$ est une courbe intègre lisse sur $k$, le faisceau $\F$ est constructible sur $X$, lisse sur l'ouvert $j\colon U\to X$ de fibre $M$. Notons $\mathscr{L}$ le faisceau lisse $j^\star\F$. Le paragraphe précédent permet de calculer un complexe représentant $\RG(X,j_\star \mathscr{L})$.  Notons $C^\bullet(G,M)$ le complexe de cochaînes usuel représentant $\RG(G,M)$. Notons encore $C^{12}(G,M)=\ker(C^1(G,M)\to C^2(G,M))$ le groupe des morphismes croisés de $G$ dans $M$. Le complexe calculé représentant $\RG(X,j_\star \mathscr{L})$ est le suivant :
\[C^0(G,M) \to C^{12}(G,M)\to \bigoplus_z \HH^1(I_z/P_z,M_{P_z}).\]
Calculons explicitement le morphisme 
\[ \RG(X,j_\star \mathscr{L})\oplus \RG(X,i_\star i^\star\F)\xrightarrow{(\id,-\phi_Z)} \RG(X,i_\star i^\star j_\star \mathscr{L}).\]
Le complexe \[ \RG(X,i_\star i^\star j_\star \mathscr{L})\simeq \HH^0(Z,i^\star j_\star \mathscr{L})[0]\simeq \bigoplus_{z\in Z}M^{I_z}[0]\]
est naturellement quasi-isomorphe au complexe
\[ \bigoplus_{z\in Z}M_{P_z} \to \bigoplus_{z\in Z}C^{12}(I_z/P_z,M_{P_z})\to \bigoplus_{z\in Z}\HH^1(I_z/P_z,M_{P_z}). \]
Le morphisme $\RG(X,j_\star \mathscr{L})\oplus\RG(X,i_\star i^\star \F)\to \RG(X,i_\star i^\star j_\star \mathscr{L})$ est donc le suivant.
\[
\begin{adjustbox}{width=\textwidth}{
\begin{tikzcd}
 0\arrow[r] &M\oplus \bigoplus_{z\in Z}\F_z \arrow[d,"\bigoplus_{z\in Z}(\id-\phi_z)"] \arrow[r,"{(\partial_G,0)}"] & C^{12}(G,M) \arrow[d,"\bigoplus_{z\in Z} \res_{I_z}^G"]\arrow[r] & \bigoplus_{z\in Z}\HH^1(I_z/P_z,M_{P_z}) \arrow[d,"\id"]\arrow[r] & 0 \\
 0 \arrow[r]& \bigoplus_{z\in Z}M_{P_z} \arrow[r,"\bigoplus_z\partial_{I_z}"] \arrow[r] & \bigoplus_{z\in Z}C^{12}(I_z/P_z,M_{P_z}) \arrow[r] & \bigoplus_{z\in Z}\HH^1(I_z/P_z,M_{P_z}) \arrow[r] & 0
\end{tikzcd}
}
\end{adjustbox}
\]
Les morphismes induits en cohomologie sont bien $\HH^0(U,\mathscr{L})\oplus \HH^0(Z,i^\star\F)\to \HH^0(X,j_\star \mathscr{L})$ en degré 0, et $\HH^i(X,j_\star \mathscr{L})\to 0$ en degré $i\geqslant 1$.
Le complexe $\RG(X,\F)[1]$ est le cône de ce morphisme. Ceci démontre le théorème \ref{th:compRG}.

\subsubsection{Fonctorialité sur $\Spec k$}

Soit $\phi\colon X'\to X$ un morphisme de courbes intègres lisses sur $k$. Soit $\F$ un faisceau constructible sur $X$, lisse sur un ouvert affine $U$. Nous allons décrire comment calculer le morphisme $\RG(X,\F)\to \RG(X',\phi^\star\F)$.
 Soit $f\colon V\to U$ un revêtement galoisien de $U$ trivialisant $\F|_U$. Soit $V_2\to V$ le revêtement de groupe $\HH^1(V,\Lambda)^\vee$, et $\overline{V_2}$ sa complétion projective lisse. Soit $V'$ la clôture galoisienne d'une composante connexe de $V\times_X X'$. Notons $U'=U\times_X X'$. Pour chaque point $z\in X-U$, choisissons un point $z_2\in \overline{V_2}$ d'image $z$. Notons $V'_2$ le revêtement de $V'$ de groupe $\HH^1(V',\Lambda)^\vee$, et $\overline{V'_2}$ sa complétion projective lisse. Pour chaque antécédent $z'$ de $z$ dans $X'-U'$, choisissons un antécédent $z'_2\in \overline{V'_2}$ de $z'$ d'image $z_2$. Alors il y a un morphisme $\Aut(V_2'|U')\subseteq \Aut(V_2'|U)\to \Aut(V_2|U)$ qui induit pour tout $z$ un morphisme $I_{z'_2}\to I_{z'}$, et par conséquent un morphisme $I_{z'_2}/P_{z'_2}\to I_{z_2}/P_{z_2}$. Le fonctorialité de la résolution bar permet d'en déduire le morphisme $\RG(X,\F)\to \RG(X',\phi^\star\F)$ cherché.

\subsubsection{Action de $\Gk$} 

Supposons désormais que $X,U,Z,\F$ proviennent de $X_0,U_0,Z_0,\F_0$ sur le corps de base $k_0$. Soit $V_0$ un revêtement galoisien de $U_0$ qui trivialise $\F|_{U_0}$.  Enfin, $\mathfrak{G}_0$ désigne toujours le groupe $\Gk$. Lorsque $V_0$ est géométriquement connexe, nous allons décrire l'action de $\mathfrak{G}_0$ sur le complexe $\RG(X,\F)$ déterminé précédemment. Nous donnerons ensuite un complexe de $\mathfrak{G}_0$-modules représentant $\RG(X,\F)$ lorsque $V_0$ n'est pas géométriquement connexe.\\

\paragraph{Si $V_0$ est géométriquement connexe} Supposons d'abord $V\coloneqq V_0\times_{k_0}k$ connexe. Soit $V_2\to V$ le revêtement de groupe $\HH^1(V,\Lambda)^\vee$, et soit $\overline{V_2}$ sa compactification lisse. Notons $G=\Aut(V_2|U)$, et $M$ le $G$-module $\HH^0(V,\F)$. Le lemme suivant montre la compatibilité de l'action de $\mathfrak{G}_0$ sur $G$ à celle de $G$ sur $M$.

\begin{lem} Soient $g\in G$ et $\sigma\in \mathfrak{G}_0$. Pour tout $m\in M$,
\[ g\cdot m=(\sigma\cdot g)\cdot m.\]
\begin{proof}
Il suffit de montrer que $g^{-1}\circ (\sigma\cdot g)$ agit trivialement sur $M$. Soit $\phi\in \Aut(V|U)$ l'image de $g$. Alors $\sigma\cdot g=\sigma g\sigma^{-1}$ a encore pour image $\phi$, car $\mathfrak{G}_0$ agit trivialement sur $\Aut(V|U)$. Rappelons que $k(V_2)=k(V)(z_1,\dots,z_r)$ où les $z_i$ sont des racines $n$-ièmes d'éléments de $k(V)$. Alors il existe $\zeta_1,\dots,\zeta_r\in \mu_n(k)$ tels que $(\sigma g)(z_i)=\zeta_i g(z_i)$, et $g^{-1}\circ (\sigma g)(z_i)=\zeta_iz_i$. Par conséquent, $g^{-1}\circ\sigma\cdot g\in \Aut(V_2|V)$ agit trivialement sur $M$.
\end{proof}
\end{lem}
Soient $z\in Z(k_0)$, et $z'\in \overline{V_2}$ un antécédent de $z$. Soient $I_{z'},P_{z'}$ les groupes d'inertie associés. Soient $\sigma\in \mathfrak{G}_0$ et $z''=\sigma(z')$. La conjugaison par $\sigma$ définit un isomorphisme $I_{z'}/P_{z'}\to I_{z''}/P_{z''}$.
Le lemme précédent assure que $M^{P_{z'}}=M^{P_{z''}}$ ; de plus, l'action de $\mathfrak{G}_0$ sur $M^{P_{z'}}$ est triviale. Comme $I_{z'}/P_{z'}$ est canoniquement isomorphe à $\mu_e(k)$, où $e$ est l'indice de ramification de $V_2\to U$ en $z'$, il est possible de calculer l'action de $\mathfrak{G}_0$ sur $\RG(\mu_e(k),M^{P_{z'}})$ grâce au diagramme commutatif suivant :
\[
\begin{tikzcd}
I_{z'}/P_{z'} \arrow[r,"\sigma"]\arrow[d,"\sim"{anchor=south, rotate=90}] & I_{z''}/P_{z''} \arrow[d,"\sim"{anchor=north, rotate=90}] \\
\mu_e(k) \arrow[r] & \mu_e(k)
\end{tikzcd}
\]
L'action sur les groupes de cohomologie de ce complexe, qui sont $(j_\star \F|_{U_0})_{\bar z}$ et $\HH^2_z(X,j_\star\F|_{U_0})$, s'en déduit immédiatement. 

Soit désormais $z\in Z(k_1)$, où $k_1$ est une extension finie de $k_0$. Appelons $W$ l'ensemble des conjugués de $z$ sous $\mathfrak{G}_0$ et calculons le complexe $\bigoplus_{w\in W} \RG(I_{w'}/P_{w'},M^{P_w'})$, où un antécédent $w'$ de chaque point $w$ de $W$ a été choisi dans une même $\mathfrak{G}_0$-orbite. L'action de $\mathfrak{G}_0$ sur ce complexe permute les différents facteurs. L'action de $\mathfrak{G}_0$ sur les groupes $\HH^0(W,j_\star \F|_{U_0})$ et $\HH^2_W(X,j_\star\F|_{U_0})$ s'en déduit. 

\paragraph{Si $V_0$ n'est pas géométriquement connexe} Considérons maintenant le cas où $V$ n'est pas connexe : le corps de fonctions de $V_0$ contient une extension galoisienne $L$ de $k_0$.
La construction de la section \ref{subsec:X2gal} permet d'obtenir un revêtement galoisien $V_{2,0}\to U_0$ dont le corps de fonctions contient une extension galoisienne $L'$ de $k_0$ suffisamment grande sur laquelle sont définis les éléments de $\HH^1(V,\Lambda)$ ainsi que les points de la compactification lisse $Y_2$ de $V_2\coloneqq V_{2,0}\times_{k_0}k$ au-dessus des points de $Z$. Ceci a l'avantage d'assurer que le corps résiduel de tout point fermé de $Y_{2,0}$ au-dessus d'un point de $Z_0$ soit exactement $L'$. 
Considérons une composante connexe $V'$ de $V_2=V_{2,0}\times_{k_0}k$ ; le choix de $V'$ n'a aucune incidence sur ce qui suit. Le groupe $G'\coloneqq\Aut(V'|U)=\ker(\Aut(V_{2,0}|U_0)\to \Gal(L|k_0))$ est le stabilisateur de $V'$ dans $G\coloneqq\Aut(V_{2,0}|U_0)$. 
Soient $M_0= \HH^0(V_{2,0},\F)$ et $M= \HH^0(V_2,\F)$ ; alors $M=\ind_{G'}^G (M_0)$, et le lemme de Shapiro \cite[Th. 4.19]{neukirch_cft} assure que \[ \tau_{\leqslant 1}\RG(G,M)=\tau_{\leqslant 1}\RG(G',M_0)=\RG(U,\F)\]
dans $\DD^b_c(\Lambda)$. En tant que groupe abélien, $M=M_0^d$ où $d$ est le nombre de composantes connexes de $V_2$. 
Le groupe $\mathfrak{G}_0$ agit naturellement sur $M$ (en permutant ses facteurs), de façon compatible à son action sur $G$. Ceci permet de calculer l'action de $\mathfrak{G}_0$ sur $\tau_{\leqslant 1}\RG(G,M)$.
Soit $z\in Z_0$ un point fermé de corps résiduel $k_1$. 
Soit $v_2$ un antécédent de $z$ dans $Y_{2,0}$. Notons $W\coloneqq \{\tau(z),\tau\in \Gal(k_1|k_0)\}$ la $\mathfrak{G}_0$-orbite de $z$. Choisissons pour tout $\tau\in \Gal(k_1|k_0)$ un antécédent $v_{2,\tau}$ de $\tau(z_0)$ dans $V_{2,0}$ ; prenons les $v_{2,\tau}$ dans une même orbite sous $\mathfrak{G}_0$.
Soit $D_\tau\triangleleft G$ le groupe de décomposition de $v_{2,\tau}$. Le morphisme $D_\tau\to \Gal(L'|k_0)$ est surjectif, de noyau le groupe d'inertie $I_\tau$ de $v_{2,\tau}$, de sorte que $M=\ind_{I_\tau}^{D_\tau}M_0$. Par le lemme de Shapiro, \[ \RG(D_\tau,M)=\RG(I_\tau,M_0)=\RG(I_\tau/P_\tau,M_0^{P_\tau}) \] dans $\DD^b_c(\Lambda)$. 
Le groupe $\mathfrak{G}_0$ agit naturellement sur le complexe $\bigoplus_{\tau}\RG(D_\tau,M)$, dont les groupes de cohomologie respectifs sont $\HH^0(W,j_\star\F|_U)$ et $\HH^2_W(X,j_\star\F|_U)$. Avec ces notations, le complexe $\RG(X,\F)[1]$ est isomorphe dans $\DD^b_c(\Lambda[\mathfrak{G}_0])$ au cône du morphisme de complexes :
\[
\begin{adjustbox}{width=\textwidth}{
\begin{tikzcd}
 0\arrow[r] & M\oplus \bigoplus_{W} \HH^0(W,\F) \arrow[d,"\bigoplus_{W,\tau}(\res_{I_\tau}^G-\phi_z)"] \arrow[r,"{(\partial_G,0)}"] & C^{12}(G,M) \arrow[d,"\bigoplus_{W,\tau} \res_{I_\tau}^G"]\arrow[r] & \bigoplus_{W}\bigoplus_{\tau}\HH^1(I_\tau/P_\tau,M_{P_\tau}) \arrow[d,"\id"]\arrow[r] & 0 \\
 0 \arrow[r]& \bigoplus_{W}\bigoplus_{\tau}M_{P_\tau} \arrow[r,"\bigoplus_z\partial_{I_\tau}"] \arrow[r] & \bigoplus_{W}\bigoplus_{\tau}C^{12}(I_\tau/P_\tau,M_{P_\tau}) \arrow[r] & \bigoplus_{W}\bigoplus_{\tau}\HH^1(I_\tau/P_\tau,M_{P_\tau}) \arrow[r] & 0
\end{tikzcd}
}
\end{adjustbox}
\]
où $W$ parcourt les $\mathfrak{G}_0$-orbites de $Z$, et $\tau$ parcourt les $k_0$-automorphismes du corps résiduel des points fermés de $W$.

\subsection{Complexité}

Cet algorithme de calcul de $\RG(X,\F)$ se compose essentiellement de deux étapes : le calcul de $G$ et de ses sous-groupes $I_z$, puis des calculs d'algèbre linéaire.
Soit $C=\Spec A$ le modèle plan de $V$ utilisé par l'algorithme. Soit $z\in X$, et $z'\in \bar V-V$ un antécédent de $z$. Soit $(f_1,\dots,f_t)$ la base de $\HH^1(V,\mu_n)$ calculée par l'algorithme ; comme indiqué dans la remarque \ref{rk:choixfct}, nous pouvons supposer que l'ordre des $f_i$ en $z$ est positif. Pour tout $i\in \{1\dots t\}$, il existe donc  des fonctions polynomiales $g_i,h_i$, où $h_i$ ne s'annule pas en $z'$, telles que \[f_i=\frac{g_i}{h_i}.\]
Le construction explicite de ces fonctions est donnée dans \cite[§1.5, p15]{shafarevich1}.
Alors $V_2$ est la normalisée de la courbe $C_2\coloneqq \Spec A[z_1,\dots,z_t]/(h_iz_i^n-g_i)$, où les $f_i$ forment une base de $\HH^1(V,\mu_n)$. Le point $z'$ a $n^{t-1}$ antécédents dans l'homogénéisée de $C_2$, et ceux-ci sont encore non singuliers, ce qui se lit sur les équations. Soit $z''$ l'un d'entre eux. Le groupe d'inertie $I_{z''}$ peut être calculé simplement en évaluant les éléments de $\Aut(V_2|U)$ en $z''$.  \\

Pour le résultat suivant, notons encore $C(g,q,n,r)$ la complexité du calcul de $\HH^1(U,\mu_n)$, où $U$ est une courbe intègre lisse sur $\overline{\FF_q}$ provenant de $\FF_q$, de compactification lisse $X$ de genre $g$, avec $r=|X-U|$. Notons également $D(g,n,r)$ le degré de la plus petite extension de $\FF_q$ sur laquelle sont définis les éléments de $\HH^1(U,\mu_n)$ obtenus. Rappelons que des majorations de ces entiers ont été obtenues dans la section \ref{subsubsec:detrep}

\begin{theorem} Soit $X_0$ une courbe lisse géométriquement connexe sur $\FF_q$. Notons $X=(X_0)_{\overline{\FF_q}}$. Soit $\bar X$ sa compactification lisse. Notons $g$ le genre de $\bar X$. Soit $\F$ un faisceau constructible de $\Lambda$-modules sur $X$, lisse sur un ouvert $U$ de $X$, de fibre générique géométrique un $\Lambda$-module $M$. Soit $d$ l'ordre de l'image du morphisme $\pi_1(X)\to\Aut_\Lambda(M)$. Notons $m$ un majorant du nombre de générateurs de $M$ et de chacun des $\F_z$, $z\in X-U$. Enfin, notons $r=|X-U|$. Supposons donné un modèle plan d'un revêtement galoisien $V$ de $U$ trivialisant $\F$, et qui n'a pas de points singuliers au-dessus de $X-U$. Il existe un algorithme probabiliste (Las Vegas) qui calcule un complexe de $\Lambda[\mathfrak{G}_0]$-modules représentant $\RG(X,\F)$ en 
\[ O\left(C(d(g+r),q,n,dr)+r(md^{3}n^{4d(g+2r)}D(d(2g+r),n,dr))^3\right)\]
opérations dans $\FF_q$.
\begin{proof} Le calcul est à peine plus compliqué dans le cas lisse (voir théorème \ref{th:complis}) : une fois déterminé le groupe $\Aut(V_2|U)$, le calcul des groupes d'inertie en découle comme décrit ci-dessus. Le nombre de calculs d'algèbre linéaire à effectuer est proportionnel au nombre de points de $\bar X-U$.
\end{proof}
\end{theorem}

La borne obtenue sur $C(g,q,n,r)$ dans la proposition \ref{th:huangrac} permet d'en déduire le résultat suivant.

\begin{cor}\label{th:complexiteRG}
Avec les notations et hypothèses du théorème, il existe un algorithme probabiliste (Las Vegas) qui calcule $\RG(X,\F)$ en
\[ \mathcal{P}(d,g,n,r,m,\log q)^{2^{O\left((d(g+r))^2\right)}}\]
opérations dans $\FF_q$, où $\mathcal{P}$ est un polynôme. Si $V$ admet un modèle plan à singularités ordinaires de degré $O(g)$, cette complexité devient \[ \mathcal{P}(d,g,n,r,m,\log q)^{O\left((d(g+r))^4\right)}.\]
\end{cor}

\subsection{Adaptation au cas des courbes nodales}

Soit $X$ une courbe intègre à singularités au pire nodales sur le corps algébriquement clos $k$. Notons $\nu\colon \tX\to X$ sa normalisation. Soit $\F$ un faisceau constructible de $\Lambda$-modules sur $X$. Soit $j\colon U\to X$ l'inclusion d'un ouvert affine lisse de $X$ tel que le faisceau $\mathscr{L}\coloneqq j^\star\F$ soit lisse, et $i\colon Z\to X$ l'inclusion du fermé réduit complémentaire. Notons $M$ la fibre générique géométrique de $\mathscr{L}$.
Comme dans le cas des courbes lisses, nous calculons $\RG(X,\F)$ à l'aide des deux triangles distingués
\[ \RG(X,\F)\to \RG(X,j_\star \mathscr{L}) \oplus \RG(X,i_\star i^\star \F)\to \RG(X,i_\star i^\star j_\star\F)\xrightarrow{+1}\]
et
\[ \RG_Z(X,j_\star\mathscr{L})\to \RG(X,j_\star \mathscr{L})\to \RG(U,\mathscr{L})\xrightarrow{+1}\]
Soit $V\to U$ un revêtement galoisien trivialisant $\mathscr{L}$. Comme $U$ est lisse, $V$ l'est également. Notons $\tilde V$ sa normalisée, et $V_2\to V$ le revêtement de groupe $\HH^1(V,\Lambda)^\vee$. Notons $G$ le groupe $\Aut(V_2|U)$. Le morphisme $\tau_{\leqslant 1}\RG(G,M)\to \RG(U,\mathscr{L})$ est toujours un quasi-isomorphisme. 
Notons $\tilde{Z}$ la préimage de $Z$ dans $\tX$. Le complexe $\RG_Z(X,j_\star\mathscr{L})$ est canoniquement quasi-isomorphe à $\RG_{\tilde Z}(\tX,\nu^\star j_\star\mathscr{L})$ d'après le corollaire \ref{cor:cohsuppnodiso}. Pour chaque point régulier $z$ de $X$ appartenant à $Z$, notons $I_z$ son stabilisateur dans $G$. Pour chaque point nodal $z$ de $X$ appartenant à $Z$, notons $Q,R$ ses antécédents dans $\tX$, et $I_Q,I_R$ leurs stabilisateurs respectifs dans $G$. Le complexe suivant représente donc $\RG(X,j_\star\mathscr{L})$ :
\[
M\to C^{12}(G,M) \to \bigoplus_{z\in \tilde Z}\HH^1(I_z,M) 
\]
Comme précédemment, il est possible de remplacer ici $\HH^1(I_z,M)$ par $\HH^1(I_z/P_z,M_{P_z})$ où $P_z$ est le sous-groupe d'inertie sauvage de $I_z$ (et de même pour les points $Q,R$). Rappelons que pour tout $z\in Z$, $(j_\star\mathscr{L})_z=\HH^0(I_Q,M)\times\HH^0(I_R,M)$. Nous noterons $\phi_z\colon \F_z\to \HH^0(I_Q,M)\times\HH^0(I_R,M)$ le morphisme de recollement. Par conséquent, $\RG(X,\F)[1]$ est représenté par le cône du morphisme de complexes :

\[
\begin{adjustbox}{width=\textwidth}{
\begin{tikzcd}
 0\arrow[r] &M\oplus \bigoplus_{z\in Z}\F_z \arrow[d,"\bigoplus_{z\in Z}(\id-\phi_z)"] \arrow[r,"{(\partial_G,0)}"] & C^{12}(G,M) \arrow[d,"\bigoplus_{\tilde{z}\in Z} \res_{I_{\tilde z}}^G"]\arrow[r] & \bigoplus_{\tilde{z}\in \tilde Z}\HH^1(I_{\tilde z}/P_{\tilde z},M_{P_{\tilde z}}) \arrow[d,"\id"]\arrow[r] & 0 \\
 0 \arrow[r]& \bigoplus_{\tilde z\in \tilde Z}M_{P_{\tilde z}} \arrow[r,"\bigoplus_{\tilde z}\partial_{I_{\tilde z}}"] \arrow[r] & \bigoplus_{\tilde z\in \tilde Z}C^{12}(I_{\tilde z}/P_{\tilde z},M_{P_{\tilde z}}) \arrow[r] & \bigoplus_{\tilde z\in \tilde Z}\HH^1(I_{\tilde z}/P_{\tilde z},M_{P_{\tilde z}}) \arrow[r] & 0
 \end{tikzcd}
}
\end{adjustbox}
\]

\subsection{Adaptation au cas des complexes}\label{subsec:algconstrcomp}

Soit $X$ une courbe intègre lisse sur $k$. Soit $j\colon U\to X$ l'inclusion d'un ouvert. Soit $i\colon Z\to X$ le fermé réduit complémentaire. Soit $\K^\bullet=[\K^0\to\cdots\to \K^b]$ un complexe de faisceaux constructibles sur $X$ tel que tous les $\K^i$ soient lisses sur $U$, et que les $j^\star\K^i$ soient trivialisés par $V\to U$. Notons $\mathscr{L}=j^\star\K$. Notons $K^\bullet$ le complexe de $G=\Aut(V_2|U)$-modules associé. Les foncteurs $j_\star,j^\star,i_\star,i^\star$ s'étendent terme à terme aux catégories (non dérivées) de complexes de faisceaux constructibles sur $X$. L'utilisation surprenante de ces foncteurs non dérivés s'explique par notre stratégie pour calculer $\RG(X,\K)$ : elle consiste à insérer $\K$ dans une suite exacte dans la catégorie des complexes de faisceaux constructibles sur $X$, puis d'en déduire un triangle distingué dans $\DD^b_c(\Lambda)$. En particulier, si $\K'$ est un autre complexe quasi-isomorphe à $\K$, les complexes obtenus représentant $\RG(X,\K)$ et $\RG(X,\K')$ seront distincts ; étant donné un quasi-isomorphisme explicite entre $\K$ et $\K'$, il sera toutefois possible de calculer le quasi-isomorphisme $\RG(X,\K)\to\RG(X,\K')$ qu'il induit.

Afin de simplifier l'exposition, supposons que $Z$ soit réduit à un seul point fermé $z$. Fixons un point $v\in V_2$ au-dessus de $z$, et notons $I_z=\{\sigma\in G\mid \sigma(v)=v \}$. Pour le cas général, il suffira de remplacer $C^j(I_z,-)$ (resp. $\HH^j(I_z,-)$) par $\bigoplus_{z\in Z}C^j(I_z,-)$ (resp. $\bigoplus_{z\in Z}\HH^j(I_z,-)$).
Il y a une suite exacte de complexes de faisceaux constructibles :
\[ 0\to \K \to j_\star \mathscr{L}\oplus i_\star i^\star \K \to i_\star Q \to 0\]
où $Q\simeq i^\star j_\star \mathscr{L}$, ce qui se déduit terme à terme de l'énoncé correspondant pour les faisceaux.
L'élément $\RG(X,\K)[1]$ de $\DD^b_c(\Lambda)$ est donc le cône de \[\RG(X,j_\star \mathscr{L})\oplus\RG(X,i_\star i^\star \K)\to \RG(X,i_\star Q).\]

\subsubsection{Calcul de $\RG(X,j_\star \mathscr{L})$}

Comme précédemment, il y a un triangle distingué dans $\DD^b_c(X,\Lambda)$ \cite[09XP]{stacks}:

\[ \RG_Z(X,j_\star \mathscr{L})\to \RG(X,j_\star \mathscr{L})\to \RG(U,\mathscr{L})\xrightarrow{+1} \]
qui montre que $\RG(X,j_\star \mathscr{L})[1]$ est le cône de $\RG(U,\mathscr{L})\to \RG_Z(X,j_\star \mathscr{L})[1]$. 

\begin{lem} Le complexe \[ 0\to 0\to \HH^1(I_z,K^0)\to \HH^1(I_z,K^1)\to \HH^1(I_z,K^2)\to\cdots \]
(où les termes ci-dessus sont placés en degrés $\geqslant 0$) représente $\RG_Z(X,j_\star \mathscr{L})$.
\begin{proof}
Soit $I^{\bullet,\bullet}$ une résolution de Cartan-Eilenberg de $j_\star \mathscr{L}$. La colonne $I^{p,\bullet}$ est alors une résolution injective de $j_\star \mathscr{L}^p$. Rappelons que pour chaque entier $p$ et chaque point $z$ de $Z$, les $\HH^i_z(X,j_\star \mathscr{L}^p)$ sont nuls dès que $i\neq 2$, et $\HH^2_z(X,j_\star \mathscr{L}^p)=\HH^1(I_z,K^p)$. Les morphismes
\[ \tau_{q\geqslant 2}\tau_{q\leqslant 2} \Gamma_Z(X,I^{\bullet,\bullet}) \leftarrow \tau_{q\leqslant 2}\Gamma_Z(X,I^{\bullet,\bullet}) \to \Gamma_Z(X,I^{\bullet,\bullet}) \]
induisent donc des isomorphismes entre les premières pages des suites spectrales associées à ces complexes doubles pour l'orientation verticale, et donc des quasi-isomorphismes entre les complexes totaux associés. Comme le complexe double de gauche a exactement un terme non nul dans chaque colonne, à hauteur 2, son complexe total associé est celui de l'énoncé du lemme.
\end{proof}
\end{lem}

Le morphisme $\RG(U,\mathscr{L})\to \RG_Z(X,j_\star \mathscr{L})[1]$ est donc représenté par le morphisme de complexes
\[
\begin{tikzcd} 
C^0(G,K^0) \arrow[r]\arrow[d]& C^{12}(G,K^0)\oplus C^0(G,K^1)\arrow[r]\arrow[d] & C^{12}(G,K^1)\oplus C^0(G,K^2)\arrow[d] \arrow[r] &\cdots\\
0 \arrow[r] & \HH^1(I_z/P_z,K^0_{P_z}) \arrow[r]&  \HH^1(I_z/P_z,K^1_{P_z}) \arrow[r] & \dots
\end{tikzcd}
\]
où les flèches verticales sont déduites du morphisme $I_z\to G$.
Le complexe $\RG(X,j_\star j^\star K)$ est donc 
\[ C^0(G,K^0) \to C^{12}(G,K^0)\oplus C^0(G,K^1) \to C^{12}(G,K^1)\oplus C^0(G,K^2)\oplus \HH^1(I_z/P_z,K^0_{P_z})\to\cdots\]
et son $i$-ième terme est 
\[ \HH^1(I_z/P_z,K^{i-2}_{P_z})\oplus C^{12}(G,K^{i-1})\oplus C^0(G,K^i).\]

\subsubsection{Calcul de $\RG(X,\K)$}

Il reste à déterminer le complexe $\RG(X,i_\star Q)$ ainsi que le morphisme \[\RG(X,j_\star \mathscr{L})\oplus\RG(X,i_\star i^\star\K)\to \RG(X,i_\star Q).\]

\begin{lem} Le complexe total associé au complexe double
\[
\begin{tikzcd}
\HH^1(I_z/P_z,K^0_{P_z}) \arrow[r]& \HH^1(I_z/P_z,K^1_{P_z}) \arrow[r] & \HH^1(I_z/P_z,K^2_{P_z})  \arrow[r] & \cdots\\
C^{12}(I_z/P_z,K^0_{P_z}) \arrow[r]\arrow[u] & C^{12}(I_z/P_z,K^1_{P_z}) \arrow[r]\arrow[u]& C^{12}(I_z/P_z,K^2_{P_z}) \arrow[u] \arrow[r] & \cdots\\
C^{0}(I_z/P_z,K^0_{P_z}) \arrow[r]\arrow[u] & C^{0}(I_z/P_z,K^1_{P_z}) \arrow[r]\arrow[u]& C^{0}(I_z/P_z,K^2_{P_z}) \arrow[u]\arrow[r] & \cdots
\end{tikzcd}
\]
représente $\RG(X,i_\star Q)=\RG(X,i_\star i^\star j_\star \mathscr{L})$. Son $i$-ième terme est
\[ L^i\coloneqq \HH^1(I_z/P_z,K^{i-2}_{P_z})\oplus C^{12}(I_z/P_z,K^{i-1}_{P_z})\oplus C^0(I_z,K^i_{P_z}).\]
\begin{proof} Pour tout complexe $A$ de groupes abéliens, $\RG(Z,A)=A$. En particulier, \[\RG(X,i_\star Q)=\RG(Z,Q)=Q=i^\star j_\star \mathscr{L}=(\HH^0(I_z/P_z,K^p_{P_z}))_{p\geqslant 0}.\] Il y a pour tout $p$ un quasi-isomorphisme de complexes \[ \HH^0(I_z/P_z,K^p_{P_z})[0]\to [C^0(I_z/P_z,K^p_{P_z})\to C^{12}(I_z/P_z,K^p_{P_z})\to \HH^1(I_z/P_z,K^p_{P_z})].\] Ce dernier étant compatible aux morphismes de transition de $K$, il induit un morphisme de complexes doubles entre le complexe ayant pour seule ligne $\HH^0(I_z/P_z,K^\bullet_{P_z})$ et le complexe de l'énoncé, qui induit un isomorphisme entre les premières pages des suites spectrales associées (pour l'orientation verticale). Par conséquent, il induit un quasi-isomorphisme entre les complexes totaux associés. 
\end{proof}
\end{lem}
Le morphisme $\RG(X,j_\star \mathscr{L})\to \RG(X,i_\star Q)$ est donc le suivant.
\[
\begin{tikzcd}
\cdots\arrow[r]&\HH^1(I_z/P_z,K^{i-2}_{P_z})\oplus C^{12}(G,K^{i-1}_{P_z})\oplus C^0(G,K^i_{P_z}) \arrow[d,"{(\id,\res_{I_z}^G,\res_{I_z}^G)}"] \arrow[r]&\cdots\\
\cdots\arrow[r]&\HH^1(I_z/P_z,K^{i-2}_{P_z})\oplus C^{12}(I_z/P_z,K^{i-1}_{P_z})\oplus C^0(I_z/P_z,K^i_{P_z})\arrow[r]&\cdots
\end{tikzcd}
\]
D'autre part, le morphisme $\RG(X,i_\star i^\star \K)\to \RG(X,i_\star Q)$ est le suivant.
\[
\begin{tikzcd}
\cdots\arrow[r]& \K^i_z \arrow[r]\arrow[d,"{(0,0,-\phi_z)}"] & \cdots \\
\cdots\arrow[r]&\HH^1(I_z,K^{i-2}_{P_z})\oplus C^{12}(I_z/P_z,K^{i-1}_{P_z})\oplus C^0(I_z/P_z,K^i_{P_z})\arrow[r]&\cdots
\end{tikzcd}
\]
Le complexe $\RG(X,\K)[1]$ est le cône de $\RG(X,j_\star \mathscr{L})\oplus\RG(X,i_\star i^\star\K)\to \RG(X,i_\star Q)$.
Nous avons démontré le théorème suivant.

\begin{theorem}
Soit $X$ une courbe intègre lisse sur $k$. Soit $\K=[\K^0\to\dots\to\K^b]$ un complexe de faisceaux constructibles de $\Lambda$-modules sur $X$. Soient $U$ un ouvert de $X$ sur lequel chaque $\K^i$ est lisse de fibre $K^i$, et $Z$ le fermé réduit complémentaire. Soient $V\to U$ un revêtement galoisien qui trivialise $\F|_U$, et $V_2\to V$ le revêtement de $V$ de groupe $\HH^1(V,\Lambda)^\vee$. Notons $G=\Aut(V_2|U)$. Pour chaque point $z\in Z$, notons $I_z$ le groupe d'inertie d'un point de la compactification lisse de $V_2$ au-dessus de $z$, et $P_z$ le groupe d'inertie sauvage correspondant. Notons $\phi_z^i\colon \K_z^i\to \HH^0({I_z},K^i)\subseteq \HH^0(P_z,K^i)\xrightarrow{\sim}K^i_{P_z}$ le morphisme de recollement en $z$. Alors $\RG(X,\F)[1]$ est le cône du morphisme de complexes suivant, où $C^{12}(G,M)$ désigne le groupe des morphismes croisés $G\to M$.
\[
\begin{tikzcd}
\cdots\arrow[r]&\HH^1(I_z/P_z,K^{i-2}_{P_z})\oplus C^{12}(G,K^{i-1}_{P_z})\oplus C^0(G,K^i_{P_z})\oplus\K^i_z \arrow[d,"{(\id,\res_{I_z}^G,\res_{I_z}^G-\phi_z)}"] \arrow[r]&\cdots\\
\cdots\arrow[r]&\HH^1(I_z/P_z,K^{i-2}_{P_z})\oplus C^{12}(I_z/P_z,K^{i-1}_{P_z})\oplus C^0(I_z/P_z,K^i_{P_z})\arrow[r]&\cdots
\end{tikzcd}
\]
\end{theorem}

\subsection{Un exemple détaillé}\label{subsec:exdetconstr}

\paragraph{Le cas étudié}
Supposons que $-1$ n'est pas un carré dans $k_0$. Posons $n=2$. Considérons les courbes $Y_0=X_0=\PP^1_{k_0}$, et le morphisme $f_0\colon Y_0\to X_0,y\mapsto y^2$. 
Soit $\F_0$ le faisceau $(f_0)_\star \Lambda$, lisse sur l'ouvert $U_0\coloneqq\GG_{m,k_0}=\Spec k[x^{\pm 1}]$ de $X_0$. Soit $V_0=\Spec k[y^{\pm 1}]$ sa préimage dans $Y_0$. 
Soient $Z_0,W_0$ les complémentaires réduits respectifs de $U_0,V_0$ dans $X_0,Y_0$. 
Le faisceau $\F_0$ est lisse sur $U_0$, trivialisé par $V_0$, de fibre $M=\Lambda^2$. L'automorphisme de $M$ induit par l'élément non trivial de $\Aut(V_0|U_0)\simeq\ZZ/2\ZZ$ échange les deux copies de $\Lambda$. 
Notons encore $X,Y,U,V,Z,W,\F$ les changements de base à $k$ respectifs de $X_0,Y_0,U_0,V_0,Z_0,W_0,\F_0$. Notons $j\colon U\to X$ l'inclusion, et $\mathscr{L}$ le faisceau lisse $j^\star\F$. Calculons $\RG(X,\F)$. 

\paragraph{Le revêtement trivialisant les $\mathscr{L}$-torseurs} 
Le revêtement $V_2\to V$ de groupe $\HH^1(V,\Lambda)$ est simplement $\GG_m\to \GG_m, z\mapsto z^2$. Par conséquent, $V_2\to U$ est le revêtement $\GG_m\to\GG_m,z\mapsto z^4$ de groupe $G=\ZZ/4\ZZ$ engendré par $\gamma\colon z\mapsto iz$, où $i$ est une racine carrée de $-1$ dans $k$. Les groupes d'inertie en $0$ et $\infty$ sont encore égaux à $G$. Par conséquent, $(j_\star j^\star)\F_0=\Lambda$ et $(j_\star \mathscr{L})_\infty=\Lambda$. Les morphismes de recollement $\F_0\to (j_\star \mathscr{L})_0$ et $\F_\infty\to (j_\star \mathscr{L})_\infty$ sont l'identité $\Lambda\to \Lambda$. 
\paragraph{Calcul de $\RG(U,\F|_U)$} 
Les morphismes croisés $G\to M$ sont uniquement déterminés par l'image $(a,b)\in\Lambda^2$ de $\gamma$. Le complexe de cochaînes usuel représentant $\RG(G,M)=\RG(U,\F|_U)$ est le suivant.
\[ \begin{array}{rcl}
\Lambda^2 &\longrightarrow& C^{12}(G,\Lambda^2)\\
(a,b)&\longmapsto& \left[ \gamma \mapsto (a+b,a+b)\right]
\end{array}\]
Le groupe $\HH^1(G,M)$ est donc isomorphe à $\Lambda$, et le morphisme $\Lambda^2\to \HH^1(G,M)$ qui associe à morphisme croisé sa classe de cohomologie a pour noyau $\langle (1,1)\rangle$ ; il s'identifie au morphisme \[ \begin{array}{rcl}\Lambda^2&\longrightarrow& \Lambda \\ (a,b)&\longmapsto& a+b.\end{array}\]

\paragraph{Calcul de $\RG(X,j_\star \mathscr{L})$} 
L'élément $\RG(X,j_\star \mathscr{L})[1]\in \DD^b_c(X,\Lambda)$ est le cône du morphisme \[\tau_{\leqslant 1}\RG(G,M)\to \HH^1(I_0,M)[-1]\oplus \HH^1(I_\infty,M)[-1]=\Lambda^2[-1].\] Par conséquent, $\RG(X,j_\star \mathscr{L})$ est représenté par le complexe 
\[ \Lambda^2\longrightarrow \Lambda^2\longrightarrow \Lambda^2 \]
où les deux morphismes sont définis par $(a,b)\mapsto (a+b,a+b)$. 

\paragraph{Calcul de $\RG(X,\F)$} 
Calculons désormais le morphisme 
\[ \RG(X,j_\star \mathscr{L})\oplus \RG(Z,i^\star\F) \to \RG(Z,i^\star j_\star \mathscr{L}).\]
D'une part, $\RG(Z,i^\star\F)=\HH^0(Z,i^\star\F)[0]=\F_0[0]\oplus \F_\infty[0]$. D'autre part, $\RG(Z,i^\star j_\star \mathscr{L})$ est représenté par le complexe 
\[ \begin{tikzcd}
\Lambda^2\oplus\Lambda^2 \arrow[r,"{\alpha'}"]& \Lambda^2\oplus\Lambda^2 \arrow[r,"{\beta'}"]& \Lambda^2 \end{tikzcd}\]
où les flèches sont $\alpha'\colon (a,b,c,d)\mapsto (a+b,a+b,c+d,c+d)$ et $\beta'\colon (a,b,c,d)\mapsto(a+b,c+d)$.
Le morphisme de complexes cherché est donc 
\[
\begin{tikzcd}
\Lambda^4 \arrow[r,"{\alpha}"]\arrow[d,"u"] & \Lambda^2 \arrow[r,"{\beta}"]\arrow[d,"v"] &\Lambda^2\arrow[d,"{\id}"] \\
\Lambda^4 \arrow[r,"{\alpha'}"] & \Lambda^4\arrow[r,"{\beta'}"] & \Lambda^2
\end{tikzcd}
\]
où, en écrivant le terme en haut à gauche comme $\F_0\oplus\F_\infty\oplus M$, \begin{itemize}[label=$\bullet$]
\item $u\colon(a,b,c,d)\mapsto (a+c,a+d,b+c,b+d)$
\item $\alpha \colon (a,b,c,d)\mapsto (c+d,c+d)$
\item $\beta \colon(a,b)\mapsto(a+b,a+b)$
\item $v\colon(a,b)\mapsto (a,b,a,b)$
\item $\alpha'\colon (a,b,c,d)\mapsto (a+b,a+b,c+d,c+d)$
\item $\beta'\colon(a,b,c,d)\mapsto(a+b,c+d)$.
\end{itemize} 
En calculant le cône de ce morphisme puis en décalant de 1, on obtient le complexe suivant, qui représente $\RG(X,\F)$.
\[
\begin{tikzcd}
\Lambda^4 \arrow[r,"\partial_0"] & \Lambda^6 \arrow[r,"\partial_1"] & \Lambda^6 \arrow[r,"\partial_2"] & \Lambda^2
\end{tikzcd}
\]
\begin{itemize}[label=$\bullet$]
\item $\partial_0\colon(a,b,c,d)\mapsto (c+d,c+d,a+c,b+c,a+d,b+d)$
\item $\partial_1 \colon(a,b,c,d,e,f)\mapsto (a+b,a+b,a+c+d,b+c+d,a+e+f,b+e+f)$
\item $\partial_2\colon(a,b,c,d,e,f)\mapsto (a+c+d,b+e+f)$. 
\end{itemize}
Ce complexe a pour groupes de cohomologie $\HH^0=\langle (1,1,1,1)\rangle$, $\HH^1=0$ et $\HH^2=\langle \overline{(1,0,1,0,0,0)}\rangle\simeq\Lambda$. C'est le résultat attendu : nous avons calculé la cohomologie du faisceau $(\PP^1\xrightarrow{x\mapsto x^2}\PP^1)_\star\Lambda$, qui est la cohomologie de $\Lambda$ sur $\PP^1$.

\paragraph{Action de Galois} L'action de $\Gk$ sur $\RG(X,\F)$ se factorise par celle de $\Gal(k_0(i)|k_0)$. Notons $\sigma\colon i\mapsto -i$ le $k_0$-automorphisme non trivial de $k_0(i)$. Le groupe $\mathfrak{G}_0$ agit sur $G=\Aut(V_2|V)$ par $\sigma\cdot \gamma=\gamma^3$, et trivialement sur $M$. Son action sur $\tau_{\leqslant 1}\RG(G,M)=\Lambda^2\to \Lambda^2$ est donc triviale sur le premier terme, et $(a,b)\mapsto (b,a)$ sur le second. En particulier, son action sur $\HH^1(G,M)=\Lambda$ est triviale.
L'action de $\sigma$ est également triviale sur $\F_0,\F_\infty$. En résumé, l'action sur le complexe 
\[ \Lambda^4\to \Lambda^6\to\Lambda^6\to \Lambda^2 \]
représentant $\RG(X,\F)$ est triviale sur le premier et le dernier terme, \[\sigma\cdot (a,b,c,d,e,f)= (b,a,c,d,e,f)\] sur le deuxième, et \[\sigma\cdot (a,b,c,d,e,f)= (a,b,d,c,f,e)\] sur le troisième. Il en découle en particulier que $\mathfrak{G}_0$ agit trivialement sur les deux groupes non nuls $\HH^0(X,\Lambda)$ et $\HH^2(X,\Lambda)$, comme attendu.

\cleartooddpage

\chapter{Cohomologie des surfaces}\label{chap:6}

Dans ce chapitre encore, $k$ désigne un corps algébriquement clos, et $n$ un entier inversible dans $k$. L'anneau $\ZZ/n\ZZ$ est toujours noté $\Lambda$. Nous indiquons comment utiliser les techniques présentées dans le chapitre \ref{chap:5} pour calculer la cohomologie d'un faisceau constant sur une $k$-surface lisse ; il reste à calculer précisément la complexité de l'algorithme obtenu. La technique utilisée est celle de la fibration en courbes projectives par le moyen d'un pinceau de Lefschetz, telle que suggérée dans \cite[Epilogue]{edixcouv}.

\section{Pinceaux de Lefschetz}\label{sec:lef}

Soit $X$ une variété projective connexe non singulière sur $k$, plongée dans $\PP^n_k$.
Un pinceau d'hyperplans dans $\PP^n$ est une droite $D$ de l'espace projectif dual $\check{\PP}^n$. Elle paramètre les hyperplans contenant un sous-espace linéaire $A$ de $\PP^n$ de codimension $2$ appelé l'axe de $D$. Pour tout $t\in D$, notons $H_t$ l'hyperplan correspondant. Soit $\tilde{X}$ le fermé réduit de $X\times D$ dont les points $(x,t)$ vérifient $x\in H_t$. Soit $\pi\colon\tilde{X}\to D$ la projection sur la deuxième coordonnée. La fibre de $\pi$ en $t$ est alors l'intersection schématique $X_t\coloneqq H_t\cap X$.

\[\begin{tikzcd}
X &\arrow[l] \tilde{X} \arrow[d,"\pi"] \\
& D
\end{tikzcd}\]

\begin{df}\cite[XVII, 2.2]{sga72} Avec ces notations, un pinceau d'hyperplans $D$ est dit de Lefschetz pour $X$ s'il vérifie les conditions suivantes.\begin{enumerate}
\item L'axe $A$ est transverse à $X$ \cite[XVII, 17.13.7]{ega44}. (Alors $\tilde{X}$ est l'éclaté de $X$ en $X\cap A$.)
\item Il existe une partie finie $S\subset D$ et pour chaque point fermé $s\in S$ un unique point $x_s\in \tilde X_s$ tels que $\pi$ soit lisse en-dehors des $x_s$.
\item Pour tout point fermé $s$ du fermé $S$ du point précédent, la courbe $X_s$ est nodale avec pour unique point singulier $x_s$.
\end{enumerate}
\end{df}

Si $D$ est un pinceau de Lefschetz pour $X$, le morphisme $\pi \colon \tilde X\to D\xrightarrow{\sim}\PP^1$ est propre, plat et admet une section. Sa fibre générique est lisse. De plus, le morphisme canonique $\OO_{\PP^1}\to \pi_\star\OO_{\tilde X}$ est un isomorphisme \cite[V, Th. 3.1]{milneEC}. \'{E}tant donné un point fermé $x$ de l'intersection de $X$ avec l'axe du pinceau, le morphisme $D\to \tilde{X}$, $t\mapsto (x,t)$ est une section de $\tilde{X}\to D$. Ceci revient à identifier $D$ avec l'une des composantes du diviseur exceptionnel de l'éclatement $\tilde{X}\to X$. Si $X$ est une surface alors, quitte à la plonger dans $\PP^{{n+3\choose 3}}$ via le plongement de Veronese de degré 3, un pinceau d'hyperplans général pour $X$ est un pinceau de Lefschetz et les fibres du morphisme $X\to\PP^1$ correspondant sont irréductibles \cite[p178, Conclusion]{igusajac}.\\

Au vu de ces résultats, la façon la plus rapide d'obtenir un pinceau de Lefschetz est probabiliste : tirer au hasard une droite de $\PP^3$, puis vérifier si elle est l'axe d'un pinceau en vérifiant les trois points de la définition, qui sont tous les trois algorithmiquement testables.

\section{Trivialisation des images directes dérivées}

\begin{prop} Soit $\pi\colon X\to S$ un morphisme propre lisse de schémas intègres normaux de type fini sur $k$. Soit $\F$ un faisceau lisse de $\Lambda$-modules sur $X$. Soient $\eta$ le point générique et $\bareta$ un point générique géométrique de $S$. Soit $i\in \mathbb{N}$. Soit $\eta'\to\eta$ un revêtement galoisien minimal par lequel se factorise l'action de $\Gal(\bareta|\eta)$ sur $\HH^i(X_\bareta,\F)$. Soit $S'$ la normalisation de $S$ dans $\eta'$. Alors le faisceau $\R^i\pi_\star\F$ est lisse, et $S'\to S$ est un revêtement étale galoisien minimal qui le trivialise.
\begin{proof} Notons $\G$ le faisceau $\R^i\pi_\star\F$.
Le théorème \ref{th:ehresmann} de changement de base propre-lisse assure que $\G$ est lisse. Le théorème \ref{th:chbp} de changement de base propre affirme que, pour tout point géométrique $\bar s$ de $S$, le morphisme canonique
\[ \G_{\bar s}\to \HH^i(X_{\bar s},\F) \]
est un isomorphisme.  Alors \[ \HH^0(\eta',\G_\bareta)=\HH^0(\Gal(\bareta|\eta'),\G_{\bareta})=\HH^0(\Gal(\bareta|\eta'),\HH^i(X_{\bareta},\F))=\HH^i(X_\bareta,\F)=\G_\bareta \]
et $\eta'\to\eta$ trivialise le faisceau $\eta^\star \G$. Notons $\mathfrak{S}$ l'image de $\pi_1 S$ dans $\Aut_\Lambda(\G_\bareta)$ : c'est le groupe de monodromie.
Montrons que $\G|_{S'}$ est constant, c'est-à-dire que le morphisme $\pi_1S'\to\mathfrak{S}$ est trivial.  Le diagramme commutatif
\[
\begin{tikzcd}
S'\arrow[d,"f"]& \arrow[l,"g'"] \eta'\arrow[d] \\
S &\arrow[l,"g"]\eta
\end{tikzcd}
\]
fournit le diagramme commutatif suivant :
\[
\begin{tikzcd}
\Gal(\bareta|\eta') \arrow[d] \arrow[r, two heads] & \pi_1S' \arrow[dr]\arrow[d]&\\
\Gal(\bareta|\eta) \arrow[r, two heads] & \pi_1S \arrow[r, two heads] & \mathfrak{S}
\end{tikzcd}
\]
Comme $g'^\star f^\star\G$ est constant, le morphisme $\Gal(\bareta|\eta')\to \mathfrak{S}$ est trivial, et il en est de même du morphisme $\pi_1S'\to\mathfrak{S}$. Il reste à montrer que le morphisme génériquement étale $S'\to S$ est étale. Soit $T\to S$ un revêtement étale galoisien minimal de $S$ trivialisant $\F$. Soit $\eta_T$ son point générique. Alors \[\Gal(\eta_T|\eta)=\frac{\Gal(\bareta|\eta)}{\ker(\Gal(\bareta|\eta)\to\mathfrak{S})}=\Gal(\eta'|\eta)\] et le théorème \ref{th:revcourb} assure que $S'\simeq T$.
\end{proof}
\end{prop}

\section{Calcul de la cohomologie de $\mu_n$ sur une surface}\label{sec:cohsurf}

Soit $X$ une surface intègre projective lisse sur $k$. Soit $D$ un pinceau de Lefschetz pour $X$. Notons $\pi\colon\tilde X\to\PP^1$ le morphisme correspondant défini dans la section \ref{sec:lef}. Soit $\F$ un faisceau constructible sur $X$, de tiré en arrière $\tilde \F$ sur $\tX$. Alors $\R\pi_\star\tF\in \DD^b_c(\PP^1,\Lambda)$ et 
\[ \RG(\tX,\tF)=\RG(\PP^1,\R\pi_\star\tF).\]
Soit $U$ un ouvert de $\PP^1$ au-dessus duquel $\pi$ est lisse et $\tF$ est lisse. La proposition précédente montre comment construire, pour tout entier naturel $i$, un revêtement galoisien de $U$ trivialisant $(\R^1\pi_\star\tF)|_U$. Soit $\bareta$ un point générique géométrique de $\PP^1$, d'image le point générique $\eta=\Spec k(t)$. Il suffit de calculer $\HH^1(\tX_\bareta,\F)$ avec son action de $\Gal(\bareta|\eta)$ : cette dernière se factorise par l'action d'un quotient $\Gal(\eta'|\eta)$, et la normalisation de $U$ dans $\eta'$ convient. Le degré de $\eta'\to \eta$ est l'ordre du groupe de monodromie $\mathfrak{S}\subseteq \Aut_\Lambda(\HH^i(\tX_\bareta,\tF))$. Lorsque $\tF=\mu_n$, il est borné par $n^{O(g^2)}$ où $g$ est le genre de $\tX_\bareta$.  \\

Le calcul des $\HH^i(\tX,\mu_n)$ est traité dans \cite[V, §3]{milneEC}. Nous allons le résumer en insistant sur la description explicite des objets et des morphismes concernés. Tout d'abord, la connexité des fibres de $\pi$ assure que $\pi_\star\mu_n=\mu_n$. De plus, le morphisme $\ZZ\to \R^2\pi_\star\mu_n$ par lequel se factorise $\R^1\pi_\star\GG_m\to \R^2\pi_\star\mu_n$ induit un isomorphisme $\Lambda\to \R^2\pi_\star\mu_n$. Pour tout point fermé $\bar s\in \PP^1$, le morphisme de spécialisation $\HH^1(\tX_\bareta,\mu_n)\to \HH^1(\tX_{\bar s},\mu_n)$ associe à la classe d'un diviseur $D\in \Div^0(X_\bareta)$ la classe de l'intersection de son adhérence dans $\tX$ avec la fibre $\tX_{\bar s}$. C'est un isomorphisme lorsque $\tX_{\bar s}$ est lisse.

\subsection{Morphismes de bord de la suite spectrale de Leray}\label{subsec:edge}

Le faisceau $\R^1\pi_\star\GG_m$ est le foncteur de Picard relatif \[ \Pic_{\tX/\PP^1}\colon \Pic(T\times_{\PP^1}\tX)/\pi_T^\star\Pic (T).\] C'est le faisceautisé du foncteur de Picard absolu $(T\to\PP^1)\mapsto \Pic(T\times_{\PP^1}\tX)$.
En particulier, le morphisme $\HH^1(\tX,\GG_m)\to \HH^0(\PP^1,\R^1\pi_\star\GG_m)$ obtenu par faisceautisation n'est autre que le quotient $\Pic (\tX)\to \Pic(\tX)/\pi^\star \Pic(\PP^1)=\colon \Pic (\tX/\PP^1)$. Comme $\pi \colon \tX\to\PP^1$ admet une section $\alpha$, la suite exacte \[ 0 \longrightarrow \Pic(\PP^1) \xrightarrow{\pi^\star}\Pic(\tX)\longrightarrow \Pic (\tX/\PP^1)\longrightarrow 0\]
est scindée, et il y a un isomorphisme
\[ \begin{array}{rcl}
\Pic(\PP^1)\oplus \Pic (\tX/\PP^1)& \overset{\sim}{\longrightarrow}& \Pic(\tX)\\
(D,[D'])&\longmapsto& D'+\pi^\star (D-\alpha^\star D').
\end{array}\]

Enfin, la platitude de $\pi$ assure \cite[I, Th. 4.2]{milneAV} que pour tout $T\to \PP^1$ et tout fibré en droites $\mathcal{L}$ sur $T$, le degré des fibres $\mathcal{L}_{T_s}$, $s\in\PP^1$ est indépendant de $s$ ; il est encore le même après un changement de base $T'\to T$. De plus, si $\mathcal{L}$ est obtenu comme tiré en arrière d'un fibré en droites sur $\PP^1$ alors $\deg\mathcal{L}_s=0$ pour tout $s\in \PP^1$. Ceci permet de définir un morphisme \[ \deg \colon \Pic_{\tX/\PP^1}\to \ZZ \]
de faisceaux sur $\PP^1$.
Le faisceau $\R^i\pi_\star \mu_n$ est le faisceau sur $\PP^1$ associé au préfaisceau \[ (T\to\PP^1)\mapsto \HH^i(T\times_{\PP^1}\tX,\mu_n).\] 
Il y a donc un morphisme canonique entre les espaces de sections globales \[ \HH^i(\tX,\mu_n)\to \HH^0(\PP^1,\R^i\pi_\star\mu_n) \] qui est le morphisme de bord correspondant de la suite spectrale de Leray
\[ E_2^{ij}=\HH^i(\PP^1,\R^j\pi_\star\mu_n)\Longrightarrow \HH^{i+j}(\tX,\mu_n).\]
Ces morphismes de bord sont décrits en général dans \cite[0, §12.2.5]{ega31}.
Les cas particuliers $i=1,2$ nous intéressent ici. Pour $i=1$,  $\R^1\pi_\star\mu_n$ est le noyau de la multiplication par $n$ sur $\R^1\pi_\star\GG_m$.
Le morphisme $\HH^1(\tX,\mu_n)\to \HH^0(\PP^1,\R^1\pi_\star \mu_n)$ est simplement le quotient $\Pic (\tX)[n]\to \Pic (\tX/\PP^1)[n]$, qui s'insère dans le diagramme commutatif à lignes exactes :
\[
\begin{tikzcd}
0 \arrow[r] & \Pic (\tX)[n]\arrow[r]\arrow[d] & \Pic (\tX) \arrow[d]\arrow[r,"n"] & \Pic(\tX) \arrow[d] \arrow[r] & 0\\
0 \arrow[r] & \Pic(\tX/\PP^1)[n]\arrow[r] & \Pic(\tX/\PP^1)\arrow[r,"n"] & \Pic(\tX/\PP^1)\arrow[r] & 0
\end{tikzcd}
\]
Pour $i=2$, le morphisme $\HH^2(\tX,\mu_n)\to \HH^0(\PP^1,\R^2\pi_\star \mu_n)$ s'insère dans le diagramme commutatif : 
\[
\begin{tikzcd}
\Pic(\tX)\arrow[r]\arrow[d] & \Pic (\tX/\PP^1) \arrow[d]\arrow[r,"\deg"] & \ZZ \arrow[d]\\
\HH^2(\tX,\mu_n)\arrow[r] & \HH^0(\PP^1,\R^2\pi_\star\mu_n)\arrow[r,"\sim"] & \Lambda
\end{tikzcd}
\]
En particulier, il est possible d'en construire une section : il suffit pour cela d'envoyer $1\in \Lambda$ sur l'image par $\HH^1(\tX,\GG_m)\to \HH^2(\tX,\mu_n)$ d'un diviseur de degré 1. Par exemple, le diviseur $E$, image de $\PP^1$ par la section $\alpha$ choisie au morphisme $\pi$, est de degré 1. En effet, par définition du degré, pour n'importe quel point fermé $\bar s\in \PP^1$, \[ \deg(E)=\deg (E\cap \tX_{\bar s})=1.\]
D'autre part, le morphisme \[ \HH^i(\PP^1,\pi_\star \mu_n)\to \HH^i(\tX,\mu_n)\]
s'obtient de la façon suivante. C'est la composée du morphisme $\HH^i(\PP^1,\pi_\star\mu_n)\to \HH^i(\tX,\pi^\star \pi_\star\mu_n)$, déduit du morphisme canonique $\HH^0(\PP^1,-)\to \HH^0(\tX,\pi^\star -)$, avec le morphisme $\HH^i(\tX,\pi^\star\pi_\star\mu_n)\to \HH^i(\tX,\mu_n)$ produit par l'adjonction $\pi^\star \dashv\pi_\star$. Dans notre cas, comme $\pi_\star\mu_n=\mu_n$, tout ceci est bien plus simple : la flèche \[ \HH^i(\PP^1,\mu_n)\to \HH^i(\tX,\mu_n)\]
est le morphisme $\pi^\star$ obtenu par fonctorialité de $\HH^i$. En particulier, pour $i=1$, c'est le tiré en arrière des diviseurs.

\subsection{Calcul des $\HH^i(\tX,\mu_n)$}
Rappelons qu'il y a des isomorphismes canoniques $\R^0\pi_\star\mu_n=\mu_n$ et $\R^2\pi_\star\mu_n=\Lambda$.

\begin{theorem}\cite[V, Th. 3.22]{milneEC}\label{th:cohsurf} Notons $\F\coloneqq \R^1\pi_\star \mu_n$. Les groupes de cohomologie de $\tX$ à valeur dans $\mu_{n}$ sont les suivants.
\[ \begin{array}{rcl}
\HH^0(\tX,\mu_n)&=& \mu_n(k) \\
\HH^1(\tX,\mu_n)&= & \HH^0(\PP^1,\F) \\
\HH^2(\tX,\mu_n)&= & \HH^1(\PP^1,\F)\oplus \HH^2(\PP^1,\mu_n) \oplus \HH^0(\PP^1,\Lambda)\\
\HH^3(\tX,\mu_n)&= & \HH^2(\PP^1,\F) \\
\HH^4(\tX,\mu_n)&= & \HH^2(\PP^1,\R^2\pi_\star\mu_n)=\mu_n(k)^\vee
\end{array}
\]
et $\HH^i(\tX,\mu_n)=0$ pour $i\geqslant 5$.
\end{theorem}
\begin{proof}
La deuxième page de la suite spectrale de Leray 
\[ E_2^{ij}=\HH^i(\PP^1,\R^j\pi_\star\mu_n)\Longrightarrow \HH^{p+q}(\tX,\mu_n) \]
s'écrit :
\[\begin{tikzpicture}
  \matrix (m) [matrix of math nodes,
    nodes in empty cells,nodes={minimum width=5ex,
    minimum height=5ex,outer sep=-5pt},
    column sep=1ex,row sep=1ex]{
    			&			   &			  & 	\\
          2     &  \Lambda    & 0 &   \mu_n(k)^\vee   \\
          1     &  \HH^0(\PP^1,\F) 		   & \HH^1(\PP^1,\F) & \HH^2(\PP^1,\F)\\
          0     &  \mu_n(k) & 0 &  \Lambda   &\\
    \quad\strut &   0  &  1  &  2  & \strut \\};
\draw[thick] (m-1-1.east) -- (m-5-1.east) ;
\draw[thick] (m-5-1.north) -- (m-5-5.north) ;
\end{tikzpicture}\]

Notons $M=\ker(\phi \colon \HH^2(\tX,\mu_n)\to \HH^0(\PP^1,\R^2\pi_\star\mu_n))$. La suite exacte de bas degré de cette suite spectrale de Leray s'écrit :
\[ 
0\longrightarrow \HH^1(\tX,\mu_n)\longrightarrow \HH^0(\PP^1,\F)\longrightarrow \HH^2(\PP^1,\mu_n)\overset{\pi^\star}{\longrightarrow} M\longrightarrow \HH^1(\PP^1,\F)\longrightarrow 0.
\]
Comme $\pi$ a une section, $\pi^\star$ a une rétraction, et la flèche $\HH^0(\PP^1,\F)\to \HH^2(\PP^1,\mu_n)$ est nulle. Ceci signifie d'une part que $\HH^1(\tX,\mu_n)\to \HH^0(\PP^1,\F)$ est un isomorphisme, et d'autre part qu'il y a une suite exacte scindée  \[ 0\longrightarrow \HH^2(\PP^1,\mu_n)\overset{\pi^\star}{\longrightarrow} M\longrightarrow \HH^1(\PP^1,\F)\longrightarrow 0\]
qui fournit un isomorphisme
\[ M\xrightarrow{\sim} \HH^1(\PP^1,\F)\oplus \HH^2(\PP^1,\mu_n).\]
De plus, l'image de $\HH^2(\tX,\mu_n)\xrightarrow{\phi} \HH^0(\PP^1,\Lambda)=\Lambda$ est $E_3^{02}=\ker(\psi \colon \HH^0(\PP^1,\R^2\pi_\star\mu_n)\to \HH^2(\PP^1,\F))$. Le conoyau de $\psi$ est $E_3^{21}=\HH^3(\tX,\mu_n)$. Il y a donc une suite exacte : 
\[ 
0\longrightarrow M \longrightarrow \HH^2(\tX,\mu_n) \overset{\phi}{\longrightarrow} \HH^0(\PP^1,\R^2\pi_\star\mu_n)\overset{\psi}{\longrightarrow}\HH^2(\PP^1,\F)\to \HH^3(\tX,\mu_n)\longrightarrow 0.
\]
La flèche $\phi$ a elle aussi une section, décrite dans la remarque \ref{subsec:edge}. Par conséquent, $\psi$ est nulle.
Il en découle d'une part que $\HH^2(\PP^1,\F)\to \HH^3(\tX,\mu_n)$ est un isomorphisme, et d'autre part qu'il y a un isomorphisme \[ \HH^2(\tX,\mu_n)\xrightarrow{\sim} M\oplus \HH^0(\PP^1,\R^2\pi_\star\mu_n). \]
Enfin, en degré 4, la suite spectrale a convergé à la deuxième page et $\HH^4(\tX,\mu_n)=\mu_n(k)^\vee$.
\end{proof}

\subsection{Cohomologie de l'éclatement}

Nous décrivons dans cette section le lien entre la cohomologie de $X$ et celle de $\tX$. Ces résultats sont bien connus et détaillés dans \cite[XVIII, §4]{sga72}. Soit $\Delta$ l'intersection de $X$ avec l'axe du pinceau, c'est-à-dire le centre de l'éclatement $\tX\to X$. Considérons le diagramme cartésien :
\[
\begin{tikzcd}
\tilde\Delta \arrow[r,"\tilde i"] \arrow[d,"g"] & \tX \arrow[d,"f"] \\
\Delta\arrow[r,"i"] & X
\end{tikzcd}
\]
Alors $\tilde\Delta\to\Delta$ est un fibré projectif, et le cup-produit par la classe de $\OO_{\tilde\Delta}(1)$ dans $\HH^2(\tilde\Delta,\mu_n)$ définit pour tout entier naturel $i$ un morphisme surjectif  
\[ \HH^i(\Delta,\Lambda)\to \HH^{i+2}(\tilde\Delta,\mu_n)\]
de noyau l'image de \[ g^\star \colon \HH^i(\Delta,\mu_n)\to \HH^i(\tilde\Delta,\mu_n).\]
Le morphisme $\HH^i(X,\mu_n)\to \HH^i(\tX,\mu_n)$ est injectif, car il admet un inverse à gauche : le morphisme de Gysin \[ f_\star\colon \HH^i(\tX,\mu_n)\to \HH^i(X,\mu_n).\]
De plus, pour tout entier naturel $i$, le morphisme $\tilde i^\star\colon \HH^i(\tX,\mu_n)\to \HH^i(\tilde\Delta,\mu_n)$ induit un isomorphisme \[ \frac{\HH^i(\tX,\mu_n)}{f^\star \HH^i(X,\mu_n)}\xrightarrow{\sim}\frac{\HH^i(\tilde\Delta,\mu_n)}{g^\star \HH^i(\Delta,\mu_n)}.\]
Le diagramme commutatif à lignes exactes
\[
\begin{tikzcd}
0 \arrow[r] & \HH^i(X,\mu_n) \arrow[r,"f^\star"]\arrow[d] & \HH^i(\tX,\mu_n)\arrow[r]\arrow[d,"\tilde i^\star"] & \HH^i_\Delta(X,\mu_n) \arrow[d,"\sim"{anchor=north, rotate=90}]\arrow[r] & 0 \\
0 \arrow[r] & \HH^i(\Delta,\mu_n) \arrow[r,"g^\star"] & \HH^i(\tilde\Delta,\mu_n)\arrow[r,"g_\star"] &\HH^{i-2}(\Delta,\Lambda) \arrow[r] & 0
\end{tikzcd}
\]
montre alors qu'il y a un isomorphisme :
\[ f_\star\oplus \tilde i^\star\colon \HH^2(\tX,\mu_n) \xrightarrow{\sim} \HH^2(X,\mu_n)\oplus \HH^2(\tilde\Delta,\mu_n)\]
d'inverse $f^\star\oplus \tilde i_\star$.
En particulier, \[ \HH^2(X,\mu_n)=\ker( \HH^2(\tX,\mu_n)\xrightarrow{\tilde i^\star}\HH^2(\tilde\Delta,\mu_n))\]
et pour tout $i\neq 2$, \[ \HH^i(X,\mu_n)=\HH^i(\tX,\mu_n).\]

\section{Algorithme et indications sur le calcul de sa complexité}

Reprenons les notations des sections précédentes ; en particulier, $\F=\R^1\pi_\star\mu_n$, et $S\subset \PP^1$ est le lieu de singularité du pinceau.
Les sections précédentes suggèrent l'algorithme suivant pour calculer $\HH^1(\tX,\mu_n)$. \begin{enumerate}
\item Calculer $\HH^1(\tX_\bareta,\mu_n)$ à l'aide de l'algorithme de Huang-Ierardi présenté dans la section \ref{sec:huang}, et en particulier une extension $K$ de $k(t)$ par laquelle se factorise l'action de $\Gal(\bareta|\eta)$ sur $\HH^1(\tX_\bareta,\mu_n)$. Soit $Y$ la normalisation de $U=\PP^1-S$ dans $K$. Le faisceau $\F$ est constructible sur $\PP^1$, lisse sur $U$, et $\F|_U$ est trivialisé par $Y\to U$.
\item Calculer les $\HH^1(\tX_{\bar{s}},\mu_n)$ pour $\bar s\in S$ à l'aide de la cohomologie de la normalisée de $\tX_{\bar s}$ et de la suite exacte \ref{prop:H1mnod}. Le calcul s'effectue comme dans le cas de la cohomologie à support décrit dans la section \ref{subsec:calculcohaff}.
\item Calculer les morphismes de spécialisation $\HH^1(\tX_\bareta,\mu_n)\to \HH^1(\tX_{\bar s},\mu_n)$ décrits au début de la section \ref{sec:cohsurf}. La description par recollement du faisceau $\F$ est maintenant complète.
\item Calculer $\RG(\PP^1,\F)$ à l'aide de l'algorithme de la section \ref{subsec:algconstr}.
\item En déduire les $\HH^i(\tX,\mu_n)$ à l'aide du théorème \ref{th:cohsurf}.
\end{enumerate}

Le calcul précis de la complexité de cet algorithme est encore à effectuer. Voici quelques indications en ce sens. Notons $U$ l'ouvert maximal de lissité de $\F$. Le revêtement minimal $V\to U$ qui le trivialise est de degré $|\HH^1(X_\bareta,\mu_n)|=n^{O(2g)}$ où $g$ est le genre de $X_\bareta$. Il est modérément ramifié au-dessus de $S=\PP^1-U$. Les fibres singulières $X_{\bar{s}}$ sont des courbes nodales de genre géométrique $g-1$. Afin d'étudier précisément la complexité de l'algorithme, il conviendrait de déterminer d'une part le genre géométrique de $Y$ (qui dépend en particulier du cardinal de $S$ et de la ramification à l'infini de $Y\to U$), et d'autre part la complexité des algorithmes calculant la $n$-torsion de la jacobienne d'une courbe définie sur $\FF_q(t)$ sans hypothèse supplémentaire.

\cleartooddpage

\chapter{Problèmes ouverts}

Pour finir, nous évoquons ici quelques questions restant sans réponse à l'issue de ce travail. Si certaines n'ont simplement pas pu être abordées par manque de temps, d'autres contiennent des difficultés réelles.

\section*{Faisceaux constructibles en dimension supérieure}

Comme esquissé à la fin du chapitre \ref{chap:3}, les représentations que nous avons données des faisceaux constructibles sur les courbes s'adaptent immédiatement à des variétés de dimension supérieure. Par contre, notre description des opérations effectuées sur ces faisceaux ne se généralise pas aussi facilement ; en particulier, il faudrait se passer du fait que le complémentaire de l'ouvert de lissité soit zéro-dimensionnel, qui a énormément facilité notre travail. 

\section*{Cohomologie des faisceaux constructibles sur les surfaces}

La question principale laissée en suspens par ce travail est le calcul de la cohomologie d'un faisceau constructible sur une surface lisse. Soit $X$ une surface intègre lisse, munie d'un morphisme propre $\pi\colon X\to\PP^1$ admettant une section. Les méthodes du chapitre \ref{chap:6}, qui utilisent une section de la fibration pour induire des scindages dans les suites exactes déduites de la suite spectrale de Leray, ne donnent des suites exactes courtes scindées que dans le cas des faisceaux constants. Afin de calculer la cohomologie d'un faisceau constructible $\F$ quelconque sur $X$, il conviendrait de calculer un complexe de faisceaux constructibles sur $\PP^1$ représentant $\R\pi_\star\F$. Une piste serait d'exploiter la forme des complexes de cohomologie $\RG(X_{\bar s},\F)$ des fibres de $\pi$ calculés dans la section \ref{subsec:cohlisG}.

\section*{Cohomologie des faisceaux sur les courbes singulières}

Nous avons montré comment calculer la cohomologie d'un faisceau constructible sur une courbe lisse ou nodale. Le calcul de la cohomologie d'une courbe $X$ avec des singularités quelconques serait un prolongement naturel de ce travail. Il y aurait essentiellement deux difficultés à surmonter : la première est la représentation des faisceaux constructibles sur une telle courbe, et la seconde est le calcul du revêtement $X_2\to X$ qui trivialise les $\Lambda$-torseurs sur $X$. Nos méthodes devraient s'adapter sans obstacle aux courbes à singularités ordinaires ; pour les autres, la construction de ce revêtement paraît en outre difficile.

\section*{Calcul de cup-produits}

La description explicite de l'accouplement de Weil pour les courbes lisses sur les corps finis a été réalisée par Bleher et Chinburg (voir section \ref{sec:bleher}). D'autre part, nous avons montré comment calculer les cup-produits \[ \HH^1\times \HH^1\xrightarrow{\cup} \HH^2\] dans la cohomologie des faisceaux lisses sur les courbes projectives lisses ou nodales sur les corps finis ou algébriquement clos. Il reste à étudier ce calcul dans le cadre plus large des faisceaux constructibles. En particulier, étant donné un faisceau lisse $\F$ sur ouvert $U$ d'une courbe projective lisse $X$, il serait intéressant d'obtenir un calcul explicite de la dualité de Poincaré 
\[ \HH^1_c(U,\F)\times \HH^1(U,\F^\vee(1))\xrightarrow{\cup} \HH^2(X,\mu_n) \]
à partir des complexes représentant $\RG(X,j_!\F)$ et $\RG(U,\F^\vee(1))$ déterminés dans le chapitre \ref{chap:5}. Notre méthode utilisant la cohomologie des groupes ne peut pas s'y appliquer immédiatement, puisque $j_!\F$ n'est pas localement constant. La difficulté de ce problème est encore difficile à évaluer.

\section*{Calcul de la $n$-torsion de la jacobienne sur les corps infinis}

Le seul algorithme dont nous disposons actuellement pour déterminer la $n$-torsion de la jacobienne d'une courbe projective lisse sur un corps infini est celui de Huang et Ierardi. Un premier travail serait d'étudier sa complexité dans le cas des corps de fonctions, ce qui paraît fastidieux mais tout à fait faisable.
S'il paraît difficile de mettre au point un algorithme plus efficace pour un corps quelconque, il serait tout de même intéressant de traiter le cas d'un corps de base fixé, par exemple $\QQ$ ou $\FF_q(t)$. Le cas de $\FF_q(t)$ revêt une importance particulière, puisqu'il est le cas de base indispensable pour le calcul de la cohomologie des surfaces sur $\FF_q$, et donc pour le comptage de points sur ces surfaces.

\section*{Six opérations}

L'objectif de calculer les 6 opérations dans $D^b_c(X,\Lambda)$, où $X$ est une courbe intègre lisse sur $k$, paraît pour l'instant hors de portée. Toutefois, certaines opérations semblent plus simples à calculer. Par exemple, pour un morphisme $f$ entre courbes lisses, les foncteurs $\R^i f_\star$ se déterminent aisément comme décrit dans la section \ref{subsec:pushforward}, et le calcul explicite du foncteur $\R f_\star$ paraît accessible. Le calcul de $\R f_!$ devrait s'en déduire. Nous n'avons pas eu le temps de nous intéresser au calcul des foncteurs $f^!$, $\R\underline{\Hom}$ et $\otimes^L$, qui paraît difficile.

\appendix

\renewcommand{\thedf}{\Alph{chapter}.\arabic{section}.\arabic{df}}
\renewcommand{\therk}{\Alph{chapter}.\arabic{section}.\arabic{rk}}
\renewcommand{\thetheorem}{\Alph{chapter}.\arabic{section}.\arabic{theorem}}
\renewcommand{\theprop}{\Alph{chapter}.\arabic{section}.\arabic{prop}}
\renewcommand{\thelem}{\Alph{chapter}.\arabic{section}.\arabic{lem}}
\renewcommand{\theex}{\Alph{chapter}.\arabic{section}.\arabic{ex}}
\renewcommand{\thecor}{\Alph{chapter}.\arabic{section}.\arabic{cor}}

\addcontentsline{toc}{chapter}{Annexes}

\cleartooddpage

\chapter{Corps calculables}\label{chap:A1}

\section{Corps calculables et complexité}\label{sec:compl}

\paragraph{Calculabilité et classes de complexité} Par fonction calculable, nous entendons toujours une fonction récursive au sens de \cite[Def. I.1.7]{odifreddi}, c'est-à-dire une fonction calculable par une machine de Turing \cite[Th. I.4.3]{odifreddi}. Le mot algorithme désignera une machine de Turing qui s'arrête pour toute entrée.\\
La classe des fonctions primitivement récursives \cite[Def. I.1.6]{odifreddi} est la plus petite classe contenant la fonction nulle, la fonction $n\mapsto n+1$, les projections $\NN^r\to\NN$, et stable par composition et récursion. De façon informelle, les algorithmes primitivement récursifs sont ceux qui s'écrivent avec une successions de boucles "for" : ils ne font pas usage de recherches non bornées.
La classe des fonctions élémentaires \cite[Def. VIII.7.1]{odifreddi2} est la plus petite classe contenant la fonction nulle, la fonction $n\mapsto n+1$, la soustraction $(x,y)\mapsto \max(0,x-y)$ et stable par composition, somme bornée et produit borné. Toute fonction élémentaire est primitivement récursive. Une fonction est élémentaire (resp. primitivement récursive) si et seulement si elle est calculable en un nombre d'opérations donné par une fonction élémentaire (resp. primitivement récursive) \cite[Th. VIII.7.6, VIII.8.8]{odifreddi2}. Intuitivement, la différence principale entre ces deux classes est la suivante : les fonctions exponentielles à nombre d'étages borné indépendamment des entrées sont élémentaires, mais la tétration $(n,x)\mapsto x\uparrow\uparrow n$ ne l'est pas.

\paragraph{Corps calculables} Nous supposons que tous les corps rencontrés sont calculables et munis d'un algorithme de factorisation. Cela signifie que l'on dispose d'une représentation des éléments de $k$, d'algorithmes calculant la somme, l'opposé, le produit et l'inverse d'éléments dans $k$, ainsi que d'un algorithme calculant, étant donné $f\in k[t]$, la décomposition de $f$ en produit de facteurs irréductibles.  Si $k$ est calculable et muni d'un algorithme de factorisation, il en est de même de tout corps de fonctions rationnelles à coefficients dans $k$ et de toute extension finie séparable de $k$ \cite[Lem. 19.2.2]{fried-jarden}. Si $k$ est calculable de caractéristique $p>0$ et est muni d'une $p$-base explicite (c'est-à-dire d'une $k^p$-base de $k$) et d'un algorithme de factorisation, il en 
est de même pour toute extension finie de $k$ \cite[Prop. 12.5]{mo}. Comme les corps $\QQ$ et $\FF_p$ vérifient ces hypothèses (voir section \ref{subsec:facpol} pour la factorisation), tous les corps globaux sont calculables et disposent d'un algorithme de factorisation. \\

\paragraph{Complexité et opérations élémentaires} Nous indiquerons les complexités en termes soit d'opérations binaires, soit d'opérations dans $\ZZ$ en mentionnant la taille des entiers utilisés, soit d'opérations dans un anneau fixé. L'addition de deux entiers de longueur binaire $n$ nécessite $O(n)$ opérations binaires. L'algorithme de multiplication de Harvey et van Der Hoeven, basé sur la transformée de Fourier discrète, multiplie deux entiers de longueur binaire $n$ en $O(n\log n)$ opérations binaires \cite[Th. 1.1]{harvey_mult}.
Une addition, multiplication ou inversion d'un élément de $\ZZ/d\ZZ$ requièrent $O(\log^2d)$ opérations dans $\ZZ$. Une telle opération dans un quotient $k[t]/(f)$ nécessite $O(\log^2\deg f)$ opérations dans le corps $k$. \'{E}tant donné un entier $m$ et un élément $x$ d'un anneau $A$, le calcul de $x^m$ nécessite $O(\log m)$ multiplications dans $A$.

\paragraph{Algorithmes probabilistes} Nous avons parfois recours à des algorithmes probabilistes. Il en existe deux grandes familles : \begin{enumerate}
\item Les algorithmes de type Las Vegas renvoient toujours une réponse correcte, quitte à effectuer si nécessaire un grand nombre de tirages aléatoires ; nous indiquerons toujours leur complexité moyenne.
\item Les algorithmes de type Monte-Carlo font un nombre fixe de tirages aléatoires, et renvoient une réponse qui est correcte avec une certaine probabilité ; nous indiquerons toujours leur complexité dans le pire cas ainsi que la probabilité que la valeur renvoyée soit correcte.
\end{enumerate}

\subsection{Algèbre linéaire}
\subsubsection{Sur un corps}\label{subsubsec:alglincorps}

\begin{df} La constante de l'algèbre linéaire, notée $\omega$, est la borne inférieure de l'ensemble des réels $\tau$ tels que pour tout anneau $A$, 
il existe un algorithme de multiplication de deux matrices de $\Mat_{n\times n}(A)$ nécessitant $O(n^\tau)$ opérations dans $A$.
\end{df}

Soit $k$ un corps. Soit $M\in \Mat_{n\times n}(k)$. Les opérations suivantes s'effectuent encore avec la même complexité que la multiplication dans $M_{n\times n}(k)$ : \begin{itemize}[label=$\bullet$]
\item calculer le déterminant de $M$ \cite[Th. 16.7]{buergisser_complexity} ;
\item calculer l'inverse de $M$ si elle est inversible \cite[Prop. 16.6]{buergisser_complexity} ;
\item échelonner $M$ \cite[Prop. 16.10]{buergisser_complexity}.
\end{itemize}

Il est clair que $\omega \geqslant 2$, puisque la matrice à calculer a $n^2$ coefficients. D'autre part, l'algorithme de multiplication naïf nécessite $O(n^3)$ opérations dans $k$ ; nous utilisons cette majoration dans le manuscrit. Les algorithmes les plus efficaces sont des raffinements de l'algorithme de Coppersmith-Winograd présenté dans \cite{coppersmith_winograd}, et il est connu que $\omega\in [2,2.3729[$ \cite{williams_mult}. Les mêmes bornes (pour la multiplication et l'échelonnement) s'appliquent à des matrices rectangulaires $M\in\Mat_{a\times b}(k)$ en prenant $n=\max(a,b)$ : il suffit d'ajouter des lignes/colonnes nulles pour se ramener au cas d'une matrice carrée.

\begin{lem}\label{rk:complinkx}\cite[§4.3.1, 4.3.2]{zisopoulos} Soient $k$ un corps et $k(t)$ le corps des fonctions rationnelles sur $k$. Soit $M\in \Mat_{n\times n}(k(x))$ une matrice dont les entrées sont des quotients de deux polynômes de degrés au plus $d$. Alors le calcul de toutes les opérations citées ci-dessus se fait en $\tilde O(n^4d)$ opérations dans $k$, où la notation $\tilde O(f(n))$ signifie $O(f(n) \log^\beta f(n))$ pour un $\beta >0$.
\end{lem}

\subsubsection{Sur $\ZZ/d\ZZ$}

Soit $\Lambda$ un anneau. La complexité $O(n^\omega)$ de la multiplication des matrices est la même que dans le cas des corps. Il n'est cependant pas évident que cette complexité soit encore celle du calcul de noyaux de morphismes de $A$-modules libres ou de la résolution de systèmes linéaires sur $A$. Dans le cas où tous les idéaux de $\Lambda$ sont principaux, la forme normale de Smith permet de réaliser cette opération ; il existe des algorithmes pour la calculer dans le cas où $\Lambda=\ZZ$ \cite[§4]{storjohann} ou $\ZZ/d\ZZ$ \cite[§3]{storjohann}. Cependant, ces algorithmes calculent seulement la forme normale $N$ d'une matrice $M$, et le calcul des matrices de transformation $P,Q$ telles que $M=PNQ$ est plus coûteux. 
Lorsque $\Lambda=\ZZ/d\ZZ$, la forme normale de Howell se prête particulièrement bien à la tâche. \`{A} chaque matrice $M\in \Mat_{n\times m}(\Lambda)$, on peut associer une unique matrice $H(M)\in \Mat_{n\times m}(\Lambda)$ échelonnée vérifiant certaines conditions décrites dans \cite[§3]{storjohann_howell} et une unique matrice $P\in \GL_n(\Lambda)$ telle que $PM=H(M)$. En particulier, $H(M)=H(N)$ si et seulement si $\ker(M)=\ker(N)$.

\begin{prop}\cite[Th. 4]{storjohann_howell}
Soient $n,m,d$ des entiers naturels non nuls. Notons $\Lambda=\ZZ/d\ZZ$. Il existe un algorithme qui, étant donné $M\in \Mat_{n\times m}(\Lambda)$, calcule sa forme normale de Howell $H(M)$ et une matrice $P\in \GL_n(\Lambda)$ telle que $PM=H(M)$ en $O(\max(n,m)^\omega)$ opérations dans $\Lambda$.
\end{prop}

Cet algorithme permet en particulier de calculer le noyau d'une matrice, une famille génératrice de la réunion ou de l'intersection de sous-modules \cite[§5, Tasks 1-3]{storjohann_howell} avec cette complexité. Un $\Lambda$-module engendré par des éléments $v_1,\dots,v_n$ vérifiant les relations linéaires $a_1,\dots,a_m$ sera décrit par la matrice $M\in \Mat_{n\times m}(\Lambda)$ dont il est le conoyau. La description des morphismes et le calcul des noyaux et conoyaux se fait exactement comme dans \cite[13.2]{mo}. Remarquons également que comme tout idéal de $\Lambda$ est principal, un sous-module de $\Lambda^n$ est toujours engendré par une famille d'au plus $n$ éléments qui se détermine simplement, ce qui fait que le nombre de relations entre les générateurs d'un $\Lambda$-module n'entre pas en compte dans le calcul de la complexité.

\section{Polynômes}

\subsection{Factorisation des polynômes univariés}\label{subsec:facpol}

\subsubsection{Sur un corps fini}

\begin{prop}\cite[Th. 1]{shoup_fact} Il existe un algorithme déterministe qui, étant donné un nombre premier $p$ et un polynôme $f\in\FF_p[t]$ de degré $d$, calcule les facteurs irréductibles de $f$ dans $\FF_p[t]$ en $O(d^{2+\epsilon} \sqrt{p}\log^2p)$ opérations dans $\FF_p$.
\end{prop}

Ceci implique encore que la factorisation d'un polynôme dans $\FF_{p^\alpha}[t]$ se calcule en un nombre d'opérations polynomial en $p$, $\alpha$ et $d$ \cite[14.40]{vzg}. Le recours aux algorithmes probabilistes permet d'obtenir une meilleure complexité.

\begin{prop}\cite[Th. 14.32]{vzg} Il existe un algorithme probabiliste (Las Vegas) qui, étant donné une puissance $q=p^\alpha$ d'un nombre premier et un polynôme $f\in\FF_q[t]$ de degré $d$,
calcule les facteurs irréductibles de $f$ dans $\FF_{q}[t]$ en $\tilde O(d^\omega\alpha^2\log^2 p)$ opérations dans $\FF_q$, avec probabilité d'échec inférieure à $\frac{1}{2}$.
\end{prop}

\begin{rk} Afin de calculer la factorisation dans $\overline{\FF_p}[t]$ d'un polynôme $f\in \FF_{p^\alpha}[t]$ de degré $d$, il suffit de le factoriser dans une extension dont le degré est le ppcm des degrés des facteurs irréductibles de $f$, qui est majoré par la fonction de Landau $\lambda(d)\leqslant \exp(d/e)$, où $e=\exp(1)$ ; pour des majorations plus précises de la fonction de Landau, voir \cite{nicolas_landau}. La complexité de la factorisation est alors encore polynomiale en $p,\alpha,\exp(d/e)$ pour l'algorithme déterministe, et $O(d^{\omega+2}exp(2\alpha/e)\log^2p)$ pour l'algorithme probabiliste.
\end{rk}

\subsubsection{Sur un corps de nombres}

Par le lemme de Gauss, la factorisation d'un polynôme de $\QQ[t]$ se ramène à celle d'un polynôme dans $\ZZ[t]$, que l'on peut supposer primitif après avoir factorisé ses coefficients dans $\ZZ$, opération dont la complexité est sous-exponentielle en la valeur absolue des coefficients. Soit donc $f\in \ZZ[t]$ de degré $d$. Notons $\Vert f\Vert_\infty$ la plus grande des valeurs absolues des coefficients de $f$.  L'algorithme "LLL" de Lenstra, Lenstra et Lov\'{a}sz calcule les facteurs premiers de $f$ dans $\ZZ[t]$ en $O(d^6+d^5\log\Vert f\Vert_\infty)$ opérations sur des entiers de longueur binaire $O(d^3+d^2\log\Vert f\Vert_\infty)$ \cite[Th. 3.6]{LLL}. \\
Cet algorithme a été adapté par A. Lenstra au cas des corps de nombres. Soit $K=\QQ[\alpha]$ un corps de nombres ; notons $n=[K:Q]$. Soit $f\in K[t]$ un polynôme de degré $d$. Soit $D$ un entier tel que $f\in \frac{1}{D}\ZZ[\alpha][t]$. On peut écrire $f=\sum_{i=0}^{n-1}\sum_{j=0}^d a_{ij}\alpha^i x^j$. Notons alors $\Vert f\Vert_2$ la norme euclidienne du vecteur $(|a_{ij}|)_{i,j}$ et $\Vert f\Vert_\infty=\max_{i,j}|a_{ij}|$.

\begin{theorem}\cite[Th. 4.5]{lenstra_lll}
Avec ces notations, il existe un algorithme déterministe qui calcule la factorisation de $f$ dans $\QQ(\alpha)[t]$ en $O(n^6d^6+n^5d^6\log(d\Vert f\Vert_2)+n^5d^5\log(d\Vert f\Vert_\infty))$ opérations arithmétiques sur des entiers de longueur binaire $O(n^3d^2+n^2d^3\log(d\Vert g\Vert_2)+n^2d^2\log(d\Vert f\Vert_\infty))$.
\end{theorem}

\subsection{Calcul de racines $n$-ièmes}\label{subsec:racn}

Dans les corps finis, le calcul d'une racine $n$-ième d'un élément par la méthode d'Adleman-Manders-Miller est bien plus rapide que la factorisation des polynômes en général. 

\begin{prop}\cite[§6]{cao_rootf} Soit $q$ une puissance d'un nombre premier. Soit $n$ un entier naturel non nul. Il existe un algorithme déterministe qui, étant donné $x\in \FF_q$, détermine si $x$ est une puissance $n$-ième dans $\FF_q$ et, le cas échéant, calcule une racine $n$-ième de $x$ en $O(\log^4q+n\log^3q)$ opérations dans $\FF_q$.
\end{prop}

\subsection{Factorisation des polynômes multivariés}

\subsubsection{En général}

Soit $k$ un corps calculable disposant d'un algorithme de factorisation des polynômes. Alors il existe un algorithme de factorisation des éléments de $k[x_1,\dots,x_m]$, basé sur l'observation suivante \cite[§11.3]{fried-jarden}. L'application \[ \begin{array}{ccl}
   k[x_1,\dots ,x_m]&\longrightarrow& k[t] \\ \sum a_ix_1^{i_1}\cdots x_m^{i_m}&\longmapsto& \sum a_it^{i_1+i_2d+i_3d^2+\cdots+i_md^{m-1}}
\end{array}\]  définit, pour chaque entier $d$, une bijection $\kappa_d$ de l'ensemble des polynômes de $k[x_1,\dots ,x_m]$ de degré $<d$ vers l'ensemble des polynômes de $k[t]$ de degré $<d^{m}$, qui vérifie $\kappa_d(fg)=\kappa_d(f)\kappa_d(g)$ dès que $\deg(fg)<d$. Son inverse se calcule explicitement par un algorithme d'écriture des entiers en base $d$. Un polynôme $f\in k[x_1,\dots,x_m]$ de degré $<d$ est irréductible si et seulement si $\kappa_d(f)$ l'est. Afin de déterminer un facteur d'un polynôme $f$ de degré $<d$, il suffit donc de calculer les facteurs irréductibles $g_1,\dots,g_s$ de $\kappa_d(f)$. Les facteurs irréductibles stricts de $f$ (s'il y en a) se trouvent parmi les images réciproques par $\kappa_d$ des facteurs de $\kappa_d(f)$, c'est-à-dire parmi les polynômes de la forme $\kappa_d^{-1}(g_{i_1}\cdots g_{i_s})$. Il ne reste plus qu'à vérifier par division euclidienne si ces éléments divisent $f$. En particulier, si l'algorithme de factorisation dans $k[t]$ a une complexité élémentaire en le degré du polynôme et la taille de la représentation de ses coefficients, il en est de même de l'algorithme de factorisation dans $k[x_1,\dots,x_n]$.

\subsubsection{Sur les corps finis}\label{susubsec:factmult}

En se basant sur l'algorithme LLL, A.K. Lenstra a également donné un algorithme de factorisation des polynômes en plusieurs variables sur les corps finis. La version qui nous est utile concerne deux variables.

\begin{theorem}\cite[Th. 2.18]{lenstra_lll_finite} Il existe un algorithme déterministe qui, étant donné une puissance $q=p^m$ d'un nombre premier $p$ et un polynôme $f\in \FF_q[x,y]$, renvoie les facteurs irréductibles de $f$ dans $\FF_q[x,y]$ en $O(\deg_x(f)^6\deg_y(f)^2+(\deg_x(f)^3+\deg_y(f)^3)pm)$ opérations dans $\FF_q$.
\end{theorem}

De même que dans le cas univarié, il existe un algorithme probabiliste de complexité polynomiale pour factoriser des polynômes en deux variables sur un corps fini.

\begin{theorem}\cite[Cor. 4.2]{wan_factoring} Il existe un algorithme probabiliste Las Vegas qui, étant donné une puissance $q$ d'un nombre premier et un polynôme $f\in \FF_q[x,y]$ de degré total $d\leqslant \sqrt{q}$, calcule les facteurs irréductibles de $f$ dans $\FF_q[x,y]$ en $O(d^{4.89}\log^2d\log q)$ opérations, avec une probabilité d'échec inférieure à $\frac{1}{\sqrt{\pi\log d}}$.
\end{theorem}

Remarquons que cet algorithme nécessite un corps fini assez grand par rapport au degré du polynôme. Lorsque $d^2\geqslant q$ avec les notations du théorème, il faut construire une extension $\FF_Q$ de $\FF_q$ de degré supérieur à $2\log_q d$, factoriser $f$ dans $\FF_Q$ puis multiplier entre eux les facteurs conjugués sous $\Aut(\FF_Q|\FF_q)$. Comme nous le verrons dans la section \ref{subsec:polirred}, la construction de cette extension $\FF_Q/\FF_q$ n'est pas plus coûteuse que la factorisation du polynôme sur $\FF_q$.

\begin{rk} Les algorithmes précédents permettent également de factoriser des polynômes en une variable sur $\FF_q(t)$. Soit $f\in \FF_q(t)[x]$. \'{E}crivons $f=\sum_i R_i(t)x^i$, où $R_i=\frac{P_i}{Q_i}\in \FF_q(t)$. Notons $d_x=\deg_xf$, et $d_t=\max_i \{\deg_t P_i,\deg_t Q_i\}$. Notons encore $Q=\ppcm_i Q_i$ ; il est de degré inférieur à $d_xd_t$. Factoriser $f$ dans $\FF_q(t)[x]$ revient à factoriser $Q(t)f\in \FF_q[t,x]$. C'est un polynôme de degré total inférieur à $d_td_x^2$. Le polynôme $f$ peut donc être factorisé par un algorithme probabiliste en $O((d_td_x^2)^{4.89}\log^2(d_td_x^2)\log q)$ opérations.
\end{rk}

\subsubsection{Sur les corps de nombres}

Soit $\QQ(\alpha)=\QQ[t]/(F)$ une extension de $\QQ$ de degré $d$. Soit $f\in \QQ(\alpha)[x_1,\dots,x_r]$. Soit $D$ un entier tel que $f\in\frac{1}{D}\ZZ[\alpha][x_1,\dots,x_r]$. Notons $n=\max(2,\min\deg_{x_i}f)$ et $N=\prod_{i=1}^r (1+\deg_{x_i}f)$. Les polynômes $f$ et $F$ sont identifiés aux vecteurs complexes de leurs coefficients.

\begin{theorem}\cite[Th. 3.26]{lenstra_factoring_multi} Avec ces notations, il existe un algorithme déterministe qui calcule la décomposition de $f$ en produit de facteurs irréductibles dans $\QQ(\alpha)[x_1,\dots,x_r]$ en \[O(n^{r-1}(dn)^5(dn+\log (d \Vert f\Vert_\infty+d\log(d\Vert F\Vert_2))\] opérations arithmétiques sur des entiers de longueur binaire \[O(n^{r-1}(dn)^2(dn+\log (d \Vert f\Vert_\infty+d\log(d\Vert F\Vert_2)).\]
\end{theorem}

\section{Extensions de corps}

Une extension finie $L$ d'un corps calculable $K$ est représentée concrètement par une $K$-base $B=(b_1,\dots,b_n)$ de $L$ et la donnée, pour tous $i,j$, de la décomposition de $b_ib_j$ dans la base $B$. Le calcul d'un produit dans $L$ nécessite donc $n^\omega$ opérations dans $K$, celui d'une somme $n$ opérations dans $k$.
Le polynôme minimal d'un élément $x\in L$ se détermine alors en calculant les puissances $1,x,x^2...$ de $x$, et en vérifiant à chaque étape si $x^i$ appartient au sous-espace vectoriel de $L$ engendré par $1,x,\dots,x^{i-1}$. Ce calcul nécessite $O(n^{3+\omega})$ opérations dans $K$.

\subsection{Calcul de polynômes irréductibles sur les corps finis}\label{subsec:polirred}

Soit $p$ un nombre premier. Soit $d$ un entier naturel non nul. La construction d'une extension de degré $d$ de $\FF_p$ nécessite le calcul d'un polynôme irréductible de $\FF_p[t]$ de degré $d$. De même que pour la factorisation, une complexité polynomiale en $\log p$ n'est pour l'instant garantie que par des algorithmes probabilistes.

\begin{prop}\label{prop:polirred}\cite[Th. 3.2]{shoup_irred} Il existe un algorithme déterministe qui, étant donné un nombre premier $p$ et un entier $d>0$, construit un polynôme irréductible dans $\FF_p[t]$ de degré $d$ en \[O(\sqrt{p}\log^3(p)d^{3+\epsilon}+\log^2(p)d^{4+\epsilon})\] opérations dans $\FF_p$ pour tout $\epsilon >0$.
\end{prop}

\begin{prop}\cite[Th. 5.1]{shoup_fast} Il existe un algorithme probabiliste (Las Vegas) qui, étant donné une puissance $q$ d'un nombre premier et un entier $d>0$, construit un polynôme irréductible dans $\FF_q[t]$ de degré $d$ en $O((d^2\log d+d\log q)\log d\log\log d)$ opérations dans $\FF_q$ en moyenne.
\end{prop}

\subsection{Extensions normales, extensions séparables}

\subsubsection{Tester si une extension est normale ou séparable}\label{subsubsec:testgal}

Soit $L/K$ une extension de corps de degré $n$. Soient $x_1,\dots,x_d$ des générateurs de $L$ comme $K$-algèbre, et $f_1,\dots,f_d\in K[t]$ leurs polynômes minimaux. La séparabilité de $L/K$ équivaut alors à celle des $f_i$, qui se teste en calculant $\pgcd(f_i,f_i')$. La normalité de $L/K$ se teste en factorisant les polynômes $f_i$ dans $L[t]$, et en vérifiant qu'ils y ont $\deg(f_i)$ racines avec multiplicité. La complexité de tester la séparabilité de $L/K$ est donc celle de la factorisation dans $L$ de $d<n$ polynômes de $K[t]$ de degré au plus $n$.

\subsubsection{Calcul d'un élément primitif}\label{subsubsec:eltprim}

Soit $K$ un corps infini. Considérons une extension finie séparable $L=K(a,b)$ de degré $d$ de $K$. Alors les éléments $\lambda\in K$ tels que $a+\lambda b$ ne soit pas un générateur de $L$ sont les racines d'un polynôme de degré $d(d-1)$ \cite[V, Th. 4.6]{lang}. Afin de déterminer un élément primitif de $L$, il suffit donc d'énumérer au plus $d(d-1)+1$ éléments $\lambda\in K$. Pour chacun de ces éléments, le polynôme minimal de $a+\lambda b$ sur $k$ se détermine en calculant ses puissances successives $(a+\lambda b)^i$ et en vérifiant par des méthodes d'algèbre linéaire si $(1,a+\lambda b,\dots,(a+\lambda b)^{d-1})$ est une base de $L/K$. Le nombre d'opérations à effectuer dans $K$ est donc polynomial en $d$. \\

Le cas général d'une extension $L=K(a_1,\dots,a_s)$ s'en déduit par récurrence sur le nombre de variables : $K(a_1,\dots,a_s)=K(a_1,\dots,a_{s-2})(a_{s-1},a_s)$. Le nombre d'opérations à effectuer dans $K$ est alors polynomial en $d^s$, puisqu'une opération dans $k(a_1,\dots,a_{s-2})$ correspond à $O(d^{s-2})$ opérations dans $K$. 

Le raisonnement suivant, adapté de \cite[§5]{yokoyama_primitive}, permet de trouver un $s$-uplet $(\lambda_1,\dots,\lambda_s)$ tel que $\lambda_1a_1+\dots+\lambda_sa_s$ soit un élément primitif de $L/K$ dans un ensemble défini à l'avance de taille $(s-1)[L:K]$. Soit $A$ un anneau principal infini de corps des fractions $K$. Soit $L=K(a_1,\dots,a_s)$ une extension finie séparable de degré $N$. Soit $S$ un sous-ensemble de $A^s$ dont toutes les familles de $s$ éléments sont linéairement indépendantes. Si $|S|>(s-1)(N-1)$ alors il contient un $s$-uplet $\lambda$ convenable \cite[Th. 4.5]{yokoyama_primitive}. En particulier, un tel ensemble peut être construit sous la forme $S=\{ (1,\lambda,\dots,\lambda^{s-1})\mid \lambda\in B\}$, où $B$ est une partie de $A$ de cardinal au moins $(s-1)(N-1)$. Pour $A=k[t]$, il suffit de prendre suffisamment d'éléments $\alpha\in k$ et de considérer les $(1,(t-\alpha),\dots,(t-\alpha)^{s-1})$.


\subsubsection{Groupe de Galois}\label{subsubsec:gpgal}

Soit $L/K$ une extension galoisienne de corps de degré $n$. Supposons avoir déjà calculé un élément primitif $x_1$ de $L/K$, ainsi que son polynôme minimal $f\in K[t]$, et les racines $x_1,x_2,\dots,x_n$ de $f$ dans $L$. Les éléments de $\Gal(L|K)$ sont alors déterminés par l'image de $x_1$, qui est l'un des autres $x_i$. La complexité du calcul de $\Gal(L|K)$ est donc dominée par celle de la factorisation dans $L$ d'un polynôme de degré $n$ à coefficients dans $K$.

\subsection{Extensions radicielles et clôture parfaite}\label{sec:insep}

Soit $k$ un corps de caractéristique $p$. Il sera souvent utile d'effectuer des calculs dans la clôture parfaite $k^\pf$ de $k$. Le principe systématiquement adopté est le suivant, proposé dans \cite[§2.3]{steel}. Il consiste à simuler une extension $k^{p^{-r}}$ assez grande pour effectuer les calculs désirés, en élevant à la puissance $p^r$ tous les éléments de $k$ rencontrés. L'isomorphisme $k\to k^{p^{r}}$ permet ainsi de remplacer $k$ par $k^{p^{r}}$, et $k^{p^{-r}}$ par $k$. L'élévation d'un élément de $k$ à la puissance $p^r$ par exponentiation rapide nécessite $O(r\log p)$ opérations dans $k$.

\cleartooddpage

\chapter{Schémas de type fini sur un corps}\label{chap:A2}

\section{Schémas et morphismes}

\subsection{Représentation des schémas et des morphismes}\label{sec:repsch}

Soit $k$ un corps. La description explicite des schémas de type fini sur $k$ telle que présentée ci-dessous est celle de \cite[§16]{mo}. Dans le cas particulier des courbes projectives, d'autres descriptions plus adaptées seront données dans la section \ref{sec:repcourb}.

\paragraph{Schémas affines} La donnée de polynômes $f_1,\dots,f_r\in k[x_1,\dots,x_m]$ définit le schéma $X=\Spec A$ où  $A=k[x_1,\dots,x_m]/(f_1,\dots,f_r)$.  Un ouvert $U$ de $X$ est défini par des polynômes $g_1,\dots,g_s$ tels que $U=\bigcup_i D(g_i)$. Soit $X'=\Spec k[x'_1,\dots,x'_p]/(f'_1,\dots,f'_t)$ un autre schéma affine. Un morphisme $X\to X'$ est défini par ses fonctions coordonnées $\phi_1,\dots,\phi_p\in k[x_1,\dots,x_m]$. \'{E}tant donné un schéma affine $X''$ et des morphismes $\psi\colon X\to X''$ et $\psi'\colon X'\to X''$, le produit fibré $X\times_{X''}X'$ est défini par la $k$-algèbre $k[x_1,\dots,x_m,x'_1,\dots,x'_p]/(f_i,f'_j,\psi_\alpha-\psi'_\alpha)$. \\

Il est possible de détecter si le schéma affine $X$ est vide : c'est le cas si et seulement s'il existe des polynômes $a_1,\dots,a_r$ de degré inférieur à une constante (que l'on sait calculer) dépendant de $m$, $s$ et des degrés des $f_i$ tels que
\[ a_1f_1+\dots +a_rf_r=1. \]
La résolution de cette équation se ramène à elle d'un système linéaire en les coefficients des $a_i$.

\paragraph{Schémas de type fini sur $k$} Par défaut, un schéma de type fini sur $k$ est représenté comme recollement de schémas affines. Reprenons les schémas $X,U,X'$ ci-dessus. Soit $U'=\bigsqcup_i D(g'_i)$ un ouvert de $X'$.
Un morphisme $U\to U'$ est défini par des morphismes \[D(g_i)=\Spec k[x_1,\dots,x_m,x]/(f_1,\dots,f_r,xg_i-1)\to X'\] qui se factorisent par $U'$ et coïncident sur $D(g_i)\cap D(g_j)=D(g_ig_j)$.
Le fait qu'ils se factorisent par $U'$ se teste en calculant, pour chaque indice $i$, le schéma $Z(g'_i)\times_{X'} X$, et en vérifiant qu'il est vide.

\subsection{Recouvrement par des voisinages de points}\label{subsec:recprimrec}

\begin{prop} Soit $X\subset \A^n_k$ un schéma de type fini sur $k$ de dimension $d$. Supposons que l'on dispose d'un algorithme primitivement récursif calculant, pour un point $x\in X$, un voisinage de Zariski de $x$ dans $X$ vérifiant une certaine propriété $(P)$. Alors il existe un algorithme primitivement récursif recouvrant $X$ par un nombre fini d'ouverts $U_1,\dots,U_r$ vérifiant $(P)$.
\begin{proof} Quitte à appliquer la procédure ci-après à chacune de ses composantes irréductibles (voir section \ref{subsec:decprim} pour leur calcul), nous pouvons supposer $X$ irréductible.
Commençons par choisir un point $x\in X$ et de lui appliquer l'algorithme pour obtenir un premier ouvert $U^{(1)}$, dont le complémentaire $F_1$ est de dimension $\leqslant d-1$. Le nombre $r$ de composantes connexes de $F_1$ est borné en fonction du degré des équations qui le définissent. L'exécution de l'algorithme pour un point de chacune de ces composantes connexes construit des ouverts $U_1^{(2)},\dots ,U_r^{(2)}$ de $X$, et ainsi de suite. Le complémentaire de la réunion des $U_i^{(j)}$, pour $1\leqslant j\leqslant r$, est de dimension au plus $d-r$, et son nombre de composantes connexes est borné par le degré de ses équations, qui ont été explicitement calculées. Le nombre de points à considérer à l'étape $j=r+1$ est donc connu. Il y a au plus $d$ étapes à effectuer.
\end{proof}
\end{prop}

\section{Bases de Gröbner et applications}
\label{subsec:dimension}
\label{subsec:réduit}
\label{subsec:decprim}

Soit $k$ un corps. Fixons un ordre monomial sur $R\coloneqq k[x_1,\dots,x_n]$, c'est-à-dire un bon ordre sur les monômes de $R$ compatible à la multiplication \cite[Def. 1.1.1]{derksen}. Soit $I$ un idéal de $R$. Une famille $(f_1,\dots,f_r)$ d'éléments de $I$ est appelée base de Gröbner si les monômes dominants des $f_i$ engendrent l'idéal de $R$ engendré par les monômes dominants de tous les éléments de $I$. Dans ce cas, l'algorithme de division multivariée permet de décider de l'appartenance à $I$ : un élément $f\in R$ appartient à $I$ si et seulement si le reste de la division de $f$ par $(f_1,\dots,f_r)$ est nul. Une base de Gröbner $(f_1,\dots,f_r)$ de $I$ est dite réduite si les $f_i$ sont unitaires et si aucun monôme d'un $f_i$ n'appartient à l'idéal engendré par les monômes dominants des $f_j\neq i$. Tout idéal admet une unique base de Gröbner réduite. 
Si l'idéal $I$ est défini par des générateurs $g_1,\dots,g_s$ de degré maximal $d$, le degré des éléments de la base de Gröbner réduite de $I$ est $O(d^{2^n})$ \cite[Th. 8.2]{dube}.
\begin{prop} \cite[Prop. 1]{bardet_faugereF5} Avec ces notations, il existe un algorithme calculant une base de Gröbner réduite de $I$ en \[ O\left(\frac{s}{n!}2^{n\omega}d^{(1+n\omega)2^n}\right) \]
opérations dans $k$.
\end{prop}
En particulier, les bases de Gröbner permettent de résoudre des systèmes par élimination. Plus précisément, pour un ordre adapté à l'élimination des variables $x_1,\dots,x_{j-1}$, si $G$ est une base de Gröbner de l'idéal $I$ alors $G\cap k[x_j,\dots,x_n]$ est une base de Gröbner de $I\cap k[x_j,\dots,x_n]$. 
Un excellent résumé des applications principales des bases de Gröbner se trouve dans \cite[Ch. 1]{derksen} ; une référence très détaillée est \cite{bw}. Supposons que le corps $k$ est calculable, muni d'une $p$-base explicite et dispose d'un algorithme de factorisation des polynômes de complexité élémentaire en le degré du polynôme et la taille de la représentation informatique de ses coefficients. Il existe alors des algorithmes de complexité élémentaire en $n,d,s$ qui utilisent les bases de Gröbner pour calculer :
\begin{enumerate}
\item le noyau d'un morphisme de $k$-algèbres de type fini ;
\item la fonction de Hilbert de $I$ \cite[Alg. 2.6]{bayer} ;
\item une normalisation de Noether de $I$ \cite[Alg. 1.13]{dickenstein} ;
\item le radical de $I$ \cite[Th. 8.99]{bw} ;
\item la décomposition primaire de $I$ \cite[Th. 8.101]{bw}, \cite[Th. 5.1]{steel} ;
\item la normalisation de $R/I$ \cite[§15.5]{mo}.
 \end{enumerate}
Le calcul de la décomposition primaire se ramène au cas zéro-dimensionnel par récurrence, et nécessite la factorisation de polynômes sur des corps de fractions rationnelles sur $k$.
La normalisation d'une $k$-algèbre $A$ de type fini est calculée de la façon suivante. Un anneau est normal si et seulement s'il est régulier en codimension 1 et vérifie la propriété $(S_2)$ de Serre \cite[IV, Th. 11]{serre}. L'algorithme de normalisation de de Jong \cite{dejong} permet de régulariser l'anneau $A$ en codimension 1 (le nombre d'étapes est alors borné par une multiplicité calculable elle-même avec complexité élémentaire \cite[§15.4, 15.5.3]{mo}). Notons $A_1$ l'anneau obtenu. La plus petite extension de $A_1$ dans $\Frac(A_1)$ qui vérifie la propriété $(S_2)$ se calcule à l'aide d'un bidual \cite[Prop. 6.21]{vasc2}, et est la normalisation de $A$. 

\section{Construction de familles d'hyperplans}

\subsection{Hyperplans en position générale}\label{subsec:hyperplans}

\'{E}tant donné un corps $k$ ayant suffisamment de points et un entier $t$ donné, comment construire un ensemble de $t$ hyperplans de $\PP^n_k$ en position générale ?

\begin{df}\label{def:posgen}
Soit $k$ un corps. Soient $n,t$ deux entiers strictement positifs. Un ensemble $S$ de $t$ hyperplans de $\PP^n_k$ est en position générale si : \begin{itemize}[label=$\bullet$]
\item pour tout $i\in \{1,\dots,n\}$ et tous $H_1,\dots,H_i\in S$ deux à deux distincts, $\dim(H_1\cap \dots\cap H_i)=n-i$ ;
\item tout point fermé de $\PP^n_k$ appartient à au plus $n$ hyperplans de $S$.
\end{itemize}
\end{df}
Le théorème suivant permet de contrôler l'intersection de deux variétés projectives.
\begin{prop}\cite[Th. 1.7.2]{hartshorne}\label{projdimthm}
Soient $Y,Z$ deux sous-variétés de $\PP^n$ de dimensions respectives $r,s$. Alors toute composante irréductible de $Y\cap Z$ est de dimension $\geqslant r+s-n$. Lorsque $r+s-n\geqslant 0$, l'intersection est non vide.
\end{prop}

\begin{lem}
Soit $S$ un ensemble de $t$ hyperplans de $\PP^n_k$ en position générale, avec $t\geqslant n$. Soit $H$ un hyperplan de $\PP^n_k$. Soient $P_1,\dots,P_s$ les points d'intersection de toutes les familles de $n$ hyperplans de $S$. Si $H$ ne contient aucun des $P_i$ alors $S\cup \{ H\}$ est en position générale.
\begin{proof}
Soit $i\leqslant n-1$. Considérons une intersection $I=H_1\cap\dots\cap H_i$ d'éléments de $S$. Alors il existe $j\in \{1,\dots,s\}$ tel que $P_j\in I$. Par conséquent, $H$ ne contient pas $I$, qui est irréductible, donc $\dim(H\cap I)=\dim(I)-1=n-i$ d'après la proposition \ref{projdimthm}.\\
Soit $P$ un point fermé de $\PP^n_k$. Si $P$ appartient à strictement moins de $n$ hyperplans de $S$, il n'y a rien à vérifier. Si $P$ appartient à $n$ hyperplans de $S$ alors c'est l'un des $P_i$ et il n'appartient pas à $H$.
\end{proof}
\end{lem}

\'Etant donné deux entiers $t,n\geqslant 1$, voici l'algorithme qui calcule un ensemble de $t$ hyperplans de $\PP^n_k$ en position générale.\\
Si $t\leqslant n+1$, l'algorithme renvoie les hyperplans définis par les coordonnées $x_0,\dots,x_{t-1}$.\\
Supposons construit un ensemble $S=\{H_1,\dots,H_t\}$ d'hyperplans en position générale, avec $t\geqslant n+1$. Voici comment construire $H_{t+1}$ tel que $\{H_1,\dots,H_{t+1}\}$ soit en position générale. Dans un premier temps, des algorithmes d'algèbre linéaire permettent de déterminer les points d'intersection $P_1,\dots,P_s$ (avec $s\leqslant {t \choose n}\leqslant \frac{t^n}{n!}$) des familles de $n$ hyperplans de $S$. Un hyperplan $H$ défini par une forme linéaire $F=\sum_i a_ix_i\in k[x_0,\dots,x_n]$ ne contient pas un point $P_i=(p_0^i:\dots:p_n^i)$ si et seulement si le point $a=(a_0,\dots,a_n)\in k^{n+1}$ n'appartient pas au noyau de la forme linéaire $F_j=\sum_j p_j^ix_j$. Par conséquent, les points $a$ convenables se trouvent en-dehors d'une réunion de $s$ hyperplans de $\A^{n+1}_k$ qui se calculent explicitement.\\

L'astuce suivante pour trouver un point $a\in \A^{n+1}(k)$ convenable provient de \cite[Lem. 2.4]{kratsch}. \'Enumérons $t^{n+1}$ éléments $b_1,\dots,b_s\in k$. Considérons la grille $R=\{b_1,\dots,b_{t^{n+1}} \}^{n+1} \subseteq\A^{n+1}(k)$. Alors tout hyperplan de $\A^{n+1}_k$ contient au plus $t^{n^2+n}$ éléments de $R$ ; en particulier, la réunion des hyperplans définis par $F_1,\dots,F_s$ contient donc moins de $\frac{1}{n!}t^{n^2+2n}$ points de $R$. Or $R$ contient $t^{n^2+2n+1}$ éléments. Il suffit donc de prendre un point $a\in R$ qui n'est pas sur l'un des hyperplans, et de définir $H_{t+1}$ comme l'hyperplan défini par la forme linéaire $\sum_i a_ix_i$.\\

Pour chaque hyperplan $H'_j$ défini par $F_j$, le calcul de $H'_j\cap R$ se fait par recherche exhaustive sur les points de $R$ : il y en a $t^{n^2+2n+1}$. Il existe donc $\alpha >0$ tel que la complexité totale de l'algorithme soit $O(t^{\alpha n^2})$. Le résultat que nous avons montré est le suivant.

\begin{prop} Soient $n,t$ deux entiers. Soit $k$ un corps ayant au moins $t^{n+1}$ éléments. Il existe $\alpha >0$ et un algorithme déterministe qui renvoie $t$ hyperplans de $\PP^n_k$ en position générale en $O(t^{\alpha n^2})$ opérations dans $k$.
\end{prop}

Si $k=\FF_q$ avec $q\leqslant t^{n+1}$, il convient de travailler sur une extension $k'$ de $\FF_q$ de degré $\lceil \log_q(t^{n+1})\rceil$, qui se calcule en temps polynomial en $q$ et $(n+1)\log_qt$. Il est alors possible de calculer avec l'algorithme ci-dessus $t$ hyperplans en position générale dans $\PP^n_{k'}$.

\begin{rk} Remarquons qu'avec un algorithme probabiliste, la construction étant donné $H_1,\dots,H_t$ en position générale, d'un hyperplan $H_{t+1}$ tel que $H_1,\dots,H_{t+1}$ soient en position générale se fait plus rapidement. En effet, en choisissant (avec les notations de la preuve de la proposition précédente) un élément aléatoire de la grille $R$, la probabilité qu'il appartienne à l'un des hyperplans $H_1,\dots,H_t$ est \[ \frac{t^{n^2+2n}}{n! t^{n^2+2n+1}}=\frac{1}{n! t}. \]
L'algorithme de type Las Vegas consistant à choisir un élément aléatoire de $R$ puis vérifier s'il vérifie l'une des équations des $H_i$ nécessite $O(tn)$ opérations dans $k$, et a une probabilité d'échec égale à $\frac{1}{n! t}$.
\end{rk}

\subsection{Hyperplans qui coupent une variété transversalement}\label{subsec:hyplis}

Soit $k$ un corps algébriquement clos. Soit $X\subset\PP^n$ un sous-schéma intègre lisse de dimension $d$. Soient $x\in X(k)$ et $r\leqslant n$. Le but de cette section est d'expliquer comment construire des hyperplans $H_1,\dots,H_{r}$ en position générale tels que l'intersection $C\coloneqq \cap_i H_i$ coupe $X$ transversalement et contienne $x$. Comme $X$ est lisse, il est localement intersection complète. Nous pouvons donc écrire $X=X_1\cup\dots \cup X_s$, où les $X_i$ sont des schémas affines d'intersection complète. Le calcul des $X_i$ est primitivement récursif : par la preuve explicite de \cite[00SC]{stacks}, il est possible pour chaque point $x\in X$ d'en calculer un voisinage ouvert qui est intersection complète ; le principe de la section \ref{subsec:recprimrec} conclut.
Remplaçons $X$ par l'un des $X_i$. \'{E}crivons $X=\Spec k[x_1,\dots,x_n]/(f_1,\dots,f_{n-d})\subset \PP^n=\Proj k[x_0,\dots ,x_n]$. Alors la preuve de \cite[XI, Th. 2.1.(i)]{sga43} montre comment déterminer explicitement des équations d'une partie constructible $Y$ de $\A^{n+1}$ paramétrant les hyperplans $H_a$ d'équation $a_0x_0+\dots+a_nx_n$ tels que $x\in H_a$ et que l'intersection $H\cap X$ ne soit pas transversale. Il suffit de trouver un $k$-point de $\PP^n$ qui ne se trouve pas dans $Y$ pour construire un premier hyperplan $H_1$. Ayant déjà construit $H_1,\dots, H_j$ en position générale, l'hyperplan $H_{j+1}$ se construit en choisissant un point de $(\PP^n)^\vee-Y$ qui vérifie également les contraintes décrites dans la section précédente (il suffit qu'il évite une réunion d'hyperplans de $\A^{n+1}$, tous distincts de celui qui paramètre les $H$ contenant $x$).\\

Le résultat ci-dessus s'étend à n'importe quel sous-schéma $X$ de $\PP^n$ dès que le point $x$ est régulier : il suffit d'appliquer la procédure à $X-X_{\rm sing}$.
Lorsque $r\geqslant d$, il est également possible d'adapter cet argument à un schéma $X$ non nécessairement lisse, pour construire des hyperplans $H_1,\dots,H_r$ qui évitent de plus $X_{\rm sing}$. \`{A} chaque étape $i\leqslant d$ de l'application de la procédure ci-dessus à $X-X_{\rm sing}$, calculons les composantes connexes de $X_{\rm sing}\cap H_1\cap \dots \cap H_i$, et déterminons un point de chacune : ceci fournit des points $P_1,\dots,P_s$. Au moment de choisir $H_{i+1}$, ajoutons à la partie $Y\subset (\PP^n)^\vee$ ci-dessus les hyperplans définis par les conditions linéaires $P_j\in H_{i+1}$. La condition $r\geqslant d$ assure que tous les points de $X_{\rm sing}$ sont évités. \\

Dans la section \ref{subsec:fibel}, nous avons besoin d'appliquer ces résultats à la situation suivante. Le schéma $X$ est normal et contient un fermé $V$ de dimension $d-1$ qui contient $X_\sing$. Comme $X$ est normal, $X_\sing$ est de dimension au plus $d-2$.
Le but est de construire des hyperplans $H_1,\dots,H_{d-1}$ dont l'intersection coupe $X$ et $V$ transversalement, et évite $X_\sing$. \`{A} chaque étape $i\in \{0,\dots,d-1\}$ sont construites des parties $Y_X^{(i)}$ et $Y_V^{(i)}$ de $(\PP^n)^\vee$ qui paramètrent les hyperplans qui ne conviennent pas. Il suffit de choisir un point de $(\PP^n)^\vee-(Y_X^{(i)}\cup Y_V^{(i)})$.\\

\section{Recherche de points}

Soit $k$ un corps algébriquement clos. Le but de cette section est le suivant : étant donné une variété $X$ plongée dans $\PP^n_k$, trouver au moins un point fermé de $X_{\bar k}$ dans chaque composante irréductible de $X$. Le fait de ne pas se préoccuper des multiplicités rend cette tâche plus facile que la décomposition primaire. Dans un premier temps, voyons comment trouver les points fermés des composantes zéro-dimensionnelles de $X$, appelés points isolés.

\subsection{Trouver des points isolés}

\begin{lem}\cite[Lem. 3.1]{ierardi_quantifier}\label{nbcomp}
Soit $k$ un corps algébriquement clos. Soient $f_1,\dots,f_m$ des polynômes homogènes de  $k[x_0,\dots,x_n]$ de degré inférieur à $d$. Alors pour tout $s\leqslant n$, le fermé $V(f_1,\dots,f_s)$ de $\PP^n$ a au plus $d^s$ composantes irréductibles de codimension inférieure à $s$. En particulier, il a au plus $d^n$ composantes irréductibles.
\end{lem}

Supposons $X$ défini par $m$ polynômes homogènes $f_1,\dots,f_m\in k[x_0,\dots,x_n]$ de degrés respectifs $d_1,\dots,d_m\leqslant d$. Quitte à ajouter des équations redondantes si $m<n$, ou à agrandir $X$ en oubliant certaines équations si $m>n$, nous pouvons supposer $n=m$. Notons $L_x(u)=\sum_{i=0}^n x_i^du_i\in k[x_0,\dots,x_n,u_0,\dots,u_n]$ et $\hat f_i(x,t)=tx_i^{d_i}+(1-t)f_i(x)\in k[x_0,\dots,x_n,t]$. Soit $\hat X \hookrightarrow \PP^n\times \A^1$ le fermé défini par les $\hat f_i$.

\begin{prop}\cite[Lem. 2.3]{ierardi_quantifier} Notons $\bar X=\overline{\hat X\cap D(t)}$ l'adhérence dans $\PP^n\times\A^1$ de l'intersection de $\hat X$ avec l'ouvert principal $D(t)$ (où $\A^1=\Spec k[t]$), et $\bar X_0$ l'intersection de $\bar X$ avec le fermé $V(t)$. (En bref : $\bar X_0 =\lim_{t\to 0} \hat X_t$) Alors $\bar X_0$ est de dimension zéro et contient tous les points fermés isolés de $X$.
\end{prop}

La procédure est décrite dans \cite[§2.2]{ierardi_quantifier}. Elle consiste à construire un polynôme $R(u_0,\dots,u_n)$ dont les zéros sont les mêmes que ceux de $\prod_{x\in \bar X_0}L_x(u_0,\dots,u_n)$. Ce $R\in k[u_0,\dots,u_n]$ est un résultant multivarié qui se calcule en temps $O (d^n)$ \cite[Prop. 2.5]{ierardi_quantifier}.

Il suffit ensuite de poser $R_i(t)=R(t,0,\dots,0,-1,0,\dots,0)$ et de factoriser $R_i(t)$ pour connaître la (puissance $d$-ième de la) projection sur le $i$-ième axe de coordonnées des points de $\bar X_0$. Ceci fournit $d^2$ coordonnées possibles sur chaque axe, et il suffit de tester les $d^{2n}$ combinaisons possibles pour savoir si elles définissent un point fermé de $X$.

\subsection{Trouver des points dans toutes les composantes}

\subsubsection{Pour un fermé de l'espace projectif}

Soit $X\subsetneq \PP^n$ donné par des équations homogènes. La stratégie pour trouver au moins un point dans chaque composante irréductible de $X$ consiste à se ramener au cas précédent en intersectant $X$ avec une succession d'hyperplans.\\
Considérons $n$ hyperplans $H_1,\dots,H_n$ en position générale (voir définition \ref{def:posgen}) ; par exemple ceux définis par $x_1,\dots,x_n$. Notons, pour $i>0$, $L_i=H_1\cap\dots\cap H_i$. Notons encore $L_0=\PP^n$.

\begin{lem} Soit $C$ une composante irréductible de $X$. Alors il existe $i\in \{0,\dots, n\}$ tel que $C\cap L_i$ soit non vide et de dimension nulle.
\begin{proof}
La suite $\dim (C\cap L_i)$ est décroissante. De plus, $\dim (C\cap L_i)\leqslant \dim(L_i)=n-i$ et $\dim (C\cap L_i)\geqslant \dim (C\cap L_{i-1})-1$. C'est donc une suite décroissante, qui décroît par pas de 0 ou 1, et qui stationne à 0. Soit $i$ le plus petit indice tel que $C\cap L_i$ soit de dimension nulle. Si $i=0$, alors $C\cap L_i=C$ est non vide et de dimension nulle. Si $i>0$, cela signifie que $\dim (C\cap L_{i-1})=1$.  D'après le théorème \ref{projdimthm}, l'intersection $C\cap L_i$ est non vide.
\end{proof}
\end{lem}

Ceci suggère l'algorithme suivant pour trouver au moins un point de chaque composante irréductible de $X$ : déterminer, pour chaque $i\in \{0,\dots,n\}$, les points isolés de $X\cap H_1\cap\dots\cap H_i$. La complexité est donc $n$ fois celle de la recherche de points isolés.

\subsubsection{Pour une partie constructible d'un espace projectif}\label{subsubsec:ptscompconstr}

\begin{lem} Soient $H_1,\dots,H_s$ des hyperplans en position générale dans $\PP^n$. Soit $C$ un fermé de $\PP^n$. Alors $C$ est contenu dans au plus $\mathrm{codim}(C)$ des hyperplans $H_i$.
\begin{proof}
La dimension de l'intersection de $\mathrm{codim}(C)+1$ de ces hyperplans, qui sont en position générale, est strictement inférieure à $\dim(C)$.
\end{proof}
\end{lem}

\begin{lem} Soit $C$ un fermé irréductible de $\PP^n$. Soit $H$ un hyperplan de $\PP^n$. Alors soit $C\subseteq H$, soit $\dim (C\cap H)=\dim(C)-1$.
\begin{proof}
Si $C$ n'est pas inclus dans $H$ alors $C\cap H$ est un fermé strict de $C$. Prenons-en une composante irréductible de dimension maximale $C'$. D'une part, $\dim(C')\geqslant \dim(C) -1$ par la proposition \ref{projdimthm} appliquée à l'intersection de $C$ et $H$. D'autre part, $\dim(C')<\dim(C)$ puisque $C'$ est un fermé irréductible strict de $C$.
\end{proof}
\end{lem}

\begin{cor}
Si $C$ est un fermé de $\PP^n$ de dimension $d$, et si $H_1,\dots,H_s$ sont des hyperplans en position générale avec $s>\dim C$, alors il existe $i_1,\dots ,i_d$ tels que l'intersection de $C$ avec $H_{i_1}\cap\dots\cap H_{i_d}$ soit de dimension nulle.
\begin{proof}
Il y a au plus $\codim C$ hyperplans $H_i$ contenant $C$. Il reste donc au moins un hyperplan ne le contenant pas : prenons-en un et notons-le $H_{i_1}$. Alors $C\cap H_{i_1}$ est inclus dans au plus $n$ des hyperplans $H_i$, dont ceux qui contiennent $C$. Il y a donc au moins un hyperplan $H_{i_2}$ qui ne le contient pas, et $\dim(C\cap H_{i_1}\cap H_{i_2})=\dim(C)-2$. Une récurrence immédiate conclut.
\end{proof}
\end{cor}

\begin{lem}
Si $C$ est un fermé de $\PP^n$ de dimension $d$, et si $H_1,\dots,H_{n-d}$ sont des hyperplans en position générale contenant $C$, alors $C=H_1\cap\dots \cap H_{n-d}$.
\begin{proof}
L'intersection est irréductible, de même dimension que $C$ et contient $C$ qui en est un fermé.
\end{proof}
\end{lem}

Adaptons l'algorithme précédent au cas d'une partie constructible d'un espace projectif. Considérons donc un fermé réduit $Y$ de $\PP^n$ défini par les équations homogènes $f_1=0,\dots,f_m=0$, et un ouvert $X$ de $Y$ défini par des inéquations homogènes $g_1\neq 0,\dots,g_s\neq 0$. Notons, pour chaque $i\in \{1\dots s\}$, $V_i$ le fermé de $\PP^n$ défini par $g_i$. Soit $d_X$ un majorant des degrés des $f_i$ et des $g_j$.  Construisons d'abord un ensemble 
$S=\{ H_1,\dots,H_D\}$ d'hyperplans de $\PP^n$ en position générale, avec $D=snd_X^n+1$.

Commençons par calculer les points isolés de $Y$ (dont font partie les points isolés de $X$). Ensuite, intersectons $Y$ avec chacun des hyperplans $H_i$. Les lemmes suivants montrent qu'alors toute composante de dimension 1 de $X$ intersecte au moins l'un des $H_i$ en un nombre fini non nul de points isolés, et que toute composante de dimension supérieure a une intersection non nulle avec au moins l'un des $H_i$.\\

Pour chaque $i$, la même procédure s'applique ensuite récursivement à $Y\cap H_i$. Elle s'arrête après une profondeur de récursion de $n$. Au total, elle nécessite le calcul de $O((sd_X^{n})^n)=O(s^nd_X^{n^2})$ intersections, et pour chacune de ces intersections, la détermination des points isolés du fermé de $\PP^n$ défini par $O(m+n)$ équations de degré $O(d_X)$. Chaque calcul d'intersection nécessite alors $O((m+n)d_X^{O(n)})$ opérations dans $k$. La complexité totale du calcul est donc de $ms^nd_X^{O(n^2)}$ opérations dans $k$.

\begin{lem}\label{lem:compirrposgen}
Soit $C$ une composante irréductible de $X$ de dimension 1. Si $S$ est un ensemble de $snd_X^n+1$ hyperplans en position générale alors il existe $H\in S$ tel que $C\cap H$ soit de dimension nulle et non vide.
\begin{proof}
Notons $\bar C$ l'adhérence de $C$ dans $Y$ : c'est un fermé de $\PP^n$. Pour tout $i$, comme $C$ n'est pas inclus dans $V_i$, l'intersection de $\bar C$ avec $V_i$ est de dimension nulle. En particulier, $C\cap V_i$ a $m_i\leqslant d_X^n$ points, et le nombre de points de $C\cap (\bigcup_i V_i)$ est $m\leqslant sd_X^n$.
Comme les hyperplans $H_i$ sont en position générale, chacun des $m_i$ points de $C\cap V_i$ est contenu dans au plus $n$ d'entre eux. Par conséquent, en prenant $nm+1\leqslant nsd_X^n+1$ hyperplans de $S$, au moins l'un d'entre eux intersecte $\bar C$ en un point qui n'est pas à l'infini.  
\end{proof}
\end{lem}

\begin{lem} Soit $C$ une composante irréductible de $X$. Notons $d_X=\dim(X)$. Soit $S$ un ensemble de $snd_X^n+1$ hyperplans en position générale. Il existe $H_1,\dots,H_d\in S$ tels que $C\cap H_{1}\cap\dots\cap H_{d}$ soit de dimension nulle et non vide.
\begin{proof} Notons $\bar C$ l'adhérence de $C$ dans $Y$. Pour chaque indice $i\in \{1,\dots,s\}$, le lemme \ref{nbcomp} assure que $\bar C\cap V_i$ a au plus $d_X^n$ composantes irréductibles. Il existe donc au plus $snd_X^n$ hyperplans de $S$ contenant l'une des composantes irréductibles de l'un des $\bar C\cap V_i$. Prenons-en un qui ne contient aucun $\bar C\cap V_i$, et notons-le $H_{1}$ ; en particulier, pour tout $i$, la dimension de $H_{1}\cap\bar{C}$ est $\dim(\bar C)-1$, et celle de $H_{1}\cap \bar{C}\cap V_i$ est égale à $\dim(\bar C)-2$. Par conséquent, $H_{1}\cap\bar{C}$ contient un point de $C$.
Il existe alors au plus $nsd_X^n$ hyperplans de $S$ contenant l'un des $\bar{C} \cap H_{1}\cap V_i$ ; leur ensemble contient celui des hyperplans précédemment écartés contenant l'un des $\bar C\cap V_i$. Ceci fournit $H_{2}$ qui ne contient aucune de ces intersections. Continuons ainsi jusqu'à obtenir $H_{1},\dots,H_{d-1}$. Le lemme \ref{lem:compirrposgen} permet alors de construire $H_d$ tel que tels que $C'\coloneqq \bar{C} \cap H_{1}\cap\dots\cap H_{d}$ soit de dimension nulle, et contienne un point qui n'est dans aucun des $V_i$, c'est-à-dire un point de $C$. Le nombre d'hyperplans en position générale à éviter est majoré par $nsd_X^n$ ; tout ensemble de $nsd_X^n+1$ hyperplans en position générale contient donc un élément convenable.
\end{proof}
\end{lem}

Par conséquent, il suffit de construire $O(nsd_X^n)$ hyperplans en position générale dans $\PP^n$. Ceci se fait en temps $(nsd_X^n)^{O(n^2)}$ avec l'algorithme décrit dans la section \ref{subsec:hyperplans}.

\subsection{Le cas équidimensionnel}\label{subsubsec:ptsequidim}

Cette méthode, bien plus élégante que celle présentée ci-avant, est décrite dans \cite[Alg. 8.5]{jinbi_jin}. Soit $X=\Spec k[x_1,\dots,x_n]/(f_1,\dots,f_m)$ un schéma affine de type fini sur $k$ dont toutes les composantes irréductibles sont de même dimension $r$. Soit $\nu\colon X\to \A^r_k$ une normalisation de Noether de $X$.

\begin{lem} La restriction de $\nu$ à toute composante irréductible de $X$ est encore surjective.
\begin{proof} Si $C$ est une composante irréductible de $X$ alors $\nu|_C$ est encore un morphisme fini, et donc $\dim(\nu(C))=\dim(C)=\dim(X)=r$ par équidimensionnalité de $X$. De plus, encore par finitude, $\nu(C)$ est un fermé de $\A_k^{r}$ de dimension $r$, il est donc égal à $\A_k^{r}$.
\end{proof}
\end{lem}
La fibre $\nu^{-1}(0)=\colon \Spec R$ intersecte donc chaque composante irréductible de $X$. Une décomposition primaire de $R$ (présenté par générateurs et relations), qui est une $k$-algèbre finie, fournit pour chaque composante primaire $C$ une base de Gröbner de son réduit. Ceci permet d'obtenir par factorisation une liste de $\bar{k}$-points, qui contient au moins un élément de chaque composante irréductible. L'algorithme est donc composé de 4 étapes : un calcul de normalisation de Noether, une décomposition primaire, des calculs de bases de Gröbner d'idéaux zéro-dimensionnels puis des factorisations. 

Nous n'avons pas précisé, dans la section sur les bases de Gröbner, la complexité des différents algorithmes; cependant, la complexité de la normalisation de Noether décrite dans \cite[Alg. 1.13]{dickenstein} est déjà, pour un idéal de $k[x_1,\dots,x_n]$ engendré par $m$ équations de degré maximal $d$, de l'ordre de $m^3d^{O(n^2)}$. Il n'y a donc pas de gain de complexité notable par rapport à la méthode par recherche de points isolés d'intersections avec des hyperplans présentée dans la section \ref{subsubsec:ptscompconstr}.

\section{Restriction de Weil}\label{sec:weilres}

\begin{df} Soit $f\colon Y\to X$ un morphisme de schémas. Soient $V$ un $Y$-schéma, et $h_V$ le faisceau représenté par $V$. Si le préfaisceau $f_\star h_V\colon U\mapsto \Hom_Y(U\times_X Y,V)$ sur $\Sch/X$ est représentable par un $X$-schéma $R$, ce dernier est appelé restriction de Weil de $V$ relativement à $f$ et noté $R_f(V)$, ou $R_{Y\to X}(V)$ si le contexte le permet.
\end{df}

\begin{prop}\cite[7.6, Th. 4]{blr} Soit $f\colon Y\to X$ un morphisme fini localement libre de schémas. Soit $V$ un $Y$-schéma tel que pour tout $x\in X$, tout sous-ensemble fini de la fibre $V_x$ soit contenu dans un ouvert affine de $V$ (c'est le cas par exemple si $V$ est affine). Alors $R_{Y\to X}(V)$ existe. 
\end{prop}

En particulier, si $X$ et $Y$ sont affines, le foncteur $R_{Y\to X}$ est l'adjoint à droite de $-\times_X Y$ sur les catégories de schémas affines sur $X$ et $Y$.
Soient
\[
\begin{array}{rcl}
A&=&k[x_1,\dots,x_n]/(f_1,\dots,f_r) \\ B&=&k[y_1,\dots,y_m]/(g_1,\dots,g_s) \\ C&=&B[v_1,\dots,v_p]/(h_1,\dots,h_t).\end{array}\] Notons $X=\Spec A$, $Y=\Spec B$, $V=\Spec C$. Supposons donné un morphisme $F\colon A\to B$ défini par $F_1,\dots,F_r\in k[y_1,\dots,y_m]$ qui fait de $B$ un $A$-module libre, et une base $(b_1,\dots,b_N)$ de $B$ comme $A$-module. Dans ce cadre, la proposition précédente assure que la restriction de Weil $R_{Y\to X}(V)$ existe et se calcule de la façon suivante (voir \cite[4.4.1]{scheiderer}).\\

Introduisons des indéterminées $w_{ij}$ où $i\in \{ 1,\dots, p\}$ et $j\in \{1,\dots N\}$, et définissons, pour $\alpha\in \{ 1,\dots, t\}$ et $\beta\in \{1,\dots,N\}$, des éléments $f_{\alpha\beta}\in A[w_{ij}]$ par :
\[ h_\alpha\left( \sum_{j=1}^N w_{1j}b_j,\dots,\sum_{j=1}^N w_{pj}b_j\right)=\sum_{\beta=1}^N f_{\alpha\beta}b_\beta \in B[w_{ij}]. \]
La restriction de Weil $R_{Y\to X}(V)$ est alors $\Spec k[w_{ij}]/(f_{\alpha\beta})$. Le morphisme d'adjonction \[V\to R_{Y\to X}(V)\times_X Y\] est donné par $k[v_i]/(h_{\gamma})\to k[w_{ij}]/(f_{\alpha\beta}), v_i\mapsto \sum_j w_{ij}b_j$. De même, si $V=U\times_X Y$, le morphisme d'adjonction $U\to R_{Y\to X}(V)$ est donné en écrivant $1=\sum_j a_{j}b_j$ avec $a_j\in A$, et en envoyant $w_{ij}$ sur $a_jv_i\in A[v_i]$.\\

La restriction $R_{Y\to X}(\phi)$ d'un morphisme $\phi\colon  V\to V'$ de $Y$-schémas se calcule par la même méthode.
De plus, si $V$ est un schéma en groupes sur $Y$, la même méthode permet de calculer la loi de groupe sur $R_{Y\to X}(V)$. En effet, comme la restriction de Weil est un adjoint à droite, elle commute aux limites, et $R_{Y\to X}(V\times_Y V)=R_{Y\to X}(V)\times_X R_{Y\to X}(V)$. Il suffit donc de déterminer le morphisme $R_{Y\to X}(V\times_Y V\to V)$.\\

\cleartooddpage

\chapter{Algorithmique des courbes}\label{chap:A3}

\section{Représentations des courbes lisses}\label{sec:repcourb}

Par défaut, une courbe sur un corps $k_0$ est représentée par un recollement de courbes affines comme décrit dans l'annexe \ref{sec:repsch}. Une courbe intègre sur $k_0$ étant toujours isomorphe à un ouvert d'un fermé d'un espace projectif \cite[0A27]{stacks}, d'autres représentations des courbes seront données ici.

\subsection{Courbes singulières, compactification lisse}

\subsubsection{Résolution des singularités d'une courbe plane}\label{subsubsec:desing}

Soit $k_0$ un corps. Soit $C=\Proj k_0[x,y,z]/(f)$ une courbe projective plane intègre de normalisée $X$. Notons $d=\deg(f)$. La courbe $C$ a au plus $\frac{(d-1)(d-2)}{2}$ points singuliers \cite[Th. 3.8]{fischer}. 
\begin{df}
Un point $P$ d'une courbe plane $C$ est dit ordinaire si le nombre de tangentes à $C$ en $P$ est égal à la multiplicité de $C$ en $P$.
\end{df}
Par moins de $d^2$ transformations quadratiques, il est possible \cite[7.4, Th. 2]{fulton} de construire une courbe $X'$ birationnelle à $X$ ayant uniquement des singularités ordinaires. Le degré de la courbe $X'$ est $O(2^{d^2})$.\\
La méthode classique de résolution des singularités sur une variété quelconque procède par éclatements successifs ; elle a l'avantage de construire directement un plongement projectif de $X$, mais l'inconvénient que l'espace projectif dans lequel $X$ est plongée est de grande dimension. En partant de $X'$ à singularités ordinaires, il suffit d'un seul éclatement par singularité. Les algorithmes n'ont souvent pas besoin de connaître un plongement projectif de $X$, mais simplement une description des points de $X$ au-dessus de chaque singularité de $C$.

\begin{prop}\cite[§5]{kozen}
Il existe un algorithme déterministe qui prend en entrée une courbe $C\subset \A^2_{\FF_q}$ (resp. $\A^2_\QQ$) définie par un polynôme $f=\sum_{i,j}a_{ij}x^iy^j\in \FF_q[x,y]$ (resp. $\ZZ[x,y]$) de degré $d$ tel que $(0,0)\in C$, et retourne l'arbre de désingularisation de $C$ en $(0,0)$, en un nombre d'opérations polynomial en $d$ et en $\log q$ (resp. en $\log \max |a_{ij}|$).
\end{prop}

\subsubsection{Compactification lisse}

Soit $k_0$ un corps parfait. Soit $X=\Spec k_0[x_1,\dots,x_m]/(f_1,\dots,f_r)$ une courbe affine lisse sur $k_0$. En homogénéisant les $f_i$, on obtient l'adhérence $Y$ de $X$ dans $\PP_{k_0}^m$. La courbe $Y$ est possiblement singulière en-dehors de son ouvert isomorphe à $X$. La normalisation $\bar X$ de $Y$ est une courbe projective lisse dont un ouvert est isomorphe à $X$. Cette courbe $\bar X$ est appelée compactification lisse de $X$ ; elle est unique à isomorphisme près.
\'{E}tant donné un morphisme de courbes affines lisses $X_1\to X_2$ sur $k_0$, le morphisme composé $X_1\to X_2\to\bar X_2$ s'étend à $\bar X_1$ par régularité \cite[I, Prop. 6.8]{hartshorne}.

\subsection{Corps de fonctions, modèle plan birationnel}\label{subsec:modplan}

Soit $k_0$ un corps. Rappelons que le foncteur qui à une courbe associe son corps de fonctions définit une équivalence entre la catégorie des courbes projectives régulières sur $k_0$ munie des morphismes non constants, et la catégorie des extensions de $k_0$ de degré de transcendance 1 \cite[0BY1]{stacks}. Une courbe projective régulière peut donc être définie par son corps de fonctions.

\begin{prop}\cite[Prop. 3.10.2.(a)]{stichtenoth}\label{prop:extsep} Soit $X$ une courbe connexe lisse sur un corps parfait $k_0$. Soient $P\in X(k_0)$ et $t$ une uniformisante en $P$. Alors le corps des fonctions $k_0(X)$ est une extension finie séparable de $k_0(t)$.
\end{prop}

\'{E}tant donné une telle courbe $X$ plongée dans $\PP^n_{k_0}$, un morphisme $X\to\PP^1_{k_0}$ fini génériquement étale se calcule donc en cherchant une uniformisante $t$ en un point $P\in X(k_0)$. Si $X(k_0)$ est vide, il suffit de remplacer $k_0$ par une extension $k_0'$ (de degré $O(\deg X)$) sur laquelle $X$ a un point pour obtenir un morphisme $X_{k_0'}\to \PP^1_{k_0'}$. Sans perte de généralité, supposons que $P$ se situe dans l'ouvert affine $U_0=\{x_0\neq 0\}$ de $\PP^n$. La courbe $X$ est donnée par des polynômes $f_1,\dots,f_m\in k[x_1,\dots,x_n]$ avec $m\geqslant n-1$. Comme $X$ est lisse sur $k_0$, la matrice $(\frac{\partial f_i}{\partial x_j})_{i,j}$ est de rang $n-1$. 
Il existe donc $j\in \{1,\dots,n\}$ tel que le morphisme $x_j-x_j(P)\colon X\cap U_0\to\A^1$ soit étale en $P$, et donc que $t\coloneqq x_j-x_j(P)$ soit une uniformisante en $P$. Le degré du morphisme $X\to\PP^1$ induit par $t$ est donc de degré au plus $\deg X$, qui est le nombre maximal de points d'intersection de $H$ avec un hyperplan de $\PP^n$. 

Afin de travailler avec un modèle birationnel plan de $X$, il suffit de chercher un élément primitif $u$ de l'extension finie séparable $k_0(X)/k_0(t)$ à l'aide de l'algorithme décrit dans la section \ref{subsubsec:eltprim}. Rappelons que si $k(X)$ est engendré par $(x_1,x_2,\dots,x_n)$, un tel élément primitif s'obtient sous la forme $u=t+\sum_{i\neq j}\lambda_i x_i$ avec $\lambda_i\in k$.
Alors $k_0(X)=k_0(t,u)= k_0(t)[u]/(f)$, où $f$ est le polynôme minimal de $u$ sur $k_0(t)$. Ceci définit un morphisme birationnel $X\to C$, où $C$ est une courbe projective plane d'équation dans une carte affine $f(t,u)=0$. Supposons, quitte à le multiplier par un élément de $k_0[t]$, que $f$ est un polynôme primitif dans $k_0[t][u]$. Le lemme de Gauss assure alors que $f$ est irréductible dans $k_0[t,u]$ et dans $k_0(u)[t]$ ; par conséquent, $\deg(f)=[k_0(X):k_0(u)]$. De même que ci-dessus, comme $u=t+\sum_i\lambda_ix_i$ est de degré 1, le degré du morphisme $X\to\PP^1$ engendré par cette fonction est majoré par $\deg X$. Or $\deg_t(f)=[k_0(X):k_0(u)]\leqslant\deg X$, et $\deg_u(f)=[k_0(X):k_0(t)]\leqslant \deg X$. Par conséquent, $\deg(C)\leqslant\deg(X)^2$.\\

De plus, le degré d'un modèle plan peut être borné en fonction du genre de la courbe $X$. Si $X$ est hyperelliptique, il est connu qu'elle admet un modèle plan de la forme \[ y^2+h(x)y=f(x) \]
avec $\deg h,\deg f\leqslant 2g+2$. Ce modèle s'obtient grâce au revêtement double $X\to \PP^1$. Si $X$ n'est pas hyperelliptique, son diviseur canonique $K$ est très ample, et le choix de $g-3$ points généraux $P_1,\dots,P_{g-3}$ de $X$ permet d'obtenir un diviseur $D=K-\sum_i P_i$ birationnellement très ample dont l'espace de Riemann-Roch $\mathscr{L}(D)$ est de dimension 3 \cite[§1]{keem_martens}. Dans tous les cas, il est possible de construire explicitement un modèle birationnel plan de $X$ de degré $O(g)$.\\
Par conséquent, toute courbe projective lisse admet un modèle plan à singularités ordinaires de degré $O(2^g)$.

\subsection{Représentation comme $\OP$-algèbre}

\subsubsection{Description}\label{subsub:OP1alg}

Cette description est celle employée par Jin dans \cite{jinbi_jin}. Soit $X$ une courbe intègre projective lisse sur un corps parfait $k_0$. Supposons-la décrite comme fermé d'un espace projectif sur $k_0$. 
D'après la proposition \ref{prop:extsep}, le calcul d'une uniformisante permet d'obtenir un morphisme $X\to \PP^1$ dont l'extension de corps de fonctions correspondante est séparable. Le morphisme $X\to \PP^1$ est donc génériquement étale par \cite[5.4.3]{mumford_oda}.

\begin{lem}Le faisceau $\phi_\star\OO_X$ est un fibré vectoriel sur $\PP^1$.
\begin{proof} Le morphisme $\phi$ est un morphisme non constant d'une courbe vers une courbe régulière, il est donc plat par \cite[0CCK]{stacks}. Un module de présentation finie est localement libre si et seulement si il est plat \cite[00NX]{stacks} ; par conséquent, $\phi_\star \OO_X$ est un $\OO_{\PP^1}$-module localement libre.
\end{proof}
\end{lem}

Notons $U_0=\Spec k_0[x]$ et $U_1=\Spec k_0[x^{-1}]$ les ouverts standard de $\PP^1$. 

\begin{lem} Si $\mathcal{F}$ est un $\OO_{\PP^1}$-module localement libre de type fini alors $\mathcal{F}(U_0)$ est un $k_0[x]$-module libre de rang fini, et $\mathcal{F}(U_1)$ est un $k_0[x^{-1}]$-module libre de rang fini.
\begin{proof}Comme $\mathcal{F}$ est un module localement libre de type fini, il est projectif \cite[00NX]{stacks}. En particulier, $\mathcal{F}_{|U_0}$ est encore projectif, il est donc libre puisque $k_0[x]$ est principal.
\end{proof}
\end{lem}

Un fibré vectoriel $\mathcal{F}$ de rang $r$ sur $\PP^1$ est défini par la donnée du $k_0[x]$-module $\mathcal{F}(U_0)$, du $k_0[x^{-1}]$-module $\mathcal{F}(U_1)$ et d'un isomorphisme de $k_0[x^{\pm 1}]$-modules \[ \mathcal{F}(U_0)\otimes_{k_0[x]}k_0[x^{\pm 1}] \to \mathcal{F}(U_1)\otimes_{k_0[x^{-1}]}k_0[x^{\pm 1}].\] Ces deux modules étant libres de rang $r$, cet isomorphisme se représente par une matrice $M_\mathcal{F}$ à coefficients dans $k_0[x^{\pm 1}]$ qui dépend des bases choisies pour ces modules libres.

\paragraph{Calcul de $\phi_\star\OO_X(U_0)$} La préimage $\phi^{-1}\Spec k_0[x]$ est la normalisation de $k_0[x]$ dans l'extension de corps de fonctions $\phi^\star \colon  k_0(x)\to k_0(C)$. Comme $k_0$ est parfait, on sait calculer un modèle plan de $C$ par le théorème de l'élément primitif, et ainsi présenter $k_0(C)$ comme $k_0(x)[y]/(f)$. Ensuite, on calcule la normalisation dans une extension de corps de la façon usuelle ; elle a dans ce cas précis une complexité plus abordable, comme décrit dans \cite[Prop. 2.127]{diem}.
D'après le théorème de Dedekind-Weber-Grothendieck \cite[Th. 11.50]{goertz_wedhorn}, un fibré vectoriel sur $\PP^1$ est isomorphe à une somme directe de fibrés en droites : \[ \mathcal{F}\simeq \OO_{\PP^1}(n_1)\oplus\cdots\oplus\OO_{\PP^1}(n_r)\] où les $n_i$ sont uniques à l'ordre près. De façon concrète, cela signifie qu'il existe des bases des modules libres $\mathcal{F}(U_0)$ et $\mathcal{F}(U_1)$ explicitement calculables telles que la matrice $M_\mathcal{F}$ soit égale à \[ \left(\begin{matrix}
x^{n_1} & & \\
& \ddots & \\
& & x^{n_r}
\end{matrix}\right).\]

Cette somme de fibrés en droites est munie d'une structure de $\OO_{\PP^1}$-algèbre de la façon suivante. D'une part, le morphisme canonique $\OO_{\PP^1}\to\phi_\star\OO_X$ définit un morphisme $\OO_{\PP^1}\to \OO_{\PP^1}(n_1)\oplus\cdots\oplus\OO_{\PP^1}(n_r)$. D'autre part, le morphisme de multiplication \[ (\bigoplus_i\OO(n_i))\otimes_{\OO_{\PP^1}}(\bigoplus_i\OO(n_i))\to\bigoplus_i\OO(n_i)\]
est défini sur $U_0$ par la multiplication dans la $k_0[x]$-algèbre $\mathcal{F}(U_0)$ écrite dans la base choisie pour $\mathcal{F}(U_0)$, et sur $U_1$ par la multiplication dans la $k_0[x^{- 1}]$-algèbre $\mathcal{F}(U_1)$ écrite dans la base choisie pour $\mathcal{F}(U_1)$. Par l'isomorphisme \[ (\bigoplus_i\OO(n_i))\otimes_{\OO_{\PP^1}}(\bigoplus_i\OO(n_i)) \xrightarrow{\sim} \bigoplus_{i,j}\OO(n_i+n_j)\] cette multiplication équivaut à la donnée d'un morphisme \[ \bigoplus_{i,j}\OO(n_i+n_j)\to \bigoplus_i\OO(n_i).\]

Comme $\underline{\Hom}_{\OP}(\OO_{\PP^1}(a),\OO_{\PP^1}(b))\simeq\OO_{\PP^1}(b-a)$, il y a des morphismes non triviaux $\OO_{\PP^1}(a)\to\OO_{\PP^1}(b)$ si et seulement si $b\geqslant a$. Dans ce cas, un tel morphisme est défini par un polynôme homogène de $k_0[X,Y]$ de degré $b-a$. La multiplication est donc définie par une matrice de taille $r\times r^2$ dont les coefficients sont des polynômes homogènes de $k_0[X,Y]$. Pour que la multiplication soit commutative, il faut et il suffit que les colonnes correspondant à $(n_i,n_j)$ et $(n_j,n_i)$ soient égales. Si $X$ est géométriquement réduite alors $n_i\leqslant 0$ pour tout $i$ \cite[Lem. 6.17]{jinbi_jin}.

Notons qu'un morphisme de $\OP$-modules $\bigoplus_{j=1}^s\OP(b_j)\to\bigoplus_{i=1}^r\OP(a_i)$ est donné par une matrice de taille $r\times s$ à coefficients dans $k_0[X,Y]$, dont le coefficient en position $(i,j)$ est nul si $b_j>a_i$, et homogène de degré $a_i-b_j$ ou nul sinon. La composition de deux tels morphismes est définie par la multiplication matricielle. Un endomorphisme d'un $\OP$-module est inversible si et seulement si la matrice qui le décrit l'est.

\paragraph{Reconstruction de la courbe} La courbe $X$ peut être reconstruite à partir de la structure de $\OO_{\PP^1}$-algèbre sur $\bigoplus_i\OO(n_i)$ induite par celle de $\phi_\star\OO_X$. En effet, le morphisme $\phi \colon  X\to\PP^1$ est un morphisme fini localement libre, il est donc affine, et par conséquent $\underline{\Spec}_{\PP^1}(\phi_\star \OO_X)\simeq X$, où $\underline{\Spec}$ désigne le spectre relatif. Il suffit donc de déterminer les courbes affines $\Spec (\oplus \OO_{\PP^1}(n_i)(U_0))$ et $\Spec (\oplus \OO_{\PP^1}(n_i)(U_1))$ et de les recoller.

\begin{rk} Certains couples (fibré, matrice de multiplication) ne définissent pas une courbe. Par exemple, le fibré $\OO\oplus\OO(-1)\oplus\OO(-2)$ avec la matrice 
\[\left(
\begin{matrix}
1 & & & &Z^3  & Z^4 \\
 &1 & &2Z & & \\
 & &1 & & & 
\end{matrix}
\right)\]
définirait au-dessus de l'ouvert $Z\neq 0$ de $\PP^1$ le schéma $\Spec k_0[x,u,v]/(u^2-2,uv-1,v^2-1)$, qui est vide. De plus, rien ne garantit a priori la lissité du schéma obtenu.

Afin de s'assurer qu'un couple (fibré, matrice de multiplication) définit une courbe lisse, il convient de tester séparément si le schéma obtenu est de dimension 1 sur $k_0$ (voir \ref{subsec:dimension}) et s'il est lisse (à l'aide du critère jacobien).
\end{rk}

\begin{rk} La même description est encore valable pour des courbes propres lisses qui ne sont pas nécessairement connexes : il suffit de faire la somme directe des fibrés obtenus pour chaque composante connexe. Cette remarque sert dans la section \ref{sec:jin}.
\end{rk}

\section{Produit fibré de courbes intègres lisses}\label{sec:prodfib}

Considérons un diagramme cartésien 
\[
\begin{tikzcd}
T \arrow[r]\arrow[d]& Y\arrow[d] \\
Z \arrow[r] & X
\end{tikzcd}
 \]
un diagramme cartésien de courbes sur un corps $k_0$ ; supposons $X,Y,Z$ intègres, et $Y\to X$ lisse. Notons $\tX,\tY,\tZ,\tT$ les normalisées respectives de $X,Y,Z,T$. L'objectif de cette section est de calculer le produit fibré $\tY\times_\tX  \tZ$, qui est encore lisse sur $k_0$ puisque $Y\to X$ l'est. Notons $S=\tY\times_{\tX}\tZ$.

\begin{lem} Il y a un isomorphisme canonique $S\to \tilde T$.
\begin{proof}
Le diagramme commutatif \[
\begin{tikzcd}
\tT \arrow[r]\arrow[d]& \tY\arrow[d] \\
\tZ \arrow[r] & \tX
\end{tikzcd}
 \]
obtenu par fonctorialité de la normalisation produit un morphisme $\tT\to S$. Le diagramme suivant
\[
\begin{tikzcd}
S \arrow[r]\arrow[d]& Y\arrow[d] \\
Z \arrow[r] & X
\end{tikzcd}
 \]
produit un morphisme $S\to T$, qui se factorise par $\tT$. La composée $S\to \tT\to S$ est l'unique morphisme $S\to S$ faisant commuter
\[
\begin{tikzcd}
S \arrow[dr] \arrow[drr] \arrow[ddr] & &\\
&S \arrow[r]\arrow[d]& Y\arrow[d] \\
&Z \arrow[r] & X
\end{tikzcd}
 \] 
\end{proof}
\end{lem}

Pour calculer $\tY\times_\tX\tZ$, qui a pour anneau total des fractions un produit de corps, il suffit donc de calculer $Y\times_X Z$, qui a ce même anneau total des fractions. Les lemmes suivants montrent comment en calculer un élément primitif sur $k(x)$.

\begin{lem} Soient $K$ un corps et $L=K(z,t)$ une $K$-algèbre réduite de dimension finie $d$ sur $L$. Dans tout ensemble $S$ de $d^{4d}$ éléments de $k$, il y a au moins un élément $\lambda$ tel que $\lambda z+ t$ engendre $L$ en tant que $K$-algèbre.
\begin{proof} Notons $L=L_1\times\dots\times L_r$, où les $L_i$ sont des corps, et $d_i=[L_i:K]$. Soit $x=(x_1,\dots,x_r)\in L$. Pour tout $i\in \{1\dots r\}$, l'élément $(x_1,\dots,x_i)$ est primitif dans $L_1\times\dots\times L_i$ si et seulement si $(x_1,\dots,x_{i-1})$ est primitif dans $L_1\times\dots\times L_{i-1}$, l'élément $x_i$ est primitif dans $L_i$ et le polynôme minimal de $x_i$ sur $K$ est premier au polynôme minimal de $(x_1,\dots,x_{i-1})$ sur $k$. Cherchons un élément primitif sous la forme $z+\lambda t$ avec $\lambda\in k$. Rappelons qu'il y a moins de $d_i^2$ valeurs de $\lambda$ pour lesquelles $z+\lambda t$ n'est pas un élément primitif de $L_i$.
Ainsi, le nombre de valeurs de $\lambda$ à éviter est inférieur à \[ d_1^2\cdot d_1 d_2^2\cdot (d_1+d_2)d_3^2\cdots (d_1+\dots+d_{r-1})d_r^2.\]
Comme $d_1+\dots+d_r=d$, ce nombre est inférieur à $d^{4d}$.
\end{proof}
\end{lem}

Dans notre situation, écrivons $Z=\Spec k_0[a,b]/u$ et $Y=\Spec k_0[s,t]/v$. Appliquons le lemme d'abord à l'extension $L=k(a,b)[s,t]$ de $k(a,b)$, en remarquant que $k(a,b)$ est bien un corps. Commençons par trouver $\lambda\in k$ tel que $L=k(a,b)[\lambda s+t]$. Appliquons ensuite le lemme à l'extension $k(a)[b,\lambda s+t]$ de $k(a)$. Le degré de la première extension est le degré $[Y:X]$ du morphisme $Y\to X$. Le degré de la deuxième extension est $\deg(Z)[Y:X]$. Il suffit donc de tester $(\deg Z)^{4\deg Z}[Y:X]^{8[Y:X]}$ éléments de la forme $\mu b+\lambda s+t$ : pour chacun d'entre eux, on calcule le polynôme minimal par des moyens d'algèbre linéaire sur $k_0(a)$ faisant intervenir des polynômes de degré majoré par $d=\max(\deg(Z),\deg(Y),\deg(Z\to X),\deg(Y\to X))$. La complexité de ces opérations d'algèbre linéaire est polynomiale en $d$. La complexité totale est donc $O(d^{12d}\mathcal{P}(d))$, où $\mathcal{P}$ est un polynôme, et c'est encore $O(d^{13d})$.

\begin{prop} Soient $Y,Z$ deux courbes intègres planes sur $k_0$. Considérons une courbe plane $X$ et des morphismes $Z\to X$, $Y\to X$ dont l'un au moins est lisse. Notons $\tX, \tY, \tZ$ les normalisées respectives de $X,Y,Z$. Soit \[d=\max \left\{\deg(Z),\deg(Y),\deg(Y\to X),\deg(Z\to X)\right\}.\] Il existe un algorithme déterministe calculant une courbe plane birationnelle à $\tZ\times_\tX\tY$  en $O(d^{13d})$ opérations dans $k$. La courbe obtenue est de degré $\deg(Z)[Y:X]^2\leqslant d^3$.
\end{prop}

Notons $W$ la courbe plane obtenue. Rappelons que les composantes connexes de $\tY\times_\tX\tZ$ sont les normalisées des composantes irréductibles de $W$ \cite[0CDV]{stacks}. Afin de déterminer les composantes connexes de $\tY\times_\tX\tZ$, il suffit donc de factoriser le polynôme définissant $Y\times_X Z$. 

\begin{prop} Soient $Y,Z$ deux courbes planes sur $k_0$. Considérons une courbe plane $X$ et des morphismes $Z\to X$, $Y\to X$ dont l'un au moins est lisse. Soit \[d=\{\max \deg(Z),\deg(Y),\deg(Y\to X),\deg(Z\to X)\}.\] Il existe un algorithme déterministe calculant les composantes connexes de $\tY\times_\tX \tZ$ en $O(d^{13d})\mathrm{Fact}_{k_0}(d^2)$ opérations dans $k_0$, où $\mathrm{Fact}_{k_0}(\cdot)$ désigne la complexité de la factorisation d'un polynôme en deux variables de degré total donné dans $k_0$.
\end{prop}

\section{Diviseurs et espaces de Riemann-Roch}\label{sec:divrr}

\subsection{Représentations des classes de diviseurs}\label{subsec:repdiv}

Soient $k_0$ un corps et $k$ une clôture algébrique de $k_0$. Soient $X_0$ une courbe projective lisse sur $k_0$ et $X=X_0\times_{k_0}k$. Supposons d'abord $X_0$ décrite par un plongement projectif. Un diviseur $D$ sur $X$ peut être représenté de deux manières : par une combinaison linéaire de points fermés de $X_0$ (chaque point étant défini par une extension $k_1$ de $k_0$ et des $k_1$-points de $X_0$) ou par sa forme de Chow.

\begin{df} Soit $Z=a_iP_1+\dots+a_rP_r$ un zéro-cycle de $\PP^n_k$. Pour tout $i$, notons $p_{ij}$ les coordonnées de $P_i$. La forme de Chow de $Z$ est le polynôme
\[ \prod_{i=1}^r\left(\sum_{j=0}^n p_{ij}u_j\right) ^{a_i}\in k[u_0,\dots,u_n]. \]
\end{df}

En particulier, si $X=\Proj k[x,y,z]/(F)$ est une courbe dans $\PP^2$, le diviseur d'un polynôme homogène $G\in k[x,y,z]$ a pour forme de Chow le $u$-résultant \cite[Lem. 15.7]{ierardi_kozen}
\[ u-\res(F,G)\coloneqq \res(F,G,xu_x+yu_y+zu_z).\]
Le calcul de la forme de Chow d'un diviseur représenté par une somme de points est évident : il suffit de calculer une extension de $k_0$ sur laquelle sont définis tous les points du diviseur.
\\

Il existe encore une autre représentation des diviseurs, qui ne travaille pas explicitement avec des équations de $X$ dans un espace projectif. D'abord proposée par Khuri-Makdisi \cite{khuri_makdisi}, elle a été utilisée dans de nombreux algorithmes probabilistes de Bruin \cite{bruin}. Supposons donné un faisceau inversible $\LL$ sur $X_0$ de degré strictement supérieur à $2g$. Il est nécessairement très ample et fournit un plongement projectif de $X$ ; soit $S$ l'anneau de coordonnées homogènes de $X$ pour ce plongement. Le morphisme 
\[ S\to \bigoplus_{i\geqslant 0}\Gamma(X,\LL^{\otimes i})\]
est un isomorphisme. Pour les calculs, il suffit de connaître le quotient $S^{(h)}$ de $S$ par l'idéal engendré par les éléments homogènes de degré strictement supérieur à un entier $h$ suffisamment grand. Un diviseur $D$ est alors représenté par l'espace des sections globales du faisceau $\LL(-D)$. Les algorithmes qui font usage de cette représentation n'ont donc pas besoin d'une réelle description de $X$, mais simplement de l'algèbre $S^{(h)}$ et du sous-espace vectoriel $\Gamma(\LL(-D))$ de $S^{(1)}=\Gamma(X,\LL)$. \'{E}tant donné un plongement projectif de $X$, un diviseur $E$ tel que $\LL=\OO_X(E)$ et un diviseur $D$ sur $X$ représenté de l'une des façons ci-dessus, l'espace $\Gamma(\LL(-D))=\Gamma(\OO_X(E-D))$ est calculé par les algorithmes de la section suivante. Réciproquement, étant donné l'espace $\Gamma(\LL(-D))$, une décomposition de $D$ en somme de diviseurs premiers est obtenue en calculant $\Gamma(D,\OO_D)$ et en en réalisant la décomposition primaire (voir \cite[Alg. 2.4]{bruin} pour les détails).

\subsection{Espaces de Riemann-Roch}\label{subsec:rr}
Soient $k_0$ un corps parfait et $k$ une clôture algébrique de $k_0$. Soient $X_0$ une courbe projective lisse sur $k_0$ et $X=X_0\times_{k_0}k$.
Soit $D\in \Div(X)$ un diviseur stable sous l'action de $\Gal(k|k_0)$. Il existe dans la littérature deux approches pour calculer l'espace de Riemann-Roch $\mathscr{L}(D)\coloneqq \HH^0(X,\OO_X(D))$ associé à $D$. 
Un état de l'art détaillé se trouve dans \cite[§1, State of the art]{aude}. Citons deux approches différentes à ce problème. 

\paragraph{Algorithmes géométriques} D'une part, les algorithmes dits géométriques s'inspirent de la méthode présentée en 1874 par Brill et Noether dans \cite{brill-noether}.  Supposons donné un modèle plan $C_0$ de $X_0$ à singularités ordinaires. Soient $P_1,\dots,P_r$ les points singuliers de $C$ et $m_1,\dots,m_r$ leurs multiplicités respectives. Soient $P_{i,1},\dots,P_{i,m_i}$ les points de $X$ au-dessus de $P_i$. Considérons le diviseur adjoint 
\[ E=\sum_{i=1}^r\sum_{j=1}^{m_i}m_i(m_i-1)P_{i,j}.\] L'algorithme de Brill-Noether calcule les éléments d'une base de $\mathscr{L}(D)$ sous la forme $\frac{f_1}{h},\dots,\frac{f_s}{h}$. Dans un premier temps, il calcule le dénominateur commun $h\in k[x,y]$. Il suffit que le diviseur de $h$ satisfasse \[ \div(h)\geqslant D+E.\] Dans un second temps, il calcule les polynômes $f_i$ de degré $\deg(h)$ ; ils satisfont \[ \div(f_i)\geqslant \div(h)-D.\]
Cette condition se traduit par un système d'équations linéaires en les coefficients des $f_i$. \\

Un algorithme inspiré de celui de Brill-Noether, qui représente les diviseurs par leur forme de Chow et utilise des résultants multivariés pour éviter le recours à la factorisation de polynômes dans des grandes extensions du corps de base, a été mis au point par Huang et Ierardi en 1994 ; c'est celui utilisé dans l'algorithme de la section \ref{sec:huang}. Les algorithmes les plus récents de cette famille sont \cite{aude} et \cite{alain}.
\begin{theorem}\cite[Th. 5.1]{huang_counting} Soit $C$ une courbe projective plane de degré $d$ sur un corps $k_0$. On suppose que tous les points singuliers de $C$ sont ordinaires et définis sur $k_0$. Soit $D\in \Div_{k_0}(C)$. \'{E}crivons $D=D^+-D^-$, avec $D^+$ et $D^-$ effectifs de degré $\leqslant m$. Il existe un algorithme déterministe qui calcule une base du $k$-espace vectoriel $\mathscr{L}(D)$ constituée d'éléments de $k_0(C)$ en $O(m^7d^{14})$ opérations dans $k_0$.
\end{theorem}

\paragraph{Algorithmes arithmétiques}
D'autre part, les algorithmes dits arithmétiques représentent les diviseurs comme des ordres dans le corps de fonctions de la courbe. L'algorithme le plus récent de cette famille, dû à Hess \cite[Algorithm 8.5]{hess}, nécessite d'avoir calculé au préalable la clôture intégrale de deux de ces ordres, ce qui peut s'avérer coûteux.

\begin{rk}\label{rk:dginf} Notons $g$ le genre de $X$, et $P_0\in X(k)$. Le théorème de Riemann-Roch assure que tout diviseur de degré zéro sur $X$ est équivalent à un diviseur de la forme $D-gP_0$, où $D\in\Div(X)$ est effectif. \'{E}tant donné un diviseur $E\in \Div^0(X)$, il est possible de calculer un diviseur effectif $D$ et une fonction $f\in k(X)$ tels que $E=D-gP_0+\div(f)$ : il suffit pour cela de calculer l'espace de Riemann-Roch (nécessairement non vide) associé au diviseur $E+gP_0$.
\end{rk}

\subsection{Diviseurs $\FF_q$-rationnels}

Dans cette section, nous fixons une clôture algébrique $\overline{\FF_q}$ de $\FF_q$ et notons $\mathfrak{G}_0$ le groupe $\Gal(\overline{\FF_q}|\FF_q)$.
\begin{lem} Soient $C_0$ une courbe projective lisse sur $\FF_q$ et $C=C_0\times_{\FF_q}\overline{\FF_q}$. Alors il existe dans $\Pic(C)$ un élément $\mathfrak{G}_0$-invariant de degré 1. 
\begin{proof}
Soit $\sigma\in \mathfrak{G}_0$ l'automorphisme de Frobenius. L'endomorphisme $\sigma-\id$ de la jacobienne $J_C$ est une isogénie de noyau $J(\FF_q)$, il est donc surjectif. Pour un point $P\in C(\overline{\FF_q})$, il existe donc une classe d'un diviseur $D$ de degré 0 tel que $(\sigma-\id)D=(\sigma-\id)P$. L'élément $P-D\in \Pic^1(C)$ est alors Galois-invariant. 
\end{proof}
\end{lem}

\begin{lem}\label{lem:PicDiv} Soient $k_0$ un corps et $k$ une clôture séparable de $k_0$. Soient $C_0$ une courbe projective lisse sur $k_0$ et $C=C_0\times_{k_0}k$. Si $\HH^2(\Gk,k^\times)=0$ alors toute classe $\Gk$-invariante dans $\Pic(C)$ contient un diviseur $\Gk$-invariant.
\begin{proof}
Soient $K_0$ et $K$ les corps des fonctions respectifs de $C_0$ et $C$.
La suite exacte courte \[ 0\to K^\times/k^\times \to \Div C\to \Pic C\to 0\]
donne la suite exacte longue en cohomologie 
\[ \HH^0(\mathfrak{G}_0,\Div C) \to \HH^0(\mathfrak{G}_0,\Pic C)\to \HH^1(\mathfrak{G}_0,K^\times/k^\times). \]
La flèche de droite se décrit de la façon suivante : à la classe d'un diviseur $D$, elle associe le cocycle $\sigma\mapsto f_{D,\sigma}$ où $D^\sigma-D=\div(f_{D,\sigma})$.
Remarquons que $\HH^1(\mathfrak{G}_0,K^\times)=\HH^1(\Gal(K|K_0),K^\times)$ est nul d'après le théorème 90 de Hilbert \cite[III.1, Lem. 1]{serre_cohomologie}. La suite exacte longue associée à \[ 0\to k^\times \to K^\times \to K^\times/k^\times\to 0\]
assure alors que $\HH^1(\Gk,K^\times/k^\times)$ s'injecte dans $\HH^2(\Gk,k^\times)$, qui est nul par hypothèse.
\end{proof}
\end{lem}

\begin{cor} Toute courbe projective lisse sur $\FF_q$ possède un diviseur $\FF_q$-rationnel de degré 1.
\begin{proof}
Le théorème de Wedderburn \cite[11.1, Th. 1]{bourbaki_alg} assure que \[\HH^2(\mathfrak{G}_0,\overline{\FF_q}^\times)=0.\] Le résultat découle alors des deux lemmes précédents.
\end{proof}
\end{cor}

Le lemme \ref{lem:PicDiv} est entièrement effectif, et se réduit à l'association du lemme du serpent avec le théorème 90 de Hilbert. Cependant, lorsque $k_0$ est fini, l'algorithme fourni par cette preuve est de complexité polynomiale en $q$, puisqu'il fait notamment intervenir le groupe des 1-cochaînes d'un certain groupe à valeurs dans $\FF_q^\times$. Les lemmes suivants fournissent une méthode plus adaptée à la pratique, suggérée par Alain Couvreur.

\begin{lem}\label{lem:divdef} Soient $k_0$ un corps et $k$ une clôture séparable de $k_0$. Soit $C_0$ une courbe intègre projective lisse sur $k_0$. Notons $C=C_0\times_{k_0}k$. Soit $D_1$ un diviseur sur $C$ défini sur une extension séparable $k_1$ de $k_0$ telle que $C_0\times_{k_0}k_1$ soit connexe. Soit $D=D_1+\div(f)$, avec $f\in k_1(C)$, un diviseur $k_0$-rationnel équivalent à $D_1$. Soit enfin $b\in k_1$ tel que $\tr_{k_1/k_0}(bf)\neq 0$. Ici, la trace relative d'une fonction rationnelle est définie coefficient par coefficient. Alors \[\tr_{k_1/k_0}\mathscr{L}_{k_1}(D_1)=\tr_{k_1/k_0}(bf)\mathscr{L}_{k_0}(D).\]
\begin{proof} Notons $V=\mathscr{L}_{k_1}(D_1)$.
D'une part, pour tout $u\in \mathscr{L}_{k_0}(D)$, $\tr_{k_1/k_0}(bf)u=\tr_{k_1/k_0}(bf u)$ appartient à $\tr_{k_1/k_0}(V)$, et $\dim_{k_0}\tr_{k_1/k_0}(V)\geqslant \dim_{k_0}\mathscr{L}_{k_0}(D)$. D'autre part, comme $V=bf\mathscr{L}_{k_1}(D)$, l'espace $\tr_{k_1/k_0}(V)$ est inclus dans $\tr_{k_1/k_0}\mathscr{L}_{k_1}(D)=\tr_{k_1/k_0}(bf)\mathscr{L}_{k_0}(D)$.
\end{proof}
\end{lem}

\begin{lem}\label{lem:divgd} Soit $C$ une courbe intègre projective lisse de genre $g$ sur un corps algébriquement clos $k$. Soit $D\in\Div(C)$ un diviseur de degré $\geqslant 2g$. Notons $V=\mathscr{L}(D)$. Alors \[D=-\sum_{P\in |C|}\min_{f\in V}v_p(f)P.\]
\begin{proof} Notons $D'=-\sum_{P\in C}\min_{f\in V}v_p(f)P$. D'une part, $D\geqslant D'$ car pour tout point $P$, $v_P(D-D')=v_P(D)+\min_{f\in V} v_P(f)\geqslant 0$. D'autre part, $V\subseteq\mathscr{L}(D')$ car pour tout $f\in V$ et tout point fermé $P\in C$, $v_P(\div(f)+D')=v_P(f)-\min_g v_P(g) \geqslant 0$. Par conséquent, $V=\mathscr{L}(D')$. Supposons que $D\neq D'$ ; il existe alors des points $P_1,\dots,P_r$ tels que $D=D'+\sum_i P_i$. Or pour tout $i$, le théorème de Riemann-Roch assure que \[\mathscr{L}(D-P_i)\neq \mathscr{L}(D)\]
car $\deg(D-P_i)>2g-2$. Comme $\mathscr{L}(D')\subseteq \mathscr{L}(D-P_i)$, c'est absurde.
\end{proof}
\end{lem}

\begin{prop}\label{prop:divrat} Soit $X_0$ une courbe intègre lisse sur $k_0=\FF_q$. Notons $X=X_0\times_{k_0}k$. Supposons donné un modèle plan de $X$ de degré $d$, ainsi qu'un diviseur $k_0$-rationnel $E$ de degré 1 sur $X$. Soient $k_1$ une extension de $\FF_q$ telle que $X_0\times_{k_0}k_1$ soit connexe et $D_1\in \Div(X)$ un diviseur défini sur $k_1$ dont la classe dans $\Pic X$ est $\Gal(k|k_0)$-invariante. Notons $D_1=D_1^+-D_1^-$ et $E=E^+-E^-$ avec $D_1^+,D_1^-, E^+,E^-$ effectifs de degré inférieur à un entier $c$. Il existe un algorithme de complexité polynomiale en $\log q$, $g$, $c$, $d$ et $[k_1:k_0]$ qui détermine un diviseur $k_0$-rationnel $D$ équivalent à $D_1$. 
\begin{proof} 
Nous savons qu'il existe un tel diviseur $D$. Soit $h$ une fonction telle que $D_1=D+\div(g)$. Par séparabilité de l'extension $k_1/k_0$, il existe $b\in k_1$ tel que $\tr_{k_1/k_0}(bh)\neq 0$.
Soit $N$ un entier tel que $\deg(D_1+NE)\geqslant 2g$. Appliquons la trace $\tr_{k_1/k_0}$ à l'espace de Riemann-Roch $\mathscr{L}_{k_1}(D_1+NE)$. L'espace obtenu est $V\coloneqq \tr(bh)\mathscr{L}_{k_0}(D+NE)$ par le lemme \ref{lem:divdef}. 
Le diviseur $D+NE$ est égal à $-\sum_{P\in X}\min_{f\in V}v_P(f)P$ par le lemme \ref{lem:divgd}. Il suffit donc de calculer une base de $V$, puis le diviseur de chacun des éléments de cette base, pour obtenir le diviseur $k_0$-rationnel $D+NE+\div(\tr bh)$ équivalent à $D_1+NE$. Pour le résultat de complexité, remarquons que $N\leqslant g+c$.
\end{proof}
\end{prop}

\subsection{Diviseur équivalent de support évitant un fermé}\label{subsec:moving}

Soit $X$ une courbe projective lisse sur un corps $k$. Soit $C$ un modèle birationnel plan de 
$X$ de degré $d$. Soient $D$ un diviseur sur $X$, et $Z$ un fermé de $X$. Il est toujours possible de calculer un diviseur $D'$ sur $X$ équivalent à $D$ et de support disjoint de $Z$ ; cette méthode est décrite par exemple dans \cite[§3.4]{couveignes_linearizing}. Soit $O$ un diviseur sur $X$ de degré inférieur à $d$ et de support disjoint de $Z$. Il suffit de calculer l'espace de Riemann-Roch $\mathscr{L}$ du diviseur $D+2gO$, qui est de dimension strictement supérieure à 1. Dans $\mathscr{L}$, les fonctions $f$ telles que $\div(f)+D$ ne soit pas disjoint de $Z$ sont contenues dans une réunion d'hyperplans, qu'il suffit d'éviter.

\begin{prop}\cite[Lem. 5]{couveignes_linearizing} Il existe un algorithme déterministe qui, étant donné une courbe plane $C$ sur $\FF_q$ et sa normalisée $X$, un diviseur $\FF_q$-rationnel $D=D^+-D^-$ de degré 0 et un diviseur effectif $Z$, calcule un diviseur $E=E^+-E^-$ linéairement équivalent à $D$ et de support disjoint de $Z$ en un nombre d'opérations dans $\FF_q$ polynomial en $d$, $\log q$, $\deg(D^+)$ et $\deg(Z)$. Le degré des diviseurs $E^+$ et $E^-$ est inférieur à $6gd(\log_q(\deg Z)+1)$.
\end{prop}

\section{Fonction zêta et comptage de points}

Soit $q=p^a$ une puissance d'un nombre premier. Soit $X$ une courbe intègre projective lisse de genre $g$ sur $\FF_q$. Il existe différentes méthodes pour compter le nombre de $\FF_q$-points ou la fonction zêta de $X$. L'une d'entre elles est le calcul de la cohomologie étale modulo $\ell$, présenté dans ce manuscrit, pour $O(g\log q)$ valeurs de $\ell$.  L'avantage des méthodes $\ell$-adiques est de proposer une complexité polynomiale en $\log q$ ; elles sont cependant exponentielles en $g$. Pour une courbe quelconque, la meilleure de ces méthodes est celle de Huang et Ierardi présentée dans la section \ref{sec:huang}. Pour les courbes hyperelliptiques, il existe des algorithmes plus efficaces, comme celui d'Adleman et Huang. \`{A} l'aide de la représentation de Mumford des diviseurs, il calcule directement les points de $\ell$-torsion de la jacobienne de la courbe en question.

\begin{prop}\cite[Th. 4.1]{adleman_huang} Soit $q$ une puissance d'un nombre premier impair. Il existe un algorithme déterministe qui, étant donné un polynôme $f\in \FF_q[t]$ de degré $2g+1$ sans racine multiple, renvoie le polynôme caractéristique de l'endomorphisme de Frobenius sur la jacobienne de la courbe hyperelliptique d'équation affine $y^2=f(x)$ en $(\log q)^{O(g^2\log g)}$ opérations.
\end{prop}

Il existe également des méthodes $p$-adiques, comme celles décrites dans \cite{wan_zeta} ou \cite{tuitman}. La complexité de la plupart de ces algorithmes est polynomiale en $p$ et en $g$. Cependant, Harvey a présenté en 2014 une méthode dont la complexité moyenne (sur $p$) est polynomiale en $\log p$.

\begin{theorem}\cite[Th. 1]{harvey} Il existe un algorithme déterministe explicite vérifiant les propriétés suivantes. \'{E}tant donné des entiers $N\geqslant 3$, $g\geqslant 1$, et un polynôme $Q\in \ZZ[x]$ définissant une courbe hyperelliptique $X$ sur $\QQ$ d'équation $y^2=Q(x)$ de genre $g$, il calcule pour tous les premiers impairs $p<N$ ne divisant pas le discriminant de $Q$ la fonction zêta de la réduction de $X$ modulo $p$. L'algorithme nécessite $g^{8+\varepsilon}N\log^2(N)\log^{1+\varepsilon}(\Vert Q\Vert_\infty N)$ opérations binaires, où $\Vert Q\Vert_\infty$ désigne le maximum des valeurs absolues des coefficients de $Q$.
\end{theorem}
Comme le nombre de nombres premiers $p<N$ est équivalent à $N/\log(N)$, le temps moyen pour chaque premier $p$ est $g^{8+\varepsilon}\log^3(N)\log^{1+\varepsilon}(\Vert Q\Vert_\infty N)$. Plus récemment, Harvey et Sutherland ont proposé et implémenté une méthode pour calculer la matrice de Hasse-Witt modulo $p$ d'une telle courbe hyperelliptique $X$, et en déduire la réduction modulo $p$ du polynôme caractéristique du Frobenius sur $X$ \cite{harvey_sutherland} ; leur algorithme est particulièrement efficace dans la pratique pour les courbes de genre 2 ou 3.

\cleartooddpage
\addcontentsline{toc}{chapter}{Bibliographie}
\printbibliography

\end{document}